\documentclass[12pt]{amsart}
\usepackage[margin=1in]{geometry}
\usepackage{amsmath}
\allowdisplaybreaks
\usepackage{mathtools}
\usepackage{amssymb}
\usepackage{amsthm}
\usepackage{enumitem}
\usepackage{mathrsfs}
\usepackage[pagebackref,colorlinks=true,linkcolor=blue,citecolor=blue]{hyperref}
\usepackage{tikz-cd}
\usepackage{tikz-3dplot}
\usepackage{xcolor}
\usepackage{pgffor}
\usepackage{ifthen}
\usepackage{asymptote}
\usepackage[margin=1in]{geometry}
\usepackage{longtable}
\usepackage[english]{babel}
\usepackage{graphicx}
\usepackage{caption}
\usepackage{float}

\usepackage{bbding} 

\tikzset{
  symbol/.style={
    draw=none,
    every to/.append style={
      edge node={node [sloped, allow upside down, auto=false]{$#1$}}}
  }
} 
\usetikzlibrary{patterns}

\newcommand{\Z}{\mathbb{Z}}

\newcommand{\R}{\mathbb{R}}
\newcommand{\BC}{\mathbb{C}}

\newcommand{\ol}{\overline}

\newcommand{\SL}{\mathrm{SL}}
\newcommand{\GL}{\mathrm{GL}}
\newcommand{\SO}{\mathrm{SO}}
\newcommand{\Sp}{\mathrm{Sp}}

\newcommand{\RG}{\mathrm{G}}

\newcommand{\Irr}{\mathrm{Irr}}

\newcommand{\Rep}{\underline{\mathrm{Rep}}}

\newcommand{\half}[1]{\frac{#1}{2}}

\newcommand{\Temp}[3]{T_{#1,#3}^{#2}}

\newcommand{\comment}[1]{}
\newcommand{\EE}{\mathcal{E}}
\newcommand{\FF}{\mathcal{F}}
\newcommand{\shrt}{\text{sla}}
\newcommand{\reg}{\text{reg}}
\newcommand{\lev}{\text{lev}}
\newcommand{\temp}{\text{temp}}

\newtheorem{thm}{Theorem}[section]

\newtheorem{lemma}[thm]{Lemma}
\newtheorem{prop}[thm]{Proposition}
\newtheorem {conj}[thm]{Conjecture}

\newtheorem {ques/conj}[thm]{Question/Conjecture}

\newtheorem{defn}[thm]{Definition}
\newtheorem{defn'}[thm]{Definition$'$}
\newtheorem{remark}[thm]{Remark}

\newtheorem{exmp}[thm]{Example}

\newtheorem{algo}[thm]{Algorithm}

\newtheorem*{globalcond*}{Global Condition}
\newtheorem*{localcond*}{Local Condition}

\newtheorem*{globalconj*}{Global Conjecture}
\newtheorem*{localconj*}{Local Conjecture}

\newtheorem*{nonzero*}{Conjecture on the non-vanishing of the normalized intertwining operators}
\newtheorem*{holo*}{Conjecture on the holomorphicity of the normalized intertwining operators}

\DeclareMathOperator{\supp}{supp}

\DeclareMathOperator{\Ind}{Ind}
\numberwithin{equation}{section}

\let\oldbullet\bullet
\renewcommand{\bullet}{{\vcenter{\hbox{\tiny$\oldbullet$}}}}

\begin{document}

\title[unitary representations of corank 4]{On corank 4 unitary representations of classical groups}

\author[Baiying Liu]{Baiying Liu}
\address{Department of Mathematics\\
Purdue University\\
West Lafayette, IN, 47907, USA}
\email{liu2053@purdue.edu}

\author[Chi-Heng Lo]{Chi-Heng Lo}
\address{Department of Mathematics\\
National University of Singapore\\
119076, Singapore}
\email{{ch\_lo@nus.edu.sg}}

\author[Brian Wen]{Brian Wen}
\address{Department of Mathematics\\
Purdue University\\
West Lafayette, IN, 47907, USA}
\email{wen190@purdue.edu}

%    General info

\date{\today}


\thanks{The research of the first-named author is partially supported by the NSF Grant DMS-1848058 and the Simons Foundation: Travel Support for Mathematicians.}

\begin{abstract}
    In this paper, we explicitly classify the corank 4 unitary representations of symplectic or split odd special orthogonal groups over non-Archimedean local fields of characteristic zero, by classifying Arthur representations of corank 4 and verifying the corresponding unitary dual conjecture recently proposed by Hazeltine-Jiang-Liu-Lo-Zhang in \cite{HJLLZ24}. 
\end{abstract}

\maketitle

\tableofcontents

\section{Introduction}

The unitary dual problem asks for the complete classification of all irreducible unitary representations of a given locally compact group $G$. Its origins trace back to the early twentieth century with the development of harmonic analysis on groups. For abelian connected Lie groups, Pontryagin duality provided a clean answer -- the unitary dual of an abelian group is itself another group (the Pontryagin dual), and Fourier analysis is simply integration over this dual. For compact groups, the Peter–Weyl theorem gave a full decomposition of $L^2(G)$ into finite-dimensional irreducibles. The challenge intensified with non-compact and non-abelian groups, especially real and $p$-adic reductive Lie groups.

The unitary dual captures all possible symmetries that act unitarily, so it is essential to harmonic analysis, number theory (via automorphic forms and the Langlands program), and quantum physics (where symmetries are represented unitarily). In the Langlands framework, the structure of the unitary dual is deeply related to spectral decompositions $L^2$-spaces on arithmetic quotients and thus to major conjectures linking representation theory and arithmetic geometry. Despite major advances, a full description of the unitary dual is still only known for certain classes of groups; for general reductive groups, the classification remains highly nontrivial and continues to drive research in modern representation theory.

In \cite{HJLLZ24}, Hazeltine, Jiang, the first two named authors, and Zhang proposed a new conjecture describing the structure of the unitary dual in terms of Arthur representations for connected reductive algebraic groups defined over any non-Archimedean local field of characteristic zero. This conjecture provides a candidate set for the unitary dual, constructed from Arthur representations. 
For classical groups, they developed an explicit algorithm to generate this candidate set. Evidence for its exhaustiveness includes compatibility with the known generic unitary dual, unramified unitary dual, and low-corank representations. As further support of the conjecture, they verified the conjecture for the unitary dual of the exceptional group of type $\RG_2$.

More precisely, let $F$ be a non-Archimedean local field of characteristic zero. The celebrated Langlands classification determines the admissible dual $\Pi(G)$ by means of the tempered dual $\Pi_\temp(M)$ of Levi subgroups $M$ of $G$, which is the subset of $\Pi(M)$ consisting of tempered members. However, there is no general approach available to determine the unitary dual $\Pi_u(G)$ based on the Langlands classification of the admissible dual $\Pi(G)$. On the other hand, the far-reaching  
endoscopic classification conjecture of J. Arthur (\cite{Art89}) produces local Arthur packets (Definition \ref{A parameters}), which are expected to be finite subsets of the unitary dual $\Pi_u(G)$. Denote by $\Pi_A(G)$ the subset of $\Pi(G)$ consisting of Arthur representations, also referred to as representations of Arthur type, which is the union of all local Arthur packets. 
It is strongly desirable that the whole unitary dual $\Pi_u(G)$ can be exhausted via various expected constructions from the Arthur representations $\Pi_A(G)$. 

In \cite{HJLLZ24}, the authors defined a set $\Pi_{\overline{A}}^{\lim}(G)$ from Arthur representations $\Pi_A(G)$, which is constructed from $\Pi_A(G)$ by applying the standard constructions of unitary representations: unitary parabolic induction, complementary series, and limits of complementary series, as explained in Definition \ref{def closure A}. Then, they conjectured that this set is exactly the unitary dual. 
Loosely speaking, the idea is that \emph{complementary series representations can be obtained by continuous Hermitian deformation of unitary representations irreducibly induced from unitary data} (see \cite[Remark 5.3(2)]{HJLLZ24} for more precise statements).
Similar ideas have been clearly reflected in the construction of the unitary dual for general linear groups (see \cite{Tad86, Vog86} and \cite[Introduction]{HJLLZ24}). The conjecture in \cite{HJLLZ24} is as follows. 

\begin{conj}[{Unitary Dual Conjecture, \cite[Conjecture 1.1]{HJLLZ24}, see also Conjecture \ref{unitary dual conjecture}}]\label{conj A+ intro}
Let $G$ be a general connected reductive group defined over $F$. Assume that the local Arthur conjecture as in \cite[Conjecture 6.1]{Art89} holds for every Levi subgroup $M$ of $G$. Then the following two sets are equal:
    $$\Pi_{\overline{A}}^{\lim}(G) = \Pi_{u}(G).$$
\end{conj}

The significance of Conjecture \ref{conj A+ intro} may be described as follows.
The unitary dual $\Pi_u(G)$ has a natural and conceptually simple definition in representation theory, yet its classification and internal structure remain largely mysterious. In contrast, according to the conjectures of J. Arthur (\cite{Art89, Art13}), the set of Arthur representations $\Pi_A(G)$ is defined via local stability and (twisted) endoscopic transfer, yielding a rich and highly structured framework; however, its intrinsic representation-theoretic meaning is not fully understood.

Conjecture \ref{conj A+ intro} serves as a bridge between these two domains by proposing a candidate subset $\Pi_{\overline{A}}^{\lim}(G)$ for the unitary dual, constructed from $\Pi_A(G)$ and designed to retain structural transparency while remaining amenable to explicit computation.  
A principal motivation of \cite{HJLLZ24} is precisely that $\Pi_{\overline{A}}^{\lim}(G)$ should be computable, thereby opening the possibility of an algorithmic determination of the unitary dual itself. Supporting this perspective, the authors provide an explicit procedure (\cite[Algorithm 8.5]{HJLLZ24}) that produces $\Pi_{\overline{A}}^{\lim}(G)$ for symplectic groups and split odd special orthogonal groups.

On the other hand, in \cite{Tad23}, Tadi{\'c} constructed the full unitary dual of classical groups for representations of corank up to 3 (see Definition \ref{def corank}). In this paper, we combine the techniques in both \cite{Tad23} and \cite{HJLLZ24}, to construct the corank 4 unitary dual of  symplectic and split odd special orthogonal groups and verify  the corresponding Conjecture \ref{conj A+ intro} (see Theorem \ref{finalconclusion}). 
In particular, we completes the classification of the unitary dual for $\Sp_8$ and split $\SO_9$.  

\begin{thm}\
\begin{enumerate}
\item Theorem \ref{finalconclusion} provides the complete list of corank 4 unitary representations of  symplectic and split odd special orthogonal groups. 
\item (Theorem \ref{main conjecture holds}) Conjecture \ref{conj A+ intro} holds for corank 4 unitary representations of symplectic and split odd special orthogonal groups. 
    \item  For $\mathrm{G}=\Sp_8$ or split $\SO_9$, Conjecture \ref{conj A+ intro} holds, i.e.,
    $\Pi_{\overline{A}}^{\lim}(G) = \Pi_{u}(G).$
\end{enumerate}
\end{thm}

As other applications, in the setting of corank 4 representations, the main result in this paper verifies some conjectures of Tadi{\'c} (\cite[Conjecture $1.1$]{Tad22}(2)(3) and \cite[Conjecture 8.15]{Tad23}) on the isolated unitary representations.

\begin{thm} 
Let $G_n$ be a symplectic or split odd special orthogonal group. Then
all corank 4 isolated unitary representations of $G_n$ are of critial (see Definition \ref{def critical type}) and Arthur type, hence are automorphic. 
\end{thm}

The above theorem is directly implied by the corank 4 unitary dual list (Theorem \ref{finalconclusion}), the corank 4 critical Arthur type list (Proposition \ref{crnk4Artcritlist}), the corank 4 critical non-Arthur type list (Proposition \ref{subquotientlist}), and Remark \ref{isolated reps}.  
Note that by the work of Arthur (\cite{Art13}),  Arthur type representations are automorphic. 

Our main methods may be summarized as follows. First, by \cite[Theorem 4.2]{Tad09a}, the
classification of the unitary dual for classical groups can be reduced to the weakly real case (see Definition \ref{def weakly real}), which is the
focus of this paper. Then, 
we compute the full good parity Arthur dual of corank $4$  using the classification of tempered representations of good parity, as well as the techniques in \cite{HJLLZ24} to determine Arthur type representations. By \cite[Theorem $1.1$]{AM25}, this is exactly equal to the unitary dual for good parity representations. Using this and the list of all possible subquotients at the critical points of corank $4$, we can determine exactly which critical point contains a non-unitarizable subquotient. Then, we use a similar reduction technique as the one used in \cite{Tad23}, in conjunction with the reduction techniques listed in \cite[Algorithm $8.5$]{HJLLZ24}, to classify all possible unitary connected components in corank $4$. This gives us all $4$ and $3$-dimensional unitary connected components for corank $4$ representations. Finally, we use an exhaustive process to compute all possible 
unitary complementary series of dimensions $1$ and $2$ to give the full corank 4 unitary dual, which we will prove to be equal to the conjectured set $\Pi_{\overline{A}}^{\lim}(G_n)$. 
We remark that our results imply that the following principles might be true for general connected reductive groups: 
\begin{enumerate}
    \item Complementary series can be obtained via
irreducible continuous Hermitian deformations from irreducible inductions of unitary data. 
\item A representation is not unitarizable if and only if it lies inside a continuous family of representations which contains an irreducible good parity subquotient that is not of Arthur type in its interior or closure.
\end{enumerate}

Here is the structure of the paper. In \S \ref{notations}, we introduce the basic notation and give some preliminary results needed in later sections, leading up to the definition of the sets $\Pi_{\overline{A}}(G_n)$ and $\Pi_{\overline{A}}^{\lim}(G_n)$. In sections \S \ref{classtempcorank3} to \S \ref{classtempcorank4}, we construct the Arthur dual for representations of corank $4$ that are of good parity. In particular, in \S \ref{classtempcorank3}, we begin by classifying good parity tempered representations of corank $3$. Then, using these results, we  classify all good parity non-tempered representations of corank $4$ in \S\ref{classnontempcorank4,1} to \S\ref{classnontempcorank4,34}, and identify which one of them are of Arthur type and of critical type. Subsequently, in \S\ref{classtempcorank4}, we classify all good-parity tempered representations of corank $4$, which gives the full Arthur dual of corank $4$ in the good parity case. We summarize our results in Appendix \ref{artlist} and \ref{nonartlist}, by giving respectively the full list of representations that are of critical type and of Arthur type, as well as the complementary list of representations of critical type, but not of Arthur type. Using these two lists, we give the full list of unitary open components, as well as their boundaries, for representations of corank $4$ in \S \ref{opencnncomponents}. To construct the full unitary dual, we append to it the unitarizable representations inside a one or two-parameter complementary series, which are constructed in \S \ref{1-parameterseries} and 
\S\ref{2-parameterseries} respectively. Finally, in \S\ref{conclusion}, we give the full unitary dual of corank $4$ and prove that it is the indeed the same as the closure of the Arthur dual $\Pi_{\overline{A}}^{\lim}(G_n)$, hence proving Conjecture \ref{conj A+ intro}.

\subsection*{Acknowledgements} 
The authors would like to thank Wee Teck Gan, Dihua Jiang, and Freydoon Shahidi for the interest and encouragement. The authors also would like to thank Alexander Hazeltine, Marko Tadi{\'c}, David Vogan, and Qing Zhang for helpful communication and suggestions. Part of this paper was written when the authors were attending the Relative Langlands Program workshop and conference at the Institute for Mathematical Sciences, the National University of Singapore (January 2026). The authors appreciate very much the warm hospitality of the Institute and the mathematics department of NUS.

\section{Notations and preliminary results} \label{notations}
Throughout the rest of this paper, let $F$ denote a non-Archimedean local field of characteristic $0$, otherwise known as a $p$-adic local field. Let $\lvert \cdot \rvert$ be the character of $\GL_{n}(F)$ obtained by composing the determinant function of $\GL_{n}(F)$ with the normalized $p$-adic absolute value of $F$. 

Fix any positive integer $n$, we consider the representations of the general linear group $\GL_{n}(F)$, or the symplectic (resp. split odd special orthogonal) groups, denoted by $G_n = $ $\mathrm{G}_n(F)$, where $\mathrm{G}_n = \Sp_{2n}$ (resp. $\SO_{2n+1}$).

\subsection{Local Langlands classification}

In this subsection, we recall the explicit form of the Langlands classification for both the general linear group and classical groups over $p$-adic local fields (\cite{Sil78,BW00,Kon03}). 

In the case of general linear groups $\GL_{n}(F)$, fix a Borel subgroup $B$  and let $P = MN$ be the standard parabolic subgroup of $\GL_n(F)$ with Levi subgroup $M \cong \GL_{n_1}(F) \times \dots \times \GL_{n_r}(F)$. Given smooth representations $\tau_i$ of $\GL_{n_i}(F)$ for $i = 1, 2, \ldots r$,  we denote the normalized parabolic induction by: 
\[\tau_1 \times \dots \times \tau_r  := \text{Ind}_{P}^{\GL_n(F)}(\tau_1 \otimes \dots \otimes \tau_r ).\]

Throughout the rest of this paper, let $\rho$ denote an irreducible supercuspidal representation of $\GL_{n}(F)$. A (Zelevinsky) segment $[x,y]_{\rho}$ is the set of supercuspidal representations of the form 
\[[x,y]_{\rho} = \{\rho\lvert \cdot \rvert^{x}, \rho\lvert \cdot \rvert^{x-1}, \ldots, \rho\lvert \cdot \rvert^{y}\}, \]
where $x,y \in \mathbb{R}$ and $x-y$ is a non-negative integer. The unique irreducible subrepresentation of $\rho\lvert \cdot \rvert^{x} \times \rho\lvert \cdot \rvert^{x-1} \times \ldots \rho\lvert \cdot \rvert^{y}$ is called the (essentially) Steinberg representation attached to the segment $[x,y]_{\rho}$ and is denoted by $\Delta_\rho[x,y]_{\rho}$. For simplicity, we treat $\Delta_{\rho}[x,y]$ as the trivial representation of $\GL_{0}(F)$ when $x < y$.

     Langlands classification provides a way to classify all equivalence classes of smooth admissible representations of $\GL_{n}(F)$, or the admissible dual $\Pi(\GL_{n}(F))$. Specifically, it states that any 
     irreducible admissible representation $\tau$ of $\GL_{n}(F)$ can be realized as the unique irreducible subrepresentation of some
     \[\Delta_{\rho_{1}} [x_1, y_1] \times \dots \times \Delta_{\rho_{r}}[x_r, y_r],\]
     where the $\rho_{i}$'s are irreducible unitary supercuspidal representations of $\GL_{n}(F)$, $[x_i, y_i]_{\rho_i}$ are segments, and $x_1 + y_1 \leq \dots \leq x_r + y_r$. In this case, we write
     \[\tau = L(\Delta_{\rho_1}[x_1, y_1], \ldots \Delta_{\rho_r}[x_r, y_r]).\]

     Let $(x_{i,j})_{1\leq i\leq s, 1\leq j \leq t}$ be real numbers such that $x_{i,j}=x_{1,1}-i+j.$ A \emph{generalized Speh representation} is an irreducible representation of  the form
\begin{equation}\label{generalized Speh representation}
\begin{pmatrix}
x_{1,1} & \cdots & x_{1,t} \\
\vdots & \ddots & \vdots \\
x_{s,1} & \cdots & x_{s,t}
\end{pmatrix}_{\rho}:=L(\Delta_{\rho}[x_{1,1},x_{s,1}],\dots,\Delta_{\rho}[x_{1,t}, x_{s,t}]).
\end{equation}
These representations will be useful in the classification of unitary dual for $\GL_{n}(F)$.

The Langlands classification of the classical groups $G_n$ can be given similarly as follows. Let $P = MN$ be a standard parabolic subgroup of $G_n$ with Levi subgroup $M \cong \GL_{n_1}(F) \times \dots \times \GL_{n_r}(F) \times G_m$ ($G_m$ a group of the same type as $G_n$ with $m \leq n$). Given smooth representations $\tau_i$ of $\GL_{n_i}(F)$ for $i = 1, 2, \ldots r$, and a smooth representation $\sigma$ of $G_m$ $(m \leq n)$, we denote the normalized parabolic induction by: 
\[\tau_1 \times \dots \tau_r \rtimes \sigma := \text{Ind}_{P}^{G_n}(\tau_1 \otimes \dots \otimes \tau_r \otimes \sigma).\]
     Then, Langlands classification for $G_n$ states that every irreducible representation $\pi$ of $G_n$ is a unique irreducible subrepresentation of 
    \[\Delta_{\rho_1}[x_1, y_1] \times \dots \times \Delta_{\rho_r}[x_r,y_r] \rtimes \pi_{temp},\]
    where the $\rho_{i}$'s are irreducible unitary supercuspidal representations of $\GL_{n}(F)$, $x_1 + y_1 \leq \ldots \leq x_r + y_r < 0$, and $\pi_{temp}$ is tempered. In this case, we write
    \[\pi = L(\Delta_{\rho_1}[x_1, y_1], \ldots \Delta_{\rho_r}[x_r, y_r]; \pi_{temp}).\]
    We call the tuple $(\Delta_{\rho_{1}}[x_1, y_1], \ldots, \Delta_{\rho_{r}}[x_r, y_r], \pi_{temp})$ the Langlands data, or $L$-data, of $\pi$, and the multi-set $\{\Delta_{\rho_{1}}[x_1, y_1], \ldots, \Delta_{\rho_{r}}[x_r, y_r]\}$ the non-tempered portion of the $L$-data of $\pi$.  

\subsection{Arthur's parameterization of the tempered spectrum}

In this subsection, for classical groups $G_n$, we recall the local Langlands correspondence (\cite[\S8.2]{Bor79} and the Arthur's parametrization of the tempered spectrum (\cite[Theorem 1.5.1]{Art13}). 

For any split reductive algebraic group $\RG$, the Langlands dual group of $\RG$ is denoted by $\RG^\vee(\BC)$. 
A local $L$-parameter of $G_n$ is a 
${\RG}^\vee_n(\BC)$-conjugacy class of an admissible homomorphism $\phi:W_F\times\SL_2(\BC)\rightarrow{\RG}^\vee_n(\BC)$ (more generally, see \cite[\S8.2]{Bor79}). In this paper, we do not distinguish $\phi$ from its conjugacy class. We say that $\phi$ is tempered if its restriction to $W_F$ has bounded image. The component group of $\phi$ is defined by
\[
\mathcal{S}_{\phi}:= \pi_0( \mathrm{Cent}_{{\RG}^\vee_n(\BC)}(\mathrm{Im}(\phi))/Z({\RG}^\vee_n(\BC))^{\Gamma}),
\]
where $\Gamma$ is the absolute Galois group of $F$.

By the local Langlands correspondence for $\GL_n(F)$ (\cite{Hen00, HT01, Sch13}), we may identify an irreducible supercuspidal representation $\rho$ of $\GL_n(F)$ with its local $L$-parameter $\phi_\rho$, which is an irreducible 
$n$-dimensional representation of the Weil group $W_F$. Explcitly, we can write 
any tempered local $L$-parameter $\phi$ in the following form:
\[
\phi=\bigoplus_{i=1}^r \rho_i\otimes S_{a_i},
\]
where $\rho_i$ are irreducible unitary supercuspidal representations of $\GL_{n_i}(F)$ 
and $S_k$ denotes the unique irreducible representation of $\SL_2(\BC)$ of dimension $k$. Let $\mathrm{Jord}(\phi)$ 
denote the multi-set consisting of all the irreducible summands occurring in $\phi$ (counting multiplicities), i.e. for $\phi$ as above, we have $\mathrm{Jord}(\phi)=\{\rho_1\otimes S_{a_1},\dots,\rho_r\otimes S_{a_r}\}.$

For the groups $G_n$, the Pontrayagin dual $\widehat{\mathcal{S}}_\phi$ of $\mathcal{S}_{\phi}$ is a finite abelian group consisting of characters which may be identified with functions $\varepsilon:\mathrm{Jord}(\phi)\rightarrow\{\pm 1\}$ such that
$\varepsilon(\rho_i\otimes S_{a_i})=\varepsilon(\rho_j\otimes S_{a_j})$ whenever $\rho_i\otimes S_{a_i}\cong\rho_j\otimes S_{a_j}$, and 
\[
\prod_{\rho\otimes S_a\in \mathrm{Jord}(\phi)} \varepsilon(\rho\otimes S_a)^{m_{\rho,a}}=1,
\]
where $m_{\rho,a}$ denotes the multiplicity of $\rho\otimes S_a$ in $\phi.$ 

From the conjectural local Langlands correspondence, we associate to each local $L$-parameter $\phi$ a finite set of irreducible admissible representation of $G$ satisfying certain properties (see \cite{Bor79}), which is called the local $L$-packet attached to $\phi$ and is denoted by $\Pi_\phi$. For $G_n$, the following theorem of Arthur shows that  $\Pi_{\phi}$ is in bijection with $\widehat{\mathcal{S}}_{\phi}$ for tempered parameters $\phi$. 

\begin{thm}[{\cite[Theorem 1.5.1]{Art13}}]\label{thm Arthur tempered}
Fix a choice of Whittaker datum for $G_n$ and let $\phi$ be a tempered local $L$-parameter. Then there is a bijective map between the tempered local $L$-packet $\Pi_{\phi}$ and $\widehat{\mathcal{S}}_\phi.$
\end{thm}
Henceforth, we fix a choice of Whittaker datum for $G_n.$ When $\phi$ is tempered, we write $\pi(\phi,\varepsilon)$ for the element of $\Pi_\phi$ corresponding to $\varepsilon\in\widehat{\mathcal{S}}_\phi$ via the bijection in Theorem \ref{thm Arthur tempered}.

\subsection{Supercuspidal representations and reducibility points}

 In this subsection, we recall M{\oe}glin's characterization of local $L$-parameters for irreducible supercuspidal representations of 
$G_n$. Let $\mathcal{C}$ (resp. $\mathcal{C}_{cl}$) be the set of supercuspidal representations of general linear groups (resp. classical groups), 
$\mathcal{C}^{u}:= \{\rho \in \mathcal{C}\ | \ \rho \text{ is unitary}\}$, and $
    \mathcal{C}^{sd}:= \{\rho \in \mathcal{C}^u\ | \ \rho \text{ is self-dual}\}.$

For any tempered $L$-parameter $ \phi= \bigoplus_{i} \rho_i \otimes S_{a_i}$, we write $\rho \otimes S_{a} \subset \phi$ if $\rho \otimes S_a$ appears as a direct summand in  $\phi$. A tempered $L$-parameter $ \phi= \bigoplus_{i} \rho_i \otimes S_{a_i}$ is called \emph{discrete} if the $\rho_i \otimes S_{a_i}$'s are pairwise non-equivalent. 
A discrete $L$-parameter $ \phi= \bigoplus_{i} \rho_i \otimes S_{a_i}$ is said to be \emph{without gaps} if for $a \geq 1$,
\[ \rho \otimes S_{a+2} \subset \phi \Rightarrow \rho \otimes S_{a} \subset \phi. \]
Applying Theorem \ref{thm Arthur tempered}, M{\oe}glin's parametrization of irreducible supercuspidal representations of $G_n$ is as follows.

\begin{thm}[{\cite[Theorem 1.5.1]{Moe11a}, \cite[Theorem 3.3]{Xu17a}}] \label{thm characterizatioin of supercuspidal}
An irreducible tempered representation $\pi(\phi,\varepsilon)$ of $G_n$ is supercuspidal if and only if the following hold.
\begin{enumerate}
    \item [$\oldbullet$] As a tempered local $L$-parameter, $\phi$ is discrete and without gaps.
    \item [$\oldbullet$] If both $\rho \otimes S_{a}\subset \phi$ and $\rho \otimes S_{a+2} \subset\phi$, then $\varepsilon(\rho \otimes S_{a}) \varepsilon(\rho \otimes S_{a+2})=-1.$
    \item [$\oldbullet$] If $\rho \otimes S_{2}\subset\phi$, then $ \varepsilon(\rho \otimes S_2)=-1$.
\end{enumerate}
\end{thm}

Let $\rho\in\mathcal{C}^u$ and $\sigma\in\mathcal{C}_{cl}$. If $\rho$ is not self-dual, then $\rho\vert\cdot\rvert^x\rtimes\sigma$ is irreducible for all $x\in\mathbb{R}.$ If $\rho$ is self-dual, then there exists a unique $\alpha_{\rho,\sigma}\in\mathbb{R}_{\geq 0}$, such that $\rho\vert\cdot\rvert^{\alpha_{\rho,\sigma}}\rtimes\sigma$ is reducible (\cite{Sil80}). The number $\alpha_{\rho, \sigma}$ is known as the reducibility point, and it is known that  $\alpha_{\rho,\sigma}\in\frac{1}{2}\mathbb{Z}_{\geq 0}$.
In fact, there is an explicit description of $\alpha_{\rho,\sigma}$ based on the local $L$-parameter $\phi_{\sigma}$ of $\sigma$ as follows, according to \cite[Remark 4.5.2]{MW06}. 

By Theorem \ref{thm characterizatioin of supercuspidal}, we may write
\begin{align}\label{eq decomp phi_sc}
    \phi_{\sigma}= \bigoplus_{\rho \in R}  \bigoplus_{i=0}^{a_{\rho}} \rho \otimes S_{2(i+\epsilon_{\rho})+1},
\end{align}
where $R$ is a finite set of $\mathcal{C}^{sd}$ and $a_{\rho} \in \Z_{\geq 0}$.
Here $\epsilon_{\rho}\in \{0,\half{1}\}$ based on the parity of $\rho$. More explicitly, if $G_n$ is an orthogonal (resp. symplectic) group, then $\epsilon_{\rho}= 0$ (resp. $\epsilon_{\rho}=\half{1}$) if the image of the homomorphism $\rho: W_{F} \to \GL(V)$ preserves a symplectic bilinear form on $V$, and $\epsilon_{\rho}= \half{1}$ (resp. $\epsilon_{\rho}=0$) if the image of $\rho$ preserves a symmetric bilinear form on $V$. Set $a_{\rho}=-1$ if $\rho$ is not in the finite set $R$. Then the reducibility point $\alpha_{\rho,\sigma}$ is given by $a_{\rho}+\epsilon_{\rho}+1$ and the decomposition \eqref{eq decomp phi_sc} can be rewritten as 
\begin{equation}\label{eq-scparameter}
\phi_{\sigma}= \bigoplus_{\rho \in R}  (\rho\otimes S_{2 \epsilon_{\rho}+1}+\rho\otimes S_{2 \epsilon_{\rho}+3}+\cdots + \rho\otimes S_{2 (\alpha_{\rho,\sigma}-1)+1}).
\end{equation}

Let us now recall the definition of corank for representations, which gives the notion of dimension in our construction of the unitary dual. 

\begin{defn}\label{def corank}
A representation $\pi\in \Pi(G_n)$ is said to be of corank $r$ if there exists an injection
$\pi \hookrightarrow \rho_1 \times \cdots \times \rho_r \times \pi_{sc},$ 
where $\rho_i \in \mathcal{C}$ and $\pi_{sc}\in \mathcal{C}_{cl}$.
\end{defn}

\subsection{Representations of good parity and critical type}
In this subsection, we recall the definition of a representation being \textit{null parity, good parity, bad parity} and \textit{critical type}. 
 We often use the unitarizability of a representation of good parity to determine whether a continuous family of representations are unitary.

\begin{defn}\label{def critical type}
    Let $\pi$ be an irreducible admissible representation of $G_n$, which is a subquotient of 
    $\rho_1\vert\cdot\rvert^{x_1} \times \cdots \times \rho_r\vert\cdot\rvert^{x_r} \rtimes \pi_{sc},$ 
    where $\rho_i \in \mathcal{C}^u$, $x_i\in \R_{\geq 0}$ and $\pi_{sc} \in \mathcal{C}_{cl}$.
    \begin{enumerate}
        \item  Define $\pi$ to be of {\emph{null parity}}
        $($named as ``ugly" in \cite{AM20}$)$ if there exists $1 \leq i \leq r$, such that either $\rho_i \not\in \mathcal{C}^{sd}$ or $x_i \not\in \half{1}\Z$.
        \item Define $\pi$ to be of \emph{good parity} if for any $1 \leq i \leq r$, $\rho_i \in \mathcal{C}^{sd}$ and $x_i \in \alpha_{\rho_i,\pi_{sc}}+ \Z$, where $\alpha_{\rho_i,\pi_{sc}} \in\mathbb{R}_{\geq 0}$ is the reducibility point of
    the pair $(\rho_i, \pi_{sc})$.
        \item Define $\pi$ to be of \emph{bad parity} if $\pi$ is of neither null-parity nor good parity.
        \item Define $\pi$ to be of \emph{critical type} if it is of good parity and for each $\rho \in \mathcal{C}^{sd}$, the set (not multi-set) 
        $\{ x_i \ | \ \rho_i \cong \rho \}$ 
        is either empty or
        forms a segment containing $\alpha_{\rho,\pi_{sc}}$. 
    \end{enumerate}
\end{defn}

We associate the notion of being good parity for representations to that of $L$-parameters as follows. A local $L$-parameter $\phi$ of $G_n$ is said to be of good parity if the local $L$-packet $\Pi_{\phi}$ contains any representation of good parity. In fact, this condition is equivalent to the one that all representations of $\Pi_{\phi}$ being good parity. 

In \cite{HJLLZ24}, the authors made the following conjecture on good parity unitary representations. 

\begin{conj}[{\cite[Conjecture 1.5]{HJLLZ24}}]\label{good parity conjecture}
    Let $G$ be a classical groups, Let $\pi$ be an irreducible unitary representation of $G$ of good parity. Then $\pi$ is unitary if and only if it is of Arthur type (see \S \ref{A packets}).
\end{conj}

This is an extension of a conjecture of M. Tadi{\'c} (\cite[Conjecture 1.1]{Tad22}, critical unitary representations are of Arthur type). In \cite{HJLLZ24}, we discussed several applications of this conjecture and verified it for unitary representations of $G_n$ with corank up to 3.  
Recently, H. Atobe and A. M{\'i}nguez proved Conjecture \ref{good parity conjecture} for symplectic and split odd orthogonal groups using a different idea (\cite{AM25}). The case for even (special) orthogonal groups is proved in \cite{HLL25} as an application of the local theta correspondence. 

\begin{thm}[\cite{AM25, HLL25}]\label{unitiffArthur}
Conjecture \ref{good parity conjecture} holds for symplectic, split odd special orthogonal, and even (special) orthogonal groups. 
\end{thm}

This result gives a full conjectural description of the unitary dual for good parity representations, and greatly simplifies our classification of the unitary dual of corank 4 in this paper. 

\subsection{Local Arthur packets and reduction to good parity}\label{A packets}

In this subsection, we recall the theory of local Arthur packets, as established in \cite{Art13}, and the reduction of the Arthur dual to the good parity case. Recall that $F$ is a $p$-adic local field, and $G_{n} = $ G$_{n}(F)$, where G$_n = Sp_{2n}$ or $SO_{2n+1}$. Let $W_{F}$ be the Weil group of $F$ and $G_{n}^{\vee}(\mathbb{C})$ be the Langlands dual group of $G_{n}$.

\begin{defn}\label{A parameters}

A \emph{local  Arthur parameter} is a homomorphism
$$\psi: W_F \times \SL_2(\mathbb{C}) \times \SL_2(\mathbb{C}) \longrightarrow {{\RG}^\vee_n(\mathbb{C})},$$
\begin{equation}\label{lap}
  \psi = \bigoplus_{i=1}^r \phi_i \otimes S_{a_i} \otimes S_{b_i},  
\end{equation}
satisfying the following conditions: 
\begin{enumerate}
    \item [(1)]$\phi_i(W_F)$ is bounded and consists of semi-simple elements, and $\dim(\phi_i)=d_i$;
    \item [(2)]the restrictions of $\psi$ to the two copies of $\SL_2(\mathbb{C})$ are algebraic, $S_k$ is the $k$-dimensional irreducible representation of $\SL_2(\mathbb{C})$, and 
    $$\sum_{i=1}^r d_ia_ib_i = N:= 
\begin{cases}
2n+1 & \text{ when } G_n=Sp_{2n},\\
2n & \text{ when } G_n=SO_{2n+1}.
\end{cases}
$$ 
\end{enumerate}
The first copy of $\SL_2(\mathbb{C})$ is called the Deligne-$\SL_2(\mathbb{C})$ and is denoted by $\SL_2^D(\mathbb{C})$. The second copy of $\SL_2(\mathbb{C})$ is called the Arthur-$\SL_2(\mathbb{C})$ and is denoted by $\SL_2^A(\mathbb{C})$. 
A local  Arthur parameter $\psi$ given in \eqref{lap} is called \emph{generic} if $b_i=1$ for $i=1, \ldots, r$. 

Given a local  Arthur parameter as in \eqref{lap}, Arthur defined a packet $\Pi_{\psi}$ in \cite[Theorem 2.2.1]{Art13}, called a \emph{local  Arthur packet}, which is a finite subset of $\Pi_u(G_n)$, satisfying certain twisted endoscopic character identities. 
Let $\Psi(G_n)$ be the subset of local Arthur parameters.
We say that a representation $\pi$ is of \emph{\textbf{Arthur type}} if $\pi\in\Pi_\psi$ for some local Arthur parameter $\psi \in \Psi(G_n)$.
Let 
$$\Pi_{A}(G_n)=\{\pi \in \Pi_{\psi} \ | \ \psi \in \Psi(G_n)\}.$$
Representations in $\Pi_{A}(G_n)$ are called \emph{\textbf{Arthur representations}}.
For $\pi\in\Pi_{A}(G_n),$ we let
\[ \Psi(\pi):= \{ \psi \in \Psi(G_n) \ | \ \pi \in \Pi_{\psi}\}.\]
\end{defn}

Now we recall the decomposition of local Arthur parameters and the reduction of the construction of local Arthur packets to the good parity case. By the Local Langlands Correspondence for $\GL_{d_i}(F)$, a bounded representation $\phi$ of $W_F$ can be identified with an irreducible unitary supercuspidal representation $\rho$ of $\GL_{d_i}(F)$ (\cite{Hen00, HT01, Sch13}). Consequently, we may write \eqref{lap} as
\begin{equation}\label{A-param decomp}
  \psi = \bigoplus_{i \in I } \rho_i\vert\cdot\rvert^{x_i} \otimes S_{a_i} \otimes S_{b_i},  
\end{equation}
where $\rho_i$'s are irreducible unitary supercuspidal representations of $\GL_{d_i}(F)$. With this decomposition, we say that $\psi$ is of \emph{good parity} if the following holds. Let $\sigma$ be any irreducible supercuspidal representation of $G_n$. Then $x_i=0$ and $\half{a_i+b_i}\in \alpha_{\rho_i,\sigma} +\Z$ for any $i \in I$.
Equivalently, $\psi$ is of good parity if and only if any representation in the local $L$-packet $\Pi_{\phi_{\psi}}$ is of good parity (Definition \ref{def critical type}(2)). We remark that by Theorem \ref{thm red from nu to gp} below and the construction of good parity local Arthur packets (see the next subsection), $\psi$ is of good parity if and only if any representation in the local Arthur packet $\Pi_{\psi}$ is of good parity.

Let $\psi \in \Psi(G_n).$ Since $\psi$ factors through $\RG_n^\vee(\BC)$, we may rewrite the decomposition \eqref{A-param decomp} as 
\begin{align*}
    \psi= & \bigoplus_{i \in I_{ngp}} (\rho_i \otimes S_{a_i} \otimes S_{b_i} +\rho_i^{\vee} \otimes S_{a_i} \otimes S_{b_i}) \oplus\bigoplus_{i \in I_{gp}} \rho_i \otimes S_{a_i} \otimes S_{b_i},
\end{align*}
where 
\begin{enumerate}
    \item [$\oldbullet$] For any $i \in I_{ngp}$, either $\rho_i\not\cong \rho_i^{\vee}$, or, $\rho_i\cong \rho_i^{\vee}$ and $\half{a_i+b_i}\not\in \alpha_{\rho_i,\sigma} +\Z$;
    \item [$\oldbullet$] For any $i \in I_{gp}$, $\rho_i \cong \rho_i^\vee$ and $\half{a_i+b_i}\in \alpha_{\rho_i,\sigma} +\Z$.
\end{enumerate}
For $\ast \in  \{ngp,\ gp  \}$, define subrepresentations $\psi_{\ast}$ of $\psi$ by 
\[ \psi_{\ast}:= \bigoplus_{i \in I_{\ast}} \rho_i \otimes S_{a_i} \otimes S_{b_i}.\]
Thus $\psi_{gp}$ is of good parity and 
\begin{align}\label{eq decomp red to gp}
    \psi=(\psi_{ngp} + \psi_{ngp}^{\vee})+ \psi_{gp}.
\end{align}

For a unitary supercuspidal representation $\rho$ of $\GL_d(F)$ and $a, b \in \Z_{>0}$, we let 
\begin{equation}\label{u rho a b}
u_{\rho}(a,b):= \begin{pmatrix}
\frac{a-b}{2} & \cdots & \frac{a+b}{2}-1 \\
\vdots & \ddots & \vdots \\
\frac{-a-b}{2}+1 & \cdots & \frac{b-a}{2}
\end{pmatrix}_{\rho}, 
\end{equation}
be the corresponding unitary generalized Speh representation (see \eqref{generalized Speh representation}). This is the unique member of the local Arthur packet $\Pi_{\rho \otimes S_{a} \otimes S_b} (\GL_{abd}(F))$. For each $i \in I_{ngp}$, define $\tau_i$ to be the generalized Speh representation
$u_{\rho_i}(a_i,b_i)$. 
Then we set 
\[ \tau_{\psi_{ngp}}:= \bigtimes_{i \in I_{ngp}} \tau_i,\]
which is irreducible. 

M{\oe}glin showed that the construction of local Arthur packets can be reduced to good parity ones as follows.

\begin{thm}[{\cite[Proposition 5.1]{Moe11b}}]\label{thm red from nu to gp}
Let $\psi\in\Psi(G_n)$ with decomposition \eqref{eq decomp red to gp}. Then, for any $\pi_{gp}\in\Pi_{\psi_{gp}},$ the parabolic induction $\tau_{\psi_{ngp}}\rtimes\pi_{gp}$ is irreducible. As a consequence, \begin{equation}\label{non-unitary A-packet}
    \Pi_\psi=\{\tau_{\psi_{ngp}}\rtimes\pi_{gp}  | ~ \pi_{gp}\in\Pi_{\psi_{gp}}\}.
\end{equation}
\end{thm}

\subsection{Extended multi-segments and operators}

In this subsection, we recall the notion of extended multi-segments and related operations, which provide  important computational tools in the construction of the Arthur dual. 

\begin{defn} [{\cite[Definition 3.1]{Ato20b}}]
(Extended multi-segments)\label{def multi-segment}

\begin{enumerate}
\item
An \emph{extended segment} is a triple $([A,B]_\rho, l, \eta)$,
where
\begin{itemize}
\item
$[A,B]_\rho = \{\rho\vert\cdot\rvert^A, \rho\vert\cdot\rvert^{A-1}, \dots, \rho\vert\cdot\rvert^B \}$ is a segment 
for an irreducible unitary supercuspidal representation $\rho$ of some $\GL_d(F)$; 
\item
$l \in \Z$ with $0 \leq l \leq \frac{b}{2}$, where $b = \#[A,B]_\rho = A-B+1$; 
\item
$\eta \in \{\pm1\}/E$, where $E=\{\pm 1\}$ if $b=2l$ and $E=\{+1\}$ if $b>2l$. 
\end{itemize}
In the statements and formulas in this section, we regard $\eta \in \{\pm 1\}$ by fixing any of its preimage except in Definition \ref{dual segment}, where the choice of the preimage is specified.
\item Consider a multi-set of extended segments of the form $\{([A_i,B_i]_{\rho},l_i,\eta_i)\}_{i \in I_{\rho}}.$
We say that a total order $>$ on $I_{\rho}$ is \emph{admissible} (or \emph{satisfies $(P)$}) if
\[ A_i< A_j, B_i< B_j\Longrightarrow i<j. \]
We say that an admissible order $>$ \emph{satisfies ($P'$)} if
$B_i< B_j\Longrightarrow i<j. $

\item
An \emph{extended multi-segment} for $G_n$ is a union of multi-sets of extended segments indexed by a collection of total ordered sets $(I_{\rho},>)$:
$\EE = \cup_{\rho}\{ ([A_i,B_i]_{\rho}, l_i, \eta_i) \}_{i \in (I_\rho,>)}$ 
such that 
\begin{enumerate}
\item
$I_\rho$ is a totally ordered finite set with a fixed total order $>$ satisfying (P);

\item
$A_i + B_i \geq 0$ for all $\rho$ and $i \in I_\rho$; 

\item
as a representation of $W_F \times \SL_2(\BC) \times \SL_2(\BC)$, 
$\psi_{\EE} = \bigoplus_\rho \bigoplus_{i \in I_\rho} \rho \otimes S_{a_i} \otimes S_{b_i},$
where $(a_i, b_i) = (A_i+B_i+1, A_i-B_i+1)$,
is a local Arthur parameter for $G_n$ of good parity. We shall denote by $\psi_{\EE}$ the local Arthur parameter associated with $\EE$. 
\item The sign condition
\begin{align*}
\prod_{\rho} \prod_{i \in I_\rho} (-1)^{[\frac{b_i}{2}]+l_i} \eta_i^{b_i} = 1
\end{align*}
holds.
\end{enumerate}
\item For each extended multi-segment $\EE$, we denote by $\pi(\EE)$ the representation associated with $\EE$ as in \cite[\S 3.2]{Ato20b}, which is either irreducible or zero. We denote by $\Rep$ the set of extended multi-segments that give nonzero representations, and by $\Rep^{(P')}$ the subset of $\Rep$ that consists of extended multi-segments whose total order on any $I_{\rho}$ satisfies $(P')$.
For a given 
$\EE = \cup_{\rho}\{ ([A_i,B_i]_{\rho}, l_i, \eta_i) \}_{i \in (I_\rho,>)} \in \Rep,$
define that 
$$\EE_{\rho}=\{([A_i,B_i]_{\rho},l_i,\eta_i)\}_{i \in (I_{\rho},>)}\quad {\rm and}\quad 
\EE^{\rho}=\cup_{\rho' \ncong \rho}
\{([A_i,B_i]_{\rho'},l_i,\eta_i)\}_{i \in (I_{\rho'},>)}.$$
Then $\EE=\EE^{\rho} \cup \EE_{\rho}$.
\item Sometimes we write $\EE_{\rho}=\{([A_1,B_1]_{\rho},l_1,\eta_1),\ldots,([A_k,B_k]_{\rho},l_k,\eta_k) \}$, implying that the elements in $\EE_{\rho}$ are listed increasingly with respect to the admissible order of $\EE_{\rho}$. Assume that 
\[ \FF=\{([A_i,B_i]_{\rho},l_i,\eta_i)\}_{i \in (I_{\rho},>)}\quad {\rm and}\quad \FF'=\{([A_i,B_i]_{\rho},l_i,\eta_i)\}_{i \in (I_{\rho}',>)}.\]
Then we let $\FF+\FF'=\{ ([A_i,B_i]_{\rho},l_i,\eta_i)\}_{i \in (I_{\rho}\sqcup I_{\rho}',\gg})$ be the extended multi-segment, 
where the admissible order $\gg$ is defined by $i\gg j$ if and only if $(i,j)\in I_{\rho}\times I_{\rho} \sqcup I_{\rho}' \times I_{\rho}' $ and $i> j$, or, $(i,j) \in I_{\rho}' \times I_{\rho}$.
\item Suppose that $\EE \in \Rep^{(P')}$ and denote $$\FF=\EE_{\rho}= \{([A_i,B_i]_{\rho},l_i,\eta_i)\}_{ i \in (I_{\rho},>)}.$$ 
Let $ x \in \R$. We define
\begin{align*}
\FF_{ >x}&:=\{ ([A_i,B_i]_{\rho},l_i,\eta_i)\}_{i \in I_{\rho}, B_i>x},\\
    \FF_{=x}&:=\{ ([A_i,B_i]_{\rho},l_i,\eta_i)\}_{i \in I_{\rho}, B_i=x},\\
    \FF_{ <x}&:=\{ ([A_i,B_i]_{\rho},l_i,\eta_i)\}_{i \in I_{\rho}, B_i<x},
\end{align*}
with the admissible order inherited from $(I_{\rho}, >)$. Note that $\FF= \FF_{<x}+\FF_{=x}+ \FF_{>x}$. We also write $\FF_{\leq x}= \FF_{<x} + \FF_{=x}$ and $\FF_{\geq x}= \FF_{=x} + \FF_{>x}$. 
\end{enumerate}
\end{defn}

Atobe attached a symbol to each extended multi-segment $\EE$ (\cite[\S3]{Ato20b}). For example, when $\EE=\{([A,B]_\rho,l,\eta)\}$ is a singleton, the symbol is as follows
\[
\EE= 
\left(
\begin{array}{rcl}
\underbrace{\overset{B}{\lhd} \lhd \cdots \overset{B+l-1}{\lhd}}_l 
&
\overset{B+l}{\odot} \odot \cdots \odot \overset{A-l} \odot 
&
\underbrace{\overset{A-l+1}{\rhd} \cdots \rhd \overset{A}{\rhd}}_l
\end{array}
\right)_\rho.
\] 
Here, $\odot\cdots\odot$ represents an alternating sequence of $\oplus$ and $\ominus$, 
starting with $\oplus$ if $\eta = 1$ (resp. $\ominus$ if $\eta = -1$). When $\EE$ is not a singleton, we stack each row vertically. See the following for an example: 

\begin{exmp}
 Let $\rho$ be the trivial representation of $\GL_1(F)$. The symbol
\[\EE= \scalebox{0.8}{\bordermatrix{
  & 0 &1 & 2 &3 \cr
  &\lhd & \oplus & \ominus & \rhd  \cr
  & &  & \lhd &\rhd  \cr
  & &  &  &  \ominus \cr
}}_{\rho}\]
corresponds to the extended multi-segment
$\EE= \{ ([ 3,0]_{\rho},1,1),([3,2]_{\rho},1,1),([3,3]_{\rho},0,-1)\}$
of $Sp_{34}$. The associated local Arthur parameter is
    $\psi_{\EE}= \rho \otimes S_{4}\otimes S_{4} + \rho \otimes S_{6}\otimes S_{2} + \rho \otimes S_7 \otimes S_1.$
\end{exmp}

Given an extended multi-segment $\mathcal{E}$, one can use the algorithms in \cite{HJLLZ24} to compute the representation $\pi(\mathcal{E})$. In particular, there are necessary conditions on $\mathcal{E}$ for $\pi(\mathcal{E})$ to be of Arthur type. This relies on the implementation of various operators on extended multi-segments as follows.

The first operator is called the \emph{row exchange}. It is used to change the admissible order of an extended multi-segment, but does not affect the corresponding local Arthur parameters. We say that $k<k+1$ are \emph{adjacent} in a total order set $(I_{\rho},>)$ if there does not exist $i \in I_{\rho}$ such that $k< i <k+1$.

\begin{defn}[{\cite[Section 4.2]{Ato20b}}, Row exchange]\label{def row exchange} 
Let $\EE\in\Rep^{(P')}$ with
$$\EE_{\rho}=\{([A_i,B_i]_{\rho},l_i,\eta_i)\}_{i \in (I_{\rho},>)}.$$
Assume that $k<k+1$ are adjacent in $(I_{\rho},>)$. Let $\gg$ be the total order on $I_\rho$ defined by $k\gg k+1$ and if $(i,j)\neq (k,k+1)$, then $ i \gg j$ if and only if $
i >j .$ 

If $\gg$ is not an admissible order on $I_{\rho}$, then we define $R_k(\EE)=\EE$. Otherwise, we define 
\[R_{k}(\EE_{\rho})=\{([A_i,B_i]_{\rho},l_i',\eta_i')\}_{i \in (I_{\rho},\gg)},\]
where $( l_i',\eta_i')=(l_i,\eta_i)$ for $i \neq k,k+1$, and $(l_k',\eta_k')$ and $(l_{k+1}', \eta_{k+1}')$ are given as follows.  Denote $\epsilon=(-1)^{A_k-B_k}\eta_k\eta_{k+1}$.
\begin{enumerate}
    \item [(1)] Assume that $[A_k,B_k]_{\rho} \supset [A_{k+1},B_{k+1}]_{\rho}$.    
    We set $(l_{k+1}',\eta_{k+1}')=(l_{k+1}, (-1)^{A_k-B_k}\eta_{k+1})$ in this case. 
    \begin{enumerate}
    \item [(a)] If $\epsilon=1$ and $b_k- 2l_k < 2(b_{k+1}-2l_{k+1})$, then
    \[ (l_k', \eta_{k}')= (b_k-(l_k+ (b_{k+1}-2l_{k+1})), (-1)^{A_{k+1}-B_{k+1}} \eta_k);  \]
    \item [(b)] If $\epsilon=1$ and $b_k- 2l_k \geq  2(b_{k+1}-2l_{k+1})$, then
    \[ (l_{k}', \eta_{k}')= (l_k+ (b_{k+1}-2l_{k+1}), (-1)^{A_{k+1}-B_{k+1}+1} \eta_k);  \]
    \item [(c)] If $\epsilon=-1$, then
    $(l_{k}', \eta_{k}')= (l_k- (b_{k+1}-2l_{k+1}), (-1)^{A_{k+1}-B_{k+1}+1} \eta_k).$
\end{enumerate}
    \item [(2)] Assume that $ [A_k,B_k]_{\rho} \subseteq [A_{k+1},B_{k+1}]_{\rho}$.     
    We set $(l_{k}',\eta_{k}')=(l_{k}, (-1)^{A_{k+1}-B_{k+1}}\eta_{k})$ in this case. 
    \begin{enumerate}
   \item [(a)] If $\epsilon=1$ and $b_{k+1}- 2l_{k+1} < 2(b_{k}-2l_{k})$, then
    \[ (l_{k+1}', \eta_{k+1}')= (b_{k+1}-(l_{k+1}+ (b_{k}-2l_{k})), (-1)^{A_{k}-B_{k}} \eta_{k+1});  \]
    \item [(b)] If $\epsilon=1$ and $b_{k+1}- 2l_{k+1} \geq  2(b_{k}-2l_{k})$,
    then
    \[ (l_{k+1}', \eta_{k+1}')= (l_{k+1}+ (b_{k}-2l_{k}), (-1)^{A_{k}-B_{k}+1} \eta_{k+1});  \]
    \item [(c)] If $\epsilon=-1$, then
    $(l_{k+1}', \eta_{k+1}')= (l_{k+1}- (b_{k}-2l_{k}), (-1)^{A_{k}-B_{k}+1} \eta_{k+1}).$
\end{enumerate}
\end{enumerate}
Finally, we define that $R_{k}(\EE)= \EE^{\rho} \cup R_{k}(\EE_{\rho})$.

If $\gg$ is another admissible order on $I_{\rho}$, then we deform  $\EE_{\rho}$ into
$\{([A_i,B_i]_{\rho},(l_i)_{\gg}, (\eta_i)_{\gg})\}_{i \in (I_{\rho}, \gg)}$ 
by applying a sequence of row exchanges on $\EE_{\rho}$. We shall denote the resulting extended multi-segment by $\EE_{\rho, \gg}$.

\end{defn}

Next, we recall the definition of the operators $sh_j^{d}, add_j^{d}$ on extended multi-segments. These operators can be useful in constructing new extended multi-segments.

\begin{defn}[Shift, Add]
Let $\EE = \cup_{\rho}\{ ([A_i,B_i]_{\rho}, l_i, \eta_i) \}_{i \in (I_\rho,>)}$ be an extended multi-segment. For $j \in I_{\rho'}$ and $d \in \Z$, we define the following operators. 

\begin{enumerate}
    \item [1.] Define $sh_j^{d}(\EE)= \cup_{\rho}\{ ([A_i',B_i']_{\rho}, l_i, \eta_i) \}_{i \in (I_\rho,>)}$ with 
    \[ [A_i',B_i']_{\rho}= \begin{cases}
    [A_i+d,B_i+d]_{\rho} & \text{ if }\rho=\rho' \text{ and } i = j,\\
     [A_i,B_i]_{\rho} & \text{ otherwise. }\end{cases} \]
    Set $sh^d_{\rho'}=\sum_{j\in I_{\rho'}} sh_j^{d}$ and $sh^d:=\sum_{\rho'} sh^d_{\rho'}$.
     \item [2.] Define $add_j^{d}(\EE)= \cup_{\rho}\{ ([A_i',B_i']_{\rho}, l_i', \eta_i) \}_{i \in (I_\rho,>)}$ with 
    \[ ([A_i',B_i']_{\rho},l_i')= \begin{cases}
    ([A_i+d,B_i-d]_{\rho},l_i+d) & \text{ if }\rho=\rho' \text{ and } i = j,\\
     ([A_i,B_i]_{\rho},l_i) & \text{ otherwise, }\end{cases} 
    \]
    Set $add^d_{\rho'}=\sum_{j\in I_{\rho'}} add_j^{d}$ and $add^d:=\sum_{\rho'} add^d_{\rho'}$.
\end{enumerate}
It is immediate that the operators commute with each other, so, we denote the composition by summation. 
We only use these notations in the case that the resulting object is still an extended multi-segment.
\end{defn}

The next operator is called the \emph{union-intersection}, which would change the corresponding local Arthur parameters if acting non-trivially. 

\begin{defn}[{\cite[Section 5.2]{Ato20b}}, Union-intersection]\label{ui def}
 Let $\EE$ be an extended multi-segment with
$\EE_{\rho}=\{([A_i,B_i]_{\rho},l_i,\eta_i)\}_{i \in (I_{\rho},>)}.$
For $k< k+1$ adjacent in $(I_{\rho},>)$, we define an operator $ui_k$, called the \emph{union-intersection}, on $\EE$ as follows. 
  Denote $\epsilon=(-1)^{A_k-B_k}\eta_k \eta_{k+1}.$ If $A_{k+1}>A_k$, $B_{k+1}>B_k$ and any of the following cases holds:
\begin{enumerate}
    \item [{Case 1}.] $ \epsilon=1$ and $A_{k+1}-l_{k+1}=A_k-l_k,$
    \item [{Case 2}.] $ \epsilon=1$ and $B_{k+1}+l_{k+1}=B_k+l_k,$
    \item [{Case 3}.] $ \epsilon=-1$ and $B_{k+1}+l_{k+1}=A_k-l_k+1,$
\end{enumerate}
we define that 
$ui_{k}(\EE_{\rho})=\{ ([A_i',B_i']_{\rho},l_i',\eta_i')\}_{i \in (I_{\rho}, >)},$
where $ ([A_i',B_i']_{\rho},l_i',\eta_i')=([A_i,B_i]_{\rho},l_i,\eta_i)$ for $i \neq k,k+1$, and $[A_k',B_k']_{\rho}=[A_{k+1},B_k]_{\rho}$, $[A_{k+1}',B_{k+1}']_{\rho}=[A_k,B_{k+1}]_{\rho}$, and $( l_k', \eta_k', l_{k+1}',\eta_{k+1}' )$ are given case by case as follows:
\begin{enumerate}
    \item[$(1)$] in Case 1, $( l_k', \eta_k', l_{k+1}',\eta_{k+1}' )= (l_k,\eta_k, l_{k+1}-(A_{k+1}-A_k), (-1)^{A_{k+1}-A_k}\eta_{k+1})$;
    \item [$(2)$] in Case 2, if $b_k-2l_k \geq A_{k+1}-A_k$, then
    \[( l_k', \eta_k', l_{k+1}',\eta_{k+1}' )= (l_k+(A_{k+1}-A_k),\eta_k, l_{k+1}, (-1)^{A_{k+1}-A_k}\eta_{k+1}),\]
    and if $b_k-2l_k < A_{k+1}-A_k$, then
    $( l_k', \eta_k', l_{k+1}',\eta_{k+1}' )= (b_k-l_k,-\eta_k, l_{k+1}, (-1)^{A_{k+1}-A_k}\eta_{k+1});$
    \item [$(3)$] in Case 3, if $l_{k+1} \leq  l_k$, then
    $( l_k', \eta_k', l_{k+1}',\eta_{k+1}' )= (l_k,\eta_k, l_{k+1}, (-1)^{A_{k+1}-A_k}\eta_{k+1})$, and 
    if $l_{k+1}> l_{k}$, then
    $( l_k', \eta_k', l_{k+1}',\eta_{k+1}' )= (l_k,\eta_k, l_{k}, (-1)^{A_{k+1}-A_k+1}\eta_{k+1});$
    \item [$(3')$] if we are in Case 3 and $l_k=l_{k+1}=0$, then we delete $ ([A_{k+1}',B_{k+1}']_{\rho},l_{k+1}',\eta_{k+1}')$ from $ui_k(\EE_{\rho})$.
\end{enumerate}
\end{defn}

The union-intersection operator can be extended to non-adjacent extended segments as follows.

\begin{defn} \label{def ui}
Let $\EE$ be an extended multi-segment with
$\EE_{\rho}=\{ ([A_i,B_i]_{\rho},l_i,\eta_i)\}_{i\in (I_{\rho,>})}.$
Given $i,j \in I_{\rho}$, we define that $ui_{i,j}(\EE_{\rho})=\EE_{\rho}$ unless the following conditions hold:
\begin{enumerate}
    \item [1.]$ A_i< A_j$, $B_i <B_j$, and $j\gg i$ are adjacent under some admissible order $\gg$ on $I_{\rho}$. 
    \item [2.] The operator $ui_i$ is applicable on $\EE_{\rho,\gg}$.
\end{enumerate}
In this case, we define that $ui_{i,j}(\EE_{\rho}):=(ui_{i}(\EE_{\rho,\gg}))_{>}$, so that the admissible order of $ui_{i,j}(\EE_{\rho})$ and $\EE_{\rho}$ are the same (if the $ui_i$ is of type $3'$, then we delete the $j$-th row). Finally, we define $ui_{i,j}(\EE):= \EE^{\rho} \cup ui_{i,j}(\EE_{\rho})$.

We say that $ui_{i,j}$ is applicable on $\EE$ if $ui_{i,j}(\EE) \neq \EE$. Furthermore, we say that $ui_{i,j}$ is of type 1, 2, 3, or $3'$ if the corresponding operator $ui_i$ is of type 1, 2, 3, or $3'$, respectively, in Definition \ref{ui def}.
\end{defn}

Let $\pi$ be an irreducible admissible representation of $G_n.$ Aubert showed that there exists $\varepsilon\in\{\pm 1\}$ such that the virtual representation defined by
$$
\hat{\pi}:=\varepsilon\sum_P (-1)^{\mathrm{dim}(A_P)}[\mathrm{Ind}_{P}^{G_n}(Jac_P(\pi))]
$$
is an irreducible representation (\cite{Aub95}). The above sum is over all standard parabolic subgroups $P$ of $G_n$, where $A_P$ denotes the maximal split torus in the center of the Levi subgroup of $P$, and $Jac_P$ denotes the Jacquet module along $P.$ We say that $\hat{\pi}$ is the Aubert-Zelevinsky dual or Aubert-Zelevinsky involution of $\pi.$ 

The next operator is known as the dual, which sends a representation corresponding to an extended multi-segment to its Aubert-Zelevinsky dual.

\begin{defn}[{\cite[Definition 6.1]{Ato20b}}, Dual]\label{dual segment}
Let $\EE= \cup_\rho \{([A_i,B_i]_{\rho},l_i,\eta_i)\}_{i\in (I_\rho, >)}$ be an extended multi-segment such that the admissible order $>$ on $I_{\rho}$ satisfies $(P')$ for all $\rho$. We define 
$$dual(\EE):=\cup_{\rho}\{([A_i,-B_i]_{\rho},l_i',\eta_i')\}_{i\in (I_\rho, >')},$$ as follows:
\begin{enumerate}
    \item The order $>'$ is defined by $i>'j$ if and only if $j>i.$ 
    \item Set \begin{align*}
l_i'=\begin{cases}
l_i+B_i  & \mathrm{if} \, B_i\in\mathbb{Z},\\
 l_i+B_i+\frac{1}{2}(-1)^{\alpha_{i}}\eta_i  & \mathrm{if} \, B_i\not\in\mathbb{Z},
\end{cases}\quad {\rm and}\quad 
\eta_i'=\begin{cases}
(-1)^{\alpha_i+\beta_i}\eta_i  & \mathrm{if} \, B_i\in\mathbb{Z},\\
 (-1)^{\alpha_i+\beta_i+1}\eta_i  & \mathrm{if} \, B_i\not\in\mathbb{Z},
\end{cases}
\end{align*}
where $\alpha_{i}=\sum_{j\in I_\rho, j<i}a_j,$ and $\beta_{i}=\sum_{j\in I_\rho, j>i}b_j,$ $a_j=A_j+B_j+1$, $b_j=A_j-B_j+1$.
\item When $B_i\not\in\mathbb{Z}$ and $l_i=\frac{b_i}{2}$, set $\eta_i=(-1)^{\alpha_i+1}.$
\end{enumerate}
Finally, we define $dual(\EE_{\rho}):= (dual(\EE))_{\rho}$.
\end{defn}

% As mentioned, Atobe showed that the operator $dual$ corresponds to the Aubert-Zelevinsky dual for the representations.

\begin{thm}[{\cite[Theorem 6.2]{Ato20b}}]\label{thm dual}
    Let $\EE \in \Rep^{(P')}$. Then $\pi(dual(\EE))= \widehat{\pi(\EE)}$ holds.
\end{thm}

The last operator we introduce is known as the partial dual. 

\begin{defn}[{\cite[Definition 6.5]{HLL22}}, Partial dual]\label{def partial dual}
Let $\EE\in \Rep^{(P')}$, and let $$\FF:=\EE_{\rho}= \{([A_i,B_i]_{\rho},l_i,\eta_i)\}_{i \in (I_{\rho},>)}.$$
For $i \in I_{\rho}$, denote 
\[\alpha_i= \sum_{ j<i} (A_j+B_j+1)\quad {\rm and}\quad {\beta_i= \sum_{j> i} (A_j-B_j+1)}.\]
Assume that there exists $k \in I_{\rho}$ such that
\begin{enumerate}
    \item [(1)] $B_k=\half{1},l_k=0$;
    \item [(2)] $(-1)^{\alpha_k}\eta_k= -1$;
    \item [(3)] for any $i < k$, $B_i < \half{1}$.
\end{enumerate}
Then we define $dual_k^{+}(\FF)$ as follows. We write the decomposition
\[ \FF= \FF_1 + \{([A_k,1/2]_{\rho},0,\eta_k)\} + \FF_2,\]
where $\FF_1=\FF_{<1/2}$ and $\FF_2=\FF_{>1/2}$, and then write
\[dual(\FF)= \widetilde{\FF_2} + \{([A_k,-1/2]_{\rho},0,(-1)^{\beta_k}) \}+ \widetilde{\FF_1},\]
where $\widetilde{\FF_2}=(dual(\FF))_{<-1/2}$ and $\widetilde{\FF_1}=(dual(\FF))_{>-1/2}$. Next, write
\[ dual(\FF'):=dual(\widetilde{\FF_2} + \{([A_k,1/2]_{\rho},0,(-1)^{\beta_k+1})\} + \widetilde{\FF_1})= \widetilde{\widetilde{\FF_1}} + \{([A_k,-1/2]_{\rho},0,-\eta_k)\}+ \widetilde{\widetilde{\FF_2}},\]
where $\widetilde{\widetilde{\FF_1}}=(dual(\FF'))_{<-1/2}$ and $\widetilde{\widetilde{\FF_2}}=(dual(\FF'))_{>-1/2}$. Then we define
\[ dual_k^{+}( \FF):= \widetilde{\widetilde{\FF_1}} + \{([A_k,-1/2]_{\rho},0,-\eta_k)\} + \FF_2, \]
and say that $dual_k^{+}$ is applicable on $\FF$.

Assume that $dual(\FF)$ satisfies above conditions (1) -- (3). Then we define
\[ dual_k^{-}(\FF):= dual \circ dual_{k}^{+} \circ dual (\FF), \]
and say that $dual_k^{-}$ is applicable on $\FF$.
We call this operators $dual_k^{+}$, $dual_k^{-}$ the \emph{partial dual}, and denote by $dual_k$ if there is no ambiguity. 

Finally, we define $dual_k(\EE):= \EE^{\rho} \cup dual_k(\EE_{\rho}).$ 
\end{defn}

Next, we introduce certain collections of the operators defined above. Note that one can easily construct the inverses of the add, shift, union-intersection, and dual operators, when they are applicable on an extended multi-segment. 

\begin{defn}\label{def raising operator}\ 
\begin{enumerate}
    \item We say that an operator $T$ defined above is a \emph{raising operator} if it is of the form  
$ ui_{i,j}^{-1}$, $dual \circ ui_{j,i} \circ dual,$ or $dual_k^{-}$.
\item We say that an extended multi-segment $\EE \in \Rep$ is \emph{absolutely maximal} if there is no raising operator applicable on $\EE$.
\end{enumerate}
\end{defn}

The definition of raising operators is used to exhaust Arthur packets in \cite{HLL22}. In later sections, they will be used to restrict the possible forms of extended multi-segments corresponding to representations of Arthur type, which allow us to exhaust the Arthur dual.

\subsection{Preservation of irreducibility and unitarizability}

In this subsection, we recall some important technical results regarding when and how irreducibility and unitarizability of a representation are preserved. 
We begin with an irreducibility criterion of Tadi\'{c}, which is useful when considering parabolic inductions with a supercuspidal representation.
For any $\pi \in \Pi(G_n)$, let $\supp(\pi)$ denote the supercuspidal support of $\pi$. We denote the set of irreducible representations of $GL_n(F)$ (resp. $G_n(F)$) by $\Irr$ (resp. $\Irr^{cl}$). Let $\mathcal{D}$ denote all irreducible essentially square integrable
representations of $\GL_{n}(F)$, $n \geq 1$, and $e(\delta)$ denote the exponent of any representation $\delta \in \mathcal{D}$. Let 
$M(\mathcal{D})$ denote all finite multi-sets of representations in $\mathcal{D}$.

\begin{thm}[{\cite[(2.23)]{Tad23}}]\label{Tadirred1}
    Let $\sigma \in \mathcal{C}_{cl}$, $\pi=L(d)$, for some $d\in M(\mathcal{D})$. Suppose now that $\supp(\pi)$ does not contain $\rho\lvert \cdot \rvert^\alpha$ or $\rho\lvert \cdot \rvert^{-\alpha}$.
% Assume that all members of $\supp(\pi)$ are contained in $\{\rho\lvert \cdot \rvert^{k+x}:k\in\Z\}$, for some fixed $x\in\tfrac12\Z$.
Denote by $d_{>0}$ (resp. $d_{<0}$) the multiset consisting of all $\delta$ in $d$
such that $e(\delta)>0$ (resp $e(\delta)<0$), counted with multiplicities.
Then, if $\pi$ is a ladder representation, or, if $\alpha\leq 1$ and all members of $\supp(\pi)$ are contained in $\{\nu^{k+\alpha}\rho;k\in\Z\}$, then the following holds
\begin{equation} \label{eq: irredcriterion,sc}
L(d)\rtimes\sigma \text{ is reducible } \iff L(d_{>0})\times L(d_{<0})\check{\ } \text{ is reducible}.
\end{equation}
\end{thm}

Next, we recall some basic notation on the support of representations, leading up to the result of Jantzen decomposition. 

\begin{defn} \label{defJantzendecomp}
    Let $X$ be a subset of $\mathcal{C}$.
\begin{enumerate}
    \item [(1)] $X$ is self-contragredient if for any $\rho \in X$, the contragredient of $\rho$ is also in $X$.
    \item [(2)] For an irreducible admissible representation $\beta$ of $\GL_d(F)$, we say that $\beta$ is supported on $X$ if the supercuspidal support of $\beta$ is contained in $X$.
    \item [(3)] Let $\pi$ be an irreducible admissible representation of $G_n$ that appears as an irreducible subquotient of 
    $\theta_1 \times \cdots \times \theta_f \rtimes \sigma,$
    where $\theta_i \in \mathcal{C}$ and $\sigma \in \mathcal{C}_{cl}$.
     We say that $\pi$ is supported on $X\cup \{\sigma\}$ for some self-contragredient $X\subseteq \mathcal{C}$ if $\theta_i \in X$ for $i=1,\dots, r$. 
    \item [(4)] Fix a supercuspidal representation $\sigma$ of $G_m$ and a self-contragredient subset $X\subseteq \mathcal{C}$. We denote by $\Irr(X;\sigma)$ the set of irreducible admissible representations $\pi$ of $G_n$ with $n\geq m$ such that $\pi$ is supported on $X \cup \{\sigma\}$. 
    \item [(5)] Suppose that $X$ is self-contragredient. Let $X= X_1 \sqcup X_2$ be a partition of $X$. 
    Such a partition is called regular if $X_1$ is self-contragredient and 
    \[ \theta \in X_1 \Longrightarrow \theta \vert\cdot\rvert^{1} \not\in X_2.  \]
    That is, there is no reducibility among $X_1$ and $X_2$. 
\end{enumerate}
\end{defn}

The following results are on the Jantzen decomposition, which give us another criterion to determine irreducibility.

\begin{thm}[{\cite[Theorem 9.3]{Jan97}}]\label{thm Jantzen} 
Let $X$ be a subset of $\mathcal{C}$ that is self-contragredient and $X=X_1 \sqcup X_2$ is a regular partition. Let $\sigma \in \mathcal{C}_{cl}$.
\begin{enumerate}
    \item [(i)] For any $\pi \in  \Irr(X;\sigma)$, there exist irreducible admissible representations $\beta_1, \beta_2$ of general linear groups and $ X_1(\pi),X_2(\pi)$ of lower rank classical groups such that 
    \begin{enumerate}
        \item [$\oldbullet$] $\beta_i$ is supported on $X_i$ and $X_i(\pi)$ is supported on $X_i \cup \{\sigma\}$; and
        \item [$\oldbullet$] there are injections
        \[ \pi \hookrightarrow \beta_1 \rtimes X_2(\pi)\quad {\rm and}\quad  \pi \hookrightarrow \beta_2 \rtimes X_1(\pi). \]
    \end{enumerate}
    The representations $X_i(\pi)$ are uniquely determined by $\pi$.
    \item [(ii)] The map 
    \begin{align*}
        \Irr(X; \sigma )\longrightarrow \Irr(X_1;\sigma) \times \Irr(X_2; \sigma ),\quad {\rm with}\ 
         \pi \mapsto  (X_1(\pi), X_2(\pi))
    \end{align*}
    is a bijection. We denote the inverse map by $\Psi_{X_1,X_2}$.
    \item [(iii)] For $i=1,2$, let $\beta_i$ be an irreducible admissible representation of a general linear group supported on $X_i$ and $\gamma_i\in \Irr(X_i; \sigma)$. Here we allow $\beta_i$ to be the trivial representation of $\GL_0(F)$. Then 
    \[ (\beta_1 \times \beta_2)\rtimes \Psi_{X_1,X_2}(\gamma_1,\gamma_2) \text{ is irreducible} \Longleftrightarrow \text{ both }\beta_i \rtimes \gamma_i \text{ are irreducible.} \]
    \item [(iv)] For any  $\pi \in \Irr(X;\sigma)$, $\pi$ is tempered (resp. square-integrable) if and only if $X_1(\pi), X_{2}(\pi)$ are both tempered (resp. square-integrable). 
\end{enumerate}
\end{thm}

Now, we recall when and how unitarizability of representations of $G_n$ are preserved. These techniques provide the basic tools for our classification of the corank 4 unitary dual and motivate the definition of $\Pi_{\overline{A}}^{\lim}(G_n)$. They can be used to prove or disprove unitarizability.

\begin{prop}[{\cite{Tad93}},{\cite[Lemma $3.3$]{LMT04}},{\cite[Section $2$]{MT11}}]\label{list} 
Fix a standard parabolic subgroup $P=MN$ of $G$ 
with $M \cong \GL_{m_1}(F) \times \dots \times \GL_{m_k}(F) \times G_0$. 
Let $\tau_i \in \Irr(\GL_{m_i}(F))$, $\pi_0 \in \Irr(G_0)$ 
and set $\pi_M = \tau_1 \otimes \dots \otimes \tau_k \otimes \pi_0 \in \Irr(M)$.
\par

\begin{enumerate}
\item(Unitary induction (UI))
If $\pi_M$ is unitary, 
then $\Ind_P^G(\pi_M)$ is a direct sum of irreducible unitary representations of $G$. 

\item(Unitary reduction (UR))
If $\pi_M$ is hermitian and if $\Ind_P^G(\pi_M)$ is irreducible and unitary, 
then $\pi_M$ is also unitary. 

\item(Complementary series (CS))
Let $x_1(t), \dots, x_r(t) \colon [0,1] \rightarrow \R$ be continuous functions, 
and set 
\[
\Pi_t = \tau_1\lvert \cdot \rvert^{x_1(t)} \times \dots \times \tau_r\lvert \cdot \rvert^{x_r(t)} \rtimes \pi_0.
\]
If $\Pi_t$ is irreducible and hermitian for $0 \leq t < 1$, and if $\Pi_0 = \Ind_P^G(\pi_M)$ is unitary, 
then all irreducible subquotients of $\Pi_t$ are unitary for $0 \leq t \leq 1$. 

\item(Beyond the first reducibility point (RP1))
Suppose that 
\begin{itemize}
\item
$k=1$; 
\item
$\tau_1 = \times_{i=1}^r\Sp(\rho_i,c_i,d_i)$ is a product of unitary Speh representations
with $\rho_i \cong \rho_i^\vee$ for $1 \leq i \leq r$; 
\item
$\pi_0$ is of Arthur type of good parity.
\end{itemize}
Let $\psi_M \in \Psi(M)$ be an $A$-parameter with $\pi_M \in \Pi_{\psi_M}$, 
and let $R_P(w,\pi_M,\psi_M)$ be the normalized intertwining operator defined by Arthur \cite[Section 2.4]{Art13} 
with $w \in W(M,G)$. 
Assume further that $R_P(w,\pi_M,\psi_M)$ is not a scalar (so in particular $\Ind_P^G(\pi_M)$ is reducible).
Then $\tau\lvert \cdot \rvert^s \rtimes \pi_0$ is not unitary for sufficiently small  $s > 0$. 
\end{enumerate}
\end{prop}

To end this subsection, we recall the definition of a weakly real representation. 

\begin{defn}\label{def weakly real}
    An representation $\pi \in \Irr^{cl}$ is called weakly real if it is a subquotient of a representation of the form 
    \[\lvert\cdot\rvert^{x_1} \rho_1 \times \ldots \lvert\cdot\rvert^{x_k} \rho_k \rtimes \sigma,\]
    where $\rho_i \in \mathcal{C}$ are self-dual,  $x_i \in \R$ and $\sigma \in \mathcal{C}_{cl}$.
\end{defn}

The key theorem by Tadi{\'c} is the following:
\begin{thm} [{\cite[Theorem $4.2$]{Tad09a}}]\label{red to weakly real}
    For any unitarizable $\pi \in \Irr^{cl}$, there exists a unitarizable $\theta \in \Irr$ and a weakly real $\pi' \in \Irr^{cl}$ such that 
    \[\pi \cong \theta \rtimes \pi'.\]
\end{thm}
Since the unitary dual of $GL_n(F)$ is well-known by \cite[Theorem $7.5$]{Tad86}, this reduces the classification of the unitary dual for classical groups to the weakly real case, which is the focus of this paper.

\subsection{\texorpdfstring{Definition of $\Pi_{\overline{A}}(G_n)$ and $\Pi_{\overline{A}}^{\lim}(G_n)$}{}}

In this subsection, we recall the definitions of $\Pi_{\overline{A}}(G_n)$ and the unitary dual candidate set $\Pi_{\overline{A}}^{\lim}(G_n)$.

\begin{defn}[{\cite[Definition 5.5]{HJLLZ24}}]\label{def closure A}\

\begin{enumerate}
    \item [(a)] For $k \in \Z_{\geq 0}$, we define subsets $ \Pi_{\overline{A}}^{(k)}(G_n)$ of $\Pi_u(G_n)$ inductively as follows. Set $ \Pi_{\overline{A}}^{(0)}(G_n):= \Pi_{A}(G_n)$. For $k \geq 1$, and any $\pi \in \Pi(G_n)$, we say that $\pi \in \Pi_{\overline{A}}^{(k)}(G_n)$ if there exist data: 
    \begin{enumerate}
        \item[$(*)$] 
        \begin{enumerate}
            \item[$\bullet$] $s \in \mathbb{Z}_{\geq 0}$; 
            \item[$\bullet$] $k_i, a_i, b_i \in \mathbb{Z}_{>0}$, for $1 \leq i \leq s$, $\sum_{i=1}^s k_ia_ib_i\leq n$;
            \item[$\bullet$] irreducible unitary supercuspidal representations $\rho_i$ of general linear groups $\GL_{k_i}(F)$ (not necessarily selfdual), for $1 \leq i \leq s$;
            \item[$\bullet$] $\pi_A \in\Pi_{A, \, gp}(G_{m})$ for some $m \leq n$, $n=m+\sum_{i=1}^s k_ia_ib_i$; \item[$\bullet$] $\underline{y}, \underline{z} \in \R^s$;
        \end{enumerate}
    \end{enumerate}
    % if there exists a triple 
    % \begin{align}\label{eq triple A bar Gn}
    %     ( \Pi_{\underline{x}}= \Pi_{x_1,\ldots, x_s}= u_{\rho_1}(a_1,b_1)\vert\cdot\rvert^{x_1}\times \cdots \times u_{\rho_s}(a_s,b_s)\vert\cdot\rvert^{x_s} \rtimes \pi_A, \underline{y} \in \R^s, \underline{z} \in \R^s),
    % \end{align}
satisfying the following conditions.
\begin{enumerate}
    \item [(1)] The point $\underline{y}$ lies in the set 
    \[U:=\{\underline{w} \in \R^s\ | \ \Pi_{\underline{w}} \text{ is irreducible and Hermitian}\},\]
    and { $\Pi_{\underline{y}}= \tau \rtimes \pi^{(k-1)}$  for some unitary representation $\tau$ of $\GL_d(F)$ and $\pi^{(k-1)}\in \Pi_{\overline{A}}^{(k-1)}(G_{n-d})$.} 
        \item [(2)] The point $\underline{z}$ lies in the unique connected component of $U$ containing $\underline{y}$ and $\pi \cong \Pi_{\underline{z}}$.
\end{enumerate}
Here, for any $\underline{x} = \{x_1,\ldots, x_s\} \in \R^s$, we denote 
\begin{equation*}\label{eq triple A bar Gn}
\Pi_{\underline{x}}= u_{\rho_1}(a_1,b_1)\vert\cdot\rvert^{x_1}\times \cdots \times u_{\rho_s}(a_s,b_s)\vert\cdot\rvert^{x_s} \rtimes \pi_A.
\end{equation*}
Finally, we let $ \Pi_{\overline{A}}(G_n):= \bigcup_{k \in \Z_{\geq 0}} \Pi_{\overline{A}}^{(k)}(G_n)$.

\item [(b)] We say that $\pi \in \Pi_{\overline{A}}^{\lim}(G_n)$ if there exist data as in $(*)$ above, satisfying 
\begin{itemize}
    \item [($1'$)] The point $\underline{y}$ lies in the set $U$ and  $\Pi_{\underline{y}} \in \Pi_{\overline{A}}(G_n)$.
    
    \item [($2'$)] The point $\underline{z}$ lies in the \textbf{closure} of the unique connected component of $U$ containing $\underline{y}$ and $\pi $ is a subquotient of $ \Pi_{\underline{z}}$.
\end{itemize}

    \item [(c)] We say that $\pi \in \Pi_{\overline{A}}^{\lim'}(G_n)$ if there exist data as in $(*)$ above, satisfying the conditions $(1)$, and $(2')$ above.
\end{enumerate}
\end{defn}

Next, we recall the following theorem of Tadi{\'c} and Vogan on a necessary condition of non-isolated unitary representations, 
which is particularly helpful to determine the isolated representations in $\Pi_{\overline{A}}^{\lim}(G_n)$. 
See \cite[Section 2.4]{HJLLZ24} for a summary of the Fell topology on $\Pi_u(G)$ and useful properties. 

\begin{thm}[{\cite[Theorem 5.6]{Tad83}, \cite[Theorem 3]{Vog07}}]\label{vogan thm}
    Let $\sigma\in \Pi_u(G)$ be a non-isolated representation, and let $\{\pi_n\}\subset \Pi_u(G)$ be a sequence  distinct from $\sigma$ such that $\pi_n$ converges to $ \sigma$ in the Fell topology. Then there is a subsequence $\{\pi_{n_j}\}$ in $\Pi_u(G)$, a parabolic subgroup $P=MN$ of $G$, an irreducible admissible representation $\rho$ of $M$, and a sequence of unramified characters $\chi_j$ of $M$ with the following properties. 
\begin{enumerate}
\item The sequence $\chi_j$ converges to the trivial character.
\item The induced representation $\Ind_{P}^G(\rho\otimes \chi_j) \cong \pi_{n_j}$.
\item The representation $\sigma$ is a composition factor of $\Ind_P^G(\rho)$.
\end{enumerate}
In this case, the limit points of $\{\pi_{n_j}\} $ are all the composition series of $\Ind_P^G(\rho)$.
\end{thm}

This theorem implies that non-isolated unitary representations are complementary series or limits of complementary series.

\subsection{The unitary dual conjecture}

In this subsection, we recall the unitary dual conjecture made by the first and second named authors, Hazeltine, Jiang, and Zhang. 

\begin{conj}[{\cite[Conjecture 5.9]{HJLLZ24}}]\label{unitary dual conjecture}
    The unitary dual of $G_n$ can be described as follows.
\begin{equation}
    \Pi_{\overline{A}}^{\lim}(G_n)=\Pi_{\ol{A}}^{\lim'}(G_n)= \Pi_{u}(G_n).
\end{equation}
\end{conj}

We remark that motivations of this unitary dual conjecture include the results of Tadi{\'c} and Vogan in Theorem \ref{vogan thm} above (i.e., non-isolated unitary representations are complementary series or limits of complementary series) and a conjecture of Tadi{\'c} on isolated representations $($\cite[Conjecture 1.1]{Tad22}$)$ (i.e., isolated representations are of Arthur type, hence in $\Pi_{\overline{A}}^{(0)}(G_n)=\Pi_A(G_n)$).

Note that in \cite[Algorithm 8.5]{HJLLZ24}, the authors developed an algorithm to generate the unitary dual candidate set $\Pi_{\overline{A}}^{\lim}(G_n)$, which will be recalled in Algorithm \ref{alg A bar} below.

\section{Classification of corank 3 tempered representations of good parity}\label{classtempcorank3}

As motivated by Definition \ref{def closure A}, the first step to construct the unitary dual of corank $4$ would be to construct the Arthur dual of corank $4$. To do so, the first natural step is to classify all representations of Arthur type that are of good parity. 
In this section, we classify all corank $3$ tempered representations of good parity, which will be used in the later sections on the classification of corank 4 non-tempered representations of good parity. 

 Let us first recall the notation in \cite[Section $7$]{HJLLZ24} regarding tempered representations of good parity, that are in particular not supercuspidal. From now on, we let $m_{\phi}(\rho \otimes S_x)$ be the multiplicity of $\rho \otimes S_x$ inside the $L$-parameter $\phi$. Recall that $\Psi(\pi)$ is the set of A-parameters $\psi$ such that $\pi \in \Pi_{\psi}$. 

\begin{thm}[{\cite[Theorem $7.1$]{HJLLZ24}}]\label{thm temp algo}
Let $\pi=\pi(\phi,\varepsilon)$ be a tempered representation of good parity. The representation $\pi$ is not supercuspidal if and only if at least one of the following holds for some $\rho$ and $x \in \half{1}\Z_{\geq 0}$. Denote $m:= m_{\phi}(\rho \otimes S_{2x+1})$ for short.
\begin{enumerate}
    \item [(I)] Suppose that $x>0$, $m_{\phi}(\rho\otimes S_{2x+1})>0$ and $ \varepsilon(\rho\otimes S_{2x+1})\varepsilon(\rho\otimes S_{2x-1})\neq -1.$ In this case, there exists a unique tempered representation $\pi_{temp}$ of smaller rank such that
   \[\pi \hookrightarrow \underbrace{\rho\lvert \cdot \rvert^{x}\times \cdots \times \rho\lvert \cdot \rvert^{x}}_{m \text{ copies}} \rtimes \pi_{temp} , \]
    and we write $\pi= \Temp{I}{x}{m}(\pi_{temp}).$ More precisely, we have $\pi_{temp}=\pi(\phi',\varepsilon')$ where 
    \[ \phi'= \phi- (\rho \otimes S_{2x+1})^{\oplus m}+(\rho \otimes S_{2x-1})^{\oplus m},\]
    and 
  \[ \varepsilon'(\rho'\otimes S_{a})=\begin{cases}
    \varepsilon(\rho\otimes S_{2x+1}) & \text{ if }\rho'\otimes S_{a} \cong \rho \otimes S_{2x-1},\\
      0 & \text{ if }\rho'\otimes S_{a} \cong \rho \otimes S_{2x+1},\\
    \varepsilon(\rho'\otimes S_{a}) & \text{ otherwise.}\end{cases}\]

    \item [(II)]Suppose that $x>0$, $m_{\phi}(\rho\otimes S_{2x+1})>1$ is \textbf{odd} and
    $ \varepsilon(\rho\otimes S_{2x+1})\varepsilon(\rho\otimes S_{2x-1})=-1 .$ In this case, there exists a unique tempered representation $\pi_{temp}$ of smaller rank such that
     \[\pi \hookrightarrow \underbrace{\rho\lvert \cdot \rvert^{x}\times \cdots \times \rho\lvert \cdot \rvert^{x}}_{m -1\text{ copies}} \rtimes \pi_{temp}, \]
    and we write $\pi= \Temp{II}{x}{m}(\pi_{temp}).$ More precisely, we have $\pi_{temp}=\pi(\phi',\varepsilon')$ where 
    \[ \phi'= \phi- (\rho \otimes S_{2x+1})^{\oplus m-1}+(\rho \otimes S_{2x-1})^{\oplus m-1},\]
    and $\varepsilon'(\rho'\otimes S_{a})=\varepsilon(\rho'\otimes S_{a})$.
    \item[(III)] Suppose that $x>0$, $m_{\phi}(\rho\otimes S_{2x+1})>1$ is \textbf{even} and $ \varepsilon(\rho\otimes S_{2x+1})\varepsilon(\rho\otimes S_{2x-1})=-1 .$ In this case, there exists a unique tempered representation $\pi_{temp}$ of smaller rank such that
    \begin{align*}
        \pi \hookrightarrow \underbrace{\rho\lvert \cdot \rvert^{x}\times \cdots \times \rho\lvert \cdot \rvert^{x}}_{m -1\text{ copies}} \times \rho\lvert \cdot \rvert^{x-1} \times \rho\lvert \cdot \rvert^{x-2} \times \cdots \times \rho\lvert \cdot \rvert^{-x} \rtimes \pi_{temp} 
    \end{align*}
    and we write $\pi= \Temp{III}{x}{m}(\pi_{temp}).$ More precisely, we have that $\pi_{temp}=\pi(\phi',\varepsilon')$ where 
    \[ \phi'= \phi- (\rho \otimes S_{2x+1})^{\oplus m}+(\rho \otimes S_{2x-1})^{\oplus m-2},\]
    and $\varepsilon'(\rho'\otimes S_{a})=\varepsilon(\rho'\otimes S_{a})$.
    \item [(IV)] Suppose that $x=0$, $m_{\phi}(\rho\otimes S_1)>1$ is odd. In this case, there exists a unique tempered representation $\pi_{temp}$ of smaller rank such that
    \[\pi \hookrightarrow \underbrace{\rho\lvert \cdot \rvert^{x}\times \cdots \times \rho\lvert \cdot \rvert^{x}}_{(m-1)/2 \text{ copies}} \rtimes \pi_{temp}, \] 
    and we write $\pi= \Temp{IV}{}{m}(\pi_{temp}).$ More precisely, we have that $\pi_{temp}=\pi(\phi',\varepsilon')$ where 
    \[ \phi'= \phi- (\rho \otimes S_{2x+1})^{\oplus m-1},\]
    and $\varepsilon'(\rho'\otimes S_{a})=\varepsilon(\rho'\otimes S_{a})$.
    \item [(V)] Suppose that $x=0$, $m_{\phi}(\rho\otimes S_1)>1$ is even.  In this case, there exists a unique tempered representation $\pi_{temp}$ of smaller rank such that
    \[\pi \hookrightarrow \underbrace{\rho\lvert \cdot \rvert^{x}\times \cdots \times \rho\lvert \cdot \rvert^{x}}_{m/2 \text{ copies}} \rtimes \pi_{temp}, \] 
    and we write $\pi= \Temp{V}{\varepsilon(\rho\otimes S_{1})}{m}(\pi_{temp}).$ More precisely, we have that $\pi_{temp}=\pi(\phi',\varepsilon')$ where 
    \[ \phi'= \phi- (\rho \otimes S_{2x+1})^{\oplus m},\]
    and 
     \[ \varepsilon'(\rho'\otimes S_{a})=\begin{cases}0 & \text{ if }\rho'\otimes S_{a} \cong \rho \otimes S_{1},\\
  \varepsilon(\rho'\otimes S_{a}) & \text{ otherwise.}\end{cases}\] 
\end{enumerate}
\end{thm}

The five cases above give us a way of describing all tempered representations. By definition, they are only well-defined up to certain restrictions, which we summarize below. 

\begin{remark}[{\cite[Remark $7.3$]{HJLLZ24}}]\label{rmk well-defined for Temp}
    Let $\pi_{temp}= \pi(\phi', \varepsilon')$.
    \begin{enumerate}
        \item [$\oldbullet$] The representation $\Temp{I}{x}{m}(\pi_{temp})$ is well-defined if and only if $m_{\phi'}(\rho\otimes S_{2x+1})=0$ and $m_{\phi'}(\rho\otimes S_{2x-1})\geq m$. 
         \item [$\oldbullet$] The representation $\Temp{II}{x}{m}(\pi_{temp})$ is well-defined if and only if $m_{\phi'}(\rho\otimes S_{2x+1})=1$, $m_{\phi'}(\rho\otimes S_{2x-1})\geq m$, and $\varepsilon'(\rho\otimes S_{2x+1})\varepsilon'(\rho\otimes S_{2x-1})=-1$.
        \item [$\oldbullet$] The representation $\Temp{III}{x}{m}(\pi_{temp})$ is well-defined if and only if $m_{\phi'}(\rho\otimes S_{2x+1})=0$ and $m_{\phi'}(\rho\otimes S_{2x-1})\geq m-1$.
         \item [$\oldbullet$] The representation $\Temp{IV}{}{m}(\pi_{temp})$ is well-defined if and only if $m_{\phi'}(\rho\otimes S_{1})=1$.
         \item [$\oldbullet$] The representation $\Temp{V}{\epsilon}{m}(\pi_{temp})$ is well-defined if and only if $m_{\phi'}(\rho\otimes S_{1})=0$.
    \end{enumerate}
\end{remark}

With these tools we can classify tempered representations of good parity of arbitrary corank. Fix some $\rho \in \mathcal{C}$. Suppose $\pi \in \Pi_{A, gp}(G_n)$ is tempered of corank $3$. Then we split into three cases: 
\begin{enumerate}[label = (\Alph*)]
    \item $\pi \hookrightarrow \rho \lvert \cdot \rvert^{x_1} \rtimes\pi_{temp}$, where $\pi_{temp}$ is tempered of corank $2$. Then $\pi$ is of the form $T_{I, 1}^{x}(\pi_{temp}), T_{IV,3}(\pi_{temp})$, or, $T_{V, 2}^{\pm}(\pi_{temp})$. 
    \item $\pi \hookrightarrow \rho\lvert \cdot \rvert^{x_1} \times \rho\lvert \cdot \rvert^{x_2} \rtimes \pi_{temp}$, where $\pi_{temp}$ is tempered of corank $1$. Then $\pi$ is of the form $T_{I, 2}^{x}(\pi_{temp}), T_{II, 3}^{x}(\pi_{temp}), T_{III, 2}^{\frac{1}{2}}(\pi_{temp}), T_{IV, 5}(\pi_{temp})$, or, $T_{V, 4}^{\pm}(\pi_{temp})$. 
    \item $\pi \hookrightarrow \rho\lvert \cdot \rvert^{x_1} \times \rho\lvert \cdot \rvert^{x_2} \times \rho\lvert \cdot \rvert^{x_3} \rtimes \pi_{sc}$, where $\pi_{sc}$ is supercuspidal. Then $\pi$ is of the form $T_{I, 3}^{x}(\pi_{sc}), T_{III, 2}^{1}(\pi_{sc}), T_{IV, 7}(\pi_{sc})$, or, $T_{V, 6}^{\pm}(\pi_{sc})$. 
\end{enumerate}

We begin with Case $(A)$. For the rest of this section,  we let $\pi_{sc} = \pi(\phi, \epsilon)$. where $\phi = \oplus_{\rho}\oplus_{i=\epsilon_{\rho}}^{\alpha-1}\rho \otimes S_{2i+1}$. Then $\pi_{sc}$ is self-dual supercuspidal with reducibility point $\alpha = \alpha_{\rho, \pi_{sc}}$.

\subsection{\texorpdfstring{Case $A:\pi \hookrightarrow \rho \lvert \cdot \rvert^{x_1} \rtimes\pi_{temp}$, where $\pi_{temp}$ is tempered of corank $2$}{}}\label{sec 3.1}
 By \cite[\S 10]{HJLLZ24}, the followings are the possible representations $\pi_{temp}$ of corank $2$: 
\begin{flalign*}
   & T_{I, 1}^{\alpha +1}(T_{I, 1}^{\alpha}(\pi_{sc})), T_{I, 1}^{\alpha-1}(T_{I, 1}^{\alpha}(\pi_{sc})), T_{IV, 3}(T_{I, 1}^{\alpha}(\pi_{sc})), \\
   & T_{V, 2}^{\pm}(T_{I, 1}^{1}(\pi_{sc})), T_{I, 1}^{1}(T_{V, 2}^{\pm}(\pi_{sc})), T_{I, 2}^{\frac{1}{2}}(\pi_{sc}), T_{II, 3}^{\frac{1}{2}}(\pi_{sc}), \\
   &T_{III, 2}^{\frac{1}{2}}(\pi_{sc}), T_{IV, 5}(\pi_{sc}), T_{V, 4}^{\pm}(\pi_{sc}). 
\end{flalign*}
There are $10$ possibilities, so there are $30$ cases in total in Case $(A)$. First we consider the case Let $\pi_{temp} = T_{I, 1}^{\alpha+1}(T_{I, 1}^{\alpha}(\pi_{sc}))$. From the classification in \cite[Proposition $10.1$]{HJLLZ24}, this is well-defined only when $\alpha > 0$. 

\begin{prop}\label{temp3,1}
Let $\pi_{temp} = T_{I, 1}^{\alpha+1}(T_{I, 1}^{\alpha}(\pi_{sc}))$ for $\alpha > 0$. 
\begin{enumerate}[label = (\roman*)]
    \item The representation $T_{I, 1}^{x}(\pi_{temp})$ is well-defined if and only if $x = \alpha +2$ for $\alpha > 0$, or $x = \alpha -1$ for $\alpha > 1$. 
    
    \item The representation $T_{I, 1}^x(\pi_{temp})$ is of critical type when $x = \alpha +2$ for $\alpha > 0$, or $x = \alpha -1$ for $\alpha > 1$.
    \item Let $x \in \{\alpha +2, \alpha -1\}$, and define 
    \begin{equation*}
        \mathcal{E}_{\alpha+2} := \{([\alpha -2, \epsilon_{\rho}]_{\rho}, 0 , \eta), ([\alpha+2, \alpha+2], 0, (-1)^{\alpha-1-\epsilon_{\rho}}\eta)\},
    \end{equation*}
    \begin{equation*}
        \mathcal{E}_{\alpha-1} :=  \{([\alpha -3, \epsilon_{\rho}]_{\rho}, 0 , \eta), ([\alpha-1, \alpha-1], 0, *), ([\alpha+1, \alpha +1],0,(-1)^{\alpha-1-\epsilon_{\rho}}\eta)\}.
    \end{equation*}
    Then we have $\pi(\mathcal{E}_x) = T_{I, 1}^{x}(\pi_{temp})$. Here are the associated symbols:
    
\[\EE_{\alpha+2}= \scalebox{0.8}{\bordermatrix{
  & \epsilon_{\rho} & \cdots & \alpha-2 &\alpha-1& \alpha & \alpha+1 & \alpha+2 \cr
  & \odot & \cdots & \odot &&&&\cr
  &&&&&&&\odot
}},\]
\[\EE_{\alpha-1}= \scalebox{0.8}{\bordermatrix{
  & \epsilon_{\rho} & \cdots & \alpha-3 &\alpha-2& \alpha-1 & \alpha & \alpha+1\cr
  & \odot & \cdots & \odot &&&&\cr
  &&&&&\odot&&\cr
  &&&&&&&\odot
}}.\]
\end{enumerate}
\end{prop}
\begin{proof} Part $(i)$ follows directly from Remark \ref{rmk well-defined for Temp}. The condition $\alpha > 1$ for $x = \alpha -1$ follow from the definition of $T_{I,1}^{x}$. Parts $(ii)$ and $(iii)$ follow from definition. \end{proof}.

\begin{prop}\label{temp3,2}
    Let $\pi_{temp} = T_{I, 1}^{\alpha+1}(T_{I, 1}^{\alpha}(\pi_{sc}))$ for $\alpha > 0$. 
    \begin{enumerate}[label = (\roman*)]
        \item The representation $T_{IV,3}(\pi_{temp})$ is well-defined if and only if $\alpha \in \mathbb{Z}_{>1}$.
        \item The representation $T_{IV, 3}(\pi_{temp})$ is not of critical type. 
        \item Define $\EE = \{([0,0]_{\rho},0,\eta),([0,0]_{\rho},0,\eta),([\alpha-2,0]_{\rho},0,\eta),([\alpha+1, \alpha+1]_{\rho}, 0,(-1)^{\alpha-1}\eta)\}$. Then $\pi(\EE) = T_{IV, 3}(\pi_{temp})$. Here is the associated symbol: 
        \[\EE= \scalebox{0.8}{\bordermatrix{ 
        &0 & \cdots &\alpha-2& \alpha-1 & \alpha &\alpha+1 \cr
        &\odot &&&&&\cr
        &\odot &&&&&\cr
        & \odot & \cdots& \odot &&&\cr 
        &&&&&&\odot 
        }}.\]
    \end{enumerate}
\end{prop}
\begin{proof} Let $\phi$ be the $L$-parameter associated with $\pi_{temp}$. By Remark \ref{rmk well-defined for Temp}, $T_{IV,3}(\pi_{temp})$ is well-defined if and only if $m_{\phi}(\rho\otimes S_1) = 1$, which happens only when $\alpha \neq 0,1$. This proves part $(i)$. Parts $(ii)$ and $(iii)$ follow from the definition. \end{proof} 

The proofs of Propositions \ref{temp3,3} to \ref{temp3,6} below are similar to that of Proposition \ref{temp3,2}, which we omit.

\begin{prop}\label{temp3,3}
    Let $\pi_{temp} = T_{I, 1}^{\alpha+1}(T_{I, 1}^{\alpha}(\pi_{sc}))$ for $\alpha > 0$. 
    \begin{enumerate}[label = (\roman*)]
        \item The representation $T_{V, 2}^{\pm}(\pi_{temp})$ is well-defined if and only if $\alpha =1$.
        \item When $\alpha =1$, the representation $T_{V, 2}^{\pm}(\pi_{temp})$ is of critical type. 
        \item Define 
        \begin{equation*}
            \EE_{\pm} = \{([0,0],0,\pm1))^2,([2,2],0,\eta)\}.
        \end{equation*}
        Then $\pi(\EE_{\pm}) = T_{V, 2}^{\pm}(\pi_{temp})$. Here are the associated symbols. 
        \[\EE_{+}= \scalebox{0.8}{\bordermatrix{ 
        &0&1&2\cr
        &\oplus&&\cr
        &\oplus&&\cr
        &&&\odot
        }}, \quad \EE_{-}= \scalebox{0.8}{\bordermatrix{ 
        &0&1&2\cr
        &\ominus&&\cr
        &\ominus&&\cr
        &&&\odot
        }}.\]         
    \end{enumerate}
\end{prop}

Now let us consider the case $\pi_{temp} = T_{I, 1}^{\alpha-1}(T_{I, 1}^{\alpha}(\pi_{sc}))$. By the classification in \cite[Proposition $10.1$]{HJLLZ24}, this is well-defined only when $\alpha > 1$. 
\begin{prop}\label{temp3,4}
    Let $\pi_{temp} =T_{I, 1}^{\alpha-1}(T_{I, 1}^{\alpha}(\pi_{sc}))$ and $\alpha > 1$. 
    \begin{enumerate}[label = (\roman*)]
        \item The representation $T_{I,1}^x(\pi_{temp})$ is well defined if and only if $x = \alpha -2$ and $\alpha > 2$, or $x = \alpha +1$. 
        \item The representation $T_{I, 1}^{x}(\pi_{temp})$ is of critical type, when $x = \alpha -2, \alpha > 2$, or $x = \alpha+1$.  When $x = \alpha+1$, the representation $T_{I,1}^{x}(\pi_{temp})$ is the same as the representation $T_{I,1}^{\alpha-1}(T_{I,1}^{\alpha+1}(T_{I,1}^{\alpha}(\pi_{sc}))$ described in Proposition \ref{temp3,1}.  
        \item  Define 
        \begin{equation*}
            \EE := \{([\alpha -4, \epsilon_{\rho}]_{\rho},0, \eta),([\alpha,\alpha-2]_{\rho},0,-\eta\}.
        \end{equation*}
        Then $\pi(\EE) = T_{I, 1}^{\alpha -2}(\pi_{temp})$. Here is the associated symbol. 
        \[\EE= \scalebox{0.8}{\bordermatrix{ 
        &\epsilon_{\rho}&\cdots&\alpha-4&\alpha-3&\alpha-2&\alpha-1&\alpha\cr
        &\odot&\cdots&\odot&&&&\cr
        &&&&&\odot&\odot&\odot
        }}.\]
    \end{enumerate}
\end{prop}

\begin{prop}\label{temp3,5}
    Let $\pi_{temp} =T_{I, 1}^{\alpha-1}(T_{I, 1}^{\alpha}(\pi_{sc}))$ and $\alpha > 1$. 
    \begin{enumerate}[label = (\roman*)]
        \item The representation $T_{IV,3}(\pi_{temp})$ is well-defined if and only if $\alpha \in \mathbb{Z}_{>2}$. 
        \item The representation $T_{IV,3}(\pi_{temp})$ is not of critical type. 
        \item Define 
        \begin{equation*}
            \EE = \{([0,0]_\rho,0,\eta),([0,0]_\rho,0,\eta),([\alpha-3, 0]_\rho,0,\eta),([\alpha,\alpha-1]_{\rho},0,(-1)^{\alpha-2}\eta)\}.
        \end{equation*}
        Then $\pi(\EE)=T_{IV,3}(\pi_{temp})$. Here are the associated symbols. 
        \[\EE= \scalebox{0.8}{\bordermatrix{ 
        &0&\cdots&\alpha-3&\alpha-2&\alpha-1&\alpha \cr
        &\odot&&&&&\cr
        &\odot&&&&&\cr
        &\odot&\cdots&\odot&&&\cr
        &&&&&\odot&\odot
        }}.\]
    \end{enumerate}
\end{prop}

\begin{prop}\label{temp3,6}
    Let $\pi_{temp} =T_{I, 1}^{\alpha-1}(T_{I, 1}^{\alpha}(\pi_{sc}))$ and $\alpha > 1$. 
    \begin{enumerate}[label = (\roman*)]
    \item The representation $T_{V,2}^{\pm}(\pi_{temp})$ is well-defined if and only if $\alpha =2$. 
    \item When $\alpha =2$, the representation $T_{V, 2}^{\pm}(\pi_{temp})$ is of critical type. 
    \item Define
    \begin{equation*}
        \EE_{\pm} = \{([0,0],0,\pm1)^2, ([2,1],0,\eta)\}.
    \end{equation*}
    Then $\pi(\EE_{\pm}) = T_{V,2}(\pi_{temp})$. Here are the associated symbols. 
    \[\EE_{+}= \scalebox{0.8}{\bordermatrix{ 
        &0&1&2\cr
        &\oplus&&\cr
        &\oplus&&\cr
        &&\odot&\odot
        }}, \quad \EE_{-}= \scalebox{0.8}{\bordermatrix{ 
        &0&1&2\cr
        &\ominus&&\cr
        &\ominus&&\cr
        &&\odot&\odot
        }}.\]    
    \end{enumerate}
\end{prop}
\begin{proof} The proof is similar to that of Proposition \ref{temp3,2}, which we omit. \end{proof} 

 Now we consider the case $\pi_{temp} = T_{IV, 3}(T_{I, 1}^{\alpha}(\pi_{sc}))$. By the classification in \cite[Proposition $10.2$]{HJLLZ24}, this is well-defined only for $\alpha \in \mathbb{Z}_{>1}$. 
\begin{prop}\label{temp3,7}
    Let $\pi_{temp} = T_{IV, 3}(T_{I, 1}^{\alpha}(\pi_{sc}))$ for $\alpha \in \mathbb{Z}_{> 1}$. 
    \begin{enumerate}[label = (\roman*)]
        \item The representation $T_{I, 1}^{x}(\pi_{temp})$ is well-defined if and only if $x = \alpha+1$ or $\alpha-1$. 
        \item When $x = \alpha+1$, the representation $T_{I, 1}^{x}(\pi_{temp})$ is not of critical type and is the same as the representation $T_{IV,3}(T_{I,1}^{\alpha+1}(T_{I,1}^{\alpha}(\pi_{sc})))$ described in Proposition \ref{temp3,2}.  If $x = \alpha-1$, and $\alpha \neq 2$, then the representation $T_{I,1}^{x}(\pi_{temp})$ is not of critical type and $T_{I,1}^{x}(\pi_{temp}) = T_{IV,3}(T_{I,1}^{\alpha-1}(T_{I,1}^{\alpha}(\pi_{sc})))$ described in Proposition \ref{temp3,5}.  
        \item  When $\alpha = 2$, the representation $\pi(\EE) = T_{I,1}^{\alpha-1}(\pi_{temp})$ is of critical type, where 
        \[\EE= \scalebox{0.8}{\bordermatrix{ 
        &0&1&2 \cr
        &\odot&&&\cr
        &\odot&&&\cr
        &&\odot&\odot
        }}.\]
        \item The representations $T_{IV,3}(\pi_{temp})$ and $T_{V,2}^{\pm}(\pi_{temp})$ are not well-defined. 
    \end{enumerate}
\end{prop}
\begin{proof} The proof of parts $(i)$ through $(iii)$ is similar to that of Proposition \ref{temp3,1}, which we omit. Part $(iv)$ follows from the fact that $m_{\phi}(\rho \otimes S_1) = 2$, where $\phi$ is the $L$-parameter associated with $\pi_{temp}$.   \end{proof} 

The next case is $\pi_{temp} = T_{V, 2}^{\pm}(T_{I,1}^{1}(\pi_{sc}))$. By \cite[Proposition $10.3$]{HJLLZ24}, this is well-defined if and only if $\alpha = 1$. The proofs of Propositions \ref{temp3,8} to \ref{temp3,10} below are similar to that of Proposition \ref{temp3,7}, which we omit. 

\begin{prop}\label{temp3,8}
    Let $\pi_{temp}^{\pm} = T_{V, 2}^{\pm}(T_{I,1}^{1}(\pi_{sc}))$ and $\alpha =1$.
    \begin{enumerate}[label = (\roman*)]
        \item The representation $T_{I, 1}^{x}(\pi_{temp}^{\pm})$ is well-defined if and only if $x = 2$. 
        \item The representation $T_{I,1}^{2}(\pi_{temp}^{\pm})$ is of critical type. 
        \item Define
        \begin{equation*}
            \EE_{\pm} = \{([0,0]_\rho,0,\pm1),([0,0]_\rho,0,\pm1),([2,2]_\rho,0,\eta)\}.
        \end{equation*}
        Then $\pi(\EE_{\pm}) = T_{I,1}^{2}(\pi_{temp}^{\pm})$. Here are the associated symbols. 
            \[\EE_{+}= \scalebox{0.8}{\bordermatrix{ 
        &0&1&2\cr
        &\oplus&&\cr
        &\oplus&&\cr
        &&&\odot
        }},\]\[\EE_{-}= \scalebox{0.8}{\bordermatrix{ 
        &0&1&2\cr
        &\ominus&&\cr
        &\ominus&&\cr
        &&&\odot
        }}.\]    
        \item The representations $T_{IV,3}(\pi_{temp})$ and $T_{V,2}^{\pm}(\pi_{temp})$ are not well-defined. 
    \end{enumerate}
\end{prop}

The next case to consider is $\pi_{temp} = T_{I, 1}^{1}(T_{V,2}^{\pm}(\pi_{sc}))$. By the classification in \cite[Proposition $10.5$]{HJLLZ24}, this is well-defined only when $\alpha = 0$.

    \begin{prop}\label{temp3,9}
        Let $\pi_{temp}^{\pm} = T_{I,1}^{1}(T_{V,2}^{\pm})(\pi_{sc})$ and $\alpha = 0$. 
        \begin{enumerate}[label = (\roman*)]
            \item The representation $T_{I,1}^{x}(\pi_{temp}^{\pm})$ is well-defined if and only if $x = 2$. 
            \item When $x =2$, the representation $T_{I,1}^{x}(\pi_{temp}^{\pm})$ is of critical type. 
            \item Define
            \begin{equation*}
                \EE_{\pm} = \{([0,0]_\rho,0,\pm1),([2,2]_\rho,0,\pm1)\},
            \end{equation*}
            then $\pi(\EE_{\pm}) = T_{I,1}^{2}(\pi_{temp}^{\pm})$. Here are the associated symbols. 
            \[\EE_{+}= \scalebox{0.8}{\bordermatrix{ 
        &0&1&2\cr
        &\oplus&&\cr
        &&&\oplus
        }}, \quad \EE_{-}= \scalebox{0.8}{\bordermatrix{ 
        &0&1&2\cr
        &\ominus&&\cr
        &&&\ominus\cr
        }}.\]    
        \end{enumerate}
    \end{prop}

    \begin{prop}\label{temp3,10}
        Let $\pi_{temp}^{\pm} = T_{I,1}^{1}(T_{V,2}^{\pm}(\pi_{sc}))$ and $\alpha = 0$.
        \begin{enumerate}[label = (\roman*)]
            \item The representation $T_{IV,3}(\pi_{temp})$ is well-defined. 
            \item The representation $T_{IV,3}(\pi_{temp})$ is of critical type. 
            \item Define
            \begin{equation*}
                \EE_{\pm}:= ([0,0]_\rho,0,\pm1)^3,([1,1]_\rho,0,\pm1)\}.
            \end{equation*}
            Then $\pi(\EE_{\pm})= T_{IV,3}(\pi_{temp})$. Here are the associated symbols. 
             \[\EE_{+}= \scalebox{0.8}{\bordermatrix{ 
        &0&1\cr
        &\oplus&\cr
        &\oplus&\cr
        &\oplus&\cr
        &&\oplus
        }}, \EE_{-}= \scalebox{0.8}{\bordermatrix{ 
        &0&1\cr
        &\ominus\cr
        &\ominus\cr
        &\ominus\cr
        &&\ominus\cr
        }}.\] 
        \item The representation $T_{V,2}^{\pm}(\pi_{temp})$ is not well-defined. 
        \end{enumerate}
    \end{prop}

    The next case to consider is $\pi_{temp} = T_{I,2}^{\frac{1}{2}}(\pi_{sc})$. By \cite[Proposition $10.6$]{HJLLZ24}, this is well-defined only when $\alpha = \frac{1}{2}$. 
    
    \begin{prop}\label{temp3,11}
        Let $\pi_{temp} = T_{I,2}^{\frac{1}{2}}(\pi_{sc})$ and $\alpha = \frac{1}{2}$. 
        \begin{enumerate}[label = (\roman*)]
            \item The representation $T_{I,1}^{x}(\pi_{temp})$ is well-defined only when $x = \frac{3}{2}$. 
            \item The representation $T_{I,1}^{\frac{3}{2}}(\pi_{temp})$ is of critical type. 
            \item The set $\Psi(T_{I,1}^{\frac{3}{2}}(\pi_{temp}))$ is a singleton, and we have $T_{I,1}^{\frac{3}{2}}(\pi_{temp}) = \pi(\EE)$, where 
            \[\EE= \scalebox{0.8}{\bordermatrix{ 
        &\frac{1}{2}&\frac{3}{2} \cr
        & \oplus \cr
        &&\oplus 
        }}.\] 
        \item The representations $T_{IV,3}(\pi_{temp})$ and $T_{V,2}^{\pm}(\pi_{temp})$ are not well-defined. 
        \end{enumerate}
    \end{prop}
    \begin{proof} The proof for parts $(i)$ and $(ii)$ is similar to that of Proposition \ref{temp3,1}, which we omit. Part $(iii)$ is easy to verify. For part $(iv)$, since $\alpha = \frac{1}{2}$, the good parity condition implies that $x \in \frac{1}{2} + \mathbb{Z}$, so by definition the two representations are not well-defined. \end{proof}
    
    Now we consider the case $\pi_{temp} = T_{II,3}^{\frac{1}{2}}(\pi_{sc})$. By \cite[Proposition $10.7$]{HJLLZ24}, this is well-defined if and only if $\alpha \in \frac{1}{2} + \mathbb{Z}_{> 0}$.
The proofs of Propositions \ref{temp3,12} to \ref{temp3,15} below are similar to that of Proposition \ref{temp3,11}, which we omit. 
    
    \begin{prop}\label{temp3,12}
        Let $\pi_{temp} = T_{II,3}^{\frac{1}{2}}(\pi_{sc})$ and let $\alpha \in \frac{1}{2} + \mathbb{Z}_{>0}$.
        \begin{enumerate}[label = (\roman*)]
            \item The representation $T_{I,1}^{x}(\pi_{temp})$ is well-defined if and only if $x = \alpha$. 
            \item The representation $T_{I,1}^{\alpha}(\pi_{temp})$ is of critical type if and only if $\alpha = \frac{3}{2}$. 
            \item Define 
            \begin{equation*}
                \EE = \{([1/2,1/2]_\rho,0,-1),([1/2,1/2]_\rho,0,-1),([\alpha -2,1/2]_\rho,0,-1),([\alpha,\alpha],0,\eta)\}.
            \end{equation*}
            Then $\pi(\EE) = T_{I,1}^{\alpha}(\pi_{temp})$. Here is the associated symbol for $\alpha > \frac{3}{2}$. 
            \[\EE= \scalebox{0.8}{\bordermatrix{ 
        &\frac{1}{2}&\cdots&\alpha-2&\alpha-1&\alpha\cr
        &\ominus&&&&\cr
        &\ominus&&&&\cr
        &\ominus&\cdots&\odot&&\cr
        &&&&&\odot
        }}.\] For $\alpha = \frac{3}{2}$ the associated symbol is
        \[\EE= \scalebox{0.8}{\bordermatrix{ 
        &\frac{1}{2}&\frac{2}{3} \cr
        &\ominus&\cr
        &\ominus&\cr
        &&\ominus\cr
        }}.\]
        \item The representations $T_{IV,3}(\pi_{temp})$ and $T_{V,2}(\pi_{temp})$ are not well-defined. 
        \end{enumerate}
    \end{prop}
    The next case to consider is $\pi_{temp} = T_{III, 2}^{\frac{1}{2}}(\pi_{sc})$. By \cite[Proposition 10.8]{HJLLZ24}, this is well-defined only when  $\alpha = \frac{1}{2}$. 
    
    \begin{prop}\label{temp3,13}
        Let $\pi_{temp} = T_{III, 2}^{\frac{1}{2}}(\pi_{sc})$ and $\alpha = \frac{1}{2}$. 
        \begin{enumerate}[label = (\roman*)]
            \item The representation $T_{I,1}^{x}(\pi_{temp})$ is well-defined if and only if $x = \frac{3}{2}$. 
            \item The representation $T_{I,1}^{\frac{3}{2}}(\pi_{temp})$ is of critical type. 
            \item The following is a complete list of extended multi-segments $\EE$ (up to row exchanges) such that $\pi(\EE) = T_{I,1}^{\frac{3}{2}}(\pi_{temp})$:   
            \[\Bigg\{\scalebox{0.8}{\bordermatrix{ 
        &\frac{1}{2}&\frac{3}{2} \cr
        &\ominus&\cr
        &&\ominus\cr
        }},{\bordermatrix{ 
        &-\frac{1}{2}&\frac{1}{2}&\frac{3}{2} \cr
        &\oplus&\ominus&\cr
        &&&\ominus\cr}}\Bigg\}.\]
        \item The representations $T_{IV,3}(\pi_{temp})$ and $T_{V,2}^{\pm}(\pi_{temp})$ are not well-defined. 
        \end{enumerate}
    \end{prop}

    Now we consider the case $\pi_{temp} = T_{IV,5}(\pi_{sc})$. From \cite[Proposition $10.9$]{HJLLZ24}, this is well-defined if and only if $\alpha \in \mathbb{Z}_{>0}$. 
    
    \begin{prop}\label{temp3,14}
        Let $\pi_{temp} = T_{IV,5}(\pi_{sc})$ and $\alpha \in \mathbb{Z}_{> 0}$. 
        \begin{enumerate}[label = (\roman*)]
            \item The representation $T_{I,1}^{x}(\pi_{temp})$ is well-defined if and only if $x = \alpha$. 
            \item The representation $T_{I,1}^{\alpha}(\pi_{temp})$ is of critical type if and only if $\alpha = 1$. 
            \item Define
            \begin{equation*}
                \EE = \{([0,0]_\rho,0,\eta)^4,\} + \{([\alpha -2,0]_\rho,0,\eta), ([\alpha,\alpha],0, -\eta)\}. 
            \end{equation*}
            Then $\pi(\EE) = T_{I,1}^{\alpha}(\pi_{temp})$. Here is the associated symbol. 
            \[\EE= \scalebox{0.8}{\bordermatrix{ 
        &0&\dots&\alpha-2&\alpha-1&\alpha \cr
        &\odot&&&&\cr
        &\odot&&&&\cr
        &\odot&&&&\cr
        &\odot&&&&\cr
        &\odot&\cdots&\odot&&\cr
        &&&&&\odot
        }}.\]
        \item The representations $T_{IV,3}(\pi_{temp})$ and $T_{V,2}^{\pm} (\pi_{temp})$ are not well-defined. 
        \end{enumerate}
    \end{prop}

    The final case in case (A) is $\pi_{temp} = T_{V,4}^{\pm}(\pi_{sc})$, which is well-defined only when $\alpha = 0$ by \cite[Proposition $10.10$]{HJLLZ24}. 
    
    \begin{prop}\label{temp3,15}
        Let $\pi_{temp}^{\pm} = T_{V,4}^{\pm}(\pi_{sc})$, and $\alpha = 0$. 
        \begin{enumerate}[label = (\roman*)]
            \item The representation $T_{I,1}^{x}(\pi_{temp}^{\pm})$ is well-defined if and only if $x = 1$. 
            \item The representation $T_{I,1}^{1}(\pi_{temp}^{\pm})$ is of critical type. 
            \item The set $\Psi(T_{V,4}^{\pm}(\pi_{sc}))$ is a singleton, and $\pi(\EE_{\pm}) = T_{I,1}^{1}(\pi_{temp}^{\pm})$, where 
            \[\EE_{+}= \scalebox{0.8}{\bordermatrix{ 
        &0&1\cr
        &\oplus&\cr
        &\oplus&\cr
        &\oplus&\cr
        &&\oplus
        }}, \quad \EE_{-}= \scalebox{0.8}{\bordermatrix{ 
        &0&1\cr
        &\ominus&\cr
        &\ominus&\cr
        &\ominus&\cr
        &&\ominus
        }}.\]    
        \item The representations $T_{IV,3}(\pi_{temp}^{\pm})$ and $T_{V,2}(\pi_{temp})^{\pm}$ are not well-defined. 
        \end{enumerate}
    \end{prop}
    
    This concludes our discussion of Case (A). Now we move on to Case $(B)$. 

    \subsection{\texorpdfstring{Case $B: \pi \hookrightarrow \rho\lvert \cdot \rvert^{x_1} \times \rho\lvert \cdot \rvert^{x_2} \rtimes \pi_{temp}$}{}}
    The three possible tempered representations of corank $1$ are $T_{I,1}^{\alpha}(\pi_{sc}), T_{IV, 3}(\pi_{sc})$ and $T_{V, 2}^{\pm}(\pi_{sc})$. First we consider the case $T_{I,1}^{\alpha}(\pi_{sc})$. By \cite[Proposition $8.1$]{HJLLZ24}, this is well-defined if and only if $\alpha > 0$. 
    
    \begin{prop}\label{temp3,16}
        Let $\pi_{temp} = T_{I,1}^{\alpha}(\pi_{sc})$, $\alpha > 0$. 
        \begin{enumerate}[label = (\roman*)]
            \item The representations $T_{I,2}^{x}(\pi_{temp})$ and $T_{II,3}^{x}(\pi_{temp})$ are not well-defined. 
            \item The representation $T_{III, 2}^{\frac{1}{2}}(\pi_{temp})$ is not well-defined. 
        \end{enumerate}
    \end{prop}
    \begin{proof} Part $(i)$ follows from the fact that the parameter $\phi$ corresponding to $\pi_{temp}$ is multiplicity free. To see part $(ii)$, suppose $T_{III,2}^{\frac{1}{2}}$ is well-defined, then we must have either $x = \alpha-1 = \frac{1}{2}$, or $x = \alpha +1 = \frac{1}{2}$. For the first case, we have $\alpha = \frac{3}{2}$, which means that the multiplicity of $\rho \otimes S_{2}$ in the parameter $\phi$ is $0$, so it is not well-defined. For the second case we have $\alpha < 0$, which is impossible. This completes the proof.   \end{proof}
    
    We check the remaining two possibilities for this tempered representation.  
    
    \begin{prop}\label{temp3,17}
        Let $\pi_{temp} = T_{I,1}^{\alpha}(\pi_{sc})$, $\alpha > 0$. 
        \begin{enumerate}[label = (\roman*)]
            \item The representation $T_{IV,5}(\pi_{temp})$ is-well-defined if and only if $\alpha \in \mathbb{Z}_{>0}$. It is the same as the representation $T_{I,1}^{\alpha}(T_{IV,5}(\pi_{sc}))$ described in Proposition \ref{temp3,14}. 
            \item The representation $T_{IV,5}(\pi_{temp})$ is of critical type if and only if $\alpha = 1$. 
        \end{enumerate}
    \end{prop}
    \begin{proof} Let $\phi$ be the $L$-parameter corresponding to $\pi_{sc}$. To ensure $m_{\phi}(\rho \otimes S_1) =1$, we need $\alpha \in \mathbb{Z}$. The rest follows from definition. 

    \begin{prop}\label{temp3,18}
        Let $\pi_{temp} = T_{I,1}^{\alpha}(\pi_{sc})$, $\alpha > 0$.
        \begin{enumerate}[label = (\roman*)]
            \item The representations $T_{V,4}^{\pm}(\pi_{temp})$ are well-defined if and only if $\alpha = 1$. 
            \item When $\alpha =1$, the representations $T_{V,4}^{\pm}(\pi_{temp})$ are of critical type. 
            \item When $\alpha =1$, define
            \begin{equation*}
                \EE_{\pm} := \{([0,0],0,\pm1)^4,  ([1,1],0,\pm1)\}.
            \end{equation*}
            Then $\pi(\EE_{\pm}) = T_{V,4}^{\pm}(\pi_{temp})$. Here are the associated symbols. 
            \[\EE= \scalebox{0.8}{\bordermatrix{ 
             &0&1\cr
        &\oplus&\cr
        &\oplus&\cr
        &\oplus&\cr
        &\oplus&\cr
        &&\oplus
        }}, \quad \EE_{-}= \scalebox{0.8}{\bordermatrix{ 
        &0&1\cr
        &\ominus&\cr
        &\ominus&\cr
        &\ominus&\cr
        &\ominus&\cr
        &&\ominus
        }}.\]    
        \end{enumerate}
    \end{prop}
    \begin{proof} To ensure $S(\rho \otimes S_1) = 1$, we must have $\alpha -1 = 0$ or $0 > \alpha$, this proves part $(i)$. The rest follows from definition. \end{proof}

    The next case to consider in Case $(B)$ is $\pi_{temp} = T_{IV, 3}(\pi_{sc})$, which is well-defined if and only if $\alpha \in \mathbb{Z}_{> 0}$ by \cite[Proposition $8.2$]{HJLLZ24}. 
    
    \begin{prop} \label{temp3,19}
        Let $\pi_{temp} = T_{IV,3}(\pi_{sc})$, and $\alpha \in \mathbb{Z}_{>0}$
        \begin{enumerate} [label = (\roman*)]
            \item The representations $T_{I,2}^{x}(\pi_{temp})$ is well-defined if and only if $\alpha = x = 1$.
            \item When $\alpha = x = 1$, the representation $T_{I,2}^{1}(\pi_{temp})$ is of critical type. 
            \item When $\alpha = x = 1$, define 
            \begin{equation*}
                \EE:= \{([0,0]_\rho,0,\eta),([1,1]_\rho,0,\eta)^2\}.
            \end{equation*}
            Then $\pi(\EE) = T_{I,2}^{1}(\pi_{temp})$. Here is the associated symbol. 
            \[\EE= \scalebox{0.8}{\bordermatrix{ 
             &0&1\cr
        &\odot&\cr
        &&\odot\cr
        &&\odot
        }}.\]
        \end{enumerate}    
    \end{prop}
    \begin{proof} Let $\phi$ be the L-parameter associated to $\pi_{temp}$.  To ensure $m_{\phi}(\rho \otimes S_{2(x-1)+1}) \geq m =2$, we must have $x = 1$. In this case if $\alpha \neq 1$, then $m_{\phi}(\rho \otimes S_{2x+1}) \neq 0$. This proves part $(i)$. The rest follows from definition. \end{proof}
    
    \begin{prop}\label{temp3,20}
        Let $ \pi_{temp} = T_{IV,3}(\pi_{sc})$ and $\alpha \in \mathbb{Z}_{>0}$. 
        \begin{enumerate}[label = (\roman*)]
            \item The representation $T_{II,3}^{x}(\pi_{temp})$ is well-defined if and only if $x = 1$. 
            \item The representation $T_{II,3}^{1}(\pi_{temp})$ is of critical type. 
            \item Define 
            \begin{equation*}
                \EE:= \{([\alpha-1,0]_\rho,0,\eta),([1,1]_\rho,0,-\eta),([1,1],0,-\eta)\}.
            \end{equation*}
            Then $\pi(\EE) = T_{II,3}^{1}(\pi_{temp})$. Here is the associated symbol. 
            \[\EE= \scalebox{0.8}{\bordermatrix{ 
             &0&1& \cdots & \alpha-1\cr
        &\odot&\odot & \cdots & \odot \cr
        &&\odot&&\cr
        &&\odot&&
        }}.\]
        \item The representations $T_{III,2}^{\frac{1}{2}}(\pi_{temp})$, $T_{IV,5}(\pi_{temp})$ and $T_{V, 4}^{\pm}(\pi_{temp})$ are not well-defined. 
        \end{enumerate}
    \end{prop}
    \begin{proof} Let $\phi$ be the L-parameter associated with $\pi_{temp}$. To ensure $m_{\phi}(\rho \otimes S_{2(x-1)+1}) \geq m = 3$, the only possibility is $x = 1$, in which case we can check it's well-defined. This proves part $(i)$. Parts $(ii)$ and $(iii)$ follow from definition. 
    For part $(iv)$, the good parity condition tells us that $T_{III, 2}^{\frac{1}{2}}$ is not well-defined, and since $m_{\phi}(\rho \otimes S_{1}) = 3$, $T_{IV,5}(\pi_{temp})$ and $T_{V, 4}^{\pm}(\pi_{temp})$ are also not well-defined. \end{proof}

    The last case to consider in Case $(B)$ is $\pi_{temp} = T_{V,4}^{\pm}(\pi_{sc})$. This is well-defined if and only if $\alpha = 0$ by \cite[Proposition $8.3$]{HJLLZ24}. 
    
    \begin{prop}\label{temp3,21}
        Let $\pi_{temp}^{\pm} = T_{V,2}^{\pm}(\pi_{sc})$ and $\alpha = 0$. 
        \begin{enumerate}[label = (\roman*)]
            \item The representations $T_{I,2}^{x}(\pi_{temp}^{\pm})$ are well-defined if and only if $x = 1$. 
            \item The representations $T_{I,2}^{1}(\pi_{temp}^{\pm})$ are of critical type. 
            \item Define 
            \begin{equation*}
                \EE_{\pm}:= \{([1,1]_\rho,0,\pm1)^2\}.
            \end{equation*}
            Then $\pi(\EE_{\pm}) = T_{I,1}^{1}(\pi_{temp}^{\pm})$. Here are the associated symbols. \[\EE_{+}= \scalebox{0.8}{\bordermatrix{ 
             &1\cr
             &\oplus\cr
             &\oplus
        }}, \quad
        \EE_{-}= \scalebox{0.8}{\bordermatrix{ 
             &1\cr
             &\ominus\cr
             &\ominus
        }}.\]
        \item The representations $T_{II,3}^{x}(\pi_{temp})$, $T_{III, 2}^{\frac{1}{2}}(\pi_{temp})$, $T_{IV, 5}(\pi_{temp})$ and $T_{V, 4}^{\pm}(\pi_{temp})$ are all not well-defined. 
        \end{enumerate}
    \end{prop}
    \begin{proof} Part $(i)$ follows from Remark \ref{rmk well-defined for Temp}.  Parts $(ii)$ and $(iii)$ follows from definition. Let $\phi$ be the L-parameter associated with $\pi_{temp}$. Then no summand of $\phi$ has multiplicity $\geq 3$, so $T_{II,3}^{x}(\pi_{temp})$ is not well-defined. By the good parity condition $T_{III,2}^{\frac{1}{2}}$ is not well-defined. Finally, since $m_{\phi}(\rho \otimes S_1) = 2$, $T_{IV,5}(\pi_{temp})$ and $T_{V,4}(\pi_{temp})$ are not well-defined. This proves part $(iv)$. \end{proof}

    This concludes our discussion of Case $(B)$. Now let's move onto Case $(C)$.
    
    \subsection{\texorpdfstring{Case $C: \pi \hookrightarrow \rho\lvert \cdot \rvert^{x_1} \times \rho\lvert \cdot \rvert^{x_2} \times \rho\lvert \cdot \rvert^{x_3} \rtimes \pi_{sc}$, where $\pi_{sc}$ is supercuspidal.}{}}
 As listed above there are $4$ total situations to consider.

    \begin{prop}\label{temp3,22}
    \begin{enumerate}[label = (\roman*)]
        \item The representation $T_{I,3}^{x}(\pi_{sc})$ is well-defined if and only if $(x,\alpha) = (\frac{1}{2}, \frac{1}{2})$. 
        \item When $\alpha = \frac{1}{2}$, the representation $T_{I,3}^{\frac{1}{2}}(\pi_{sc})$ is of critical type. 
        \item When $\alpha = \frac{1}{2}$, we have that $\pi(\EE)= T_{I,3}^{\frac{1}{2}}(\pi_{sc})$, where 
        \[\EE= \scalebox{0.8}{\bordermatrix{ 
        &\frac{1}{2}\cr
        &\oplus\cr
        &\oplus\cr
        &\oplus\cr
        }}.\]
    \end{enumerate}
    \end{prop}
    \begin{proof}
        Since $\pi_{sc}$ is supercuspidal, $\phi$ must be multiplicity free. For $T_{I,3}^{x}(\pi_{sc})$ to be well-defined we must have $m_{\phi}(\rho \otimes S_{2x -1}) \geq 3$, which can only happen when $x = \frac{1}{2}$. We also require $m_{\phi}(\rho \otimes S_{2x +1}) = 0$, which additionally force $\alpha = \frac{1}{2}$ by the good parity condition. This proves part $(i)$. Parts $(ii)$ and $(iii)$ follow from definition. 
    \end{proof}

        The proofs of Propositions \ref{temp3,25} to \ref{temp3,27} below are similar to that of Proposition \ref{temp3,22}, which we omit.   
    
    \begin{prop}\label{temp3,25}
    \begin{enumerate} [label = (\roman*)]
        \item The representation $T_{III,2}^{1}(\pi_{sc})$ is well-defined if and only if $\alpha = 1$. 
        \item When $\alpha = 1$, the representation $T_{III,2}^{1}(\pi_{sc})$ is of critical type. 
        \item When $\alpha = 1$, define: 
        \begin{equation*}
                \EE:= \{([0,0]_\rho,0,\eta),([1,1]_\rho,0,-\eta),([1,1]_\rho,0,-\eta),
        \end{equation*}
        then $\pi(\EE) = T_{III,2}^{1}(\pi_{sc})$. Here is the associated symbol. 
        \[\EE= \scalebox{0.8}{\bordermatrix{ 
             &0&1\cr
        &\odot&\cr
        &&\odot\cr
        &&\odot
        }}.\]
    \end{enumerate}
    \end{prop}
    Note that this is different than the symbol in Proposition \ref{temp3,19} since by definition, the signs in the first and second column of the symbol of $T_{III,2}^{1}(\pi_{sc})$ must be different, whereas in Proposition \ref{temp3,19} they are the same. 
    \\
    \begin{proof} Part $(i)$ follows from Remark \ref{rmk well-defined for Temp}. The rest follows from definition.
    \end{proof}

    \begin{prop}\label{temp3,26}
        \begin{enumerate}[label = (\roman*)]
            \item The representation $T_{IV,7}(\pi_{sc})$ is well-defined if and only if $\alpha \in \mathbb{Z}_{>0}$. 
            \item The representation $T_{IV,7}(\pi_{sc})$ is not of critical type. 
            \item When $\alpha \in \mathbb{Z}_{>0}$, define 
            \begin{equation*}
                \EE = \{([0,0]_\rho,0,\eta)^{6}\} + \{[\alpha -1,0]_\rho,0,\eta)\}.
            \end{equation*}
            Then $\pi(\EE) = T_{IV,7}(\pi_{sc})$. Here is the associated symbol. 
            \[\EE= \scalebox{0.8}{\bordermatrix{ 
             &0 &\cdots &\alpha -1\cr
            &\odot &&\cr
            &\odot &&\cr
            &\odot &&\cr
            &\odot &&\cr
            &\odot &&\cr
            &\odot &&\cr
            &\odot &\dots&\odot\cr
        }}.\]
        \end{enumerate}
    \end{prop}

    \begin{prop} \label{temp3,27}
        \begin{enumerate} [label = (\roman*)]
            \item The representation $T_{V,6}^{\pm}(\pi_{sc})$ is well-defined if and only if $\alpha = 0$. 
            \item When $\alpha = 0$, the representation $T_{V,6}^{\pm}(\pi_{sc})$ is of critical type.
            \item When $\alpha = 0$, the set $\Psi(T_{V,6}^{\pm})$ is a singleton and $\pi(\EE_{\pm}) = \pi_{V,6}^{\pm}(\pi_{sc})$, where 
            \[\EE_{+}= \scalebox{0.8}{\bordermatrix{ 
             &0\cr
             &\oplus\cr
             &\oplus\cr
             &\oplus\cr
             &\oplus\cr
             &\oplus\cr
             &\oplus
        }}, \quad \EE_{-}= \scalebox{0.8}{\bordermatrix{ 
             &0\cr
             &\ominus\cr
             &\ominus\cr
             &\ominus\cr
             &\ominus\cr
             &\ominus\cr
             &\ominus
        }}.\]
        \end{enumerate}
    \end{prop}

This concludes our classification of all corank 3 tempered representations of good parity. 

\section{\texorpdfstring{Classification of corank 4 non-tempered  representations of good parity $(f(\pi) = 1$)}{}}\label{classnontempcorank4,1}

 In the last section, we have the list of all tempered representations of corank $3$, in this section, we classify all $\pi \in \Pi_{gp}(G_n)$, which are non-tempered of corank $4$.  In particular, we  identify those representations which are of Arthur type. 

To exhaust all such representations $\pi$, we classify them using the number of segments in their $L$-data, denoted by $f(\pi)$. By definition, the maximum value for $f(\pi)$ is $4$ and the case $f(\pi) = 0$ corresponds to the case when $\pi$ is tempered, so it suffices to consider the case where $f(\pi) = 1,2,3,4$. In this section, we focus on the case $f(\pi) = 1$.  
To determine exactly when a non-tempered representation $\pi$ is of Arthur type, we reduce the number of segments in their $L$-data to obtain a representation $\pi^{\rho,-}$ of smaller corank, and use the known algorithms in \cite[Section 6]{HJLLZ24}. We recall the relevant definitions below. 
\begin{defn}[{\cite[Definition $6.2$]{HJLLZ24}}]\label{rhominus}
    Suppose $\pi \in \Pi(G_n)$ is non-tempered. Write $\pi = L(\Delta_{\rho_1}[x_1, -y_1], \ldots, \Delta_{\rho_f}[x_f, -y_f]; \pi(\phi, \epsilon)). $
\begin{enumerate}
    \item We define $\pi^{\rho, -}$ to be the representation whose $L$-data is obtained by removing all copies of $\Delta_{\rho}[x,-y]$ from $\pi$ such that 
    \begin{equation*}
        x = \min\{x_i | x_i - y_i = \min\{x_j - y_j | \rho_j \cong \rho\}\}
    \end{equation*}
    \begin{equation*}
        y = x - \min\{x_j - y_j | \rho_j \cong \rho\}
    \end{equation*}
    We write $\pi = \pi^{\rho, -} + \{(\Delta_{\rho}[x, -y])^r\}$, where $r$ denotes the multiplicity.
    % \item We define $\pi_{\rho, -}$ to be the rep whose $L$-data is obtained by removing all $\Delta_{\rho_i}[x_i, -y_i]$'s from $\pi$ such that $\rho_i \cong \rho$ and $x_i = \min\{x_j | \rho_j \cong \rho\}$
\end{enumerate}
\end{defn}

\begin{defn} [{\cite[Definition $6.8$]{HJLLZ24}}]\label{psiminus}
    Suppose $\pi$ is a rep of $G_n$ of good parity, and $\pi = \pi^{\rho, -}+ \{(\Delta_{\rho}[x,-y]^r)\}$. 
\begin{enumerate}
    \item We denote by $\Psi(\pi^{\rho, -}; \Delta_{\rho}[x,-y],r)$ the set of local Arthur parameters such that 
    \begin{itemize}
        \item $\pi^{\rho, -} \in \Pi_{\psi}$
        \item If $y-x-1> 0$, then $\psi$ contains $r$ copies of $\rho \otimes S_{x+y+1} \otimes S_{y-x-1}$; 
        \item Any summand of $\psi$ of the form $\rho \otimes S_a \otimes S_b$ satisfies $b \leq y-x+1$, and $a > x+y+1$ if $b = y-x+1$
    \end{itemize}
    For any $\psi \in \Psi(\pi^{\rho, -}, \Delta_{\rho}[x,-y], r)$, we define 
    \begin{equation*}
        \psi^+ := \psi - (\rho \otimes S_{x+y+1} \otimes S_{y-x-1})^{\oplus r}+ (\rho \otimes S_{x+y+1}\otimes S_{y-x+1})^{\oplus r} 
    \end{equation*}
    \item We denote by $\mathscr{E}(\pi^{\rho, -}; \Delta_{\rho}[x, -y], r)$ the set of extended multi-segments $\mathcal{E} \in \Rep^{(P')}$ such that $\pi(\mathcal{E}) = \pi^{\rho, -}$ and $\psi_{\mathcal{E}}\in \Psi(\pi^{\rho, -}; \Delta_{\rho}[x,-y],r)$. For each $\mathcal{E} \in \mathscr{E}(\pi^{\rho, -}; \Delta_{\rho}[x,-y], r)$ we define $\mathcal{E}^{\rho, +}$ as follows. 
    \begin{enumerate}
        \item If $y-x =1$. then we define $\mathcal{E}^{\rho, +}$ by inserting $r$ copies of $([x, x-1]_{\rho}, 1, 1)$ in $\mathcal{E}$ with admissible order $\gg$ on $I_{\rho} \sqcup \{j_1, \ldots, j_r\}$, where $j_k$ corresponds to the $k$-th copy we inserted, as follows: 
        $\begin{cases}
            j_r \gg j_{r-1} \gg \ldots \gg j_1, \\
            \alpha \gg \beta \Longleftrightarrow \alpha > \beta & \text{for } \alpha, \beta \in I_{\rho} \\
            \alpha \gg j_k \Longleftrightarrow B_{\alpha} > x-1 & \text{for } \alpha \in I_{\rho}
        \end{cases}$
        \item If $y-x > 1$, then change the admissible order if necessary so that there exists $j_1 \ldots j_r \in I_{\rho}$ such that 
        \begin{itemize}
            \item $j_1 < \dots < j_r$ are adjacent under the admissible order on $I_\rho$, 
            \item $[A_{j_1}, B_{j_1}]_{\rho} = \dots = [A_{j_r}, B_{j_r}]_{\rho} = [y-1, x+1]_{\rho}$
            \item $j_1 = \min\{i \in I_\rho | B_j = B_{j_1}\}$, 
            \item $A_{j_1} - B_{j_1} + 3 \geq A_i - B_i + 1$ for all $i \in I_\rho$, and the equality does not hold for $i < j_1$. 
        \end{itemize}
        Then we define $\mathcal{E}^{\rho, +} := \sum_{k=1}^{r} add_{j_k}^{1}(\mathcal{E}) \in \Rep^{(P')}$. 
    \end{enumerate}
\end{enumerate}
\end{defn}

We mainly concern the reduction from $\pi$ to $\pi^{\rho,-}$, so the following result is crucial. Explicitly, it gives us a way to construct an extended multi-segment corresponding to $\pi$ if we know the representation $\pi^{\rho,-}$ is of Arthur type. 

\begin{thm} [{\cite[Theorem $6.10$]{HJLLZ24}}] \label{nontemp red}
    An irreducible representation $\pi$ is of Arthur type if and only if  $\Psi(\pi^{\rho, -}; \Delta_{\rho}[x,-y],r)$ is nonempty. Moreover, for any $\mathcal{E} \in \mathscr{E}(\pi^{\rho, -}; \Delta_{\rho}[x,-y],r)$, we have $\pi = \pi(\mathcal{E}^{\rho, +})$. In other words, for any $\psi \in \Psi(\pi^{\rho, -}; \Delta_{\rho}[x, -y], r)$, we have $\psi^+ \in \Psi(\pi)$. 
\end{thm}

The following theorem provides certain restrictions on the $L$-data of $\pi^{\rho,-}$ for the original representation $\pi$ to be of Arthur type. Here, we make no distinction between a representation and its $L$-data. 
\begin{lemma}[{\cite[Lemma $7.8$]{HJLLZ24}}]\label{absolutevalue}
    Suppose $\pi \in \Pi_{gp}(G_n)$ and $\pi = \pi^{\rho,-}+ \{\Delta_\rho[x,-y]^{s}\}$. Then $\pi$ is of Arthur type only if $y-x =1$ or there exists at least $s$ copies of $\{\rho\lvert \cdot \rvert^{z}\}_{z = |x+1|}^{y-1}$ in $|\Omega|(\pi^{\rho,-})_{\rho}$. 
\end{lemma}

To show that certain representations are not of Arthur type, we state the following lemma, which is a generalization of \cite[Lemma $4.23$]{HJLLZ24}. 
\begin{lemma} \label{generalconstraint}
    Suppose $\mathcal{E} \in \Rep$ with $\mathcal{E}_{\rho} = \{([A_i, B_i]_\rho,l_i,\eta_i)\}_{i \in (I_{\rho},>)} + \{([X_k, X_k], 0, *)^{s_k}\}_{k =1}^{n}$  satisfying the following conditions: 
    \begin{enumerate}
    \item $X_k > \frac{1}{2}$ for $1 \leq k \leq n$
    \item For all $1 \leq k \leq n$ and $i \in I_{\rho}$, $X_k > A_i +1$.
    \end{enumerate}
    Then any extended multi-segment $\mathcal{E}'$ such that $\pi(\mathcal{E}') \cong \pi(\mathcal{E})$ must also contain $s_k$ copies of $([X_k, X_k]_\rho,0,*)$ and does not contain any extended segment of the form $([X_k,-X_k]_\rho, *, *)$.
\end{lemma}
\begin{proof} It suffices to check that if $T$ or $T^{-1}$ is a raising operator, then $T(\mathcal{E})$ or $T^{-1}(\mathcal{E})$ satisfies the given conditions. We use induction on $n$. When $n = 1$, this is exactly the same statement of \cite[Lemma $4.23$]{HJLLZ24}. 

Now assume the statement holds for $n-1$ segments of the form $([X_k, X_k],0,*)$, then applying the raising operators $ui_{i,j}^{-1}$, $dual \circ ui_{j,i} \circ dual$, or $dual_{k}^{-}$ or their inverses to the segment $([X_n, X_n],0,*)$ gives back another segment of the form $([X_n, X_n]_\rho, 0, *)$, and there are no segments of the form $([X_n, -X_n]_\rho, *,*)$, since the two conditions guarantee that the segment $[X_n, X_n]$ is separated from the other segments. This completes the proof. \end{proof}

When $f(\pi) = 1$, we have the following four cases: 
\begin{enumerate}[label = (\Alph*)]
    \item $\pi = L(\Delta_{\rho}[-x,-x]; \pi_{temp})$, where $x > 0$ and $\pi_{temp}$ is tempered of corank $3$. 
    \item $\pi = L(\Delta_\rho[-x,-x-1];\pi_{temp})$, where $x > -\frac{1}{2}$ and $\pi_{temp}$ is tempered of corank $2$. 
    \item $\pi = L(\Delta_\rho[-x,-x-2];\pi_{temp})$, where $x > -1$ and $\pi_{temp}$ is tempered of corank $1$. 
    \item $\pi = L(\Delta_\rho[-x,-x-3];\pi_{sc})$, where $x> -\frac{3}{2}$ and $\pi_{sc}$ is supercuspidal. 
\end{enumerate}

\subsection{\texorpdfstring{Case $(A)$: $\pi = L(\Delta_{\rho}[-x,-x]; \pi_{temp})$}{}}
Let us begin with Case $(A)$. Following the classification of corank 3 tempered representations of good parity in \S \ref{classtempcorank3}, we consider case by case as follows. 

\begin{prop} \label{nontempA1}
    Consider the representation
    \begin{equation*}
        \pi_x = L(\Delta_{\rho}[-x,-x]; \pi_{temp}),
    \end{equation*}
    where $\pi_{temp} = T_{I,1}^{\alpha+2}(T_{I,1}^{\alpha+1}(T_{I,1}^{\alpha}(\pi_{sc})))$ for $x > 0$ and $\alpha > 0$ as described in Proposition \ref{temp3,1}. 
    \begin{enumerate} [label = (\roman*)]
        \item The representation $\pi_x$ is of Arthur type if and only if $\frac{1}{2} \leq x \leq \alpha-1$ when $\alpha \geq \frac{3}{2}$, or $x = \alpha = \frac{1}{2}$. 
        \item The representation $\pi_x$ is of critical type when $x \in \{\alpha-1, \alpha, \alpha +1, \alpha +2, \alpha +3\}$. 
        \item Define $\EE_x$ in various cases as follows. Then $\pi_x = \pi(\EE_x)$. When $\alpha \geq 2$, set 
        \begin{flalign*}
            \EE_x &:= \{([x,-x]_\rho,\lfloor x \rfloor, (-1)^{2x}\eta), ([x-2, \epsilon_\rho]_\rho,0,-\eta)\} \\
            &+ \{([z,z]_\rho,0,(-1)^{z - \epsilon_\rho}\eta)\}_{z = x}^{\alpha-2} + \{([\alpha+2,\alpha+2]_\rho,0,(-1)^{\alpha-1-\epsilon_\rho}\eta)\}. 
        \end{flalign*}
        Here is the associated symbol. 
            \[\EE_{x}= \scalebox{0.8}{\bordermatrix{ 
             &-x &-(x-1) &\cdots & \epsilon_\rho & \cdots & x-2 & x-1 & x & \cdots & \alpha-2 & \alpha-1 & \alpha & \alpha+1 & \alpha+2\cr
            &\lhd & \lhd & \cdots & \cdots & \dots & \rhd & \rhd &\rhd \cr
            &&&&\odot & \cdots & \odot \cr
            &&&&&&&&\odot \cr
            &&&&&&&&&\ddots \cr
            &&&&&&&&&& \odot \cr
            &&&&&&&&&&&&&&\odot 
        }}.\]
        When $\alpha = \frac{3}{2}$ (resp. $\frac{1}{2}$), $\EE_{\frac{1}{2}}$ are given as follows respectively, 
        \[\EE_{\frac{1}{2}}= \scalebox{0.8}{\bordermatrix{ 
             &-\frac{1}{2} & \frac{1}{2} & \frac{3}{2} & \frac{5}{2} & \frac{7}{2} \cr
             & \lhd & \rhd \cr
             &&&&&\oplus
        }}, \quad \Bigg(\text{resp. }\EE_{\frac{1}{2}}= \scalebox{0.8}{\bordermatrix{ 
             &-\frac{1}{2} & \frac{1}{2} & \frac{3}{2} & \frac{5}{2} \cr
             & \lhd & \rhd \cr
             &&&&\oplus
        }}\Bigg).\]
    \end{enumerate}
\end{prop}
{\it Proof:}
We first show part $(iii)$. When $\alpha < 2$, the assertion can be easily verified. For $\alpha \geq 2$, Proposition \ref{temp3,1} tells us that $\pi(add_{1}^{-1}(\EE_x)) = \pi_{temp}$ after applying  a sequence of operators (e.g. take the dual of the first row and take the union-intersection of the first two rows) on $\pi(add_{1}^{-1}(\EE_x))$. Therefore by Theorem \ref{nontemp red}, we obtain: 
\begin{equation*}
    \pi(\EE_x) = \pi(add_{1}^{-1}(\EE_x) + \{\Delta_{\rho}[-x,-x]\}) = \pi_x.
\end{equation*} This proves part $(iii)$ and the sufficient direction of part $(i)$. 

Now let's show the necessary direction of part $(i)$. Since $|\Omega| (\pi_{temp}) = \{\rho \lvert \cdot \rvert^{y}: \epsilon_\rho \leq y \leq \alpha -2 \} + \{\rho \lvert \cdot \rvert^{\alpha+2}\}$. By Lemma \ref{absolutevalue}, the representation $\pi_x$ is of Arthur type only in one of the following situations: 
\begin{enumerate}
    \item $x = \frac{1}{2}$ and $\alpha \in \mathbb{Z}_{> 0}$,
    \item $\frac{1}{2} \leq x \leq \alpha -1 $ and $\alpha \geq \frac{3}{2}$,
    \item $x = \alpha + 3$ with $\alpha \in x + \mathbb{Z}$. 
\end{enumerate}
    Therefore it suffices to show that $\pi_x$ is not of Arthur type in situation $(3)$. Let $x = \alpha + 3$ and $\alpha \in x + \mathbb{Z}$. A priori, we have an extended multi-segment $\EE_{temp}$ listed in Proposition \ref{temp3,1} (iii) such that $\pi(\EE_{temp}) = \pi_{temp}$. Then by Lemma \ref{generalconstraint}, there does not exist any multi-segment $\EE$ that contains $([\alpha +2 ,-\alpha -2]_\rho,*,*)$ such that $\pi(\EE) = \pi_{temp}$. 

However, by Definition \ref{psiminus}, any element in $\Psi(\pi_{temp}, \Delta_{\rho}[-x,-x],1)$ must contain $1$ copy of $\rho \otimes S_1 \otimes S_{2\alpha +5}$, which is equivalent to saying that the corresponding extended multi-segment must contain some element of the form $([\alpha+2, -\alpha-2]_\rho,*,*)$. This implies that the set $\Psi(\pi_{temp}; \Delta_{\rho}[-x,-x],1)$ must be empty, so $\pi_x$ is not of Arthur type by Theorem \ref{nontemp red}. This completes the proof. 
\end{proof}

Next we consider the case $\pi_{temp} = T_{I,1}^{\alpha-1}(T_{I,1}^{\alpha+1}(T_{I,1}^{\alpha}(\pi_{sc})))$. The proofs of Propositions \ref{nontempA2} to \ref{nontempA4} below are similar to that of Proposition \ref{nontempA1}, which we omit. 

\begin{prop} \label{nontempA2}
    Consider the representation
    \begin{equation*}
        \pi_x = L(\Delta_{\rho}[-x,-x]; \pi_{temp}),
    \end{equation*}
    where $\pi_{temp} = T_{I,1}^{\alpha-1}(T_{I,1}^{\alpha+1}(T_{I,1}^{\alpha}(\pi_{sc})))$ for $x > 0$ and $\alpha > 1$ as described in Proposition \ref{temp3,1}. 
    \begin{enumerate}[label = (\roman*)]
        \item The representation $\pi_x$  is of Arthur type if and only if $\frac{1}{2} \leq x \leq \alpha -2$  when $\alpha \geq \frac{5}{2}$, or when $x = \alpha = \frac{3}{2}$. 
        \item The representation $\pi_x$ is of critical type when $x \in \{\alpha-2, \alpha-1, \alpha, \alpha +1, \alpha +2\}$.
        \item When $\frac{1}{2} \leq x \leq \alpha -2$ for $\alpha \geq \frac{5}{2}$, or $x = \alpha = \frac{3}{2}$, define
        \begin{flalign*}
            \EE_x &:= \{([x,-x]_\rho,\lfloor x \rfloor, (-1)^{2x}\eta), ([x-2, \epsilon_\rho]_\rho,0,-\eta)\} +
             \{([z,z]_\rho,0,(-1)^{z - \epsilon_\rho}\eta)\}_{z = x}^{\alpha-3} \\
             &+ \{([\alpha-1,\alpha-1]_\rho,0,(-1)^{\alpha-1-\epsilon_\rho}\eta),([\alpha+1,\alpha+1]_\rho,0,(-1)^{\alpha-2-
            \epsilon_\rho}\eta\}. 
        \end{flalign*}
        Then $\pi_x = \pi(\EE_x)$. Here is the associated symbol. 
        \[\EE= \scalebox{0.8}{\bordermatrix{ 
             &-x &-(x-1) &\cdots & \epsilon_\rho & \cdots & x-2 & x-1 & x & \cdots & \alpha-3 & \alpha -2& \alpha-1 & \alpha & \alpha+1 \cr
            &\lhd & \lhd & \cdots & \odot & \cdots & \rhd & \rhd &\rhd \cr
            &&&&\odot & \dots & \odot \cr
            &&&&&&&&\odot \cr
            &&&&&&&&&\ddots \cr
            &&&&&&&&&& \odot \cr
            &&&&&&&&&&&&&&\odot 
        }}.\]
    \end{enumerate}
\end{prop}

The next case to consider is $\pi_{temp} = T_{IV,3}(T_{I,1}^{\alpha + 1}(T_{I,1}^{\alpha}(\pi_{sc})))$ for $\alpha \in \mathbb{Z}_{>1}$. 

\begin{prop}\label{nontempA3}
Consider the representation 
\begin{equation*}
    \pi_x = L(\Delta_{\rho}[-x,-x]; \pi_{temp})
\end{equation*}
where 
     $\pi_{temp} = T_{IV,3}(T_{I,1}^{\alpha + 1}(T_{I,1}^{\alpha}(\pi_{sc})))$ for $x > 0$ and  $\alpha \in \mathbb{Z}_{>1}$, as described in Proposition \ref{temp3,2}. 
    \begin{enumerate}[label = (\roman*)]
        \item The representation $\pi_x$ is of Arthur type if and only if $1 \leq x \leq \alpha -1$. 
        \item The representation $\pi_x$ is of critical type when $(x,\alpha) = (1,2)$. 
        \item When $1 \leq x \leq \alpha -1$, define 
        \begin{flalign*}
            \EE_x &:= \{([x,-x]_\rho,\lfloor x \rfloor, (-1)^{2x}\eta)\}  + \{([0,0]_\rho,0,\eta))^2\} + ([x-2, 0]_\rho,0,\eta)\}  \\
            &+ \{([z,z]_\rho,0,(-1)^{z - \epsilon_\rho}\eta)\}_{z = x}^{\alpha-2} 
              + \{([\alpha+1,\alpha+1],0,(-1)^{\alpha-2-
            \epsilon_\rho}\eta\}. 
        \end{flalign*}
        Then $\pi_x = \pi(\EE_x)$. Here is the associated symbol. 
        \[\EE= \scalebox{0.8}{\bordermatrix{ 
             &-x  & -(x-1)& \cdots & 0 & \cdots  & x-2 & x-1 &x & \cdots & \alpha -2& \alpha-1 & \alpha & \alpha+1 \cr
            &\lhd  &\lhd & \cdots & \odot & \cdots & \rhd & \rhd&\rhd \cr
            &&&&\odot\cr
            &&&&\odot \cr
            &&&&\odot & \cdots & \odot \cr
            &&&&&&&&\odot \cr
            &&&&&&&&&\ddots \cr
            &&&&&&&&&& \odot \cr
            &&&&&&&&&&&&&\odot 
        }}.\]
    \end{enumerate}
\end{prop}

Now we consider the case $\pi_{temp} = T_{V,2}^{\pm}(T_{I,1}^{2}(T_{I,1}^{1}(\pi_{sc})))$. 
\begin{prop}\label{nontempA4}
    Consider the representation 
\begin{equation*}
    \pi_x^{\pm} = L(\Delta_{\rho}[-x,-x]; \pi_{temp}^{\pm}),
\end{equation*}
where $\pi_{temp}^{\pm} = T_{V,2}^{\pm}(T_{I,1}^{2}(T_{I,1}^{1}(\pi_{sc})))$ and $\alpha = 1$,  as described in Proposition \ref{temp3,3}. 
\begin{enumerate}[label = (\roman*)]
    \item The representation $\pi_x^{\pm}$ is of Arthur type if and only if $x = 1$. 
    \item The representation $\pi_x^{\pm}$ is of critical type for $x \in \{1,2,3\}$.
    \item Define 
    \begin{equation*}
        \EE_{\pm} = \{([-1,-1],1,\eta)\}\{([0,0],0,\pm1))^3,\} + \{([2,2],0,\eta)\}.
    \end{equation*}
     Then $\pi_1^{\pm} = \pi(\EE_{\pm})$. Here are the associated symbols. 
     \[\EE_{+}= \scalebox{0.8}{\bordermatrix{ 
        &0&1&2\cr
        &\lhd & \odot & \rhd \cr
        &\oplus&&\cr
        &\oplus&&\cr
        &\oplus&&\cr
        &&&\odot
        }}, \quad \EE_{-}= \scalebox{0.8}{\bordermatrix{ 
        &0&1&2\cr
        &\lhd & \odot & \rhd \cr
        &\ominus&&\cr
        &\ominus&&\cr
        &\ominus&&\cr
        &&&\odot
        }}.\]         
\end{enumerate}
\end{prop}

The next case is $\pi_{temp} = T_{I,1}^{\alpha-2}(T_{I,1}^{\alpha-1}(T_{I,1}^{\alpha}(\pi_{sc})))$. 
\begin{prop}\label{nontempA5}
    Consider the representation 
\begin{equation*}
    \pi_x = L(\Delta_{\rho}[-x,-x]; \pi_{temp})
\end{equation*}
where $\pi_{temp} = T_{I,1}^{\alpha-2}(T_{I,1}^{\alpha-1}(T_{I,1}^{\alpha}(\pi_{sc})))$ for $x> 0$ and $\alpha > 2$ as described in Proposition \ref{temp3,4}. 
\begin{enumerate}[label = (\roman*)]
    \item The representation $\pi_x$ is of Arthur type if and only if $0 < x \leq \alpha -3$ for $\alpha \geq \frac{7}{2}$. 
    \item The representation $\pi_x$ is of critical type when $x \in \{\alpha -3, \alpha -2, \alpha -1, \alpha, \alpha+1\}$. 
    \item Define 
    \begin{equation*}
        \EE_x = \{([x,-x]_\rho,\lfloor x \rfloor, (-1)^{2x}\eta), ([x-2, \epsilon_\rho]_\rho,0,\eta),([\alpha -4, x]_\rho, 0, (-1)^{x - \epsilon_
        \rho}\eta),([\alpha, \alpha-2]_\rho, 0, -\eta)\}.
    \end{equation*}
    Then $\pi_x = \pi(\EE_x)$. Here is the associated symbol. 
    \[\EE_x= \scalebox{0.8}{\bordermatrix{ 
             &-x  & -(x-1)& \cdots & \epsilon_\rho & \cdots  & x-2 & x-1 &x & \cdots & \alpha -4& \alpha-3 & \alpha-2 & \alpha-1 &\alpha \cr
            &\lhd  &\lhd & \cdots & \odot  & \cdots & \rhd & \rhd&\rhd \cr
            &&&&\odot & \cdots  & \odot \cr
            &&&&&&&&&\odot \cr
            &&&&&&&&&&\ddots \cr
            &&&&&&&&&&& \odot \cr
            &&&&&&&&&&&&\odot&\odot &\odot
        }}.\]
\end{enumerate}
\begin{proof} By Lemma \ref{absolutevalue}, it suffices to prove that $\pi_x$ is not of Arthur type when $x = \alpha-1, \alpha, \alpha + 1$. If $x = \alpha -1$, then by Proposition \ref{temp3,4} and Lemma \ref{generalconstraint}, there does not exist any extended multi-segment $\mathcal{E}$ that contains $([\alpha-2, -\alpha +2]_{\rho}, *, *)$ such that $\pi(\mathcal{E}) = \pi_{temp}$. 

However, by Definition \ref{psiminus}, any element in $\Psi(\pi_{temp}, \Delta_{\rho}[-x,-x],1)$ must contain $1$ copy of $\rho \otimes S_1 \otimes S_{2\alpha -1}$, which means that the corresponding extended multi-segment must contain a segment of the form $([\alpha-2, -\alpha +2]_{\rho}, *, *)$. Therefore the set $\Psi(\pi_{temp}, \Delta_{\rho}[-x,-x],1)$ must be empty, so $\pi_x$ is not of Arthur 
type by Theorem \ref{nontemp red}. This completes the proof for $x = \alpha -1$. The cases are $x = \alpha, \alpha +1$ are exactly the same. 
\end{proof}
\end{prop}

Now we consider the case $\pi_{temp} = T_{IV, 3}(T_{I,1}^{\alpha -1}(T_{I,1}^{\alpha}(\pi_{sc})))$.
The proofs of Propositions \ref{nontempA6} to \ref{nontempA11} below are similar to that of Proposition \ref{nontempA4}, which we omit. 

\begin{prop}\label{nontempA6}
     Consider the representation 
\begin{equation*}
    \pi_x = L(\Delta_{\rho}[-x,-x]; \pi_{temp}),
\end{equation*}
where $\pi_{temp} = T_{IV, 3}(T_{I,1}^{\alpha -1}(T_{I,1}^{\alpha}(\pi_{sc})))$ for $x > 0$ and $\alpha \in \mathbb{Z}_{>2}$, as described in Proposition \ref{temp3,5}. 
\begin{enumerate}[label = (\roman*)]
    \item The representation $\pi_x$ is of Arthur type if and only if $0 < x \leq \alpha -2$. 
    \item The representation $\pi_x$ is of critical type when $(x,\alpha) = (1,3)$. 
    \item When $0 \leq x \leq \alpha -2$, define
    \begin{flalign*}
        \EE_x &= \{([x,-x]_\rho,\lfloor x \rfloor, (-1)^{2x}\eta), ([0,0]_\rho,0,\eta)^2 ,([x-2, 0]_\rho,0,\eta),\\
        &([\alpha -3, x]_\rho, 0, (-1)^{x}\eta),([\alpha-1, \alpha-1]_\rho, 0, -\eta),([\alpha,\alpha]_\rho,0,\eta)\}.
    \end{flalign*}
    Then $\pi(\EE_x) = \pi_x$. Here is the associated symbol.
    \[\EE_x= \scalebox{0.8}{\bordermatrix{ 
             &-x  & -(x-1)& \cdots & 0 & \dots  & x-2 & x-1 &x & \cdots & \alpha -3 &\alpha-2 & \alpha-1 &\alpha \cr
            &\lhd  &\lhd & \cdots & \odot & \cdots & \rhd & \rhd&\rhd \cr
            &&&&\odot \cr
            &&&&\odot \cr
            &&&&\odot & \cdots  & \odot \cr
            &&&&&&&&\odot \cr
            &&&&&&&&&\ddots \cr
            &&&&&&&&&& \odot \cr
            &&&&&&&&&&&&\odot\cr
            &&&&&&&&&&&&&\odot\cr
        }}.\]
\end{enumerate}
\end{prop}

The next case is $\pi_{temp} = T_{V,2}^{\pm} (T_{I,1}^{1}(T_{I,1}^{2}(\pi_{sc})))$.
\begin{prop}\label{nontempA7}
    Consider the representation 
\begin{equation*}
    \pi_x^{\pm} = L(\Delta_{\rho}[-x,-x]; \pi_{temp}^{\pm}),
\end{equation*}
where $\pi_{temp}^{\pm} = T_{V, 2}^{\pm}(T_{I,1}^{1}(T_{I,1}^{2}(\pi_{sc})))$ for $x> 0$, $\alpha = 2$, as described in Proposition \ref{temp3,6}. 
\begin{enumerate}
    \item The representation $\pi_{x}^{\pm}$ is of Arthur type if and only if $x = 1$.
    \item The representation $\pi_{x}^{\pm}$ is of critical type for $x \in \{1,2,3\}$.
    \item Define
    \begin{equation*}
        \EE_{\pm} := \{([1,-1]_\rho,1,\pm1),([0,0]_\rho,0,\pm1)^2,([2,1]_\rho,0,\eta)\}.
    \end{equation*}
    Then $\pi(\EE_{\pm}) = \pi_{1}^{\pm}$. Here are the associated symbols.
    \[\EE_{+}= \scalebox{0.8}{\bordermatrix{ 
        &-1 & 0&1&2\cr
        &\lhd & \oplus & \rhd \cr
        &&\oplus&&\cr
        &&&\oplus&\ominus\cr
        }}, \quad \EE_{-}= \scalebox{0.8}{\bordermatrix{ 
        &-1 & 0&1&2\cr
        &\lhd & \ominus & \rhd \cr
        &&\ominus&&\cr
        &&&\ominus&\oplus\cr
        }}.\]       
\end{enumerate}
\end{prop}

% As described in Proposition \ref{temp3,7}, there is one extra case not included in the previous propositions due to the restraints on $\alpha$. 
% \begin{prop}\label{nontempA8}
%      Consider the representation 
% \begin{equation*}
%     \pi_x = L(\Delta_{\rho}[-x,-x]; \pi_{temp}),
% \end{equation*}
% where $\pi_{temp} = T_{I,1}^{1}(T_{IV,3}(T_{I,1}^{2}(\pi_{sc})))$ for $x>0$ and $\alpha = 2$, as described in Proposition \ref{temp3,7}. 
% \begin{enumerate}[label = (\roman*)]
%     \item The representation $\pi_x$ is of Arthur type if and only if $x = 1$. 
%     \item The representation $\pi_1$  is of critical type when $x \in \{1,2,3\}$. 
%     \item Define
%     \begin{equation*}
%         \EE := \{([1,-1]_\rho,1,\eta),([2,0]_\rho,0,\eta)\}.
%     \end{equation*}
%     Then $\pi(\EE) = \pi_1$. Here is the associated symbol. 
%     \[\EE= \scalebox{0.8}{\bordermatrix{ 
%         &-1&0&1&2\cr
%         &\lhd & \odot & \rhd \cr
%         &&\odot&\odot&\odot\cr
%         }}.\]
% \end{enumerate}
% \end{prop}
% \begin{proof} The proof is similar to Proposition \ref{nontempA7}, which we omit. \end{proof}
Now we move onto the case $\pi_{temp}^{\pm} = T_{I,1}^{1}(T_{V,2}^{\pm}(T_{I,1}^{1}(\pi_{sc})))$.
\begin{prop}\label{nontempA9}
    Consider the representation 
\begin{equation*}
    \pi_x^{\pm} = L(\Delta_{\rho}[-x,-x]; \pi_{temp}^{\pm}),
\end{equation*}
where $\pi_{temp}^{\pm} = T_{I,1}^{2}(T_{V,2}^{\pm}(T_{I,1}^{1}(\pi_{sc})))$ for $x > 0$ and $\alpha =1$, as described in Proposition \ref{temp3,8}. 
\begin{enumerate}[label = (\roman*)]
    \item The representation $\pi_x$ is of Arthur type if and only if $x = 1$. 
    \item The representation $\pi_x$ is of critical type when $x \in \{1,2,3\}$.
    \item Define 
    \begin{equation*}
        \EE_{\pm} := \{([1,-1]_\rho,1,\pm1),([0,0]_\rho,0,\pm1)^2,([2,1]_\rho,0,\eta)\}.
    \end{equation*}
    Then $\pi(\EE_{\pm}) = \pi_{1}^{\pm}$. Here are the associated symbols. 
    \[\EE_{+}= \scalebox{0.8}{\bordermatrix{ 
        &-1 & 0&1&2\cr
        &\lhd & \oplus & \rhd \cr
        &&\oplus&&\cr
        &&&&\odot
        }}, \quad \EE_{-}= \scalebox{0.8}{\bordermatrix{ 
        &-1 & 0&1&2\cr
        &\lhd & \ominus & \rhd \cr
        &&\ominus&&\cr
        &&&&\odot
        }}.\]       
\end{enumerate}
\end{prop}

Now we consider the case $\pi_{temp} = T_{I,1}^{2}(T_{I,1}^{1}(T_{V,2}^{\pm}(\pi_{sc})))$. 
\begin{prop}\label{nontempA10}
    Consider the representation 
\begin{equation*}
    \pi_x^{\pm} = L(\Delta_{\rho}[-x,-x]; \pi_{temp}^{\pm}),
\end{equation*}
where $\pi_{temp}^{\pm} = T_{I,1}^{2}(T_{I,1}^{1}(T_{V,2}^{\pm}(\pi_{sc})))$ for $x > 0$ and $\alpha = 0$, as described in Proposition \ref{temp3,9}. 
\begin{enumerate}[label = (\roman*)]
    \item The representation $\pi_x^{\pm}$ is of Arthur type if and only if $x = 1$. 
    \item The representation $\pi_x^{\pm}$ is of critical type when $x \in \{1,2,3\}$.
    \item Define 
    \begin{equation*}
        \EE_{\pm} := \{([1,-1]_\rho,1,\pm1),([2,2]_\rho,0,\pm1)\}.
    \end{equation*}
    Then $\pi(\EE_{\pm}) = \pi_{1}^{\pm}$. Here are the associated symbols. 
    \[\EE_{+}= \scalebox{0.8}{\bordermatrix{ 
        &-1 & 0&1&2\cr
        &\lhd & \oplus & \rhd \cr
        &&&&\oplus
        }}, \quad \EE_{-}= \scalebox{0.8}{\bordermatrix{ 
        &-1 & 0&1&2\cr
        &\lhd & \ominus & \rhd \cr
        &&&&\ominus
        }}.\]
\end{enumerate}
\end{prop}

The next case is $\pi_{temp} = T_{IV,3}(T_{I,1}^{1}(T_{V,2}^{\pm}(\pi_{sc})))$. 
\begin{prop}\label{nontempA11}
     Consider the representation 
\begin{equation*}
    \pi_x^{\pm} = L(\Delta_{\rho}[-x,-x]; \pi_{temp}^{\pm}),
\end{equation*}
where $\pi_{temp}^{\pm} =T_{IV,3}(T_{I,1}^{1}(T_{V,2}^{\pm}(\pi_{sc})))$ for $\alpha = 0$, as described in Proposition \ref{temp3,10}. 
\begin{enumerate}[label = (\roman*)]
    \item The representation $\pi_x^{\pm}$ is of Arthur type if and only if $x = 1$. 
    \item The representation $\pi_x^{\pm}$ is of critical type when $x \in \{1,2\}$.
    \item Define 
    \begin{equation*}
        \EE_{\pm} := \{([1,-1]_\rho,1,\pm1), ([0,0]_\rho,0,\pm1), ([1,1]_\rho,0,\pm1)\}.
    \end{equation*}
    Then $\pi(\EE_{\pm}) = \pi_1^{\pm}$. Here are the associated symbols. 
     \[\EE_{+}= \scalebox{0.8}{\bordermatrix{ 
        &-1 & 0&1\cr
        &\lhd & \odot & \rhd \cr
        &&\odot \cr
        &&\oplus \cr
        &&&\oplus
        }}, \quad \EE_{-}= \scalebox{0.8}{\bordermatrix{ 
        &-1 & 0&1\cr
        &\lhd & \odot & \rhd \cr
        &&\odot \cr
        &&\ominus \cr
        &&&\ominus
        }}.\]
\end{enumerate}
\end{prop}

Now let's look at $\pi_{temp} = T_{I,1}^{\frac{3}{2}}(T_{I,2}^{\frac{1}{2}}(\pi_{sc}))$. 
\begin{prop}\label{nontempA12}
    Consider the representation 
\begin{equation*}
    \pi_x = L(\Delta_{\rho}[-x,-x]; \pi_{temp}),
\end{equation*}
where $\pi_{temp} = T_{I,1}^{\frac{3}{2}}(T_{I,2}^{\frac{1}{2}}(\pi_{sc}))$ for $x > 0$ and $\alpha = \frac{1}{2}$, as described in Proposition \ref{temp3,11}. 
\begin{enumerate}[label = (\roman*)]
    \item The representation $\pi_x$ is not of Arthur type for any $x$. 
    \item The representation $\pi_x$ of critical type when $x \in \{\frac{1}{2}, \frac{3}{2}, \frac{5}{2}\}$.
\end{enumerate}
\end{prop}
\begin{proof} By Lemma \ref{absolutevalue}, $\pi_x$ can only be of Arthur type if $x = \frac{5}{2}$. When $x = \frac{5}{2}$, we see from Proposition \ref{temp3,11} that the support of the extended multi-segment corresponding to $\pi_x$ does not include $\frac{1}{2}$, so by Definition \ref{psiminus} and Theorem \ref{nontemp red}, it cannot be of Arthur type. \end{proof}

The next case we have is $\pi_{temp} = T_{I,1}^{\alpha}(T_{II,3}^{\frac{1}{2}}(\pi_{sc}))$. 
\begin{prop}\label{nontempA13}
      Consider the representation 
\begin{equation*}
    \pi_x = L(\Delta_{\rho}[-x,-x]; \pi_{temp}),
\end{equation*}
where $\pi_{temp} = T_{I,1}^{\alpha}(T_{II,3}^{\frac{1}{2}}(\pi_{sc}))$ for $x > 0$ and $\alpha \in \frac{1}{2} + \mathbb{Z}_{>0}$, as described in Proposition \ref{temp3,12}. 
\begin{enumerate}[label = (\roman*)]
    \item The representation $\pi_x$ is of Arthur type if and only if $\frac{1}{2} \leq x \leq \alpha -1$.  
    \item The representation $\pi_x$ is of critical type when $(x,\alpha) = (\frac{1}{2},\frac{3}{2}),(\frac{3}{2}, \frac{3}{2}),(\frac{5}{2}, \frac{3}{2})$ or $(\frac{3}{2}, \frac{5}{2})$. 
    \item Define 
    \begin{flalign*}
        \EE_x := &\{([x,-x]_\rho,\lfloor x \rfloor, (-1)^{x + \frac{1}{2}}),([\frac{1}{2}, \frac{1}{2}]_\rho,0,-1)^2,\\
        &([x -2, \frac{3}{2}]_\rho, 0,1),([\alpha-2,x]_\rho,0,(-1)^{x+\frac{1}{2}}),([\alpha,\alpha]_\rho,0,(-1)^{\alpha - \frac{1}{2}})\}.
    \end{flalign*}
    Then $\pi(\EE_x) = \pi_x$. Here is the associated symbol. 
    \[\EE_x= \scalebox{0.8}{\bordermatrix{ 
             &-x  & -(x-1)& \cdots & \frac{1}{2} & \frac{3}{2}& \cdots  & x-2 & x-1 &x & \cdots & \alpha -2 &\alpha-1  &\alpha \cr
            &\lhd  &\lhd & \cdots & \ominus & \cdots & \rhd & \rhd&\rhd \cr
            &&&&\ominus \cr
            &&&&\ominus \cr
            &&&&& \oplus &\cdots& \odot \cr
            &&&&&&&&&\odot \cr
            &&&&&&&&&&\ddots \cr
            &&&&&&&&&&& \odot \cr
            &&&&&&&&&&&&&\odot\cr
        }}.\]
\end{enumerate}
\end{prop}
\begin{proof}
The proof is similar to that of Proposition \ref{nontempA12}, which we omit. \end{proof}

Now we move onto the case $\pi_{temp} = T_{I,1}^{\frac{3}{2}}(T_{III,2}^{\frac{1}{2}}(\pi_{sc}))$. 
\begin{prop}\label{nontempA14}
     Consider the representation 
\begin{equation*}
    \pi_x = L(\Delta_{\rho}[-x,-x]; \pi_{temp}),
\end{equation*}
where $\pi_{temp} =T_{I,1}^{\frac{3}{2}}(T_{III,2}^{\frac{1}{2}}(\pi_{sc})) $, for $x > 0$ and $\alpha = \frac{1}{2}$, as described in Proposition \ref{temp3,13}. 
\begin{enumerate}[label = (\roman*)]
    \item The representation $\pi_x$ is of Arthur type if and only if $x = \frac{1}{2}$ or $\frac{3}{2}$. 
    \item The representation $\pi_{x}$ is of critical type when $x \in \{\frac{1}{2}, \frac{3}{2}, \frac{5}{2}\}$. 
    \item For $x \in \{\frac{1}{2}, \frac{3}{2}\}$, we have $\pi(\EE_x) = \pi_x$, where $\EE_x$ are given as follows. 
    \[\EE_{\frac{1}{2}}= \scalebox{0.8}{\bordermatrix{ 
        &-\frac{1}{2} & \frac{1}{2}&\frac{3}{2}\cr
        &\lhd & \ominus & \rhd \cr
        &&\ominus \cr
        &&&\ominus
        }}, \quad \EE_{\frac{3}{2}}= \scalebox{0.8}{\bordermatrix{ 
        &-\frac{3}{2} & -\frac{1}{2}&\frac{1}{2}&\frac{3}{2}\cr
        &\lhd & \oplus & \ominus &\rhd \cr
        &&&&\ominus\cr
        }}.\]   
\end{enumerate}
\end{prop}
\begin{proof} By Lemma \ref{absolutevalue} and Proposition \ref{temp3,13}, the only possible $x$ for $\pi_x$ to be of Arthur type are $x = \frac{1}{2}, \frac{3}{2}, \frac{5}{2}$. By Definition \ref{rhominus}, the set $\Psi(\pi_{temp}, [-x,-x],1)$ is nonempty if and only if it contains some parameter $\psi$ containing $\rho \otimes S_1 \otimes S_{2x-1}$, or that the corresponding extended multisegment contains the segment $([x-1,-(x-1)]_\rho,*,*)$. 

Looking at the possible extended multisegments in Proposition \ref{temp3,13}, this only happens when $x = \frac{1}{2}, \frac{3}{2}$. By Theorem \ref{nontemp red}, this is exactly when $\pi_x$ is of Arthur type. This proves part $(i)$. Part $(ii)$ follows from definition, and part $(iii)$ follows from Proposition \ref{temp3,13} $(iii)$. \end{proof}

The next case is $\pi_{temp} = T_{I,1}^{\alpha}(T_{IV,5}(\pi_{sc}))$. The proofs of Propositions \ref{nontempA15} to \ref{nontempA19} below are similar to that of Proposition \ref{nontempA14}, which we omit.

\begin{prop}\label{nontempA15}
     Consider the representation 
\begin{equation*}
    \pi_x = L(\Delta_{\rho}[-x,-x]; \pi_{temp}),
\end{equation*}
where $\pi_{temp} = T_{I,1}^{\alpha}(T_{IV,5}(\pi_{sc}))$ for $x> 0$ and $\alpha \in \mathbb{Z}_{> 0}$, as described in Proposition \ref{temp3,14}. 
\begin{enumerate}[label = (\roman*)]
    \item The representation $\pi_x$ if of Arthur type if and only if $1 \leq x \leq \alpha -1$ for $\alpha \geq 2$, or when $x = \alpha =1$. 
    \item The representation $\pi_x$ is of critical type when $(x,\alpha) = (1,1)$ or $(1,2)$. 
    \item Define 
    \begin{flalign*}
        \EE_x := &\{([x,-x]_\rho,\lfloor x \rfloor, (-1)^{x}\eta),([0, 0]_\rho,0,\eta)^5,\\
        &([x -2, 1]_\rho, 0,-\eta),([\alpha-2,x]_\rho,0,(-1)^{x}\eta),([\alpha,\alpha]_\rho,0,(-1)^{\alpha-1}\eta)\}.
    \end{flalign*}
    Then $\pi(\EE_x) = \pi_x$. Here is the associated symbol. 
     \[\EE_x= \scalebox{0.8}{\bordermatrix{ 
             &-x  & -(x-1)& \cdots & 0 & 1& \cdots  & x-2 & x-1 &x & \cdots & \alpha -2 &\alpha-1  &\alpha \cr
            &\lhd  &\lhd & \cdots & \ominus & \oplus & \cdots & \rhd & \rhd&\rhd \cr
            &&&&\odot \cr
            &&&&\odot \cr
            &&&&& \odot &\cdots& \odot \cr
            &&&&&&&&&\odot \cr
            &&&&&&&&&&\ddots \cr
            &&&&&&&&&&& \odot \cr
            &&&&&&&&&&&&&\odot\cr
        }}.\]
\end{enumerate}
\end{prop}

Now we move onto the case $\pi_{temp} = T_{I,1}^{1}(T_{V,4}^{\pm}(\pi_{sc}))$. 
\begin{prop}\label{nontempA16}
    Consider the representation 
\begin{equation*}
    \pi_x^{\pm} = L(\Delta_{\rho}[-x,-x]; \pi_{temp}^{\pm}),
\end{equation*}
where $\pi_{temp}^{\pm} =T_{I,1}^{1}(T_{V,4}^{\pm}(\pi_{sc}))$ for $x > 0$ and $\alpha = 0$, as described in Proposition \ref{temp3,15}. 
\begin{enumerate}[label = (\roman*)]
    \item The representation $\pi_{x}^{\pm}$ is of Arthur type if and only if $x = 1$. 
    \item The representation $\pi_x^{\pm}$ is of critical type when $x \in \{1,2\}$. 
    \item We have that $\pi(\EE_{\pm}) = \pi_{x}^{\pm}$, where 
     \[\EE_{+}= \scalebox{0.8}{\bordermatrix{ 
        &-1 & 0&1\cr
        &\lhd & \oplus & \rhd \cr
        &&\oplus \cr
        &&\oplus \cr
        &&&\oplus
        }}, \quad \EE_{-}= \scalebox{0.8}{\bordermatrix{ 
        &-1 & 0&1\cr
        &\lhd & \ominus & \rhd \cr
        &&\ominus \cr
        &&\ominus \cr
        &&&\ominus
        }}.\]
\end{enumerate}
\end{prop}

The next relevant case is $\pi_{temp} = T_{V,4}^{\pm}(T_{I,1}^{1}(\pi_{sc}))$. 
\begin{prop}\label{nontempA17}
    Consider the representation 
\begin{equation*}
    \pi_x^{\pm} = L(\Delta_{\rho}[-x,-x]; \pi_{temp}^{\pm}),
\end{equation*}
where $\pi_{temp}^{\pm} = T_{V,4}^{\pm}(T_{I,1}^{1}(\pi_{sc}))$ for $x > 0$ and $\alpha = 1$, as described in Proposition \ref{temp3,18}. 
\begin{enumerate}[label = (\roman*)]
    \item The representation $\pi_{x}^{\pm}$ is of Arthur type if and only if $x = 1$. 
    \item The representation $\pi_{x}^{\pm}$ is of critical type for $x \in \{1,2\}$. 
    \item We have $\pi(\EE_{\pm}) = \pi_x^{\pm}$, where 
    \[\EE_{+}= \scalebox{0.8}{\bordermatrix{ 
        &-1 & 0&1\cr
        &\lhd & \oplus & \rhd \cr
        &&\oplus \cr
        &&\oplus \cr
        &&\oplus \cr
        &&&\oplus
        }}, \quad \EE_{-}= \scalebox{0.8}{\bordermatrix{ 
        &-1 & 0&1\cr
        &\lhd & \ominus & \rhd \cr
        &&\ominus \cr
        &&\ominus \cr
        &&\ominus \cr
        &&&\ominus
        }}.\]
\end{enumerate}
\end{prop}

Now we move onto the case $\pi_{temp} = T_{I,2}^{1}(T_{IV,3}(\pi_{sc}))$. 
\begin{prop}\label{nontempA18}
    Consider the representation 
\begin{equation*}
    \pi_x = L(\Delta_{\rho}[-x,-x]; \pi_{temp}),
\end{equation*}
where $\pi_{temp} = T_{I,2}^{1}(T_{IV,3}(\pi_{sc}))$ for $x > 0$ and $\alpha = 1$, as described in Proposition \ref{temp3,19}. 
\begin{enumerate}[label = (\roman*)]
    \item The representation $\pi_x$ is of Arthur type if and only if $x = 1$. 
    \item The representation $\pi_x$ is of critical type when $x \in \{1,2\}$.  
    \item We have $\pi(\EE) = \pi_1$, where 
    \[\EE= \scalebox{0.8}{\bordermatrix{ 
        &-1 & 0&1\cr
        &\lhd & \odot & \rhd \cr
        &&&\odot \cr
        &&&\odot
        }}.\]
\end{enumerate}
\end{prop}

The next case to consider is $\pi_{temp} = T_{II,3}^{1}(T_{IV,3}(\pi_{sc}))$. 
\begin{prop}\label{nontempA19}
     Consider the representation 
\begin{equation*}
    \pi_x = L(\Delta_{\rho}[-x,-x]; \pi_{temp}),
\end{equation*}
where $\pi_{temp} = T_{II,3}^{1}(T_{IV,3}(\pi_{sc}))$, for $x > 0$ and $\alpha \in \mathbb{Z}_{>0}$, as described in Proposition \ref{temp3,20}. 
\end{prop}
\begin{enumerate}[label = (\roman*)]
    \item The representation $\pi_x$ is of Arthur type if and only if $1 \leq x \leq \alpha$. 
    \item The representation $\pi_x$ is of critical type when $(x,\alpha) = (1,1),(2,1),(2,2)$. 
    \item Define 
    \begin{flalign*}
        \EE_x := &\{([x,-x]_\rho,\lfloor x \rfloor, (-1)^{x}\eta),
        ([x -2, 0]_\rho, 0,-\eta), \\
        &([1,1]_\rho,0,-\eta)^2, ([\alpha-1,x]_\rho,0,(-1)^{x}\eta)\}.
    \end{flalign*}
    Then $\pi(\EE_x) = \pi_x$. Here is the associated symbol. 
     \[\EE_x= \scalebox{0.8}{\bordermatrix{ 
             &-x  & -(x-1)& \cdots & 0 & 1& \cdots  & x-2 & x-1 &x & \cdots  &\alpha-1  \cr
            &\lhd  &\lhd & \cdots & \odot &\odot & \cdots & \rhd & \rhd&\rhd \cr
            &&&&\odot& \odot &\cdots& \odot \cr
            &&&&&\odot \cr
            &&&&&\odot \cr
            &&&&&&&&&\odot \cr
            &&&&&&&&&&\ddots \cr
            &&&&&&&&&&& \odot \cr
        }}.\]
\end{enumerate}

There are four more cases to consider in Case $(A)$. Next we have $\pi_{temp} = T_{I,2}^{1}(T_{V,2}^{\pm}(\pi_{sc}))$.
\begin{prop}\label{nontempA20}
    Consider the representation 
\begin{equation*}
    \pi_x^{\pm} = L(\Delta_{\rho}[-x,-x]; \pi_{temp}^{\pm}),
\end{equation*}
where $\pi_{temp}^{\pm} = T_{I,2}^{1}(T_{V,2}^{\pm}(\pi_{sc}))$ for $x > 0$ and $\alpha = 0$, as described in Proposition \ref{temp3,21}.
\begin{enumerate}[label = (\roman*)]
    \item The representation $\pi_x^{\pm}$ is not of Arthur type for any $x$. 
    \item The representation $\pi_{x}^{\pm}$ is of critical type when $x \in \{1,2\}$. 
\end{enumerate}
\end{prop}
\begin{proof} This follows from the fact that the set $\Psi(\pi_{temp}^{\pm})$ is a singleton. \end{proof}

The next case is $\pi_{temp} = T_{III,2}^{1}(\pi_{sc})$. 
\begin{prop}\label{nontempA21}
    Consider the representation 
\begin{equation*}
    \pi_x = L(\Delta_{\rho}[-x,-x]; \pi_{temp}),
\end{equation*}
where $\pi_{temp} = T_{III,2}^{1}(\pi_{sc})$ for $\alpha = 1$, as described in Proposition \ref{temp3,22}. 
\begin{enumerate}[label = (\roman*)]
    \item The representation $\pi_x$ is of Arthur type if and only if $x \in \{1,2\}$. 
    \item The representation $\pi_x$ is of critical type when $x \in \{1,2\}$. 
    \item We have that $\pi(\EE_x) = \pi_x$, where
    \[\EE_{1}= \scalebox{0.8}{\bordermatrix{ 
        &-1 & 0&1\cr
        &\lhd & \odot & \rhd \cr
        &&&\odot \cr
        &&&\odot \cr
        }}, \quad \EE_{2}= \scalebox{0.8}{\bordermatrix{ 
        &-2&-1 & 0&1 &2\cr
        &\lhd &\lhd & \odot& \rhd &\rhd\cr
        &&&\odot \cr
        &&&&\odot \cr
        }}.\]
\end{enumerate}
\end{prop}
\begin{proof} By Lemma \ref{absolutevalue}, the only possible $x$ for $\pi_x$ to be of Arthur type are $x = 0$ or $1$. It is easy to verify that both cases hold. The rest follows from definition. \end{proof}

We continue to the case $\pi_{temp} = T_{I,3}^{\frac{1}{2}}(\pi_{sc})$, which is well-defined if and only if $\alpha = \frac{1}{2}$ by Proposition \ref{temp3,22}. 

\begin{prop}
    Consider the representation 
\begin{equation*}
    \pi_x = L(\Delta_{\rho}[-x,-x]; \pi_{temp}),
\end{equation*}
where $\pi_{temp} = T_{I,3}^{\frac{1}{2}}(\pi_{sc})$ and $\alpha = \frac{1}{2}$, as described in Proposition \ref{temp3,22}. 
\begin{enumerate}[label = (\roman*)]
    \item The representation $\pi_{x}$ is of Arthur type if and only if $x = \frac{1}{2}$.
    \item The representation $\pi_{x}$ is of critical type if and only if $x = \frac{1}{2}$  or $x = \frac{3}{2}$. 
    \item 
    W have $\pi(\EE_{\frac{1}{2}}) = \pi_{\frac{1}{2}}$, where 
    \[\EE_{\frac{1}{2}}= \scalebox{0.8}{\bordermatrix{ 
             &-\frac{1}{2} & \frac{1}{2} \cr
             &\lhd & \rhd \cr
             &&\ominus \cr
             &&\ominus \cr
             &&\ominus \cr
        }}.\]
\end{enumerate}
\end{prop}
\begin{proof}
    Part $(i)$ follows directly from Theorem \ref{nontemp red}. Parts $(ii)$ and $(iii)$ follow from definition. 
\end{proof}

Two more cases remain in Case $(A)$. The next one is $\pi_{temp} = T_{IV,7}(\pi_{sc})$. The proofs of Propositions \ref{nontempA22} to \ref{nontempA23} below are similar to that of Proposition \ref{nontempA21}, which we omit.

\begin{prop}\label{nontempA22}
     Consider the representation 
\begin{equation*}
    \pi_x = L(\Delta_{\rho}[-x,-x]; \pi_{temp}),
\end{equation*}
where $\pi_{temp} = T_{IV,7}(\pi_{sc})$ for $x > 0$ and $\alpha \in \mathbb{Z}_{>0}$, as described in Proposition \ref{temp3,26}. 
\begin{enumerate}[label = (\roman*)]
    \item The representation $\pi_x$ is of Arthur type if and only if $1 \leq x \leq \alpha$.
    \item The representation is of critical type when $(x,\alpha) = (1,1)$. 
    \item Define 
    \begin{flalign*}
        \EE_x := &\{([x,-x]_\rho,\lfloor x \rfloor, (-1)^{x}\eta), ([0,0]_\rho,0,\eta)^6,
        ([x -2, 0]_\rho, 0,-\eta), \\
        &([\alpha-1,x]_\rho,0,(-1)^{x}\eta)\}.
    \end{flalign*}
    Then $\pi(\EE_x) = \pi_x$. Here is the associated symbol. 
     \[\EE_x= \scalebox{0.8}{\bordermatrix{ 
             &-x  & -(x-1)& \cdots & 0 & \cdots  & x-2 & x-1 &x & \cdots  &\alpha-1  \cr
            &\lhd  &\lhd & \cdots & \odot &\odot & \cdots & \rhd & \rhd&\rhd \cr
            &&&&\odot \cr
            &&&&\odot \cr
            &&&&\odot \cr
            &&&&\odot \cr
            &&&&\odot \cr
            &&&&\odot \cr
            &&&&\odot &\dots& \odot \cr
            &&&&&&&&\odot &\cdots &\odot
        }}.\]
\end{enumerate}
\end{prop}

The final case remaining in Case $(A)$ is $\pi_{temp} = T_{V,6}^{\pm}(\pi_{sc})$. 
\begin{prop}\label{nontempA23}
     Consider the representation 
\begin{equation*}
    \pi_x^{\pm} = L(\Delta_{\rho}[-x,-x]; \pi_{temp}^{\pm}),
\end{equation*}
where $\pi_{temp}^{\pm} = T_{V,6}^{\pm}(\pi_{sc})$, for $x > 0$ and $\alpha = 0$, as described in Proposition \ref{temp3,27}. 
\begin{enumerate}[label = (\roman*)]
    \item The representation $\pi_x^{\pm}$ is of Arthur type if and only if $x = 1$. 
    \item The representation $\pi_x^{\pm}$ is of critical type when $x = 1$. 
    \item We have $\pi(\EE_{\pm}) = \pi_x^{\pm}$, where 
    \[\EE_{+}= \scalebox{0.8}{\bordermatrix{ 
        &-1 & 0&1\cr
        &\lhd & \oplus & \rhd \cr
        &&\oplus \cr
        &&\oplus \cr
        &&\oplus \cr
        &&\oplus\cr
        &&\oplus\cr
        }}, \quad \EE_{-}= \scalebox{0.8}{\bordermatrix{ 
        &-1 & 0&1\cr
        &\lhd & \ominus & \rhd \cr
        &&\ominus \cr
        &&\ominus \cr
        &&\ominus \cr
        &&\ominus \cr
        &&\ominus \cr
        }}.\]
\end{enumerate}
\end{prop}

This concludes our discussion of Case $(A)$. Now let 's move onto Case (B).

\subsection{\texorpdfstring{Case $(B): \pi = L(\Delta_{\rho}[-x,-x-1]; \pi_{temp})$}{}}
In this case, we are interested in the representation $\pi_{[-x,-x-1]} := L(\Delta_{\rho}[-x,-x-1]; \pi_{temp})$ where $\pi_{temp}$ is a corank 2 tempered representation of good parity. By Langlands classification, we have the restriction $x > -\frac{1}{2}$. 

Similarly to the classification of tempered representations in \S \ref{sec 3.1}, there are $10$ total cases to consider. We start with the case $\pi_{temp} = T_{I,1}^{\alpha +1}(T_{I,1}^{\alpha}(\pi_{sc}))$. 

\begin{prop}\label{nontempB1}
    Consider the representation 
    \begin{equation*}
        \pi_{[-x,-x-1]} = L(\Delta_{\rho}[-x, -x-1];\pi_{temp}),
    \end{equation*}
    for $x > -\frac{1}{2}$, where $\pi_{temp} = T_{I,1}^{\alpha +1}(T_{I,1}^{\alpha}(\pi_{sc})) $ for $\alpha > 0$. 
    \begin{enumerate}[label = (\roman*)]
        \item The representation $\pi_{[-x,-x-1]}$ is of Arthur type if and only if $\epsilon_\rho \leq x \leq \alpha -2$ when $\alpha \geq 2$, or when $(x,\alpha) = (0,1)$.
        \item The representation $\pi_{[-x,-x-1]}$ is of critical type if and only if $x \in \{\alpha -2, \alpha -1, \alpha, \alpha +1, \alpha +2\}$. 
        \item Define $\EE_{[-x,-x-1]}$ in various cases as follows. When $\alpha \geq 2$, let 
        \begin{flalign*}
            \EE_{[-x,-x-1]} := &\{([x+1,-x]_\rho, \lfloor x \rfloor, (-1)^x \eta),([x-2,\epsilon_\rho]_\rho,0,\eta), \\
            &([\alpha -2,x+1]_\rho,0,(-1)^{x+1 - \epsilon_\rho}\eta),([\alpha +1, \alpha+1]_\rho,0,(-1)^{\alpha - \epsilon_\rho} \eta)\}.
        \end{flalign*}
        Here is the associated symbol. 
        \[\EE_{[-x,-x-1]}= \scalebox{0.8}{\bordermatrix{ 
             &-x  & -(x-1)&\cdots & \epsilon_\rho & \cdots  & x-2 & x-1 &x &x+1 &\cdots  &\alpha-2 & \alpha -1 & \alpha & \alpha +1  \cr
            &\lhd  &\lhd & \cdots & \odot &\cdots & \rhd & \rhd & \rhd&\rhd \cr
            &&&&\odot& \cdots &\odot \cr
            &&&&&&&&&\odot&\cdots & \odot  \cr
            &&&&&&&&&&&&&& \odot \cr
        }}.\]
        When $(x,\alpha) = (0,1)$ , let 
        \[\EE_{[0,-1]}= \scalebox{0.8}{\bordermatrix{ 
        &0 & 1&2\cr
        &\lhd &\odot & \rhd \cr
        &&&\odot \cr
        }}.\]
        Then we have $\pi(\EE_{[-x,-x-1]}) = \pi_{[-x,-x-1]}(\pi_{temp})$.
        
    \end{enumerate}
\end{prop}
\begin{proof} The sufficient direction for part $(i)$ can be proven in a similar way as Proposition \ref{nontempA1}, which we omit. For the necessary direction, we see that \cite[Lemma $7.8$]{HJLLZ24} gives us the restriction that $\pi_{[-x,-x-1]}$ is of Arthur type only if
\begin{enumerate}
    \item $x = 0$,
    \item $x = \frac{1}{2}$ and $|\Omega|(\pi_{temp})_{\rho}$ contains $\{\rho\lvert \cdot \rvert^{\frac{1}{2}}\}$,
    \item $x > \frac{1}{2}$ and $|\Omega|(\pi_{temp})_{\rho}$ contains $\{\rho\lvert \cdot \rvert^{x-1}, \rho\lvert \cdot \rvert^{x}\}$.
\end{enumerate}
Since $|\Omega|(\pi_{temp})_{\rho}  = \{\rho\lvert \cdot \rvert^y: \epsilon_\rho \leq y \leq \alpha -2\} \cup \{\rho\lvert \cdot \rvert^{\alpha+1}\}$, we can easily check that the three cases combined gives us the criterion in $(i)$. This proves part $(i)$. Part $(ii)$ and $(iii)$ follow from definition. This proves the proposition. \end{proof}

Now we move onto the case $\pi_{temp} = T_{I,1}^{\alpha -1}(T_{I,1}^{\alpha}(\pi_{sc}))$. 
\begin{prop}\label{nontempB2}
    Consider the representation 
    \begin{equation*}
        \pi_{[-x,-x-1]} = L(\Delta_{\rho}[-x, -x-1];\pi_{temp}),
    \end{equation*}
    for $x > -\frac{1}{2}$, where $\pi_{temp} = T_{I,1}^{\alpha -1}(T_{I,1}^{\alpha}(\pi_{sc}))$ and $\alpha > 1$. 
    \begin{enumerate}[label = (\roman*)]
        \item The representation $\pi_{[-x,-x-1]}$ is of Arthur type if and only if $\epsilon_\rho \leq x \leq \alpha -3$ for $\alpha \geq 3$, or when $(x,\alpha) = (0,2)$.
        \item The representation $\pi_{[-x,-x-1]}$ is of critical type if and only if $x \in \{\alpha -2, \alpha -1, \alpha, \alpha +1\}$. 
        \item Define $\EE_{[-x,-x-1]}$ in various cases as follows. When $\alpha \geq 3$, let 
        \begin{flalign*}
            \EE_{[-x,-x-1]} := &\{([x+1,-x]_\rho, \lfloor x \rfloor, (-1)^x \eta),([x-2,\epsilon_\rho]_\rho,0,\eta), \\
            &([\alpha -3,x+1]_\rho,0,(-1)^{x+1 - \epsilon_\rho}\eta),([\alpha -1, \alpha-1]_\rho,0,(-1)^{\alpha - \epsilon_\rho }\eta), 
            \\&([\alpha,\alpha]_\rho,0,(-1)^{\alpha -1 -\epsilon_\rho}\eta)\}.
        \end{flalign*}
         Here is the associated symbol. 
        \[\EE_{[-x,-x-1]}= \scalebox{0.8}{\bordermatrix{ 
             &-x  & -(x-1)&\cdots & \epsilon_\rho & \cdots  & x-2 & x-1 &x &x+1 &\cdots  &\alpha-3 & \alpha-2 & \alpha -1 & \alpha  \cr
            &\lhd  &\lhd & \cdots & \odot &\dots & \rhd & \rhd & \rhd&\rhd \cr
            &&&&\odot& \cdots &\odot \cr
            &&&&&&&&&\odot&\cdots & \odot  \cr
            &&&&&&&&&&&&& \odot \cr
            &&&&&&&&&&&&&& \odot \cr
        }}.\]
        When $(x,\alpha) = (0,2)$ , let 
        \[\EE_{[0,-1]}= \scalebox{0.8}{\bordermatrix{ 
        &0 & 1&2\cr
        &\lhd &\odot & \rhd \cr
        &&\odot \cr
        &&&\odot \cr
        }}.\]
        Then we have $\pi(\EE_{[-x,-x-1]}) = \pi_{[-x,-x-1]}(\pi_{temp})$.
    \end{enumerate}
\end{prop}
\begin{proof} The proof is similar to that of  Proposition \ref{nontempB1}, except that we also have to eliminate the case where $x = \alpha$ for $\alpha \geq 3$. This can be done in the same way as the proof of Proposition \ref{nontempA1}, which we omit. \end{proof}

The next case is $\pi_{temp} = T_{IV,3}(T_{I,1}^{\alpha}(\pi_{sc}))$. 
\begin{prop}\label{nontempB3}
    Consider the representation 
    \begin{equation*}
        \pi_{[-x,-x-1]} = L(\Delta_{\rho}[-x, -x-1];\pi_{temp}),
    \end{equation*}
    for $x > -\frac{1}{2}$, where $\pi_{temp} =T_{IV,3}(T_{I,1}^{\alpha}(\pi_{sc})) $ and $\alpha \in \mathbb{Z}_{>1}$.
    \begin{enumerate}[label = (\roman*)]
         \item The representation $\pi_{[-x,-x-1]}$ is of Arthur type if and only if $0 \leq x \leq \alpha -2$. 
         \item The representation $\pi_{[-x,-x-1]}$ is of critical type if and only if $(x,\alpha) = (0,2)$ or $(1,3)$. 
         \item Define 
        \begin{flalign*}
            \EE_{[-x,-x-1]} := &\{([x+1,-x]_\rho, \lfloor x \rfloor, (-1)^x \eta), ([0,0]_\rho,0,\eta)^2,([x-2,0]_\rho,0,\eta), \\
            &([\alpha -2,x+1],0,(-1)^{x+1 - \epsilon_\rho}\eta), 
            ([\alpha,\alpha]_\rho,0,(-1)^{\alpha -1 -\epsilon_\rho}\eta)\}.
        \end{flalign*}
        Then $\pi(\EE_{[-x,-x-1]}) = \pi_{[-x,-x-1]}(\pi_{temp})$. Here is the associated symbol.
        \[\EE_{[-x,-x-1]}= \scalebox{0.8}{\bordermatrix{ 
             &-x  & -(x-1)&\cdots & 0 & \cdots  & x-2 & x-1 &x &x+1 &\cdots  & \alpha-2 & \alpha -1 & \alpha  \cr
            &\lhd  &\lhd & \cdots & \odot &\cdots & \rhd & \rhd & \rhd&\rhd \cr
            &&&&\odot \cr
            &&&&\odot \cr
            &&&&\odot& \cdots &\odot \cr
            &&&&&&&&&\odot&\cdots & \odot  \cr
            &&&&&&&&&&&&& \odot \cr
        }}.\]
    \end{enumerate}
\end{prop}
\begin{proof} The proof is similar to that of Proposition \ref{nontempB1}, which we omit. \end{proof}

Let's look at the next case, which is $\pi_{temp} = T_{V,2}^{\pm}(T_{I,1}^{1}(\pi_{sc}))$. 
\begin{prop}\label{nontempB4}
    Consider the representation 
    \begin{equation*}
        \pi_{[-x,-x-1]}^{\pm} = L(\Delta_{\rho}[-x, -x-1];\pi_{temp}^{\pm}),
    \end{equation*}
    for $x > -\frac{1}{2}$, where $\pi_{temp}^{\pm} = T_{V,2}^{\pm}(T_{I,1}^{1}(\pi_{sc}))$ and $\alpha =1$. 
    \begin{enumerate}[label = (\roman*)]
        \item The representation $\pi_{[-x,-x-1]}^{\pm}$ is of Arthur type if and only if $x = 0$, or $x = 1$ and $\epsilon_{sc}(\rho \otimes S_3) = \mp 1$. 
        \item The representation $\pi_{[-x,-x-1]}^{\pm}$ of critical type if and only if $x \in \{0,1,2\}$. 
        \item When $x = 0$, we have $\pi(\EE_{[0,-1]}^{\pm}) = \pi_{[0,-1]}^{\pm}$, where 
        \[\EE^{+}_{[0,-1]}= \scalebox{0.8}{\bordermatrix{ 
         &0&1\cr
        &\oplus \cr
        &\oplus &\ominus \cr
        &\ominus&\oplus\cr
        }}, \quad \EE^{-}_{[0,-1]}= \scalebox{0.8}{\bordermatrix{ 
        & 0&1\cr
        &\ominus \cr
        &\ominus &\oplus\cr
        &\oplus&\ominus\cr
        }}.\]
        When $x = 1$ and $\epsilon_{sc}(\rho\otimes S_3) = \mp 1$, we have $\pi(\EE_{[-1,-2]}^{\pm}) = \pi_{[-1,-2]}^{\pm}$, where 
        \[\EE^{+}_{[-1,-2]}= \scalebox{0.8}{\bordermatrix{ 
         &-1&0&1&2\cr
        &\lhd&\lhd&\rhd&\rhd \cr
        &&\ominus 
        }}, \quad \EE^{-}_{[-1,-2]}= \scalebox{0.8}{\bordermatrix{ 
        &-1&0&1&2\cr
        &\lhd&\lhd&\rhd&\rhd \cr
        &&\oplus 
        }}.\]
    \end{enumerate}
\end{prop}
\begin{proof} By Lemma \ref{absolutevalue}, the representation $\pi_{[-x,-x-1]}$ can only be of Arthur type when $x = 0$ or $1$. When $x = 0$, it's easy to verify that the representation is of Arthur type. 
      When $x = 1$, Definition \ref{psiminus} and Theorem \ref{nontemp red} tells us that $\pi_{[-1,-2]}$ is of Arthur type if and only if the segment $[1,0]_\rho$ is contained in the extended multi-segment corresponding to $\pi_{temp}$. 
    By \cite[Proposition $10.3$]{HJLLZ24}, this can only happen when you can apply the $ui$ operator on the second and third row, and by \cite[Definition $3.23$]{HLL22}, we have the condition stated. \end{proof}

The next case is $\pi_{temp} = T_{I,1}^{1}(T_{V,2}(\pi_{sc}))$.
\begin{prop}\label{nontempB5}
    Consider the representation 
    \begin{equation*}
        \pi_{[-x,-x-1]}^{\pm} = L(\Delta_{\rho}[-x, -x-1];\pi_{temp}^{\pm}),
    \end{equation*}
    for $x > -\frac{1}{2}$, where $\pi_{temp}^{\pm} = T_{I,1}^{1}(T_{V,2}^{\pm}(\pi_{sc}))$ and $\alpha = 0$. 
    \begin{enumerate}[label = (\roman*)]
        \item The representation $\pi_{[-x,-x-1]}^{\pm}$ is of Arthur type if and only if $x = 0$. 
        \item The representation $\pi_{[-x,-x-1]}^{\pm}$ is of critical type if and only if $x \in \{0,1,2\}$. 
        \item We have $\pi(\EE_{[0,-1]}^{\pm}) = \pi_{[0,-1]}^{\pm}$, where 
         \[\EE^{+}_{[0,-1]}= \scalebox{0.8}{\bordermatrix{ 
         &0&1\cr
        &\oplus&\ominus \cr
        &\ominus &\oplus \cr
        }}, \quad \EE^{-}_{[0,-1]}= \scalebox{0.8}{\bordermatrix{ 
        & 0&1\cr
        &\ominus&\oplus \cr
        &\oplus &\ominus\cr
        }}.\]
    \end{enumerate}
\end{prop}
\begin{proof} By \cite[Proposition 10.5]{HJLLZ24}, the set $\Psi(\pi_{temp}^{\pm})$ is a singleton, and it does not contain the segment $([1,0]_\rho,*,*)$, so when $x = 1$, the representation $\pi_{[-x,-x-1]}$ is not of Arthur type. The rest of the proof is similar to that of Proposition \ref{nontempB4}, which we omit. \end{proof}

Now we move onto the case $\pi_{temp} = T_{I,2}^{\frac{1}{2}}(\pi_{sc})$.
The proofs of Propositions \ref{nontempB6} to \ref{nontempB10} below are similar to that of Proposition \ref{nontempB5}, which we omit.

\begin{prop}\label{nontempB6}
    Consider the representation 
    \begin{equation*}
        \pi_{[-x,-x-1]} = L(\Delta_{\rho}[-x, -x-1];\pi_{temp}),
    \end{equation*}
    for $x > -\frac{1}{2}$, where $\pi_{temp} = T_{I,2}^{\frac{1}{2}}(\pi_{sc})$ and $\alpha = \frac{1}{2}$. 
    \begin{enumerate}[label = (\roman*)]
        \item The representation $\pi_{[-x,-x-1]}$ is of Arthur type if and only if $x = \frac{1}{2}$. 
        \item The representation $\pi_{[-x,-x-1]}$ is of critical type if and only if $x \in \{\frac{1}{2}, \frac{3}{2}\}$.
        \item We have that $\pi(\EE_{[-\frac{1}{2},-\frac{3}{2}]}) = \pi_{[-\frac{1}{2},-\frac{3}{2}]}$, where 
         \[\EE_{[-\frac{1}{2},-\frac{3}{2}]}= \scalebox{0.8}{\bordermatrix{ 
         &-\frac{1}{2}&\frac{1}{2}&\frac{3}{2}\cr
        &\lhd&\oplus&\rhd \cr
        & &\oplus  \cr
        }}.
        \]
    \end{enumerate}
\end{prop}

The next case to consider is $\pi_{temp} = T_{II,3}^{\frac{1}{2}}(\pi_{sc})$. 
\begin{prop}\label{nontempB7}
    Consider the representation 
    \begin{equation*}
        \pi_{[-x,-x-1]} = L(\Delta_{\rho}[-x, -x-1];\pi_{temp}),
    \end{equation*}
    for $x > -\frac{1}{2}$, where $\pi_{temp} = T_{II,3}^{\frac{1}{2}}(\pi_{sc})$ and $\alpha \in \frac{1}{2} + \mathbb{Z}_{>0}$. 
    \begin{enumerate}[label = (\roman*)]
        \item The representation $\pi_{[-x,-x-1]}$ is of Arthur type if and only if $\frac{1}{2} \leq x \leq \alpha -1$. 
        \item The representation $\pi_{[-x,-x-1]}$ is of critical type if and only if $(x,\alpha) = (\frac{1}{2}, \frac{3}{2}),(\frac{3}{2}, \frac{3}{2}),(\frac{3}{2}, \frac{5}{2})$.
        \item Define 
        \begin{flalign*}
            \EE_{[-x,-x-1]} := &\{([x+1,-x]_\rho, \lfloor x \rfloor, (-1)^{x + \frac{1}{2}}),([\frac{1}{2},\frac{1}{2}]_\rho,0,-1)^2 ,([x-2,\frac{1}{2}]_\rho,0,-1), \\
            &([\alpha -1,x+1]_\rho,0,(-1)^{x-\frac{1}{2}})\}.
        \end{flalign*}
        Then $\pi(\EE_{[-x,-x-1]}) = \pi_{[-x,-x-1]}$. Here is the associated symbol. 
        \[\EE_{[-x,-x-1]}= \scalebox{0.8}{\bordermatrix{ 
             &-x  & -(x-1)&\dots & \frac{1}{2} & \cdots  & x-2 & x-1 &x &x+1 &\cdots  &  \alpha -1   \cr
            &\lhd  &\lhd & \cdots & \odot &\cdots & \rhd & \rhd & \rhd&\rhd \cr
            &&&&\odot \cr
            &&&&\odot \cr
            &&&&\odot& \cdots &\odot \cr
            &&&&&&&&&\odot&\cdots & \odot  \cr
        }}.\]
    \end{enumerate}
\end{prop}

Three cases remain in Case $(B)$. The next one is $\pi_{temp} = T_{III,2}^{\frac{1}{2}}(\pi_{sc})$. 
\begin{prop}\label{nontempB8}
    Consider the representation 
    \begin{equation*}
        \pi_{[-x,-x-1]} = L(\Delta_{\rho}[-x, -x-1];\pi_{temp}),
    \end{equation*}
    for $x > -\frac{1}{2}$, where $\pi_{temp} = T_{III,2}^{\frac{1}{2}}(\pi_{sc})$ and $\alpha = \frac{1}{2}$. 
    \begin{enumerate}[label = (\roman*)]
        \item The representation $\pi_{[-x,-x-1]}$ is of Arthur type if and only if $x = \frac{1}{2}$. 
        \item The representation $\pi_{[-x,-x-1]}$ is of critical type if and only if $x \in \{\frac{1}{2}, \frac{3}{2}\}$. 
        \item We have $\pi(\EE_{[-\frac{1}{2},-\frac{3}{2}]}) = \pi_{[-\frac{1}{2},-\frac{3}{2}]}$, where 
         \[\EE_{[-\frac{1}{2},-\frac{3}{2}]}= \scalebox{0.8}{\bordermatrix{ 
        &-\frac{1}{2} & \frac{1}{2}&\frac{3}{2}\cr
        &\lhd &\ominus & \rhd \cr
        &&\ominus \cr
        }}.\]
    \end{enumerate}
\end{prop}

The next remaining case is $\pi_{temp} = T_{IV,5}(\pi_{sc})$. 
\begin{prop}\label{nontempB9}
    Consider the representation 
    \begin{equation*}
        \pi_{[-x,-x-1]} = L(\Delta_{\rho}[-x, -x-1];\pi_{temp}),
    \end{equation*}
    for $x > -\frac{1}{2}$, where $\pi_{temp} = T_{IV,5}(\pi_{sc})$ and $\alpha \in \mathbb{Z}_{>0}$. 
    \begin{enumerate}[label = (\roman*)]
        \item The representation $\pi_{[-x,-x-1]} $ is of Arthur type if and only if $0 \leq x \leq \alpha -1$. 
        \item The representation $\pi_{[-x,-x-1]} $ is of critical type if and only if $(x,\alpha) = (0,1), (1,1)$, or $(1,2)$. 
        \item Define
        \begin{flalign*}
            \EE_{[-x,-x-1]} := &\{([x+1,-x]_\rho, \lfloor x \rfloor, (-1)^{x}\eta),([0,0]_\rho,0,\eta)^4 ,([x-2,0]_\rho,0,\eta), \\
            &([\alpha -1,x+1]_\rho,0,(-1)^{x+1}\eta)\}.
        \end{flalign*}
        Then $\pi(\EE_{[-x,-x-1]}) = \pi_{[-x,-x-1]}$. Here is the associated symbol. 
        \[\EE_{[-x,-x-1]}= \scalebox{0.8}{\bordermatrix{ 
             &-x  & -(x-1)&\cdots & 0 & \cdots  & x-2 & x-1 &x &x+1 &\cdots  &  \alpha -1   \cr
            &\lhd  &\lhd & \cdots & \odot &\cdots & \rhd & \rhd & \rhd&\rhd \cr
            &&&&\odot \cr
            &&&&\odot \cr
            &&&&\odot \cr
            &&&&\odot \cr
            &&&&\odot& \cdots &\odot \cr
            &&&&&&&&&\odot&\cdots & \odot  \cr
        }}.\]
    \end{enumerate}
\end{prop}

The final case to consider in Case $(B)$ is $\pi_{temp} = T_{V,4}^{\pm}(\pi_{sc})$. 
\begin{prop}\label{nontempB10}
     Consider the representation 
    \begin{equation*}
        \pi_{[-x,-x-1]}^{\pm} = L(\Delta_{\rho}[-x, -x-1];\pi_{temp}^{\pm}),
    \end{equation*}
    for $x > -\frac{1}{2}$, where $\pi_{temp
    }^{\pm} =T_{V,4}^{\pm}(\pi_{sc}) $ and $\alpha = 0$. 
    \begin{enumerate}[label = (\roman*)]
        \item The representation $\pi_{[-x,-x-1]}^{\pm}$ is of Arthur type if and only if $x = 0$. 
        \item The representation $\pi_{[-x,-x-1]}^{\pm}$ is of critical type when $x \in \{0,1\}$. 
        \item We have $\pi(\EE_{[0,-1]}^{\pm}) = \pi_{[0,-1]}^{\pm}$, where 
        \[\EE^{+}_{[0,-1]}= \scalebox{0.8}{\bordermatrix{ 
         &0&1\cr
        &\lhd&\rhd \cr
        &\oplus \cr
        &\oplus \cr
        &\oplus \cr
        &\oplus \cr
        }}, \quad \EE^{-}_{[0,-1]}= \scalebox{0.8}{\bordermatrix{ 
        &0&1\cr
        &\lhd&\rhd \cr
        &\ominus \cr
        &\ominus \cr
        &\ominus \cr
        &\ominus \cr
        }}.\]
    \end{enumerate}
\end{prop}

This concludes our discussion of Case (B). Now we move onto Case $(C)$. 

\subsection{\texorpdfstring{Case $(C)$: $\pi = L(\Delta_{\rho}[-x,-x-2];\pi_{temp})$.}{}}

In this section we consider representations of the form $\pi_{[-x,-x-2]} := L(\Delta_{\rho}[-x,-x-2]; \pi_{temp})$, where $\pi_{temp}$ is a corank 1 tempered representation of good parity. By the classification in \cite[Section $8$]{HJLLZ24}, there are three cases to consider. By Langlands classification, we have a natural restriction that $x > -1$. 

The first case is $\pi_{temp} = T_{I,1}^{\alpha}(\pi_{sc})$. 

\begin{prop}\label{nontempC1}
    Consider the representation 
    \begin{equation*}
        \pi_{[-x,-x-2]} = L(\Delta_{\rho}[-x,-x-2]; \pi_{temp}),
    \end{equation*}
    for $x > -1$, where $\pi_{temp} = T_{I,1}^{\alpha}(\pi_{sc})$ and $\alpha > 0$. 
    \begin{enumerate}[label = (\roman*)]
        \item The representation $\pi_{[-x,-x-2]}$ is of Arthur type if and only if one of the following cases hold: 
        \begin{enumerate}
            \item $x = -\frac{1}{2}$,
            \item $x = 0$ and $\alpha \in \mathbb{Z}_{>0} \setminus \{2\}$,
            \item $x = \frac{1}{2}$ and $\alpha \in \frac{1}{2} + \mathbb{Z}_{>2}$,
            \item $\epsilon_\rho + 1 \leq x \leq \alpha -3$ for $\alpha \geq 1$.
        \end{enumerate}
        \item The representation $\pi_{[-x,-x-2]}$ is of critical type if and only if $x \in \{\alpha -3, \alpha -2, \alpha -1, \alpha, \alpha +1\}$. 
        \item 
        Define $\EE_[-x,-x-2]$ in various cases as follows. For $x = -\frac{1}{2}$, let 
        \[\EE_{[\frac{1}{2},-\frac{3}{2}]}= \scalebox{0.8}{\bordermatrix{ 
        &\frac{1}{2}&\frac{3}{2}& \dots &\alpha -2 & \alpha -1 & \alpha\cr
        &\odot \cr
        &\lhd&\rhd \cr
        &&\odot & \dots & \odot  \cr
        &&&&&&\odot \cr
        }}.\]
        For $x = 0$, let 
        \[\EE_{[0,-2]}= \scalebox{0.8}{\bordermatrix{ 
        &0&1&2& \dots &\alpha -2 & \alpha -1 & \alpha\cr
        &\odot \cr
        &\lhd&\odot &\rhd \cr
        &&&\odot & \dots & \odot  \cr
        &&&&&&&\odot \cr
        }}.\]
        For $x = \frac{1}{2}$, let 
        \[\EE_{[-\frac{1}{2},-\frac{5}{2}]}= \scalebox{0.8}{\bordermatrix{ 
        &-\frac{1}{2}&\frac{1}{2}&\frac{3}{2}&\frac{5}{2}& \cdots &\alpha -2 & \alpha -1 & \alpha\cr
        &\lhd &\odot &\odot &\rhd \cr
        &&&&\odot & \cdots & \odot  \cr
        &&&&&&&&\odot \cr
        }}.\]
        For $\epsilon_\rho + 1 \leq x = \alpha -3$  and $\alpha \geq 1$, define
        \begin{flalign*}
            \EE_{[-x,-x-2]} = &\{([x+2,-x]_\rho,0,(-1)^{x - \epsilon_\rho}\eta),([\epsilon_\rho,\epsilon_\rho]_\rho,0,\eta),\\
            &([\alpha -2,x+2]_\rho,0,(-1)^{x-\epsilon_\rho}\eta),([\alpha,\alpha]_\rho,0,(-1)^{\alpha - \epsilon_\rho}\eta)\}.
        \end{flalign*}
        Then we have $\pi(\EE_{[-x,-x-2]}) := \pi_{[-x,-x-2]}$. Here is the associated symbol. 
        \[\EE_{[-x,-x-2]}= \scalebox{0.8}{\bordermatrix{ 
        &-x&\cdots&\epsilon_\rho & \cdots & x+2 & \dots  &\alpha -2 & \alpha -1 & \alpha\cr
        &\lhd &\cdots &\odot &\cdots &\rhd \cr
        &&&\odot   \cr
        &&&&&\odot &\cdots &\odot \cr
        &&&&&&&&&\odot \cr
        }}.\]
        
    \end{enumerate}
\end{prop}
\begin{proof}
Part $(iii)$ and the sufficient direction of part $(i)$ can be proven in a similar way as Proposition \ref{nontempB1}, which we omit. For the necessary condition, we have the following restriction from Lemma \ref{absolutevalue}: 
    $\pi$ \text{is of Arthur type only if} $x = -\frac{1}{2}$ \text{or the set} $\{\rho \lvert \cdot \rvert^{z}\}_{z = |-x+1|}^{x+1}$ \text{lies in} $|\Omega|(\pi_{temp})$. 
    For $\pi_{temp} = T_{I,1}^{\alpha}(\pi_{sc})$ we have that $|\Omega|(\pi_{temp}) = \{\rho\lvert \cdot \rvert^{y}\}_{y = \epsilon_\rho}^{\alpha -2} \cup \{\rho \lvert \cdot \rvert^{\alpha}\}$. Therefore we reduce to the four cases listed in the proposition. Part $(ii)$ follows from definition. 
    \end{proof}
 
The second case in Case $(C)$ is $\pi_{temp} = T_{IV, 3}(\pi_{sc})$. 
\begin{prop}\label{nontempC2}
     Consider the representation 
    \begin{equation*}
        \pi_{[-x,-x-2]} = L(\Delta_{\rho}[-x,-x-2]; \pi_{temp}),
    \end{equation*}
    for $x > -1$, where $\pi_{temp} = T_{IV,3}(\pi_{sc})$ and $\alpha \in \mathbb{Z}_{>0}$.
    \begin{enumerate}[label = (\roman*)]
        \item The representation $\pi_{[-x,-x-2]}$ is of Arthur type if and only if one of the following cases hold: 
        \begin{enumerate}
            \item $x = 0$ and $\alpha \in \mathbb{Z}_{>1}$,
            \item $1 \leq x \leq \alpha -2$ for $\alpha \geq 3$.
        \end{enumerate}
        \item The representation $\pi_{[-x,-x-2]}$ is of critical type if and only if 
        \begin{equation*}
            (x,\alpha) = (0,1),(0,2),(1,1),(1,2),(1,3).
        \end{equation*}
        \item Define $\EE_{[-x,-x-2]}$ in various cases as follows. For $x = 0$ and $\alpha \in \mathbb{Z}_{>1}$, let
        \[\EE_{[0,-2]}= \scalebox{0.8}{\bordermatrix{ 
        &0&1&2& \cdots  & \alpha -1 \cr
        &\odot \cr
        &\odot \cr
        &\odot \cr
        &\lhd&\odot &\rhd \cr
        &&&\odot & \cdots & \odot  \cr
        }}.\]
        For $1 \leq x \leq \alpha -2$ and $\alpha \geq 3$, define
        \begin{flalign*}
            \EE_{[-x,-x-2]} := &\{([x+2,-x]_\rho,0,(-1)^{x}\eta),([0,0]_\rho,0,\eta)^2,\\
            &([\alpha -1,x+2]_\rho,0,(-1)^{x}\eta)\}.
        \end{flalign*}
        Then we have $\pi(\EE_{[-x,-x-2]}) = \pi_{[-x,-x-2]}$. Here is the associated symbol. 
        \[\EE_{[-x,-x-2]}= \scalebox{0.8}{\bordermatrix{ 
        &-x&\cdots&0 & \dots & x+2 & \cdots  & \alpha -1 \cr
        &\lhd &\cdots &\odot &\cdots &\rhd \cr
        &&&\odot   \cr
        &&&\odot   \cr
        &&&&&\odot &\cdots &\odot \cr
        }}.\]
    \end{enumerate}
\end{prop}
\begin{proof} The proof is similar to that of Proposition \ref{nontempC1}, which we omit.  \end{proof}

The final case in Case $(C)$ is $\pi_{temp} = T_{V,2}^{\pm}(\pi_{sc})$. 
\begin{prop}\label{nontempC3}
     Consider the representation 
    \begin{equation*}
        \pi_{[-x,-x-2]}^{\pm} = L(\Delta_{\rho}[-x,-x-2]; \pi_{temp}^{\pm}),
    \end{equation*}
    where $\pi_{temp}^{\pm} = T_{V,2}^{\pm}(\pi_{sc})$ for $x > -1$ and $\alpha = 0$. 
    \begin{enumerate}[label = (\roman*)]
        \item The representation $\pi_{[-x,-x-2]}^{\pm}$ is not of Arthur type for any $x$. 
        \item The representation $\pi_{[-x,-x-2]}^{\pm}$ is of critical type if and only if $x \in \{0,1\}$. 
    \end{enumerate}
\end{prop}
\begin{proof} This follows from Lemma \ref{absolutevalue} and the fact that $|\Omega|(\pi_{temp}^{\pm})$ is the singleton $\{\rho\}$. \end{proof}

This concludes our discussion of Case $(C)$. Let us now move onto the final case when $f(\pi) = 1$. 

\subsection{\texorpdfstring{Case $(D): \pi = L(\Delta_{\rho}[-x,-x-3]; \pi_{sc})$}{}}
In this section, we are interested in representations of the form $\pi_{[-x,-x-3]} := L(\Delta_{\rho}[-x,-x-3]; \pi_{sc})$. The natural restriction from Langlands classification in this case is that $x > -\frac{3}{2}$. 

Since we are working with supercuspidal representations, we introduce the following result to simplify our computations. 
\begin{prop}[{\cite[Proposition $8.4$]{HJLLZ24}}]\label{scArthur}
    The representation (of good parity)
    \begin{equation*}
        \pi = L(\Delta_{\rho}[x,-y]; \pi_{sc})
    \end{equation*}
    is of Arthur type if and only if $x - y = -1$ or $y \leq \alpha$. Note that $x$ is not arbitrary since we require $x - y < 0$ for the Langlands classification. 
\end{prop}

From this, we can easily derive the following result. 
\begin{prop}\label{nontempD1}
    Consider the representation 
    \begin{equation*}
        \pi_{[-x,-x-3]} = L(\Delta_{\rho}[-x,-x-3]; \pi_{sc}), 
    \end{equation*}
    where $x > -\frac{3}{2}$. 
    \begin{enumerate}
        \item The representation $\pi_{[-x,-x-3]}$ is of Arthur type if and only if $x = -1$ or $-\frac{1}{2} \leq x \leq \alpha -3$.
        \item The representation $\pi_{[-x,-x-3]}$ is of critical type if and only if $x \in \{\alpha -3, \alpha -2, \alpha -1, \alpha\}$. 
        \item Define $\EE_{sc}$ and $\EE_{[-x,-x-3]}$ in various cases as follows. Then $\pi(\EE_{sc}) = \pi_{sc}$ and $\pi(\EE_{[-x,-x-3]}) = \pi_{[-x,-x-3]}$. When $x = -1$, let
        \[\EE_{sc}= \scalebox{0.8}{\bordermatrix{ 
        &0 & 1 &2 & \dots & \alpha -1 \cr
        &\odot &\odot \cr
        &&&\odot&\dots &\odot   \cr
        }},\]

        \[\EE_{[-x,-x-3]}= \scalebox{0.8}{\bordermatrix{ 
        &-1&0 &1 & 2 & \cdots & \alpha -1 \cr
        &\lhd &\odot &\odot &\rhd \cr
        &&&&\odot   \cr
        &&&&&\ddots   \cr
        &&&&&&\odot  \cr
        }}.\]
        When $x = -\frac{1}{2}$, let 
         \[\EE_{sc}= \scalebox{0.8}{\bordermatrix{ 
        &\frac{1}{2} & \dots  & \alpha -1 \cr
        &\odot &\dots &\odot
        }},\]

        \[\EE_{[-x,-x-3]}= \scalebox{0.8}{\bordermatrix{ 
        &\frac{1}{2} & \frac{3}{2} &\dots &\alpha-1 \cr
        &\odot \cr
        &\lhd  &\rhd \cr
        &&\odot   \cr
        &&&\ddots   \cr
        &&&&\odot  \cr
        }}.\]
        When $0 \leq x \leq \alpha -3$, define: 
        \begin{flalign*}
            \EE_{sc} := &\{([x+2,-x +1]_\rho, \lfloor x-1 \rfloor, (-1)^{x-\epsilon_\rho} \eta),([x - 4,\epsilon_\rho]_\rho,0,\eta),\\
            &([\alpha -1,x+1]_\rho,0,(-1)^{x+1 - \epsilon_\rho}\eta)\},\\
            \EE_{[-x,-x-3]} := &\{([x+3,-x]_\rho, \lfloor x \rfloor, (-1)^{x-\epsilon_\rho} \eta),([x - 4,\epsilon_\rho]_\rho,0,\eta),\\
            &([\alpha -1,x+1]_\rho,0,(-1)^{x+1 - \epsilon_\rho}\eta)\}.
        \end{flalign*}
        Here are the associated symbols. 
        \[\EE_{sc}= \scalebox{0.8}{\bordermatrix{ 
        &-x + 1 & \cdots & \epsilon_\rho & \cdots & x-4 & x-3 & x-2 & x-1 & x & x+ 1 &x + 2 & \cdots & \alpha -1\cr
        &\lhd & \cdots & \odot & \cdots & \cdots &\dots & \cdots &\cdots &\cdots & \cdots &\rhd \cr
        &&&\odot & \cdots & \odot\cr
        &&&&&&&&&&\odot &\odot & \cdots &\odot 
        }},\]

        \[\EE_{[-x,-x-3]}= \scalebox{0.8}{\bordermatrix{ 
        &-x & \cdots & \epsilon_\rho & \cdots & x-4 & x-3 & x-2 & x-1 & x & x+ 1 &x + 2 &x+3 & \cdots & \alpha -1\cr
        &\lhd & \cdots & \odot & \cdots & \cdots &\cdots & \cdots &\cdots &\cdots & \cdots &\cdots &\rhd \cr
        &&&\odot & \cdots & \odot\cr
        &&&&&&&&&&\odot &\odot & \odot &\cdots &\odot 
        }}.\]
    \end{enumerate}
\end{prop}
\begin{proof} From Proposition \ref{scArthur}, we see directly that $\pi_{[-x,-x-3]}$ is of Arthur type if and only if $x = -1$ or $x + 3 \leq \alpha$. This proves part $(i)$. Part $(ii)$ follows from definition and part $(iii)$ can be proven using the similar kinds of calculation as Proposition \ref{nontempA1} for supercuspidal representations. We omit the details. 
\end{proof}
This concludes our discussion for the case $f(\pi) = 1$. 

\section{Classification of corank 4 non-tempered representations of good parity \texorpdfstring{($f(\pi) = 2$)}{}}\label{classnontempcorank4,2}
In this section we classify corank 4 non-tempered representations with two segments in their $L$-data. Under the restriction imposed by Langlands classification, the following are the  different cases we need to consider: 
\begin{enumerate}[label = (\Alph*)]
    \item $\pi = L(\Delta_{\rho}[-x_1,-x_1],\Delta_{\rho}[-x_2,-x_2]; \pi_{temp})$, where $x_1 > x_2 \geq \frac{1}{2}$ and $\pi_{temp}$ is tempered of corank $2$. 
    \item $\pi = L(\Delta_{\rho}[-x,-x],\Delta_{\rho}[-x,-x]; \pi_{temp})$, where $x \geq \frac{1}{2}$ and $\pi_{temp}$ is tempered of corank $2$. 
    \item $\pi = L(\Delta_{\rho}[-x_1,-x_1 -1], \Delta_{\rho}[-x_2,-x_2]; \pi_{temp})$, where $x_1 > -\frac{1}{2}$, $ x_2 >0$,$x_2 - x_1 \leq \frac{1}{2}$ and $\pi_{temp}$ is tempered of corank $1$. 
    \item $\pi = L(\Delta_{\rho}[-x_1,-x_1], \Delta_{\rho}[-x_2,-x_2-1]; \pi_{temp})$, where $x_1 \geq \frac{1}{2}, x_2 > -\frac{1}{2}$, $x_1 - x_2 \geq \frac{1}{2}$ and $\pi_{temp}$ is tempered of corank $1$. 
    \item $\pi = L(\Delta_{\rho}[-x_1, -x_1 -1], \Delta_{\rho}[-x_2,-x_2-1]; \pi_{sc})$, where $x_1, x_2 > -\frac{1}{2}$ and $x_1 \geq x_2$. 
    \item $\pi = L(\Delta_{\rho}[-x_1,-x_1-2], \Delta_{\rho}[-x_2,-x_2]; \pi_{sc})$, where $x_1 > -1$, $x_2 \geq \frac{1}{2}$, and $x_2 - x_1 \leq 1$. 
    \item $\pi = L(\Delta_{\rho}[-x_1,-x_1],\Delta_{\rho}[-x_2,-x_2-2]; \pi_{sc})$, where $x_1 \geq \frac{1}{2}$, $x_2 > -1$ and $x_1 - x_2 \geq 1$. 
\end{enumerate}

We begin with Case $(A)$. 

\subsection{\texorpdfstring{Case $(A): \pi = L(\Delta_{\rho}[-x_1,-x_1], \Delta_{\rho}[-x_2,-x_2]; \pi_{temp})$}{}}
Note that in this case, $x_1 > x_2 \geq \frac{1}{2}$ and $\pi_{temp}$ is a corank 2 tempered representation of good parity. 
Therefore, similar as before, there are a total of $10$ subcases to consider. Let us start with $\pi_{temp} = T_{I,1}^{\alpha +1}(T_{I,1}^{\alpha}(\pi_{sc}))$. 
\begin{prop}\label{nontemp2A1}
    Consider the representation 
    \begin{equation*}
        \pi_{x_1, x_2} = L(\Delta_{\rho}[-x_1,-x_1], \Delta_{\rho}[-x_2,-x_2]; \pi_{temp}),
    \end{equation*}
    where $x_1 > x_2 \geq \frac{1}{2}$ and $\pi_{temp} = T_{I,1}^{\alpha +1}(T_{I,1}^{\alpha}(\pi_{sc}))$ for $\alpha > 0$. 
    \begin{enumerate}[label = (\roman*)]
        \item The representation $\pi_{x_1,x_2}$ is of Arthur type if and only if $\frac{1}{2} \leq x_2 < x_1 \leq \alpha -1$ or $(x_1,x_2) = (\alpha,\alpha -1)$ when $\alpha \geq \frac{3}{2}$, or when $(x_1, x_2, \alpha) = (\frac{3}{2}, \frac{1}{2}, \frac{1}{2})$.
        \item The representation $\pi_{x_1,x_2}$ if of critical type in the following cases: 
        \begin{itemize}
            \item $(x_1,x_2) = (\alpha-1, \alpha-2)$ for $\alpha \geq \frac{3}{2}$,
            \item $(x_1, x_2) = (\alpha, \alpha -1)$ when $\alpha \geq \frac{3}{2}$,
            \item $(x_1, x_2) = (\alpha+1, \alpha)$, $(\alpha +2, \alpha +1)$, $(\alpha +3, \alpha +2), (\alpha+1, \alpha-1),(\alpha+2, \alpha),(\alpha+2, \alpha-1)$ for $\alpha \geq \frac{1}{2}$.  
        \end{itemize}
        \item Define $\EE_{x_1, x_2}$ in various cases as follows. Then $\pi(\EE_{x_1, x_2}) = \pi_{x_1, x_2}$. When $\alpha \geq \frac{3}{2}$, and $\frac{1}{2} \leq x_2 < x_1 \leq \alpha -1$, define 
        \begin{flalign*}
            \EE_{x_1, x_2}:= &\{([x_1,-x_1]_\rho,\lfloor x_1 \rfloor, (-1)^{x_1 +1 - \epsilon_\rho}\eta),([x_2 -2, \epsilon_\rho]_\rho,0,-\eta),
            \\
            &([\alpha-2,x_2]_\rho,0,(-1)^{x_2+1 - \epsilon_\rho}\eta), ([\alpha +1, \alpha +1]_\rho,0,(-1)^{\alpha - \epsilon_\rho}\eta)\}.
        \end{flalign*}
        Here is the associated symbol. 
        \[\EE_{x_1, x_2}= \scalebox{0.8}{\bordermatrix{ 
        &-x_1  & \cdots & \epsilon_\rho & \cdots & x_2 -2 & x_2 -1 & x_2 & \cdots & x_1 & \cdots &\alpha -2 & \alpha -1 & \alpha & \alpha +1\cr
        &\lhd & \cdots & \odot & \cdots & \cdots &\cdots & \cdots &\cdots  &\rhd \cr
        &&&\odot & \cdots & \odot\cr
        &&&&&&&\odot & \cdots & \odot & \cdots & \odot \cr
        &&&&&&&&&&&&&&\odot 
        }}.\]
        When $\alpha \geq \frac{3}{2}$ and $(x_1,x_2) = (\alpha, \alpha -1)$, define
        \begin{flalign*}
            \EE_{x_1, x_2}:= &\{([x_1,-x_1]_\rho,\lfloor x_1 \rfloor, (-1)^{x_1 +1 - \epsilon_\rho}\eta),([\alpha -3, \epsilon_\rho]_\rho,0,-\eta),
            \\
            &([\alpha +1, \alpha +1]_\rho,0,(-1)^{\alpha  - \epsilon_\rho}\eta)\}.
        \end{flalign*}
        Here is the associated symbol. 
        \[\EE_{x_1, x_2}= \scalebox{0.8}{\bordermatrix{ 
        &-\alpha  & \cdots & \epsilon_\rho & \cdots & \alpha -3 & \alpha -2 &\alpha -1 & \alpha & \alpha +1 \cr
        &\lhd & \cdots & \odot & \cdots & \cdots  & \cdots &\cdots  &\rhd \cr
        &&&\odot & \cdots & \odot\cr
        &&&&&&&&&\odot 
        }}.\]
        Finally, when $(x_1, x_2, \alpha) = (\frac{3}{2}, \frac{1}{2}, \frac{1}{2})$, let 
        \[\EE_{\frac{3}{2}, \frac{1}{2}}= \scalebox{0.8}{\bordermatrix{ 
        &-\frac{3}{2} & -\frac{1}{2} & \frac{1}{2} & \frac{3}{2} \cr
        &\lhd & \lhd & \rhd & \rhd\cr
        &&&&\odot 
        }}.\]
    \end{enumerate}
\end{prop}
\begin{proof} Let us first show the sufficient direction for $(i)$. The case $(x_1,x_2, \alpha) = (\frac{3}{2}, \frac{1}{2}, \frac{1}{2})$ is easy to verify. Now consider the case when $\alpha \geq \frac{3}{2}$ and $ x_2 + 1 = x_1$, $x_2 \leq \alpha -1$. Here we can construct an extended multi-segment $\EE_{x_2}$ of the form: 
\[\EE_{x_2}= \scalebox{0.8}{\bordermatrix{ 
        &-x_2  & \cdots & \epsilon_\rho & \cdots & x_2-2 & x_2-1 & x_2 & \dots &\alpha -2 &\alpha -1 & \alpha & \alpha +1 \cr
        &\lhd & \cdots & \odot & \cdots  & \dots &\cdots  &\rhd \cr
        &&&\odot & \cdots & \odot\cr
        &&&&&&&\odot&\cdots&\odot  \cr
        &&&&&&&&&&&&\odot
        }},\]
such that $\EE_{x_2}$ lies in $\mathscr{E}(\pi_{x_1,x_2}^{\rho,-}; \Delta_{\rho}[-x_1,-x_1], 1)$, and thus by Theorem \ref{nontemp red}, $\pi_{x_1,x_2}$ is of Arthur type with $\pi_{x_1,x_2} = \pi(\EE_{x_2}^{\rho,+})$. Now when $x_2 +1 < x_1 \leq \alpha -1$, we see from the discussion in \cite[Section 9]{HJLLZ24}, that there exists a multi-segment $\EE'$ containing $([\alpha -1,\alpha -1]_\rho,0,*)$ in its bottom row such that 
\begin{equation*}
    \pi(\EE') = L(\Delta_{\rho}[-x_1,-x_1],\Delta_{\rho}[-x_2,-x_2]; \pi_{sc}).
\end{equation*}

Define a new multi-segment $\EE$ which is identical to $\EE'$ except that the bottom row is now replaced with $([\alpha +1,\alpha+1]_{\rho},0,*)$. It follows that $\pi(\EE) = \pi_{x_1,x_2}$, and thus $\pi_{x_1,x_2}$ is of Arthur type. This concludes the proof of the sufficient direction. 

Now we show the necessary direction. By definition we have 
\begin{equation*}
    \pi_{x_1,x_2}^{\rho,-} = L(\Delta_{\rho}[-x_2,-x_2]; \pi_{temp}).
\end{equation*}
Then Theorem \ref{nontemp red} suggests that $\pi_{x_1, x_2}$ is of Arthur type only if $\pi_{x_1,x_2}^{\rho,-}$ is also of Arthur type. By \cite[Proposition $11.1$]{HJLLZ24}, this happens when $\frac{1}{2} \leq x_2 \leq \alpha -1$ when $\alpha \geq \frac{3}{2}$, or when $x_2 = \frac{1}{2} = \alpha$. 

Furthermore, we have
\begin{equation*}
    |\Omega|(\pi_{x_1,x_2}^{\rho, -}) \subseteq \{\rho\lvert \cdot \rvert^y: \epsilon_\rho \leq y \leq \alpha -2\} \cup \{\rho \lvert \cdot \rvert^{\alpha+1}\} \cup \{\rho\lvert \cdot \rvert^{x_2}, \rho\lvert \cdot \rvert^{-x_2}\}.
\end{equation*}

Lemma \ref{absolutevalue} gives the constriction that $\pi_{x_1,x_2}$ is of Arthur type only if $x_1 = \frac{1}{2}$ or $\rho\lvert \cdot \rvert^{x_1-1} \in |\Omega|(\pi^{\rho,-})$. Since $x_1 > \frac{1}{2}$ by definition, we may conclude that $\pi^{\rho,-}$ is of Arthur type only if one of the following holds: 
\begin{itemize}
    \item $x_2 + 1 \leq x_1 \leq \alpha -1$, or
    \item $x_1 = \alpha +2$, or
    \item $x_2 = \alpha -1$, $x_1 = \alpha$ 
\end{itemize}
It suffices to show that $\pi^{\rho,-}$ is not of Arthur type when $x_1 = \alpha +2$. We resolve this case in the same way as we did before. When $x_1 = \alpha +2$, by Lemma \ref{generalconstraint}, there does not exist any extended multi-segment $\EE$ containing the segment $([\alpha+1, -(\alpha+1)]_\rho,*,*)$ such that $\pi(\EE) = \pi_{x_1,x_2}^{\rho -}$. By Definition \ref{psiminus}, this means that the set $\Psi(\pi_{x_1,x_2}^{\rho,-}; \Delta_{\rho}[x_1,-x_1], 1)$ is empty, so $\pi_{x_1,x_2}$ is not of Arthur type by Proposition \ref{nontemp red}. The rest follows from definition. \end{proof}

Now we move onto the case $\pi_{temp} = T_{I,1}^{\alpha -1}(T_{I,1}^{\alpha}(\pi_{sc}))$. 
\begin{prop}\label{nontemp2A2}
    Consider the representation 
    \begin{equation*}
        \pi_{x_1, x_2} = L(\Delta_{\rho}[-x_1,-x_1], \Delta_{\rho}[-x_2,-x_2]; \pi_{temp}),
    \end{equation*}
    where $x_1 > x_2 \geq \frac{1}{2}$ and $\pi_{temp} = T_{I,1}^{\alpha -1}(T_{I,1}^{\alpha}(\pi_{sc}))$ for $\alpha > 1$. 
    \begin{enumerate}[label = (\roman*)]
        \item The representation $\pi_{x_1, x_2}$ is of Arthur type if and only if $\frac{1}{2} < x_2 < x_1 \leq \alpha -2$ or $(x_1,x_2) = (\alpha -1,\alpha -2)$ when $\alpha \geq \frac{5}{2}$, or when $(x_1,x_2,\alpha) = (\frac{3}{2}, \frac{1}{2}, \frac{3}{2})$. 
        \item The representation is of critical type in the following cases: 
        \begin{itemize}
            \item $(x_1,x_2) = (\alpha -2,\alpha -3)$ for $\alpha \geq \frac{7}{2}$,
            \item $(x_1,x_2) = (\alpha-1, \alpha -2)$ for $\alpha \geq \frac{5}{2}$,
            \item $(x_1, x_2) = (\alpha, \alpha -1), (\alpha +1, \alpha), (\alpha +2, \alpha +1),(\alpha+1, \alpha-1),(\alpha, \alpha-2), (\alpha+1, \alpha-2)$. 
        \end{itemize}
        \item Define $\EE_{x_1, x_2}$ in various cases as follows. Then $\pi(\EE_{x_1, x_2}) = \pi_{x_1, x_2}$. When $\alpha \geq \frac{3}{2}$, and $\frac{1}{2} \leq x_2 < x_1 \leq \alpha -1$, define 
        \begin{flalign*}
            \EE_{x_1, x_2}:= &\{([-x_1,x_1]_\rho,\lfloor x_1 \rfloor, (-1)^{x_1 +1 - \epsilon_\rho}\eta),([x_2 -2, \epsilon_\rho]_\rho,0,-\eta),
            \\
            &([\alpha-3,x_2]_\rho,0,(-1)^{x_2+1 - \epsilon_\rho}\eta), ([\alpha , \alpha ]_\rho,0,(-1)^{\alpha +1 - \epsilon_\rho}\eta)\}.
        \end{flalign*}
        Here is the associated symbol: 
         \[\EE_{x_1, x_2}= \scalebox{0.8}{\bordermatrix{ 
        &-x_1  & \cdots & \epsilon_\rho & \cdots & x_2 -2 & x_2 -1 & x_2 & \cdots & x_1 & \cdots &\alpha -3 &\alpha -2 & \alpha -1 & \alpha \cr
        &\lhd & \cdots & \odot & \cdots & \cdots &\dots & \cdots &\cdots  &\rhd \cr
        &&&\odot & \cdots & \odot\cr
        &&&&&&&\odot & \cdots & \odot & \cdots & \odot \cr
        &&&&&&&&&&&&&&\odot 
        }}.\]
        When $\alpha \geq \frac{3}{2}$ and $(x_1, x_2) = (\alpha -1, \alpha -2)$, define 
        \begin{flalign*}
            \EE_{x_1, x_2}:= &\{([-x_1,x_1]_\rho,\lfloor x_1 \rfloor, (-1)^{x_1 +1 - \epsilon_\rho}\eta),([\alpha -4, \epsilon_\rho]_\rho,0,-\eta),
            \\
            &([\alpha , \alpha ]_\rho,0,(-1)^{\alpha +1  - \epsilon_\rho}\eta)\}.
        \end{flalign*}
        Here is the associated symbol: 
        \[\EE_{x_1, x_2}= \scalebox{0.8}{\bordermatrix{ 
        &-(\alpha-1)  & \cdots & \epsilon_\rho & \cdots &\alpha -4 & \alpha -3 & \alpha -2 &\alpha -1 & \alpha  \cr
        &\lhd & \cdots & \odot & \cdots & \cdots  & \dots &\cdots  &\rhd \cr
        &&&\odot & \cdots & \odot\cr
        &&&&&&&&&\odot 
        }}.\]
        When $(x_1, x_2, \alpha) = (\frac{3}{2}, \frac{1}{2},\frac{3}{2})$, then define $\EE_{x_1,x_2}$ to be the same extended multi-segment defined in Proposition \ref{nontemp2A1}, when $(x_1, x_2, \alpha) = (\frac{3}{2}, \frac{1}{2}, \frac{1}{2})$
    \end{enumerate}
\end{prop}
\begin{proof} To prove the necessary direction, we note that $\pi_{x_1,x_2}^{\rho,-}$ is of Arthur type if and only if $\frac{1}{2} \leq x_2 \leq \alpha -2$ for $\alpha \geq \frac{3}{2}$, or $(x,\alpha) = (\frac{1}{2}, \frac{3}{2})$ by \cite[Proposition $11.3$]{HJLLZ24}. The rest of the proof is similar to that of Proposition \ref{nontemp2A1}, which we omit. \end{proof}

The next case is $\pi_{temp} = T_{IV, 3}(T_{I,1}^{\alpha}(\pi_{sc}))$. 
\begin{prop}\label{nontemp2A3}
    Consider the representation 
    \begin{equation*}
        \pi_{x_1, x_2} = L(\Delta_{\rho}[-x_1,-x_1], \Delta_{\rho}[-x_2,-x_2]; \pi_{temp}),
    \end{equation*}
    where $x_1 > x_2 \geq \frac{1}{2}$ and $\pi_{temp} = T_{IV, 3}(T_{I,1}^{\alpha}(\pi_{sc}))$ for $\alpha \in \mathbb{Z}_{>1}$. 
    \begin{enumerate}[label = (\roman*)]
        \item The representation $\pi_{[x_1,x_2]}$ is of Arthur type if and only if $1 \leq x_2 < x_1 \leq \alpha -1$, or $(x_1,x_2) = (\alpha, \alpha -1)$. 
        \item The representation $\pi_{[x_1,x_2]}$ is of critical type when $(x_1,x_2, \alpha) = (2,1,2)$, $(2,1,3)$, $(3,1,2)$. 
        \item Define $\EE_{x_1, x_2}$ in various cases as follows. Then $\pi(\EE_{x_1,x_2}) = \pi_{x_1, x_2}$. 
        When $1 \leq x_2 < x_1 \leq \alpha -1$, define 
        \begin{flalign*}
            \EE_{x_1, x_2}:= &\{([x_1,-x_1]_\rho,\lfloor x_1 \rfloor, (-1)^{x_1 +1 - \epsilon_\rho}\eta),([0,0]_{\rho},0,-\eta)^2, ([x_2 -2, \epsilon_\rho]_\rho,0,-\eta),\\
            &([\alpha-2,x_2]_\rho,0,(-1)^{x_2+1 - \epsilon_\rho}\eta), ([\alpha , \alpha ]_\rho,0,(-1)^{\alpha +1 - \epsilon_\rho}\eta)\}.
        \end{flalign*}
        Here is the associated symbol. 
        \[\EE_{x_1, x_2}= \scalebox{0.8}{\bordermatrix{ 
        &-x_1  & \cdots & 0 & \cdots & x_2 -2 & x_2 -1 & x_2 & \cdots & x_1 & \cdots &\alpha -2 & \alpha -1 & \alpha \cr
        &\lhd & \cdots & \odot & \cdots & \cdots &\cdots & \cdots &\cdots  &\rhd \cr
        &&&\odot \cr
        &&&\odot \cr
        &&&\odot & \cdots & \odot\cr
        &&&&&&&\odot & \cdots & \odot & \cdots & \odot \cr
        &&&&&&&&&&&&&&\odot 
        }}.\]
        When $(x_1, x_2) = (\alpha, \alpha -1)$, define 
        \begin{flalign*}
            \EE_{x_1, x_2}:= &\{([-x_1,x_1]_\rho,\lfloor x_1 \rfloor, (-1)^{x_1 +1 - \epsilon_\rho}\eta), ([0,0]_\rho,0,-\eta)^2, \\
            &([\alpha -3, \epsilon_\rho]_\rho,0,-\eta),([\alpha , \alpha ]_\rho,0,(-1)^{\alpha +1  - \epsilon_\rho}\eta)\}.
        \end{flalign*}
        Here is the associated symbol. 
        \[\EE_{x_1, x_2}= \scalebox{0.8}{\bordermatrix{ 
        &-\alpha  & \cdots & \epsilon_\rho & \cdots & \alpha -3 & \alpha -2 &\alpha -1 & \alpha  \cr
        &\lhd & \cdots & \odot & \cdots & \cdots  & \cdots &\cdots  &\rhd \cr
        &&&\odot & \cdots & \odot\cr
        &&&&&&&&\odot 
        }}.\]
    \end{enumerate}
\end{prop}
\begin{proof} To prove the necessary direction, we note that $\pi_{x_1,x_2}^{\rho,-}$ is of Arthur type if and only if $1 \leq x_2 \leq \alpha -1$ by \cite[Proposition $11.3$]{HJLLZ24}. The rest of the proof is similar to that of Proposition \ref{nontemp2A1}, which we omit. \end{proof}

Now we move onto the case $\pi_{temp} = T_{V,2}^{\pm}(T_{I,1}^{1}(\pi_{sc}))$. 
\begin{prop}\label{nontemp2A4}
     Consider the representation 
    \begin{equation*}
        \pi_{x_1, x_2}^{\pm} = L(\Delta_{\rho}[-x_1,-x_1], \Delta_{\rho}[-x_2,-x_2]; \pi_{temp}^{\pm}),
    \end{equation*}
    where $x_1 > x_2 \geq \frac{1}{2}$ and $\pi_{temp}^{\pm} = T_{V,2}^{\pm}(T_{I,1}^{1}(\pi_{sc}))$ for $\alpha = 0$.
    \begin{enumerate}[label = (\roman*)]
        \item The representation $\pi_{[x_1,x_2]}^{\pm}$ is of Arthur type if and only if $(x_1, x_2) = (2,1)$. 
        \item The representation $\pi_{[x_1,x_2]}^{\pm}$ is of critical type if and only if  $(x_1, x_2) = (2,1)$ or $(3,2)$. 
        \item We have $\pi(\EE^{\pm}) = \pi_{2,1}^{\pm}$, , where 
         \[\EE^{+}= \scalebox{0.8}{\bordermatrix{ 
         &-2&-1&0&1&2\cr
        &\lhd & \lhd & \oplus & \rhd & \rhd \cr
        &&&\oplus &\ominus \cr
        }}, \quad \EE^{-}= \scalebox{0.8}{\bordermatrix{ 
         &-2&-1&0&1&2\cr
        &\lhd & \lhd & \ominus & \rhd & \rhd\cr
        &&&\ominus &\oplus \cr
        }}.\]
    \end{enumerate}
\end{prop}
\begin{proof} By \cite[Proposition $11.4$]{HJLLZ24}, the representation $(\pi_{x_1, x_2}^{\pm})^{\rho,-} = L(\Delta_{\rho}[-x_2,-x_2]; \pi_{temp})$ is of Arthur type if and only if $x_2 = 1$. The additional requirement from Lemma \ref{absolutevalue} tells us that we must have $x_1 = 2$. By Theorem \ref{nontemp red}, this is the only case where $\pi_{x_1, x_2}^{\pm}$ is of Arthur type. The rest follows from definition. \end{proof}

The next relevant case is $\pi_{temp} = T_{I,1}^{1}(T_{V,2}^{\pm}(\pi_{sc}))$. 
\begin{prop}\label{nontemp2A5}
    Consider the representation 
    \begin{equation*}
        \pi_{x_1, x_2}^{\pm} = L(\Delta_{\rho}[-x_1,-x_1], \Delta_{\rho}[-x_2,-x_2]; \pi_{temp}^{\pm}),
    \end{equation*}
    where $x_1 > x_2 \geq \frac{1}{2}$ and $\pi_{temp}^{\pm} = T_{I,1}^{1}(T_{V,2}^{\pm}(\pi_{sc}))$ for $\alpha = 0$. 
    \begin{enumerate}[label = (\roman*)]
        \item The representation $\pi_{x_1, x_2}^{\pm}$ is of Arthur type if and only if $(x_1, x_2) = (2,1)$. 
        \item The representation $\pi_{[x_1,x_2]}^{\pm}$ is of critical type if and only if $(x_1, x_2) = (2,1)$ or $(3,2)$.  
        \item We have $\pi(\EE^{\pm}) = \pi_{2,1}^{\pm}$, , where 
         \[\EE^{+}= \scalebox{0.8}{\bordermatrix{ 
         &-2&-1&0&1&2\cr
        &\lhd & \lhd & \oplus & \rhd & \rhd \cr
        && &&\oplus \cr
        }}, \quad \EE^{-}= \scalebox{0.8}{\bordermatrix{ 
         &-2&-1&0&1&2\cr
        &\lhd & \lhd & \ominus & \rhd & \rhd\cr
        &&& &\ominus \cr
        }}.\]
    \end{enumerate}
\end{prop}
\begin{proof} By \cite[Proposition $11.5$]{HJLLZ24}, the representation $(\pi_{x_1, x_2}^{\pm})^{\rho,-} = L(\Delta_{\rho}[-x_2,-x_2]; \pi_{temp})$ is of Arthur type if and only if $x_2 = 1$. The rest follows from the same argument as Proposition \ref{nontemp2A4}, which we omit. \end{proof}

We move onto the next case, which is $\pi_{temp} = T_{I,2}^{\frac{1}{2}}(\pi_{sc})$. 
\begin{prop}\label{nontemp2A6}
    Consider the representation 
    \begin{equation*}
        \pi_{x_1, x_2} = L(\Delta_{\rho}[-x_1,-x_1], \Delta_{\rho}[-x_2,-x_2]; \pi_{temp}),
    \end{equation*}
    where $x_1 > x_2 \geq \frac{1}{2}$ and $\pi_{temp} = T_{I,2}^{\frac{1}{2}}(\pi_{sc})$ for $\alpha = \frac{1}{2}$.  
    \begin{enumerate}[label = (\roman*)]
        \item The representation $\pi_{x_1, x_2}$  is of Arthur type if and only if $(x_1, x_2) = (\frac{3}{2}, \frac{1}{2})$. 
        \item The representation $\pi_{x_1, x_2}$ is of critical type if and only if $(x_1, x_2) = (\frac{3}{2}, \frac{1}{2})$ or $(\frac{5}{2}, \frac{3}{2})$. 
        \item We have $\pi(\EE) = \pi_{\frac{3}{2}, \frac{1}{2}}$, where 
         \[\EE= \scalebox{0.8}{\bordermatrix{ 
         &-\frac{1}{2} & \frac{1}{2}\cr
        &\lhd & \rhd \cr
        && \oplus \cr
        && \oplus \cr
        }}.\]
    \end{enumerate}
\end{prop}
\begin{proof} By \cite[Proposition $11.6$]{HJLLZ24}, the representation $(\pi_{x_1, x_2})^{\rho,-} = L(\Delta_{\rho}[-x_2,-x_2]; \pi_{temp})$ is of Arthur type if and only if $x_2 = \frac{1}{2}$. The rest follows from the same argument as Proposition \ref{nontemp2A4}, which we omit. \end{proof}

Now let us consider $\pi_{temp} = T_{II,3}^{\frac{1}{2}}(\pi_{sc})$. 
\begin{prop}\label{nontemp2A7}
     Consider the representation 
    \begin{equation*}
        \pi_{x_1, x_2} = L(\Delta_{\rho}[-x_1,-x_1], \Delta_{\rho}[-x_2,-x_2]; \pi_{temp}),
    \end{equation*}
    where $x_1 > x_2 \geq \frac{1}{2}$ and $\pi_{temp} = T_{II,3}^{\frac{1}{2}}(\pi_{sc})$ and $\alpha \in \frac{1}{2} + \mathbb{Z}_{>0}$. 
    \begin{enumerate}[label = (\roman*)]
        \item The representation $\pi_{x_1, x_2}$ is of Arthur type if and only if $\frac{1}{2} \leq x_2 < x_1 \leq \alpha$ or $(x_1, x_2) = (\alpha +1, \alpha)$.
        \item The representation $\pi_{x_1, x_2}$ is of critical type if and only if $(x_1, x_2, \alpha) = (\frac{3}{2}, \frac{1}{2}, \frac{3}{2}), (\frac{5}{2}, \frac{3}{2}, \frac{3}{2})$ or $(\frac{5}{2}, \frac{3}{2}, \frac{5}{2})$. 
        \item Define $\EE_{x_1, x_2}$ in various cases as follows. Then $\pi(\EE_{x_1, x_2}) = \pi_{x_1, x_2}$. When $\frac{1}{2} \leq x_2 < x_1 \leq \alpha$, define 
        \begin{flalign*}
            \EE_{x_1, x_2}:= &\{([x_1,-x_1]_\rho,\lfloor x_1 \rfloor, 1),([\frac{1}{2}, \frac{1}{2}]_{\rho},0,1)^2, ([x_2 -2, \frac{1}{2}]_\rho,0,1),\\
            &([\alpha-1,x_2]_\rho,0,(-1)^{x_2+ \frac{3}{2}}).
        \end{flalign*}
        Here is the associated symbol. 
        \[\EE_{x_1, x_2}= \scalebox{0.8}{\bordermatrix{ 
        &-x_1  & \cdots & \frac{1}{2} & \cdots & x_2 -2 & x_2 -1 & x_2 & \cdots & x_1 & \cdots & \alpha -1  \cr
        &\lhd & \cdots & \oplus & \cdots & \cdots &\dots & \cdots &\cdots  &\rhd \cr
        &&&\oplus \cr
        &&&\oplus \cr
        &&&\oplus & \cdots & \odot\cr
        &&&&&&&\odot & \cdots & \cdots & \cdots & \odot \cr
        }}.\]
         When $(x_1, x_2) = (\alpha +1, \alpha)$, define 
        \begin{flalign*}
            \EE_{x_1, x_2}:= &\{([-x_1,x_1]_\rho,\lfloor x_1 \rfloor, 1), ([\frac{1}{2}, \frac{1}{2}]_\rho,0,1)^2, ([\alpha - 2, \frac{1}{2}]_\rho,0,1)\}.
        \end{flalign*}
        Here is the associated symbol. 
        \[\EE_{x_1, x_2}= \scalebox{0.8}{\bordermatrix{ 
        &-(\alpha+1)  & \cdots & \frac{1}{2} & \cdots  & \alpha -2 &\alpha -1 & \alpha+1  \cr
        &\lhd & \cdots & \oplus & \cdots  & \dots &\cdots  &\rhd \cr
        &&&\oplus \cr
        &&&\oplus \cr
        &&&\oplus & \cdots & \odot\cr
        }}.\]
    \end{enumerate}
\end{prop}
\begin{proof} By \cite[Proposition $11.7$]{HJLLZ24}, the representation $\pi_{x_1, x_2}^{\rho,-} = L(\Delta_{\rho}[-x_2,-x_2]; \pi_{temp})$ is of Arthur type if and only if $\frac{1}{2} \leq x_2 \leq \alpha$. Since 
\begin{equation*}
    |\Omega|(\pi_{x_1, x_2}^{\rho,-})_\rho \subseteq \{\rho \lvert \cdot \rvert^y: \frac{1}{2} \leq y \leq \alpha -1\} \cup \{\rho\lvert \cdot \rvert^{x_2}, \rho\lvert \cdot \rvert^{-x_2}\},
\end{equation*}
the necessary direction follows from Lemma \ref{absolutevalue} and Proposition \ref{nontemp red}. The rest of the proof is similar to that of Proposition \ref{nontemp2A1}, which we omit. \end{proof}

There are three remaining cases we need to consider. The next one is $\pi_{temp} = T_{III, 2}^{\frac{1}{2}}(\pi_{sc})$. 
\begin{prop}\label{nontemp2A8}
    Consider the representation 
    \begin{equation*}
        \pi_{x_1, x_2} = L(\Delta_{\rho}[-x_1,-x_1], \Delta_{\rho}[-x_2,-x_2]; \pi_{temp}),
    \end{equation*}
    where $x_1 > x_2 \geq \frac{1}{2}$ and $\pi_{temp} = T_{III, 2}^{\frac{1}{2}}(\pi_{sc})$ for $\alpha = \frac{1}{2}$.
    \begin{enumerate}[label = (\roman*)]
        \item The representation $\pi_{x_1, x_2}$ is of Arthur type if and only if $(x_1, x_2) = (\frac{3}{2}, \frac{1}{2})$ or $(\frac{5}{2}, \frac{3}{2})$. 
        \item The representation $\pi_{x_1, x_2}$ is of critical type if and only if $(x_1, x_2) = (\frac{3}{2}, \frac{1}{2})$ or $(\frac{5}{2}, \frac{3}{2})$. 
        \item We have $\pi(\EE) = \pi_{x_1, x_2}$, where 
        \[\EE_{\frac{3}{2}, \frac{1}{2}}= \scalebox{0.8}{\bordermatrix{ 
        &-\frac{3}{2} & -\frac{1}{2} & \frac{1}{2} & \frac{3}{2}  \cr
        &\lhd & \lhd & \rhd &\rhd \cr
        &&&\ominus \cr
        &&&\ominus \cr
        }},\]
        \[\EE_{\frac{5}{2}, \frac{3}{2}}= \scalebox{0.8}{\bordermatrix{ 
        &-\frac{5}{2} & -\frac{3}{2} & -\frac{1}{2} & \frac{1}{2} & \frac{3}{2} &\frac{5}{2} \cr
        &\lhd & \lhd & \oplus & \ominus & \rhd &\rhd \cr
        &&&&\ominus \cr
        }}.\]
    \end{enumerate}
\end{prop}
\begin{proof} By \cite[Proposition $10.8$]{HJLLZ24}, the representation $\pi_{x_1, x_2}^{\rho,-} = L(\Delta_{\rho}[-x_2,-x_2]; \pi_{temp})$ is of Arthur type if and only if $x_2 = \frac{1}{2}$ or $\frac{3}{2}$. The rest of the proof is similar to that of Proposition \ref{nontemp2A4}, which we omit.  \end{proof}

The second last case is $\pi_{temp} = T_{IV, 5}(\pi_{sc})$. 
\begin{prop}\label{nontemp2A9}
    Consider the representation 
    \begin{equation*}
        \pi_{x_1, x_2} = L(\Delta_{\rho}[-x_1,-x_1], \Delta_{\rho}[-x_2,-x_2]; \pi_{temp}),
    \end{equation*}
    where $x_1 > x_2 \geq \frac{1}{2}$ and $\pi_{temp} = T_{IV,5}(\pi_{sc})$ for $\alpha \in \mathbb{Z}_{>0}$. 
    \begin{enumerate}[label = (\roman*)]
        \item The representation $\pi_{x_1, x_2}$ is of Arthur type if and only if $1 \leq x_2 < x_1 \leq \alpha$, or when $(x_1, x_2) = (\alpha +1, \alpha)$. 
        \item The $\pi_{x_1, x_2}$ is of critical type if and only if $(x_1, x_2,\alpha) = (2,1,1)$ or $(2,1,2)$. 
        \item Define $\EE_{x_1, x_2}$ in various cases as follows. Then $\pi(\EE_{x_1, x_2}) = \pi_{x_1, x_2}$. When $1 \leq x_2 < x_1 \leq \alpha$, define 
        \begin{flalign*}
            \EE_{x_1, x_2}:= &\{([x_1,-x_1]_\rho,\lfloor x_1 \rfloor, (-1)^{x_1 +1} \eta),([0,0]_{\rho},0,-\eta)^4, ([x_2 -2, 0]_\rho,0,(-1)^{x_2 -1}\eta),\\
            &([\alpha-1,x_2]_\rho,0,(-1)^{x_2 +1}\eta)\}.
        \end{flalign*}
        Here is the associated symbol. 
        \[\EE_{x_1, x_2}= \scalebox{0.8}{\bordermatrix{ 
        &-x_1  & \cdots & 0 & \cdots & x_2 -2 & x_2 -1 & x_2 & \cdots & x_1 & \cdots & \alpha -1  \cr
        &\lhd & \cdots & \odot & \cdots & \cdots &\cdots & \cdots &\cdots  &\rhd \cr
        &&&\odot \cr
        &&&\odot \cr
        &&&\odot \cr
        &&&\odot \cr
        &&&\odot & \cdots & \odot\cr
        &&&&&&&\odot & \cdots & \cdots & \cdots & \odot \cr
        }}.\]
        When $(x_1, x_2) = (\alpha +1, \alpha)$, define 
        \begin{flalign*}
            \EE_{x_1, x_2}:= &\{([-x_1,x_1]_\rho,\lfloor x_1 \rfloor, (-1)^{\alpha} \eta),([0,0]_{\rho},0,-\eta)^4, ([\alpha -2, 0]_\rho,0,(-1)^{\alpha -1}\eta)\}.
        \end{flalign*}
        Here is the associated symbol. 
        \[\EE_{x_1, x_2}= \scalebox{0.8}{\bordermatrix{ 
        &-(\alpha+1)  & \cdots & 0 & \cdots  & \alpha -2 &\alpha -1 & \alpha+1  \cr
        &\lhd & \cdots & \odot & \cdots  & \cdots &\cdots  &\rhd \cr
        &&&\odot \cr
        &&&\odot \cr
        &&&\odot \cr
        &&&\odot \cr
        &&&\odot & \dots & \odot\cr
        }}.\]
    \end{enumerate}
\end{prop}
\begin{proof} By \cite[Proposition $11.9$]{HJLLZ24}, the representation $\pi_{x_1, x_2}^{\rho,-} = L(\Delta_{\rho}[-x_2,-x_2]; \pi_{temp})$ is of Arthur type if and only if $1 \leq x_2 \leq \alpha$. The rest of the proof is similar to that of Proposition \ref{nontemp2A7}, which we omit. 
\end{proof}

The final case we need to consider in Case $(A)$ is $\pi_{temp} = T_{V,4}^{\pm}(\pi_{sc})$. 
\begin{prop}\label{nontemp2A10}
    Consider the representation 
    \begin{equation*}
        \pi_{x_1, x_2}^{\pm} = L(\Delta_{\rho}[-x_1,-x_1], \Delta_{\rho}[-x_2,-x_2]; \pi_{temp}^{\pm}),
    \end{equation*}
    where $x_1 > x_2 \geq \frac{1}{2}$ and $\pi_{temp}^{\pm} = T_{V,4}^{\pm}(\pi_{sc})$ for $\alpha = 0$. 
    \begin{enumerate}[label = (\roman*)]
        \item The representation $\pi_{x_1, x_2}^{\pm}$ is of Arthur type if and only if $(x_1, x_2) = (2,1)$. 
        \item The representation $\pi_{x_1, x_2}^{\pm}$ is of critical type if and only if $(x_1, x_2) = (2,1)$. 
        \item We have $\pi(\EE^{\pm}) = \pi_{2,1}^{\pm}$, where 
         \[\EE^{+}= \scalebox{0.8}{\bordermatrix{ 
         &-2&-1&0&1&2\cr
        &\lhd & \lhd & \oplus & \rhd & \rhd \cr
        && &\oplus \cr
        && &\oplus \cr
        && &\oplus \cr
        }}, \quad \EE^{-}= \scalebox{0.8}{\bordermatrix{ 
         &-2&-1&0&1&2\cr
        &\lhd & \lhd & \ominus & \rhd & \rhd\cr
        && &\ominus \cr
        && &\ominus \cr
        && &\ominus \cr
        }}.\]
    \end{enumerate}
\end{prop}
\begin{proof} By \cite[Proposition $11.10$]{HJLLZ24}, the representation $(\pi_{x_1, x_2}^{\pm})^{\rho,-} = L(\Delta_{\rho}[-x_2,-x_2]; \pi_{temp}^{\pm})$ is of Arthur type if and only if $x_2 = 1$. The rest of the proof is similar to that of Proposition \ref{nontemp2A5}, which we omit. \end{proof}

This concludes our discussion of Case $(A)$. 

\subsection{\texorpdfstring{Case $(B): \pi = L(\Delta_{\rho}[-x,-x],\Delta_{\rho}[-x,-x]; \pi_{temp})$}{}}
In this subsection, we consider representations of the form $\pi_x = L(\Delta_{\rho}[-x,-x],\Delta_{\rho}[-x,-x]; \pi_{temp})$, where $x \geq \frac{1}{2}$ and $\pi_{temp}$  is a corank 2 tempered representation of good parityy. 
Therefore, as in Case $(A)$, 
there are a total of 10 subcases to consider. 
The first case is $\pi_{temp} = T_{I,1}^{\alpha +1}(T_{I,1}^{\alpha}(\pi_{sc}))$.

\begin{prop}\label{nontemp2B1}
    Consider the representation 
    \begin{equation*}
    \pi_x = L(\Delta_{\rho}[-x,-x],\Delta_{\rho}[-x,-x]; \pi_{temp}),
    \end{equation*}
    for $x \geq \frac{1}{2}$, where $\pi_{temp} =T_{I,1}^{\alpha +1}(T_{I,1}^{\alpha}(\pi_{sc}))$ and $\alpha > 0$. 
    \begin{enumerate}[label = (\roman*)]
        \item The representation $\pi_x$ is of Arthur type if and only if $x = \frac{1}{2}$. 
        \item The representation $\pi_{x}$ is of critical type if and only if $x \in \{\alpha-1, \alpha, \alpha +1, \alpha +2\}$. 
        \item We have $\pi(\EE) = \pi_{\frac{1}{2}}$ where 
        \[\EE= \scalebox{0.8}{\bordermatrix{ 
         &-\frac{1}{2}&\frac{1}{2}& \cdots &\alpha-2 & \alpha -1 & \alpha & \alpha +1\cr
        &\lhd &\rhd \cr
        &\lhd &\rhd \cr
        &&\odot &\cdots &\odot \cr
        &&&&&&&\odot
        }}.\]
    \end{enumerate}
\end{prop}
\begin{proof} The sufficient direction is easy to check. For the necessary direction, by Lemma \ref{absolutevalue}, the representation $\pi_x$ is of Arthur type only if $x = \frac{1}{2}$ or $|\Omega|(\pi_{temp})$ contains $2$ copies of $\rho \lvert \cdot \rvert^{x-1}$, but $|\Omega|(\pi_{temp})$ is multiplicity free, so the conclusion follows. The rest of the proposition follows from definition. \end{proof}

The next case is $\pi_{temp} = T_{I,1}^{\alpha -1}(T_{I,1}^{\alpha}(\pi_{sc}))$, which is similar. 
\begin{prop}\label{nontemp2B2}
     Consider the representation 
    \begin{equation*}
    \pi_x = L(\Delta_{\rho}[-x,-x],\Delta_{\rho}[-x,-x]; \pi_{temp}),
    \end{equation*}
    for $x \geq \frac{1}{2}$, where $\pi_{temp} =T_{I,1}^{\alpha -1}(T_{I,1}^{\alpha}(\pi_{sc}))$ and $\alpha > 1$. 
    \begin{enumerate}[label = (\roman*)]
        \item The representation $\pi_x$ is of Arthur type if and only if $x = \frac{1}{2}$. 
        \item The representation $\pi_{x}$ is of critical type if and only if $x \in \{\alpha -2, \alpha -1, \alpha , \alpha +1\}$. 
        \item We have $\pi(\EE) = \pi_{\frac{1}{2}}$, where 
        \[\EE= \scalebox{0.8}{\bordermatrix{ 
         &-\frac{1}{2}&\frac{1}{2}& \cdots &\alpha-3 & \alpha -2 & \alpha-1 & \alpha \cr
        &\lhd &\rhd \cr
        &\lhd &\rhd \cr
        &&\odot &\cdots &\odot \cr
        &&&&&&\odot \cr
        &&&&&&&\odot 
        }}.\]
    \end{enumerate}
\end{prop}

The next case we have is $\pi_{temp} = T_{IV,3}(T_{I,1}^{\alpha}(\pi_{sc}))$. 
\begin{prop}\label{nontemp2B3}
    Consider the representation 
    \begin{equation*}
    \pi_x = L(\Delta_{\rho}[-x,-x],\Delta_{\rho}[-x,-x]; \pi_{temp}),
    \end{equation*}
    for $x \geq \frac{1}{2}$, where $\pi_{temp} = T_{IV,3}(T_{I,1}^{\alpha}(\pi_{sc}))$ and $\alpha \in \mathbb{Z}_{>1}$. 
    \begin{enumerate}[label = (\roman*)]
        \item The representation $\pi_x$ is of Arthur type if and only if $x = 1$. 
        \item The representation $\pi_x$ is of critical type if and only if $(x,\alpha) = (1,2)$. 
        \item We have $\pi(\EE) = \pi_{1}$ where 
        \[\EE= \scalebox{0.8}{\bordermatrix{ 
         &-1 &0 & 1& \cdots & \alpha -2 & \alpha-1 & \alpha \cr
        &\lhd &\odot &\rhd \cr
        &\lhd &\odot &\rhd \cr
        &&\odot &\odot&\cdots &\odot \cr
        &&&&&&&\odot 
        }}.\]
    \end{enumerate}
\end{prop}

Now we move onto the case $\pi_{temp} = T_{V,2}^{\pm}(T_{I,1}^{1}(\pi_{sc}))$. 
\begin{prop}\label{nontemp2B4}
    Consider the representation 
    \begin{equation*}
    \pi_x^{\pm} = L(\Delta_{\rho}[-x,-x],\Delta_{\rho}[-x,-x]; \pi_{temp}^{\pm}),
    \end{equation*}
    for $x \geq \frac{1}{2}$, where $\pi_{temp}^{\pm} = T_{V,2}^{\pm}(T_{I,1}^{1}(\pi_{sc}))$ and $\alpha = 1$.
    \begin{enumerate}[label = (\roman*)]
        \item The representation $\pi_x$ is of Arthur type if and only if $x = 1$. 
        \item The representation $\pi_x$ is of critical type if and only if $x \in \{1,2\}$. 
        \item We have $\pi(\EE^{\pm}) = \pi_{1}^{\pm}$, where 
        \[\EE^{+}= \scalebox{0.8}{\bordermatrix{ 
         &-1 &0 & 1  \cr
        &\lhd &\oplus &\rhd \cr
        &\lhd &\oplus &\rhd \cr
        &&&\odot \cr
        }}, \quad \EE^{-}= \scalebox{0.8}{\bordermatrix{ 
         &-1 &0 & 1  \cr
        &\lhd &\ominus &\rhd \cr
        &\lhd &\ominus &\rhd \cr
        &&&\odot \cr
        }}.\]
    \end{enumerate}
\end{prop}

The next case to consider is $\pi_{temp} = T_{I,1}^{1}(T_{V,2}^{\pm}(\pi_{sc}))$.
\begin{prop}\label{nontemp2B5}
    Consider the representation 
    \begin{equation*}
    \pi_x^{\pm} = L(\Delta_{\rho}[-x,-x],\Delta_{\rho}[-x,-x]; \pi_{temp}^{\pm}),
    \end{equation*}
    for $x \geq \frac{1}{2}$, where $\pi_{temp}^{\pm} = T_{I,1}(T_{V,2}^{\pm}(\pi_{sc}))$ and $\alpha = 0$. 
    \begin{enumerate}[label = (\roman*)]
        \item The representation $\pi_x$ is not of Arthur type for any $x$. 
        \item The representation $\pi_x$ is of critical type if and only if $x \in \{1,2\}$. 
    \end{enumerate}
\end{prop}

We move onto the case $\pi_{temp} = T_{I,2}^{\frac{1}{2}}(\pi_{sc})$. 
\begin{prop}\label{nontemp2B6}
     Consider the representation 
    \begin{equation*}
    \pi_x = L(\Delta_{\rho}[-x,-x],\Delta_{\rho}[-x,-x]; \pi_{temp}),
    \end{equation*}
    for $x \geq \frac{1}{2}$, where $\pi_{temp} = T_{I,2}^{\frac{1}{2}}(\pi_{sc})$ and $\alpha = \frac{1}{2}$. 
    \begin{enumerate}[label = (\roman*)]
        \item The representation $\pi_x$ is of Arthur type if and only if $x = \frac{1}{2}$. 
        \item The representation $\pi_{x}$ is of critical type if and only if $x \in \{\frac{1}{2}, \frac{3}{2}\}$. 
        \item We have that $\pi(\EE)= \pi_{\frac{1}{2}}$, where 
        \[\EE= \scalebox{0.8}{\bordermatrix{ 
         &-\frac{1}{2} & \frac{1}{2} \cr
        &\lhd  &\rhd \cr
        &\lhd  &\rhd \cr
        &&\oplus \cr
        &&\oplus 
        }}.\]
    \end{enumerate}
\end{prop}
The proofs of Propositions \ref{nontemp2B2} to \ref{nontemp2B6} are similar to that of Proposition \ref{nontemp2B1}, which we omit. The next case is $\pi_{temp} = T_{II,3}^{\frac{1}{2}}(\pi_{sc})$. 
\begin{prop}\label{nontemp2B7}
    Consider the representation 
    \begin{equation*}
    \pi_x = L(\Delta_{\rho}[-x,-x],\Delta_{\rho}[-x,-x]; \pi_{temp}),
    \end{equation*}
    for $x \geq \frac{1}{2}$, where $\pi_{temp} = T_{II,3}^{\frac{1}{2}}(\pi_{sc})$ and $\alpha \in \frac{1}{2} + \mathbb{Z}_{> 0}$. 
    \begin{enumerate}[label = (\roman*)]
        \item The representation $\pi_x$ is of Arthur type if and only if $x = \frac{1}{2}$. 
        \item The representation $\pi_{x}$ is of critical type if and only if  $(x,\alpha) = (\frac{3}{2}, \frac{3}{2})$. 
        \item We have that $\pi(\EE) = \pi_{\frac{1}{2}}$, where 
        \[\EE= \scalebox{0.8}{\bordermatrix{ 
         &-\frac{1}{2} & \frac{1}{2}&\cdots & \alpha -1 \cr
        &\lhd  &\rhd \cr
        &\lhd  &\rhd \cr
        &&\ominus \cr
        &&\ominus \cr
        &&\ominus & \dots & \odot 
        }}.\]
    \end{enumerate}
\end{prop}
\begin{proof} 
The sufficient direction is easy to verify. By Lemma \ref{absolutevalue}, the representation $\pi_x$ can only be of Arthur type if $x = \frac{1}{2}$ or $x = \frac{3}{2}$. Therefore, to prove the necessary direction, it suffices to show that the case $x = \frac{3}{2}$ is not of Arthur type. 

Suppose it is, then by Definition \ref{psiminus}, the extended multi-segment corresponding to $\pi_{temp}$ must contain at least two segments of the form $([\frac{1}{2}, -\frac{1}{2}]_\rho, *,*)$. If this is true, then by looking at $\Omega(\pi_{temp})$, we see that there can be at most $1$ more segment of the form $([A_i, B_i]_\rho,*,*)$ with $B_i = \frac{1}{2}$. This contradicts \cite[Proposition $3.5$]{HLL22} since $\pi_{temp}$ has an order $2$ nonzero derivative at $x = \frac{1}{2}$. This proves part $(i)$. The rest follows from definition. \end{proof}

Three more cases remain. The next one is $\pi_{temp} = T_{III, 2}^{\frac{1}{2}}(\pi_{sc})$. 
\begin{prop}\label{nontemp2B8}
     Consider the representation 
    \begin{equation*}
    \pi_x = L(\Delta_{\rho}[-x,-x],\Delta_{\rho}[-x,-x]; \pi_{temp}),
    \end{equation*}
    for $x \geq \frac{1}{2}$, where $\pi_{temp}= T_{III, 2}^{\frac{1}{2}}(\pi_{sc})$ and $\alpha = \frac{1}{2}$. 
    \begin{enumerate}
        \item The representation $\pi_x$ is of Arthur type if and only if $x = \frac{1}{2}$.
        \item The representation $\pi_{x}$ is of critical type if and only if $x \in \{\frac{1}{2}, \frac{3}{2}\}$.
        \item We have $\pi(\EE) = \pi_{\frac{1}{2}}$, where 
        \[\EE= \scalebox{0.8}{\bordermatrix{ 
         &-\frac{1}{2} & \frac{1}{2} \cr
        &\lhd  &\rhd \cr
        &\lhd  &\rhd \cr
        &&\ominus \cr
        &&\ominus \cr
        }}.\]
    \end{enumerate}
\end{prop}

The next tempered representation of corank $2$ is $\pi_{temp} = T_{IV,5}(\pi_{sc})$. 
\begin{prop}\label{nontemp2B9}
     Consider the representation 
    \begin{equation*}
    \pi_x = L(\Delta_{\rho}[-x,-x],\Delta_{\rho}[-x,-x]; \pi_{temp}),
    \end{equation*}
    for $x \geq \frac{1}{2}$, where $\pi_{temp} = T_{IV,5}(\pi_{sc})$ and $\alpha \in \mathbb{Z}_{>0}$. 
    \begin{enumerate}[label = (\roman*)]
        \item The representation $\pi_x$ is of Arthur type if and only if $x =1$. 
        \item The representation $\pi_x$ is of critical type if and only if $(x,\alpha) = (1,1)$. 
        \item We have that $\pi(\EE) = \pi_1$, where 
        \[\EE= \scalebox{0.8}{\bordermatrix{ 
         &-1 & 0 & 1&\cdots & \alpha -1 \cr
        &\lhd  & \odot &\rhd \cr
        &\lhd  &\odot &\rhd \cr
        &&\odot \cr
        &&\odot \cr
        &&\odot & \odot &\cdots & \odot
        }}.\]
    \end{enumerate}
\end{prop}

The final case in this section is $\pi_{temp} = T_{V,4}^{\pm}(\pi_{sc})$. 
\begin{prop}\label{nontemp2B10}
    Consider the representation 
    \begin{equation*}
    \pi_x^{\pm} = L(\Delta_{\rho}[-x,-x],\Delta_{\rho}[-x,-x]; \pi_{temp}^{\pm}),
    \end{equation*}
    for $x \geq \frac{1}{2}$, where $\pi_{temp}^{\pm} = T_{V,4}^{\pm}(\pi_{sc})$ and $\alpha = 0$. 
    \begin{enumerate}
        \item The representation $\pi_x$ is of Arthur type if and only if $x = 1$. 
        \item The representation $\pi_x$ is of critical type if and only if $x = 1$.
        \item We have $\pi(\EE^{\pm}) = \pi_1^{\pm}$, where 
        \[\EE^{+}= \scalebox{0.8}{\bordermatrix{ 
         &-1 &0 &1  \cr
        &\lhd &\oplus &\rhd \cr
        &\lhd &\oplus &\rhd \cr
        &&\oplus \cr
        &&\oplus \cr
        }}, \quad \EE^{-}= \scalebox{0.8}{\bordermatrix{ 
         &-1 &0 &1  \cr
        &\lhd &\ominus &\rhd \cr
        &\lhd &\ominus &\rhd \cr
        &&\ominus \cr
         &&\ominus \cr
        }}.\]
    \end{enumerate}
\end{prop}
The proofs of Propositions \ref{nontemp2B8} to \ref{nontemp2B10} is similar to that of Proposition \ref{nontemp2B7}, which we omit. 
This concludes our work in Case $(B)$. Now we move onto Case $(C)$.

\subsection{\texorpdfstring{Case $(C)$: $\pi = L(\Delta_{\rho}[-x_1, -x_1-1], \Delta_{\rho}[-x_2, -x_2]; \pi_{temp})$}{}}
In this subsection we are concerned with representations of the form $\pi_{x_1, x_2} = L(\Delta_{\rho}[-x_1, -x_1-1], \Delta_{\rho}[-x_2, -x_2]; \pi_{temp})$, where $\pi_{temp}$ is a corank 1 tempered of good parity. Here we require $x_1 \geq 0, x_2 \geq \frac{1}{2}$ and $x_2 - x_1 \leq \frac{1}{2}$ by Langlands classification of classical groups. We assume our representations are of good parity, so the last condition can be rewritten as $x_2 \leq x_1$.

Since we're dealing with tempered representations of corank $1$, there are three cases to consider. Let us begin with $\pi_{temp} = T_{I,1}^{\alpha}(\pi_{sc})$. 
\begin{prop}\label{nontemp2C1}
    Consider the representation 
    \begin{equation*}
        \pi_{x_1, x_2} = L(\Delta_{\rho}[-x_1, -x_1-1], \Delta_{\rho}[-x_2, -x_2]; \pi_{temp}),
    \end{equation*}
    where $x_1 \geq 0, x_2 \geq \frac{1}{2}$,  $x_2 \leq x_1$, and $\pi_{temp} = T_{I,1}^{\alpha}(\pi_{sc})$ for $\alpha > 0$. 
    \begin{enumerate}[label = (\roman*)]
        \item The representation $\pi_{x_1, x_2}$ is of Arthur type if and only if $\frac{1}{2} \leq x_2 \leq \alpha -3$, $\epsilon_\rho +1 \leq x_1 \leq \alpha -2$ and $x_1 > x_2$ for $\alpha \geq \frac{7}{2}$. 
        \item $\pi_{x_1, x_2}$ is of critical type when $(x_1, x_2) = (\alpha-2, \alpha-3),(\alpha -2, \alpha -2),(\alpha-1, \alpha -1), (\alpha -1, \alpha -2),(\alpha, \alpha), (\alpha, \alpha -1), (\alpha +1, \alpha +1), (\alpha +1, \alpha), (\alpha+1, \alpha-1)$ or $(\alpha +2, \alpha +1)$.  
        \item Define $\EE_{x_1, x_2}$ in various cases as follows, then we have $\pi(\EE_{x_1, x_2}) = \pi_{x_1, x_2}$.  When  $\frac{1}{2} \leq x_2 \leq \alpha -3$, $\epsilon_\rho +1 \leq x_1 \leq \alpha -2$ and $x_1 = x_2 + 1$, define
        \begin{flalign*}
            \EE_{x_1, x_2}:= &\{([x_1+1,-x_1]_\rho,\lfloor x_1 \rfloor, (-1)^{x_1 +1} \eta),([x_1-3, \epsilon_\rho]_\rho,0,\eta), ([x_1 -1, x_1 -1]_\rho,0,(-1)^{x_1 - \epsilon_\rho}\eta),\\
            &([\alpha -2, x_1 +1]_\rho, 0, (-1)^{x_1 +1- \epsilon_\rho}\eta),([\alpha, \alpha]_\rho,0,(-1)^{x_2 +1}\eta)\}.
        \end{flalign*}
        Here is the associated symbol. 
        \[\EE_{x_1, x_2}= \scalebox{0.8}{\bordermatrix{ 
        &-x_1  & \cdots & \epsilon_\rho & \cdots & x_1-3 & x_1-2 & x_1-1 &x_1 & x_1+1 & \cdots & \alpha -2 &\alpha -1 & \alpha \cr
        &\lhd & \cdots & \odot & \cdots & \cdots &\cdots & \cdots &\cdots  &\rhd \cr
        &&&\odot &\cdots & \odot \cr
        &&&&&&&\odot \cr
        &&&&&&&&&\odot & \cdots & \odot \cr
        &&&&&&&&&&&&&\odot
        }}.\]
        When $\frac{1}{2} \leq x_2 \leq \alpha -3$, $\epsilon_\rho +1 \leq x_1 \leq \alpha -2$ and $x_1 - x_1 > 1$, define 
        \begin{flalign*}
            \EE_{x_1, x_2}:= &\{([x_1+1,-x_1]_\rho,\lfloor x_1 \rfloor, (-1)^{x_1 +1} \eta),([x_2, -x_2]_\rho,\lfloor x_2 \rfloor, (-1)^{x_2 - \epsilon_\rho} \eta),\\
            &([x_2-2, 0]_\rho,0,-\eta), ([x_1 -2, x_2]_\rho,0,(-1)^{x_2}\eta),\\
            &([\alpha -2, x_1 +1]_\rho, 0, (-1)^{x_1 +1}\eta),([\alpha, \alpha]_\rho,0,(-1)^{x_2 +1}\eta)\}.
        \end{flalign*}
        Here is the associated symbol for $\EE_{x_1, x_2}:$
        \[ \scalebox{0.8}{\bordermatrix{ 
        &-x_1  & \cdots &-x_2&\cdots& \epsilon_\rho & \cdots & x_2 -2& x_2 -1 & x_2 & \cdots & x_1 -2& x_1-1 &x_1 & x_1+1 & \cdots & \alpha -2 &\alpha -1 & \alpha \cr
        &\lhd & \cdots & \lhd&\cdots & \odot & \cdots & \cdots &\cdots  & \rhd &\cdots &\rhd &\rhd &\rhd&\rhd\cr
        &&&\lhd & \cdots & \odot & \cdots & \cdots & \cdots& \rhd&\cr
        &&&&&\odot &\cdots & \odot \cr
        &&&&&&&\odot \cr
        &&&&&&&&&\odot & \cdots & \odot \cr
        &&&&&&&&&&&&&&\odot & \cdots & \odot\cr
        &&&&&&&&&&&&&&&&&&\odot
        }}.\]
        
    \end{enumerate}
\end{prop}
\begin{proof} The sufficient condition can be proven in a similar way as Proposition \ref{nontemp2B1}, which we omit. Let us now show the necessary direction. First by definition we must have: 
\begin{equation*}
    \pi_{x_1, x_2}^{\rho,-} = L(\Delta_{\rho}[-x_2,-x_2]; \pi_{temp}),
\end{equation*} 
which is of Arthur type if and only if $\frac{1}{2} \leq x_2 \leq \alpha - 1$ for $\alpha \geq \frac{3}{2}$, or $x_2 = \alpha = \frac{1}{2}$, which is not possible in this case. Furthermore, since 
\begin{equation*}
    |\Omega|(\pi_{x_1, x_2}^{\rho,-}) \subseteq \{\rho\lvert \cdot \rvert^{y}: \epsilon_\rho \leq y \leq \alpha -2\}\cup \{\rho\lvert \cdot \rvert^{\alpha}\} \cup \{\rho \lvert \cdot \rvert^{x_2}, \rho\lvert \cdot \rvert^{-x_2}\},
\end{equation*}
we must have either
\begin{itemize}
    \item $\epsilon_\rho+1 \leq x_1 \leq \alpha -2$, or 
    \item $x_1 = x_2$ and $\epsilon_\rho \leq x_2 -1 \leq \alpha -2$
\end{itemize}
It suffices to show that the case $x_1 = x_2$ does not give a representation of Arthur type. Suppose the contrary, then there exists some segment $([x_1, -x_1+1]_\rho,*,*)$ inside the extended multi-segment $\EE_{x_1, x_2}$ corresponding to $\pi_{x_1, x_2}$ by Definition \ref{psiminus}. This is impossible since the first segment inside $\EE_{x_1, x_2}$ is of the form $([x_1, -x_1], \lfloor x_1 \rfloor, *)$, and applying the raising operators does not produce the desired segment.  The other segments contain a gap at $x_1 -1$, so producing a segment of the form $([x_1, -x_1+1]_\rho,*,*)$ is impossible. This proves the sufficient direction. Considering our restrictions on $x_1, x_2$ gives the desired conditions. This proves part $(i)$. Part $(ii)$ follows from definition.  \end{proof}

The next case is $\pi_{temp} = T_{IV, 3}(\pi_{sc})$. 
\begin{prop}\label{nontemp2C2}
    Consider the representation 
    \begin{equation*}
        \pi_{x_1, x_2} = L(\Delta_{\rho}[-x_1, -x_1-1], \Delta_{\rho}[-x_2, -x_2]; \pi_{temp}),
    \end{equation*}
    where $x_1 \geq 0, x_2\geq \frac{1}{2}$,  $x_2 \leq x_1$, and $\pi_{temp} = T_{IV, 3}(\pi_{sc})$ for $\alpha \in \mathbb{Z}_{>0}$. 
    \begin{enumerate}[label = (\roman*)]
        \item The representation $\pi_{x_1, x_2}$ is of Arthur type if and only if $1 < x_1 \leq \alpha -1$ and $1 \leq x_2 \leq \alpha -2$ with $x_2 < x_1$, or when $x_1 = x_2 = 1$. 
        \item The representation $\pi_{x_1, x_2}$ is of critical type if and only if $(x_1, x_2, \alpha) = (2,1,3), (2,1,2), (2,1,1)$.
        \item Define $\pi_{x_1, x_2}$ in various cases as follows, then we have $\pi(\EE_{x_1, x_2}) = \pi_{x_1, x_2}$.
        When  $1 < x_1 \leq \alpha -1$, $1 \leq x_2 \leq \alpha -2$ and $x_1 = x_2 + 1$, define
        \begin{flalign*}
            \EE_{x_1, x_2}:= &\{([x_1+1,-x_1]_\rho,\lfloor x_1 \rfloor, (-1)^{x_1 +1} \eta), ([0,0]_\rho,0,\eta)^2, ([x_1-3, 0]_\rho,0,\eta), \\
            &([x_1 -1, x_1 -1]_\rho,0,(-1)^{x_1 }\eta),([\alpha -1, x_1 +1]_\rho, 0, (-1)^{x_1 +1}\eta).
        \end{flalign*}
        Here is the associated symbol. 
        \[\EE_{x_1, x_2}= \scalebox{0.8}{\bordermatrix{ 
        &-x_1  & \cdots & 0 & \cdots & x_1-3 & x_1-2 & x_1-1 &x_1 & x_1+1 & \cdots & \alpha -2 &\alpha -1 \cr
        &\lhd & \cdots & \odot & \cdots &\cdots & \rhd & \rhd &\rhd &\rhd \cr
        &&&\odot \cr
        &&&\odot \cr
        &&&\odot &\cdots & \odot  \cr
        &&&&&&&\odot \cr
        &&&&&&&&&\odot & \cdots & \odot &\odot 
        }}.\]
        When $\frac{1}{2} \leq x_2 \leq \alpha -3$, $\epsilon_\rho +1 \leq x_1 \leq \alpha -2$ and $x_1 - x_1 > 1$, define 
        \begin{flalign*}
            \EE_{x_1, x_2}:= &\{([x_1+1,-x_1]_\rho,\lfloor x_1 \rfloor, (-1)^{x_1 +1} \eta),([x_2, -x_2]_\rho,\lfloor x_2 \rfloor, (-1)^{x_2 } \eta),\\
            &([0,0]_\rho,0,-\eta)^2,([x_2-2, 0]_\rho,0,\eta), ([x_1 -2, x_2]_\rho,0,(-1)^{x_2}\eta),\\
            &([\alpha -1, x_1 +1]_\rho, 0, (-1)^{x_1 +1}\eta)\}.
        \end{flalign*}
        Here is the associated symbol for $\EE_{x_1, x_2}$:
        \[\scalebox{0.8}{\bordermatrix{ 
        &-x_1  & \cdots &-x_2&\cdots& 0 & \cdots & x_2 -2& x_2 -1 & x_2 & \cdots & x_1 -2& x_1-1 &x_1 & x_1+1 & \cdots & \alpha -2 &\alpha -1 \cr
        &\lhd & \cdots & \lhd&\cdots & \odot & \cdots & \cdots &\cdots  & \rhd &\cdots &\rhd &\rhd &\rhd&\rhd\cr
        &&&\lhd & \cdots & \odot & \cdots & \cdots & \cdots& \rhd&\cr
        &&&&&\odot &\cdots & \odot \cr
        &&&&&&&\odot \cr
        &&&&&&&&&\odot & \cdots & \odot \cr
        &&&&&&&&&&&&&&\odot & \cdots & \odot&\odot \cr
        }}.\]
        When $x_1 = x_2 = 1$, we have that $\pi_{1,1} = \pi(\EE_{1,1})$, where 
        \[\EE_{1,1}= \scalebox{0.8}{\bordermatrix{ 
        & -1 & 0 & 1 & 2 & \ldots & \alpha -1& \alpha \cr
        &\lhd & \odot & \rhd \cr
        & \lhd & \lhd & \odot&\odot & \ldots & \rhd & \rhd \cr
        &&\odot \cr
        &&&&\odot & \ldots &\odot & \odot 
        }}.\]
    \end{enumerate}
\end{prop}
\begin{proof} From \cite[Proposition $10.13$]{HJLLZ24}, we see that the representation 
\begin{equation*}
    \pi_{x_1, x_2}^{\rho,-} = L(\Delta_{\rho}[-x_2, -x_2]; \pi_{temp})
\end{equation*}
is of Arthur type if and only if $1 \leq x_2 \leq \alpha$. The rest of the proof is similar to that of Proposition \ref{nontemp2C1}, which we omit. \end{proof}

The last tempered representation of corank $1$ is $\pi_{temp} =T_{V,2}^{\pm}(\pi_{sc})$. 
\begin{prop}\label{nontemp2C3}
     Consider the representation 
    \begin{equation*}
        \pi_{x_1, x_2}^{\pm} = L(\Delta_{\rho}[-x_1, -x_1-1], \Delta_{\rho}[-x_2, -x_2]; \pi_{temp}^{\pm}),
    \end{equation*}
    where $x_1 \geq 0, x_2\geq \frac{1}{2}$,  $x_2 \leq x_1$, and $\pi_{temp}^{\pm} = T_{V,2}^{\pm}(\pi_{sc})$ for $\alpha = 0$. 
    \begin{enumerate}[label = (\roman*)]
        \item The representation $\pi_{x_1, x_2}^{\pm}(\pi_{sc})$ is not of Arthur type. 
        \item The representation $\pi_{x_1, x_2}^{\pm}$ is of critical type if and only if  $(x_1, x_2)= (2,1)$ or $(1,1)$. 
    \end{enumerate}
\end{prop}
\begin{proof} By \cite[Proposition $10.14$]{HJLLZ24}, the representation 
\begin{equation*}
    (\pi_{x_1, x_2}^{\pm})^{\rho,-} = L(\Delta_{\rho}[-x_2, -x_2]; \pi_{temp}^{\pm})
\end{equation*}
is of Arthur type if and only if $x_2 =1$. In this case, $\Psi((\pi_{x_1, x_2}^{\pm})^{\rho,-})$ are singletons with
\begin{equation*}
    |\Omega|(L(\Delta_{\rho}[-x_2, -x_2]; \pi_{temp}^{\pm})) \subseteq \{\rho \lvert \cdot \rvert^{y}, -1 \leq y \leq 1\} \cup \{\rho \lvert \cdot \rvert^{x_2}, \rho\lvert \cdot \rvert^{-x_2}\}.
\end{equation*}
From Lemma \ref{absolutevalue}, the only possible $x_1$ for $\pi_{x_1, x_2}^{\pm}$ to be of Arthur type is $x_1 = x_2 = 1$, but this case can also be eliminated by Definition \ref{psiminus} and Theorem \ref{nontemp red}. This proves part $(i)$. Part $(ii)$ follows from definition. \end{proof}. 

This concludes our discussion of Case $(C)$. Now let's move onto Case $(D)$. 
\subsection{\texorpdfstring{Case $(D): \pi = L(\Delta_{\rho}[-x_1, -x_1], \Delta_{\rho}[-x_2, -x_2 -1]; \pi_{temp})$}{}}
In this subsection we consider representations of the form $\pi_{x_1, x_2} = L(\Delta_{\rho}[-x_1, -x_1], \Delta_{\rho}[-x_2, -x_2 -1]; \pi_{temp})$, where $\pi_{temp}$ is tempered of corank $1$. Here the situation is similar to that of Case $(C)$ except the orders of the segments are exchanged. Just like in Case $(C)$, we have $x_1 > 0, x_2 > -\frac{1}{2}$ and $x_1 - x_2 \geq \frac{1}{2}$ by the Langlands classification of classical groups. By the good parity condition, the last statement can be translated to $x_1 > x_2$.

We consider the same three cases as before, starting with $\pi_{temp} = T_{I,1}^{\alpha}(\pi_{sc})$. 
\begin{prop}\label{nontemp2D1}
    Consider the representation 
    \begin{equation*}
        \pi_{x_1, x_2} = L(\Delta_{\rho}[-x_1,-x_1], \Delta_\rho[-x_2, -x_2-1]; \pi_{temp}),
    \end{equation*}
    where $x_1 > x_2 \geq 0$ and $\pi_{temp} = T_{I,1}^{\alpha}(\pi_{sc})$ for $\alpha > 0$. 
    \begin{enumerate}[label = (\roman*)]
        \item The representation $\pi_{x_1, x_2}$ is of Arthur type if and only if $\epsilon_\rho +1 \leq x_1 \leq \alpha$, $0 \leq x_2 \leq \alpha -2$ and $x_1 - x_2 \geq 2$ for $\alpha \geq 3$, or $(x_1, x_2) = (1,0)$ for $\alpha \in \mathbb{Z}_{>1}$, or $(x_1, x_2, \alpha) = (\frac{5}{2}, \frac{1}{2}, \frac{1}{2})$ and $\epsilon_{sc}(\rho \otimes S_2) = -1$.  
        \item The representation $\pi_{x_1, x_2}$ is of critical type if and only if $(x_1, x_2) = (\alpha -1, \alpha-3), (\alpha -1, \alpha -2), (\alpha, \alpha -2), (\alpha, \alpha -1), (\alpha +1, \alpha-1), (\alpha +1, \alpha), (\alpha +2, \alpha), (\alpha +2, \alpha +1), (\alpha +3, \alpha +1)$. 
        \item Define $\mathcal{E}_{x_1,x_2}$ in various cases as follows, then $\pi(\EE_{x_1, x_2}) = \pi_{x_1, x_2}$. When $\epsilon_\rho + 1 \leq x_1 \leq \alpha -1$, $0 \leq x_2 \leq \alpha -2$ and $x_1 - x_2 \geq 2$ for $\alpha \geq 3$, define
        \begin{flalign*}
            \EE_{x_1, x_2} := &\{([x_1,-x_1]_\rho,\lfloor x_1 \rfloor, \eta), ([x_2+1, -x_2]_\rho, \lfloor x_2 \rfloor, -\eta), 
            ([x_2-2,\epsilon_\rho]_\rho,0,\eta), \\
            &([\alpha-2, x_2+1]_\rho,0,(-1)^{x_2 +1 - \epsilon_\rho}\eta), ([\alpha,\alpha]_\rho,0,(-1)^{\alpha - \epsilon_\rho}\eta)\}.
        \end{flalign*}
        
        Here is the associated symbol. 
        \[\EE_{x_1, x_2}= \scalebox{0.8}{\bordermatrix{ 
        &-x_1  & \cdots & -x_2 & \cdots &\epsilon_\rho & \cdots & x_2 -2 & x_2 -1 & x_2 & x_2 +1& \cdots & x_1 & \cdots &\alpha -2 & \alpha -1 & \alpha \cr
        &\lhd & \cdots & \lhd &\cdots & \odot & \cdots & \rhd &\rhd & \rhd &\rhd  &\cdots &\rhd \cr
        &&&\lhd &\cdots & \odot & \cdots & \rhd & \rhd & \rhd&\rhd \cr
        &&&&&\odot & \cdots & \odot\cr
        &&&&&&& && & \odot & \cdots & \cdots &\cdots &\odot \cr
        &&&&&&&&&&&&&&&&\odot 
        }}.\]
        When $(x_1, x_2) = (1,0)$ and $\alpha \in \mathbb{Z}_{>1}$, define 
        \begin{flalign*}
            \EE_{x_1, x_2} := &\{([1,-1]_\rho,1, \eta), ([1,0]_\rho, 1, -\eta), 
            ([\alpha -2, 1]_\rho,0,\eta), \\
            &([\alpha, \alpha]_\rho,0,(-1)^{\alpha -1}\eta)\}.
        \end{flalign*}
        Here is the associated symbol. 
        \[\EE_{x_1, x_2}= \scalebox{0.8}{\bordermatrix{ 
        &-1 & 0 & 1 & \cdots & \alpha -2 & \alpha -1 & \alpha \cr
        &\lhd & \odot & \rhd \cr
        &&\lhd & \rhd \cr
        &&&\odot & \cdots & \odot \cr
        &&&&&&&\odot 
        }}.\]
        When $(x_1, x_2, \alpha) = (\frac{5}{2}, \frac{1}{2}, \frac{1}{2})$ and $\epsilon_{sc}(\rho \otimes S_2)= -1$, let 
        \[\EE_{x_1, x_2}= \scalebox{0.8}{\bordermatrix{ 
        &-\frac{5}{2} & -\frac{3}{2} & -\frac{1}{2} & \frac{1}{2} & \frac{3}{2} & \frac{5}{2} \cr
        &\lhd & \lhd & \oplus & \ominus &\rhd & \rhd \cr
        &&&\lhd & \rhd 
        }}.\]
    \end{enumerate}
\end{prop}
\begin{proof}The sufficient direction can be proven in the same way as Proposition \ref{nontemp2A1}, which we omit. Now we show the necessary direction. By Proposition \cite[Proposition $11.11$]{HJLLZ24}, the representation 
\begin{equation*}
    \pi_{x_1, x_2}^{\rho,-} = L(\Delta_\rho[-x_2, -x_2-1]; \pi_{temp})
\end{equation*}
is of Arthur type if and only if $0 \leq x_2 \leq \alpha-2$ for $\alpha \geq 2$, or when $(x_2, \alpha) = (0,1)$ or $(x_2, \alpha) = (\frac{1}{2}, \frac{1}{2})$. Matching this with our restrictions gives the condition on $x_2$. For the restrictions on $x_1$, we see from Lemma \ref{absolutevalue} that either $\epsilon_\rho +1 \leq x_1 \leq \alpha -1$, or $x_1 = \alpha$. Therefore, it remains to show that the case $x_2 = x_1 -1$ is not of Arthur type, when $x_1 > 1, x_2 \geq 0$ satisfies the restrictions above. 

Suppose $\pi_{x_1, x_2}$ is of Arthur type in this case. By Definition \ref{psiminus}, there exists an $\EE$ with $\pi(\EE) = L(\Delta_\rho[-x_2, -x_2-1]; \pi_{temp})$ that contains a segment of the form $([x_1-1, -(x_1-1)]_\rho, *, *)$. However, any such $\EE$ must have the first segment (by the admissible order) be of the form $([x_1, -(x_1-1)]_\rho,*,*)$. This strictly contains our desired segment, and we cannot reduce the segment any further by definition, which gives a contradiction. This proves the necessary direction and part $(i)$. Parts $(ii)$ and $(iii)$ follow from definition. \end{proof}

The next case is $\pi_{temp} = T_{IV, 3}(\pi_{sc})$. 
\begin{prop}\label{nontemp2D2}
    Consider the representation 
    \begin{equation*}
        \pi_{x_1, x_2} = L(\Delta_{\rho}[-x_1,-x_1], \Delta_\rho[-x_2, -x_2-1]; \pi_{temp}),
    \end{equation*}
    where $x_1 > x_2 \geq 0$ and $\pi_{temp} = T_{IV, 3}(\pi_{sc})$ for $\alpha \in \mathbb{Z}_{> 0}$. 
    \begin{enumerate}[label = (\roman*)]
        \item The representation $\pi_{x_1, x_2}$ is of Arthur type if and only if $2 \leq x_1 \leq \alpha$, $0 \leq x_2 \leq \alpha -2$, and $x_1 - x_2 \geq 2$, or when $(x_1, x_2) = (1,0)$ or $(\alpha+1, \alpha-1)$.
        \item The representation $\pi_{x_1, x_2}$ is of critical type if and only if 
        \begin{equation*}
            (x_1, x_2, \alpha) = (1,0,1),(2,1,1),(3,1,1),(2,1,2),(3,1,2),(3,1,3). 
        \end{equation*} 
        \item Define $\mathcal{E}_{x_1, x_2}$ in various cases as follows, then $\pi(\EE_{x_1, x_2}) = \pi_{x_1, x_2}$. When $2 \leq x_1 \leq \alpha, 0 \leq x_2 \leq \alpha -2$ and $x_1 - x_2 \geq 2$, define
        \begin{flalign*}
            \EE_{x_1, x_2} := &\{([x_1,-x_1]_\rho,\lfloor x_1 \rfloor, \eta), ([x_2, -x_2]_\rho, \lfloor x_2 \rfloor, -\eta), ([0,0]_\rho,0,\eta)^2 \\
            &([x_2-2,0]_\rho,0,\eta), 
            ([\alpha-1, x_2+1]_\rho,0,(-1)^{x_2 +1}\eta)\}.
        \end{flalign*}
         Here is the associated symbol. 
        \[\EE_{x_1, x_2}= \scalebox{0.8}{\bordermatrix{ 
        &-x_1  & \cdots & -x_2 & \cdots &0 & \cdots & x_2 -2 & x_2 -1 & x_2 & x_2 +1& \cdots & x_1 & \cdots & \alpha -1 &  \cr
        &\lhd & \cdots & \lhd &\cdots & \odot & \cdots & \rhd &\rhd & \rhd &\rhd  &\cdots &\rhd \cr
        &&&\lhd &\cdots & \odot & \cdots & \rhd & \rhd & \rhd \cr
        &&&&&\odot \cr
        &&&&&\odot \cr
        &&&&&\odot & \cdots & \odot \cr
        &&&&&&& &  && \odot & \cdots & \cdots &\cdots  &\odot \cr
        }}.\]
        When $(x_1, x_2) = (\alpha+1, \alpha-1)$, define
        \begin{flalign*}
            \EE_{x_1, x_2} := &\{([x_1,-x_1]_\rho,\lfloor x_1 \rfloor, \eta), ([x_2, -x_2]_\rho, \lfloor x_2 \rfloor, -\eta), ([0,0]_\rho,0,\eta)^2 \\
            &([x_2-2,0]_\rho,0,\eta)\}.
        \end{flalign*}
         Here is the associated symbol. 
        \[\EE_{x_1, x_2}= \scalebox{0.8}{\bordermatrix{ 
        &-x_1  & -\alpha & -x_2 & \cdots &0 & \cdots & x_2 -2 & x_2 -1 & x_2 & \alpha & x_1 &   \cr
        &\lhd & \lhd & \lhd &\cdots & \odot & \cdots & \rhd &\rhd  &\rhd  &\rhd &\rhd \cr
        &&&\lhd &\cdots & \odot & \cdots & \rhd & \rhd & \rhd \cr
        &&&&&\odot \cr
        &&&&&\odot \cr
        &&&&&\odot & \cdots & \odot \cr
        }}.\]
        When $(x_1, x_2) = (1,0)$, define 
        \begin{flalign*}
            \EE_{x_1, x_2} := &\{([1,-1]_\rho,1, \eta), ([0,0]_\rho,0,\eta)^2, ([\alpha-1, 0]_\rho,1,-\eta),([1,1]_\rho,0,\eta)\}.
        \end{flalign*}
        Here is the associated symbol. 
        \[\EE_{x_1, x_2}= \scalebox{0.8}{\bordermatrix{ 
        &-1 & 0 & 1 & \cdots  & \alpha -1 & \alpha \cr
        &\lhd & \odot & \rhd \cr
        &&\odot \cr
        &&\odot \cr
        &&\lhd &\odot & \cdots &\rhd \cr
        &&&\odot 
        }}.\]
    \end{enumerate}
\end{prop}
\begin{proof} From \cite[Proposition $11.12$]{HJLLZ24}, we see that the representation 
\begin{equation*}
    \pi_{x_1, x_2}^{\rho,-} = L(\Delta_\rho[-x_2, -x_2-1]; \pi_{temp})
\end{equation*}
is of Arthur type if and only if $0 \leq x_2 \leq \alpha-1$. From here, the rest of the proof follows the exact same way as Proposition \ref{nontemp2D1}, which we omit. \end{proof}

The last tempered representation of corank $1$ is $\pi_{temp}= T_{V,2}^{\pm}(\pi_{sc})$. 
\begin{prop}\label{nontemp2D3}
    Consider the representation 
    \begin{equation*}
        \pi_{x_1, x_2}^{\pm} = L(\Delta_{\rho}[-x_1,-x_1], \Delta_\rho[-x_2, -x_2-1]; \pi_{temp}^{\pm}),
    \end{equation*}
    where $x_1 > x_2 \geq 0$ and $\pi_{temp}^{\pm} =  T_{V,2}^{\pm}(\pi_{sc})$ for $\alpha = 0$. 
    \begin{enumerate}[label = (\roman*)]
        \item The representation $\pi_{x_1, x_2}^{\pm}$ is of Arthur type if and only if $(x_1, x_2) = (1,0)$.
        \item The representation $\pi_{x_1, x_2}^{\pm}$ is of critical type if and only if $(x_1, x_2) = (1,0),(2,0)$. 
        \item We have $\pi(\EE^{\pm}) = \pi_{0,1}^{\pm}$, where 
        \[\EE^{+}= \scalebox{0.8}{\bordermatrix{ 
         &-1 &0 &1  \cr
        &\lhd &\oplus &\rhd \cr
        & &\oplus \cr
        &&\oplus &\ominus\cr
        }}, \quad \EE^{-}= \scalebox{0.8}{\bordermatrix{ 
         &-1 &0 &1  \cr
        &\lhd &\ominus &\rhd \cr
        & &\ominus \cr
        &&\ominus &\oplus\cr
        }}.\]
    \end{enumerate}
\end{prop}
\begin{proof} This follows immediately from \cite[Proposition $11.13$]{HJLLZ24}, which states that the representation 
\begin{equation*}
    (\pi_{x_1, x_2}^{\pm})^{\rho,-} = L(\Delta_\rho[-x_2, -x_2-1]; \pi_{temp}^{\pm})
\end{equation*}
is of Arthur type if and only if $x_2 = 0$. \end{proof}

This concludes our discussion in Case $(D)$. In the next subsection, we consider Cases $(E), (F), (G)$ together. 
\subsection{\texorpdfstring{Cases $(E),(F),(G)$}{} involving supercuspidal representations}
In this subsection, we look at the remaining three cases together, since they all involve supercuspidal representations. The first case we'll examine is Case $(E)$, where $\pi_{x_1, x_2} = L(\Delta_{\rho}[-x_1, -x_1-1], \Delta_{\rho}[-x_2, -x_2 -1]; \pi_{sc})$. 

From Langlands classification we obtain the bounds $x_1, x_2> -\frac{1}{2}$ and $x_1 \geq x_2$. 
\begin{prop}\label{nontemp2E}
    Consider the representation 
    \begin{equation*}
        \pi_{x_1, x_2} = L(\Delta_{\rho}[-x_1, -x_1-1], \Delta_{\rho}[-x_2, -x_2 -1]; \pi_{sc}),
    \end{equation*}
    for $x_1, x_2> -\frac{1}{2}$ and $x_1 \geq x_2$.
    \begin{enumerate}
        \item The representation $\pi_{x_1, x_2}$ is of Arthur type if and only if $1 \leq x_1 \leq \alpha -1$ and $x_1 > x_2$, or when $(x_1, x_2) = (0,0)$ or $(\alpha, \alpha-1)$.  
        \item The representation $\pi_{x_1, x_2}$ is of critical type if and only if  $(x_1, x_2) = (\alpha-1, \alpha -3),(\alpha-1, \alpha -2),(\alpha-1, \alpha -1),(\alpha, \alpha -2),(\alpha, \alpha -1),(\alpha, \alpha),(\alpha +1, \alpha -1),(\alpha +1, \alpha),(\alpha +2,\alpha)$. 
        \item When  $1 \leq x_1 \leq \alpha -1$ and $x_1 > x_2$, define 
         \begin{flalign*}
            \EE_{x_1, x_2} := &\{([x_1,-x_1]_\rho,\lfloor x_1 \rfloor, \eta), 
            ([x_2-2,\epsilon_\rho]_\rho,0,\eta), 
            ([\alpha-1, x_2+1]_\rho,0,(-1)^{x_2 +1}\eta).
        \end{flalign*}
        Then $\pi(\EE_{x_1, x_2}) = \pi_{x_1, x_2}$. Here is the associated symbol. 
        \[\EE_{x_1, x_2}= \scalebox{0.8}{\bordermatrix{ 
        &-x_1  & \cdots & -x_2 & \cdots &\epsilon_\rho & \cdots & x_2 -2 & x_2 -1 & x_2 & x_2 +1& \cdots & x_1 & \cdots & \alpha -1 &  \cr
        &\lhd & \cdots & \lhd &\cdots & \odot & \cdots & \rhd &\rhd & \rhd &\rhd  &\cdots &\rhd \cr
        &&&&&\odot & \cdots & \odot \cr
        &&&&&&& &  && \odot & \cdots & \cdots &\cdots  &\odot \cr
        }}.\]
    \end{enumerate}
\end{prop}
\begin{proof} The sufficient direction can be proved in the same way as Proposition \ref{nontemp2A1}, which we omit. For the necessary direction, first note that from Proposition \ref{scArthur}, we have that the representation
\begin{equation*}
    \pi_{x_1, x_2}^{\rho,-} = L(\Delta_{\rho}[-x_2, -x_2-1]; \pi_{sc})
\end{equation*}
is of Arthur type if and only if $x_2 \leq \alpha -1$. Furthermore, Lemma \ref{absolutevalue} gives us the constraint that either $(x_1, x_2) = (\alpha, \alpha-1), (0,0)$ or $x_1 \leq \alpha -1$ since we have $|\Omega|(\pi_{sc}) = \{\rho \lvert \cdot \rvert^{y}: \epsilon_\rho \leq y \leq \alpha-1\}$. Therefore, it suffices to show that the case $x_1 = x_2$ does not give a representation of Arthur type. 

Suppose on the contrary that $\pi_{x_1, x_2}$ is of Arthur type when $x_1 = x_2$. Then by Theorem \ref{nontemp red}, there exists an extended multi-segment $\EE$ containing $2$ copies of the segment
$([x_1, -x_1 +1],*,*)$ such that $\pi(\EE) = \pi_{sc}$. However, we know that $|\Omega|(\pi_{sc})$ is multiplicity-free, which gives a contradiction. This proves part $(i)$. The rest follows from definition. \end{proof}

Now we move onto Case $(F)$, where $\pi_{x_1, x_2} = L(\Delta_{\rho}[-x_1, -x_1-2], \Delta_{\rho}[-x_2, -x_2]; \pi_{sc})$, where the Langlands classification gives us the natural constraint $x_1 > -1, x_2 \geq \frac{1}{2}$ and $x_2 - x_1 \leq 1$. 
\begin{prop}\label{nontemp2F}
    Consider the representation 
    \begin{equation*}
        \pi_{x_1, x_2} = L(\Delta_{\rho}[-x_1, -x_1-2], \Delta_{\rho}[-x_2, -x_2]; \pi_{sc}),
    \end{equation*}
    where $x_1 > -1, x_2 \geq \frac{1}{2}$ and $x_2 - x_1 \leq 1$. 
    \begin{enumerate}[label = (\roman*)]
        \item The representation $\pi_{x_1, x_2}$ of Arthur type if and only if $\epsilon_\rho+1 \leq x_1 \leq \alpha -2$ and $x_2 < x_1$. 
        \item The representation $\pi_{x_1, x_2}$ is of critical type if and only if $(x_1, x_2) = (\alpha -2, \alpha -3),(\alpha -2, \alpha-2),(\alpha-2,\alpha-1), (\alpha -1, \alpha -2),(\alpha -1, \alpha -1),(\alpha-1,\alpha),(\alpha, \alpha -1),(\alpha, \alpha),(\alpha,\alpha+1),(\alpha +1, \alpha)$. 
        \item Define $\EE_{x_1, x_2}$ in various cases as follows, then $\pi(\EE_{x_1, x_2}) = \pi_{x_1, x_2}$. When $\epsilon_\rho +1 \leq x_1 \leq \alpha -2$ and $x_2 = x_1 -1$, define
        \begin{flalign*}
            \EE_{x_1, x_2} := &\{([x_1+2,-x_1]_\rho,\lfloor x_1 \rfloor, \eta),
            ([x_2-2,\epsilon_\rho])_\rho,0,\eta), 
            ([x_1 -1,x_2]_\rho,0,(-1)^{x_2+1 -\epsilon_\rho}\eta), \\
            &([\alpha -1, x_1 +2]_\rho,0,(-1)^{x_1+1- \epsilon_\rho}\eta).
        \end{flalign*}
        Here is the associated symbol for $\EE_{x_1, x_2}:$
        \[ \scalebox{0.8}{\bordermatrix{ 
        &-x_1  & \cdots & -x_2 & \cdots &\epsilon_\rho & \cdots &x_2-2&x_2 -1&x_2 & \cdots & x_1 -1&x_1 &x_1 +1 & x_1 +2 & \cdots & \alpha -1 &  \cr
        &\lhd & \cdots & \lhd &\cdots & \odot & \cdots & \rhd &\rhd & \rhd   &\cdots &\rhd &\rhd &\rhd &\rhd \cr
        &&&&&\odot & \cdots & \odot \cr
        &&&&&&& &  & \odot & \cdots  &\odot \cr
        &&&&&&&&&&&&&&\odot & \cdots & \odot 
        }}.\]
        When $\epsilon_\rho +1 \leq x_1 \leq \alpha -2$ and $x_1 - x_2 > 1$, define 
        \begin{flalign*}
            \EE_{x_1, x_2}:= &\{([x_1+2,-x_1]_\rho,\lfloor x_1 \rfloor, (-1)^{x_1 - \epsilon_\rho} \eta),([x_2, -x_2]_\rho,\lfloor x_2 \rfloor, (-1)^{x_2 - \epsilon_\rho } \eta),\\
            &([x_2-2, \epsilon_\rho]_\rho,0,\eta), ([x_1 -2, x_2]_\rho,0,(-1)^{x_2- \epsilon_\rho}\eta),\\
            &([\alpha -1, x_1 +2]_\rho, 0, (-1)^{x_1 +1 - \epsilon_\rho}\eta)\},
        \end{flalign*}
        Here is the associated symbol for $\EE_{x_1, x_2}:$ 
        \[\scalebox{0.8}{\bordermatrix{ 
        &-x_1  & \cdots &-x_2&\cdots& \epsilon_\rho & \cdots & x_2 -2& x_2 -1 & x_2 & \cdots & x_1 -2& x_1-1 &x_1 & x_1+1 & x_1 + 2\cdots &\alpha -1 \cr
        &\lhd & \cdots & \lhd&\cdots & \odot & \cdots & \cdots &\cdots  & \rhd &\cdots &\rhd &\rhd &\rhd&\rhd&\rhd\cr
        &&&\lhd & \cdots & \odot & \cdots & \cdots & \cdots& \rhd&\cr
        &&&&&\odot &\cdots & \odot \cr
        &&&&&&&&&\odot & \cdots & \odot \cr
        &&&&&&&&&&&&&&\odot & \cdots & \odot\cr
        }}.\]
    \end{enumerate}
\end{prop}
\begin{proof} The sufficient direction can be proven in a similar way as Proposition \ref{nontemp2E}, which we omit. Now we show the necessary direction. Assume $x_2 = x_1 + 1$, and $\pi_{x_1, x_2}$ is of Arthur type, then by Definition \ref{psiminus} and Theorem \ref{nontemp red}, there exists an extended multi-segment $\EE$ that contains two segments of the form $([x_1, -x_1]_\rho, *,*)$ and $([x_1 +1, -(x_1 -1)]_\rho,*,*)$ such that $\pi(\EE) = \pi_{sc}$. This is impossible, since the $L$-parameter of $\pi_{sc}$ must be multiplicity free. Similarly, we see that the representation $\pi_{x_1, x_2}$ where $x_1 = x_2$ cannot be of Arthur type. Therefore, it follows that $x_2 < x_1$. Furthermore, by Lemma \ref{absolutevalue} and Theorem \ref{scArthur}, we have that $\pi_{x_1, x_2}$ is of Arthur type only if $\epsilon_\rho +1 \leq x_1 \leq \alpha -2$. This proves the necessary direction. The rest follows from definition.  \end{proof}

To wrap up the case where $f(\pi) = 2$, we look at the final subcase, Case $(G)$, where $\pi_{x_1, x_2} = L(\Delta_{\rho}[-x_1, -x_1], \Delta_{\rho}[-x_2, -x_2-2]; \pi_{sc})$. The natural restriction from Langlands classification gives us $x_1 \geq \frac{1}{2}, x_2 > -1$ and $x_1 - x_2 \geq 1$. 
\begin{prop}\label{nontemp2G}
    Consider the representation 
    \begin{equation*}
        \pi_{x_1, x_2} = L(\Delta_{\rho}[-x_1, -x_1], \Delta_{\rho}[-x_2, -x_2-2]; \pi_{sc}),
    \end{equation*}
    where $x_1 \geq \frac{1}{2}, x_2 > -1$ and $x_1 - x_2 \geq 1$.
    \begin{enumerate}
        \item The representation $\pi_{x_1, x_2}$ is of Arthur type if and only if $\epsilon_\rho +1 \leq x_1 \leq \alpha$ and $x_2 \leq x_1 -3$, or when $(x_1, x_2) = (\alpha+1, \alpha-2)$ or $(1,0)$. 
        \item The representation $\pi_{x_1, x_2}$ is of critical type if and only if $(x_1,x_2) = (\alpha -1, \alpha -2),(\alpha, \alpha -2),(\alpha +1, \alpha -2),(\alpha, \alpha -1),(\alpha +1,\alpha -1),(\alpha +2,\alpha -1),(\alpha +1, \alpha),(\alpha +2,\alpha),(\alpha +3,\alpha),(\alpha, \alpha -3)$.
        \item Define $\EE_{x_1, x_2}$ in various cases as follows, then $\pi(\EE_{x_1, x_2}) = \pi_{x_1, x_2}$. When $\epsilon_\rho +1 \leq x_1 \leq \alpha$ and $x_2 = x_1 -3$, define
        \begin{flalign*}
            \EE_{x_1, x_2} := &\{([x_1,-x_1]_\rho,\lfloor x_1 \rfloor, \eta),
            ([x_1 -2,-(x_1 - 3)])_\rho,0,-\eta), 
            ([x_1- 5, \epsilon_\rho]_\rho,0,-\eta), \\
            &([\alpha -1, x_1 -1]_\rho,0,(-1)^{x_2+1- \epsilon_\rho}\eta).
        \end{flalign*}
        Here is the associated symbol for $\EE_{x_1, x_2}:$
        \[ \scalebox{0.8}{\bordermatrix{ 
        &-x_1  & -(x_1 -1) & -(x_1 - 2)& -(x_1 -3) & \cdots &\epsilon_\rho & \cdots &x_1 -5&\cdots &x_1 -2 & x_1 -1 & x_1 & \cdots & \alpha -1 \cr
        &\lhd & \lhd & \lhd &\lhd &\cdots & \odot & \cdots & \rhd &\cdots & \rhd & \rhd   &\rhd\cr
        &&&&\lhd & \cdots & \odot & \cdots & \rhd & \cdots & \rhd\cr
        &&&&&&\odot & \cdots & \odot \cr
        &&&&&&& &   \odot & \cdots  &\odot \cr
        &&&&&&&&&&&\odot & \odot & \cdots & \odot 
        }}.\]
        When $\epsilon_\rho +1 \leq x_1 \leq \alpha $ and $x_1 - x_2 > 3$, define 
        \begin{flalign*}
            \EE_{x_1, x_2}:= &\{([x_1,-x_1]_\rho,\lfloor x_1 \rfloor, (-1)^{x_1 - \epsilon_\rho} \eta),([x_2+2, -x_2]_\rho,\lfloor x_2 \rfloor, (-1)^{x_2 - \epsilon_\rho } \eta),\\
            &([x_2-2, \epsilon_\rho]_\rho,0,\eta), ([x_1 -2, x_2+2]_\rho,0,(-1)^{x_2- \epsilon_\rho}\eta),\\
            &([\alpha -1, x_1]_\rho, 0, (-1)^{x_1 +1 - \epsilon_\rho}\eta)\}.
        \end{flalign*}
        Here is the associated symbol. 
        \[\EE_{x_1, x_2}= \scalebox{0.8}{\bordermatrix{ 
        &-x_1  & \cdots &-x_2&\cdots& \epsilon_\rho & \cdots & x_2 -2&  \cdots & x_2 +2& \cdots & x_1-2 &x_1 -1  & x_1 & \cdots &\alpha -1 \cr
        &\lhd & \cdots & \lhd&\cdots & \odot & \cdots & \odot &\cdots  & \rhd &\cdots &\rhd &\rhd &\rhd\cr
        &&&\lhd & \cdots & \odot & \cdots & \cdots & \cdots& \rhd&\cr
        &&&&&\odot &\cdots & \odot \cr
        &&&&&&&&&\odot & \cdots & \odot \cr
        &&&&&&&&&&&&&\odot & \cdots & \odot\cr
        }}.\]
        When $(x_1, x_2) = (\alpha+1, \alpha-2)$,
        define
        \begin{flalign*}
            \EE_{x_1, x_2} :=\{([x_1,-x_1]_\rho,\lfloor x_1 \rfloor, \eta),
            ([x_1 -2,-(x_1 - 3)])_\rho,0,-\eta), 
            ([x_1- 4, \epsilon_\rho]_\rho,0,-\eta)\}.\\
        \end{flalign*}
        Here is the associated symbol for $\EE_{x_1, x_2}$:
        \[ \scalebox{0.8}{\bordermatrix{ 
        &-x_1  & -(x_1 -1) & -(x_1 - 2)& -(x_1 -3) & \cdots &\epsilon_\rho & \cdots &x_1 -5& x_1-4 & x_1 -3 &x_1 -2 & x_1 -1 & x_1 \cr
        &\lhd & \lhd & \lhd &\lhd &\cdots & \odot & \cdots & \rhd &\rhd & \rhd & \rhd   &\rhd &\rhd\cr
        &&&&\lhd & \cdots & \odot & \cdots & \rhd & \rhd & \rhd  & \rhd\cr
        &&&&&&\odot & \cdots & \odot &\odot \cr
        }}.\]
        When $(x_1, x_2) = (1,0)$, define 
        \[\EE = \{([1,-1]_\rho, 1,\eta),([\alpha-1, 0]_\rho,1,-\eta), ([2,2]_\rho, 0,\eta)\}.\]
        Then $\pi_{1,0} = \pi(\EE)$. Here is the associated symbol. 
        \[\EE_{x_1, x_2}= \scalebox{0.8}{\bordermatrix{ 
        &-1 & 0 & 1 & 2 & \ldots & \alpha -1 \cr
        & \lhd & \odot & \rhd \cr
        &&\lhd & \odot &\odot & \dots & \rhd \cr
        &&&&\odot
        }}.\]
    \end{enumerate}
\end{prop}
\begin{proof} The proof is similar to that of Proposition $\ref{nontemp2E}$, which we omit. \end{proof} 

With this, we are done with the classification of non-tempered representations $\pi$ of good parity with $f(\pi) = 2$. We move onto the cases $f(\pi) = 3$ and $f(\pi) = 4$ now. 

 \section{Classification of corank 4 non-tempered representations of good parity \texorpdfstring{($f(\pi) = 3,4$)}{}}\label{classnontempcorank4,34}
 
In this section, we consider non-tempered representations $\pi$ of corank $4$, where there are exactly $3$ or $4$ segments in the $L$-data of $\pi$, i.e. $f(\pi) = 3, 4$.

\subsection{The case \texorpdfstring{$f(\pi) = 3$}{}}
We begin with the case where $f(\pi) = 3$.  In this case, we can group them into the following subcases: 
\begin{enumerate}[label = (\Alph*)]
    \item $\pi = L(\Delta_{\rho}[-x_1, -x_1], \Delta_{\rho}[-x_2, -x_2], \Delta_{\rho}[-x_3, -x_3]; \pi_{temp})$, where $x_1 > x_2 > x_3 \geq \frac{1}{2}$ and $\pi_{temp}$ is tempered of corank $1$. 
    \item $\pi = L(\Delta_{\rho}[-x_1, -x_1], \Delta_{\rho}[-x_1, -x_1], \Delta_{\rho}[-x_2, -x_2]; \pi_{temp})$, where $x_1 > x_2 \geq \frac{1}{2}$, and $\pi_{temp}$ is tempered of corank $1$. 
    \item $\pi = L(\Delta_{\rho}[-x_1, -x_1], \Delta_{\rho}[-x_2, -x_2], \Delta_{\rho}[-x_2, -x_2]; \pi_{temp})$, where $x_1 > x_2 \geq \frac{1}{2}$, and $\pi_{temp}$ is tempered of corank $1$.
    \item  $\pi = L(\Delta_{\rho}[-x, -x], \Delta_{\rho}[-x, -x], \Delta_{\rho}[-x, -x]; \pi_{temp})$, where $x \geq \frac{1}{2}$ and $\pi_{temp}$ is tempered of corank $1$. 
    \item $\pi = L(\Delta_{\rho}[-x_1, -x_1 -1], \Delta_{\rho}[-x_2, -x_2], \Delta_{\rho}[-x_3, -x_3]; \pi_{sc})$, where $x_1 \geq x_2 \geq x_3$. 
    \item $\pi = L(\Delta_{\rho}[-x_1, -x_1 ], \Delta_{\rho}[-x_2, -x_2-1], \Delta_{\rho}[-x_3, -x_3]; \pi_{sc})$, where $x_1 > x_2 \geq x_3$. 
    \item $\pi = L(\Delta_{\rho}[-x_1, -x_1 ], \Delta_{\rho}[-x_2, -x_2], \Delta_{\rho}[-x_3, -x_3-1]; \pi_{sc})$, where $x_1 \geq x_2 > x_3$. 
\end{enumerate}

\subsection{\texorpdfstring{Case $(A):\pi = L(\Delta_{\rho}[-x_2, -x_2], \Delta_{\rho}[-x_1, -x_1], \Delta_{\rho}[-x_3, -x_3]; \pi_{temp})$}{}}

We begin with Case $(A)$, where we denote $\pi_{x_1, x_2, x_3} = L(\Delta_{\rho}[-x_1, -x_1], \Delta_{\rho}[-x_2, -x_2], \Delta_{\rho}[-x_3, -x_3]; \pi_{temp})$ 
for $x_1 > x_2 > x_3 \geq \frac{1}{2}$ and $\pi_{temp}$ is a corank 1 tempered representation of good parity. There are $3$ cases to consider. The first one is $\pi_{temp} = T_{I,1}^{\alpha}(\pi_{sc})$
\begin{prop}\label{nontemp3A1}
    Consider the representation 
    \begin{equation*}
        \pi_{x_1, x_2, x_3} = L(\Delta_{\rho}[-x_1, -x_1], \Delta_{\rho}[-x_2, -x_2], \Delta_{\rho}[-x_3, -x_3]; \pi_{temp}),
    \end{equation*}
    for $x_1 > x_2 > x_3 \geq \frac{1}{2}$, where $\pi_{temp} = T_{I,1}^{\alpha}(\pi_{sc})$ and $\alpha > 0$. 
    \begin{enumerate}[label = (\roman*)]
        \item The representation $\pi_{x_1, x_2, x_3}$ is of Arthur type if and only if one of the following holds: 
        \begin{itemize}
            \item $\epsilon_\rho +1 \leq x_1 \leq \alpha -1$ and $\frac{1}{2} \leq x_3 < x_2 \leq \alpha -2$.
            \item $(x_1, x_2) = (\alpha, \alpha -1)$ and $\frac{1}{2} < x_3 \leq \alpha -2$.
            \item $(x_1, x_2, x_3) = (\alpha +1, \alpha, \alpha -1)$.
            \item $(x_1, x_2, x_3, \alpha) = (\frac{5}{2}, \frac{3}{2}, \frac{1}{2}, \frac{1}{2})$.            
        \end{itemize}
        \item The representation $\pi_{x_1, x_2, x_3}$ is of critical type in the following cases: 
        \begin{itemize}
            \item $(x_1, x_2,x_3)= (\alpha-1, \alpha -2, \alpha-3)$. 
            \item $(x_1, x_2, x_3) = (\alpha, \alpha -1, \alpha -2)$. 
            \item $(x_1, x_2, x_3) = (\alpha+1, \alpha, \alpha -1)$.
            \item $(x_1, x_2, x_3) = (\alpha+2, \alpha +1, \alpha)$ or $(\alpha +3, \alpha +2, \alpha +1)$.
        \end{itemize}
        \item Define $\EE_{x_1, x_2, x_3}$ in various cases as follows. Then $\pi(\EE_{x_1, x_2, x_3}) = \pi_{x_1, x_2, x_3}$. When $\epsilon_\rho +1 \leq x_1 \leq \alpha+1$ and $\frac{1}{2} \leq x_3 < x_2 \leq \alpha$, with $x_1 = x_2 +1 = x_3 +2$, define 
        \begin{flalign*}
            \EE_{x_1, x_2, x_3}:= &\{[x_1, -x_1]_\rho,\lfloor x_1 \rfloor, \eta),([x_3 -2, \epsilon_\rho]_\rho,0,-\eta), \\&([\alpha -2, x_2]_\rho,0,(-1)^{x_2 - \epsilon_\rho} \eta),([\alpha,\alpha]_\rho,0,(-1)^{\alpha - \epsilon_\rho}\eta)\}.
        \end{flalign*}
        Here is the associated symbol. 
        \[\EE_{x_1, x_2,x_3}= \scalebox{0.8}{\bordermatrix{ 
        &-x_1  & \cdots & \epsilon_\rho & \cdots & x_3 -2&  \cdots & x_2 &x_1 & \cdots & \alpha -2 &\alpha-1 & \alpha \cr
        &\lhd & \cdots & \odot  & \cdots& \rhd & \cdots & \rhd &\rhd\cr
        &&&\odot & \cdots & \odot\cr
        &&&&&&&\odot & \odot & \cdots & \odot \cr
        &&&&&&&&&&&&\odot\cr
        }}.\]
         When $\epsilon_\rho +1 \leq x_1 \leq \alpha +1$ and $\frac{1}{2} \leq x_3 < x_2 \leq \alpha$, with $x_1 = x_2 +1  > x_3 +2$, define 
        \begin{flalign*}
            \EE_{x_1, x_2, x_3}:= &\{[x_1, -x_1]_\rho,\lfloor x_1 \rfloor, \eta),
            ([x_3, -x_3]_\rho,\lfloor x_3 \rfloor, - \eta),([x_3 -2, \epsilon_\rho]_\rho,0,-\eta), \\&([x_2 -2, x_3]_\rho,0,(-1)^{x_3+1 - \epsilon_\rho} \eta),([\alpha-2, x_2]_\rho,0,(-1)^{x_2 - \epsilon_\rho}\eta),([\alpha,\alpha]_\rho,0,(-1)^{\alpha - \epsilon_\rho}\eta)\}.
        \end{flalign*}
        Here is the associated symbol for  $\EE_{x_1, x_2,x_3}$: 
        \[\scalebox{0.8}{\bordermatrix{ 
        &-x_1  & \cdots &-x_3 & \cdots & \epsilon_\rho & \cdots & x_3 -2&  x_3 -1 & x_3 & \cdots & x_2 -2 & x_2 -1 & x_2 &x_1 & \cdots & \alpha -2 &\alpha-1 & \alpha \cr
        &\lhd & \cdots & \lhd & \cdots & \odot & \cdots & \rhd & \rhd & \rhd & \cdots & \rhd & \rhd & \rhd & \rhd & \rhd \cr
        &&&\lhd & \cdots & \odot & \cdots & \rhd & \rhd & \rhd\cr
        &&&&&\odot & \cdots & \odot \cr
        &&&&&&&&&\odot & \cdots & \odot \cr
        &&&&&&&&&&&&&\odot & \odot & \cdots & \odot \cr
        &&&&&&&&&&&&&&&&&&\odot\cr
        }}.\]
        When $\epsilon_\rho +1 \leq x_1 \leq \alpha -1$ and $\frac{1}{2} \leq x_3 < x_2 \leq \alpha-2$, with $x_1 > x_2 +1 = x_3 + 2$, define 
        \begin{flalign*}
             \EE_{x_1, x_2, x_3} := &\{[x_1, -x_1]_\rho,\lfloor x_1 \rfloor, \eta),
            ([x_2, -x_2]_\rho,\lfloor x_2 \rfloor, - \eta),([x_3 -2, \epsilon_\rho]_\rho,0,-\eta), \\
            &([x_1-2, x_3]_\rho,0,(-1)^{x_3+1 - \epsilon_\rho} \eta), ([\alpha-2, x_1]_\rho,0,(-1)^{x_1 - \epsilon_\rho}\eta),([\alpha,\alpha]_\rho,0,(-1)^{\alpha - \epsilon_\rho}\eta)\}.
        \end{flalign*}
        Here is the associated symbol for $\EE_{x_1, x_2,x_3}$: 
        \[ \scalebox{0.8}{\bordermatrix{ 
        &-x_1  & \cdots &-x_2 & \cdots & \epsilon_\rho & \cdots & x_3 -2&  x_3 -1 & x_3 & x_2 &\cdots & x_1 -2 & x_1 -1  &x_1 & \cdots & \alpha -2 &\alpha-1 & \alpha \cr
        &\lhd & \cdots & \lhd & \cdots & \odot & \cdots & \rhd & \rhd & \rhd & \rhd& \cdots & \rhd & \rhd & \rhd  \cr
        &&&\lhd & \cdots & \odot & \cdots & \rhd & \rhd & \rhd&\rhd\cr
        &&&&&\odot & \cdots & \odot \cr
        &&&&&&&&&\odot &\odot & \cdots & \odot \cr
        &&&&&&&&&&&&&&\odot & \cdots & \odot \cr
        &&&&&&&&&&&&&&&&&&\odot\cr
        }}.\]
        When $\epsilon_\rho +1 \leq x_1 \leq \alpha -1$ and $\frac{1}{2} \leq x_3 < x_2 \leq \alpha -2$, with $x_1 > x_2 +1 > x_3 + 2$, define 
        \begin{flalign*}
            \EE_{x_1, x_2, x_3}:= &\{[x_1, -x_1]_\rho,\lfloor x_1 \rfloor, \eta),
            ([x_2, -x_2]_\rho,\lfloor x_2 \rfloor, -\eta),([x_3, -x_3]_\rho,\lfloor x_3 \rfloor,\eta),([x_3 -2, \epsilon_\rho]_\rho,0,-\eta)\\
            &([x_1-2, x_3]_\rho,0,(-1)^{x_3+1 - \epsilon_\rho} \eta), 
            ([x_2 -2, x_3]_\rho,0,(-1)^{x_3 -\epsilon_\rho}), ([x_1 -2, x_2]_\rho,0,(-1)^{x_2 - \epsilon_\rho}\eta),\\
            &([\alpha -2, x_1]_\rho,0,(-1)^{x_1 - \epsilon_\rho}\eta), 
            ([\alpha,\alpha]_\rho,0,(-1)^{\alpha - \epsilon_\rho}\eta)\}.
        \end{flalign*}
        Here is the associated symbol for $\EE_{x_1, x_2,x_3}$: 
        \[ \scalebox{0.6}{\bordermatrix{ 
        &-x_1  & \cdots &-x_2 & \cdots & -x_3 & \cdots &\epsilon_\rho & \cdots & x_3 -2&  x_3 -1 & x_3 & \cdots &x_2 -2& x_2 -1 & x_2 & \cdots & x_1 -2 & x_1 -1  &x_1 & \cdots & \alpha -2 &\alpha-1 & \alpha \cr
        &\lhd & \cdots & \lhd & \cdots & \lhd & \cdots & \odot & \cdots & \rhd & \rhd & \rhd & \cdots & \rhd & \rhd & \rhd&\cdots & \rhd & \rhd & \rhd  \cr
        &&&\lhd & \cdots & \lhd & \cdots & \odot & \cdots & \rhd & \rhd & \rhd&\cdots&\rhd & \rhd &\rhd\cr
        &&&&&\lhd & \cdots & \odot & \cdots & \rhd & \rhd & \rhd \cr
        &&&&&&&\odot & \cdots & \odot \cr
        &&&&&&&&&&&\odot & \cdots & \odot \cr
        &&&&&&&&&&&&&&&\odot  & \cdots & \odot \cr
        &&&&&&&&&&&&&&&&&&&\odot & \cdots & \odot \cr
        &&&&&&&&&&&&&&&&&&&&&&&\odot\cr
        }}.\]
        When $(x_1, x_2, x_3, \alpha) = (\frac{5}{2}, \frac{3}{2}, \frac{1}{2}, \frac{1}{2})$, define 
        \[\EE_{\frac{5}{2}, \frac{3}{2}, \frac{1}{2}}= \scalebox{0.8}{\bordermatrix{ 
         &-\frac{5}{2}&-\frac{3}{2}&-\frac{1}{2}&\frac{1}{2}&\frac{3}{2}&\frac{5}{2}\cr
        &\lhd &\lhd &\lhd  & \rhd &\rhd &\rhd \cr
        &&&&\oplus
        }}.\]
    \end{enumerate}
\end{prop}
\begin{proof}The sufficient direction can be proven in a similar way as Proposition \ref{nontemp2A1}, which we omit. Now let's show the necessary direction. 

First, by \cite[Proposition $11.15$]{HJLLZ24}, we have that: 
\begin{equation*}
    \pi_{x_1, x_2, x_3}^{\rho,-} = L(\Delta_{\rho}[-x_2, -x_2], \Delta_{\rho}[-x_3, -x_3]; \pi_{temp})
\end{equation*}
is of Arthur type if and only if $\frac{1}{2} \leq x_3 < x_2 \leq \alpha -1$, $(x_2, x_3) = (\alpha, \alpha -1)$ or $(x_2, x_3, \alpha) = (\frac{3}{2}, \frac{1}{2}, \frac{1}{2})$. Furthermore, if $\pi_{x_1, x_2, x_3}$ is of Arthur type, then by Lemma \ref{absolutevalue}, we require that 
\begin{equation*}
    x_1 - 1 \in |\Omega|(\pi_{x_1, x_2, x_3}(\pi_{sc})^{\rho,-}) \subseteq \{\rho\lvert \cdot \rvert^{y}: \epsilon_\rho \leq y \leq \alpha -2\} \cup \{\rho\lvert \cdot \rvert^{\alpha}\} \cup \{\rho\lvert \cdot \rvert^{x_2}. \rho\lvert \cdot \rvert^{-x_2}, \rho\lvert \cdot \rvert^{x_3}, \rho\lvert \cdot \rvert^{-x_3}\}
\end{equation*}
We conclude that $\pi_{x_1, x_2, x_3}$ is of Arthur type only if one of the following holds: 
\begin{itemize}
    \item $\epsilon_\rho +1 \leq x_1 \leq \alpha -1$ and $\frac{1}{2} \leq x_3 < x_2 \leq \alpha -2$.
    \item $(x_1, x_2) = (\alpha, \alpha -1)$ and $\frac{1}{2} \leq x_3 \leq \alpha -2$
    \item $(x_1, x_2, x_3) = (\alpha +1, \alpha, \alpha -1)$
    \item $x_1 = \alpha +1$ and $\frac{1}{2} \leq x_3 < x_2 \leq \alpha -1$. 
\end{itemize}
Therefore it suffices to eliminate the last case. This follows from the exact same argument used in Proposition \ref{nontemp2A1}, which we omit. This proves part $(i)$ of the proposition. Part $(ii)$ follows from definition. \end{proof}

The next case to consider is $\pi_{temp} = T_{IV, 3}(\pi_{sc})$. 
\begin{prop}\label{nontemp3A2}
    Consider the representation 
    \begin{equation*}
        \pi_{x_1, x_2, x_3} = L(\Delta_{\rho}[-x_1, -x_1], \Delta_{\rho}[-x_2, -x_2], \Delta_{\rho}[-x_3, -x_3]; \pi_{temp}),
    \end{equation*}
    for $x_1 > x_2 > x_3 \geq \frac{1}{2}$, where $\pi_{temp} = T_{IV, 3}(\pi_{sc})$ and $\alpha \in \mathbb{Z}_{>0}$. 
    \begin{enumerate}[label = (\roman*)]
        \item The representation $\pi_{x_1, x_2, x_3}(\pi_{sc})$ is of Arthur type if and only if one of the following holds: 
        \begin{itemize}
            \item $1 \leq x_1 \leq \alpha$ and $1 \leq x_3 < x_2 \leq \alpha -1$.
            \item $(x_1, x_2) = (\alpha +1, \alpha)$ and $1 \leq x_3 \leq \alpha -1$. 
            \item $(x_1, x_2, x_3) = (\alpha +2, \alpha +1, \alpha)$. 
        \end{itemize}
        \item The representation $\pi_{x_1, x_2, x_3}(\pi_{sc})$ is of critical type in the following cases: 
        \begin{itemize}
            \item $(x_1, x_2, x_3) = (\alpha, \alpha -1, \alpha -2)$ for $\alpha = 3$.
            \item $(x_1, x_2, x_3) = (\alpha+1, \alpha, \alpha-1)$ for $\alpha = 2$.
            \item $(x_1, x_2, x_3) = (\alpha +2, \alpha+1, \alpha)$ for $\alpha = 1$.
        \end{itemize}
        \item 
        Define $\EE_{x_1, x_2, x_3}$ in various cases as follows. Then $\pi(\EE_{x_1, x_2, x_3}) = \pi_{x_1, x_2, x_3}$. When $1 \leq x_1 \leq \alpha+2$ and $1 \leq x_3 < x_2 \leq \alpha+1$, with $x_1 = x_2 +1 = x_3 +2$, define 
        \begin{flalign*}
            \EE_{x_1, x_2, x_3}:= &\{[x_1, -x_1]_\rho,\lfloor x_1 \rfloor, \eta),
            ([0,0]_\rho,0,\eta)^2,([x_3 -2, 0]_\rho,0,-\eta), \\&([\alpha -1, x_2]_\rho,0,(-1)^{x_2} \eta).
        \end{flalign*}
        Here is the associated symbol. 
        \[\EE_{x_1, x_2,x_3}= \scalebox{0.8}{\bordermatrix{ 
        &-x_1  & \cdots & 0 & \cdots & x_3 -2&  \cdots & x_2 &x_1 & \cdots  &\alpha-1 \cr
        &\lhd & \cdots & \odot &\cdots & \rhd & \cdots & \rhd & \rhd \cr
        &&&\odot \cr
        &&&\odot \cr
        &&&\odot & \cdots & \odot\cr
        &&&&&&&\odot & \odot & \cdots & \odot \cr
        }}.\]
         When $ 1 \leq x_1 \leq \alpha +1$ and $1 \leq x_3 < x_2 \leq \alpha$, with $x_1 = x_2 +1  > x_3 +2$, define 
        \begin{flalign*}
            \EE_{x_1, x_2, x_3}:= &\{[x_1, -x_1]_\rho,\lfloor x_1 \rfloor, \eta),
            ([x_3, -x_3]_\rho,\lfloor x_3 \rfloor, - \eta),([0,0]_\rho,0,-\eta)^2,([x_3 -2, 0]_\rho,0,-\eta), \\&([x_2 -2, x_3]_\rho,0,(-1)^{x_3+1} \eta),([\alpha-1, x_2]_\rho,0,(-1)^{x_2}\eta).
        \end{flalign*}
        Here is the associated symbol for  $\EE_{x_1, x_2,x_3}$: 
        \[ \scalebox{0.8}{\bordermatrix{ 
        &-x_1  & \cdots &-x_3 & \cdots & 0 & \cdots & x_3 -2&  x_3 -1 & x_3 & \cdots & x_2 -2 & x_2 -1 & x_2 &x_1 & \cdots &\alpha-1 \cr
        &\lhd & \cdots & \lhd & \cdots & \odot & \cdots & \rhd & \rhd & \rhd & \cdots & \rhd & \rhd & \rhd & \rhd \cr
        &&&\lhd & \cdots & \odot & \cdots & \rhd & \rhd & \rhd\cr
        &&&&&\odot \cr
        &&&&&\odot \cr
        &&&&&\odot & \cdots & \odot \cr
        &&&&&&&&&\odot & \cdots & \odot \cr
        &&&&&&&&&&&&&\odot & \odot & \cdots & \odot \cr
        }}.\]
        When $1 \leq x_1 \leq \alpha -1$ and $1 \leq x_3 < x_2 \leq \alpha-2$, with $x_1 > x_2 +1 = x_3 + 2$, define 
        \begin{flalign*}
            \EE_{x_1, x_2, x_3}:= &\{[x_1, -x_1]_\rho,\lfloor x_1 \rfloor, \eta),
            ([x_2, -x_2]_\rho,\lfloor x_2 \rfloor, - \eta),([0,0]_\rho,0,-\eta)^2,([x_3 -2, 0]_\rho,0,-\eta),\\
            &([x_1-2, x_3]_\rho,0,(-1)^{x_3+1 } \eta),
             ([\alpha-1, x_1]_\rho,0,(-1)^{x_1}\eta)\}.
        \end{flalign*}
        Here is the associated symbol for $\EE_{x_1, x_2,x_3}$:
        \[ \scalebox{0.8}{\bordermatrix{ 
        &-x_1  & \cdots &-x_2 & \cdots & 0 & \cdots & x_3 -2&  x_3 -1 & x_3 & x_2 &\cdots & x_1 -2 & x_1 -1  &x_1 & \cdots &\alpha-1  \cr
        &\lhd & \cdots & \lhd & \cdots & \odot & \cdots & \rhd & \rhd & \rhd & \rhd& \cdots & \rhd & \rhd & \rhd  \cr
        &&&\lhd & \cdots & \odot & \cdots & \rhd & \rhd & \rhd&\rhd\cr
        &&&&&\odot \cr
        &&&&&\odot\cr
        &&&&&\odot & \cdots & \odot \cr
        &&&&&&&&&\odot &\odot & \cdots & \odot \cr
        &&&&&&&&&&&&&&\odot & \cdots & \odot \cr
        }}.\]
        When $1 \leq x_1 \leq \alpha -1$ and $1 \leq x_3 < x_2 \leq \alpha -2$, with $x_1 > x_2 +1 > x_3 + 2$, define 
        \begin{flalign*}
            \EE_{x_1, x_2, x_3}:= &\{[x_1, -x_1]_\rho,\lfloor x_1 \rfloor, \eta),
            ([x_2, -x_2]_\rho,\lfloor x_2 \rfloor, -\eta),([x_3, -x_3]_\rho,\lfloor x_3 \rfloor,\eta),([x_3 -2, 0]_\rho,0,-\eta)\\
            &([x_1-2, x_3]_\rho,0,(-1)^{x_3+1} \eta), 
            ([x_2 -2, x_3]_\rho,0,(-1)^{x_3 -\epsilon_\rho}), ([x_1 -2, x_2]_\rho,0,(-1)^{x_2 }\eta),\\
            &([\alpha -1, x_1]_\rho,0,(-1)^{x_1}\eta).
        \end{flalign*}
        Here is the associated symbol for $\EE_{x_1, x_2,x_3}$:  
        \[\scalebox{0.6}{\bordermatrix{ 
        &-x_1  & \cdots &-x_2 & \cdots & -x_3 & \cdots &\epsilon_\rho & \cdots & x_3 -2&  x_3 -1 & x_3 & \cdots &x_2 -2& x_2 -1 & x_2 & \cdots & x_1 -2 & x_1 -1  &x_1 & \cdots&\alpha-1\cr
        &\lhd & \cdots & \lhd & \cdots & \lhd & \cdots & \odot & \cdots & \rhd & \rhd & \rhd & \cdots & \rhd & \rhd & \rhd&\cdots & \rhd & \rhd & \rhd  \cr
        &&&\lhd & \cdots & \lhd & \cdots & \odot & \cdots & \rhd & \rhd & \rhd&\cdots&\rhd & \rhd &\rhd\cr
        &&&&&\lhd & \cdots & \odot & \cdots & \rhd & \rhd & \rhd \cr
        &&&&&&&\odot & \cdots & \odot \cr
        &&&&&&&&&&&\odot & \cdots & \odot \cr
        &&&&&&&&&&&&&&&\odot  & \cdots & \odot \cr
        &&&&&&&&&&&&&&&&&&&\odot & \cdots & \odot \cr
        }}.\]
    \end{enumerate}
\end{prop}
\begin{proof} By \cite[Proposition $11.16$]{HJLLZ24}, we know that the representation 
\begin{equation*}
    \pi_{x_1, x_2, x_3}^{\rho,-} = L(\Delta_{\rho}[-x_2, -x_2], \Delta_{\rho}[-x_3, -x_3]; \pi_{temp})
\end{equation*}
is of Arthur type if and only if $1 \leq x_3 < x_2 \leq \alpha$ or $(x_2, x_3) = (\alpha +1, \alpha)$. The rest of the proof follows in a similar way as Proposition \ref{nontemp3A1}, which we omit. \end{proof}

The last case to consider in Case $(A)$ is $\pi_{temp} = T_{V,2}^{\pm}(\pi_{sc})$. 
\begin{prop}\label{nontemp3A3}
     Consider the representation 
    \begin{equation*}
        \pi_{x_1, x_2, x_3}^{\pm} = L(\Delta_{\rho}[-x_1, -x_1], \Delta_{\rho}[-x_2, -x_2], \Delta_{\rho}[-x_3, -x_3]; \pi_{temp}^{\pm}),
    \end{equation*}
    for $x_1 > x_2 > x_3 \geq \frac{1}{2}$, where $\pi_{temp}^{\pm} = T_{V,2}(\pi_{sc})^{\pm}$ and $\alpha = 0$. 
    \begin{enumerate}
        \item The representation $\pi_{x_1, x_2, x_3}^{\pm}$ is of Arthur type if and only if $(x_1, x_2, x_3) = (3,2,1)$. 
        \item The representation $\pi_{3,2,1}^{\pm}$ is of critical type if and only if $(x_1, x_2, x_3) = (3,2,1)$. 
        \item We have $\pi(\EE^{\pm}) = \pi_{x_1, x_2, x_3}^{\pm}$, where 
        \[\EE^{+}= \scalebox{0.8}{\bordermatrix{ 
         &-3 &-2 &-1 &0 & 1 &2 &3  \cr
        &\lhd &\lhd &\lhd &\oplus &\rhd &\rhd &\rhd\cr
        &&&&\oplus
        }}, \quad \EE^{-}= \scalebox{0.8}{\bordermatrix{ 
         &-3 &-2 &-1 &0 & 1 &2 &3  \cr
        &\lhd &\lhd &\lhd &\oplus &\rhd &\rhd &\rhd\cr
        &&&&\ominus
        }}.\]
    \end{enumerate}
\end{prop}
\begin{proof} This follows immediately from Lemma \ref{absolutevalue} and \cite[Proposition $11.17$]{HJLLZ24}, which state that the representation 
\begin{equation*}
    (\pi_{x_1, x_2, x_3}^{\pm})^{\rho,-} = L(\Delta_{\rho}[-x_2, -x_2], \Delta_{\rho}[-x_3, -x_3]; \pi_{temp}^{\pm})
\end{equation*}
is of Arthur type if and only if $(x_2, x_3) = (2,1)$. \end{proof}

This concludes all the cases in Case $(A)$. Now we move onto Case $(B)$. 
\subsection{\texorpdfstring{Case $(B): \pi = L(\Delta_{\rho}[-x_1, -x_1], \Delta_{\rho}[-x_1, -x_1], \Delta_{\rho}[-x_2, -x_2]; \pi_{temp})$}{}}
In this subsection we investigate representations of the form \[\pi_{x_1, x_2} = L(\Delta_{\rho}[-x_1, -x_1], \Delta_{\rho}[-x_1, -x_1], \Delta_{\rho}[-x_2, -x_2]; \pi_{temp})\] for $x_1 > x_2 \geq \frac{1}{2}$ where $\pi_{temp}$ is a corank 1 tempered representation of good parity. Again, there are three cases to consider. The first one is $\pi_{temp} = T_{I,1}^{\alpha}(\pi_{sc})$. 
\begin{prop}\label{nontemp3B1}
    Consider the representation 
    \begin{equation*}
        \pi_{x_1, x_2} = L(\Delta_{\rho}[-x_1, -x_1], \Delta_{\rho}[-x_1, -x_1], \Delta_{\rho}[-x_2, -x_2]; \pi_{temp}),
    \end{equation*}
    for $x_1 > x_2 \geq \frac{1}{2}$, where $\pi_{temp} = T_{I,1}^{\alpha}(\pi_{sc})$ and $\alpha > 0$. 
    \begin{enumerate}[label = (\roman*)]
        \item The representation $\pi_{x_1, x_2}$ is of Arthur type if and only if $(x_1, x_2, \alpha) = (\frac{3}{2}, \frac{1}{2}, \frac{1}{2})$. 
        \item The representation $\pi_{x_1, x_2}$ is of critical type if and only if $(x_1, x_2) = (\alpha +2, \alpha +1),(\alpha +1, \alpha), (\alpha +1, \alpha -1), (\alpha, \alpha -1)$ or $ (\alpha -1, \alpha -2)$. 
        \item When $\alpha = \frac{1}{2}$, we have $\pi(\EE) =\pi_{\frac{3}{2}, \frac{1}{2}}$, where
        \[\EE= \scalebox{0.8}{\bordermatrix{ 
         &-\frac{3}{2} &-\frac{1}{2} &\frac{1}{2}&\frac{3}{2}  \cr
        &\lhd &\odot &\odot &\rhd\cr
        &\lhd &\odot &\odot &\rhd\cr
        }}.\]
    \end{enumerate}
\end{prop}
\begin{proof} We only need to show the necessary direction. By \cite[Proposition $10.12$]{HJLLZ24}, the representation 
\begin{equation*}
    \pi_{x_1, x_2}^{\rho,-} = L(\Delta_{\rho}[-x_2, -x_2]; \pi_{temp})
\end{equation*}
is of Arthur type if and only if $\frac{1}{2} \leq x_2 \leq \alpha -1$ when $\alpha \geq \frac{3}{2}$, or $x_2 = \frac{1}{2} = \alpha$. On the other hand, Lemma \ref{absolutevalue} tells us that $\pi_{x_1, x_2}$ is of Arthur type only if $|\Omega|(\pi_{x_1, x_2}^{\rho,-})$ contains two copies of $\rho\lvert \cdot \rvert^{x_1 -1}$. We know that: 
\begin{equation*}
    |\Omega|(\pi_{x_1, x_2}^{\rho,-}) \subset \{\rho\lvert \cdot \rvert^{y}: 0 \leq y \leq \alpha -2\} \cup \{\rho\lvert \cdot \rvert^{\alpha}\} \cup \{\rho\lvert \cdot \rvert^{x_2}, \rho\lvert \cdot \rvert^{-x_2}\}.
\end{equation*}
By our assumption we cannot have $x_1= 1$. Therefore, from this we may conclude that $\pi_{x_1, x_2}$ can only be of Arthur type in the following cases: 
\begin{itemize}
    \item $1 \leq x_2 \leq \alpha -2$ and $x_1 = x_2 + 1$ for $\alpha \geq \frac{3}{2}$
    \item $(x_1, x_2, \alpha) = (\frac{3}{2}, \frac{1}{2}, \frac{1}{2})$
\end{itemize}
We show that the first case fails. Suppose $1\leq x_2 \leq \alpha -2$ and $x_1 = x_2 + 1$ for $\alpha \geq \frac{3}{2}$ and $\pi_{x_1, x_2}$ is of Arthur type. Then by Theorem \ref{nontemp red} and Definition \ref{psiminus}, there exists an extended multi-segment $\EE$ that contains two segments of the form $([x_2, -x_2]_\rho,*,*)$ such that $\pi(\EE) = L(\Delta_{\rho}[-x_2, -x_2]; \pi_{temp})$. However, by \cite[Proposition $10.12$]{HJLLZ24}, the multiplicity of $x_2 -1$ inside the extended multi-segment $\EE$ can only be $1$, so we reach a contradiction. This proves part $(i)$. Part $(ii)$ follows from definition. \end{proof}

The second case we need to consider is $\pi_{temp} = T_{IV, 3}(\pi_{sc})$. 
\begin{prop}\label{nontemp3B2}
     Consider the representation 
    \begin{equation*}
        \pi_{x_1, x_2} = L(\Delta_{\rho}[-x_1, -x_1], \Delta_{\rho}[-x_1, -x_1], \Delta_{\rho}[-x_2, -x_2]; \pi_{temp}),
    \end{equation*}
    for $x_1 > x_2 \geq 1$, where $\pi_{temp} = T_{IV, 3}(\pi_{sc})$ and $\alpha \in \mathbb{Z}_{> 0}$. 
    \begin{enumerate}[label = (\roman*)]
        \item The representation $\pi_{x_1, x_2}$ is not of Arthur type for any $x_1, x_2$. 
        \item The representation $\pi_{x_1, x_2}$ is of critical type if and only if one of the following holds: 
        \begin{itemize}
            \item $(x_1, x_2) = (\alpha +1, \alpha)$ and $\alpha = 1$.
            \item $(x_1, x_2) = (\alpha, \alpha -1)$ and $\alpha = 2$.
        \end{itemize} 
    \end{enumerate}
\end{prop}
\begin{proof} From \cite[Proposition $10.13$]{HJLLZ24}, the representation 
\begin{equation*}
    \pi_{x_1, x_2}^{\rho,-} = L(\Delta_{\rho}[-x_2, -x_2]; \pi_{temp})
\end{equation*}
can only be of Arthur type when $1 \leq x_2 \leq \alpha$. On the other hand, Lemma \ref{absolutevalue} tells us that $\pi_{x_1, x_2}$ is of Arthur type only if $|\Omega|(\pi_{x_1, x_2}^{\rho,-})$ contains two copies of $\rho\lvert \cdot \rvert^{x_1 -1}$. 
We know that: 
\begin{equation*}
    |\Omega|(\pi_{x_1, x_2}^{\rho,-}) \subset \{\rho\lvert \cdot \rvert^{y}: 0 \leq y \leq \alpha -1\} \cup \{\rho, \rho\} \cup \{\rho\lvert \cdot \rvert^{x_2}, \rho\lvert \cdot \rvert^{-x_2}\}.
\end{equation*}
By assumption we cannot have $x_1 = 1$, which means that $\pi_{x_1, x_2}$ can only be of Arthur type if $1 \leq x_2 \leq \alpha -1$ and $x_1 = x_2 +1$. However, by the same argument as in the proof of Proposition \ref{nontemp3B1}, we can disregard this case as well. Therefore $\pi_{x_1, x_2}$ is not of Arthur type for any $x_1, x_2$. \end{proof}

The last case to consider is $\pi_{temp} = T_{V,2}^{\pm}(\pi_{sc})$. 
\begin{prop}\label{nontemp3B3}
     Consider the representation 
    \begin{equation*}
        \pi_{x_1, x_2}^{\pm} = L(\Delta_{\rho}[-x_1, -x_1], \Delta_{\rho}[-x_1, -x_1], \Delta_{\rho}[-x_2, -x_2]; \pi_{temp}^{\pm}),
    \end{equation*}
    for $x_1 > x_2 \geq 1$, where $\pi_{temp}^{\pm}= T_{V,2}^{\pm}(\pi_{sc})$ and $\alpha =0$.
    \begin{enumerate}
        \item The representation $\pi_{x_1, x_2}^{\pm}$ is not of Arthur type for any $x_1, x_2$. 
        \item The representation $\pi_{x_1, x_2}$ is of critical type if and only if  $(x_1, x_2) = (2,1)$. 
    \end{enumerate}
\end{prop}
\begin{proof} By Lemma \ref{absolutevalue} and \cite[Proposition $10.14$]{HJLLZ24}, the representation $\pi_{x_1, x_2}^{\pm}$ can only be of Arthur type if $x_1 = 2, x_2 = 1$. It's easy to verify that this case does not work. This proves the proposition.  \end{proof}

This concludes our discussion in Case $(B)$. We move onto Case $(C)$. 
\subsection{\texorpdfstring{Case $(C): \pi = L(\Delta_{\rho}[-x_1, -x_1], \Delta_{\rho}[-x_2,-x_2],\Delta_{\rho}[-x_2,-x_2]; \pi_{temp})$}{}}
In this subsection, we consider representations of the form 
\[\pi_{x_1, x_2} = L(\Delta_{\rho}[-x_1, -x_1], \Delta_{\rho}[-x_2,-x_2],\Delta_{\rho}[-x_2,-x_2]; \pi_{temp}),\] where $x_1 > x_2 \geq \frac{1}{2}$ and $\pi_{temp}$ is a corank 1 tempered representation of good parity. This may look similar to the previous case, but since the Steinberg segments do not commute in general, we need to deal with them separately. The first case to consider is $\pi_{temp} = T_{I,1}^{\alpha}(\pi_{sc})$. 
\begin{prop}\label{nontemp3C1}
    Consider the representation 
    \begin{equation*}
        \pi_{x_1, x_2} = L(\Delta_{\rho}[-x_1,-x_1], \Delta_{\rho}[-x_2,-x_2], \Delta_{\rho}[-x_2,-x_2]; \pi_{temp}),
    \end{equation*}
    where $x_1 > x_2 \geq \frac{1}{2}$ and $\pi_{temp} = T_{I,1}^{\alpha}(\pi_{sc})$ for $\alpha > 0$. 
    \begin{enumerate}[label = (\roman*)]
        \item The representation $\pi_{x_1, x_2}$ is of Arthur type if and only if $(x_1, x_2) = (\frac{3}{2}, \frac{1}{2})$, or if $x_2 = \frac{1}{2}, \frac{3}{2} \leq x_1 \leq \alpha-1$. 
        \item The representation $\pi_{x_1, x_2}$ is of critical type if and only if 
        \begin{equation*}
            (x_1, x_2) = (\alpha +2, \alpha +1), (\alpha +1, \alpha), (\alpha +1, \alpha -1), (\alpha, \alpha -1), (\alpha -1, \alpha -2).
        \end{equation*}
        \item We have $\pi(\EE_{\frac{3}{2}, \frac{1}{2}}) = \pi_{\frac{3}{2}, \frac{1}{2}}$, where 
        \[\EE_{\frac{3}{2}, \frac{1}{2}}= \scalebox{0.8}{\bordermatrix{ 
         &-\frac{3}{2} &-\frac{1}{2} &\frac{1}{2}&\frac{3}{2} & \cdots & \alpha -1  \cr
        &\lhd &\lhd & \rhd &\rhd\cr
        &&\lhd &\rhd\cr
        &&&\odot & \odot & \cdots & \odot 
        }}.\]
        
    \end{enumerate}
\end{prop}
\begin{proof} The sufficient direction is easy to verify. Now we show the necessary direction. By \cite[Proposition $11.18$]{HJLLZ24}, the representation 
\begin{equation*}
    \pi_{x_1, x_2}^{\rho,-} = L(\Delta_\rho[-x_2, -x_2], \Delta_\rho[-x_2, -x_2]; \pi_{temp})
\end{equation*}
is of Arthur type if and only if $x_2 = \frac{1}{2}$. By Lemma \ref{absolutevalue}, the representation $\pi_{x_1, x_2}$ is only of Arthur type when the set $|\Omega|(\pi_{x_1, x_2}^{\rho,-})$ contains $\rho \lvert \cdot \rvert^{x_1}$. This happens only when $x_1 = \frac{3}{2}$, $x_1 \leq \alpha-1$, or $x_1 = \alpha+1$. The last case may be eliminated via a similar argument as before. 
\end{proof}

The next case to consider is $\pi_{temp} = T_{IV, 3}(\pi_{sc})$. 
\begin{prop}\label{nontemp3C2}
     Consider the representation 
    \begin{equation*}
        \pi_{x_1, x_2} = L(\Delta_{\rho}[-x_1,-x_1], \Delta_{\rho}[-x_2,-x_2], \Delta_{\rho}[-x_2,-x_2]; \pi_{temp}),
    \end{equation*}
    where $x_1 > x_2 \geq \frac{1}{2}$ and $\pi_{temp} = T_{IV,3}(\pi_{sc})$ for $\alpha \in \mathbb{Z}_{>0}$. 
    \begin{enumerate}[label = (\roman*)]
        \item The representation $\pi_{x_1, x_2}$ is of Arthur type if and only if $1 < x_1 \leq \alpha -1$ and $x_2 = 1$, or when $(x_1, x_2) = (2,1)$. 
        \item The representation $\pi_{x_1, x_2}$ is of critical type when $(x_1, x_2, \alpha) = (2,1,1)$ or $(2,1,2)$. 
        \item Define $\EE_{x_1, x_2}$ in various cases as follows. Then $\pi(\EE_{x_1, x_2}) = \pi_{x_1, x_2}$. When $(x_1, x_2) = (2,1)$, define 
       \[\EE_{2,1}= \scalebox{0.8}{\bordermatrix{ 
         &-2&-1&0&1&2 & \cdots & \alpha -1  \cr
        &\lhd &\lhd & \odot & \rhd &\rhd\cr
        &&\lhd &\odot & \rhd\cr
        &&&\odot &\odot &\odot & \cdots & \odot 
        }}.\]
        When $x_2 = 1, x_1 > 2$, define 
        \begin{flalign*}
            \EE_{x_1, x_2}:= \{([x_1, -x_1]_\rho,\lfloor x_1 \rfloor,\eta),([1,-1]_\rho,1, \eta)^2, ([x_1 -2, 0]_\rho,0,\eta),([\alpha-1, x_1]_\rho,0,(-1)^{x_1}\eta).
        \end{flalign*}
        Here is the associated symbol. 
        \[\EE_{x_1, x_2}= \scalebox{0.8}{\bordermatrix{ 
         &-x_1 & \cdots & -1 & 0& 1& \cdots & x_1 -2 & x_1 -1 & x_1 & \cdots &\alpha -1  \cr
        &\lhd & \cdots & \lhd & \odot & \rhd & \cdots & \rhd &\rhd & \rhd \cr
        &&&\lhd &\odot & \rhd\cr
        &&&\lhd &\odot & \rhd\cr
        &&&&\odot &\odot  & \cdots & \odot &\odot &\odot \cr&&&&&&&&&\odot & \cdots & \odot 
        }}.\]
    \end{enumerate}
\end{prop}
\begin{proof}
The sufficient direction is easy to verify. The necessary direction follows from Lemma \ref{absolutevalue} and \cite[Proposition $11.19$]{HJLLZ24}. The rest follows from definition. \end{proof}

The final case we need to consider in Case $(C)$ is $\pi_{temp} = T_{V,2}^{\pm}(\pi_{sc})$
\begin{prop}
     Consider the representation 
    \begin{equation*}
        \pi_{x_1, x_2}^{\pm} = L(\Delta_{\rho}[-x_1,-x_1], \Delta_{\rho}[-x_2,-x_2], \Delta_{\rho}[-x_2,-x_2]; \pi_{temp}^{\pm}),
    \end{equation*}
    where $x_1 > x_2 \geq \frac{1}{2}$ and $\pi_{temp}^{\pm} = T_{V,2}^{\pm}(\pi_{sc})$ for $\alpha = 0$. 
    \begin{enumerate}[label = (\roman*)]
        \item The representation $\pi_{x_1, x_2}^{\pm}$ is of Arthur type if and only if $(x_1, x_2) = (2,1)$. 
        \item The representation $\pi_{x_1, x_2}^{\pm}$ is of critical type if and only if $(x_1, x_2) = (2,1)$. 
        \item We have $\pi(\EE^{\pm}) = \pi_{2,1}^{\pm}$, where 
        \[\EE_{+}= \scalebox{0.8}{\bordermatrix{ 
         &-2 & -1& 0 & 1 & 2  \cr
        &\lhd &\lhd & \oplus & \rhd &\rhd\cr
        &&\lhd &\oplus &\rhd\cr
        }}, \quad
        \EE_{-}= \scalebox{0.8}{\bordermatrix{ 
         &-2 & -1& 0 & 1 & 2  \cr
        &\lhd &\lhd & \ominus & \rhd &\rhd\cr
        &&\lhd &\ominus &\rhd\cr
        }}.\]
    \end{enumerate}
\end{prop}
\begin{proof} This follows directly from \cite[Proposition $11.20$]{HJLLZ24} and Lemma \ref{absolutevalue}. \end{proof}

This concludes our discussion of Case $(C)$. Now we proceed to Case $(D)$. 
\subsection{\texorpdfstring{Case $(D):\pi = L(\Delta_\rho[-x,-x],\Delta_\rho[-x,-x],\Delta_\rho[-x,-x]; \pi_{temp})$}{}}
In this subsection we consider the situation where all three segments inside the $L$-data of $\pi$ are identical, i.e. when $\pi_{x} = L(\Delta_\rho[-x,-x],\Delta_\rho[-x,-x],\Delta_\rho[-x,-x]; \pi_{temp})$ for $x \geq \frac{1}{2}$ and $\pi_{temp}$ is tempered of corank $1$. As usual, there are three cases to consider. First is $\pi_{temp} = T_{I,1}^{\alpha}(\pi_{sc})$.
\begin{prop}\label{nontemp3D1}
    Consider the representation 
    \begin{equation*}
        \pi_x = L(\Delta_\rho[-x,-x],\Delta_\rho[-x,-x],\Delta_\rho[-x,-x]; \pi_{temp}),
    \end{equation*}
    where $x \geq \frac{1}{2}$ and $\pi_{temp} = T_{I,1}^{\alpha}(\pi_{sc})$ for $\alpha > 0$. 
    \begin{enumerate}[label = (\roman*)]
        \item The representation $\pi_x$ is not of Arthur type for any $x$. 
        \item The representation $\pi_x$ is of critical type when $x \in \{\alpha -1, \alpha, \alpha +1\}$. 
    \end{enumerate}
\end{prop}
\begin{proof} Part $(i)$ follows directly from Lemma \ref{absolutevalue}. The key point is that $\pi_{temp}$ is multiplicity-free. Part $(ii)$ follows from definition. \end{proof}

The next case is $\pi_{temp} = T_{IV,3}(\pi_{sc})$. 
\begin{prop}\label{nontemp3D2}
    Consider the representation 
    \begin{equation*}
        \pi_x = L(\Delta_\rho[-x,-x],\Delta_\rho[-x,-x],\Delta_\rho[-x,-x]; \pi_{temp}),
    \end{equation*}
    where $x \geq \frac{1}{2}$ and $\pi_{temp} = T_{IV,3}(\pi_{sc})$ for $\alpha \in \mathbb{Z}_{>0}$. 
    \begin{enumerate}[label = (\roman*)]
        \item The representation $\pi_{x}$ is of Arthur type if and only if $x = 1$. 
        \item The representation $\pi_x$ is of critical type if and only if $x = \alpha = 1$. 
        \item We have $\pi(\EE) = \pi_1$, where 
        \[\EE= \scalebox{0.8}{\bordermatrix{ 
         &-1&0&1 & \cdots & \alpha -1  \cr
        &\lhd  & \odot & \rhd \cr
        &\lhd  & \odot & \rhd \cr
        &\lhd  & \odot & \rhd \cr
        &&&\odot&\cdots & \odot  
        }}.\]
    \end{enumerate}
\end{prop}
\begin{proof} By Lemma \ref{absolutevalue}, the representation $\pi_{x}$ is of Arthur type only if there exists three copies of $\rho\lvert \cdot \rvert^{x-1}$ inside $|\Omega|(\pi_x^{\rho,-})$. This can only happen when $x = 1$. This proves the necessary direction. The sufficient direction is easy to verify and the rest follows from definition. \end{proof}

The last case we have is $\pi_{temp} = T_{V,2}^{\pm}(\pi_{sc})$. 
\begin{prop}\label{nontemp3D3}
     Consider the representation 
    \begin{equation*}
        \pi_x^{\pm} = L(\Delta_\rho[-x,-x],\Delta_\rho[-x,-x],\Delta_\rho[-x,-x]; \pi_{temp}^{\pm}),
    \end{equation*}
    where $x \geq \frac{1}{2}$ and $\pi_{temp}^{\pm} = T_{V,2}^{\pm}(\pi_{sc})$ for $\alpha = 0$.
    \begin{enumerate}[label = (\roman*)]
        \item The representation $\pi_{x}$ is not of Arthur type for any $x$. 
        \item The representation $\pi_x$ is of critical type if and only if $x = 1$. 
    \end{enumerate}
\end{prop}
\begin{proof} This follows directly from Lemma \ref{absolutevalue}. \end{proof}

This concludes our discussion of Case $(D)$. Now we are ready to move on to Cases $(E), (F), (G)$, where the three segments inside the $L$-data of $\pi$ exhaust the corank of $\pi$ by $4$. 

\subsection{\texorpdfstring{Cases $(E),(F),(G)$ involving supercuspidal representations}{}}
The last three subcases in this section all involve supercuspidal representations. We start with Case $(E)$, where $\pi = L(\Delta_\rho[-x_1, -x_1-1], \Delta_\rho[-x_2, -x_2], \Delta_\rho[-x_3,-x_3]; \pi_{sc})$, where $x_1 \geq x_2 \geq x_3$ by Langlands classification and the good parity condition. 

\begin{prop}\label{nontemp3E}
    Consider the representation 
    \begin{equation*}
        \pi_{x_1, x_2, x_3} =  L(\Delta_\rho[-x_1, -x_1-1], \Delta_\rho[-x_2, -x_2], \Delta_\rho[-x_3,-x_3]; \pi_{sc}),
    \end{equation*}
    where $x_1 \geq x_2 \geq x_3 \geq \frac{1}{2}$. 
    \begin{enumerate}[label = (\roman*)]
        \item The representation $\pi_{x_1, x_2, x_3}$ is of Arthur type if and only if $\frac{1}{2} \leq x_3 < x_2 < x_1 \leq \alpha -1$, or $(x_1, x_2, x_3) = (\frac{3}{2}, \frac{1}{2}, \frac{1}{2})$ or $(\frac{1}{2}, \frac{1}{2}, \frac{1}{2})$. 
        \item The representation $\pi_{x_1, x_2, x_3}$ is of critical type if and only if $(x_1, x_2, x_3)$ is one of the following points: $(\alpha, \alpha, \alpha),(\alpha,\alpha, \alpha -1), (\alpha, \alpha -1, \alpha -1), (\alpha,\alpha -1,\alpha-2), (\alpha -1,\alpha-1, \alpha -1),(\alpha -1, \alpha -1, \alpha -2),(\alpha -1,\alpha -2, \alpha -2),(\alpha -1, \alpha -2,\alpha -3),(\alpha +1, \alpha, \alpha),(\alpha +1, \alpha, \alpha -1),(\alpha+1, \alpha+1, \alpha),(\alpha+2, \alpha+1, \alpha)$. 
         \item Define $\EE_{x_1, x_2, x_3}$ in various cases as follows. Then $\pi(\EE_{x_1, x_2, x_3}) = \pi_{x_1, x_2, x_3}$. When $\frac{1}{2} \leq x_3 < x_2 < x_1 \leq \alpha -1$, and $x_1 = x_2 +1 = x_3 + 2$, define
            \begin{flalign*}
                \EE_{x_1, x_2, x_3}:= &\{([x_1 +1, -x_1]_\rho,\lfloor x_1 \rfloor, \eta),([x_3 -2, \epsilon_\rho]_\rho,0,\eta),([x_1 -1, x_3]_\rho,0,(-1)^{x_3 - \epsilon_\rho}\eta),\\
                &([\alpha -1, x_1+1]_\rho,0,(-1)^{x_1 +1 - \epsilon_\rho}\eta)\}.
            \end{flalign*}
            Here is the associated symbol. 
        
        \[\EE_{x_1, x_2,x_3}= \scalebox{0.8}{\bordermatrix{ 
        &-x_1   & \cdots & \epsilon_\rho & \cdots & x_3 -2&  x_3 -1 & x_3  &\cdots &x_1 -1 &x_1 &x_1 +1& \cdots &\alpha-1 \cr
        &\lhd  & \cdots & \odot & \cdots & \rhd & \rhd & \rhd & \cdots & \rhd & \rhd & \rhd\cr
        &&&\odot & \cdots & \odot \cr
        &&&&&&&\odot&\cdots &\odot  \cr
        &&&&&&&&&& &\odot & \cdots &\odot \cr
        }}.\]
        When $\frac{1}{2} \leq x_3 < x_2 < x_1 \leq \alpha -1$, and $x_1 > x_2 +1  = x_3 + 2$, define
            \begin{flalign*}
                \EE_{x_1, x_2, x_3}:= &\{([x_1 +1, -x_1]_\rho,\lfloor x_1 \rfloor, \eta),
                ([x_2, -x_2]_\rho,\lfloor x_2 \rfloor, \eta), ([x_3 -2, \epsilon_\rho]_\rho,0,\eta),\\
                &([x_1 -2, x_3]_\rho,0,(-1)^{x_3 - \epsilon_\rho}\eta),
                ([\alpha -1, x_1+1]_\rho,0,(-1)^{x_1 +1 - \epsilon_\rho}\eta)\}.
            \end{flalign*}
            Here is the associated symbol for $\EE_{x_1, x_2,x_3}$: 
        
        \[ \scalebox{0.8}{\bordermatrix{ 
        &-x_1   & \cdots &-x_2& \cdots &\epsilon_\rho & \cdots & x_3 -2&  x_3 -1 & x_3 & x_2 &\cdots &x_1 -2 &x_1 -1 & x_1 &x_1 +1& \cdots &\alpha-1 \cr
        &\lhd  & \cdots & \lhd &\cdots & \odot & \cdots & \rhd & \rhd & \rhd &\rhd & \cdots & \rhd & \rhd & \rhd &\rhd\cr
        &&&\lhd & \cdots & \odot & \cdots & \rhd & \rhd & \rhd & \rhd \cr 
        &&&&&\odot  & \cdots & \odot \cr
        &&&&&&&&&\odot&\odot &\cdots &\odot  \cr
        &&&&&&&&&&&&&& &\odot & \cdots &\odot \cr
        }}.\]
        When $\frac{1}{2} \leq x_3 < x_2 < x_1 \leq \alpha -1$, and $x_1 > x_2 +1 > x_3 + 1$, define $\EE_{x_1, x_2, x_3}:=$
            \begin{flalign*}
                 &\{([x_1 +1, -x_1]_\rho,\lfloor x_1 \rfloor, \eta),
                ([x_2, -x_2]_\rho,\lfloor x_2 \rfloor, \eta), ([x_3, -x_3]_\rho, \lfloor x_3 \rfloor, \eta),([x_3 -2, \epsilon_\rho]_\rho,0,\eta),\\
                &([x_2 -2, x_3]_\rho,0,(-1)^{x_3 - \epsilon_\rho}\eta), ([x_1 -2, x_2]_\rho,0,(-1)^{x_2 - \epsilon_\rho}\eta),
                ([\alpha -1, x_1+1]_\rho,0,(-1)^{x_1 +1 - \epsilon_\rho}\eta)\}.
            \end{flalign*}
            Here is the associated symbol for $\EE_{x_1, x_2,x_3}$: 
        
        \[ \scalebox{0.6}{\bordermatrix{ 
        &-x_1   & \cdots &-x_2& \cdots &-x_3 & \cdots &\epsilon_\rho & \cdots & x_3 -2&  x_3 -1 & x_3 & \cdots &x_2 -2 & x_2 -1 & x_2 & \cdots &x_1 -2 &x_1 -1 & x_1 &x_1 +1& \cdots &\alpha-1 \cr
        &\lhd  & \cdots & \lhd &\cdots &\lhd & \cdots& \odot & \cdots & \rhd & \rhd & \rhd &\cdots & \rhd & \rhd & \rhd & \cdots & \rhd & \rhd & \rhd &\rhd\cr
        &&&\lhd & \cdots & \lhd & \cdots & \odot & \cdots & \rhd &\rhd & \rhd &\cdots & \rhd & \rhd & \rhd \cr
        &&&&&\lhd & \cdots & \odot & \cdots & \rhd & \rhd &\rhd \cr
        &&&&&&&\odot & \cdots & \odot \cr
        &&&&&&&&&&&\odot&\cdots &\odot  \cr
        &&&&&&&&&&&&&& &\odot & \cdots &\odot \cr
        &&&&&&&&&&&&&&&&&&&&\odot &\cdots & \odot
        }}.\]
        When $(x_1, x_2, x_3) = (\frac{1}{2}, \frac{1}{2}, \frac{1}{2})$, let 
        \[\EE_{x_1, x_2,x_3}= \scalebox{0.8}{\bordermatrix{ 
         &-\frac{1}{2} & \frac{1}{2}&\frac{3}{2}&\cdots & \alpha -1  \cr
        &\lhd  & \rhd \cr
        &\lhd   & \rhd \cr
        &\lhd  &\odot  & \rhd \cr
        &&&\odot & \cdots &\odot
        }}.\]
         When $(x_1, x_2, x_3) = (\frac{3}{2}, \frac{1}{2}, \frac{1}{2})$, let 
        \[\EE_{x_1, x_2,x_3}= \scalebox{0.8}{\bordermatrix{ 
         &-\frac{3}{2} &-\frac{1}{2} & \frac{1}{2}&\frac{3}{2}&\frac{5}{2}&\cdots & \alpha -1  \cr
         &\lhd & \lhd & \odot & \rhd & \rhd \cr
        &&\lhd  & \rhd \cr
        &&&\odot \cr
        &&&&&\odot & \cdots &\odot
        }}.\]
    \end{enumerate}
\end{prop}
\begin{proof} The sufficient direction can be proved in a similar way as Proposition \ref{nontemp3A1}, which we omit. Now we show the necessary direction. First assume $x_1 = x_2$ and $\pi_{x_1, x_2, x_3}$ is of Arthur type. Then by Definition \ref{psiminus} and Theorem \ref{nontemp red}, there exists an extended multi-segment $\EE$ containing a segment of the form $([x, -(x-1)]_\rho,*,*)$ and another segment of the form $[x-1, -(x-1)]_\rho$ such that $\pi(\EE) = \pi_{sc})$. This cannot happen since the $L$-parameter of $\pi_{sc}$ is multiplicity-free. Thus we may conclude that $x_1 > x_2$. 

Then by \cite[Proposition $9.12$]{HJLLZ24}, the representation 
\begin{equation*}
    \pi_{x_1, x_2, x_3}^{\rho,-} = L(\Delta_\rho[-x_2, -x_2], \Delta_\rho[-x_3, -x_3]; \pi_{sc})
\end{equation*}
is of Arthur type if and only if one of the following holds: 
\begin{itemize}
    \item $\frac{1}{2} \leq x_3 < x_2 \leq \alpha$, 
    \item $(x_2, x_3) = (\alpha +1, \alpha)$,
    \item $(x_2, x_3) = (\alpha, \alpha -1)$,
    \item $x_2 = x_3 = \frac{1}{2}$.
\end{itemize}
Combining these conditions with Lemma \ref{absolutevalue}, we may conclude that $\pi_{x_1, x_2, x_3}$ is of Arthur type only when $\frac{1}{2} \leq x_3 < x_2 < x_1 \leq \alpha -1$ and when $x_2 = x_3 = \frac{1}{2}$. In the second case, we must have $x_1 = \frac{1}{2}$ or $\frac{3}{2}$, since otherwise there exists an extended multi-segment $\EE$ containing a segment of the form $([x_1, -(x_1 -1)]_\rho$ such that $\pi(\EE) = L(\Delta_\rho[-\frac{1}{2}, -\frac{1}{2}], \Delta_\rho[-\frac{1}{2}, -\frac{1}{2}]; \pi_{sc})$ with $x_1 > \frac{3}{2}$, which cannot happen. This proves part $(i)$. The rest follows from definition \end{proof}

  This concludes our discussion of Case $(E)$. Now we move on to Case $(F)$, where 
  \begin{equation*}
      \pi = L(\Delta_\rho[-x_1, -x_1], \Delta_\rho[-x_2, -x_2 -1], \Delta_\rho[-x_3, -x_3] ;\pi_{sc}),
  \end{equation*}
  for $x_1 > x_2 \geq x_3$. 

\begin{prop}\label{nontemp3F}
     Consider the representation 
    \begin{equation*}
        \pi_{x_1, x_2, x_3} = L(\Delta_\rho[-x_1, -x_1], \Delta_\rho[-x_2, -x_2 -1], \Delta_\rho[-x_3, -x_3]; \pi_{sc}),
    \end{equation*}
    where $x_1 > x_2 \geq x_3 \geq \frac{1}{2}$. 
    \begin{enumerate}[label = (\roman*)]
        \item The representation $\pi_{x_1, x_2, x_3}$ is of Arthur type if and only if 
        one of the following holds: 
        \begin{itemize}
            \item $\epsilon_\rho +1 < x_1 \leq \alpha$, $\frac{1}{2} \leq x_3 < x_2 \leq \alpha -2$ and $x_1 - x_2 > 1$.
            \item $(x_1, x_2) = (\alpha+1, \alpha-1)$, and $\frac{1}{2} \leq x_3 < x_2$.
            \item $(x_1, x_2, x_3) = (\frac{3}{2}, \frac{1}{2}, \frac{1}{2})$ and $\alpha \in \frac{1}{2} + \mathbb{Z}_{>0}$. 
        \end{itemize}
        \item The representation $\pi_{x_1, x_2, x_3}$ is of critical type if and only if $(x_1, x_2,x_3) = (\alpha+3, \alpha+1, \alpha),(\alpha+2, \alpha+1, \alpha),(\alpha+2, \alpha, \alpha),(\alpha+2, \alpha, \alpha-1),(\alpha+1, \alpha, \alpha),(\alpha+1, \alpha, \alpha-1),(\alpha+1, \alpha-1, \alpha-1),(\alpha+1, \alpha-1,\alpha-2),(\alpha, \alpha-1, \alpha-1),(\alpha, \alpha-1, \alpha-2),(\alpha, \alpha-2,\alpha-2),(\alpha, \alpha-2,\alpha-3)$. 
        \item Define $\EE_{x_1, x_2, x_3}$ in various cases as follows. Then $\pi(\EE_{x_1, x_2, x_3}) = \pi_{x_1, x_2, x_3}$. When $\epsilon_\rho +1 < x_1 \leq \alpha, x_1 - x_2 = 2$, and $\frac{1}{2}\leq x_3 = x_2-1 \leq \alpha -3$, define
        \begin{flalign*}
                \EE_{x_1, x_2, x_3}:= &\{([x_1, -x_1]_\rho,\lfloor x_1 \rfloor, \eta),
                ([x_2, -x_2]_\rho,\lfloor x_2 \rfloor, -\eta),
                ([x_3 -2, \epsilon_\rho]_\rho,0,\eta),\\
                &([x_3, x_3]_\rho,0,(-1)^{x_3 - \epsilon_\rho}\eta),([\alpha-1, x_2 +1]_\rho, 0, (-1)^{x_2 +1 - \epsilon_\rho}\eta)\}.
            \end{flalign*}
            Here is the associated symbol. 
        
        \[\EE_{x_1, x_2,x_3}= \scalebox{0.8}{\bordermatrix{ 
        &-x_1   & \cdots &-x_2& \cdots &\epsilon_\rho & \cdots & x_3 -2&  x_3 -1 & x_3 & x_2 &x_1 -1 & x_1 & \cdots &\alpha-1 \cr
        &\lhd  & \cdots & \lhd &\cdots & \odot & \cdots & \rhd & \rhd & \rhd &\rhd &\rhd &\rhd \cr
        &&&\lhd & \cdots & \odot & \cdots & \rhd & \rhd & \rhd & \rhd \cr 
        &&&&&\odot & \cdots &\odot \cr
        &&&&&&&&&\odot \cr
        &&&&&&&&&&&\odot&\odot &\cdots &\odot  \cr
        }}.\]
        When $\epsilon_\rho +1 < x_1 \leq \alpha, x_1 - x_2 > 2$, and $\frac{1}{2}\leq x_3 = x_2-1 \leq \alpha -3$, define 
        \begin{flalign*}
                \EE_{x_1, x_2, x_3}:= &\{([x_1, -x_1]_\rho,\lfloor x_1 \rfloor, \eta),
                ([x_2+1, -x_2]_\rho,\lfloor x_2 \rfloor, -\eta),
                ([x_3 -2, \epsilon_\rho]_\rho,0,-\eta),\\
                &([x_3, x_3]_\rho,0,(-1)^{x_3 - \epsilon_\rho}\eta),([x_1 -2, x_2 +1]_\rho, 0, (-1)^{x_2 - \epsilon_\rho}\eta),([\alpha -1, x_1]_\rho, 0, (-1)^{x_1 - \epsilon_\rho} \eta)\}.
            \end{flalign*}
            Here is the associated symbol for $\EE_{x_1, x_2,x_3}$: 
        \[ \scalebox{0.8}{\bordermatrix{ 
        &-x_1   & \cdots &-x_2& \cdots &\epsilon_\rho & \cdots & x_3 -2&  x_3 -1 & x_3 & x_2 &x_2 +1 & \cdots &x_1 -2 & x_1 -1& x_1 & \cdots &\alpha-1 \cr
        &\lhd  & \cdots & \lhd &\cdots & \odot & \cdots & \rhd & \rhd & \rhd &\rhd &\rhd &\cdots & \rhd & \rhd & \rhd  \cr
        &&&\lhd & \cdots & \odot & \cdots & \rhd & \rhd & \rhd & \rhd &\rhd\cr 
        &&&&&\odot & \cdots & \odot \cr
        &&&&&&&&&\odot \cr
        &&&&&&&&&&&\odot& \cdots & \odot   \cr
        &&&&&&&&&&&&&&&\odot & \cdots & \odot 
        }}.\]
         When $\epsilon_\rho +1 < x_1 \leq \alpha, x_1 - x_2 > 2$, and $\frac{1}{2}\leq x_3 < x_2-1 \leq \alpha -3$, define 
        \begin{flalign*}
                \EE_{x_1, x_2, x_3}:= &\{([x_1, -x_1]_\rho,\lfloor x_1 \rfloor, \eta),
                ([x_2+1, -x_2]_\rho,\lfloor x_2 \rfloor, -\eta), ([x_2-2, -x_3]_\rho,0,(-1)^{x_3 - \epsilon_\rho}\eta), \\
                &([x_3 -2, \epsilon_\rho]_\rho,0,-\eta),
                 ([x_2-2, x_3]_\rho,0,(-1)^{x_3 - \epsilon_\rho}\eta),([x_1 -2, x_2 +1]_\rho, 0, (-1)^{x_2 - \epsilon_\rho}\eta), \\
                 &([\alpha -1, x_1]_\rho, 0, (-1)^{x_1 - \epsilon_\rho} \eta)\}.
            \end{flalign*}
            Here is the associated symbol for $\EE_{x_1, x_2,x_3}$: 
        \[ \scalebox{0.7}{\bordermatrix{ 
        &-x_1   & \cdots &-x_2& \cdots &\epsilon_\rho & \cdots & x_3 -2&  x_3 -1 & x_3 &\cdots &x_2 -2 & x_2 -1 & x_2 &x_2 +1 & \cdots &x_1 -2 & x_1 -1& x_1 & \cdots &\alpha-1 \cr
        &\lhd  & \cdots & \lhd &\cdots & \odot & \cdots & \rhd & \rhd & \rhd &\cdots &\rhd &\rhd &\rhd &\rhd &\cdots & \rhd & \rhd & \rhd  \cr
        &&&\lhd & \cdots & \odot & \cdots & \rhd & \rhd & \rhd & \rhd &\rhd\cr 
        &&&&&\odot & \cdots & \odot \cr
        &&&&&&&&&\odot&\cdots & \odot  \cr
        &&&&&&&&&&&&&&\odot& \cdots & \odot   \cr
        &&&&&&&&&&&&&&&&&&\odot & \cdots & \odot 
        }}.\]
        When $(x_1, x_2, x_3) = (\frac{3}{2}, \frac{1}{2}, \frac{1}{2})$ and $\alpha \in \frac{1}{2} + \mathbb{Z}_{>0}$, let 
        \[\EE_{x_1, x_2,x_3}= \scalebox{0.8}{\bordermatrix{ 
         &-\frac{3}{2} & -\frac{1}{2} & \frac{1}{2}&\frac{3}{2}&\cdots & \alpha -1  \cr
        &\lhd  &\lhd &\rhd & \rhd \cr
        &&\lhd   &\odot & \rhd \cr
        &&&&\odot & \cdots &\odot
        }}.\]
        
    \end{enumerate}
\end{prop}
\begin{proof} The sufficient direction can be proved in a similar way as Proposition \ref{nontemp3E}, which we omit. Now we show the necessary direction. 
First, for $x_2 > \frac{1}{2}$, we assume that $x_1 = x_2 +1$, and $\pi_{x_1, x_2, x_3}$ is of Arthur type. Then by Definition \ref{psiminus} and Theorem \ref{nontemp red}, there exists an extended multi-segment $\EE$ that contains a segment of the form $\{[x_1 -1, -(x_1-1)]_\rho,*,*\}$ and another segment of the form $([x_2, -(x_2 -1)]_\rho,*,*)$, such that $\pi(\EE) = \pi_{sc}$. This gives a contradiction since the $L$-parameter of $\pi_{sc}$ is multiplicity-free. Thus, we may conclude that $x_1 - x_2 > 1$. 

Now by \cite[Proposition $11.21$ and $11.22$]{HJLLZ24}, the representation 
\begin{flalign*}
    \pi_{x_1, x_2, x_3}^{\rho,-} = L(\Delta_{\rho}[-x_2, -x_2 -1], \Delta_\rho[-x_3, -x_3]; \pi_{sc})
\end{flalign*}
is of Arthur type if and only if one of the following holds: 
\begin{itemize}
    \item $\frac{1}{2} \leq x_3 < x_2 \leq \alpha -1$,
    \item $x_2 = x_3 = \frac{1}{2}$ and $\alpha > \frac{1}{2}$.
\end{itemize}
In order for $\pi_{x_1, x_2, x_3}$ to be of Arthur type, we must also have $\epsilon_\rho +1 \leq x_1 \leq \alpha$ for the first case, or $(x_1,x_2) = (\alpha+1, \alpha-1)$ by Lemma \ref{absolutevalue}. In the second case, if $x_1 \geq \frac{5}{2}$, then $\pi_{x_1, x_2, x_3}$ is of Arthur type only if there exists an extended multi-segment $\EE$ containing the segment $([x_1 -1, -(x_1-1)]_\rho,*,*)$ such that $\pi(\EE) = L(\Delta_\rho[-\frac{1}{2}, -\frac{3}{2}], \Delta_\rho[-\frac{1}{2}, -\frac{1}{2}])$. You can show that this cannot be the case by exhausting the possible Arthur packets the representation can appear in, using the algorithm described in \cite{HLL22}. This proves the necessary direction. The rest follows from definition. \end{proof}

This takes care of Case $(F)$. Now we tackle the final case, Case $(G)$, where 
\begin{equation*}
    \pi = L(\Delta_\rho[-x_1,-x_1], \Delta_\rho[-x_2, -x_2], \Delta_\rho[-x_3, -x_3-1]; \pi_{sc}).
\end{equation*}
\begin{prop}
    Consider the representation 
    \begin{equation*}
    \pi_{x_1,x_2, x_3} = L(\Delta_\rho[-x_1,-x_1], \Delta_\rho[-x_2, -x_2], \Delta_\rho[-x_3, -x_3-1]; \pi_{sc}),
\end{equation*}
where $x_1 \geq x_2 > x_3 \geq 0$. 
\begin{enumerate}[label  =(\roman*)]
    \item The representation $\pi_{x_1, x_2, x_3}$ is of Arthur type if and only if one of the following holds: 
    \begin{itemize}
        \item $1 \leq x_3 +1 < x_2 < x_1 \leq \alpha$, 
        \item $(x_1, x_2) = (\alpha +1, \alpha)$ and $0 \leq x_3 < \alpha -1$,
        \item $(x_1, x_2, x_3) = (\alpha +2, \alpha +1, \alpha -1)$,
        \item $(x_1, x_2, x_3, \alpha) = (2,1,0,1)$. 
    \end{itemize}
    \item The representation $\pi_{x_1, x_2, x_3}$ is of critical type if and only if $(x_1, x_2, x_3)$ is one of the following points: $(\alpha +1, \alpha +1, \alpha),(\alpha +2, \alpha +1, \alpha),(\alpha +2, \alpha +2, \alpha),(\alpha +3, \alpha +2, \alpha), (\alpha, \alpha, \alpha -1),(\alpha +1, \alpha, \alpha -1),(\alpha+1, \alpha +1, \alpha -1),(\alpha +2, \alpha +1, \alpha-1)
    ,(\alpha, \alpha, \alpha -2),(\alpha, \alpha-1, \alpha -2),(\alpha, \alpha -1,\alpha -3),(\alpha +1, \alpha, \alpha -2)$.  
    \item Define $\EE_{x_1, x_2, x_3}$ in various cases as follows, then we have $\pi(\EE_{x_1, x_2, x_3}) = \pi_{x_1, x_2, x_3}$. When $1 \leq x_3 +2 =  x_2 = x_1-1 \leq \alpha +1$, define
    \begin{flalign*}
        \EE_{x_1, x_2, x_3} := &\{([x_1, -x_1]_\rho, \lfloor x_1 \rfloor, \eta), ([x_3, -x_3]_\rho, \lfloor x_3 \rfloor, -\eta), ([x_3 -2, \epsilon_\rho]_\rho,0,\eta)\\
        &([\alpha -1, x_2]_\rho, 0, (-1)^{x_2 - \epsilon_\rho}\eta).
    \end{flalign*}
    Here is the associated symbol. 
    \[\EE_{x_1, x_2,x_3}= \scalebox{0.8}{\bordermatrix{ 
        &-x_1   & -x_2 & -x_3& \cdots &\epsilon_\rho & \cdots & x_3 -2&  x_3 -1 & x_3 & x_2 & x_1 & \cdots &\alpha-1 \cr
        &\lhd  & \lhd & \lhd &\cdots & \odot & \cdots & \rhd & \rhd & \rhd &\rhd &\rhd \cr
        &&&\lhd & \cdots & \odot & \cdots & \rhd & \rhd & \rhd \cr 
        &&&&&\odot & \cdots &\odot \cr
        &&&&&&&&&&\odot&\odot &\cdots &\odot  \cr
        }}.\]
         When $1 \leq x_3 +2 <  x_2  = x_1-1 \leq \alpha $, define
    \begin{flalign*}
        \EE_{x_1, x_2, x_3} := &\{([x_1, -x_1]_\rho, \lfloor x_1 \rfloor, \eta), ([x_3+1, -x_3]_\rho, \lfloor x_3 \rfloor, -\eta), ([x_3 -2, \epsilon_\rho]_\rho,0,\eta)\\
        &([x_2-2, x_3+1]_\rho, 0, (-1)^{x_3 +1 - \epsilon_\rho}\eta), ([\alpha -1, x_2]_\rho, 0, (-1)^{x_2 - \epsilon_\rho}\eta)\}.
    \end{flalign*}
    Here is the associated symbol for $\EE_{x_1, x_2,x_3}$: 
    \[\scalebox{0.8}{\bordermatrix{ 
        &-x_1   & \cdots &-x_3& \cdots &\epsilon_\rho & \cdots & x_3 -2&  x_3 -1 & x_3 & x_3+1 & \cdots &x_2 -2 & x_2 -1 & x_2 & x_1 & \cdots &\alpha-1 \cr
        &\lhd  & \cdots & \lhd &\cdots & \odot & \cdots & \rhd & \rhd & \rhd &\rhd &\rhd &\cdots & \rhd & \rhd & \rhd  \cr
        &&&\lhd & \cdots & \odot & \cdots & \rhd & \rhd & \rhd & \rhd \cr 
        &&&&&\odot & \cdots & \odot \cr
        &&&&&&&&&&\odot &\cdots & \odot\cr
        &&&&&&&&&&&&&&\odot&\odot & \cdots & \odot 
        }}.\]
    Finally, when $1 \leq x_3 +2 < x_2 < x_1 -1 \leq \alpha$, define $\EE_{x_1, x_2, x_3} :=$
    \begin{flalign*}
         &\{([x_1, -x_1]_\rho, \lfloor x_1 \rfloor, \eta), ([x_2, -x_2]_\rho, \lfloor x_2 \rfloor, -\eta), ([x_3+1, -x_3]_\rho, \lfloor x_3 \rfloor, -\eta), ([x_3 -2, \epsilon_\rho]_\rho,0,\eta)\\
        &([x_2-2, x_3+1]_\rho, 0, (-1)^{x_3 +1 - \epsilon_\rho}\eta), ([x_1 -2, x_2]_\rho, 0, (-1)^{x_2 +1 - \epsilon_\rho}\eta), ([\alpha-1, x_1]_\rho, 0, (-1)^{x_1 - \epsilon_\rho} \eta).
    \end{flalign*}
    Here is the associated symbol for $\EE_{x_1, x_2,x_3}$:  
    \[ \scalebox{0.6}{\bordermatrix{ 
        &-x_1   &\cdots & -x_2 & \cdots &-x_3& \cdots &\epsilon_\rho & \cdots & x_3 -2&  x_3 -1 & x_3 & x_3+1 & \cdots &x_2 -2 & x_2 -1 & x_2 &   \cdots &x_1 -2 & x_1 -1 & x_1 & \cdots &\alpha-1 \cr
        &\lhd  & \cdots & \lhd &\cdots & \lhd & \cdots & \odot & \cdots & \rhd & \rhd & \rhd &\rhd &\cdots & \rhd & \rhd & \rhd  & \cdots & \rhd &\rhd &\rhd \cr
        &&&\lhd & \cdots &\lhd & \cdots& \odot & \cdots & \rhd & \rhd & \rhd & \rhd  &\cdots & \rhd &\rhd &\rhd\cr 
        &&&&&\lhd & \cdots &\odot & \cdots & \rhd &\rhd&\rhd\cr
        &&&&&&&\odot & \cdots & \odot \cr
        &&&&&&&&&&&&\odot &\cdots & \odot\cr
        &&&&&&&&&&&&&&&&\odot & \cdots & \odot \cr
        &&&&&&&&&&&&&&&&&&&&\odot & \cdots &\odot
        }}.\]
\end{enumerate}
\end{prop}
\begin{proof} The sufficient direction can be proven in a similar way as Proposition \ref{nontemp3E}, which we omit. Now we prove the necessary direction. 

Suppose $\pi_{x_1, x_2, x_3}$ is of Arthur type, then we must have $x_1 > x_2$. Otherwise, by Definition \ref{psiminus} and Theorem \ref{nontemp red}, there exists an extended multi-segment $\EE$ containing two copies of $([x_1 -1, -(x_1-1)]_\rho$ such that $\pi(\EE)$ such that $\pi(\EE) = \pi_{sc}$. This is impossible since the $L$-parameter of $\pi_{sc}$ is multiplicity-free. 

Next, by \cite[Proposition $11.23$ and $11.24$]{HJLLZ24}, the representation
\begin{equation*}
    \pi_{x_1, x_2, x_3}^{\rho,-} = L(\Delta_\rho[-x_2, -x_2], \Delta_\rho[-x_3, -x_3 -1])
\end{equation*}
is of Arthur type if and only if $1 \leq x_3 +1 < x_2 \leq \alpha$ or $(x_2, x_3) = (\alpha +1, \alpha -1)$. Combining this with Lemma \ref{absolutevalue}, we see that $\pi_{x_1, x_2, x_3}$ is of Arthur type only if one of the following holds: 
\begin{itemize}
    \item $1 < x_3 +1 < x_2 < x_1 \leq \alpha$,
    \item $(x_1, x_2) = (\alpha +1, \alpha)$ and $0 \leq x_3 < \alpha -1$,
    \item $(x_1, x_2, x_3) = (\alpha +2, \alpha +1, \alpha -1)$,
    \item $(x_1, x_2, x_3) = (2,1,0)$ and $\alpha = 1$.
\end{itemize}
This gives the desired condition and proves the necessary direction. The rest follows from definition. \end{proof}
This concludes our discussion of Case $(F)$ and the classification of all the corank 4 non-tempered representations $\pi$ of good parity with $f(\pi) = 3$. 

\subsection{The case \texorpdfstring{$f(\pi) = 4$}{}}
We move onto the final case in this section, where there are $4$ segments in the $L$-data of $\pi$. This situation is a lot easier to deal with compared to the previous cases, since the number of segments match the corank of $\pi$, we must have 
\begin{equation*}
    \pi_{x_1, x_2, x_3, x_4} = L(\Delta_{\rho}[-x_1, -x_1], \Delta_{\rho}[-x_2, -x_2], \Delta_{\rho}[-x_3, -x_3], \Delta_{\rho}[-x_4, -x_4]; \pi_{sc})
\end{equation*}
There is a way to directly characterize which one of  such representations lie in an Arthur packet, introduced as Condition ($\mathcal{A}$) in \cite{HJLLZ24}. We reference their result below. In general, we consider representations of the form 
\begin{align}\label{eq f=r}
     \pi= L(\Delta_{\rho}[-x_1,-x_1]^{m_{x_1}},\dots, \Delta_{\rho}[-x_f,-x_f]^{m_{x_f}}; \pi_{sc}),
\end{align} 
where $\half{1} \leq x_f <  \cdots< x_1$. Here, the $m_{x_i}$ denotes the multiplicity of a given segment inside the $L$-data. For continuity, we set $m_x = 0$ for $x \in \frac{1}{2}\mathbb{Z}$ if $x \neq x_i$ for any $1 \leq i \leq f$. 
\begin{defn}[{\cite[Definition $9.1$]{HJLLZ24}}]\label{conditionA}
    Suppose the $L$-data of $\pi$ is of the form \eqref{eq f=r}. We say the $L$-data of $\pi$ satisfies condition $(\mathcal{A})$ if the following holds. 
\begin{enumerate}
    \item [$\oldbullet$]If $\alpha \geq 1$, then for any $1 \leq  x \leq \alpha$, we have
    $\left\lfloor \half{m_x}\right\rfloor \leq \left\lfloor \half{m_{x-1}}\right\rfloor $.
    \item [$\oldbullet$] For any $x > \alpha$, we have $m_x \leq m_{x-1}$.
\end{enumerate}
\end{defn}

The characterization is given as follows. 
\begin{thm}[{\cite[Theorem $9.2$]{HJLLZ24}}]\label{thm f=r}
A representation $\pi$ of the form \eqref{eq f=r} is of Arthur type (hence unitary) if and only if its $L$-data satisfies condition $(\mathcal{A})$.
\end{thm}
This theorem handles all representations $\pi$ where $f(\pi) = \text{corank}(\pi)$. In particular, we can apply it to the case where $\text{corank}(\pi) = 4$.
\begin{prop}\label{nontemp4}
    Let 
\begin{equation*}
    \pi_{x_1, x_2, x_3, x_4} = L(\Delta_{\rho}[-x_1, -x_1],\Delta_{\rho}[-x_2, -x_2],\Delta_{\rho}[-x_3, -x_3],\Delta_{\rho}[-x_4, -x_4]; \pi_{sc})
\end{equation*}
with $x_1 \geq x_2 \geq x_3 \geq x_4 > 0$ be a representation of good parity. 
\begin{enumerate}
    \item When $x_1 > x_2 > x_3 > x_4$, $\pi_{x_1, x_2, x_3, x_4}$ is of Arthur type if and only if one of the following holds:
\begin{itemize}
    \item $x_1 \leq \alpha$,
    \item $(x_1, x_2) = (\alpha +1, \alpha)$,
    \item $(x_1, x_2, x_3) = (\alpha +2, \alpha+1, \alpha)$,
    \item $(x_1, x_2, x_3, x_4) = (\alpha +3, \alpha +2, \alpha +1, \alpha)$.
\end{itemize}
It's of critical type when $(x_1, x_2, x_3, x_4) = (\alpha +3, \alpha +2, \alpha +1, \alpha), (\alpha +2, \alpha +1, \alpha , \alpha-1), (\alpha +1, \alpha , \alpha -1, \alpha-2), (\alpha , \alpha -1, \alpha -2, \alpha -3)$. 
    \item When $x_1 = x_2 > x_3 > x_4$, $\pi_{x_1, x_2, x_3, x_4}$ is not of Arthur type. It is of critical type when $(x_1, x_2, x_3, x_4) =(\alpha +2, \alpha +2, \alpha +1, \alpha), (\alpha+1, \alpha +1, \alpha, \alpha -1),(\alpha, \alpha, \alpha -1, \alpha -2)$.
    \item When $x_1 > x_2 = x_3 > x_4$, $\pi_{x_1, x_2, x_3, x_4}$ is not of Arthur type. It is of critical type if and only if $(x_1, x_2, x_3, x_4) = (\alpha +2, \alpha +1, \alpha +1, \alpha), (\alpha +1, \alpha , \alpha, \alpha -1), (\alpha, \alpha -1, \alpha -1, \alpha -2)$. 
    \item When $x_1 > x_2 > x_3 = x_4$, $\pi_{x_1, x_2, x_3, x_4}$ is of Arthur type in the following cases: 
    \begin{itemize}
        \item $\alpha > 1, \frac{1}{2} = x_3 = x_4 < x_2 < x_1 \leq \alpha+1$,
        \item $\alpha = \frac{1}{2}, (x_1, x_2, x_3, x_4) = (\frac{5}{2}, \frac{3}{2}, \frac{1}{2}, \frac{1}{2})$.
    \end{itemize}
    It is of critical type if and only if $(x_1, x_2, x_3, x_4) = (\alpha +2, \alpha +1, \alpha, \alpha), (\alpha +1, \alpha , \alpha -1, \alpha -1), (\alpha, \alpha -1, \alpha -2, \alpha -2)$. 
    \item When $x_1 = x_2 = x_3 > x_4$, $\pi_{x_1, x_2, x_3, x_4}$ is not of Arthur type. It is of critical type when $(x_1, x_2, x_3, x_4) = (\alpha, \alpha, \alpha, \alpha -1)$ or $(\alpha +1, \alpha +1, \alpha +1, \alpha)$. 
    \item When $x_1 = x_2 > x_3 = x_4$, $\pi_{x_1, x_2, x_3, x_4}$ is not of Arthur type. It is of critical type when $(x_1, x_2, x_3, x_4) = (\alpha, \alpha, \alpha -1, \alpha -1)$ or $(\alpha +1, \alpha +1, \alpha , \alpha)$. 
    \item When $x_1 > x_2 = x_3 = x_4$, $\pi_{x_1, x_2, x_3, x_4}$ is not of Arthur type. It is of critical type when $(x_1, x_2, x_3, x_4) = (\alpha, \alpha -1, \alpha -1, \alpha -1),(\alpha +1, \alpha, \alpha, \alpha)$. 
    \item When $x_1 = x_2 = x_3 = x_4$, $\pi_{x_1, x_2, x_3, x_4}$ is not of Arthur type. It is of critical type when $(x_1, x_2, x_3, x_4) = (\alpha, \alpha, \alpha, \alpha)$. 
    \end{enumerate}
\end{prop}
\begin{proof} This follows directly from Theorem \ref{thm f=r}. \end{proof}

With this, we are done with the characterization of all corank 4 non-tempered representations that are of Arthur type.  
Towards the full corank $4$ unitary dual and the corresponding Conjecture \ref{unitary dual conjecture}, we need a complete list of all critical type representations of corank $4$ that are of Arthur type. With results in \S \ref{classnontempcorank4,1} to \S \ref{classnontempcorank4,34}, we only need to classify all corank 4 tempered representations of good parity, which we'll do in the next section.

\section{Classification of corank 4 tempered representations of good parity} \label{classtempcorank4}

In this section, we classify all corank 4 tempered representations of good parity, following a similar procedure as in Section \ref{classtempcorank3}. Combining results in the previous sections, this gives us the full list of good parity representations of corank $4$.
In particular, we can identify independently within the admissible dual which representations are of Arthur type and which representations are of critical type. 
By Theorem \ref{thm red from nu to gp}, this allows us to construct the full Arthur dual of corank $4$. This will form the basis of our main results later. 

Again, we fix some self-dual $\rho \in \mathcal{C}$. Let $\alpha = \alpha_{\rho, \sigma}$ and suppose that $\pi \in \Pi_{A,gp}(G_n)$ is tempered of corank $4$. Then there are four cases to consider. 

\begin{enumerate}[label = (\Alph*)]
    \item $\pi \hookrightarrow \rho\lvert \cdot \rvert^{x_1} \rtimes \pi_{temp}$, where $\pi_{temp}$ is tempered of corank $3$. Then $\pi$ is of the form $T_{I,1}^{x}(\pi_{temp}), T_{IV,3}(\pi_{temp})$ or $T_{V,2}^{\pm}(\pi_{temp})$, 
    \item $\pi \hookrightarrow \rho\lvert \cdot \rvert^{x_1} \rtimes \lvert \cdot \rvert^{x_2} \rtimes \pi_{temp}$, where $\pi_{temp}$ is tempered of corank $2$. Then $\pi$ is of the form $T_{I, 2}^{x}(\pi_{temp}), T_{II, 3}^{x}(\pi_{temp}), T_{III, 2}^{\frac{1}{2}}(\pi_{temp}), T_{IV, 5}(\pi_{temp})$ or $T_{V, 4}^{\pm}(\pi_{temp})$, 
    \item $\pi \hookrightarrow \rho\lvert \cdot \rvert^{x_1} \rtimes \lvert \cdot \rvert^{x_2} \times \lvert \cdot \rvert^{x_3} \rtimes \pi_{temp}$, where $\pi_{temp}$ is tempered of corank $1$. Then $\pi$ is of the form $T_{I,3}^{x}(\pi_{temp}), T_{III,2}^{1}(\pi_{temp}), T_{IV,5}(\pi_{temp}), T_{V,6}^{\pm}(\pi_{temp})$, 
    \item $\pi \hookrightarrow \rho\lvert \cdot \rvert^{x_1} \rtimes \lvert \cdot \rvert^{x_2} \times \lvert \cdot \rvert^{x_3} \times \lvert \cdot \rvert^{x_4} \rtimes \pi_{sc}$, where $\pi_{sc}$ is supercuspidal. Then $\pi$ is of the form $T_{I,4}^{x}(\pi_{sc}), T_{II,5}^{x}(\pi_{sc}), T_{III,2}^{\frac{3}{2}}(\pi_{sc}), T_{III, 4}^{\frac{1}{2}}(\pi_{sc}), T_{IV,9}(\pi_{sc}), T_{V,8}^{\pm}(\pi_{sc})$. 
\end{enumerate}
\subsection{\texorpdfstring{Case $(A): \pi \hookrightarrow \rho\lvert \cdot \rvert^{x_1} \rtimes \pi_{temp}$, where $\pi_{temp}$ is tempered of corank $3$}{}}

We begin with Case $(A)$. In Section \ref{classtempcorank3}, we classified all good parity tempered representations of corank $3$. These are given below: 
\begin{flalign*}
    & T_{I,1}^{\alpha +2}(T_{I,1}^{\alpha+1}(T_{I,1}^{\alpha}(\pi_{sc}))), T_{I,1}^{\alpha-1}(T_{I,1}^{\alpha+1}(T_{I,1}^{\alpha}(\pi_{sc}))), T_{IV,3}(T_{I,1}^{\alpha+1}(T_{I,1}^{\alpha}(\pi_{sc}))), 
    T_{V,2}^{\pm}(T_{I,1}^{2}(T_{I,1}^{1}(\pi_{sc}))), \\
    &T_{I,1}^{\alpha -2}(T_{I,1}^{\alpha -1}(T_{I,1}^{\alpha}(\pi_{sc}))), T_{IV,3}(T_{I,1}^{\alpha -1}(T_{I,1}^{\alpha}(\pi_{sc}))), T_{V,2}^{\pm}(T_{I,1}^{1}(T_{I,1}^{2}(\pi_{sc})))), T_{I,1}^{2}(T_{V,2}^{\pm}(T_{I,1}^{1}(\pi_{sc}))), \\ 
    &T_{I,1}^{2}(T_{I,1}^{1}(T_{V,2}^{\pm}(\pi_{sc}))) 
    T_{IV,3}(T_{I,1}^{1}(T_{V,2}^{\pm}(\pi_{sc}))), T_{I,1}^{\frac{3}{2}}(T_{I,2}^{\frac{1}{2}}(\pi_{sc})), T_{I,1}^{\alpha}(T_{II,3}^{\frac{1}{2}}(\pi_{sc})), T_{I,1}^{\frac{3}{2}}(T_{III,2}^{\frac{1}{2}}(\pi_{sc})), \\ 
    & T_{I,1}^{\alpha}(T_{IV,5}(\pi_{sc})), 
    T_{I,1}^{1}(T_{V,4}^{\pm}(\pi_{sc})), 
    T_{V,4}^{\pm}(T_{I,1}^{1}(\pi_{sc})), T_{I,2}^{1}(T_{IV,3}(\pi_{sc})), T_{II,3}^{1}(T_{IV,3}(\pi_{sc})), T_{I,2}^{1}(T_{V,2}^{\pm}(\pi_{sc})), \\
    &T_{III,2}^{1}(\pi_{sc}), 
     T_{IV,7}(\pi_{sc}), T_{V,6}^{\pm}(\pi_{sc}).
\end{flalign*}

There are $22$ possible tempered representations of corank $3$, so there are $66$ cases to consider in Case $(A)$. We begin with $\pi_{temp} = T_{I,1}^{\alpha +2}(T_{I,1}^{\alpha+1}(T_{I,1}^{\alpha}(\pi_{sc})))$.
\begin{prop}\label{temp4A,1}
    Let $\pi_{temp} = T_{I,1}^{\alpha +2}(T_{I,1}^{\alpha+1}(T_{I,1}^{\alpha}(\pi_{sc})))$ for $\alpha > 0$. 
    \begin{enumerate}[label = (\roman*)]
        \item The representation $T_{I,1}^{x}(\pi_{temp})$ is well-defined if and only if $x = \alpha +3$ for $\alpha > 0$, or $x = \alpha -1$ for $\alpha > 1$. 
        \item The representation $T_{I,1}^{x}(\pi_{temp})$ is of critical type when $x = \alpha +3$ for $\alpha > 0$, or $x = \alpha -1$ for $\alpha > 1$. 
        \item Let $x \in \{\alpha +3, \alpha-1\}$, and define
        \begin{equation*}
            \mathcal{E}_{\alpha +3} := \{([\alpha -2, \epsilon_\rho]_\rho,0,\eta),([\alpha +3, \alpha +3]_\rho, 0, (-1)^{\alpha -1 - \epsilon_\rho}\eta)\}, 
        \end{equation*}
        \begin{equation*}
            \mathcal{E}_{\alpha -1} := \{([\alpha -3, \epsilon_\rho]_\rho,0,\eta),([\alpha -1, \alpha -1]_\rho, 0, (-1)^{\alpha  - \epsilon_\rho}\eta), ([\alpha +1, \alpha +1]_\rho, 0, (-1)^{\alpha  - \epsilon_\rho}\eta)\}.
        \end{equation*}
        Then we have $\pi(\mathcal{E}_x) = T_{I,1}^{x}(\pi_{temp})$. Here are the associated symbols. .
        \[\EE_{\alpha+3}= \scalebox{0.8}{\bordermatrix{
  & \epsilon_{\rho} & \cdots & \alpha-2 &\alpha-1& \alpha & \alpha+1 & \alpha+2 &\alpha +3\cr
  & \odot & \cdots & \odot &&&&\cr
  &&&&&&&&\odot
}},\]
\[\EE_{\alpha-1}= \scalebox{0.8}{\bordermatrix{
  & \epsilon_{\rho} & \cdots & \alpha-3 &\alpha-2& \alpha-1 & \alpha & \alpha+1&\alpha +2\cr
  & \odot & \cdots & \odot &&&&\cr
  &&&&&\odot&&\cr
  &&&&&&&&\odot
}}.\]
    \end{enumerate}
\end{prop}
\begin{proof}
    Part $(i)$ follows from Remark \ref{rmk well-defined for Temp} and definition of $T_{I,1}^{x}(\pi_{temp})$. Parts $(ii)$ and $(iii)$ follow from definition. 
\end{proof}

The proofs of Propositions \ref{temp4A,2} to \ref{temp4A,32} below are similar to that of Proposition \ref{temp4A,1}, which we mostly omit. 

\begin{prop}\label{temp4A,2}
Let $\pi_{temp} = T_{I,1}^{\alpha +2}(T_{I,1}^{\alpha+1}(T_{I,1}^{\alpha}(\pi_{sc})))$ for $\alpha > 0$. 
\begin{enumerate}[label = (\roman*)]
    \item The representation $T_{IV,3}(\pi_{temp})$ is well-defined if and only if $x \in \mathbb{Z}_{>1}$. 
    \item The representation $T_{IV,3}(\pi_{temp})$ is of critical type if and only if $\alpha = 1$. 
    \item Define
    \begin{equation*}
            \mathcal{E} := \{([0,0]_\rho,0,\eta)^2,([\alpha -2, 0]_\rho,0,\eta), ([\alpha +2, \alpha +2]_\rho, 0, (-1)^{\alpha -1}\eta)\}. 
        \end{equation*}
    Then $\pi(\EE) = T_{IV,3}(\pi_{temp})$. Here is the associated symbol. 
    \[\EE= \scalebox{0.8}{\bordermatrix{
  & 0 & \cdots  &\alpha-2& \alpha-1 & \alpha & \alpha+1&\alpha +2\cr
    &\odot \cr
  &\odot \cr
  & \odot & \cdots & \odot &&&&\cr
  &&&&&&&\odot
}}.\]
\end{enumerate}
\end{prop}

\begin{prop}\label{temp4A,3}
    Let $\pi_{temp} = T_{I,1}^{\alpha +2}(T_{I,1}^{\alpha+1}(T_{I,1}^{\alpha}(\pi_{sc})))$ for $\alpha > 0$. 
    \begin{enumerate}[label = (\roman*)]
        \item The representation $T_{V,2}^{\pm}(\pi_{temp})$ is well-defined if and only if $\alpha = 1$. 
        \item The representation $T_{V,2}^{\pm}(\pi_{temp})$ is of critical type when $\alpha = 1$. 
        \item When $\alpha = 1$, define
        \begin{equation*}
            \mathcal{E}_{\pm} := \{([0, 0]_\rho,0,\pm)^2, ([3,3]_\rho, 0, \eta)\}. 
        \end{equation*}
    \end{enumerate}
    Then $\pi(\EE_{\pm}) = T_{V,2}^{\pm}(\pi_{temp})$. Here are the associated symbols. 
    \[\EE_{+}= \scalebox{0.8}{\bordermatrix{
  & 0 & 1& 2 & 3\cr
  & \oplus \cr
  &\oplus \cr
  &&&&\odot \cr
}}, \quad
\EE_{-}= \scalebox{0.8}{\bordermatrix{
  & 0 & 1& 2 & 3\cr
  & \ominus\cr
  &\ominus \cr
  &&&&\odot \cr
}}.\]
\end{prop}

We move onto the second case, which is $\pi_{temp} = T_{I,1}^{\alpha-1}(T_{I,1}^{\alpha +1}(T_{I,1}^{\alpha}(\pi_{sc})))$. 
\begin{prop}\label{temp4A,4}
    Let $\pi_{temp} = T_{I,1}^{\alpha-1}(T_{I,1}^{\alpha +1}(T_{I,1}^{\alpha}(\pi_{sc})))$ for $\alpha > 1$. 
    \begin{enumerate}[label = (\roman*)]
        \item The representation $T_{I,1}^{x}(\pi_{temp})$ is well-defined if and only if one of the following holds: 
        \begin{itemize}
            \item $x = \alpha +2$,  
            \item $x = \alpha $,
            \item $x = \alpha -2$ for $\alpha > 2$. 
        \end{itemize}
        \item When $x = \alpha +3$, $x = \alpha$, or $x = \alpha -2$ for $\alpha > 2$, the representation $T_{I,1}^{x}(\pi_{temp})$ is of critical type. 
        \item Let $x \in \{\alpha +2, \alpha, \alpha -2\}$. Define 
        \begin{equation*}
            \mathcal{E}_{\alpha +2} := \{([\alpha -3, \epsilon_\rho]_\rho,0,\eta),([\alpha -1, \alpha -1]_\rho, 0, (-1)^{\alpha  - \epsilon_\rho}\eta), ([\alpha +2, \alpha+2]_\rho, 0,(-1)^{\alpha -1 - \epsilon_\rho}\eta)\},
        \end{equation*}
        \begin{equation*}
            \mathcal{E}_{\alpha} := \{([\alpha -3, \epsilon_\rho]_\rho,0,\eta),([\alpha , \alpha ]_\rho, 0, (-1)^{\alpha  - \epsilon_\rho}\eta), ([\alpha +1, \alpha+1]_\rho, 0,(-1)^{\alpha -1 - \epsilon_\rho}\eta)\}.
        \end{equation*}
        When $\alpha > 2$, define
        \begin{equation*}
            \mathcal{E}_{\alpha -2} := \{([\alpha -4, \epsilon_\rho]_\rho,0,\eta),([\alpha -2, \alpha -2]_\rho, 0, (-1)^{\alpha  - 1 - \epsilon_\rho}\eta), ([\alpha +1, \alpha+1]_\rho, 0,(-1)^{\alpha -1 - \epsilon_\rho}\eta)\}.
        \end{equation*}
        Then we have $\pi(\EE_{x}) = T_{I,1}^{x}(\pi_{temp})$. Here are the associated symbols. 
        \[\EE_{\alpha+2}= \scalebox{0.8}{\bordermatrix{
  & \epsilon_{\rho} & \cdots &\alpha -3& \alpha-2 &\alpha-1& \alpha & \alpha+1 & \alpha+2 \cr
  & \odot & \cdots & \odot &&&&\cr
&&&&&\odot \cr
  &&&&&&&&\odot
}},\]
\[\EE_{\alpha}= \scalebox{0.8}{\bordermatrix{
  & \epsilon_{\rho} & \cdots &\alpha -3& \alpha-2 &\alpha-1& \alpha & \alpha+1 \cr
  & \odot & \cdots & \odot &&&&\cr
&&&&&&\odot \cr
  &&&&&&&\odot
}},\]
\[\EE_{\alpha -2}= \scalebox{0.8}{\bordermatrix{
  & \epsilon_{\rho} & \cdots &\alpha -4 &\alpha -3& \alpha-2 &\alpha-1& \alpha & \alpha+1 \cr
  & \odot & \cdots & \odot &&&&\cr
&&&&&\odot \cr
  &&&&&&&&\odot
}}.\]
    \end{enumerate}
\end{prop}

\begin{prop}\label{temp4A,5}
    Let $\pi_{temp} = T_{I,1}^{\alpha-1}(T_{I,1}^{\alpha +1}(T_{I,1}^{\alpha}(\pi_{sc})))$ for $\alpha > 1$.
    \begin{enumerate}[label = (\roman*)]
        \item The representation $T_{IV,3}(\pi_{temp})$ is well-defined if and only if $\alpha \in \mathbb{Z}_{>2}$. 
        \item The representation $T_{IV,3}(\pi_{temp})$ is not of critical type. 
        \item Define 
        \begin{equation*}
            \mathcal{E} := \{([0,0]_\rho,0,\eta)^2,([\alpha -3, 0]_\rho,0,\eta),  ([\alpha -1, \alpha -1]_\rho, 0, (-1)^{\alpha  - 1}\eta), ([\alpha +1, \alpha+1]_\rho, 0,(-1)^{\alpha -1}\eta)\}.
        \end{equation*}
    \end{enumerate}
    Then we have $\pi(\EE)= T_{IV,5}(\pi_{temp})$. Here is the associated symbol. 
    \[\EE= \scalebox{0.8}{\bordermatrix{
  & 0 & \cdots  &\alpha -3 &\alpha-2& \alpha-1 & \alpha & \alpha+1\cr
    &\odot \cr
  &\odot \cr
  & \odot & \cdots & \odot &&&&\cr
  &&&&&\odot \cr
  &&&&&&&\odot
}.}\]
\end{prop}

\begin{prop}\label{temp4A,6}
     Let $\pi_{temp} = T_{I,1}^{\alpha-1}(T_{I,1}^{\alpha +1}(T_{I,1}^{\alpha}(\pi_{sc})))$ for $\alpha > 1$.
     \begin{enumerate}[label = (\roman*)]
         \item The representation $T_{V,2}^{\pm}(\pi_{temp})$ is well-defined if and only if $\alpha = 2$. 
         \item When $\alpha = 2$, the representation $T_{V,2}^{\pm}(\pi_{temp})$ is of critical type. 
         \item When $\alpha = 2$, define 
         \begin{equation*}
            \mathcal{E}_{\pm} := \{([0, 0]_\rho,0,\pm1)^2, ([1,1]_\rho,0,\eta),([3,3]_\rho, 0, -\eta)\}. 
        \end{equation*}. 
    Then $\pi(\EE_{\pm}) = T_{V,2}^{\pm}(\pi_{temp})$. Here are the associated symbols. 
    \[\EE_{+}= \scalebox{0.8}{\bordermatrix{
  & 0 & 1& 2 & 3\cr
  & \oplus \cr
  &\oplus \cr
  &&\odot \cr
  &&&&\odot \cr
}}, \quad
\EE_{-}= \scalebox{0.8}{\bordermatrix{
  & 0 & 1& 2 & 3\cr
  & \ominus\cr
  &\ominus \cr
  &&\odot \cr
  &&&&\odot \cr
}}.\]
\end{enumerate}
\end{prop}
 The next case to consider is $\pi_{temp} = T_{IV,3}(T_{I,1}^{\alpha+1}(T_{I,1}^{\alpha}(\pi_{sc})))$. 
\begin{prop}\label{temp4A,7}
    Let $\pi_{temp} = T_{IV,3}(T_{I,1}^{\alpha+1}(T_{I,1}^{\alpha}(\pi_{sc})))$ for $\alpha \in \mathbb{Z}_{> 1}$. 
    \begin{enumerate}[label = (\roman*)]
        \item The representation $T_{I,1}^{x}(\pi_{temp})$ is well-defined if and only if $x = \alpha -1$ or $\alpha +2$. When $x = \alpha -1$ and $\alpha \in \mathbb{Z}_{>2}$, then $T_{I,1}^{x}(\pi_{temp})$ is the same as the representation $T_{IV,3}(T_{I,1}^{\alpha -1}(T_{I,1}^{\alpha +1}(T_{I,1}^{\alpha}(\pi_{sc}))))$ defined in Proposition \ref{temp4A,5}.  
        \item The representation $T_{I,1}^{x}(\pi_{temp})$ is of critical type only when $(x, \alpha) = (1,2)$. 
        \item Let $x \in \{\alpha -1, \alpha +2\}$. Define
        \begin{equation*}
            \mathcal{E}_{\alpha -1} := \{([0,0]_\rho,0,\eta)^2,([\alpha -3, 0]_\rho,0,\eta),([\alpha -1, \alpha -1]_\rho, 0, (-1)^{\alpha }\eta), ([\alpha +1, \alpha+1]_\rho, 0,(-1)^{\alpha -1 }\eta)\},
        \end{equation*}
        \begin{equation*}
            \mathcal{E}_{\alpha +2} := \{([0,0]_\rho,0,\eta)^2,([\alpha -2, 0]_\rho,0,\eta), ([\alpha +2, \alpha+2]_\rho, 0,(-1)^{\alpha -1 }\eta)\}.
        \end{equation*}
        Then we have $\pi(\EE_x) = T_{I,1}^{x}(\pi_{temp})$. Here are the associated symbols. 
         \[\EE_{\alpha -1}= \scalebox{0.8}{\bordermatrix{
  & 0 & \cdots  &\alpha -3 &\alpha-2& \alpha-1 & \alpha & \alpha+1\cr
    &\odot \cr
  &\odot \cr
  & \odot & \cdots & \odot &&&&\cr
  &&&&&\odot \cr
  &&&&&&&\odot
},}\]
\[\EE_{\alpha +2}= \scalebox{0.8}{\bordermatrix{
  & 0 & \cdots   &\alpha-2& \alpha-1 & \alpha & \alpha+1&\alpha +2\cr
    &\odot \cr
  &\odot \cr
  & \odot & \cdots & \odot &&&&\cr
  &&&&&&&\odot
}.}\]
\item The representations $T_{IV,3}(\pi_{temp})$ and $T_{V,2}^{\pm}(\pi_{temp})$ are not well-defined. 
    \end{enumerate}
\end{prop}
Now we move onto the case $\pi_{temp} = T_{V,2}^{\pm}(T_{I,1}^{\alpha+1}(T_{I,1}^{\alpha}(\pi_{sc})))$, which is well-defined only when $\alpha =1$, by Proposition \ref{temp3,3}. 
\begin{prop}\label{temp4A,8}
    Let $\pi_{temp}^{\pm} = T_{V,2}^{\pm}(T_{I,1}^{2}(T_{I,1}^{1}(\pi_{sc})))$ for $\alpha = 1$. 
    \begin{enumerate}[label = (\roman*)]
        \item The representation $T_{I,1}^{x}(\pi_{temp}^{\pm})$ is well-defined if and only if $x = 1$ or $x = 3$. 
        \item When $x = 1$ or $3$, the representation $T_{I,1}^{x}(\pi_{temp}^{\pm})$ is of critical type. 
        \item Let $x \in \{1,3\}$. Define
        \begin{equation*}
            \mathcal{E}_{1,\pm} := \{([0,0]_\rho,0,\pm1),([1,1]_\rho,0,\pm1), ([2,2]_\rho, 0,\eta)\},
        \end{equation*}
        \begin{equation*}
            \mathcal{E}_{3,\pm} := \{([0,0]_\rho,0,\pm1)^2, ([3,3]_\rho, 0,\eta)\}.
        \end{equation*}
        Then we have $\pi(\EE_{x,\pm}) = T_{I,1}^{x}(\pi_{temp}^{\pm})$. Here are the associated symbols. 
        \[\EE_{1,+}= \scalebox{0.8}{\bordermatrix{
  & 0 & 1& 2\cr
  & \oplus \cr
  &&\oplus \cr
  &&&\odot \cr
}}, \quad
\EE_{1,-}= \scalebox{0.8}{\bordermatrix{
  & 0 & 1& 2 \cr
  & \ominus\cr
  &&\ominus \cr
  &&&\odot \cr
}}\]
\[\EE_{3,+}= \scalebox{0.8}{\bordermatrix{
  & 0 & 1& 2&3\cr
  & \oplus \cr
  &\oplus \cr
  &&&&\odot \cr
}}, \quad
\EE_{3,-}= \scalebox{0.8}{\bordermatrix{
  & 0 & 1& 2 &3\cr
  & \ominus\cr
  &\ominus \cr
  &&&&\odot \cr
}}.\]
        \item The representations $T_{IV,3}(\pi_{temp})$ and $T_{V,2}^{\pm}(\pi_{temp})$ are not well-defined.
    \end{enumerate}
    
\end{prop}

The next case we want to consider is $\pi_{temp} = T_{I,1}^{\alpha -2}(T_{I,1}^{\alpha -1}(T_{I,1}^{\alpha}(\pi_{sc})))$. 
\begin{prop}\label{temp4A,9}
    Let $\pi_{temp} = T_{I,1}^{\alpha -2}(T_{I,1}^{\alpha -1}(T_{I,1}^{\alpha}(\pi_{sc})))$ for $\alpha > 2$. 
    \begin{enumerate}[label = (\roman*)]
        \item The representation $T_{I,1}^{x}(\pi_{temp})$ is well-defined if and only if $x = \alpha -3$ and $\alpha > 3$, or $x = \alpha +1$. 
        \item When $x = \alpha -3$ and $\alpha > 3$, or when $x = \alpha +1$, the representation $T_{I,1}^{x}(\pi_{temp})$ is of critical type. 
        \item Let $x \in \{\alpha-3, \alpha +1\}$. When $\alpha > 3$, define
        \begin{equation*}
            \mathcal{E}_{\alpha -3} := \{([\alpha -5, \epsilon_\rho]_\rho,0,\eta),([\alpha -3, \alpha -3]_\rho,0,(-1)^{\alpha - \epsilon_\rho}\eta),([\alpha , \alpha -2]_\rho, 0, (-1)^{\alpha -1 -\epsilon_\rho}\eta)\},
        \end{equation*}
         \begin{equation*}
            \EE_{\alpha +1} := \{([\alpha -4, \epsilon_\rho]_\rho,0,\eta),([\alpha -1, \alpha -2]_\rho,0,(-1)^{\alpha - 1-\epsilon_\rho}\eta),([\alpha+1, \alpha +1]_\rho, 0, (-1)^{\alpha -1 -\epsilon_\rho}\eta)\}.
        \end{equation*}
        Then we have $\pi(\EE_{x}) = T_{I,1}^{x}(\pi_{temp})$. Here are the associated symbols. 
        \[\EE_{\alpha -3}= \scalebox{0.8}{\bordermatrix{
  & \epsilon_\rho & \cdots  &\alpha-5 & \alpha -4 &\alpha -3 &\alpha-2& \alpha-1 & \alpha \cr
  & \odot & \cdots & \odot &&&&\cr
  &&&&&\odot \cr
  &&&&&&\odot&\odot & \odot
}},\]
\[\EE_{\alpha +1}= \scalebox{0.8}{\bordermatrix{
  & \epsilon_\rho & \cdots  &\alpha-4& \alpha -3 &\alpha-2& \alpha-1 & \alpha & \alpha+1\cr
    &\odot & \cdots & \odot \cr
  &&&&&\odot&\odot \cr
  &&&&&&&&\odot
}}.\]
    \end{enumerate}
\end{prop}

\begin{prop}\label{temp4A,10}
    Let $\pi_{temp} = T_{I,1}^{\alpha -2}(T_{I,1}^{\alpha -1}(T_{I,1}^{\alpha}(\pi_{sc})))$ for $\alpha > 2$.
    \begin{enumerate}[label = (\roman*)]
        \item The representation $T_{IV,3}(\pi_{temp})$ is well-defined if and only if $\alpha \in \mathbb{Z}_{>3}$. 
        \item The representation $T_{IV,3}(\pi_{temp})$ is not of critical type. 
        \item When $\alpha \in \mathbb{Z}_{>3}$, define
        \begin{equation*}
            \EE := \{([0,0]_\rho, 0,\eta)^2,[\alpha -4, 0]_\rho,0,\eta),([\alpha , \alpha -2]_\rho,0,(-1)^{\alpha - 1}\eta).
        \end{equation*}
        Then $\pi(\EE) = T_{IV,3}(\pi_{temp})$. Here is the associated symbol. 
        \[\EE  = \scalebox{0.8}{\bordermatrix{
  & 0 & \cdots  &\alpha-4& \alpha -3 &\alpha-2& \alpha-1 & \alpha \cr
  &\odot \cr
  &\odot \cr
    &\odot & \cdots & \odot \cr
  &&&&&\odot&\odot&\odot \cr
}}.\]
    \end{enumerate}
\end{prop}

\begin{prop}\label{temp4A,11}
     Let $\pi_{temp} = T_{I,1}^{\alpha -2}(T_{I,1}^{\alpha -1}(T_{I,1}^{\alpha}(\pi_{sc})))$ for $\alpha > 2$.
     \begin{enumerate}[label = (\roman*)]
         \item The representation $T_{V,2}^{\pm}(\pi_{temp})$ is well-defined if and only if $\alpha = 3$. 
         \item When $\alpha =3$, the representation $T_{V,2}^{\pm}(\pi_{temp})$ is of critical type. 
         \item When $\alpha = 3$, define 
         \begin{equation*}
            \EE_{\pm} := \{([0,0]_\rho, 0,\pm1)^2,([3,1]_\rho,0,-\eta)\}.
        \end{equation*}
        Then we have $\pi(\EE_{\pm}) = T_{V,2}^{\pm}(\pi_{temp})$. Here are the associated symbols. 
         \[\EE_{+}= \scalebox{0.8}{\bordermatrix{
  & 0 & 1& 2&3\cr
  & \oplus \cr
  &\oplus \cr
  &&\odot&\odot &\odot \cr
}}, \quad
\EE_{-}= \scalebox{0.8}{\bordermatrix{
  & 0 & 1& 2&3\cr
  & \ominus \cr
  &\ominus \cr
  &&\odot&\odot &\odot \cr
}}.\]
     \end{enumerate}
\end{prop}

Now we move onto the case $\pi_{temp} = T_{IV,3}(T_{I,1}^{\alpha-1}(T_{I,1}^{\alpha}(\pi_{sc})))$. By Proposition \ref{temp3,5}, this is well -defined if and only if $\alpha \in \mathbb{Z}_{>2}$. 
\begin{prop}\label{temp4A,12}
    Let  $\pi_{temp} = T_{IV,3}(T_{I,1}^{\alpha-1}(T_{I,1}^{\alpha}(\pi_{sc})))$ for $\alpha \in \mathbb{Z}_{>2}$. 
    \begin{enumerate}[label = (\roman*)]
        \item The representation $T_{I,1}^{x}(\pi_{temp})$ is well-defined if and only if $x = \alpha -2$ or $\alpha +1$. 
        \item The representation $T_{I,1}^{x}(\pi_{temp})$ is of critical type if and only if $(x,\alpha) = (1,3)$.
        \item Let $x \in \{\alpha-2, \alpha +1\}$. Define
        \begin{equation*}
            \mathcal{E}_{\alpha -2} := \{(([0,0]_\rho, 0,\eta)^2, [\alpha -4, 0]_\rho,0,\eta),([\alpha -2, \alpha -2]_\rho,0,(-1)^{\alpha - 1}\eta),([\alpha , \alpha -1]_\rho, 0, (-1)^{\alpha}\eta)\},
        \end{equation*}
         \begin{equation*}
            \mathcal{E}_{\alpha +1} := \{(([0,0]_\rho, 0,\eta)^2, [\alpha -3, 0]_\rho,0,\eta),([\alpha -1, \alpha -1]_\rho,0,(-1)^{\alpha}\eta),([\alpha +1, \alpha +1]_\rho, 0, (-1)^{\alpha-1}\eta)\}.
        \end{equation*}
        Then we have $\pi(\EE_{x}) = T_{I,1}^{x}(\pi_{temp})$. Here are the associated symbols. 
        \[\EE_{\alpha -2}= \scalebox{0.8}{\bordermatrix{
  & 0 & \cdots  & \alpha -4 &\alpha -3 &\alpha-2& \alpha-1 & \alpha \cr
  &\odot \cr
  &\odot \cr
  & \odot & \cdots & \odot &&&&\cr
  &&&&&\odot \cr
  &&&&&&\odot
}},\]
\[\EE_{\alpha +1}= \scalebox{0.8}{\bordermatrix{
  & 0 & \cdots  &\alpha -3 &\alpha-2& \alpha-1 & \alpha &\alpha +1\cr
  &\odot \cr
  &\odot \cr
  & \odot & \cdots & \odot &&&&\cr
  &&&&&\odot \cr
  &&&&&&&\odot
}}.\]
    \end{enumerate}
\end{prop}

The next case we'll look at is $\pi_{temp} = T_{V,2}^{\pm}(T_{I,1}^{1}(T_{I,1}^{2}(\pi_{sc})))$, which is well-defined only for $\alpha =2$, by Proposition \ref{temp3,6}. 
\begin{prop}\label{temp4A,13}
    Let $\pi_{temp}^{\pm} = T_{V,2}^{\pm}(T_{I,1}^{1}(T_{I,1}^{2}(\pi_{sc})))$, and $\alpha = 2$. 
    \begin{enumerate}[label = (\roman*)]
        \item The representation $T_{I,1}^{x}(\pi_{temp}^{\pm})$ is well-defined if and only if $x = 3$. 
        \item When $x = 3$, the representation $T_{I,1}^{x}(\pi_{temp}^{\pm})$ is of critical type. In this case, it is the same as the representations $T_{V,2}^{\pm}(T_{I,1}^{1}(T_{I,1}^{3}(T_{I,1}^{3}(\pi_{sc}))))$, as described in Proposition \ref{temp4A,6}. 
\item The representations $T_{IV,3}(\pi_{temp})$ and $T_{V,2}^{\pm}(\pi_{temp})$ are not well-defined.
    \end{enumerate}
\end{prop}

A similar case to the one above is $\pi_{temp} = T_{I,1}^{2}(T_{V,2}^{\pm}(T_{I,1}^{1}(\pi_{sc})))$, which is well-defined only for $\alpha = 1$. 
\begin{prop}\label{temp4A,14}
    Let $\pi_{temp}^{\pm} = T_{I,1}^{2}((T_{V,2}^{\pm}T_{I,1}^{1}(\pi_{sc})))$ and $\alpha = 1$. 
    \begin{enumerate}[label = (\roman*)]
        \item The representation $T_{I,1}^{x}(\pi_{temp}^{\pm})$ is well-defined if and only if $x = 1$ or $x = 3$. When $x = 1$, then $T_{I,1}^{x}(\pi_{temp}^{\pm})$ is the same as the representation $T_{I,1}^{x}(T_{V,2}^{\pm}(T_{I,1}^{2}(T_{I,1}^{1}(\pi_{sc}))))$, as described in Proposition \ref{temp4A,9}. 
        \item When $x = 1$ or $x = 3$, the representation $T_{I,1}^{x}(\pi_{temp}^{\pm})$ is of critical type. 
        \item The representations $T_{IV,3}(\pi_{temp})$ and $T_{V,2}^{\pm}(\pi_{temp})$ are not well-defined.
    \end{enumerate}
\end{prop}
 Another similar case would be $\pi_{temp} = T_{I,1}^{2}(T_{I,1}^{1}(T_{V,2}^{\pm}(\pi_{sc})))$, which is well-defined if and only if $\alpha = 0$.  
\begin{prop}\label{temp4A,15}
    Let  $\pi_{temp}^{\pm} = T_{I,1}^{2}(T_{I,1}^{1}(T_{V,2}^{\pm}(\pi_{sc})))$ and $\alpha = 0$. 
    \begin{enumerate}[label = (\roman*)]
        \item The representation $T_{I,1}^{x}(\pi_{temp}^{\pm})$ is well-defined if and only if $x = 1$ or $x = 3$. 
        \item When $x = 1$ or $x = 3$, the representation $T_{I,1}^{x}(\pi_{temp}^{\pm})$ is of critical type. 
        \item Let $x \in \{1,3\}$. Define
        \begin{flalign*}  
            \EE_{1,\pm} &:= \{([1,1]_\rho, 0,\pm1),([2,2]_\rho,0,\pm1)\}, \\
            \EE_{3,\pm} &:= \{([0,0]_\rho, 0,\pm1),([3,3]_\rho,0,\pm1)\}.
        \end{flalign*}
        Then we have that $\pi(\EE_{x,\pm}) = T_{I,1}^{x}(\pi_{temp}^{\pm})$. Here are the associated symbols. 
         \[\EE_{1,+}= \scalebox{0.8}{\bordermatrix{
   & 1& 2\cr
  &\oplus \cr
  &&\oplus \cr
}}, \quad
\EE_{1,-}= \scalebox{0.8}{\bordermatrix{
  & 1& 2 \cr
  &\ominus \cr
  &&\ominus \cr
}},\]
\[\EE_{3,+}= \scalebox{0.8}{\bordermatrix{
  & 0 & 1& 2&3\cr
  & \oplus \cr
  &&&&\oplus \cr
}}, \quad
\EE_{3,-}= \scalebox{0.8}{\bordermatrix{
  & 0 & 1& 2 &3\cr
  & \ominus\cr
  &&&&\ominus \cr
}}.\]
    \end{enumerate}
\end{prop}

\begin{prop}\label{temp4A,16}
    Let  $\pi_{temp}^{\pm} = T_{I,1}^{2}(T_{I,1}^{1}(T_{V,2}^{\pm}(\pi_{sc})))$ and $\alpha = 0$. 
    \begin{enumerate}[label = (\roman*)]
        \item The representation $T_{IV,3}(\pi_{temp}^{\pm})$ is well-defined and of critical type. 
        \item The representation $T_{V,2}^{\pm}(\pi_{temp}^{\pm})$ is not well-defined. 
        \item Define
        \begin{flalign*}  
            \EE_{\pm} &:= \{([0,0]_\rho, 0,\pm1)^3,([2,2]_\rho,0,\pm1)\}. 
        \end{flalign*}
        Then $\pi(\EE_{\pm}) = T_{IV,3}(\pi_{temp}^{\pm})$. Here are the associated symbols. 
        \[\EE_{+}= \scalebox{0.8}{\bordermatrix{
   &0& 1& 2\cr
  &\oplus \cr
  &\oplus \cr
  &\oplus \cr
  &&&\oplus \cr
}}, \quad
\EE_{-}= \scalebox{0.8}{\bordermatrix{
  & 0&1& 2 \cr
  &\ominus \cr
  &\ominus \cr
  &\ominus \cr
  &&&\ominus \cr
}}.\]
    \end{enumerate}
\end{prop}

Another tempered representation of corank $3$ is $\pi_{temp} = T_{IV,3}(T_{I,1}^{1}(T_{V,2}^{\pm}(\pi_{sc})))$, which is also well-defined only when $\alpha = 0$. 
\begin{prop}\label{temp4A,17}
    Let $\pi_{temp}^{\pm} = T_{IV,3}(T_{I,1}^{1}(T_{V,2}^{\pm}(\pi_{sc})))$ and $\alpha = 0$. 
    \begin{enumerate}[label = (\roman*)]
        \item The representation $T_{I,1}^{x}(\pi_{temp}^{\pm})$ is well-defined if and only if $x = 2$. The representations $T_{I,1}^{2}(\pi_{temp}^{\pm})$ are the same as the representations $T_{IV,3}(T_{I,1}^{2}(T_{I,1}^{1}(T_{V,2}^{\pm}(\pi_{sc}))))$, as described in Proposition \ref{temp4A,16}. 
        \item The representation $T_{I,1}^{2}(\pi_{temp}^{\pm})$ is of critical type. 
        \item The representations $T_{IV,3}(\pi_{temp})$ and $T_{V,2}^{\pm}(\pi_{temp})$ are not well-defined.
     \end{enumerate}
\end{prop}

We now proceed to the case $\pi_{temp} = T_{I,1}^{\frac{3}{2}}(T_{I,2}^{\frac{1}{2}}(\pi_{sc}))$, defined only when $\alpha = \frac{1}{2}$. 
\begin{prop}\label{temp4A,18}
    Let $\pi_{temp} = T_{I,1}^{\frac{3}{2}}(T_{I,2}^{\frac{1}{2}}(\pi_{sc}))$ and $\alpha = \frac{1}{2}$. 
    \begin{enumerate}[label = (\roman*)]
        \item The representation $T_{I,1}^{x}(\pi_{temp})$ is well-defined if and only if $x = \frac{5}{2}$. 
        \item The representation $T_{I,1}^{\frac{5}{2}}(\pi_{temp})$ is of critical type. 
        \item The set $\Psi(T_{I,1}^{\frac{5}{2}}(\pi_{temp}))$ is a singleton, and we have $T_{I,1}^{\frac{5}{2}}(\pi_{temp}) = \pi(\EE)$, where 
        \[\EE= \scalebox{0.8}{\bordermatrix{
   &\frac{1}{2}& \frac{3}{2}& \frac{5}{2}\cr
&\oplus \cr
&&&\oplus
}}.\]
  \item The representations $T_{IV,3}(\pi_{temp})$ and $T_{V,2}^{\pm}(\pi_{temp})$ are not well-defined. 
    \end{enumerate}
\end{prop}

The next case is $\pi_{temp} = T_{I,1}^{\alpha}(T_{II,3}^{\frac{1}{2}}(\pi_{sc}))$, defined for $\alpha \in \frac{1}{2} + \mathbb{Z}_{> 0}$. 
\begin{prop}\label{temp4A,19}
    Let $\pi_{temp} = T_{I,1}^{\alpha}(T_{II,3}^{\frac{1}{2}}(\pi_{sc}))$ for $\alpha \in \frac{1}{2} + \mathbb{Z}$. 
    \begin{enumerate}[label = (\roman*)]
        \item The representation $T_{I,1}^{x}(\pi_{temp})$ is well-defined if and only if $x = \alpha +1$, or $x = \alpha -1$, and $\alpha \geq \frac{5}{2}$. 
        \item The representation $T_{I,1}^{\alpha +1}(\pi_{temp})$ if of critical type when $\alpha \in \{\frac{3}{2}, \frac{5}{2}\}$. The representation $T_{I,1}^{\alpha -1}(\pi_{temp})$ is of critical type only when $\alpha = \frac{5}{2}$. 
        \item Let $x \in \{\alpha +1, \alpha -1\}$. Define
        \begin{flalign*}  
            \EE_{\alpha +1} &:= \{([\frac{1}{2}, \frac{1}{2}]_\rho, 0,-1)^2,([\alpha -2, \frac{1}{2}]_\rho, 0, -1), ([\alpha +1, \alpha+1]_\rho, 0, (-1)^{\alpha - \frac{1}{2}})\}. 
        \end{flalign*}
        When $\alpha \in \frac{3}{2} + \mathbb{Z}_{>0}$, define
        \begin{flalign*}  
            \EE_{\alpha -1} &:= \{([\frac{1}{2}, \frac{1}{2}]_\rho, 0,-1)^2,([\alpha -3, \frac{1}{2}]_\rho, 0, -1), ([\alpha , \alpha-1]_\rho, 0, (-1)^{\alpha + \frac{1}{2}})\}. 
        \end{flalign*}
        Then $\pi(\EE_x) = T_{I,1}^{x}(\pi_{temp})$. Here are the associated symbols. 
        \[\EE_{\alpha +1}= \scalebox{0.8}{\bordermatrix{
  & \frac{1}{2} & \cdots  & \alpha -2 &\alpha -1 &\alpha & \alpha+1  \cr
  &\ominus \cr
  &\ominus \cr
  & \ominus & \cdots & \odot &&&&\cr
  &&&&&&\odot \cr
}},\]
\[\EE_{\alpha -1}= \scalebox{0.8}{\bordermatrix{
  & \frac{1}{2} & \cdots &\alpha -3 &\alpha-2& \alpha-1 & \alpha \cr
  &\ominus \cr
  &\ominus \cr
  & \ominus & \cdots & \odot &&&&\cr
  &&&&&\odot&\odot \cr
}}.\]
\item The representations $T_{IV,3}(\pi_{temp})$ and $T_{V,2}^{\pm}(\pi_{temp})$ are not well-defined. 
    \end{enumerate}
\end{prop}

Another case to consider is $\pi_{temp} = T_{I,1}^{\frac{3}{2}}(T_{III,2}^{\frac{1}{2}})$, which is well-defined only when $\alpha = \frac{1}{2}$, by Proposition \ref{temp3,13}. This will not appear in our Arthur type list later on, since it involves negative supercuspidal support, but we include it here for completion.  

\begin{prop}\label{temp4A,20}
    Let $\pi_{temp} = T_{I,1}^{\frac{3}{2}}(T_{III,2}^{\frac{1}{2}})$ and $\alpha = \frac{1}{2}$. 
    \begin{enumerate}[label = (\roman*)]
        \item The representation $T_{I,1}^{x}(\pi_{temp})$ is well-defined if and only if $x = \frac{5}{2}$. 
        \item The representation $T_{I,1}^{\frac{5}{2}}(\pi_{temp})$ is of critical type. 
        \item The following is the set of extended multi-segments $\EE$ (up to row exchanges) such that $\pi(\EE) = T_{I,1}^{\frac{3}{2}}(\pi_{temp})$: 
            \[\Bigg\{\scalebox{0.8}{\bordermatrix{ 
        &\frac{1}{2}&\frac{3}{2} & \frac{5}{2} \cr
        &\ominus&\cr
        &&&\ominus\cr
        }},{\bordermatrix{ 
        &-\frac{1}{2}&\frac{1}{2}&\frac{3}{2}&\frac{5}{2} \cr
        &\oplus&\ominus&\cr
        &&&&\ominus\cr}}\Bigg\}.\]
        \item The representations $T_{IV,3}(\pi_{temp})$ and $T_{V,2}^{\pm}(\pi_{temp})$ are not well-defined. 
    \end{enumerate}
\end{prop}
\begin{proof}
    Parts $(i)$ through $(iii)$ follows from Proposition \ref{temp3,13} and Remark \ref{rmk well-defined for Temp}. Part $(iv)$ follows from the good parity condition. 
\end{proof}

The next case to look at is $\pi_{temp} = T_{I,1}^{\alpha}(T_{IV,5}(\pi_{sc}))$, which is well-defined if and only if $\alpha \in \mathbb{Z}_{>0}$, by Proposition \ref{temp3,14}. 
\begin{prop}\label{temp4A,21}
    Let $\pi_{temp} = T_{I,1}^{\alpha}(T_{IV,5}(\pi_{sc}))$ and $\alpha \in \mathbb{Z}_{>0} $. 
    \begin{enumerate}[label = (\roman*)]
        \item The representation $T_{I,1}^{x}(\pi_{temp})$ is well-defined if and only if $x = \alpha +1$, or $x = \alpha -1$, and $\alpha > 1$. 
        \item The representation $T_{I,1}^{x}(\pi_{temp})$ is of critical type if and only if $(x,\alpha) = (1,2)$ or $(2,1)$. 
        \item  Let $\{x \in \alpha +1, \alpha -1\}$. Define
        \begin{flalign*}  
            \EE_{\alpha +1} &:= \{([0,0]_\rho,0,\eta)^4, ([\alpha -2, 0]_\rho, 0, \eta), ([\alpha +1, \alpha+1]_\rho, 0, (-1)^{\alpha - 1}\eta)\}. 
        \end{flalign*}
        When $\alpha > 1$, define 
        \begin{flalign*}  
            \EE_{\alpha -1} &:= \{([0,0]_\rho,0,\eta)^4, ([\alpha -3, 0]_\rho, 0, \eta), ([\alpha -1, \alpha-1]_\rho, 0, (-1)^{\alpha }\eta), ([\alpha, \alpha]_\rho, 0, (-1)^{\alpha -1}\eta)\}. 
        \end{flalign*}
        Then we have $\pi(\EE_{x}) = T_{I,1}^{x}(\pi_{temp})$. Here are the associated symbols. 
        \[\EE_{\alpha +1}= \scalebox{0.8}{\bordermatrix{
  & 0 & \cdots  & \alpha -2 &\alpha -1 &\alpha& \alpha+1  \cr
  &\odot \cr
  &\odot \cr
  &\odot \cr
  &\odot \cr
  & \odot & \cdots & \odot &&&&\cr
  &&&&&&\odot \cr
}},\]
\[\EE_{\alpha -1}= \scalebox{0.8}{\bordermatrix{
  & 0 & \cdots &\alpha -3 &\alpha-2& \alpha-1 & \alpha \cr
  &\odot \cr
  &\odot \cr
  &\odot \cr
  &\odot \cr
  & \odot & \cdots & \odot &&&&\cr
  &&&&&\odot \cr
  &&&&&&\odot
}}.\]
\item The representations $T_{IV,3}(\pi_{temp})$ and $T_{V,2}^{\pm}(\pi_{temp})$ are not well-defined. 
    \end{enumerate}
\end{prop}

We now move onto the case $\pi_{temp} = T_{I,1}^{1}(T_{V,4}^{\pm}(\pi_{sc}))$, which is well-defined if and only if $\alpha = 0$. 
\begin{prop}\label{temp4A,22}
    Let $\pi_{temp}^{\pm} = T_{I,1}^{1}(T_{V,4}^{\pm}(\pi_{sc}))$ and $\alpha = 0$. 
    \begin{enumerate}[label = (\roman*)]
        \item The representation $T_{I,1}^{x}(\pi_{temp}^{\pm})$ is well-defined if and only if $x = 2$. In this case, it is the same as the representation $T_{IV,3}(T_{I,1}^{2}(T_{I,1}^{1}(T_{V,2}^{\pm}(\pi_{sc}))))$ described in Proposition \ref{temp4A,16}. 
        \item The representation $T_{I,1}^{2}(\pi_{temp}^{\pm})$ is of critical type. 
        \item The representations $T_{IV,3}(\pi_{temp})$ and $T_{V,2}^{\pm}(\pi_{temp})$ are not well-defined. 
    \end{enumerate}
\end{prop}

The next case is $\pi_{temp} = T_{V,4}^{\pm}(T_{I,1}^{1}(\pi_{sc}))$, which is well-defined only for $\alpha = 1$, by Proposition \ref{temp3,18}.
\begin{prop}\label{temp4A,23}
    Let $\pi_{temp}^{\pm} = T_{V,4}^{\pm}(T_{I,1}^{1}(\pi_{sc}))$ and $\alpha =1$. 
    \begin{enumerate}[label = (\roman*)]
        \item The representation $T_{I,1}^{x}(\pi_{temp}^{\pm})$ is well-defined if and only if $x = 2$. 
        \item The representation $T_{I,1}^{2}(\pi_{temp}^{\pm})$ is of critical type. 
        \item Define 
        \begin{flalign*}  
            \EE_{\pm} &:= \{([0,0]_\rho, 0,\pm1)^4,([2,2]_\rho,0,\pm1)\}. 
        \end{flalign*}
        Then we have $\pi(\EE_{\pm}) = T_{I,1}^{2}(\pi_{temp}^{\pm})$. Here are the associated symbols. 
         \[\EE_{+}= \scalebox{0.8}{\bordermatrix{
   &0& 1& 2\cr
  &\oplus \cr
  &\oplus \cr
  &\oplus \cr
  &\oplus \cr
  &&&\oplus \cr
}}, \quad
\EE_{-}= \scalebox{0.8}{\bordermatrix{
  & 0&1& 2 \cr
  &\ominus \cr
  &\ominus \cr
  &\ominus \cr
  &\ominus \cr
  &&&\ominus \cr
}}.\]
\item The representations $T_{IV,3}(\pi_{temp}^{\pm})$ and $T_{V,2}^{\pm}(\pi_{temp}^{\pm})$ are not well-defined. 
    \end{enumerate}
\end{prop}

By Proposition \ref{temp3,20}, we have another tempered representation of corank $3$ which is only well-defined at $\alpha = 1$, namely $\pi_{temp} = T_{I,2}^{1}(T_{IV,3}(\pi_{sc}))$. 
\begin{prop}\label{temp4A,24}
    Let $\pi_{temp} = T_{I,2}^{1}(T_{IV,3}(\pi_{sc}))$ and $\alpha = 1$. 
    \begin{enumerate}[label = (\roman*)]
        \item The representation $T_{I,1}^{x}(\pi_{temp})$ is well-defined if and only if $x = 2$. 
        \item The representation $T_{I,1}^{2}(\pi_{temp})$ is of critical type. 
        \item Define 
\begin{flalign*}  
            \EE &:= \{([0,0]_\rho, 0,\eta),([1,1]_\rho,0,\eta),([2,2]_\rho,0,\eta)\}. 
        \end{flalign*}        
        Then we have $\pi(\EE) = T_{I,1}^{2}(\pi_{temp})$. Here is the associated symbol. 
        \[\EE= \scalebox{0.8}{\bordermatrix{
  & 0&1& 2 \cr
  &\odot \cr
  &&\odot \cr
  &&&\odot \cr
}}.\]
    \end{enumerate}
\end{prop}

\begin{remark}
    Note that despite their similarities, the extended multi-segment in the proposition above is not the same as the extended multi-segment of $T_{I,1}^{1}(T_{V,2}^{\pm}(T_{I,1}^{2}(T_{I,1}^{1}(\pi_{sc}))))$, as shown in Proposition \ref{temp4A,8}, due to the sign restrictions. 
\end{remark}

\begin{prop}\label{temp4A,26}
    Let $\pi_{temp} = T_{I,2}^{1}(T_{IV,3}(\pi_{sc}))$ and $\alpha = 1$. 
    \begin{enumerate}[label = (\roman*)]
        \item The representation $T_{IV,3}(\pi_{temp})$ is well-defined and of critical type. 
        \item Define
        \begin{flalign*}  
            \EE &:= \{([0,0]_\rho, 0,\eta)^3,([1,1]_\rho,0,\eta)^2\}.
        \end{flalign*}      
        Then $\pi(\EE) = T_{IV,3}(\pi_{temp})$. Here is the associated symbol. 
        \[\EE= \scalebox{0.8}{\bordermatrix{
  & 0&1 \cr
  &\odot \cr
   &\odot \cr
   &\odot \cr
  &&\odot \cr
  &&\odot \cr
}}.\]
\item The representation $T_{V,2}^{\pm}(\pi_{temp})$ is not well-defined. 
    \end{enumerate}
\end{prop}

Five more cases remain to wrap up Case $(A)$. The next one is $\pi_{temp} = T_{II,3}^{1}(T_{IV,3}(\pi_{sc}))$, which is well-defined only when $\alpha \in \mathbb{Z}_{>0}$, by Proposition \ref{temp3,20}. 
\begin{prop}\label{temp4A,27}
    Let $\pi_{temp} = T_{II,3}^{1}(T_{IV,3}(\pi_{sc}))$, and $\alpha \in \mathbb{Z}_{>0}$. 
    \begin{enumerate}[label = (\roman*)]
        \item The representation $T_{I,1}^{x}(\pi_{temp})$ is well-defined if and only if $x = \alpha$. 
        \item The representation $T_{I,1}^{\alpha}(\pi_{temp})$ is of critical type only when $\alpha \in \{1,2\}$. 
        \item Define 
        \begin{flalign*}  
            \EE &:= \{([\alpha -2,0]_\rho, 0,\eta),([1,1]_\rho,0,-\eta)^2,([\alpha, \alpha]_\rho,0, *)\}.
        \end{flalign*}      
        Then we have $\pi(\EE) = T_{I,1}^{\alpha}(\pi_{temp})$. Here is the associated symbol. 
        \[\EE= \scalebox{0.8}{\bordermatrix{
  & 0&1&\cdots & \alpha -2 & \alpha -1 & \alpha \cr
  &\odot &\cdots & \cdots & \odot \cr
   &&\odot \cr
   &&\odot \cr
  &&&&&&\odot \cr
}}.\]
    \end{enumerate}
\end{prop}

\begin{prop}\label{temp4A,28}
     Let $\pi_{temp} = T_{II,3}^{1}(T_{IV,3}(\pi_{sc}))$, and $\alpha \in \mathbb{Z}_{>0}$. 
     \begin{enumerate}[label = (\roman*)]
         \item The representation $T_{IV,3}(\pi_{temp})$ is well-defined.
         \item The representation $T_{IV,3}(\pi_{temp})$ is of critical type if and only if $\alpha = 1$. 
         \item Define 
        \begin{flalign*}  
            \EE &:= \{([0,0]_\rho,0,\eta)^2,([\alpha -1,0]_\rho, 0,\eta),([1,1]_\rho,0,-\eta)^2\}.
        \end{flalign*}
        Then we have $\pi(\EE)= T_{IV,3}(\pi_{temp})$. Here is the associated symbol. 
        \[\EE= \scalebox{0.8}{\bordermatrix{
  & 0&1&\cdots & \alpha -2 & \alpha -1  \cr
  &\odot \cr
  &\odot \cr
  &\odot &\cdots & \cdots & \cdots & \odot \cr
   &&\odot \cr
   &&\odot \cr
}}.\]
\item The representation $T_{V,2}^{\pm}(\pi_{temp})$ is not well-defined. 
     \end{enumerate}
\end{prop}

The next case we need to consider is $\pi_{temp} = T_{I,2}^{1}(T_{V,2}^{\pm}(\pi_{sc}))$. This is well-defined if and only if $\alpha = 0$, by Proposition \ref{temp3,21}. 
\begin{prop}\label{temp4A,29}
    Let $\pi_{temp}^{\pm} = T_{I,2}^{1}(T_{V,2}^{\pm}(\pi_{sc}))$ and $\alpha = 0$. 
    \begin{enumerate}
        \item The representation $T_{I,1}^{x}(\pi_{temp})$ is well-defined if and only if $x = 2$. When $x = 2$, these are the same as the representations $T_{I,1}^{1}(T_{I,1}^{2}(T_{I,1}^{1}(T_{V,2}^{\pm}(\pi_{sc}))))$, as described in Proposition \ref{temp4A,15}. 
        \item The representation $T_{I,1}^{2}(\pi_{temp})$ is of critical type. 
        \item The representations $T_{IV,3}(\pi_{temp})$ and $T_{V,2}^{\pm}(\pi_{temp})$ are not well-defined. 
    \end{enumerate}
\end{prop}

Three more cases remain in Case $(A)$. The next one to consider is $\pi_{temp} = T_{III,2}^{1}(\pi_{sc})$, which is defined only when $\alpha = 1$, by Proposition \ref{temp3,25}. 
\begin{prop}\label{temp4A,30}
    Let $\pi_{temp} = T_{III,2}^{1}(\pi_{sc})$ and $\alpha = 1$. 
    \begin{enumerate}[label = (\roman*)]
        \item The representation $T_{I,1}^{x}(\pi_{temp})$ is well-defined if and only if $x = 2$. When $x = 2$, the representation $T_{I,1}^{2}(\pi_{temp})$ is of critical type, and is the same as the representation $T_{I,1}^{2}(T_{I,2}^{1}(T_{IV,3}(\pi_{sc})))$, as described in Proposition \ref{temp4A,24}. 
        \item The representation $T_{IV,3}(\pi_{temp})$ is well-defined and of critical type. It is the same as the representation $T_{IV,3}(T_{I,2}^{1}(T_{IV,3}(\pi_{sc})))$, as described in Proposition \ref{temp4A,26}. 
        \item The representation $T_{V,2}^{\pm}(\pi_{temp})$ is not well-defined. 
    \end{enumerate} 
\end{prop}

The second last case in Case $(A)$ is $\pi_{temp} = T_{IV,7}(\pi_{sc})$, which is defined only when $\alpha \in \mathbb{Z}_{>0}$, by Proposition \ref{temp3,26}. 
\begin{prop}\label{temp4A,31}
    Let $\pi_{temp} = T_{IV,7}(\pi_{sc})$ and $\alpha \in \mathbb{Z}_{>0}$. 
    \begin{enumerate}[label = (\roman*)]
        \item The representation $T_{I,1}^{x}(\pi_{temp})$ is well-defined if and only if $x = \alpha$. 
        \item The representation $T_{I,1}^{\alpha}(\pi_{temp})$ is of critical type if and only if $\alpha = 1$. 
        \item Define 
        \begin{flalign*}  
            \EE &:= \{([0,0]_\rho,0,\eta)^6,([\alpha -2,0]_\rho, 0,\eta),([\alpha,\alpha]_\rho,0,(-1)^{\alpha-1}\eta)\}.
        \end{flalign*}
        Then we have $\pi(\EE) = T_{I,1}^{\alpha}(\pi_{temp})$. Here is the associated symbol. 
        \[\EE= \scalebox{0.8}{\bordermatrix{
  & 0&\cdots & \alpha -2 & \alpha -1  &\alpha\cr
  &\odot \cr
  &\odot \cr
  &\odot \cr
  &\odot \cr
  &\odot \cr
  &\odot \cr
  &\odot &\cdots  & \odot \cr
  &&&&&\odot
}}.\]
\item The representations $T_{IV,3}(\pi_{temp})$ and $T_{V,2}^{\pm}(\pi_{temp})$ are not well-defined. 
    \end{enumerate}
\end{prop}

Finally, the last case in Case $(A)$ is $\pi_{temp} = T_{V,6}^{\pm}(\pi_{sc})$, which is defined only when $\alpha = 0$, by Proposition \ref{temp3,25}. 
\begin{prop}\label{temp4A,32}
    Let $\pi_{temp}^{\pm} = T_{V,6}^{\pm}(\pi_{sc})$ and $\alpha = 0$.
    \begin{enumerate}[label = (\roman*)]
        \item The representations $T_{I,1}^{x}(\pi_{temp}^{\pm})$ are well-defined if and only if $x = 1$. 
        \item The representations $T_{I,1}^{1}(\pi_{temp}^{\pm})$ are of critical type. 
        \item Define 
        \begin{flalign*}  
            \EE_{\pm} &:= \{([0,0]_\rho,0,\pm1)^5,([1,1]_\rho,0,\pm1)\}.
        \end{flalign*}
        Then $\pi(\EE_{\pm}) = T_{I,1}^{1}(\pi_{temp}^{\pm})$. Here are the associated symbols. 
        \[\EE_{+}= \scalebox{0.8}{\bordermatrix{
  & 0&1\cr
  &\oplus \cr
  &\oplus \cr
  &\oplus \cr
  &\oplus \cr
  &\oplus \cr
  &&\oplus   
}}, \quad
\EE_{-}= \scalebox{0.8}{\bordermatrix{
  & 0&1\cr
  &\ominus \cr
  &\ominus \cr
  &\ominus \cr
  &\ominus \cr
  &\ominus \cr
  &&\ominus   
}},\]
\item The representations $T_{IV,3}(\pi_{temp}^{\pm})$ and $T_{V,2}^{\pm}(\pi_{temp}^{\pm})$ are not well-defined. 
    \end{enumerate}
\end{prop}

This concludes all our discussion of Case $(A)$. We move onto Case $(B)$. 
\subsection{\texorpdfstring{Case $(B): \pi \hookrightarrow \rho\lvert \cdot \rvert^{x_1} \times \rho \lvert \cdot \rvert^{x_2} \rtimes \pi_{temp}$, where $\pi_{temp}$ is a corank 2 tempered representation of good parity.}{}}
In this section, we consider the possibility when  $\pi \hookrightarrow \rho\lvert \cdot \rvert^{x_1} \times \rho \lvert \cdot \rvert^{x_2} \rtimes \pi_{temp}$, where $\pi_{temp}$ is tempered of corank $2$. Here there are five possibilities for the form of $\pi$, namely $\pi = T_{I,2}^{x}(\pi_{temp}), T_{II,3}^{x}(\pi_{temp}), T_{III,2}^{\frac{1}{2}}(\pi_{temp}), T_{IV,5}(\pi_{temp})$, and $T_{V,4}^{\pm}(\pi_{temp})$. We use the same list of tempered, corank $2$ representations that is provided at the beginning of Section \ref{classtempcorank3}. 

We begin with the case $T_{I,1}^{\alpha+1}(T_{I,1}^{\alpha}(\pi_{sc}))$. 
\begin{prop}\label{temp4B,1}
    Let $\pi_{temp} = T_{I,1}^{\alpha+1}(T_{I,1}^{\alpha}(\pi_{sc}))$ for $\alpha > 0$.
    \begin{enumerate}[label = (\roman*)]
    \item The representation $T_{I,2}^{x}(\pi_{temp})$ is well-defined if and only if $x = \frac{1}{2}$, and $\alpha \in \{\frac{1}{2}, \frac{3}{2}\}$. 
    \item When $\alpha \in \{\frac{1}{2}, \frac{3}{2}\}$, the representation $T_{I,2}^{\frac{1}{2}}(\pi_{temp})$ is of critical type. 
    \item When $\alpha = \frac{1}{2}$, we have that $\pi(\EE_{I,\frac{1}{2}}) = T_{I,2}^{\frac{1}{2}}(\pi_{temp})$, where 
    \[\EE_{I,\frac{1}{2}}= \scalebox{0.8}{\bordermatrix{
  & \frac{1}{2} & \frac{3}{2} \cr
  &\oplus \cr
  &\oplus \cr
  &&\oplus
}}.\]
When $\alpha = \frac{3}{2}$, we have that $\pi(\EE_{I,\frac{3}{2}}) = T_{I,2}^{\frac{1}{2}}(\pi_{temp})$, where 
\[\EE_{I,\frac{3}{2}}= \scalebox{0.8}{\bordermatrix{
  & \frac{1}{2} & \frac{3}{2} &\frac{5}{2}\cr
  &\oplus \cr
  &\oplus \cr
  &&&\ominus
}}.\]
\item The representation $T_{II,3}^{x}(\pi_{temp})$ is well-defined if and only if $x = \frac{1}{2}$, and $\alpha \in \frac{1}{2} + \mathbb{Z}_{>1}$. In this case, it is the same as the representation $T_{I,1}^{\alpha+1}(T_{I,1}^{\alpha}(T_{II,3}^{\frac{1}{2}}(\pi_{sc})))$, as described in Proposition \ref{temp4A,19}. 
\item The representation $T_{III,2}^{\frac{1}{2}}(\pi_{temp})$ is well-defined if and only if $\alpha \in \{\frac{1}{2}, \frac{3}{2}\}$. 
\item When $\alpha \in \{\frac{1}{2}, \frac{3}{2}\}$, the representation $T_{III,2}^{\frac{1}{2}}(\pi_{temp})$ is of critical type. 
\item  When $\alpha = \frac{1}{2}$, we have that $\pi(\EE_{III,\frac{1}{2}}) = T_{III,2}^{\frac{1}{2}}(\pi_{temp})$, where 
    \[\EE_{III,\frac{1}{2}}= \scalebox{0.8}{\bordermatrix{
  & \frac{1}{2} & \frac{3}{2} \cr
  &\ominus \cr
  &\ominus \cr
  &&\oplus
}}.\] 
When $\alpha = \frac{3}{2}$, we have that $\pi(\EE_{III,\frac{3}{2}}) = T_{I,2}^{\frac{1}{2}}(\pi_{temp})$, where 
\[\EE_{III,\frac{3}{2}}= \scalebox{0.8}{\bordermatrix{
  & \frac{1}{2} & \frac{3}{2} &\frac{5}{2}\cr
  &\ominus \cr
  &\ominus \cr
  &&&\ominus
}}.\]
        
        \item The representation $T_{IV,5}(\pi_{temp})$ is well-defined if and only if $\alpha \in \mathbb{Z}_{>1}$. In this case, the representation is the same as the representation $T_{I,1}^{\alpha +1}(T_{I,1}^{\alpha}(T_{IV,5}(\pi_{sc})))$, as described in Proposition \ref{temp4A,21}. 
        \item The representation $T_{V,4}^{\pm}(\pi_{temp})$ is well-defined if and only if $\alpha = 1$. In this case, the representation is the same as the representation $T_{V,4}^{\pm}(T_{I,1}^{1}(\pi_{sc}))$, as described in Proposition \ref{temp4A,23}. 
    \end{enumerate}
\end{prop}
\begin{proof}
    Define 
    \begin{flalign*}  
            \EE &:= \{([\alpha -2, \epsilon_\rho]_\rho,0,\eta),([\alpha+1, \alpha+1]_\rho,0,(-1)^{\alpha -1 - \epsilon_\rho} \eta)\}.
    \end{flalign*}
    Then we have $\pi(\EE) = \pi_{temp}$. Here is the associated symbol. 
     \[\EE= \scalebox{0.8}{\bordermatrix{
  & \epsilon_\rho&\cdots & \alpha -2 & \alpha -1  &\alpha & \alpha +1\cr
  &\odot & \cdots & \odot \cr
  &&&&&&\odot
}}.\] As shown, the extended multi-segment is multiplicity free. In order for $T_{I,2}^{x}$ to be well-defined, we need to have $x = \frac{1}{2}$, since by convention the multiplicity $m_{\phi}(\rho \otimes S_0) = \infty$. We also need that $m_{\phi}(\rho \otimes S_{2}) = 0$, which can only happen when $\alpha \in \{\frac{1}{2}, \frac{3}{2}\}$. This proves part $(i)$. Parts $(ii)$ and $(iii)$ follow from definition. Following a similar argument, we can prove parts $(v)$ to $(vii)$. 

Parts $(iv)$, $(viii)$ and $(ix)$ follow from comparing the extended multi-segment above to the one given in previous propositions. This completes the proof. 
\end{proof}

The second case is $\pi_{temp} = T_{I,1}^{\alpha-1}(T_{I,1}^{\alpha}(\pi_{sc}))$, defined for $\alpha > 1$. 

\begin{prop}\label{temp4B,2}
    Let $\pi_{temp} = T_{I,1}^{\alpha-1}(T_{I,1}^{\alpha}(\pi_{sc}))$ and $\alpha > 1$. 
    \begin{enumerate}[label = (\roman*)]
    \item The representation $T_{I,2}^{x}(\pi_{temp})$ is well-defined if and only if $(x,\alpha) = (\frac{1}{2}, \frac{5}{2})$. 
    \item When $\alpha = \frac{5}{2}$, the representation $T_{I,2}^{\frac{1}{2}}(\pi_{temp})$ is of critical type. 
    \item Let $\alpha = \frac{5}{2}$. Then we have $\pi(\EE_{I}) = T_{III,2}^{\frac{1}{2}}(\pi_{temp})$, where 
    \[\EE_{I}= \scalebox{0.8}{\bordermatrix{
  &\frac{1}{2} & \frac{3}{2} & \frac{5}{2} \cr
  &\oplus \cr
  &\oplus \cr
  &&\odot & \odot
}}.\]
\item The representation $T_{II,3}^{x}(\pi_{temp})$ is well-defined if and only if $x = \frac{1}{2}$ and $\alpha \in \frac{1}{2} + \mathbb{Z}_{>2}$. In this case, it is the same as the representation $T_{I,1}^{\alpha-1}(T_{I,1}^{\alpha}(T_{II,3}^{\frac{1}{2}}(\pi_{sc})))$, as defined in Proposition \ref{temp4A,19}. 
    \item The representation $T_{III,2}^{\frac{1}{2}}(\pi_{temp})$ is well-defined if and only if $\alpha =  \frac{5}{2}$. 
    \item When $\alpha = \frac{5}{2}$, the representation $T_{III,2}^{\frac{1}{2}}(\pi_{temp})$ is of critical type. 
    \item Let $\alpha = \frac{5}{2}$. Then we have $\pi(\EE_{III}) = T_{III,2}^{\frac{1}{2}}(\pi_{temp})$, where 
    \[\EE_{III}= \scalebox{0.8}{\bordermatrix{
  &\frac{1}{2} & \frac{3}{2} & \frac{5}{2} \cr
  &\ominus \cr
  &\ominus \cr
  &&\odot & \odot
}}.\]
        \item The representation $T_{IV,5}(\pi_{temp})$ is well-defined if and only if $\alpha \in \mathbb{Z}_{\geq 3}$. In this case, it's the same as the representation $T_{I,1}^{\alpha-1}(T_{I,1}^{\alpha}(T_{IV,5}(\pi_{sc})))$, as described in Proposition \ref{temp4A,21}. 
        \item The representation $T_{V,4}^{\pm}(\pi_{temp})$ is well-defined if and only if $\alpha = 2$. In this case, define 
        \begin{flalign*}  
            \EE_{\pm} &:= \{([0,0]_\rho,0,\pm1)^5,([2,1]_\rho,0,\eta)\}.
    \end{flalign*}
    Then we have $\pi(\EE_{\pm}) = T_{V,4}^{\pm}(\pi_{temp})$. Here are the associated symbols. 
    \[\EE_{+}= \scalebox{0.8}{\bordermatrix{
  & 0&1&2\cr
  &\oplus \cr
  &\oplus \cr
  &\oplus \cr
  &\oplus \cr
  &&\odot &\odot
}},\quad
\EE_{-}= \scalebox{0.8}{\bordermatrix{
  & 0&1&2\cr
  &\ominus \cr
  &\ominus \cr
  &\ominus \cr
  &\ominus \cr
  &&\odot &\odot
}}.\]
\item The representations $T_{I,2}^{x}(\pi_{temp}), T_{II,3}^{x}(\pi_{temp})$ are not well-defined. 
    \end{enumerate}
\end{prop}
\begin{proof}
    The proof is similar to that of  Proposition \ref{temp4B,1}, which we omit. 
\end{proof}
The next case to consider is $\pi_{temp} = T_{IV,3}(T_{I,1}^{\alpha}(\pi_{sc}))$, defined for $\alpha \in \mathbb{Z}_{>0}$. 
\begin{prop}\label{temp4B,3}
    Let  $\pi_{temp} = T_{IV,3}(T_{I,1}^{\alpha}(\pi_{sc}))$ and $\alpha \in \mathbb{Z}_{>0}$. 
    \begin{enumerate}[label = (\roman*)]
        \item The representation $T_{I,2}^{x}(\pi_{temp})$ is well-defined if and only if $(x,\alpha) = (1,2)$.
        \item When $\alpha = 2$, the representation $T_{I,2}^{1}(\pi_{temp})$ is of critical type. 
        \item Define 
        \begin{flalign*}  
            \EE &:= \{([0,0]_\rho,0,\eta)([1,1]_\rho,0,\eta)^2, ([2,2]_\rho,0,\eta)\}.
    \end{flalign*}
    Then $\pi(\EE)= T_{I,1}^{1}(\pi_{temp})$. Here is the associated symbol. 
    \[\EE= \scalebox{0.8}{\bordermatrix{
  & 0&1&2\cr
  &\odot\cr
  &&\odot \cr
  &&\odot \cr
  &&&\odot \cr
}}.\]
\item The representation $T_{II,3}^{x}(\pi_{temp})$ is well-defined if and only if $x = 1$. In this case, it is the same as the representation $T_{I,1}^{\alpha}(T_{II,3}^{1}(T_{IV,3}(\pi_{sc})))$, as described in Proposition \ref{temp4A,27}. 
\item The representations $T_{III,2}^{\frac{1}{2}}(\pi_{temp}), T_{IV,5}(\pi_{temp}), T_{V,4}^{\pm}(\pi_{temp})$ are not well-defined. 
    \end{enumerate}
\end{prop}
\begin{proof}
    Let $\phi$ be the $L$-parameter associated to $\pi_{temp}$. For $T_{I,2}^{x}(\pi_{temp})$ to be well-defined, we need $m_{\phi}(\rho \otimes S_{2x-1}) \geq 2$ and $S(\rho \otimes S_{2x+1}) = 0$. By examining the extended multi-segment, we see that this is only possible when $(x,\alpha) = (1,2)$. This proves part $(i)$. Parts $(ii), (iii)$ and $(v)$ follows from definition. 
    For part $(iv)$, we see that $T_{II,3}^{x}$ is only well-defined when $m_{\phi}(\rho\otimes S_{2x-1}) \geq 3$, which means that $x = 1$. The result follows by comparing the corresponding extended multi-segments. This concludes the proof. 
\end{proof}

We now move onto the case $\pi_{temp} = T_{V,2}^{\pm}(T_{I,1}^{1}(\pi_{sc}))$, which is well-defined only for $\alpha = 1$. 
\begin{prop}\label{temp4B,4}
    Let $\pi_{temp} = T_{V,2}^{\pm}(T_{I,1}^{1}(\pi_{sc}))$ and $\alpha = 1$, then the representations $T_{I,2}^{x}(\pi_{temp})$ \\
    $, T_{II,3}^{x}(\pi_{temp}), T_{III,2}^{\frac{1}{2}}(\pi_{temp}), T_{IV,5}(\pi_{temp})$ and $T_{V,4}^{\pm}(\pi_{temp})$ are all not well-defined. 
\end{prop}
\begin{proof}
    Let $\phi$ be the $L$-parameter corresponding to $\pi_{temp}$, then we have that $m_{\phi}(\rho \otimes S_1) = 2, m_{\phi}(\rho \otimes S_3) = 1$, and $m_{\phi}(\rho \otimes S_{a}) = 0$ otherwise for $a \in \mathbb{Z}$. This proves that 
    \[T_{I,2}^{x}(\pi_{temp}), T_{II,3}^{x}(\pi_{temp}), T_{IV,5}(\pi_{temp}), T_{V,2}^{\pm}(\pi_{temp})\]
    are not well-defined. 

    Finally, $T_{III,2}^{\frac{1}{2}}(\pi_{temp})$ is not well-defined by the good parity condition. 
\end{proof}

A similar case to the one above is $\pi_{temp} = T_{I,1}^{1}(T_{V,2}^{\pm}(\pi_{sc}))$, which is well-defined if and only if $\alpha = 0$. The proofs of Propositions \ref{temp4B,5} to \ref{temp4B,9} below are similar to that of Propositions \ref{temp4B,3} and \ref{temp4B,4}, which we omit. 
\begin{prop}\label{temp4B,5}
    Let $\pi_{temp} = T_{I,1}^{1}(T_{V,2}^{\pm}(\pi_{sc}))$ and $\alpha =0$. 
    \begin{enumerate}[label = (\roman*)]
        \item The representation $T_{IV,5}(\pi_{temp})$ is well-defined and of critical type. It is the same as the representation $T_{I,1}^{1}(T_{V,6}^{\pm}(\pi_{sc}))$, as described in Proposition \ref{temp4A,32}. 
        \item The representations $T_{I,2}^{x}(\pi_{temp}), T_{II,3}^{x}(\pi_{temp}), T_{III,2}^{\frac{1}{2}}(\pi_{temp})$ and $T_{V,4}^{\pm}(\pi_{temp})$ are not well-defined. 
    \end{enumerate}
\end{prop}

The next case we have is $\pi_{temp} = T_{I,2}^{\frac{1}{2}}(\pi_{sc})$, which is defined only when $\alpha = \frac{1}{2}$. 
\begin{prop}\label{temp4B,6}
    Let $\pi_{temp} = T_{I,2}^{\frac{1}{2}}(\pi_{sc})$ and $\alpha = \frac{1}{2}$. 
    \begin{enumerate}[label = (\roman*)]
        \item The representation $T_{I,2}^{x}(\pi_{temp})$ is well-defined if and only if $x = \frac{3}{2}$. 
        \item The representation $T_{I,2}^{\frac{3}{2}}(\pi_{temp})$ is of critical type. 
        \item We have $\pi(\EE) = T_{I,2}^{\frac{3}{2}}(\pi_{temp})$, where 
        \[\EE= \scalebox{0.8}{\bordermatrix{
  & \frac{3}{2}\cr
  &\oplus\cr
  &\oplus 
}}.\]
\item The representations $T_{II,3}^{x}(\pi_{temp}), T_{III,2}^{\frac{1}{2}}(\pi_{temp}), T_{IV,5}(\pi_{temp})$, and $T_{V,4}^{\pm}(\pi_{temp})$ are all not well-defined. 
    \end{enumerate}
\end{prop}

The next case to consider is $\pi_{temp} = T_{II,3}^{\frac{1}{2}}(\pi_{sc})$, which is well-defined only when $\alpha \in \frac{1}{2} + \mathbb{Z}_{>0}$. 
\begin{prop}\label{temp4B,7}
    Let $\pi_{temp} = T_{II,3}^{\frac{1}{2}}(\pi_{sc})$, and $\alpha \in \frac{1}{2} + \mathbb{Z}_{>0}$. 
    \begin{enumerate}[label = (\roman*)]
        \item The representation $T_{II,3}^{x}(\pi_{temp})$ is well-defined if and only if $x = \frac{3}{2}$. 
        \item The representation $T_{II,3}^{\frac{3}{2}}(\pi_{temp})$ is of critical type. 
        \item Define 
        \begin{flalign*}  
            \EE &:= \{([\alpha -1, \frac{1}{2}]_\rho,0,-1),[\frac{3}{2}, \frac{3}{2}]_\rho,0,1)^2\}.
        \end{flalign*}
    Then $\pi(\EE)= T_{II,3}^{\frac{1}{2}}(\pi_{temp})$. Here is the associated symbol. 
    \[\EE= \scalebox{0.8}{\bordermatrix{
  & \frac{1}{2} &\frac{3}{2}&\cdots & \alpha - 1\cr
  &\ominus &\oplus & \cdots & \odot\cr
  &&\oplus \cr
  &&\oplus \cr
}}.\]
\item The representations $T_{I,2}^{x}(\pi_{temp}), T_{III,2}^{\frac{1}{2}}(\pi_{temp}), T_{IV,5}(\pi_{temp})$ and $T_{V,4}^{\pm}(\pi_{temp})$ are all not well-defined. 
    \end{enumerate}
\end{prop}

There are three more cases to examine in Case $(B)$. The next one is $\pi_{temp} = T_{III,2}^{\frac{1}{2}}(\pi_{sc})$, which is well-defined only when $\alpha = \frac{1}{2}$. 
\begin{prop}\label{temp4B,8}
    Let $\pi_{temp} = T_{III,2}^{\frac{1}{2}}(\pi_{sc})$ and $\alpha = \frac{1}{2}$. 
    \begin{enumerate}[label = (\roman*)]
        \item The representation $T_{I,2}^{x}(\pi_{temp})$ is well-defined if and only if $x = \frac{3}{2}$. 
        \item The representation $T_{I,2}^{\frac{3}{2}}(\pi_{temp})$ is of critical type. 
        \item We have that $\pi(\EE) = T_{I,2}^{\frac{3}{2}}(\pi_{temp})$, where 
        \[\EE= \scalebox{0.8}{\bordermatrix{
  & \frac{3}{2}\cr
  &\ominus\cr
  &\ominus 
}}.\]
\item The representations $T_{II,3}^{x}(\pi_{temp}), T_{III,2}^{\frac{1}{2}}(\pi_{temp}), T_{IV,5}(\pi_{temp})$ and $T_{V,4}^{\pm}(\pi_{temp})$ are all not well-defined. 
    \end{enumerate}
\end{prop}

The second last case in Case $(B)$ is $\pi_{temp} = T_{IV,5}(\pi_{sc})$, which is well-defined if and only if $\alpha \in \mathbb{Z}_{>0}$. 
\begin{prop}\label{temp4B,9}
    Let $\pi_{temp}= T_{IV,5}(\pi_{sc})$, and $\alpha \in \mathbb{Z}_{>0}$. 
    \begin{enumerate}[label = (\roman*)]
        \item The representation $T_{I,2}^{x}(\pi_{temp})$ is well-defined if and only if $(x,\alpha) = (1,1)$. In this case, $T_{I,2}^{1}(\pi_{temp})$ is of critical type and is the same as the representation $T_{I,2}^{1}(T_{IV,3}(\pi_{sc}))$, as described in Proposition \ref{temp4A,26}. 
        \item The representation $T_{II,3}^{x}(\pi_{temp})$ is well-defined if and only if $x = 1$, and $\alpha \in \mathbb{Z}_{>1}$. In this case, the representation $T_{II,3}^{1}(\pi_{temp})$ is not of critical type, and it is the same as the representation $T_{IV,3}(T_{II,3}^{1}(T_{IV,3}(\pi_{sc})))$, as described in Proposition \ref{temp4A,28}. 
        \item The representations $T_{III,2}^{\frac{1}{2}}(\pi_{temp}), T_{IV,5}(\pi_{temp})$, and $T_{V,4}^{\pm}(\pi_{temp})$ are all not well-defined. 
    \end{enumerate}
\end{prop}

The final case in Case $(B)$ is $\pi_{temp} = T_{V,4}^{\pm}(\pi_{sc})$, which is well-defined if and only if $\alpha = 0$. 
\begin{prop}\label{temp4B,10}
    Let $\pi_{temp}^{\pm} = T_{V,4}^{\pm}(\pi_{sc})$ and $\alpha = 0$. 
    \begin{enumerate}[label = (\roman*)]
        \item The representations $T_{I,2}^{x}(\pi_{temp}^{\pm})$ are well-defined if and only if $x = 1$. 
        \item The representations $T_{I,2}^{1}(\pi_{temp}^{\pm})$ are of critical type. 
        \item  Define 
        \begin{flalign*}  
            \EE_{\pm} &:= \{([0,0]_\rho,0,\pm1)^2,([1,1]_\rho,0,\pm1)^2\}.
        \end{flalign*}
    Then $\pi(\EE_{\pm})=T_{I,2}^{1}(\pi_{temp}^{\pm})$. Here are the associated symbols. 
    \[\EE_{+}= \scalebox{0.8}{\bordermatrix{
  & 0& 1\cr
  &\oplus \cr
  &\oplus \cr
  &&\oplus \cr
  &&\oplus \cr
}},\quad
\EE_{-}= \scalebox{0.8}{\bordermatrix{
  & 0& 1\cr
  &\ominus \cr
  &\ominus \cr
  &&\ominus \cr
  &&\ominus \cr
}}.\]
    \end{enumerate}
\end{prop}

This concludes our discussion of Case $(B)$. We now move onto Case $(C)$. 
\subsection{\texorpdfstring{Case $(C): \pi \hookrightarrow \rho\lvert \cdot \rvert^{x_1} \times \rho\lvert \cdot \rvert^{x_1} \times \rho\lvert \cdot \rvert^{x_3} \rtimes \pi_{temp}$, where $\pi_{temp}$ is tempered of corank $1$}{}}
In this subsection, we consider representations $\pi$ of the form 
\begin{equation*}
    T_{I,3}^{x}(\pi_{temp}), T_{III,2}^{1}(\pi_{temp}), T_{IV,7}(\pi_{temp}),T_{V,6}^{\pm}(\pi_{temp})
\end{equation*} where $\pi_{temp}$ is tempered of corank $1$. Since there are three tempered representations of corank $1$, namely $T_{I,1}^{\alpha}(\pi_{sc}), T_{IV,5}(\pi_{sc})$ and $T_{V,2}^{\pm}(\pi_{sc})$. There are $12$ total cases to consider. We begin with $\pi_{temp} = T_{I,1}^{\alpha}(\pi_{sc})$, which is defined if and only when $\alpha > 0$.

\begin{prop}\label{temp4C,1}
    Let $\pi_{temp} = T_{I,1}^{\alpha}(\pi_{sc})$, for $\alpha > 0$. 
    \begin{enumerate}[label = (\roman*)]
        \item The representation $T_{III,2}^{1}(\pi_{temp})$ is well-defined if and only if $\alpha = 2$. When $\alpha = 2$, the representation $T_{III,2}^{1}(\pi_{temp})$ is of critical type. It is the same as the representation $T_{I,1}^{1}(T_{IV,3}(T_{I,1}^{2}(\pi_{sc})))$, as described in Proposition \ref{temp4B,3}. 
        \item The representation $T_{IV,7}(\pi_{temp})$ is well-defined if and only if $\alpha \in \mathbb{Z}_{>0}$. It is the same as the representation $T_{I,1}^{\alpha}(T_{IV,7}(\pi_{sc}))$, as described in Proposition \ref{temp4A,31}, and it is of critical type only when $\alpha = 1$. 
        \item The representations $T_{I,3}^{x}(\pi_{temp}),  T_{V,6}^{\pm}(\pi_{temp})$ are not well-defined. 
    \end{enumerate}
\end{prop}
\begin{proof}
    Parts $(i)$ and $(ii)$ follows from definition, and by comparing the resulting extended multi-segments to the ones given in previous propositions. Part $(iii)$ follows from the fact that $\alpha > 0$, and $\phi$ is multiplicity free, where $\phi$ is the $L$-parameter corresponding to $\pi_{temp}$. 
\end{proof}

The second tempered representation of corank $1$ is $\pi_{temp} = T_{IV,3}(\pi_{sc})$, which is well-defined only when $\alpha \in \mathbb{Z}_{>0}$. 
\begin{prop}\label{temp4C,2}
    Let $\pi_{temp} = T_{IV,3}(\pi_{sc})$, and $\alpha \in \mathbb{Z}_{>0}$. 
    \begin{enumerate}[label = (\roman*)]
        \item The representation $T_{I,3}^{x}(\pi_{temp})$ is well-defined if and only if $(x,\alpha) = (1,1)$. 
        \item When $\alpha = 1$, the representation $T_{I,3}^{1}(\pi_{temp})$ is of critical type. 
        \item When $\alpha = 1$, we have $\pi(\EE) = T_{I,3}^{1}(\pi_{temp})$, where 
        \[\EE_= \scalebox{0.8}{\bordermatrix{
  &  1\cr
  &\odot \cr
  &\odot \cr
  &\odot \cr
}}.\]
 \item The representation $T_{III,2}^{1}(\pi_{temp})$ is well-defined if and only if $\alpha = 1$. In this case, it is of critical type, and is the same as the representation $T_{IV,3}(T_{I,2}^{1}(T_{IV,3}(\pi_{sc})))$, as described in Proposition \ref{temp4A,26}. 
    \item The representations $T_{IV,7}(\pi_{temp})$ and $T_{V,6}^{\pm}(\pi_{temp})$ are not well-defined. 
    \end{enumerate}
\end{prop}
\begin{proof}
    Let $\phi$ be the $L$-parameter corresponding to $\pi_{temp}$. Then $m_{\phi}(\rho \otimes S_{2x+1}) \geq 3$ only when $x = 0$. If $T_{I,3}^{x}(\pi_{temp})$ is well-defined then we also need $m_{\phi}(\rho \otimes S_{2x+1}) = 0$, which is only possible when $\alpha = 1$. This proves part $(i)$. Parts $(ii)$ and $(iii)$ follows from definition. 

    For $T_{III,2}^{1}(\pi_{temp})$ to be well-defined, we need $m_{\phi}(\rho \otimes S_{3})= 0$, which can only happen when $\alpha =1$. By comparing the corresponding extended multi-segments in this case, we prove $(iv)$. Part $(v)$ follows from definition and the good parity condition. 
\end{proof}

The last tempered representation of corank $1$ is $\pi_{temp} = T_{V,2}^{\pm}(\pi_{sc})$, which is well-defined if and only if $\alpha = 0$. 
\begin{prop}\label{temp4C,3}
    Let $\pi_{temp}^{\pm} = T_{V,2}^{\pm}(\pi_{sc})$ and $\alpha = 0$. 
    \begin{enumerate}[label = (\roman*)]
        \item The representations $T_{III,2}^{1}(\pi_{temp}^{\pm})$ are well-defined and of critical type. 
        \item Define 
        \begin{flalign*}  
            \EE_{\pm} &:= \{([0,0]_\rho,0,\pm1)^2,([1,1]_\rho,0,\mp1)^2\}.
        \end{flalign*}
        Then we have that $\pi(\EE_{\pm}) = T_{III,2}^{1}(\pi_{temp}^{\pm})$, where 
        \[\EE_{+}= \scalebox{0.8}{\bordermatrix{
  &  0 & 1\cr
  &\oplus \cr
  &\oplus \cr
  &&\ominus \cr
  &&\ominus \cr
}}, \quad \EE_{-}= \scalebox{0.8}{\bordermatrix{
  &  0 & 1\cr
  &\ominus \cr
  &\ominus \cr
  &&\oplus \cr
  &&\oplus \cr
}}.  \]
\item The representations $T_{I,3}^{x}(\pi_{temp}), T_{IV,7}(\pi_{temp})$ and $T_{V,6}^{\pm}(\pi_{temp})$ are not well-defined. 
    \end{enumerate}
\end{prop}
\begin{proof}
    The proof is similar to that of Proposition \ref{temp4C,2}, which we omit.  
\end{proof}

This concludes Case $(C)$. We now move onto the final case in our classification of tempered representations of corank $4$. 
\subsection{\texorpdfstring{Case $(D): \pi \hookrightarrow \rho\lvert \cdot \rvert^{x_1} \times \rho\lvert \cdot \rvert^{x_1} \times \rho\lvert \cdot \rvert^{x_3} \times \rho\lvert \cdot \rvert^{x_4} \rtimes\pi_{sc}$, where $\pi_{sc}$ is supercuspidal}{}}

By Theorem \ref{thm temp algo}, there are $6$ total possibilities to consider in this case, as listed in the beginning of the section. Let $\pi_{sc} = \pi(\phi, \epsilon)$, then $\phi$ must be multiplicity free.

The first case is $T_{I,4}^{x}(\pi_{sc})$. 
\begin{prop}\label{temp4D,ex1}
    Let $\pi_{sc} = \pi(\phi, \epsilon)$ and $\alpha = \alpha_{\rho, \epsilon}$. 
    \begin{enumerate}[label = (\roman*)]
        \item The representation $T_{I,4}^{x}(\pi_{sc})$ is well-defined if and only if $(x,\alpha) = (\frac{1}{2}, \frac{1}{2})$. 
        \item When $\alpha = \frac{1}{2}$, the representation $T_{I,4}^{x}(\pi_{sc})$ is of critical type. 
        \item When $\alpha = \frac{1}{2}$, we have that $\pi(\EE) = T_{I,4}^{\frac{1}{2}}(\pi_{sc})$, where 
        \[\EE= \scalebox{0.8}{\bordermatrix{
  &  \frac{1}{2} \cr
  &\oplus \cr
  &\oplus \cr
  &\oplus \cr
  &\oplus \cr
}}.\]
    \end{enumerate}
\end{prop}
\begin{proof}
    $T_{I,4}^{x}(\pi_{temp})$ is only well-defined when $m_{\phi}(\rho \otimes S_{2x -1}) \geq 4$. Since $\phi$ is multiplicity free, this can only occur when $x = \frac{1}{2}$. Additionally we need that $m_{\phi}(\rho \otimes S_{2x +1}) = 0$, which means that $\alpha$ must be $\frac{1}{2}$. This proves part $(i)$. Parts $(ii)$ and $(iii)$ follow from definition. 
\end{proof}

The second case is $\pi_{temp} = T_{II,5}^{x}(\pi_{sc})$. The proofs of Propositions \ref{temp4D,ex2} to \ref{temp4D,3} below are similar to that of Proposition \ref{temp4D,ex1}, which we omit. 

\begin{prop}\label{temp4D,ex2}
     Let $\pi_{sc} = \pi(\phi, \epsilon)$ and $\alpha = \alpha_{\rho, \epsilon}$. 
     \begin{enumerate}[label = (\roman*)]
         \item The representation $T_{II,5}^{x}(\pi_{sc})$ is well-defined if and only if $x = \frac{1}{2}$ and $\alpha \in \frac{1}{2} + \mathbb{Z}_{>0}$. 
         \item When $\alpha \in \frac{1}{2} + \mathbb{Z}_{>0}$, the representation $T_{II,5}^{\frac{1}{2}}(\pi_{sc})$ is not of critical type. 
         \item When $\alpha \in \frac{1}{2} + \mathbb{Z}_{>0}$, define 
         \begin{flalign*}  
            \EE &:= \{([\frac{1}{2}, \frac{1}{2}]_\rho,0,-1)^4,([\alpha -1, \frac{1}{2}]_\rho, 0, -1)\}.
        \end{flalign*}
        Then we have $\pi(\EE) = T_{II,5}^{\frac{1}{2}}(\pi_{sc})$. Here is the associated symbol. 
        \[\EE= \scalebox{0.8}{\bordermatrix{
  &  \frac{1}{2} & \cdots & \alpha-1 \cr
  &\ominus \cr
  &\ominus \cr
  &\ominus \cr
  &\ominus \cr
  &\ominus & \cdots & \odot \cr
}}.\]
     \end{enumerate}
\end{prop}

The third case is $\pi_{temp} = T_{III,2}^{\frac{3}{2}}(\pi_{sc})$. 

\begin{prop}\label{temp4D,1}
    \begin{enumerate}[label = (\roman*)]
    Let $\pi_{sc} = \pi(\phi, \epsilon)$ and $\alpha = \alpha_{\rho, \epsilon}$. 
        \item The representation $T_{III,2}^{\frac{3}{2}}(\pi_{sc})$ is well-defined if and only if $\alpha = \frac{3}{2}$. 
        \item When $\alpha = \frac{3}{2}$, the representation $T_{III,2}^{\frac{3}{2}}(\pi_{temp})$ is of critical type. 
        \item Define 
        \begin{flalign*}  
            \EE &:= \{([\frac{1}{2}, \frac{1}{2}]_\rho,0,\eta)^2,([\frac{3}{2}, \frac{3}{2}]_\rho, 0, -\eta)^2\}.
        \end{flalign*}
        Then we have $\pi(\EE) = T_{III,2}^{\frac{3}{2}}(\pi_{sc})$. Here is the associated symbol. 
        \[\EE= \scalebox{0.8}{\bordermatrix{
  &  \frac{1}{2} & \frac{3}{2}\cr
  &\odot \cr
  &&\odot \cr
  &&\odot 
}}.\]
    \end{enumerate}
\end{prop}

We move onto the case $T_{III,4}^{\frac{1}{2}}(\pi_{sc})$. 
\begin{prop}\label{temp4D,ex3}
    Let $\pi_{sc} = \pi(\phi, \epsilon)$ and $\alpha = \alpha_{\rho, \epsilon}$. 
    \begin{enumerate}
        \item The representation $T_{III,4}^{\frac{1}{2}}(\pi_{sc})$ is well-defined if and only if $\alpha = \frac{1}{2}$. 
        \item When $\alpha = \frac{1}{2}$, the representation $T_{III,4}^{\frac{1}{2}}(\pi_{sc})$ is of critical type. 
        \item We have that $\pi(\EE) = T_{III,4}^{\frac{1}{2}}(\pi_{sc})$, where 
        \[\EE= \scalebox{0.8}{\bordermatrix{
  &  \frac{1}{2} \cr
  &\ominus \cr
  &\ominus \cr
  &\ominus \cr
  &\ominus \cr
}}.\]
    \end{enumerate}
\end{prop}

The next case to consider is the representation $\pi_{temp} = T_{IV,9}(\pi_{sc})$. 
\begin{prop}\label{temp4D,2}
    Let $\pi_{sc} = \pi(\phi, \epsilon)$  and $\alpha = \alpha_{\rho, \epsilon}$.
    \begin{enumerate}[label = (\roman*)]
        \item The representation $T_{IV,9}(\pi_{sc})$ is well-defined if and only if $\alpha \in \mathbb{Z}_{>0}$. 
        \item The representation $T_{IV,9}(\pi_{sc})$ is not of critical type. 
        \item      Define 
        \begin{flalign*}  
            \EE &:= \{([0,0]_\rho,0,\eta)^8,([\alpha -1, 0]_\rho, 0, \eta)\}.
        \end{flalign*}
        Then we have $\pi(\EE) = T_{IV,9}(\pi_{sc})$. Here is the associated symbol. 
        \[\EE= \scalebox{0.8}{\bordermatrix{
  &  0 & \cdots & \alpha -1\cr
  &\odot \cr
  &\odot \cr
  &\odot \cr
  &\odot \cr
  &\odot \cr
  &\odot \cr
  &\odot \cr
  &\odot \cr
  &\odot & \cdots &\odot \cr
}}.\]
    \end{enumerate}
\end{prop}

The final case to consider is the representation 
$\pi_{temp} = T_{V,8}^{\pm}(\pi_{sc})$. 
\begin{prop}\label{temp4D,3}
    Let $\pi_{sc} = \pi(\phi, \epsilon)$  and $\alpha = \alpha_{\rho, \epsilon}$.
    \begin{enumerate}[label = (\roman*)]
        \item The representation $\pi_{temp}^{\pm} = T_{V,8}^{\pm}(\pi_{sc})$ is well-defined if and only if $\alpha = 0$.
        \item When $\alpha = 0$, the representation $\pi_{temp}^{\pm}$ is of critical type. 
        \item We have that $\pi(\EE_{\pm}) = \pi_{temp}^{\pm}$, where 
        \[\EE_{+}= \scalebox{0.8}{\bordermatrix{
  &  0 \cr
  &\oplus \cr
  &\oplus \cr
  &\oplus \cr
  &\oplus \cr
  &\oplus \cr
  &\oplus \cr
  &\oplus \cr
  &\oplus \cr
}}, \quad \EE_{-}= \scalebox{0.8}{\bordermatrix{
  &  0\cr
  &\ominus \cr
  &\ominus \cr
  &\ominus \cr
  &\ominus \cr
  &\ominus \cr
  &\ominus \cr
  &\ominus \cr
  &\ominus \cr
}}.\]
    \end{enumerate}
\end{prop}

With this we have classified all the corank 4 tempered representations of good parity. Combined with the non-tempered representations of corank $4$ which we classified in \S\ref{classnontempcorank4,1} to \S\ref{classnontempcorank4,34}, we can give a complete list of all representations of corank $4$, which are both of Arthur type and of critical type. This will be given in the next section.

By Theorem \ref{unitiffArthur}, we can conclude that that all of the representations listed in Proposition \ref{crnk4Artcritlist} are unitarizable. Furthermore, the list contains all representations of $G(F)$ of corank $4$ that are unitarizable and of critical type. 

\section{\texorpdfstring{Open connected components in the unitary dual  of corank $4$}{}} \label{opencnncomponents}

In this section, we use the lists given in Appendices \S\ref{artlist} and \S\ref{nonartlist} to give the full list of open unitary connected components in corank $4$. This will be the first step to constructing the full unitary dual. To this end, we  use the algorithm to generate the unitary dual candidate set  $\Pi_{\overline{A}}^{\lim}(G_n)$ introduced in \cite[\S 8]{HJLLZ24} to determine inductively the unitarity of a given open connected component. 

Using the technique of unitary reduction in Step $3$ of Algorithm \ref{alg A bar} below, we later on prove that certain connected components are non-unitary. In fact, we highlight that there are five main methods to prove that a certain connected component in $C \subset \R^n$ is non-unitary. Four of them are given in steps $(3$-$1)$ to $(3$-$4)$ of Algorithm \ref{alg A bar}, respectively. The last method is to show directly that there exists a point on the boundary of $C$ which contains non-unitarizable subquotients.

The steps to construct the full unitary dual of corank $4$ is as follows: First, in Proposition \ref{connectedisunitary} below, we give a list of open connected components, which we prove to be unitary. Then, in subsections \ref{ssregularcomponents} to \ref{levelsection}, we show that all other connected components are non-unitary, using Algorithm \ref{alg A bar} together with Tadi{\'c}'s techniques of dimensionality reduction in \cite{Tad23}. Explicitly, we characterize all possible connected components of dimension $3$ that can possibly appear on the boundary of an open unitary connected component of dimension $4$ (Propositions  \ref{finalslavregunit} and \ref{finallevelunit}). This provides the necessary conditions for the open unitary connected components to  satisfy, which enable us to show that the list of open connected components in Proposition \ref{connectedisunitary} is complete. 

Later, in \S\ref{1-parameterseries} and \S\ref{2-parameterseries}, we append to the unitary dual all representations which appear as part of a lower-dimensional unitarizable family (i.e., those with one or two parameters), as described in Step $2$ of Algorithm \ref{alg A bar}. This will give us the list of all possible unitarizable representations of corank $4$.

\subsection{\texorpdfstring{Algorithm for generating $\Pi_{\overline{A}}^{\lim}(G_n)$}{}}
\begin{algo}  [{\cite[Algorithm $8.5$]{HJLLZ24}}] \label{alg A bar}
   Let $\pi_{sc}$ be an irreducible supercuspidal representation of $G_n$, $\rho$ be an irreducible self-dual supercuspidal representation of $\GL_d(F)$, and let $r \in \Z_{\geq 0}$. In this algorithm, we output the set
\[ \Pi_{\overline{A}}(\pi_{sc}, \rho, r):= \Pi_{\overline{A}}(G_{n+rd}) \cap \Irr(X_{\rho}; \pi_{sc}).\]
Write $\alpha=\alpha_{\pi_{sc},\rho}$ for short.
\begin{enumerate}
    \item [Step 1.] Compute $\Pi_{A, \, gp}(\pi_{sc},\rho,s)$ for $0 \leq s \leq r$. Let $\Omega_{A,gp}:= \cup_{0 \leq s \leq r} \Pi_{A, \, gp}(\pi_{sc},\rho,s).$ 
    \item [Step 2.] For $ 0 \leq s \leq r' \leq r$, initiate $R(\pi_{sc},\rho,r',s)$ to be the empty set. Repeat the following steps for each $\pi_A \in \Pi_{A, \, gp}(\pi_{sc},\rho,s)$: 
    \begin{enumerate}
        \item [(2-1)] Let $\Psi(r'-s)$ be the collection of ordered tuples of pairs $\psi=((a_1,b_1),\ldots, (a_{l(\psi)}, b_{l(\psi)}))\in(\mathbb{Z}_{\geq 1}^2)^{l(\psi)}$, ${l(\psi)} \in \mathbb{Z}_{>0}$, such that $\sum_{i=1}^{l(\psi)} a_ib_i=r'-s$. By convention, let $\Psi(0):=\{\emptyset\}$. 
        
        \item [(2-2)] Suppose that $r'-s>0$. For each $\psi \in \Psi(r'-s)$, let $\mathcal{H}(\psi, \pi_A)$ denote the set consisting of the following reducibility hyperplanes in $\R^{{l(\psi)}}$: 
        \begin{itemize}
           
            \item $\{x_i= t \ | \ u_{\rho}(a_i,b_i)\lvert\cdot\rvert^{t} \rtimes \pi_A \text{ is reducible}\}$,
            \item $\{x_i\pm x_j = t\ | \ u_{\rho}(a_i,b_i)\lvert\cdot\rvert^{t} \times u_{\rho}(a_j,b_j) \text{ is reducible}, 1 \leq i < j \leq {l(\psi)}\}$.
        \end{itemize}
       
        \item [(2-3)] If $r'-s>0$, let $R(\psi,\pi_A)$ denote the (finite) set of connected components of $\R^{l(\psi)} \setminus (\cup_{H \in \mathcal{H}(\psi, \pi_A)} H)$.        
        For each $C \in R(\psi,\pi_A)$, append the triple $(C, \psi, \pi_A)$ into $R(\pi_{sc},\rho,r',s)$. If $r'-s=0$, then append $\{ (\R^0, \emptyset , \pi_A) \}$ into $R(\pi_{sc},\rho,r',s)$. 
    \end{enumerate}
    \item [Step 3.] Let $R:= \cup_{0 \leq s \leq r' \leq r} R(\pi_{sc},\rho,r',s) \sqcup \{-1\}$. Define an equivalence relation $\sim$ on $R$ as follows. Let $(C,\psi, \pi_A) \in R(\pi_{sc},\rho,r',s)$ where $C \subseteq \R^{l(\psi)}$.
    \begin{itemize}
    \item [(3-1)] Suppose that $C$ is unbounded. Then we define $ (C,\psi, \pi_A) \sim -1$.
    \item [(3-2)]Suppose that $C \cap \{x_i=0\} \neq \emptyset$. Define $\psi^{-}$ by removing $(a_i,b_i)$ from $\psi$. Take any point $\underline{y} \in C \cap \{x_i=0\}$ and define $\underline{y}^{-}\in \R^{l(\psi)-1}$ by removing the $i$-th coordinate $($if ${l(\psi)}=1$, then set $\underline{y}^-=0$$)$. Let $(C^-, \psi^{-}, \pi_A)$ be the unique element in $R$ such that $\underline{y}^- \in C^- $. Then we define $(C,\psi,\pi_A) \sim (C^-, \psi^{-}, \pi_A)$.
    \item [(3-3)] Suppose that $C \cap \{x_i=t\} \neq \emptyset$ for some $t\in (\alpha+ \half{a_i+b_i})+\Z$. Let $\pi_{A}^+:= u_{\rho}(a_i,b_i)\lvert\cdot\rvert^{t} \rtimes \pi_A$, which is irreducible and of good parity. If $\pi_A^+ $ is not in $\Omega_{A,gp}$, then define $(C,\psi, \pi_A) \sim -1$. If $\pi_A^+ \in \Omega_{A,gp}$, then  define $\psi^{-}$ and take a point $\underline{y}^- \in \R^{{l(\psi)}-1}$ as in the previous case. Let $(C^-, \psi^{-}, \pi_A^+)$ be the unique element in $R$ such that $\underline{y}^- \in C^- $. Then we define $(C,\psi,\pi_A) \sim (C^-, \psi^{-}, \pi_A^+)$.
    \item [(3-4)] Suppose that $(a_i,b_i)=(a_j,b_j)$ with $i \neq j$ and $C \cap \{x_i={ \pm}x_j\}\neq \emptyset$. Define $\psi^-$ by removing both $(a_i,b_i)$ and $(a_j,b_j)$ from $\psi$. Take any point $\underline{y}=(y_1,\ldots, y_{l}) \in C  \cap \{x_i={ \pm}x_j\}$. If $|y_i|>\half{1}$, then define $(C,\psi, \pi_A) \sim -1$. If $|y_i|< 1/2$, then define $\underline{y}^{-} \in \R^{{l(\psi)}-2}$ by removing the $i$-th and $j$-th coordinates $($if ${l(\psi)}=2$, then set $\underline{y}^-=0$$)$. Let $(C^-, \psi^{-}, \pi_A)$ be the unique element in $R$ such that $\underline{y}^- \in C^- $. Then we define $(C,\psi,\pi_A) \sim (C^-, \psi^{-}, \pi_A)$.
    \end{itemize}
     Let $R_{\overline{A}}$ be the collection of $(C,\psi, \pi_A) \in R$ such that $(C, \psi, \pi_A) \sim ( \R^0, \emptyset, \pi_A')$ for some $\pi_A' \in \Omega_{A,gp}$. 
    \item [Step 4.] For each $(C,\psi=((a_1,b_1),\ldots,(a_{l(\psi)},b_{l(\psi)})),\pi_A) \in R$, let 
   \[ \Pi (C, \psi, \pi_A):= \left\{ \bigtimes_{i=1}^{l(\psi)} u_{\rho}(a_i,b_i)\lvert\cdot\rvert^{y_i} \rtimes \pi_A \ | \ (y_1,\ldots, y_{l(\psi)})\in C\right\}, \] 
    if $l(\psi)>0,$ and $\Pi(\R^0, \emptyset , \pi_A):= \{\pi_A\}$. Then 
    \[ \Pi_{\overline{A}}(\pi_{sc}, \rho, r)= \left(\bigcup_{(C, \psi, \pi_A) \in R_{\overline{A}}} \Pi (C, \psi, \pi_A)\right) \cap \Pi(G_{n+rd}).\]
\end{enumerate}

Then, we can describe the set $\Pi_{\overline{A}}^{\lim}(\pi_{sc}, \rho, r):= \Pi_{\overline{A}}^{\lim}(G_{n+rd})\cap \Irr(X_{\rho};\pi_{sc})$ abstractly from $\Pi_{\overline{A}}(\pi_{sc}, \rho, r)$ as follows. For each $(C,\psi,\pi_A)\in R_{\overline{A}}$, define
\[ \Pi^{\lim} ({C}, \psi, \pi_A):= \left\{ \pi \ \middle| \ \pi \leq \bigtimes_{i=1}^{l(\psi)} u_{\rho}(a_i,b_i)\lvert\cdot\rvert^{y_i} \rtimes \pi_A \text{ for some } (y_1,\ldots, y_{l(\psi)})\in \overline{C}\right\}. \]
Then by definition,
\[  \Pi_{\overline{A}}^{\lim}(\pi_{sc}, \rho, r)= \left(\bigcup_{(C, \psi, \pi_A) \in R_{\overline{A}}} {\Pi^{\lim} ({C}, \psi, \pi_A)}\right) \cap \Pi(G_{n+rd}).\]
\end{algo}
\begin{remark}
    Note that to compute the set $\Pi_{A,gp}(\pi_{sc}, \rho, s)$, we repeat the steps from Sections \ref{classtempcorank3} to \ref{classtempcorank4}, which can be done by induction for any arbitrary corank $r$. 
\end{remark}

\subsection{Unitarizability for the regular components}\label{ssregularcomponents}

For the rest of this section, we follow Tadi{\'c}'s notation in \cite[Chapter $8$]{Tad23}. For brevity, let $\underline{x} = (x_1, x_2, x_3, x_4)$ and $\Pi_{\underline{x}} = \Pi_{x_1, x_2, x_3, x_4}$.  Let 
\[\mathbb{R}_{++}^{4} = \{\underline{x} \in \mathbb{R}^4: 0 \leq x_1 \leq x_2 \leq x_3 \leq x_4\}.\]
Then $\Pi_{\underline{x}}$ is reducible if and only if $\underline{x}$ lies on one of the following singular affine hyperplanes: 

\begin{flalign}\label{4 d hyperplanes}
\begin{split}
    x_i &= \pm \alpha, \quad i = 1,2,3,4, \\
    x_i \pm x_j &= \pm1, \quad  1 \leq i < j \leq 4.
    \end{split}
\end{flalign}
We say that $\underline{x} \in \mathbb{R}^4$ is regular if it does not lie on any of these planes, and denote the set of such elements as $\mathbb{R}_{\text{reg}}^{4}$. Define 
\[\mathbb{R}_{\text{reg},++}^{4}  = \mathbb{R}_{\text{reg}}^{4} \cap \mathbb{R}_{++}^{4}.  \]
We also say that a point $\underline{x} \in \mathbb{R}^4$ is {\it strongly unitary} (resp.{ \it strongly non-unitary}) if all irreducible subquotients of $\Pi_{\underline{x}}$ are unitarizable (resp. non-unitarizable).

Similarly, let $\mathbb{R}^3_{\text{reg}}$ be the set of $(x_1, x_2, x_3) \in \mathbb{R}^3$ that lies in the complement of the hyperplanes
\begin{flalign}\label{3 d hyperplanes}
\begin{split}
    x_i &= \pm \alpha, \quad i = 1,2,3, \\
    x_i \pm x_j &= \pm1, \quad 1 \leq i < j \leq 3,
    \end{split}
\end{flalign}
and $\mathbb{R}^3_{\text{reg},++} = \mathbb{R}^3_{++} \cap \mathbb{R}^3_{\text{reg}}$.

We wish to classify all unitary representations of corank $4$. To begin, we consider the case where $\underline{x}$ is regular, in which case $\Pi_{\underline{x}}$ is either strongly unitary or strongly non-unitary. Denote the group of signed permutations by $W$. Since the property of being strongly unitary/non-unitary depends only on the $W$-orbit of $\underline{x}$, it suffices to consider the region $\mathbb{R}_{++}^{4}$. 

We say that an open connected component $\Omega \subset \mathbb{R}_{\text{reg}, ++}^{4}$ is unitary if for all points $\underline{x} \in \Omega$, $\Pi_{\underline{x}}$ is strongly unitary. For convenience, we do not distinguish between the connected component and the inequalities defining the connected component. 

\begin{prop}\label{connectedisunitary}
    The following open connected components of $\mathbb{R}_{\text{reg},++}^4$ are unitary. 
    For $\alpha \geq 1$,
    \begin{gather}
        x_3 + x_4 < 1. \tag{8.1a}\label{ciu1}\\
        x_2 + x_3 < 1, \quad x_4 - x_3 > 1, \quad x_4 < \alpha \quad (\alpha > 1). \tag{8.1b} \label{ciu2}\\
        x_2 + x_3 < 1, \ \ x_4 - x_2 < 1, \ \ x_4 - x_1 > 1, \ \ x_4 < \alpha \ \ (\alpha > 1). \tag{8.1c}\label{ciu3}\\
        x_2 - x_1 > 1, \quad x_3 - x_2 > 1, \quad x_4 - x_3 > 1, \quad x_4 < \alpha,  \quad (\alpha > 3). \tag{8.1d}\label{ciu4}\\      
        x_1 + x_2 < 1, \quad x_3 - x_2 > 1, \quad x_4 - x_3 > 1, \quad x_4 < \alpha, \quad (\alpha > 2). \tag{8.1e} \label{ciu5} \\
        x_1 + x_2 < 1, \quad x_1 + x_3 > 1, \quad x_3 - x_1 < 1, \quad x_4 - x_3 > 1, \quad x_4 < \alpha, \quad (\alpha > \frac{3}{2}). \tag{8.1f} \label{ciu6} \\
        x_2 + x_3 < 1, \quad x_1 + x_4 < 1, \quad x_2 + x_4 > 1,  \ \ x_4 < \alpha. \tag{8.1g}\label{ciu7} \\
        %x_2 - x_1 > 1, \quad x_4 - x_2 < 1 \quad x_4 < \alpha, \quad (\alpha > 1) \label{ciu5}\\
        %x_1 + x_2 < 1, \quad x_4 - x_2 < 1, \quad x_4 < \alpha \label{ciu6}
    \end{gather}
    For $\alpha =  \frac{1}{2}$: 
    \begin{gather}
        x_4 < \frac{1}{2}. \tag{8.1i} \label{ciu8}
    \end{gather}
    For $\alpha = 0$: 
    \begin{gather}
        x_3 + x_4 < 1. \tag{8.1j} \label{ciu9}
    \end{gather}
    Consequently, any point $\underline{x} \in \mathbb{R}^4_{++}$ that lies in the closure of the above region is also strongly unitary. 
\end{prop}

\begin{proof}
    For $\alpha \geq 1$, note that the non-empty connected components defined by (\ref{ciu1}) to (\ref{ciu7}) are all mutually disjoint. When $\alpha \neq 0$, $\Pi_{\vec{0}}$ is unitarizable and thus so is its connected component, which is given by (\ref{ciu1}) when $\alpha \geq 1$ and (\ref{ciu8}) when $\alpha = \frac{1}{2}$. 

    For $\alpha = 0$, we want to show that the connected component (\ref{ciu9}) is unitary. In this case, note that the condition 
    \[x_3 + x_4 < 1\]
    implies $0 \leq x_1 \leq x_2 \leq x_3 \leq x_4 \leq \frac{1}{2}$. We take the intersection of (\ref{ciu9}) with the hyperplane $x_3 = x_4$, which is nonempty. Applying Step ($3$-$4$) of Algorithm \ref{alg A bar}, we reduce to the connected component 
    \[0 \leq x_1 \leq x_2 < \frac{1}{2}\]
    in corank $2$, and using Theorem \ref{unitcrnk2}, one can conclude that (\ref{ciu9}) is indeed unitary. 

    Assume now that $\alpha \geq 1$. Let us consider the complementary series
    \[\pi_{x_4} = \Delta_\rho[-x_4, -x_4] \rtimes \pi_{sc},\]
    which is irreducible and unitarizable for $0 \leq x_4 < \alpha$. 
    We use it to prove that the connected components (\ref{ciu2}) to (\ref{ciu7}) are unitary one by one. 

    \begin{enumerate}
        \item Assume $\alpha > 1$ and fix $1 < x_4 < \alpha$. Then $(0,0,0,x_4) \in \mathbb{R}_{\text{reg}, ++}^{4}$ is strongly unitary by considering $\pi_{x_4}$. Since (\ref{ciu2}) is the connected component of $(0,0,0,x_4)$, it must be unitary. 
        \item  Let $\alpha > 1$. Fix $1 < x_4 < \frac{3}{2}$. Then $\pi_{x_4}$ is unitarizable. Now fix $x_4 - 1 < x_3 < \frac{1}{2}$, then the complementary series 
        \[\Delta_\rho[-x_3,-x_3] \times \Delta_\rho[-x_3, -x_3] \rtimes \pi_{x_4}\]
        is irreducible and unitarizable. Finally, if we take $x_1  <x_4 -1$, then the complementary series 
        \[\Delta_\rho[-x_1, -x_1] \times \Delta_\rho[-x_3,-x_3] \times \Delta_\rho[-x_3, -x_3] \rtimes \pi_{x_4}\]
        is irreducible and unitarizable. The set of such points $(x_1, x_3, x_3, x_4)$ lies in the connected component (\ref{ciu3}), so the component must be unitary. 
        \item Let $\alpha > 3$. Fix $3 < x_4 < \alpha$ and construct the complementary series
        \[ \Delta_\rho[-x_3, -x_3] \rtimes \pi_{x_4},\]
        for $0 \leq x_3 < x_4 -1$. Fixing $2 \leq x_3 < x_4-1$, we can construct the complementary series
        \[\Delta_\rho[-x_2, -x_2] \times \Delta_\rho[-x_3, -x_3] \rtimes \pi_{x_4},\]
        for $0 \leq x_2 < x_3 -1$. Now fix $1 < x_2 < 1 - x_3$ and construct the complementary series 
        \[\Delta_\rho[-x_1, -x_1] \times \Delta_\rho[-x_2, -x_2] \times \Delta_\rho[-x_3, -x_3] \rtimes \pi_{x_4},\]
        for $0 \leq x_1< x_2-1$. The corresponding $(x_1, x_2, x_3, x_4)$ is unitary and is contained in the connected component (\ref{ciu4}). Therefore (\ref{ciu4}) must be unitary. 
        \item Let $\alpha > 2$. Fixing $2 < x_4 < \alpha$, we can construct the complementary series 
        \[\Delta_\rho[-x_3, -x_3] \rtimes \pi_{x_4},\]
        which is irreducible and unitarizable for $0 \leq x_3 < x_4 -1$. Now fix $1 < x_3 < x_4 -1$ and consider the complementary series
        \[\Delta_\rho[-x_2, -x_2] \times \Delta_\rho[-x_3, -x_3] \rtimes \pi_{x_4}.\]
        This is irreducible and unitarizable for $0 \leq x_2 < x_3 -1$. Subsequently, fix $0 < x_2 < \min(\frac{1}{2}, x_3 -1)$, the complementary series 
        \[\Delta_\rho[-x_1, -x_1] \times \Delta_\rho[-x_2, -x_2] \times  \Delta_\rho[-x_3, -x_3] \rtimes \pi_{x_4},\]
        is irreducible and unitarizable for $0 \leq x_1 < x_2$. The connected component containing such $(x_1, x_2, x_3, x_4)$ is (\ref{ciu5}), so (\ref{ciu5}) must be unitary.

        \item Let $\alpha > \frac{3}{2}$. Fixing $\frac{3}{2} < x_4 < \alpha$ we can construct the complementary series 
        \[\Delta_\rho[-x_3, -x_3] \rtimes \pi_{x_4},\]
        which is irreducible and unitarizable for $0 \leq x_3 < x_4 -1$. Now fix $\frac{1}{2} < x_3 < x_4 -1$. The point $(0,0,x_3, x_4)$ is now strongly unitary. This point lies in the connected component (\ref{ciu6}), so (\ref{ciu6}) must be unitary. 

        \item Let $\alpha \geq 1$. Then the connected component (\ref{ciu7}) has nontrivial $3$-dimensional intersection with the hyperplane $x_1 = 0$. The intersection can be described by the inequalities
        \[x_2 + x_3 < 1, \ \ x_4 < 1, \ \ x_2 + x_4 > 1.\]
        For $(x_1, x_2, x_3, x_4) \in \R^4_{\reg,++}$, the representation 
        \[\Pi_{0,x_2,x_3,x_4} = \rho \times \rho \lvert \cdot \rvert^{x_2} \times \rho \lvert \cdot \rvert^{x_3} \times \rho \lvert \cdot \rvert^{x_4}\]
        is unitarizable if and only if the representation 
        \[\Pi_{x_2, x_3, x_4} = \lvert \cdot \rvert^{x_2} \times \rho \lvert \cdot \rvert^{x_3} \times \rho \lvert \cdot \rvert^{x_4}\]
        is unitarizable. Since above $3$-dimensional region is contained in one of the connected components listed in Proposition \ref{crnk3unitcc}, it follows that (\ref{ciu7}) is unitary. 

    \end{enumerate}
    
\end{proof}

In \S \ref{exhaustion} below, we show that the other open connected components in $\R^4_{++}$ are not unitary. To this end, we need to consider the cases where $\underline{x}$ is not regular, beginning with the singular affine hyperplanes.

\subsection{Unitarizability for the irregular components - slanted hyperplanes}\label{slantedsection}

As described in Step $(3$-$3)$ and $(3$-$4)$ of Algorithm \ref{alg A bar}, there are two types of reducibility hyperplanes, which we will classify in this and the following subsection. In this subsection, we aim to to perform the inductive process described in Step $(3$-$4)$ of Algorithm \ref{alg A bar}. To begin, let us consider the hyperplane 
\[
H_{\text{sla}} = \{\underline{x} \in \mathbb{R}^4: x_2 - x_1 = 1\}.
\]
We call the $W$-orbit of $H_{\text{sla}}$ the set of {\it slanted hyperplanes}, which is precisely the set of hyperplanes defined by equations of the form

\[
x_i \pm x_j = \pm1, \quad 1 \leq i < j \leq 4.
\]
Let $x = x_1 + \frac{1}{2}, y = x_3, z = x_4$, then clearly $\iota(x,y,z): \mathbb{R}^3 \to H_{\text{sla}}$ defined by 
\[
\iota(x,y,z) = (x - \frac{1}{2}, x + \frac{1}{2}, y,z)
\]
is an affine isomorphism.

Consider all reducibility hyperplanes which intersect $H_{\shrt}$ non-trivially. These are precisely the set of hyperplanes described in \eqref{4 d hyperplanes} other than $x_2-x_1=\pm1$.
We define $\R^3_{v\reg,\shrt} \subset \R^3$ to be the complement of the hyperplanes
\begin{gather*}
x\pm y=\pm\tfrac12,\ \ x \pm z = \pm\tfrac12, \ \ x=\pm\tfrac12,\\
x\pm y=\pm\tfrac32, \ \ x \pm z = \pm\tfrac32, \ \  \ \ y\pm z = \pm 1, \\
 y=\pm\alpha,\ \ z = \pm \alpha, \ \ x=\pm(\alpha-\tfrac12),\ \ x=\pm(\alpha+\tfrac12).
\end{gather*}
Note that $\R^3_{v\reg,\shrt,++} := \{(x,y,z) \in \R^3_{v\reg,\shrt} \, | \, 0 \leq x \leq y \leq z\}$ is the complement of the following hyperplanes: 
\begin{gather*}
    y \pm x = \tfrac12,  \ \  z \pm x = \tfrac12,   \ \ x = \tfrac12, \\
    y \pm x = \tfrac32, \ \ z \pm x = \tfrac32, \ \ z \pm y = 1, \\
    y = \alpha, \ \ z = \alpha, \ \ x = |\alpha - \tfrac12|, \ \ x = \alpha + \tfrac12.
\end{gather*}
Then $\iota(\R^3_{v\reg, \shrt}) = H_{\shrt} \cap \R^4_{\reg}$. 

For $(x,y,z) \in \mathbb{R}^3$, we can decompose $\Pi_{\iota(x,y,z)}$ in the Grothendieck group as 
\[
\Pi_{\iota(x,y,z)} = \pi_{(x,y,z)}^{+} + \pi_{(x,y,z)}^{-},
\]
where
\begin{flalign*}
    \pi_{(x,y,z)}^{+} &= \lvert \cdot \rvert^{x}\Delta_\rho[\frac{1}{2}, -\frac{1}{2}] \times \Delta_\rho[-y,-y] \times \Delta_\rho[-z,-z] \rtimes \pi_{sc}, \\
    \pi_{(x,y,z)}^{-} &= \lvert \cdot \rvert^{x}L(\Delta_\rho[-\frac{1}{2}, -\frac{1}{2}], \Delta_\rho[\frac{1}{2}, \frac{1}{2}]) \times \Delta_\rho[-y,-y] \times \Delta_\rho[-z,-z] \rtimes \pi_{sc}.
\end{flalign*}
Now, let  $\R^3_{\reg,\shrt} \subset \R^3$ be the complement of the hyperplanes
\begin{gather*}
    x\pm y = \pm\tfrac32, \ \ x\pm z = \pm\tfrac32, \ \  y \pm z = \pm1,\\
     y = \pm \alpha, \ \ z = \pm \alpha, \ \ x = \pm(\alpha - \frac{1}{2}), \ \ x = \pm(\alpha + \frac{1}{2}).
\end{gather*}
Then, $\pi_{(x,y,z)}^{+}$ and $\pi_{(x,y,z)}^{-}$ are both irreducible if and only if $(x,y,z) \in \mathbb{R}^3_{\reg, \shrt}$, since the hyperplanes above describe all the possible reducibility conditions for the parabolic induction $\pi^{\pm}(x,y,z)$. In other words, for $(x,y,z) \in \mathbb{R}^3_{\reg, \shrt}$, $\iota(x,y,z)$ is strongly unitary (resp., strongly non-unitary) if and only if both $\pi_{(x,y,z)}^{+}$ and $\pi_{(x,y,z)}^{-}$ are unitarizable (resp., non-unitarizable). 

Denote
\[
W_{\text{sla}} = \{w \in W: w(H_{\text{sla}}) = H_{\text{sla}}\}.
\]
Then $W_{\text{sla}} \cong \{\pm1\}^{3}$, and $W_\shrt$ can be generated by the signed permutations
\begin{flalign*}
    (x_1, x_2, x_3, x_4) &\mapsto (-x_2, -x_1, x_3, x_4), \\
    (x_1, x_2, x_3, x_4) &\mapsto (x_1, x_2, -x_3, x_4), \\
    (x_1, x_2, x_3, x_4) &\mapsto (x_1, x_2, x_3, -x_4).
\end{flalign*}

We say that a point $(x,y,z) \in \mathbb{R}^3_{\reg,\shrt}$ is unitary$^{+}$ (resp., unitary $^{-}$) if $\pi_{(x,y,z)}^{+}$ (resp., $\pi_{(x,y,z)}^{-}$) is unitarizable. We say that a point $(x,y,z) \in \mathbb{R}^3_{\reg,\shrt}$ is unitary$^{\pm}$ if it is both unitary$^{+}$ and unitary$^{-}$. We say that a connected component is unitary$^{+}$ (resp., unitary$^-$, unitary$^{\pm}$) if every point in it is unitary$^+$ (resp., unitary$^{-}$, unitary$^\pm$). 

Recall that $\R^3_{+} = \{(x,y,z) \in \R^3: x,y,z \geq 0\}$, $\R^3_{++} = \{(x,y,z) \in \R^3: 0\leq x \leq y \leq z\}$, and let $\R^3_{\reg, \shrt, +}= \R^3_{+} \cap \R^3_{\reg, \shrt}, \R^3_{\reg, \shrt, ++}= \R^3_{++} \cap \R^3_{\reg, \shrt}$. We start by looking at the preimage of the slanted hyperplane boundaries of the regions described in Proposition \ref{connectedisunitary}, under the isomorphism $\iota$. In the proposition below, we denote the connected components (\ref{ciu1}) to (\ref{ciu7}) defined in Proposition \ref{connectedisunitary} by $C_a, C_b, C_c, C_d, C_e, C_f, C_g, C_h$. 

\begin{prop}\label{3dslaunitcut}
For $\alpha \geq 1$, any connected component of $\mathbb{R}^{3}_{v\reg, \shrt, ++}$ is either fully contained in one of the following regions in $\R^3_{++}$, or has empty intersection with each of them. Moreover, in the former case, that component is unitary$^{\pm}$: 
    \begin{gather*}
          x+z < \frac{1}{2}, \tag{$L_1'$} \label{3dLgeq1,1} \\ 
          x+y > \frac{1}{2}, \ \ z-x < \frac{1}{2}, \ \  y+z < 1, \ \ x < \frac{1}{2}. \tag{$L_2'$} \label{3dLgeq1,2} \\
     \end{gather*}
     For $\alpha > 1$: 
     \begin{gather*}
         x+y < \frac{1}{2}, \ \ z-x > \frac{3}{2}, \ \ z < \alpha. \tag{$L_3'$}\label{3dLg1,1}\\
     \end{gather*}
     For $\alpha > \frac{3}{2}$: 
     \begin{gather*}
      x+y < \frac{3}{2}, \ \ z-y > 1, \ \ x < \frac{1}{2}, \ \ z < \alpha, \tag{$L_4'$}\label{3dLg3/2,1} \\
     x+y > \frac{3}{2}, \ \ z-x > \frac{3}{2}, \ \ x < \frac{1}{2}, \ \ z < \alpha, \tag{$L_5'$} \label{3dLg3/2,2} \\
    x+y < \frac{3}{2}, \ \ z-x > \frac{3}{2}, \ \ x > \frac{1}{2}, \ \ z < \alpha. \tag{$L_6'$}\label{3dLg3/2,3}
     \end{gather*}

     For $\alpha > 3$: 
     \begin{gather*}
        y-x > \frac{3}{2}, \ \ z-y > 1, \ \ z < \alpha. \tag{$L_7'$}\label{3dLg3,1}\\
     \end{gather*}

    Furthermore, for each of the connected components $L_i \in \{$\ref{3dLgeq1,1}, $\ldots,$ \ref{3dLg3,1}$\}$, there exists some $w_i \in W$ and $X_i \in \{C_a, C_b, C_c, C_d, C_e, C_f, C_g, C_h\}$ such that $\iota(L_i')$ is contained in $\partial(w_i(X_i))$. 
\end{prop}

\begin{proof}
    For $\alpha \geq 1$, consider the region $L_1'$, defined by 
    \begin{equation*}
         x+z < \frac{1}{2}, \ \ 0 \leq x \leq y \leq z.
    \end{equation*}
     Take $X_1 = C_a$, which is defined by the equation 
    \begin{equation*}
        x_3 + x_4 < 1, \ \ 0 \leq x_1 \leq x_2 \leq x_3 \leq x_4.
    \end{equation*}
    The boundary of the connected component $X_1$ is the region defined by 
    \begin{equation*}
        x_3 + x_4 = 1, \ \ 0 \leq x_1 \leq x_2 \leq x_3 \leq x_4.
    \end{equation*}
    Now let $w_1$ be the map defined by 
    \[w_1(x_1, x_2, x_3, x_4) = (-x_3, x_4, x_1, x_2) := (x_1', x_2', x_3', x_4').\]
    Then the boundary $\partial(w_1(X_1))$ becomes 
    \[x_2' - x_1' = 1, \ \ 0 \leq x_3' \leq x_4' \leq -x_1' \leq x_2'.\]
    Introducing the transformation $(x_1', x_2', x_3', x_4') = (x-\frac{1}{2}, x + \frac{1}{2}, y, z)$ denoted as $\iota$, the condition on the coordinates translates to 
    \[0 \leq y \leq z \leq -x + \frac{1}{2} \leq x + \frac{1}{2}.\]
    We see then clearly the region $L_1'$ satisfies these constraints, and thus is contained in the boundary $\partial (w_1(X_1))$. It follows that any connected component of $\R^3_{v\reg, \shrt}$ that is fully contained in $L_1'$ is unitary$^{\pm}$. This proves the statement for $i = 1$. Now consider 
    \begin{flalign*}
        X_2 &= C_g, \ \ w_2: (x_1, x_2, x_3, x_4) \mapsto (-x_1, x_3, x_4, x_2), \\
        X_3 &= C_b, \ \ w_3: (x_1, x_2, x_3, x_4) \mapsto (-x_1, x_3, x_1, x_4), \\
        X_4 &= C_f, \ \ w_4: (x_1, x_2, x_3, x_4) \mapsto (-x_1, x_3, x_2, x_4), \\
        X_5 &= C_f, \ \ w_5: (x_1, x_2, x_3, x_4) \mapsto (-x_1, x_3, -x_2, x_4), \\
        X_6 &= C_f, \ \ w_6: (x_1, x_2, x_3, x_4) \mapsto (x_1, x_3, x_2, x_4), \\
        X_7 &= C_d, \ \ w_7: (x_1, x_2, x_3, x_4) \mapsto (x_1, x_2, x_3, x_4). \\
    \end{flalign*}
    Using these $X_i$s and $w_i$s, we can prove the statement for $i = 2, \ldots 9$ in an analogous way.    
\end{proof}

Using this result, we can find all of the unitary connected components of $\mathbb{R}^3_{\reg,\shrt, ++}$. To do this, we need to invoke our classification of critical points in the previous sections. 
Let us first recall Tadi{\'c}'s result on the unitary dual of corank $2$: 
\begin{prop}[{\cite[Proposition $7.2$]{Tad23}}]\label{unitcrnk2}
    The irreducible unitarizable subquotients of $\Pi_{\underline{x}}$, where $\underline{x} = (x_1, x_2) \in \R^2_{++}$, are the following: 
    \begin{enumerate}
        \item $(\alpha > 1)$ All irreducible subquotients when $x_1 +1 \leq x_2 \leq \alpha$. 
        \item $(\alpha \neq \frac{1}{2})$ All irreducible subquotients when $x_1 + x_2 \leq 1$. 
        \item $(\alpha = \frac{1}{2})$ All irreducible subquotients when $x_2 \leq \frac{1}{2}$. 
        \item $(\alpha > 0)$ The irreducible representations $T_{I,1}^{\alpha+1}(T_{I,1}^{\alpha}(\pi_{sc}))$ and $L([-\alpha, -\alpha], [-\alpha-1, -\alpha-1]; \pi_{sc})$.
    \end{enumerate}
\end{prop}

Let us now give the list of the unitary components inside $\R^3_{\reg, \shrt, ++}$, which we prove to be the complete list using a method of exhaustion. Note that one can use various computer algorithms, such as Sage, to determine the possible bounded connected components in a region, given a set of hyperplanes separating the components. 

\begin{prop}\label{finalslantunit}
     The unitary$^{\pm}$ connected components of $\R^3_{\reg, \shrt, ++}$ are precisely given by 
    \begin{gather}
        y+z < 1, \ \ (\alpha \geq 1,  \quad \alpha = 0)\tag{\ref{finalslantunit}a}\label{finalslantunit1}\\
        x+z < \frac{3}{2}, \ \ z-y > 1 \ \ z < \alpha \tag{\ref{finalslantunit}b}\label{finalslantunit2}, \ \ (\alpha > 1)\\\
        x+y < \frac{3}{2}, \ \ x+z > \frac{3}{2}, \ \ z-x < \frac{3}{2}, \ \ z-y > 1 \label{finalslantunit3}, \tag{\ref{finalslantunit}c}\ \ (\alpha > 1) \\\
        x+y < \frac{3}{2}, \ \ z-x > \frac{3}{2}, \ \ z-y > 1, \ \ z < \alpha \label{finalslantunit4},\tag{\ref{finalslantunit}d} \ \ (\alpha > \frac{3}{2}) \\\
        y-x > \frac{3}{2}, \ \ z-y > 1, \ \ z < \alpha, \ \ (\alpha > \frac{5}{2}) \tag{\ref{finalslantunit}e}\label{finalslantunit5} \\
        z < \frac{1}{2} \ \ (\alpha = \frac{1}{2}).  \tag{\ref{finalslantunit}f}\label{finalslantunit6}
    \end{gather}
\end{prop}
\begin{proof}
    For $\alpha > 3$, the following are all the possible bounded connected components of $\R^3_{\reg, \shrt, ++}$:
    \begin{gather*}
        y+z < 1, \tag{$C_{>3, 1}$} \label{regslaunitg3,1}\\
        x+z < \frac{3}{2}, \ \ y+z > 1, \ \ z-y < 1, \tag{$C_{>3, 2}$} \label{regslaunitg3,2}\\
        x+z < \frac{3}{2}, \ \ z-y > 1, \tag{$C_{>3, 3}$} \label{regslaunitg3,3}\\
        x+y < \frac{3}{2}, \ \ x+z > \frac{3}{2}, \ \ z-x < \frac{3}{2}, \ \ z-y < 1, \tag{$C_{>3, 4}$} \label{regslaunitg3,4}\\
        x+y < \frac{3}{2}, \ \ x+z > \frac{3}{2}, \ \ z-x < \frac{3}{2}, \ \ z-y > 1, \tag{$C_{>3, 5}$} \label{regslaunitg3,5} \\
        x+y < \frac{3}{2}, \ \ z-x > \frac{3}{2}, \ \ z-y < 1, \tag{$C_{>3, 6}$} \label{regslaunitg3,6}\\
        x+y < \frac{3}{2}, \ \ z-x > \frac{3}{2}, \ \ z-y > 1, \ \ z < \alpha, \tag{$C_{>3, 7}$} \label{regslaunitg3,7}\\
        x < \alpha - \frac{1}{2}, \ \ z < \alpha, \ \ x+y > \frac{3}{2}, \ \ z-y < 1, \ \ z-x < \frac{3}{2}, \tag{$C_{>3, 8}$} \label{regslaunitg3,8}\\
        x < \alpha - \frac{1}{2}, \ \ y < \alpha, \ \ z-x < \frac{3}{2}, \ \ z > \alpha, \ \ z-y < 1, \tag{$C_{>3, 9}$} \label{regslaunitg3,9}\\
        %x < \alpha - \frac{1}{2}, \ \ z-x < \frac{3}{2}, \ \ y > \alpha \tag{$C_{>3, 10}$} \label{regslaunitg3,10} \\
        x > \alpha -\frac{1}{2}, \ \ z < \alpha, \tag{$C_{>3, 10}$} \label{regslaunitg3,10} \\
        x > \alpha - \frac{1}{2}, \ \ y < \alpha, \ \ z > \alpha, \ \ z-y < 1, \tag{$C_{>3, 11}$} \label{regslaunitg3,11}\\
        \alpha - \frac{1}{2} < x < \alpha + \frac{1}{2}, \ \ z-x < \frac{3}{2}, \ \ z-y < 1, \ \ y > \alpha, \tag{$C_{>3, 12}$} \label{regslaunitg3,12}\\
        z-y > 1, \ \ z < \alpha, \ \ x+y > \frac{3}{2}, \ \ z-x < \frac{3}{2}, \tag{$C_{>3, 13}$} \label{regslaunitg3,13}\\
       % x < \alpha - \frac{1}{2}, \ \ z-y > 1, \ \ z-x < \frac{3}{2}, \ \ z > \alpha \tag{$C_{>3, 15}$} \label{regslaunitg3,15}\\
        y < \alpha, \ \ z-y > 1, \ \ z-x < \frac{3}{2}, \ \ x > \alpha - \frac{1}{2}, \tag{$C_{>3, 14}$} \label{regslaunitg3,14}\\
        %x < \alpha + \frac{1}{2}, \ \ z-y > 1, \ \ z-x < \frac{3}{2}, \ \ y > \alpha \tag{$C_{>3, 17}$} \label{regslaunitg3,17}\\
        z-x > \frac{3}{2}, \ \ z < \alpha, \ \ x+y > \frac{3}{2}, \ \ z-y > 1, \ \ y-x < \frac{3}{2}, \tag{$C_{>3, 15}$} \label{regslaunitg3,15}\\
        z-x > \frac{3}{2}, \ \ y < \alpha, \ \ y-x < \frac{3}{2}, \ \ z > \alpha, \ \ z-y < 1, \tag{$C_{>3, 16}$} \label{regslaunitg3,16} \\
        z-x > \frac{3}{2}, \ \ x < \alpha - \frac{1}{2}, \ \ y-x < \frac{3}{2}, \ \ z-y < 1, \ \ y > \alpha, \tag{$C_{>3, 17}$} \label{regslaunitg3,17} \\
        z-x > \frac{3}{2}, \ \ \alpha - \frac{1}{2} < x < \alpha + \frac{1}{2}, \ \ y-x < \frac{3}{2}, \ \ z-y < 1, \tag{$C_{>3, 18}$} \label{regslaunitg3,18}\\
        z-x > \frac{3}{2}, \ \ y-x < \frac{3}{2}, \ \ z-y < 1, \ \ x > \alpha + \frac{1}{2}, \tag{$C_{>3, 19}$} \label{regslaunitg3,19}\\
        z-x > \frac{3}{2}, \ \ z-y > 1, \ \ z < \alpha, \ \ x+y > \frac{3}{2}, \ \ y-x < \frac{3}{2}, \tag{$C_{>3, 20}$} \label{regslaunitg3,20}\\
        y-x > \frac{3}{2}, \ \ z < \alpha, \ \ z -y < 1, \tag{$C_{>3, 21}$} \label{regslaunitg3,21}\\
        y-x > \frac{3}{2}, \ \ y < \alpha, \ \ z > \alpha, \ \ z-y < 1, \tag{$C_{>3, 22}$} \label{regslaunitg3,22}\\
        y-x > \frac{3}{2}, \ \ z-y > 1 , \ \ z < \alpha. \tag{$C_{>3, 23}$} \label{regslaunitg3,23}
    \end{gather*}
    One can verify through direct computation that $C_{>3,i} \subset \cup_{j=1}^{7}L_{j}'$ for $i = 1,3,5,7,23$. We now show that the other connected components are not unitary. In the case $\alpha > 3$, there are $6$ points that gets mapped to not strongly unitary points in $\R^4_{++}$ under $\iota$: 
    \[(\alpha +\frac{1}{2}, \alpha+2, \alpha+3), \ \ (\alpha- \frac{1}{2}, \alpha+1, \alpha+2), \ \ (\alpha-\frac{1}{2}, \alpha, \alpha+1),\] 
    \[(\alpha-\frac{3}{2}, \alpha, \alpha+1), \ \ (\alpha -\frac{3}{2}, \alpha, \alpha), \ \ (\alpha - \frac{3}{2}, \alpha-1, \alpha).\]
    One can easily verify that $C_{>3,i}$ contains at least one of the above critical points on its boundary for $i = 8,9,11,12,13\ldots,23$. The only components remaining are $C_{>3,i}$ for $i = 2,4,6,10$. 

    For $i = 2,4,6$, the connected component $C_{>3,i}$ has nonempty $3$-dimensional intersection with the hyperplane $x = 0$. Applying Step $(3$-$2)$ of Algorithm \ref{alg A bar}, one can see that the representations
    \begin{flalign*}
        \pi_{(0,y,z)}^{+} &= \Delta_\rho[\frac{1}{2}, -\frac{1}{2}] \times \Delta_\rho[-y,-y] \times \Delta_\rho[-z,-z] \rtimes \pi_{sc}, \\
        \pi_{(0,y,z)}^{-} &= L(\Delta_\rho[-\frac{1}{2}, -\frac{1}{2}], \Delta_\rho[\frac{1}{2}, \frac{1}{2}]) \times \Delta_\rho[-y,-y] \times \Delta_\rho[-z,-z ] \rtimes \pi_{sc}.
    \end{flalign*}
    are unitarizable if and only if $\pi_{(x,y)} = \Delta_\rho[-x,-x] \times \Delta_\rho[-y,-y] \rtimes \pi_{sc}$ contains only unitarizable subquotients. 
    By Proposition \ref{unitcrnk2}, the representations $\pi_{(0,y,z)}$ is not unitary for $(y,z) = (\frac{2}{3}, \frac{2}{3})$ or $(\frac{4}{3}, \frac{5}{3})$. However, the points $(0,\frac{2}{3}, \frac{2}{3})$ lies in $C_{>3,2}$ and the point $(0,\frac{4}{3}, \frac{5}{3})$ lies in $C_{>3,6}$. Therefore these two connected components cannot be unitary$^{\pm}$. 

    For $C_{>3,10}$, consider the critical point $(\alpha-\frac{1}{2}, \alpha, \alpha)$, which lies on the boundary. By Proposition \ref{subquotientlist}, the point 
    \[\iota(\alpha - \frac{1}{2}, \alpha, \alpha) = (\alpha-1, \alpha, \alpha, \alpha)\]
    is not strongly unitary for $\alpha > 3$. Thus the connected component $C_{>3,11}$ is not unitary$^{\pm}$. 

    Finally, for $C_{>3,4}$, we split into two cases. When $\alpha \in \frac{1}{2} + \mathbb{Z}$, the point $(0,\frac{3}{2}, \frac{3}{2})$ is a critical point that lies on the boundary of $C_{>3,4}$. By Proposition \ref{unitcrnk2}, $(\frac{3}{2}, \frac{3}{2})$ is not in the unitary dual of corank $2$, so the representations $(0,\frac{3}{2}, \frac{3}{2})$ is not strongly unitary, which means that $C_{>3,4}$ is not unitary$^{\pm}$. 

    Now consider the case $\alpha \in \Z_{>3}$. In this case, the connected component $C_{>3,4}$ contains no points of good parity. Instead, we consider the point $(\frac{1}{2},1,1)$ which lies its boundary. The representations
    \begin{flalign*}
        \pi_{(\frac{1}{2},1,1)}^{+} &= \lvert \cdot \rvert^{\frac{1}{2}}\Delta_\rho[\frac{1}{2}, -\frac{1}{2}] \times \rho \lvert \cdot \rvert^{1} \times \rho \lvert \cdot \rvert^{1} \rtimes \pi_{sc}, \\
        \pi_{(\frac{1}{2},1,1)}^{-} &= \lvert \cdot \rvert^{\frac{1}{2}}L(\Delta_\rho[-\frac{1}{2}, -\frac{1}{2}], \Delta_\rho[-\frac{1}{2}, \frac{1}{2}]) \times \rho \lvert \cdot \rvert^{1} \times \rho \lvert \cdot \rvert^{1} \rtimes \pi_{sc},
    \end{flalign*}
    are irreducible. Consider now the family
    \begin{flalign*}
        \tau_{(x,1,1)}^{+} &= \lvert \cdot \rvert^{x}\Delta_\rho[\frac{1}{2}, -\frac{1}{2}] \times \rho \lvert \cdot \rvert^{1} \times \rho \lvert \cdot \rvert^{1} \rtimes \pi_{sc}, \\
        \tau_{(x,1,1)}^{-} &= \lvert \cdot \rvert^{x}(\Delta_\rho[-\frac{1}{2}, -\frac{1}{2}], \Delta_\rho[-\frac{1}{2}, \frac{1}{2}]) \times \rho \lvert \cdot \rvert^{1} \times \rho \lvert \cdot \rvert^{1} \rtimes \pi_{sc},
    \end{flalign*}
    which is irreducible for $0 \leq x \leq \frac{1}{2}$. When $x = 0$, we can apply Step $(3$-$2)$ of Algorithm \ref{alg A bar} again to conclude that $\tau_{(0,1,1)}$ is non-unitarizable, since $(1,1)$ does not lie in the unitary dual of corank $2$. Therefore it follows that $\tau_{\frac{1}{2},1,1}^{\pm}$ is non-unitarizable, which means that $C_{>3,4}$ is not unitary$^{\pm}$. This proves the claim for $\alpha > 3$. When $\alpha = 3$, the number of connected components and their descriptions are completely identical to the case $\alpha > 3$. The critical points for $\alpha = 3$ also matches the description for the critical points for $\alpha > 3$. The proof in this case is exactly the same as above, so we omit it. 

    Now consider the case $\alpha = \frac{5}{2}$. In this case the list of possible bounded connected components in $\R^3_{\reg, \shrt, ++}$ is as follows: 
    \begin{gather*}
        y+z <1, \tag{$C_{\frac{5}{2}, 1}$} \label{regslaunit5/2,1}\\
        x+z < \frac{3}{2}, \ \ y+z > 1, \ \ z-y < 1, \tag{$C_{\frac{5}{2}, 2}$} \label{regslaunit5/2,2}\\
        x+z < \frac{3}{2}, \ \ z-y > 1, \tag{$C_{\frac{5}{2}, 3}$} \label{regslaunit5/2,3}\\
        x+y < \frac{3}{2}, \ \ x+z > \frac{3}{2}, \ \ z-y < 1, \ \ z-x < \frac{3}{2}, \tag{$C_{\frac{5}{2}, 4}$} \label{regslaunit5/2,4}\\
        x+y < \frac{3}{2},  \ \ x+z > \frac{3}{2}, \ \ z-y > 1, \ \ z-x < \frac{3}{2}, \tag{$C_{\frac{5}{2}, 5}$} \label{regslaunit5/2,5} \\
        x+y < \frac{3}{2}, \ \ z-x > \frac{3}{2}, \ \ z-y < 1, \tag{$C_{\frac{5}{2}, 6}$} \label{regslaunit5/2,6}\\
        x+y < \frac{3}{2}, \ \ z-x > \frac{3}{2}, \ \ z-y > 1, \ \ z < \alpha, \tag{$C_{\frac{5}{2}, 7}$} \label{regslaunit5/2,7}\\
        x < \alpha - \frac{1}{2}, \ \ z < \alpha, \ \ x+y > \frac{3}{2}, \ \ z-y < 1, \ \ z-x < \frac{3}{2}, \tag{$C_{\frac{5}{2}, 8}$} \label{regslaunit5/2,8}\\
        x < \alpha - \frac{1}{2}, \ \ y < \alpha, \ \ z-x < \frac{3}{2}, \ \ z > \alpha, \ \ z-y < 1, \tag{$C_{\frac{5}{2}, 9}$} \label{regslaunit5/2,9}\\
        x > \alpha - \frac{1}{2}, \ \ z < \alpha, \tag{$C_{\frac{5}{2}, 10}$} \label{regslaunit5/2,10} \\
        x > \alpha -\frac{1}{2}, \ \ y<\alpha, \ \ z > \alpha \ \ z-y < 1, \tag{$C_{\frac{5}{2}, 11}$} \label{regslaunit5/2,11}\\
        \alpha - \frac{1}{2} < x < \alpha + \frac{1}{2}, \ \ z-x < \frac{3}{2}, \ \ z-y < 1, \ \ y > \alpha, \tag{$C_{\frac{5}{2}, 12}$} \label{regslaunit5/2,12}\\
        z-y > 1, \ \ z < \alpha \ \ x+y > \frac{3}{2}, \ \ z-x < \frac{3}{2}, \tag{$C_{\frac{5}{2}, 13}$} \label{regslaunit5/2,13} \\
        z-y > 1,  \ \ y < \alpha, \ \ z-x < \frac{3}{2}, \ \ x > \alpha - \frac{1}{2},\tag{$C_{\frac{5}{2}, 14}$} \label{regslaunit5/2,14}\\
        z-x > \frac{3}{2}, \ \ z < \alpha, \ \ x+y > \frac{3}{2}, \ \ z-y < 1, \ \ y-x < \frac{3}{2}, \tag{$C_{\frac{5}{2}, 15}$} \label{regslaunit5/2,15}\\
        z-x > \frac{3}{2}, \ \ y < \alpha, \ \ y-x < \frac{3}{2}, \ \ z > \alpha, \ \ z-y < 1, \tag{$C_{\frac{5}{2}, 16}$} \label{regslaunit5/2,16}\\
        z-x > \frac{3}{2}, \ \ x < \alpha - \frac{1}{2}, \ \ y-x < \frac{3}{2},  \ \ z-y < 1, \ \ y > \alpha, \tag{$C_{\frac{5}{2}, 17}$} \label{regslaunit5/2,17}\\
        z-x > \frac{3}{2}, \ \ \alpha - \frac{1}{2} < x < \alpha + \frac{1}{2}, \ \ y-x < \frac{3}{2}, \ \ z-y < 1, \tag{$C_{\frac{5}{2}, 18}$} \label{regslaunit5/2,18}\\
        z-x > \frac{3}{2}, \ \ z-y > 1, \ \ z < \alpha, \ \ x+y > \frac{3}{2}, \tag{$C_{\frac{5}{2}, 19}$} \label{regslaunit5/2,19}\\
        y-x > \frac{3}{2}, \ \ z < \alpha, \tag{$C_{\frac{5}{2}, 20}$} \label{regslaunit5/2,20}\\
        y-x > \frac{3}{2}, \ \ y < \alpha, \ \ z > \alpha, \ \ z-y < 1. \tag{$C_{\frac{5}{2}, 21}$} \label{regslaunit5/2,21}\\
    \end{gather*}
In this case, one can verify that $C_{\frac{5}{2}, i} \subset \cup_{j=1}^{6}L_{j'}$ for $i = 1,3,5,7$. By Proposition \ref{subquotientlist}, there are $6$ points that get mapped to not strongly unitary points in $\R^4_{++}$ under $\iota$, which are 
\[(\alpha +\frac{1}{2}, \alpha+2, \alpha+3), \ \ (\alpha- \frac{1}{2}, \alpha+1, \alpha+2), \ \ (\alpha-\frac{1}{2}, \alpha, \alpha+1),\] 
    \[(\alpha-\frac{3}{2}, \alpha, \alpha+1), \ \ (\alpha -\frac{3}{2}, \alpha, \alpha), \ \ (\alpha - \frac{3}{2}, \alpha-1, \alpha).\]
The proof of non-unitarity for the other regions is identical to the case $\alpha > 3$, which we omit. 

The next case to consider is $\alpha = 2$. In this case, the list of possible bounded  connected components is as follows: 
\begin{gather*}
    y+z < 1, \tag{$C_{2, 1}$} \label{regslaunit2,1}\\
    x+z < \frac{3}{2}, \ \ z-y < 1, \ \ y+z > 1, \tag{$C_{2, 2}$} \label{regslaunit2,2}\\
    x+z < \frac{3}{2}, \ \ z -y > 1, \tag{$C_{2, 3}$} \label{regslaunit2,3}\\
    x+y < \frac{3}{2}, \ \ x+z > \frac{3}{2}, \ \ z-y < 1, \ \ z-x < \frac{3}{2}, \tag{$C_{2, 4}$} \label{regslaunit2,4}\\
    x+y < \frac{3}{2}, \ \ z-y > 1, \ \ z < \alpha, \ \ x+z > \frac{3}{2}, \ \ z-x < \frac{3}{2}, \tag{$C_{2, 5}$} \label{regslaunit2,5}\\
    x+y < \frac{3}{2}, \ \ z > \alpha, \ \ z-y < 1, \tag{$C_{2, 6}$} \label{regslaunit2,6}\\
    z-x > \frac{3}{2}, \ \ z-y > 1, \ \ z < \alpha, \tag{$C_{2, 7}$} \label{regslaunit2,7}\\
    x < \alpha - \frac{1}{2}, \ \ z < \alpha, \ \ x+y > \frac{3}{2}, \ \ z-y <1 , \ \ z-x < \frac{3}{2}, \tag{$C_{2, 8}$} \label{regslaunit2,8}\\
    x < \alpha - \frac{1}{2}, \ \ y < \alpha, \ \ z-x < \frac{3}{2}, \ \ z > \alpha, \ \ z-y < 1, \tag{$C_{2, 9}$} \label{regslaunit2,9}\\
     x > \alpha - \frac{1}{2}, \ \ z < \alpha,  \tag{$C_{2, 10}$} \label{regslaunit2,10}\\
    x > \alpha - \frac{1}{2}, \ \ y < \alpha, \ \ z > \alpha , \ \ z- y < 1, \tag{$C_{2, 11}$} \label{regslaunit2,11}\\
    \alpha - \frac{1}{2} < x < \alpha + \frac{1}{2}, \ \ z-y < 1, \ \ y > \alpha, \ \ z-x < \frac{3}{2}, \tag{$C_{2, 12}$} \label{regslaunit2,12}\\
    z-y > 1, \ \ z < \alpha, \ \ x+y > \frac{3}{2}, \tag{$C_{2, 13}$} \label{regslaunit2,13}\\
    x > \alpha -\frac{1}{2}, \ \ y  <\alpha, \ \ z-y > 1, \ \ z-x < \frac{3}{2}, \tag{$C_{2, 14}$} \label{regslaunit2,14}\\
    z-x > \frac{3}{2}, \ \ z < \alpha, \ \ x+y > \frac{3}{2}, \ \ y-x < \frac{3}{2}, \tag{$C_{2, 15}$} \label{regslaunit2,15}\\
    z-x > \frac{3}{2}, \ \ y < \alpha, \ \ x+y > \frac{3}{2}, \ \ z > \alpha, \ \ z-y < 1, \ \ y-x < \frac{3}{2}, \tag{$C_{2, 16}$} \label{regslaunit2,16}\\
    x < \alpha -\frac{1}{2}, \ \ z-x > \frac{3}{2}, \ \ y-x < \frac{3}{2}, \ \ z-y < 1, \ \ y > \alpha, \tag{$C_{2, 17}$} \label{regslaunit2,17}\\
    \alpha -\frac{1}{2} < x < \alpha + \frac{1}{2}, \ \ z-x > \frac{1}{2}, \ \ z-y < 1, \ \ y-x < \frac{3}{2}, \tag{$C_{2, 18}$} \label{regslaunit2,18}\\
    y-x > \frac{3}{2}, \ \ z < \alpha, \tag{$C_{2, 19}$} \label{regslaunit2,19}\\
    y-x > \frac{3}{2}, \ \ y < \alpha, \ \ z > \alpha, \ \ z-y < 1. \tag{$C_{2, 20}$} \label{regslaunit2,20}
\end{gather*}
One can verify that $C_{2,i} \subset \cup_{j=1}^{6} L_{j}'$ for $i = 1,3,5,7$. Note that (\ref{regslaunit2,7}) is exactly the same as (\ref{finalslantunit4}) when $\alpha = 2$. By considering the critical points described in Proposition \ref{subquotientlist} and using Step $(3$-$2)$ of Algorithm \ref{alg A bar}, we can prove the non-unitarity for all regions $C_{2,i}$, except when $i = 4$. 

For $i = 4$, we consider the point $(\frac{1}{2}, 1,1)$, which lies on the boundary of (\ref{regslaunit2,4}). When $\alpha = 2$, the point 
\[\iota(\frac{1}{2}, 1,1) = (0,1,1,1) \in \R^4_{++}\]
is not strongly unitary, since it contains the representation 
\[L(\Delta_\rho[-1,-1], \Delta_\rho[-1,-1], \Delta_\rho[0,-1]; \pi_{sc}),\]
which is not of Arthur type and hence not unitarizable by Proposition \ref{unitiffArthur}. Suppose (\ref{regslaunit2,4}) is a unitary connected component in $\R^3_{\reg, \shrt, ++}$, then $\iota$(\ref{regslaunit2,4}) must be lie in the $W$-orbit of a unitary connected component in $R^4$. However, this is impossible since the point $(0,1,1,1)$ lies on the boundary of (\ref{regslaunit2,4}). Therefore (\ref{regslaunit2,4}) is non-unitary. This proves the Proposition for $\alpha = 2$. 

Two cases remain. Let us now consider the case $\alpha = \frac{3}{2}$. The list of possible bounded connected components in $\R^3_{\reg,\shrt,++}$ in this case is as follows: 
\begin{gather*}
    y+z < 1, \tag{$C_{\frac{3}{2}, 1}$} \label{regslaunit3/2,1}\\
    x+z < \frac{3}{2}, \ \ z-y < 1, \ \ y+z > 1, \tag{$C_{\frac{3}{2}, 2}$} \label{regslaunit3/2,2}\\
    x+z < \frac{3}{2}, \ \ z-y > 1, \tag{$C_{\frac{3}{2}, 3}$} \label{regslaunit3/2,3}\\
    x+y < \frac{3}{2}, \ \ z < \alpha, \ \ x+z > \frac{3}{2}, \ \ z-y < 1, \tag{$C_{\frac{3}{2}, 4}$} \label{regslaunit3/2,4}\\
    x+y < \frac{3}{2}, \ \ z-x < \frac{3}{2}, \ \ z > \alpha, \ \ z-y < 1, \tag{$C_{\frac{3}{2}, 5}$} \label{regslaunit3/2,5} \\
    z-y > 1, \ \ z < \alpha, \ \ x+z > \frac{3}{2}, \tag{$C_{\frac{3}{2}, 6}$} \label{regslaunit3/2,6}\\
    x+y < \frac{3}{2}, \ \ z-y > 1, \ \ z-x < \frac{3}{2}, \ \ z > \alpha, \tag{$C_{\frac{3}{2}, 7}$} \label{regslaunit3/2,7}\\
    x+y < \frac{3}{2}, \ \ z-x > \frac{3}{2}, \ \ z-y < 1, \tag{$C_{\frac{3}{2}, 8}$} \label{regslaunit3/2,8}\\
    x < \alpha - \frac{1}{2}, \ \ z < \alpha, \ \ x+y > \frac{3}{2}, \tag{$C_{\frac{3}{2}, 9}$} \label{regslaunit3/2,9}\\
    x < \alpha - \frac{1}{2}, \ \ y < \alpha, \ \ x+y > \frac{3}{2}, \ \ z > \alpha, \ \ z-y < 1, \ \ z-x < \frac{3}{2}, \tag{$C_{\frac{3}{2}, 10}$} \label{regslaunit3/2,10}\\
    x > \alpha - \frac{1}{2}, \ \ z < \alpha, \tag{$C_{\frac{3}{2}, 11}$} \label{regslaunit3/2,11}\\
    x > \alpha - \frac{1}{2}, \ \ y < \alpha, \ \ z > \alpha, \ \ z-y < 1, \tag{$C_{\frac{3}{2}, 12}$} \label{regslaunit3/2,12}\\
    \alpha - \frac{1}{2} < x < \alpha + \frac{1}{2}, \ \ z-x < \frac{3}{2}, \ \ z-y < 1, \ \ y > \alpha, \tag{$C_{\frac{3}{2}, 13}$} \label{regslaunit3/2,13}\\
    x > \alpha - \frac{1}{2}, \ \ y < \alpha, \ \ z-y > 1, \ \ z-x < \frac{3}{2}, \tag{$C_{\frac{3}{2}, 14}$} \label{regslaunit3/2,14}\\
    z-x > \frac{3}{2}, \ \ y < \alpha, \ \ x+y > \frac{3}{2}, \ \ z-y < 1, \tag{$C_{\frac{3}{2}, 15}$} \label{regslaunit3/2,15}\\
    z-x > \frac{3}{2}, \ \ x < \alpha - \frac{1}{2}, \ \ y-x < \frac{3}{2}, \ \ z-y < 1, \ \ y > \alpha, \tag{$C_{\frac{3}{2}, 16}$} \label{regslaunit3/2,16}\\
    z-x > \frac{3}{2}, \ \ \alpha -\frac{1}{2} < x < \alpha + \frac{1}{2}, \ \ y-x < \frac{3}{2}, \ \ z-y < 1. \tag{$C_{\frac{3}{2}, 17}$} \label{regslaunit3/2,17}\\
\end{gather*}
In this case, we have $C_{\frac{3}{2},i} \in \cup_{j=1}^{3}L_{j}'$ if and only if $i=1$. From Proposition \ref{subquotientlist}, there are $3$ points that get mapped to not strongly unitary points in $\R^4_{++}$ under $\iota$: 
\[(\alpha + \frac{1}{2}, \alpha +2, \alpha +3), \ \ (\alpha - \frac{1}{2}, \alpha +1, \alpha+2), \ \ (\alpha - \frac{1}{2}, \alpha, \alpha+1).\]

These points lie on the boundary of $C_{\frac{3}{2},i}$ for $i = 10, 12, 13, \ldots 17$. For $i = 3,6$, $C_{\frac{3}{2},i}$ has nonempty intersection with the hyperplane $x = 0$. Applying unitary reduction, one can conclude that (\ref{regslaunit3/2,3}) and (\ref{regslaunit3/2,6}) are unitary$^{\pm}$. Similarly one can also show that (\ref{regslaunit3/2,2})),(\ref{regslaunit3/2,4}), (\ref{regslaunit3/2,8}) are not unitary$^{\pm}$. It remains to show that (\ref{regslaunit3/2,4})), (\ref{regslaunit3/2,7})), (\ref{regslaunit3/2,9})) are not unitary$^{\pm}$. 

For (\ref{regslaunit3/2,4})), consider the point $(0,\frac{3}{2}, \frac{3}{2})$ which lies on its boundary. We have that 
\[\iota(0,\frac{3}{2}, \frac{3}{2}) = (-\frac{1}{2}, \frac{1}{2}, \frac{3}{2}, \frac{3}{2}),\]
and 
\[\Pi_{(-\frac{1}{2}, \frac{1}{2}, \frac{3}{2}, \frac{3}{2})} \cong \Pi_{(\frac{1}{2}, \frac{1}{2}, \frac{3}{2}, \frac{3}{2})}.\]
contains non-unitarizable subquotients by Proposition \ref{subquotientlist} and \ref{unitiffArthur}. Therefore (\ref{regslaunit3/2,4}) must be non-unitary. 

Similarly, for (\ref{regslaunit3/2,7})), consider the point $(0,\frac{1}{2}, \frac{3}{2})$ which lies on its boundary and gets mapped to $(-\frac{1}{2}, \frac{1}{2}, \frac{1}{2}, \frac{3}{2})$ under $\iota$. 
We have that 
\[\Pi_{(-\frac{1}{2}, \frac{1}{2}, \frac{1}{2}, \frac{3}{2})} \cong \Pi_{(\frac{1}{2}, \frac{1}{2}, \frac{1}{2}, \frac{3}{2})},\]
which contains non-unitarizable subquotients by Proposition \ref{subquotientlist} and \ref{unitiffArthur}. Therefore (\ref{regslaunit3/2,7}) must be non-unitary. 

Now consider the point $(1, \frac{3}{2}, \frac{3}{2})$. This point lies on the boundary of (\ref{regslaunit3/2,9}) but 
\[\iota(1,\frac{3}{2},\frac{3}{2}) = (\frac{1}{2}, \frac{3}{2}, \frac{3}{2}, \frac{3}{2})\]
contains non-unitarizable subquotients by Proposition \ref{subquotientlist} and \ref{unitiffArthur}. Therefore (\ref{regslaunit3/2,9}) is not unitary$^{\pm}$. 

This proves the claim for $\alpha = \frac{3}{2}$. Now suppose $\alpha = 1$. In this case, the list of possible bounded connected components of $\R^3_{\reg,\shrt,++}$ is as follows: 
\begin{gather*}
    y+z < 1, \tag{$C_{1, 1}$} \label{regslaunit1,1}\\
    x < \alpha -\frac{1}{2}, \ \ z < \alpha, \ \ y+z > 1, \tag{$C_{1, 2}$} \label{regslaunit1,2}\\
    x+z < \frac{3}{2}, \ \ y < \alpha, \ \ z > \alpha, \ \ z-y < 1, \tag{$C_{1, 3}$} \label{regslaunit1,3}\\
    x+z < \frac{3}{2}, \ \ y > \alpha, \tag{$C_{1, 4}$} \label{regslaunit1,4}\\
    x +z < \frac{3}{2}, \ \ x > \alpha - \frac{1}{2}, \tag{$C_{1, 5}$} \label{regslaunit1,5}\\
    x+z < \frac{3}{2}, \ \ z-y > 1, \tag{$C_{1, 6}$} \label{regslaunit1,6}\\
    x < \alpha -\frac{1}{2}, \ \ y < a, \ \ x+z > \frac{3}{2}, \ \ z-y < 1, \ \ z-x < \frac{3}{2}, \tag{$C_{1, 7}$} \label{regslaunit1,7}\\
    x+y < \frac{3}{2}, \ \ x+z > \frac{3}{2}, \ \ z-x < \frac{3}{2}, \tag{$C_{1, 8}$} \label{regslaunit1,8}\\
    x+y < \frac{3}{2}, \ \ z < \alpha, \ \ x+z > \frac{3}{2}, \tag{$C_{1, 9}$} \label{regslaunit1,9}\\
    x+y < \frac{3}{2}, \ \ x > \alpha -\frac{1}{2}, \ \ z > \alpha, \ \ z-y < 1, \tag{$C_{1, 10}$} \label{regslaunit1,10}\\
    x+y < \frac{3}{2}, \ \ z-y > 1, \ \ x > \alpha - \frac{1}{2}, \ \ z-x < \frac{3}{2}, \tag{$C_{1, 11}$} \label{regslaunit1,11} \\
    z-x > \frac{3}{2}, \ \ y < \alpha, \ \ z-y <1, \tag{$C_{1, 12}$} \label{regslaunit1,12}\\
    x+y < \frac{3}{2}, \ \ z-x > \frac{3}{2}, \ \ z-y < 1, \ \ y> \alpha, \tag{$C_{1, 13}$} \label{regslaunit1,13} \\
    x+y > \frac{3}{2}, \ \ z < \alpha, \tag{$C_{1, 14}$} \label{regslaunit1,14}\\
    y < \alpha, \ \ x+y > \frac{3}{2}, \ \ z > \alpha, \ \ z-y < 1, \tag{$C_{1, 15}$} \label{regslaunit1,15}\\
    \alpha -\frac{1}{2} < x < \alpha + \frac{1}{2}, \ \ z-y < 1, \ \ y > \alpha, \ \ z-x < \frac{3}{2},\tag{$C_{1, 16}$} \label{regslaunit1,16} \\
    y < \alpha, \ \ z-y > 1, \ \ x+y > \frac{3}{2}, \ \ z-x < \frac{3}{2},\tag{$C_{1, 17}$} \label{regslaunit1,17}\\
    x < \alpha - \frac{1}{2}, \ \ z-x > \frac{3}{2}, \ \ x+y > \frac{3}{2}, \ \ z-y < 1, \ \ y-x < \frac{3}{2}, \tag{$C_{1, 18}$} \label{regslaunit1,18}\\
    \alpha - \frac{1}{2} < x < \alpha + \frac{1}{2}, \ \ z-x > \frac{3}{2}, \ \ z-y < 1, \ \ y-x < \frac{3}{2}, \tag{$C_{1, 19}$} \label{regslaunit1,19}\\
    z-x > \frac{3}{2}, \ \ x > \alpha + \frac{1}{2}, \ \ z-y < 1, \ \ y-x < \frac{3}{2}. \tag{$C_{1, 20}$} \label{regslaunit1,20}\\
\end{gather*}
By Proposition \ref{subquotientlist}, the following critical points are not strongly unitary:
\[(\alpha + \frac{1}{2}, \alpha +2, \alpha +3), (\alpha - \frac{1}{2}, \alpha +1, \alpha+2),
(\alpha - \frac{1}{2}, \alpha, \alpha+1), (\alpha -\frac{1}{2}, \alpha, \alpha).\]
One can easily verify that all the other regions contain at least one of the critical points above in their boundaries except for the region (\ref{regslaunit1,6}). 

This region has nonempty $3$-dimensional intersection with the hyperplane $x = 0$ inside $\R^{3}_{\reg,\shrt,++}$. Applying Step $(3$-$2)$ of Algorithm \ref{alg A bar}, we can conclude that the region is not unitary$^{\pm}$, since for example, the point $(0,\frac{4}{3})$ does not lie in the unitary dual of corank $2$, but $(0,0,\frac{4}{3}) \in$ (\ref{regslaunit1,6}). 

Now let $\alpha = \frac{1}{2}$. By Proposition \ref{subquotientlist}, the following points are not strongly unitary: 
\[(\alpha + \frac{1}{2}, \alpha+1, \alpha+1), (\alpha+\frac{1}{2}, \alpha+1, \alpha+2),(\alpha + \frac{1}{2}, \alpha+1, \alpha+2),(\alpha + \frac{1}{2}, \alpha+2, \alpha+3).\]
Using these points, by the same method as before, we can eliminate all connected components except for the following: 
\begin{gather*}
    z < \frac{1}{2}, \tag{$C_{1/2, 1}$} \label{regslaunit1/2,1}\\   
    y + z < 1, \ \ z > \frac{1}{2}, \tag{$C_{1/2, 2}$}\label{regslaunit1/2,2}\\
    x+z < \frac{3}{2}, \ \ y < \frac{1}{2}, \ \ z-y < 1, \ \ y+z >1, \tag{$C_{1/2, 3}$}\label{regslaunit1/2,3}\\
    x+z < \frac{3}{2}, \ \ y > \frac{1}{2}, \tag{$C_{1/2, 4}$}\label{regslaunit1/2,4}\\
    x+ z < \frac{3}{2}, \ \ z-y > 1, \tag{$C_{1/2, 5}$}\label{regslaunit1/2,5}\\
   y < \frac{1}{2}, \ \ x+ z > \frac{3}{2}, \ \ z-y < 1, \tag{$C_{1/2, 6}$}\label{regslaunit1/2,6}\\
   x+y < \frac{3}{2}, \ \ z-y < 1, \ \ y > \frac{1}{2}, \ \ z-x > \frac{3}{2},\tag{$C_{1/2, 7}$}\label{regslaunit1/2,7} \\
    y < \frac{1}{2}, \ \ z-y > 1, \ \ x+z > \frac{3}{2}, \ \ z-x < \frac{3}{2},\tag{$C_{1/2, 8}$}\label{regslaunit1/2,8} \\
    x+ y< \frac{3}{2}, \ \ z-y > 1, \ \ z- x < \frac{3}{2}, \ \ y> \frac{1}{2}, \tag{$C_{1/2, 9}$} \label{regslaunit1/2,9}\\
    x+y < \frac{3}{2}, \ \ z-x > \frac{3}{2}, \ \ z-y < 1. \tag{$C_{1/2, 10}$} \label{regslaunit1/2,10}
\end{gather*}
Using Step ($3$-$2$) of Algorithm \ref{alg A bar} and Proposition \ref{unitcrnk2}, we can conclude that ($C_{1/2,i})$ is not unitary$^{\pm}$ for $i = 2,3,4,5,7,10$. 

For $i = 6,8,9$, consider the point $(0,\frac{1}{2}, \frac{3}{2})$ which lies on its boundary of ($C_{1/2,i}$). Under $\iota$, this gets mapped to 
\[\Pi_{(-\frac{1}{2}, \frac{1}{2}, \frac{1}{2}, \frac{3}{2})} \cong \Pi_{(\frac{1}{2}, \frac{1}{2}, \frac{1}{2}, \frac{3}{2})}.\]
By Proposition \ref{subquotientlist}, this point is not strongly unitary, and thus ($C_{1/2,i}$) are not unitary$^{\pm}$ for $i = 6,8,9$. One can easily prove that (\ref{regslaunit1/2,1}) is unitary by mapping it to the boundary of (\ref{ciu8}). 

Finally, let $\alpha = 0$. In this case, by Proposition \ref{subquotientlist}, the points
\[(\alpha + \frac{1}{2}, \alpha+1, \alpha+1), (\alpha + \frac{1}{2}, \alpha+1, \alpha+2), (\alpha + \frac{1}{2}, \alpha+2, \alpha+2), (\alpha + \frac{1}{2}, \alpha+2, \alpha+3)\]
are not strongly unitary. Using these points, we can eliminate all connected components other than 
\begin{gather*}
    y + z < 1, \tag{$C_{0,1}$} \label{regslaunit0,1} \\
    x + z < \frac{3}{2}, \ \ z - y > 1 \tag{$C_{0,2}$} \label{regslaunit0,2}
\end{gather*}
By applying Step ($3$-$2$) of Algorithm \ref{alg A bar} and using Proposition \ref{unitcrnk2}, one can show that (\ref{regslaunit0,2}) is not unitary$^{\pm}$ and (\ref{regslaunit0,1}) is unitary$^{\pm}$. This concludes the proof.
\end{proof}

Finally, we can list all the unitary$^{\pm}$ connected components of $\R^3_{v\reg,\shrt,++}$ as follows. Note that when $\alpha = 0$ or $\frac{1}{2}$, the unitary$^{\pm}$ connected components of $\R^3_{v\reg,\shrt,++}$ are exactly the unitary$^{\pm}$ components listed in Proposition \ref{finalslantunit}. 

\begin{prop}\label{finalslavregunit}
    For $\alpha \geq 1$, the unitary$^{\pm}$ connected components of $\R^3_{v\reg,\shrt,++}$ are exactly as follows: 
    \begin{gather*}
        x + z< \frac{1}{2}, \tag{$L_{1}$} \label{finalvregslaunit1}\\
        y + z < 1, \ \ x+z > \frac{1}{2}, \ \ z-x < \frac{1}{2}, \tag{$L_{2}$} \label{finalvregslaunit2}\\
        y+z < 1, \ \ z-x > \frac{1}{2}, \tag{$L_{3}$} \label{finalvregslaunit3} \\
        x+ z < \frac{3}{2}, \ \ z-y > 1, \ \ (\alpha > 1) \tag{$L_{4}$} \label{finalvregslaunit4} \\
        x < \frac{1}{2}, \ \ z-y > 1, \ \ x+z > \frac{3}{2}, \ \ z-x < \frac{3}{2}, \ \ (\alpha > 1) \tag{$L_{5}$} \label{finalvregslaunit5}\\
        x > \frac{1}{2}, \ \ x+y < \frac{3}{2}, \ \ z-y > 1, \ \ z < \alpha, \ \  (\alpha > \frac{3}{2}) \tag{$L_{6}$} \label{finalvregslaunit6}\\
        x < \frac{1}{2}, \ \ z-x > \frac{3}{2}, \ \ y-x < \frac{1}{2}, \ \ z < \alpha, \ \ (\alpha > \frac{3}{2}) \tag{$L_{7}$} \label{finalvregslaunit7}\\
        x + y < \frac{3}{2}, \ \ y-x > \frac{1}{2}, \ \ z-y > 1, \ \ z < \alpha, \ \  (\alpha > \frac{3}{2})\tag{$L_{8}$} \label{finalvregslaunit8}\\
        x+y < \frac{3}{2}, \ \ z -x > \frac{3}{2}, \ \ z < \alpha, \ \ x > \frac{1}{2}, \ \ (\alpha > 2) \tag{$L_{9}$} \label{finalvregslaunit9}\\
        x < \frac{1}{2}, \ \ y-x > \frac{3}{2}, \ \ z-y > 1, \ \ z < \alpha, \ \ (\alpha > \frac{5}{2}) \tag{$L_{10}$} \label{finalvregslaunit10}\\
        x > \frac{1}{2}, \ \ y-x > \frac{3}{2}, \ \ z-y > 1, \ \ z < \alpha \ \ (\alpha > 3). \tag{$L_{11}$} \label{finalvregslaunit11}
    \end{gather*}
\end{prop}
\begin{proof}
    The regions above are obtained by slicing up the regions (\ref{finalslantunit1}) to (\ref{finalslantunit5}) using the reducibility hyperplanes inside $R^{3}_{v\reg,\shrt,++}$. The result follows directly from Proposition \ref{finalslantunit}. 
\end{proof}

\subsection{Unitarizability for the irregular components - level hyperplanes}\label{levelsection}
In this subsection we perform the inductive process described in Step $(3$-$3)$ of Algorithm \ref{alg A bar}. In particular, we consider reducibility hyperplanes which are contained in the $W$-orbit of $H_{\lev}$, where
\[
H_{\lev} = \{\underline{x} \in \R^4: x_4 = \alpha\}.
\]
We call these the set of {\it{level hyperplanes}}. Similar to the slanted hyperplanes case, for any $\underline{x} \in H_{\lev}$, we can decompose $\Pi_{(x,y,z,\alpha)}$ in the Grothendieck group as $\tau_{(x,y,z)}^{+} + \tau_{(x,y,z)}^{-}$, where 
\begin{flalign*}
    \tau_{(x,y,z)}^{+} &= \Delta_\rho[-x,-x] \times \Delta_\rho[-y,-y] \times \Delta_\rho[-z, -z] \rtimes T_{I,1}^{\alpha}(\pi_{sc}), \\
    \tau_{(x,y,z)}^{-} &= \Delta_\rho[-x,-x] \times \Delta_\rho[-y,-y] \times \Delta_\rho[-z, -z] \rtimes L(\Delta_\rho[-\alpha, -\alpha]; \pi_{sc}),
\end{flalign*}
when $\alpha \neq 0$, and
\[\tau_{(x,y,z)}^{\pm} = \Delta_\rho[-x,-x] \times \Delta_\rho[-y,-y] \times \Delta_\rho[-z, -z] \rtimes T_{V,2}^{\pm}(\pi_{sc})\] when $\alpha = 0$. 

Similar to before, define $\R^3_{\reg,\lev}$ to be the complement of the singular affine hyperplanes other than $H_{\lev}$. In other words, if we identify $(x,y,z)$ with $(x_1, x_2, x_3)$, it is the complement of the hyperplanes
\begin{flalign*}
    x_i \pm x_j = \epsilon, \ x_i = \pm(\alpha + \epsilon), \ \ \epsilon = \pm1, \ \ i = 1,2,3, \ \ i< j.
\end{flalign*}

We say that a point $(x,y,z) \in \mathbb{R}^3_{\reg, \lev}$ is unitary$^{+}$ (resp., unitary $^{-}$) if $\pi_{(x,y,z)}^{+}$ (resp., $\pi_{(x,y,z)}^{-}$ is unitarizable). We say that a point $(x,y,z) \in \mathbb{R}^3_{\reg,\lev}$ is unitary$^{\pm}$ if it is both unitary$^{+}$ and unitary$^{-}$. We say that a connected component is unitary$^{+}$ (resp., unitary$^-$, unitary$^{\pm}$) if every point in it is unitary$^+$ (resp., unitary$^{-}$, unitary$^\pm$). Then, $\tau^{+}_{(x,y,z)}$ and $\tau^{-}_{(x,y,z)}$ are irreducible precisely when $(x,y,z) \in \R^{3}_{\reg, \lev}$. Therefore, for all $(x,y,z) \in \R^3_{\reg, ++}$, $\iota(x,y,z)$ is strongly unitary (resp., strongly non-unitary) if and only if both $\tau^{+}_{(x,y,z)}$ and $\tau^{-}_{(x,y,z})$ are unitarizable (resp., non-unitarizable). 

Let $\R^3_{\reg, \lev, +}= \R^3_{+} \cap \R^3_{\reg, \lev}, \R^3_{\reg, \lev, ++}= \R^3_{++} \cap \R^3_{\reg, \lev}$. The following proposition describes the unitary$^{\pm}$ connected components of $\R^3_{\reg,\lev, ++}$. 

\begin{prop}\label{finallevelunit}
    The unitary$^{\pm}$ connected components of $\R^3_{\reg, \lev, ++}$ are as follows: 
    \begin{gather}
    (\alpha > 3) \ \ y-x > 1, \ \ z - y > 1, \ \ z < \alpha -1, \label{level3++ucc,1}\\
    (\alpha > 2) \ \ x+y < 1, \ \ z-y > 1, \ \ z < \alpha -1, \label{level3++ucc,2} \\
    (\alpha \geq 1) \ \ y+z < 1, \ \ z < \alpha-1, \label{level3++ucc,3}\\
    (\alpha > \frac{3}{2}) \ \ x+y < 1, \ \ x+z > 1, \ \ z-x < 1, \ \ z < \alpha -1, \label{level3++ucc,4}\\
    (\alpha = \frac{1}{2}) \ \  z < \frac{1}{2}, \label{level3++ucc,5}\\
    (\alpha = 2) \ \ x+y < 1, \ \ z-x < 1, \ \ z > \alpha -1, \label{level3++ucc,6}\\
    (\alpha = \frac{3}{2}) \ \ y+ z < 1, \ \ z > \alpha -1. \label{level3++ucc,7}\\
    (\alpha = 0) \ \ y+z < 1 \label{level3++ucc,8}
    \end{gather}  
\end{prop}

\begin{proof}
    The connected components (\ref{level3++ucc,3}),  (\ref{level3++ucc,5}), (\ref{level3++ucc,8}) are precisely the connected components that contain the origin, in the cases $\alpha \geq 1, \alpha = \frac{1}{2}, \alpha = 0$ respectively, so they are unitary$^{\pm}$. 

    Now let us assume $\alpha > \frac{3}{2}$. Define the representations
    \begin{flalign*}
        \sigma^{+} &= T_{I,1}^{\alpha}(\pi_{sc}), \\
        \sigma^{-} &= L(\Delta_\rho[-\alpha, -\alpha]; \pi_{sc}).
    \end{flalign*}
    Consider the complementary series
    \[\pi_{z} = \Delta_\rho[-z,-z] \rtimes \sigma^{\pm},\]
    which is irreducible and unitary for $0 \leq z < \alpha-1$. 
    First let us look at the complementary series 
    \[
    \Delta_\rho[-x,-x] \times \Delta_\rho[x,x] \rtimes \pi_{z},
    \]
    for $\frac{1}{2} < z < \min(\frac{3}{2}, \alpha-1)$ and $|1-z| < x < \frac{1}{2}$. This gives us an irreducible unitarizable representation $\Pi_{(x,-x,z)} \cong \Pi_{(x,x,z)}$. Since (\ref{level3++ucc,4}) is the connected component of 
    \[\{(x,x,z): |1-z| < x < \frac{1}{2}, \ \ \frac{1}{2} < z < \min(\frac{3}{2}, \alpha-1) \} .\]
    It must be unitary$^{\pm}$.

    Now assume $\alpha > 2$. Fix $1 < z < \alpha-1$. Clearly the representation $\rho  \times \rho \rtimes \pi_z$ is irreducible and unitary. It follows that $(0,0,z)$ is unitary$^{\pm}$. Since (\ref{level3++ucc,2}) is the connected component of $(0,0,z)$ for $1 < z < \alpha-1$, it must be unitary$^{\pm}$. 

    Now assume $\alpha > 3$. First fix $2 < z < \alpha-1$. Construct the complementary series  
    \[\Delta_\rho[-y,-y] \rtimes \pi_z,\]
    for $0 \leq y < z-1$, which is irreducible and unitarizable. Now fix $1 < y < z-1$ and construct the complementary series 
    \[
    \Delta_\rho[-x,-x] \rtimes \Delta_\rho[-y,-y] \rtimes \pi_z,
    \]
    which is irreducible and unitarizable for $0 \leq x < y-1$. 
    Since (\ref{level3++ucc,1}) is the connected component containing such $(x,y,z)$, we may conclude that it is unitary$^{\pm}$. 

    Let $\alpha = 2$. Fix $1 < z < \frac{3}{2}$. Then the complementary series $\pi_z$ is irreducible and unitarizable. Now pick $z-1 < x < \frac{1}{2}$. Then the complementary series 
    \[\rho\lvert \cdot \rvert^{x} \times \rho \lvert \cdot \rvert^{x} \rtimes \pi_z\]
    is irreducible and unitarizable. Since such $(x,x,z)$ is contained in (\ref{level3++ucc,6}), the connected component must be unitary$^{\pm}$.

    Let $\alpha = \frac{3}{2}$. By applying Step $(3$-$2)$ of Algorithm \ref{alg A bar} and using Proposition \ref{unitcrnk2}, one sees that (\ref{level3++ucc,7}) is unitary$^{\pm}$. This proves that all connected components listed in the proposition are unitary$^{\pm}$. 

    The proof that these are the only unitary$^{\pm}$ components is similar to that of Proposition \ref{finalslantunit}, which we omit. 
\end{proof}

\subsection{\texorpdfstring{Final list of open unitary connected components}{}}\label{exhaustion}

In the previous two subsections, we've completed Step $3$ of Algorithm \ref{alg A bar}. In this subsection, we show that the list of open unitary connected components in Proposition \ref{connectedisunitary} is in fact the full list. First, let us recall Tadi{\'c}'s result on the unitary dual of corank $3$: 

\begin{prop} [{\cite[Proposition $8.3$]{Tad23}}] \label{crnk3unitcc}
The following connected components of $\R^3_{\reg,++}$ are unitary.

For $\alpha\ge1:$
\begin{subequations}
\begin{gather}
\label{eq: xx2x31} x_2+x_3<1,\\
\label{eq: xx1x212} x_1+x_2<1, \ \ , x_3-x_2>1, \ \ , x_3<\alpha,\ \ \ (\alpha>1)\\
\label{eq: xx1x213} x_1+x_2<1, \ \ , x_1+x_3>1, \ \, x_3-x_1<1,\ \ , x_3<\alpha,\\
\label{eq: xx1x2x14} x_2-x_1>1, \ \ , x_3-x_2>1, \ \, x_3<\alpha\ \ \ (\alpha>2).
\end{gather}
\end{subequations}
(The constraint $x_3<\alpha$ in \eqref{eq: xx1x213} is redundant unless $\alpha=1$.)

For $\alpha=\tfrac12:$
\begin{equation} \label{eq: x312}
x_3<\tfrac12.
\end{equation}

Consequently, any $\mathbf{x}\in\R^3_{++}$ in the closure of the above regions (i.e., changing strict inequalities
to non-strict ones) is strongly unitary.
\end{prop}

By \cite[Proposition $8.12$]{Tad23}, this list is exhaustive. Now we are ready to prove the same for our list in corank $4$. 

\begin{prop} \label{exhaustconclusion}
    The open unitary connected components of $\R^4_{++}$ are exactly those listed in Proposition \ref{connectedisunitary}.
\end{prop}

\begin{proof}
    Let $C$ be an open unitary connected component of $\R^4_{\reg, ++}$. Let $\mathscr{C}$ be the unitary connected component of $\R^4_{++}$ containing $\mathscr{C}$. $C$ and $\mathscr{C}$ must be bounded, and their boundary must be contained in the union of all reducibility hyperplanes. In other words, any three dimensional volume in the boundary of $\mathscr{C}$ must be contained in a $W$-translate of $H_{\lev}$ or $H_{\shrt}$, since all slanted/level hyperplanes are in the same $W$-orbit of $H_{\shrt}/H_{\lev}$. 

    Consider the case $\alpha = 0$. If the boundary of $\mathscr{C}$ is contained entirely in the union of level hyperplanes, then $\mathscr{C}$ cannot be bounded. Therefore there exists some unitary$^{\pm}$ $3$-dimensional volume in the boundary of $\mathscr{C}$ contained in a slanted hyperplane. By Proposition \ref{finalslantunit},  the boundary must lie inside $W$-orbit of the component 
    \[x_2 - x_1 = 1, \ \  x_3 + x_4 < 1.\]
    We take the unique $w \in W$ such that the image of the above component under $w$ has nonempty intersection with $\R^{4}_{++}$, which is 
    \[w: (x_1, x_2, x_3, x_4) \to (x_3,x_4,-x_1,x_2) = (y_1, y_2, y_3, y_4),\]
    under which the transformed boundary becomes the region 
    \[y_1 + y_2 < 1, \ \ y_3 + y_4 = 1.\]
    The two connected components in $\R^4$ sharing this boundary are
    (\ref{ciu9}) and some unbounded region. This proves that when $\alpha = 0$, the only unitary connected component in $\R^4_{++}$ is (\ref{ciu9}). 
    
    For all $\alpha > 0$, using  Proposition \ref{finalslantunit} and \ref{finallevelunit}, we can conclude that 
    \[\mathscr{C} \subset \{x \in \R^4_{\reg}: |x_i| < \alpha, i = 1,2,3,4\}\]
    and 
    \[C \subset \{x \in \R^4_{\reg,++}: x_4< \alpha\},\]
    by passing them to the boundary. 
    For $\alpha = \frac{1}{2}$, we can conclude that $C$ is exactly the component given in Proposition \ref{connectedisunitary}. 

    For the rest of the proof we consider the case $\alpha \geq 1$. Now suppose some $3$-dimensional volume on the boundary of $\mathscr{C}$ lies in a level hyperplane. Up to $W$-translation,  we may assume this to be $H_{\lev}$. It suffices to consider the unitary$^{\pm}$ components of $R^3_{\reg,\lev,++}$. 

     By Proposition \ref{finallevelunit}, there are $6$ possible unitary$^{\pm}$ components in this case. Let us first consider the case $\alpha > 2$. The first possibility is the region 
     \[y-x > 1, \ \ z - y > 1, \ \ z < \alpha-1, \]
     for $\alpha > 3$. One can easily show that this is contained in the closure of the region (\ref{ciu4}), by considering the affine isomorphism $(x_1, x_2, x_3) \mapsto (x_1, x_2, x_3, \alpha)$. Moreover, (\ref{ciu4}) is the only connected component that contains (\ref{level3++ucc,1}) and is contained in the region $\{\underline{x} \in \R^4_{++}: x_4 < \alpha\}$. Thus we can conclude that $C$ is (\ref{ciu4}). Similarly, the region (\ref{level3++ucc,2}) described by 
     \[x+y < 1,  \ \ z-y > 1, \ \ z < \alpha -1, \]
     for $\alpha > 2$, is contained in the closure of (\ref{ciu5}). By the same reasoning, we may conclude that $C$ is (\ref{ciu5}). 

     Now let's look at the region (\ref{level3++ucc,3}), described by 
     \[y+ z< 1, \ \ z < \alpha -1, \]
     for $\alpha \geq 1$. For $\alpha > 1$, this is contained in the closure of (\ref{ciu2}). 

     Lastly, we consider the region (\ref{level3++ucc,4}), described by 
     \[x+y < 1, \ \ x+z > 1, \ \ z-x < 1, \ \ z < \alpha-1, \]
     for $\alpha > \frac{3}{2}$. This is contained in the closure of the connected component (\ref{ciu6}). By the same reasoning as above, we can conclude that $C$ can only be the connected component (\ref{ciu6}). 

     When $\alpha = 2$, (\ref{level3++ucc,6}) is contained entirely in the boundary of (\ref{ciu7}). When $\alpha = \frac{3}{2}$, (\ref{level3++ucc,7}) is contained entirely in the boundary of (\ref{ciu1}). In both cases, we can conclude that they are the only possibilities for $C$. 

     It remains to consider the case where the boundary of $C$ is contained entirely in slanted hyperplanes. Take $H$ to be a slanted hyperplane that contains a $3$-dimensional volume in the boundary of $C$. Fix $w_C \in W$ such that $w_C(H) = H_{\shrt}$. It follows that $w_C(\partial(C)) \in H_{\shrt}$ and 
     \[w_C \partial(C) = \partial (w_C(C)) \supset \iota(w_i L_i),\]
     for some $1 \leq i \leq 11$ and $w_i \in W$. This implies that 
     \[\partial w_C(C) \cap \R^4_{++} \supset \iota(w_i L_i) \cap \R^4_{++}.\]

     It suffices now to show the following: For all $1 \leq i \leq 11$, for all $w_i \in W$ such that $\iota(w_i L_i) \cap \R^4_{++}$ is nonempty, out of the two connected components in $\R^4_{++}$ which share the common boundary that contains $\iota(w_i L_i) \cap \R^4_{++}$, one is unitary if and only if it is listed in Proposition \ref{connectedisunitary}. In fact, by definition of $\R^4_{++}$, one can show that such $w_i$ is unique.  
     
     We begin with the case $\alpha \geq 1$. For the component (\ref{finalvregslaunit1}), the only possible $w_1 \in W$ such that 
     \[\iota(w_1 L_1) \cap \R^4_{++} \neq \varnothing,\]
     is the map
     \[w_1: (x_1, x_2, x_3, x_4) \mapsto (x_3,x_4, -x_1, x_2) = (y_1, y_2, y_3, y_4).\]
     Under this map, the resulting boundary can be described by 
     \[y_3 + y_4 = 1, \ \ 0 \leq \frac{1}{2} - y_3 \leq y_1 \leq y_2. \]
     Clearly this is the subset of the boundary of the connected component (\ref{ciu1}). Denote the only other connected component in $\R^4_{++}$ sharing this boundary by $C_{a,opp}$. Then $C_{a,opp}$ can be described by 
     \[x_2 + x_4 < 1, \ \ x_3 + x_4 > 1.\]
     This connected component has $4$-dimensional intersection with the hyperplanes $x_1 = 0$ and $x_2 = 0$. Applying Step $(3$-$2)$ of Algorithm \ref{alg A bar} twice, one can conclude that $C_{a,opp}$ is non-unitary, since, for example $(0,0,\frac{2}{3}, \frac{2}{3})$ is contained in $C_{a, opp}$, but $(\frac{2}{3}, \frac{2}{3})$ is not contained in the unitary dual of corank $2$. 

     The only possible $w_2$ is    
        \[w_2: (x_1, x_2, x_3, x_4) \mapsto (x_1, x_3, x_4, x_2) = (y_1, y_2, y_3, y_4).\]
     In this case, the transformed region in $\R^4_{++}$  can be described by 
    \[y_2 + y_3 < 1, \ \ y_4 - y_1 = 1, \ \ y_3 - y_1 < 1, \ \ 0 \leq y_1 + \frac{1}{2} \leq y_2 \leq y_3.\]
     This region has empty intersection with $\R^4_{++}$, so we may ignore it. 
     
     For $i = 3$, we must have: 
     \[w_3: (x_1, x_2, x_3, x_4) \mapsto (x_3, -x_1, x_2, x_4) = (y_1, y_2, y_3, y_4).\]
     Then the transformed region in $\R^4_{++}$ is 
     \[y_2 + y_3 = 1, \ \ y_1 + y_4 < 1, \ \ y_2 + y_4 > 1, \ \ 0 \leq y_3 - \frac{1}{2} \leq y_1 \leq y_4.\]
     There are two possible connected components with this boundary. Call them $C_{3,1}$ and $C_{3,2}$, defined by
     \[C_{3,1}: y_2 + y_3 < 1, \ \ y_1 + y_4 < 1, \ \ y_2 + y_4 > 1,\]
     \[C_{3,2}: y_2 + y_3 > 1, \ \ y_1 + y_4 < 1.\]

$C_{3,1}$ is the same as (\ref{ciu7}), which is unitary. $C_{3,2}$ has nontrivial $3$-dimensional intersection with the hyperplane $y_1 = 0$. The intersection can be described by 
\[y_2 + y_3 > 1, \ \ y_4 < 1.\]
One can easily verify that this region is not contained in any of the connected components listed in Proposition \ref{crnk3unitcc}, so $C_{3,2}$
is non-unitary by Step $(3$-$2)$ of Algorithm \ref{alg A bar}. 

Now let $\alpha > 1$. For $i = 4$, the only possible $w_4$ is 
\[w_4: (x_1, x_2, x_3, x_4) \mapsto (x_3, -x_1, x_2, x_4) = (y_1, y_2, y_3, y_4).\]
The transformed region in $\R^4_{++}$ is now
\[y_2 + y_3 = 1, \ \ y_4 - y_2 < 1, \ \ y_4 - y_1 > 1 \ \ 0 
\leq y_3 - \frac{1}{2} \leq y_1 \leq y_4.\]
This is contained in the boundary of (\ref{ciu2}). The other possible component containing this region in its boundary is 
\[x_1 + x_3 < 1, \ \ x_4 - x_1 > 1, \ \ x_4 - x_2 < 1, \ \ x_2 + x_3 > 1, \ \ x_4 < \alpha.\]
Using Step $(3$-$2)$ of Algorithm \ref{alg A bar} and Proposition \ref{crnk3unitcc}, this connected component is non-unitary. 

For $i = 5$, the only possibility is 
\[w_5: (x_1, x_2, x_3, x_4) \mapsto (x_3, -x_1, x_2, x_4) = (y_1, y_2, y_3, y_4).\]
The transformed region in $\R^4_{++}$  is 
\[y_2 + y_3 = 1, \ \ y_3 - y_2 < 1, \ \ y_4 - y_2 > 1, \ \ y_4 - y_3  < 1. \]
This is contained in the boundary of (\ref{ciu7}). The other possible component can be described by 
\[x_2 + x_3 > 1, \ \ x_3 - x_2 < 1, \ \ x_4 - x_2 > 1, \ \ x_4 - x_3 < 1.\]
By applying Step $(3$-$2)$ of Algorithm \ref{alg A bar}, one can conclude that this connected component is non-unitary. 

Now consider $\alpha > \frac{3}{2}$. We must have
\[w_6: (x_1, x_2, x_3, x_4) \mapsto (x_1, x_3, x_2, x_4) = (y_1, y_2, y_3, y_4).\]
The transformed region in $\R^4_{++}$  is 
\[y_3 - y_1 = 1, \ \ y_1 + y_2 < 1, \ \ y_4 - y_2 > 1, \ \ y_4 < \alpha, \ \ 0 \leq y_3 - \frac{1}{2} \leq y_1 \leq y_4.\]
This is contained in the boundary of (\ref{ciu6}). The other possible component can be described by 
\[x_1 + x_2 < 1, \ \ x_3 - x_1 > 1, \ \ x_4 - x_2 > 1, \ \ x_3 - x_2 < 1, \ \ x_4 < \alpha.\]
This is non-unitary by Step $(3$-$2)$ of Algorithm \ref{alg A bar} and Proposition \ref{crnk3unitcc}.

For $i = 7$, we have
\[w_7: (x_1, x_2, x_3, x_4) \mapsto (-x_1, x_3, x_2, x_4) = (y_1, y_2, y_3, y_4).\]
The transformed region in $\R^4_{++}$  is 
\[y_1 + y_3 = 1, \ \ y_4 - y_3 > 1, \ \ y_1 + y_2 < 1 \ \ y_4 < \alpha, 0 \leq y_3 - \frac{1}{2} \leq y_2 \leq y_4 .\]
This is contained in the boundary of (\ref{ciu6}). The other component containing with boundary containing the above region is 
\[x_1 + x_3 < 1, \ \ x_4 - x_3 > 1, \ \ x_2 + x_3 > 1, \ \ x_4 < \alpha.\]
which is non-unitary by the same reasoning as before. 

For $i = 8$, we have 
\[w_8: (x_1, x_2, x_3, x_4) \mapsto (-x_1, x_2, x_3, x_4) = (y_1, y_2, y_3, y_4).\]
The transformed region in $\R^4_{++}$  is 
\[y_1 + y_2 = 1, \ \ y_3 - y_1 < 1, \ \ y_1 + y_3 > 1, \ \ y_4 - y_3 > 1 \ \ y_4 < \alpha, \ \ 0 \leq y_2 - \frac{1}{2} \leq y_3 \leq y_4.\]
This is contained in the boundary of (\ref{ciu6}). The other component  with boundary containing the above region is 
\[x_1 + x_2 > 1, \ \ x_3 - x_1 < 1, \ \ x_1 + x_3 > 1, \ \ x_4 - x_3 > 1, \ \ x_4 < \alpha.\]
When $\alpha \in \mathbb{Z}$, then $(0,1,1,2)$ is a point of good parity contained in the boundary of the connected component above. By Proposition \ref{nontemp4}, this point contains subquotients that are not of Arthur type and hence it is not strongly unitary by Proposition \ref{unitiffArthur}. When $\alpha \in \frac{1}{2} + \mathbb{Z}$, then $(\frac{1}{2}, \frac{1}{2}, \frac{1}{2},\frac{3}{2})$ is a point of good parity in the boundary of the above component. It also contains subquotients that are not of Arthur type by Proposition \ref{nontemp4}. Thus,in both cases, we can conclude that the above component is non-unitary. 

Now we move onto the case $\alpha > 2$. For $i = 9$, we have 
\[w_9: (x_1, x_2, x_3, x_4) \mapsto (x_1, x_3, x_2, x_4) = (y_1, y_2, y_3, y_4).\]
The transformed region in $\R^4_{++}$  is 
\[y_3 - y_1 = 1, \ \ y_1 + y_2 < 1, \ \ y_4 - y_3 > 1, \ \ y_1 + y_3 > 1, \ \ y_4 < \alpha, \ \ 0 \leq y_3 - \frac{1}{2} \leq y_2 \leq y_4. \]
This is contained in the boundary of (\ref{ciu6}). The other component with boundary containing the above region is 
\[x_1 + x_2 < 1, \ \ x_3 - x_1 > 1, \ \ x_3 - x_2 < 1 \ \ x_4 - x_2 > 1, \ \ x_4 - x_3 < 1, \ \ x_4 < \alpha.\]
Using Step $(3$-$2)$ of Algorithm \ref{alg A bar} and Proposition \ref{crnk3unitcc}, we can conclude that the above component is non-unitary. 

Now let $\alpha > \frac{5}{2}$. For $i = 10$, we have
\[w_{10}: (x_1, x_2, x_3, x_4) \mapsto (-x_1, x_2, x_3, x_4) = (y_1, y_2, y_3, y_4).\]
The transformed region in $\R^4_{++}$  is 
\[y_1 + y_2 = 1, \ \ y_2 - y_1 < 1, \ \ y_3 - y_2 > 1, \ \ y_4 - y_3 > 1, \ \ 0 \leq y_2 - \frac{1}{2} \leq y_3 \leq y_4.\]
This is contained in the boundary of (\ref{ciu5}). The other component with boundary containing the above region is 
\[x_4 - x_3 > 1, \ \ x_1 + x_2 > 1, \ \ x_3 - x_1 < 1, \ \ x_4 < \alpha.\]
Using the same method as the case $i = 8$, we can prove show that the above component is non-unitary. 

Finally, let $\alpha > 3$. For $i = 11$, we have
\[w_{10}: (x_1, x_2, x_3, x_4) \mapsto (x_1, x_2, x_3, x_4).\]
In this case, the transformed region in $\R^4_{++}$ is 
\[x_2 - x_1 = 1, \ \ x_1 + x_2 > 1, \ \ x_3 - x_2 > 1, \ \ x_4 - x_3 > 1, \ \ x_4 < \alpha, \ \ 0 \leq x_2 - \frac{1}{2} \leq x_3 \leq x_4. \]
This is contained in the boundary of (\ref{ciu4})). The other component with boundary containing the above region is 
\[x_1 + x_2 > 1, \ \ x_2 - x_1 < 1, \ \ x_3 - x_2 > 1, \ \ x_4 - x_3 > 1, \ \ x_4 < \alpha.\]
When $\alpha \in \mathbb{Z}$, the point $(1,1,2,3)$ is a point of good parity contained in the boundary of the above component. By Proposition \ref{nontemp4}, the point contains subquotients that are not of Arthur type and hence, by Proposition \ref{unitiffArthur}, it is not strongly unitary. Similarly, when $\alpha \in \frac{1}{2} + \mathbb{Z}$, the point $(\frac{3}{2}, \frac{3}{2}, \frac{5}{2}, \frac{7}{2})$ is a point of good parity contained in the boundary. It contains subquotients that are not of Arthur type and hence is not strongly unitary. In both cases, we can conclude that the above component is non-unitary. This concludes the proof. 
\end{proof}

\begin{remark} \label{n-1dcompseries}
    It is clear from our proof that any 3-parameter complementary series within the corank $4$ unitary dual is contained in the boundary of some open unitary connected component in $\R^4$. We expect this to hold for any arbitrary corank $n$, i.e., within the unitary dual of corank $n$, any $n-1$-dimensional complementary series should be fully contained in the boundary of some open unitary connected component in $\R^n$.  
\end{remark}

\section{One-parameter unitarizable families}\label{1-parameterseries}

From the last section, we have all unitarizable representations of corank $4$ which appear as part of the $4$ or $3$-dimensional complementary series. To construct the full unitary dual, we also need to include unitarizable representations of corank $4$ that appear as part of a unitarizable family with one or two parameters, as described in Step $2$ of Algorithm \ref{alg A bar}. In this section, we classify all one-parameter continuous families of irreducible, unitarizable representations of corank $4$. 

 \subsection{One-parameter complementary series }
We begin by listing representations inside a one-parameter complementary series, which is constructed from a critical type, unitarizable, irreducible representation of corank up to $3$ (see \cite{Tad23} for the complete list). From now on, we say that an induced representation is unitary if all of its irreducible subquotients are unitary. We begin by considering all complementary series of the form 
\[u_{\rho}(a,b)\lvert \cdot \rvert^{x} \rtimes \pi_{A}.\]

\begin{prop}\label{1dcompleseries}
    In Table \ref{tab: redpnt} below, we list all possible one parameter complementary series of the form 
    \[\pi_{x} = u_{\rho}(a,b)\lvert \cdot \rvert^{x} \rtimes \pi_{A},\]
    where 
    $\pi_{A} \in \Pi_{A,gp}$ is an irreducible critical unitarizable representation of a smaller rank group $G_m$. For each $\pi_{A}$ and its dual (denoted by case $'$), we list all the corresponding families $\pi_x$ for $x \geq 0$ and all the reducibility points. When irreducible, $\pi_{x}$, $x \geq 0$ is unitarizable up to the first nonzero reducibility point. When reducibility occurs at $0$, there are no such complementary series. 
\end{prop}

\begin{proof}
   The statement is clear when there is only one reducibility point for $\pi_x$. When there are $2$ or $3$ reducibility points, then it suffices to show that $\pi_{x_0}$ contains a non-unitarizable subquotient, when $x_0$ is the second nonzero reducibility point. When there are $4$ reducibility points, then it suffices to show that $\pi_{x}$ contains a non-unitarizable subquotient at the second and third nonzero reducibility point, or the second and fourth nonzero reducibility point. All the relevant information is given in Table \ref{tab: nonunitsub1} below, where the particular non-unitarizable subquotient is denoted by $\pi$. For the dual case, we take the Aubert-Zelevinsky dual of the given non-unitarizable representation $\pi$ in Table \ref{tab: nonunitsub1}, which one can verify to be a non Arthur type, hence non-unitarizable subquotient of the dual one-parameter family. This proves the statement for all cases except for when $N^{\circ} = 21$. 

   When $N^{\circ} = 21$, all irreducible subquotients are unitarizable at both the first and second reducibility point. In this case, when $\frac{1}{2} < x < \frac{3}{2}$, the one-parameter complementary series falls within a non-unitary connected component of the two-parameter family   
   \[\rho\lvert \cdot \rvert^{x} \times \rho\lvert \cdot \rvert^{y} \rtimes L(\Delta_\rho[-\frac{1}{2}, -\frac{1}{2}], T_{I,1}^{\frac{1}{2}}(\pi_{sc}))\]
   with $y = \frac{1}{2}$. This is proved in Proposition \ref{2dcompseries5}. This proves the claim for $N^{\circ} = 21$. 
   
   The fact that the list is exhaustive follows from \cite{Tad23}. It suffices now to show that the unitarizable representations we provide  in Table \ref{tab: nonunitsub1} are indeed non-unitarizable subquotients of the fully induced representation in Table \ref{tab: redpnt}. By exploiting Aubert duality and Remark \ref{rmknonunitsub}, it suffices to show this when $N^{\circ} = 30$. 

   When $N^{\circ} = 30$, it suffices to consider the case $\alpha \geq \frac{1}{2}$. In this case, we need to find a non-unitarizable subquotient of 
   \[\Pi_{\alpha} = L(\Delta_\rho[0,-1], \Delta_\rho[1,0])\lvert \cdot \rvert^{\alpha} \rtimes \pi_{sc} = L(\Delta_\rho[\alpha, |\alpha-1|], \Delta_\rho[\alpha+1, \alpha]) \rtimes \pi_{sc}.\]
   First, note that 
   \[D_{\rho\lvert \cdot \rvert^{\alpha}}^{(k_{\alpha})} D_{\rho\lvert \cdot \rvert^{|\alpha-1|}}^{(k_{(|\alpha-1|})}(\Pi_{\alpha}) =\Delta_\rho[\alpha+1, \alpha] \rtimes \pi_{sc},\]
   where $k_{i}$ is the highest order derivative of $\Pi_{\alpha}$ with respect to $\rho\lvert \cdot \rvert^{i}$. In particular, we have $k_{\alpha} = k_{(|\alpha-1|)} = 1$. The above fully induced representation has a subquotient 
   \[L(\Delta_\rho[-\alpha+1, -\alpha+1]; T_{I,1}^{\alpha}(\pi_{sc})),\]
   which is of Arthur type and unitarizable. Further more, we have
   \begin{flalign*}
       S_{\rho\lvert \cdot \rvert^{\alpha}}^{(k_{\alpha})} S_{\rho\lvert \cdot \rvert^{|\alpha-1|}}^{(k_{(|\alpha-1|})}(L(\Delta_\rho[-\alpha+1, -\alpha+1]; T_{I,1}^{\alpha}(\pi_{sc}))) &= L(\Delta_\rho[-\alpha, -\alpha-1], \Delta_\rho[-|\alpha-1|, -\alpha]; \pi_{sc}).\\
   \end{flalign*}
   This implies that $L(\Delta_\rho[-\alpha, -\alpha-1], \Delta_\rho[-|\alpha-1|, -\alpha]; \pi_{sc})$ is a subquotient of $\Pi_{\alpha}$, as desired. It is non-unitarizable since it is not of Arthur type. 
\end{proof}

\begin{center}
{\footnotesize
\begin{longtable}{|@{}c@{}|l@{}|l@{}|l@{}|@{}c@{}|}
\caption{Reducibility data for one parameter families induced from unitarizable representations} \label{tab: redpnt}\\
\hline
\text{\small{\rm N}}$^\circ$ & $\pi_x$ & Reducibility points & Cases \\
\hline
\endfirsthead

\hline
\text{\small{\rm N}}$^\circ$ & $\pi_x$ & Reducibility points & Cases \\
\hline
\endhead

\hline
\endfoot

\hline
\endlastfoot

\hline
1. & $\rho \lvert \cdot \rvert^{x} \rtimes T_{I,1}^{\alpha+2}(T_{I,1}^{\alpha+1}(T_{I,1}^{\alpha}(\pi_{sc})))$ & $\alpha-1, \alpha+3$ & $\alpha \geq \tfrac{1}{2}$ \\
1'. & $\rho \lvert \cdot \rvert^{x} \rtimes L(\Delta_\rho[-\alpha-2,-\alpha-2], \Delta_\rho[-\alpha-1, -\alpha-1], \Delta_\rho[-\alpha, -\alpha]; \pi_{sc})$ & $\alpha-1, \alpha+3$ & $\alpha \geq \tfrac{1}{2}$ \\
\hline
2. & $\rho \lvert \cdot \rvert^{x} \rtimes L(\Delta_\rho[-\alpha-1, -\alpha-1], \Delta_\rho[-\alpha+1, -\alpha]; \pi_{sc})$ & $|\alpha-2|, \alpha, \alpha+2$ & $\alpha > 1$ \\
2'. & $\rho \lvert \cdot \rvert^{x} \rtimes T_{I,1}^{\alpha-1}(T_{I,1}^{\alpha+1}(T_{I,1}^{\alpha}(\pi_{sc})))$ & $|\alpha-2|, \alpha, \alpha+2$ & $\alpha > 1$ \\
3. & $\rho \lvert \cdot \rvert^{x} \rtimes L(\Delta_\rho[-\alpha+1, -\alpha+1]; T_{I,1}^{\alpha+1}(T_{I,1}^{\alpha}(\pi_{sc})))$ & $|\alpha-2|, \alpha, \alpha+2$ & $\alpha > 1$ \\
3'. & $\rho \lvert \cdot \rvert^{x} \rtimes L(\Delta_\rho[-\alpha-1, -\alpha-1], \Delta_\rho[-\alpha, -\alpha], \Delta_\rho[-\alpha+1, -\alpha+1]; \pi_{sc})$ & $|\alpha-2|, \alpha, \alpha+2$ & $\alpha > 1$ \\
\hline
4. & $\rho \lvert \cdot \rvert^{x} \rtimes T_{I,1}^{\alpha-2}(T_{I,1}^{\alpha-1}(T_{I,1}^{\alpha}(\pi_{sc})))$ & $|\alpha-3|, \alpha+1$ & $\alpha \geq 2$ \\
4'. & $\rho \lvert \cdot \rvert^{x} \rtimes L(\Delta_\rho[-\alpha+2, -\alpha]; \pi_{sc})$ & $|\alpha-3|, \alpha+1$ & $\alpha \geq 2$ \\
5. & $\rho \lvert \cdot \rvert^{x} \rtimes L(\Delta_\rho[-\alpha+2, -\alpha+2]; T_{I,1}^{\alpha-1}(T_{I,1}^{\alpha}(\pi_{sc})))$ & $|\alpha-3|, \alpha-1, \alpha+1$ & $\alpha \geq 2$ \\
5'. & $\rho \lvert \cdot \rvert^{x} \rtimes L(\Delta_\rho[-\alpha+1, -\alpha], \Delta_\rho[-\alpha+2, -\alpha+2]; \pi_{sc})$ & $|\alpha-3|, \alpha-1, \alpha+1$ & $\alpha \geq 2$ \\
6. & $\rho \lvert \cdot \rvert^{x} \rtimes L(\Delta_\rho[-\alpha +1, -\alpha+1], \Delta_\rho[-\alpha+2, -\alpha+2]; T_{I,1}^{\alpha}(\pi_{sc}))$ & $|\alpha-3|,\alpha-1, \alpha, \alpha+1$ & $\alpha \geq 2$ \\
6'. & $\rho \lvert \cdot \rvert^{x} \rtimes L(\Delta_\rho[-\alpha, -\alpha]; \Delta_\rho[-\alpha+2, -\alpha+1]; \pi_{sc})$ & $|\alpha-3|,\alpha-1, \alpha, \alpha+1$ & $\alpha \geq 2$ \\
7. & $\rho \lvert \cdot \rvert^{x} \rtimes L(\Delta_\rho[-\alpha+2, -\alpha+1]; T_{I,1}^{\alpha}(\pi_{sc}))$ & $|\alpha-3|, \alpha, \alpha+1$ & $\alpha \geq 2$ \\
7'. & $\rho \lvert \cdot \rvert^{x} \rtimes L(\Delta_\rho[-\alpha, -\alpha], \Delta_\rho[-\alpha+1, -\alpha+1], \Delta_\rho[-\alpha+2, -\alpha+2]; \pi_{sc})$ & $|\alpha-3|, \alpha, \alpha+1$ & $\alpha \geq 2$ \\
\hline
8. & $\rho \lvert \cdot \rvert^{x} \rtimes L(\Delta_\rho[-2,-2], \Delta_\rho[-1,-1]; T_{IV,3}(\pi_{sc}))$ & $0,1,3$ & $\alpha = 1$ \\
8'. & $\rho \lvert \cdot \rvert^{x} \rtimes T_{V,2}^{+}(T_{I,1}^{2}(T_{I,1}^{1}(\pi_{sc})))$ & $0,1,3$ & $\alpha = 1$ \\
9. & $\rho \lvert \cdot \rvert^{x} \rtimes L(\Delta_\rho[-2,-2], \Delta_\rho[0,-1]; \pi_{sc})$ & $0,1,3$ & $\alpha = 1$ \\
9'. & $\rho \lvert \cdot \rvert^{x} \rtimes T_{V,2}^{-}(T_{I,1}^{2}(T_{I,1}^{1}(\pi_{sc})))$ & $0,1,3$ & $\alpha = 1$ \\
10. & $\rho \lvert \cdot \rvert^{x} \rtimes L(\Delta_\rho[-1,-1], \Delta_\rho[0,-1]; \pi_{sc})$ & $0,2$ & $\alpha = 1$ \\
10'. & $\rho \lvert \cdot \rvert^{x} \rtimes  L(\Delta_\rho[0,-1]; T_{I,1}^{1}(\pi_{sc}))$ & $0,2$ & $\alpha = 1$ \\
11. & $\rho \lvert \cdot \rvert^{x} \rtimes L(\Delta_\rho[-1,-1], \Delta_\rho[-1,-1]; T_{IV,3}(\pi_{sc}))$ & $0,1,2$ & $\alpha = 1$ \\
11'. & $\rho \lvert \cdot \rvert^{x} \rtimes T_{I,2}^{1}(T_{IV,3}(\pi_{sc}))$ & $0,1,2$ & $\alpha = 1$ \\
12. & $\rho \lvert \cdot \rvert^{x} \rtimes L(\Delta_\rho[-1,-1]; T_{V,2}^{-}(T_{I,1}^{1}(\pi_{sc})))$ & $0,2$ & $\alpha = 1$ \\
12'. & $\rho \lvert \cdot \rvert^{x} \rtimes T_{III,2}^{1}(T_{IV,3}(\pi_{sc}))$ & $0,2$ & $\alpha = 1$ \\
13. & $\rho \lvert \cdot \rvert^{x} \rtimes L(\Delta_\rho[-1,-1]; T_{V,2}^{+}(T_{I,1}^{1}(\pi_{sc})))$ & $0,1,2$ & $\alpha = 1$ \\
\hline
14. & $\rho \lvert \cdot \rvert^{x} \rtimes L(\Delta_\rho[-1,-1]; T_{IV,5}(\pi_{sc}))$ & $0,1,2$ & $\alpha = 1$ \\
14'. & $\rho \lvert \cdot \rvert^{x} \rtimes T_{V,4}^{+}(T_{I,1}^{1}(\pi_{sc}))$ & $0,1,2$ & $\alpha = 1$ \\
15. & $\rho \lvert \cdot \rvert^{x} \rtimes L(\Delta_\rho[0,-1]; T_{IV,3}(\pi_{sc}))$ & $0,1,2$ & $\alpha = 1$ \\
15'. & $\rho \lvert \cdot \rvert^{x} \rtimes T_{V,4}^{-}(T_{I,1}^{1}(\pi_{sc}))$ & $0,1,2$ & $\alpha = 1$ \\
\hline
16. & $\rho \lvert \cdot \rvert^{x} \rtimes L(\Delta_\rho[-\tfrac{3}{2}, -\tfrac{3}{2}],\Delta_\rho[-\tfrac{1}{2}, -\tfrac{1}{2}], \Delta_\rho[-\tfrac{1}{2},-\tfrac{1}{2}]; \pi_{sc})$ & $\tfrac{1}{2}, \tfrac{3}{2}, \tfrac{5}{2}$ & $\alpha = \tfrac{1}{2}$ \\
16'. & $\rho \lvert \cdot \rvert^{x} \rtimes T_{I,1}^{\tfrac{3}{2}}(T_{II,3}^{\tfrac{1}{2}}(\pi_{sc}))$ & $\tfrac{1}{2}, \tfrac{3}{2}, \tfrac{5}{2}$ & $\alpha = \tfrac{1}{2}$ \\
17. & $\rho \lvert \cdot \rvert^{x} \rtimes L(\Delta_\rho[\tfrac{1}{2}, -\tfrac{3}{2}]; \pi_{sc})$ & $\tfrac{1}{2}$ & $\alpha = \tfrac{1}{2}$ \\
17'. & $\rho \lvert \cdot \rvert^{x} \rtimes L(\Delta_\rho[-\tfrac{1}{2}, -\tfrac{3}{2}]; T_{I,1}^{\tfrac{1}{2}}(\pi_{sc}))$ & $\tfrac{1}{2}$ & $\alpha = \tfrac{1}{2}$ \\
18. & $\rho \lvert \cdot \rvert^{x} \rtimes L(\Delta_\rho[-\tfrac{3}{2}, -\tfrac{3}{2}], \Delta_\rho[-\tfrac{1}{2}, -\tfrac{1}{2}]; T_{I,1}^{\tfrac{1}{2}}(\pi_{sc}))$ & $\tfrac{3}{2}, \tfrac{5}{2}$ & $\alpha = \tfrac{1}{2}$ \\
18'. & $\rho \lvert \cdot \rvert^{x} \rtimes T_{I,1}^{\tfrac{3}{2}}(T_{I,2}^{\tfrac{1}{2}}(\pi_{sc}))$ & $\tfrac{3}{2}, \tfrac{5}{2}$ & $\alpha = \tfrac{1}{2}$ \\
19. & $\rho \lvert \cdot \rvert^{x} \rtimes L(\Delta_\rho[-\tfrac{1}{2}, -\tfrac{1}{2}]; T_{I,1}^{\tfrac{3}{2}}(T_{I,1}^{\tfrac{1}{2}}(\pi_{sc})))$ & $\tfrac{1}{2}, \tfrac{3}{2}, \tfrac{5}{2}$ & $\alpha = \tfrac{1}{2}$ \\
19'. & $\rho \lvert \cdot \rvert^{x} \rtimes L(\Delta_\rho[-\tfrac{3}{2}, -\tfrac{3}{2}]; T_{I,2}^{\tfrac{1}{2}}(\pi_{sc}))$ & $\tfrac{1}{2}, \tfrac{3}{2}, \tfrac{5}{2}$ & $\alpha = \tfrac{1}{2}$ \\
\hline
20. & $\rho \lvert \cdot \rvert^{x} \rtimes T_{I,3}^{\tfrac{1}{2}}(\pi_{sc})$ & $\tfrac{1}{2}, \tfrac{3}{2}$ & $\alpha = \tfrac{1}{2}$ \\
20'. & $\rho \lvert \cdot \rvert^{x} \rtimes L(\Delta_\rho[-\tfrac{1}{2}, -\tfrac{1}{2}], \Delta_\rho[-\tfrac{1}{2}, -\tfrac{1}{2}], \Delta_\rho[-\tfrac{1}{2}, -\tfrac{1}{2}]; \pi_{sc})$ & $\tfrac{1}{2}, \tfrac{3}{2}$ & $\alpha = \tfrac{1}{2}$ \\
21. & $\rho \lvert \cdot \rvert^{x} \rtimes L(\Delta_\rho[-\tfrac{1}{2}, -\tfrac{1}{2}], \Delta_\rho[-\tfrac{1}{2}, -\tfrac{1}{2}]; T_{I,1}^{\tfrac{1}{2}}(\pi_{sc}))$ & $\tfrac{1}{2}, \tfrac{3}{2}$ & $\alpha = \tfrac{1}{2}$ \\
21'. & $\rho \lvert \cdot \rvert^{x} \rtimes L(\Delta_\rho[-\tfrac{1}{2}, -\tfrac{1}{2}]; T_{II,3}^{\tfrac{1}{2}}(\pi_{sc}))$ & $\tfrac{1}{2}, \tfrac{3}{2}$ & $\alpha = \tfrac{1}{2}$ \\
22. & $\rho \lvert \cdot \rvert^{x} \rtimes L(\Delta_\rho[-\tfrac{1}{2}, -\tfrac{1}{2}]; T_{I,2}^{\tfrac{1}{2}}(\pi_{sc}))$ & $\tfrac{3}{2}$ & $\alpha = \tfrac{1}{2}$ \\
\hline
23. & $\rho\lvert \cdot \rvert^{x} \rtimes L(\Delta_\rho[-2,-2],\Delta_\rho[-1,-1]; T_{V,2}^{+}(\pi_{sc}))$ & $0,1,3$ & $\alpha = 0$ \\
23'. & $\rho\lvert \cdot \rvert^{x} \rtimes T_{I,1}^{2}(T_{I,1}^{1}(T_{V,2}^{-}(\pi_{sc})))$ & $0,1,3$ & $\alpha = 0$ \\
24. & $\rho\lvert \cdot \rvert^{x} \rtimes L(\Delta_\rho[-2,-2],\Delta_\rho[-1,-1]; T_{V,2}^{-}(\pi_{sc}))$ & $0,1,3$ & $\alpha = 0$ \\
24'. & $\rho\lvert \cdot \rvert^{x} \rtimes T_{I,1}^{2}(T_{I,1}^{1}(T_{V,2}^{+}(\pi_{sc})))$ & $0,1,3$ & $\alpha = 0$ \\
25. & $\rho\lvert \cdot \rvert^{x} \rtimes L(\Delta_\rho[-1,-1], \Delta_\rho[-1,-1]; T_{V,2}^{\pm}(\pi_{sc}))$ & $0,2$ & $\alpha = 0$ \\
25'. & $\rho\lvert \cdot \rvert^{x} \rtimes T_{I,2}^{\pm}(T_{V,2}^{\pm}(\pi_{sc}))$ & $0,2$ & $\alpha = 0$ \\
26. & $\rho\lvert \cdot \rvert^{x} \rtimes L(\Delta_\rho[-1,-1], \Delta_\rho[-1,-1]; T_{I,1}^{1}(T_{V,2}^{\pm}(\pi_{sc})))$ & $0,2$ & $\alpha = 0$ \\
27. & $\rho\lvert \cdot \rvert^{x} \rtimes L(\Delta_\rho[-1,-1]; T_{V,4}^{\pm}(\pi_{sc}))$ & $0,1,2$ & $\alpha = 0$ \\
27'. & $\rho\lvert \cdot \rvert^{x} \rtimes T_{I,1}^{1}(T_{V,4}^{\pm}(\pi_{sc}))$ & $0,1,2$ & $\alpha = 0$ \\
28. & $\rho\lvert \cdot \rvert^{x} \rtimes L(\Delta_\rho[0,-1]; T_{V,2}^{\pm}(\pi_{sc}))$ & $0,1,2$ & $\alpha = 0$ \\
29. & $\rho\lvert \cdot \rvert^{x} \rtimes T_{V,6}^{\pm}(\pi_{sc})$ & $0,1$ & $\alpha = 0$ \\
\hline
30. & $L(\Delta_\rho[0,-1],\Delta_\rho[1,0])\lvert \cdot \rvert^{x}  \rtimes \pi_{sc}$ & $|\alpha -1|, \alpha, \alpha+1$ & all \\
31. & $ \Delta_\rho[1,-1] \lvert \cdot \rvert^{x} \rtimes T_{I,1}^{\alpha}(\pi_{sc})$ & $|\alpha-2|, \alpha, \alpha+1, \alpha+2$ & $\alpha \neq 0$ \\
31'. & $ L(\Delta_\rho[-1,-1], \Delta_\rho[0,0], \Delta_\rho[1,1]) \lvert \cdot \rvert^{x} \rtimes L(\Delta_\rho[-\alpha, -\alpha]; \pi_{sc})$ & $|\alpha-2|, \alpha, \alpha+1, \alpha+2$ & $\alpha \neq 0$ \\
32. & $ \Delta_\rho[1,-1]\lvert \cdot \rvert^{x} \rtimes L(\Delta_\rho[-\alpha, -\alpha]; \pi_{sc}) $ & $|\alpha-2|, |\alpha -1|, \alpha, \alpha+2$ & $\alpha \neq 0$ \\
32'. & $L(\Delta_\rho[-1,-1], \Delta_\rho[0,0], \Delta_\rho[1,1])\lvert \cdot \rvert^{x} \rtimes T_{I,1}^{\alpha}(\pi_{sc})$ &$|\alpha-2|, |\alpha -1|, \alpha, \alpha+2$ & $\alpha \neq 0$ \\
33. & $ \Delta_\rho[1,-1] \lvert \cdot \rvert^{x}\rtimes T_{V,2}^{\pm}(\pi_{sc})$ & $0,1,2$ & $\alpha = 0$ \\
33'. & $ L(\Delta_\rho[-1,-1], \Delta_\rho[0,0], \Delta_\rho[1,1]) \lvert \cdot \rvert^{x}\rtimes T_{V,2}^{\mp}(\pi_{sc})$ & $0,1,2$ & $\alpha = 0$ \\
\hline

34. & $\Delta_\rho[\frac{1}{2}, -\frac{1}{2}] \lvert \cdot \rvert^{x}\rtimes L(\Delta_\rho[-\alpha-1, -\alpha-1], \Delta_\rho[-\alpha, -\alpha]; \pi_{sc})$ & $|\alpha - \frac{3}{2}|, \alpha - \frac{1}{2}, \alpha + \frac{5}{2}$ & $\alpha \geq \frac{1}{2}$\\
34'. & $ L(\Delta_\rho[-\frac{1}{2}, -\frac{1}{2}], \Delta_\rho[\frac{1}{2}, \frac{1}{2}])\lvert \cdot \rvert^{x} \rtimes T_{I,1}^{\alpha+1}(T_{I,1}^{\alpha}(\pi_{sc}))$ & $|\alpha - \frac{3}{2}|, \alpha - \frac{1}{2}, \alpha + \frac{5}{2}$  & $\alpha \geq \frac{1}{2}$ \\
35. & $ \Delta_\rho[\frac{1}{2}, -\frac{1}{2}] \lvert \cdot \rvert^{x} \rtimes T_{I,1}^{\alpha+1}(T_{I,1}^{\alpha}(\pi_{sc}))$ & $|\alpha - \frac{3}{2}|, \alpha - \frac{1}{2}, \alpha + \frac{3}{2}, \alpha + \frac{5}{2}$ & $\alpha \geq \frac{1}{2}$\\
35'. & $  L(\Delta_\rho[-\frac{1}{2}, -\frac{1}{2}], \Delta_\rho[\frac{1}{2}, \frac{1}{2}])\lvert \cdot \rvert^{x} \rtimes L(\Delta_\rho[-\alpha-1, -\alpha-1], \Delta_\rho[-\alpha, -\alpha]; \pi_{sc})$ &$|\alpha - \frac{3}{2}|, \alpha - \frac{1}{2}, \alpha + \frac{3}{2}, \alpha + \frac{5}{2}$ & $\alpha \geq \frac{1}{2}$ \\

36. & $ \Delta_\rho[\frac{1}{2}, -\frac{1}{2}] \lvert \cdot \rvert^{x}\rtimes L(\Delta_\rho[-\alpha+1, -\alpha]; \pi_{sc})$ & $|\alpha - \frac{5}{2}|, \alpha - \frac{3}{2}, \alpha + \frac{1}{2}, \alpha + \frac{3}{2}$& $\alpha \geq \frac{3}{2}$\\
36'. & $  L(\Delta_\rho[-\frac{1}{2}, -\frac{1}{2}], \Delta_\rho[\frac{1}{2}, \frac{1}{2}]) \lvert \cdot \rvert^{x}\rtimes T_{I,1}^{\alpha-1}(T_{I,1}^{\alpha}(\pi_{sc}))$ &$|\alpha - \frac{5}{2}|, \alpha - \frac{3}{2}, \alpha + \frac{1}{2}, \alpha + \frac{3}{2}$& $\alpha \geq \frac{3}{2}$ \\

37. & $ \Delta_\rho[\frac{1}{2}, -\frac{1}{2}]\lvert \cdot \rvert^{x} \rtimes  T_{I,1}^{\alpha-1}(T_{I,1}^{\alpha}(\pi_{sc}))$ & $|\alpha - \frac{5}{2}|, \alpha + \frac{1}{2}, \alpha + \frac{3}{2}$& $\alpha \geq \frac{3}{2}$\\
37'. & $  L(\Delta_\rho[-\frac{1}{2}, -\frac{1}{2}], \Delta_\rho[\frac{1}{2}, \frac{1}{2}]) \lvert \cdot \rvert^{x}\rtimes L(\Delta_\rho[-\alpha+1, -\alpha]; \pi_{sc})$ & $|\alpha - \frac{5}{2}|, \alpha + \frac{1}{2}, \alpha + \frac{3}{2}$ & $\alpha \geq \frac{3}{2}$ \\

38. & $ \Delta_\rho[\frac{1}{2}, -\frac{1}{2}]\lvert \cdot \rvert^{x} \rtimes L(\Delta_\rho[-\alpha+1, -\alpha+1]; T_{I,1}^{\alpha}(\pi_{sc}))$ & $|\alpha - \frac{5}{2}|, \alpha - \frac{3}{2}, \alpha + \frac{1}{2}, \alpha + \frac{3}{2}$ & $\alpha \geq \frac{3}{2}$\\
38'. & $  L(\Delta_\rho[-\frac{1}{2}, -\frac{1}{2}], \Delta_\rho[\frac{1}{2}, \frac{1}{2}]) \lvert \cdot \rvert^{x} \rtimes L(\Delta_\rho[-\alpha+1, -\alpha+1], \Delta_\rho[-\alpha, -\alpha]; \pi_{sc})$ &$|\alpha - \frac{5}{2}|, \alpha - \frac{3}{2}, \alpha + \frac{1}{2}, \alpha + \frac{3}{2}$ & $\alpha \geq \frac{3}{2}$ \\
39. & $ \Delta_\rho[\frac{1}{2}, -\frac{1}{2}] \lvert \cdot \rvert^{x}\rtimes  L(\Delta_\rho[-\alpha+1, -\alpha+1], \Delta_\rho[-\alpha, -\alpha]; \pi_{sc})$ & $|\alpha - \frac{5}{2}|, \alpha - \frac{1}{2}, \alpha + \frac{1}{2}, \alpha + \frac{3}{2}$ & $\alpha \geq \frac{3}{2}$\\
39'. & $ L(\Delta_\rho[-\frac{1}{2}, -\frac{1}{2}], \Delta_\rho[\frac{1}{2}, \frac{1}{2}]) \lvert \cdot \rvert^{x}\rtimes L(\Delta_\rho[-\alpha+1, -\alpha+1]; T_{I,1}^{\alpha}(\pi_{sc}))$ & $|\alpha - \frac{5}{2}|, \alpha - \frac{1}{2}, \alpha + \frac{1}{2}, \alpha + \frac{3}{2}$& $\alpha \geq \frac{3}{2}$ \\
\hline
40. & $ \Delta_\rho[\frac{1}{2}, -\frac{1}{2}] \lvert \cdot \rvert^{x}\rtimes L(\Delta_\rho[-1,-1]; T_{IV,3}(\pi_{sc}))$ &$\frac{1}{2}, \frac{3}{2}, \frac{5}{2}$ & $\alpha = 1$\\
40'. & $L(\Delta_\rho[-\frac{1}{2}, -\frac{1}{2}], \Delta_\rho[\frac{1}{2}, \frac{1}{2}])\lvert \cdot \rvert^{x} \rtimes T_{V,2}^{+}(T_{I,1}^{1}(\pi_{sc}))$ &$\frac{1}{2}, \frac{3}{2}, \frac{5}{2}$ & $\alpha = 1$ \\
41. & $ \Delta_\rho[\frac{1}{2}, -\frac{1}{2}] \lvert \cdot \rvert^{x}\rtimes T_{V,2}^{+}(T_{I,1}^{1}(\pi_{sc}))$ & $\frac{1}{2}, \frac{3}{2}, \frac{5}{2}$& $\alpha = 1$\\
41'. & $  L(\Delta_\rho[-\frac{1}{2}, -\frac{1}{2}], \Delta_\rho[\frac{1}{2}, \frac{1}{2}])\lvert \cdot \rvert^{x} \rtimes L(\Delta_\rho[-1,-1]; T_{IV,3}(\pi_{sc}))$ &$\frac{1}{2}, \frac{3}{2}, \frac{5}{2}$ & $\alpha = 1$ \\

42. & $ \Delta_\rho[\frac{1}{2}, -\frac{1}{2}] \lvert \cdot \rvert^{x}\rtimes L(\Delta_\rho[0,-1]; \pi_{sc})$ & $\frac{1}{2}, \frac{3}{2}, \frac{5}{2}$& $\alpha = 1$\\
42'. & $  L(\Delta_\rho[-\frac{1}{2}, -\frac{1}{2}], \Delta_\rho[\frac{1}{2}, \frac{1}{2}]) \lvert \cdot \rvert^{x} \rtimes T_{V,2}^{-}(T_{I,1}^{1}(\pi_{sc}))$ &$\frac{1}{2}, \frac{3}{2}, \frac{5}{2}$ & $\alpha = 1$ \\
43. & $ \Delta_\rho[\frac{1}{2}, -\frac{1}{2}] \lvert \cdot \rvert^{x} \rtimes T_{V,2}^{-}(T_{I,1}^{1}(\pi_{sc}))$ & $\frac{1}{2}, \frac{3}{2}, \frac{5}{2}$& $\alpha = 1$\\
43'. & $  L(\Delta_\rho[-\frac{1}{2}, -\frac{1}{2}], \Delta_\rho[\frac{1}{2}, \frac{1}{2}]) \lvert \cdot \rvert^{x} \rtimes L(\Delta_\rho[0,-1]; \pi_{sc})$ &$\frac{1}{2}, \frac{3}{2}, \frac{5}{2}$ & $\alpha = 1$ \\
\hline
44. & $ \Delta_\rho[\frac{1}{2}, -\frac{1}{2}] \lvert \cdot \rvert^{x} \rtimes L(\Delta_\rho[-\frac{1}{2}, -\frac{1}{2}], \Delta_\rho[-\frac{1}{2}, -\frac{1}{2}]; \pi_{sc})$ & $0,1,2$& $\alpha = \frac{1}{2}$\\
44'. & $  L(\Delta_\rho[-\frac{1}{2}, -\frac{1}{2}], \Delta_\rho[\frac{1}{2}, \frac{1}{2}]) \lvert \cdot \rvert^{x} \rtimes T_{I,2}^{\frac{1}{2}}(\pi_{sc})$ & $0,1,2$& $\alpha = \frac{1}{2}$ \\
45. & $ \Delta_\rho[\frac{1}{2}, -\frac{1}{2}] \lvert \cdot \rvert^{x} \rtimes T_{I,2}^{\frac{1}{2}}(\pi_{sc})$ & $1,2$& $\alpha = \frac{1}{2}$\\
45'. & $ L(\Delta_\rho[-\frac{1}{2}, -\frac{1}{2}], \Delta_\rho[\frac{1}{2}, \frac{1}{2}]) \lvert \cdot \rvert^{x} \rtimes L(\Delta_\rho[-\frac{1}{2}, -\frac{1}{2}], \Delta_\rho[-\frac{1}{2}, -\frac{1}{2}]; \pi_{sc})$ & $1,2$& $\alpha = \frac{1}{2}$ \\

46. & $ \Delta_\rho[\frac{1}{2}, -\frac{1}{2}] \lvert \cdot \rvert^{x} \rtimes L(\Delta_\rho[-\frac{1}{2}, -\frac{1}{2}]; T_{I,1}^{\frac{1}{2}}(\pi_{sc}))$ & $1,2$ & $\alpha = \frac{1}{2}$\\
46'. & $ L(\Delta_\rho[-\frac{1}{2}, -\frac{1}{2}], \Delta_\rho[\frac{1}{2}, \frac{1}{2}]) \lvert \cdot \rvert^{x} \rtimes T_{III,2}^{\frac{1}{2}}(\pi_{sc})$ & $1,2$& $\alpha = \frac{1}{2}$ \\
47. & $ \Delta_\rho[\frac{1}{2}, -\frac{1}{2}] \lvert \cdot \rvert^{x} \rtimes T_{III,2}^{\frac{1}{2}}(\pi_{sc})$ &$1,2$ & $\alpha = \frac{1}{2}$\\
47'. & $  L(\Delta_\rho[-\frac{1}{2}, -\frac{1}{2}], \Delta_\rho[\frac{1}{2}, \frac{1}{2}]) \lvert \cdot \rvert^{x} \rtimes T_{I,1}^{\frac{1}{2}}(\pi_{sc}))$ &$1,2$  & $\alpha = \frac{1}{2}$ \\
\hline

48. & $ \Delta_\rho[\frac{1}{2}, -\frac{1}{2}] \lvert \cdot \rvert^{x}\rtimes L(\Delta_\rho[0,-1]; \pi_{sc})$ & $\frac{1}{2}, \frac{3}{2}, \frac{5}{2}$& $\alpha = 0$\\
48'. & $  L(\Delta_\rho[-\frac{1}{2}, -\frac{1}{2}], \Delta_\rho[\frac{1}{2}, \frac{1}{2}]) \lvert \cdot \rvert^{x} \rtimes L(\Delta_\rho[0,-1]; \pi_{sc})$ & $\frac{1}{2}, \frac{3}{2}, \frac{5}{2}$& $\alpha = 0$ \\

49. & $ \Delta_\rho[\frac{1}{2}, -\frac{1}{2}] \lvert \cdot \rvert^{x} \rtimes L(\Delta_\rho[-1,-1]; T_{V,2}^{\pm}(\pi_{sc}))$ &$\frac{1}{2}, \frac{3}{2}, \frac{5}{2}$ & $\alpha = 0$\\
49'. & $ L(\Delta_\rho[-\frac{1}{2}, -\frac{1}{2}], \Delta_\rho[\frac{1}{2}, \frac{1}{2}])\lvert \cdot \rvert^{x} \rtimes T_{I,1}^{1}(T_{V,2}^{\mp}(\pi_{sc}))$ &$\frac{1}{2}, \frac{3}{2}, \frac{5}{2}$ & $\alpha = 0$ \\
50. & $ \Delta_\rho[\frac{1}{2}, -\frac{1}{2}] \lvert \cdot \rvert^{x}\rtimes T_{I,1}^{1}(T_{V,2}^{\mp}(\pi_{sc}))$ & $\frac{1}{2}, \frac{3}{2}, \frac{5}{2}$& $\alpha = 0$\\
50'. & $  L(\Delta_\rho[-\frac{1}{2}, -\frac{1}{2}], \Delta_\rho[\frac{1}{2}, \frac{1}{2}]) \lvert \cdot \rvert^{x} \rtimes L(\Delta_\rho[-1,-1]; T_{V,2}^{\pm}(\pi_{sc}))$ & $\frac{1}{2}, \frac{3}{2}, \frac{5}{2}$ & $\alpha = 0$ \\
% Continue the rest of the table here, same format
% Use \\ at the end of each row, and wrap math in $...$
% You can resume pasting here...
\end{longtable}
}
\end{center}

\begin{remark}\label{rmknonunitsub}
 In general, when computing subquotients of an induced representation of the form 
\[u_{\rho}(a,b) \lvert \cdot \rvert^{x} \rtimes \pi_{A},\]
where $\pi_{A} \in \Pi_{A,gp}$ and $x$ is a reducibility point of the parabolic induction, the difficulty arises when there are more than $1$ Steinberg segments in the $L$-data of $u_{\rho}(a,b)$. When there is exactly $1$ Steinberg segment,  the induced representation is of the form 

\begin{flalign*}
    \Delta_\rho[x_1, y_1] \lvert \cdot \rvert^{x} \rtimes \pi_{u} &= \Delta_\rho[x_1+x, y_1+x]  \rtimes \pi_{u} \\
 &\leq  \rho\lvert \cdot \rvert^{x_1 +x} \rtimes (\ldots \rtimes (\rho\lvert \cdot \rvert^{y_1+x} \rtimes \pi_{u})).
\end{flalign*}

All subquotients of the induced representation above have a fixed supercuspidal support. By induction, one can compute all subquotients of the above parabolic induction if one can compute all subquotients of $\rho\lvert \cdot \rvert^{x'} \rtimes \pi_{u}$ for a given point $x'$. This can be done using the process illustrated in Sections \ref{classtempcorank3} to \ref{classtempcorank4}, namely, identifying those representations with the given supercuspidal support and the corresponding derivatives. After that, by Theorem \ref{unitiffArthur}, one can easily verify whether a given subquotient is unitarizable by verifying whether it is of Arthur type. 

When there are more than $1$ Steinberg segments which are linked, each case must be treated individually, due to the more complicated derivatives. 
\end{remark}

\setlength{\tabcolsep}{4pt} % Adjust column spacing if needed
\renewcommand{\arraystretch}{1.2} % Adjust row spacing

\setlength{\tabcolsep}{4pt} % Adjust column spacing
\renewcommand{\arraystretch}{1.2} % Adjust row spacing

\begin{center}
\begin{longtable}{|@{}c@{}|l@{}|l@{}|l@{}|@{}c@{}|}
\caption{Non-unitarizable irreducible subquotients at the reducibility point $x_0$} \label{tab: nonunitsub1} \\
\hline
\text{\small{\rm N}}$^\circ$ & $\pi$ & $x_0$ & Cases \\
\hline
\endfirsthead

\hline
\text{\small{\rm N}}$^\circ$ & $\pi$ & $x_0$ & Cases \\
\hline
\endhead

\hline
\endfoot

\hline
\endlastfoot

\hline
1. & $L(\Delta_\rho[-\alpha-3, -\alpha-3], T_{I,1}^{\alpha+2}(T_{I,1}^{\alpha+1}(T_{I,1}^{\alpha}(\pi_{sc}))))$ & $\alpha+3$ & $\alpha \geq \tfrac{1}{2}$ \\
\hline
2. & $L(\Delta_\rho[-\alpha-1, -\alpha-1], \Delta_\rho[-\alpha, -\alpha], \Delta_\rho[-\alpha+1, -\alpha]; \pi_{sc})$ & $\alpha$ & $\alpha > 1$ \\
3. & $L(\Delta_\rho[-\alpha+1, -\alpha]; T_{I,1}^{\alpha+1}(T_{I,1}^{\alpha}(\pi_{sc})))$ & $\alpha$ & $\alpha > 1$ \\
\hline
4. & $L(\Delta_\rho[-\alpha-1, -\alpha-1]; T_{I,1}^{\alpha-2}(T_{I,1}^{\alpha-1}(T_{I,1}^{\alpha}(\pi_{sc}))))$ & $\alpha+1$ & $\alpha \geq 2$ \\
5. & $L(\Delta_\rho[-\alpha-2, -\alpha-1]; T_{I,1}^{\alpha-1}(T_{I,1}^{\alpha}(\pi_{sc})))$ & $\alpha-1$ & $\alpha \geq 2$ \\
6. & $L(\Delta_\rho[-\alpha+1, -\alpha+1], \Delta_\rho[-\alpha+1, -\alpha+1], \Delta_\rho[-\alpha+2, -\alpha+2]; T_{I,1}^{\alpha}(\pi_{sc}))$ & $\alpha-1$ & $\alpha \geq 2$ \\
6. & $L(\Delta_\rho[-\alpha-1, -\alpha-1], \Delta_\rho[-\alpha+1, -\alpha+1], \Delta_\rho[-\alpha+2, -\alpha+2]; T_{I,1}^{\alpha}(\pi_{sc}))$ & $\alpha+1$ & $\alpha \geq 2$ \\
7. & $L(\Delta_\rho[-\alpha+2, -\alpha]; T_{I,1}^{\alpha}(\pi_{sc}))$ & $\alpha$ & $\alpha \geq 2$ \\
\hline
8. & $L(\Delta_\rho[-2,-2]; T_{V,4}^{\pm}(T_{I,1}^{1}(\pi_{sc}))$ & $0$ &$\alpha = 1$ \\
8. & $L(\Delta_\rho[-2,-2], \Delta_\rho[-2,-2], \Delta_\rho[-1,-1]; T_{IV,3}(\pi_{sc})$ & $3$ &$\alpha = 1$ \\
9. & $L(\Delta_\rho[-2,-2]; T_{V,4}^{\pm}(T_{I,1}^{1}(\pi_{sc}))$ & $0$ &$\alpha = 1$ \\
9. & $L(\Delta_\rho[-2,-3], \Delta_\rho[0,-1];\pi_{sc})$ & $3$ & $\alpha = 1$ \\
10. & $L(\Delta_\rho[-1,-1],\Delta_\rho[0,-2]; \pi_{sc})$ & $2$ & $\alpha = 1$ \\
11. & $L(\Delta_\rho[-1,-2],\Delta_\rho[-1,-1]; T_{IV,3}(\pi_{sc}))$& $2$ & $\alpha = 1$ \\
13. & $L(\Delta_\rho[-1,-2]; T_{I,1}^{1}(T_{IV,3}(\pi_{sc}))$  & $2$ &$\alpha = 1$\\
14. & $L(\Delta_\rho[-1,-2]; T_{I,1}^{1}(T_{IV,3}(\pi_{sc}))$  & $2$ &$\alpha = 1$\\
15. & $L(\Delta_\rho[0,-2]; T_{IV,3}(\pi_{sc}))$  & $2$ &$\alpha = 1$\\
\hline
16. & $L(\Delta_\rho[-\tfrac{3}{2}, -\tfrac{3}{2}], \Delta_\rho[-\tfrac{1}{2}, -\tfrac{3}{2}], \Delta_\rho[-\tfrac{1}{2}, -\tfrac{1}{2}]; \pi_{sc})$ & $\tfrac{3}{2}$ & $\alpha = \tfrac{1}{2}$ \\
18. & $L(\Delta_\rho[-\tfrac{3}{2}, -\tfrac{5}{2}], \Delta_\rho[-\tfrac{1}{2}, -\tfrac{1}{2}]; T_{I,1}^{\tfrac{1}{2}}(\pi_{sc}))$ & $\tfrac{5}{2}$ & $\alpha = \tfrac{1}{2}$ \\
19. & $L(\Delta_\rho[-\tfrac{1}{2}, -\tfrac{3}{2}]; T_{I,1}^{\tfrac{3}{2}}(T_{I,1}^{\tfrac{1}{2}}(\pi_{sc})))$ & $\tfrac{3}{2}$ & $\alpha = \tfrac{1}{2}$ \\
20. & $L(\Delta_\rho[-\tfrac{3}{2}, -\tfrac{3}{2}]; T_{I,3}^{\tfrac{1}{2}}(\pi_{sc}))$ & $\tfrac{3}{2}$ & $\alpha = \tfrac{1}{2}$ \\
\hline
23. & $L(\Delta_\rho[-2,-3],\Delta_\rho[-1,-1]; T_{V,2}^{+}(\pi_{sc}))$ & $3$ &$\alpha = 0$ \\
24. & $L(\Delta_\rho[-2,-3],\Delta_\rho[-1,-1]; T_{V,2}^{-}(\pi_{sc}))$ & $3$ &$\alpha = 0$ \\
27. & $L(\Delta_\rho[-1,-2]; T_{V,4}^{\pm}(\pi_{sc}))$ & $2$ &$\alpha = 0$ \\
28. & $L(\Delta_\rho[-2,-2], \Delta_\rho[0,-1]; T_{V,2}^{\pm}(\pi_{sc}))$ & $2$ &$\alpha = 0$ \\
\hline
30. & $L(\Delta_\rho[-\alpha-1, -\alpha], \Delta_\rho[-|\alpha-1|, -|\alpha-1|]; T_{I,1}^{\alpha}(\pi_{sc}))$ & $\alpha$ & $\alpha \neq 0$ \\
31. & $L(\Delta_\rho[-\alpha+1, -\alpha-1]; T_{I,1}^{\alpha}(\pi_{sc}))$ & $\alpha$ & $\alpha > 1$\\
31. & $L(\Delta_\rho[-\alpha, -\alpha-2]; T_{I,1}^{\alpha}(\pi_{sc}))$ & $\alpha+1$ & $\alpha \neq 0$\\
32. & $L(\Delta_\rho[-\alpha, -\alpha-2], \Delta_\rho[-\alpha, -\alpha]; \pi_{sc})$ & $\alpha -1$ &  $\alpha \geq \frac{3}{2}$\\
32. & $L(\Delta_\rho[-\alpha+1, -\alpha-1], \Delta_\rho[-\alpha, -\alpha]; \pi_{sc})$ & $\alpha$ & $\alpha \geq 1$\\
32. & $L(\Delta_\rho[-\frac{1}{2}, -\frac{5}{2}], \Delta_\rho[-\frac{1}{2}, -\frac{1}{2}]; \pi_{sc})$ & $\frac{3}{2}$ & $\alpha = \frac{1}{2}$\\
33. & $L(\Delta_\rho[-1,-3]; T_{V,2}^{\pm}(\pi_{sc}))$ & $2$ & $\alpha = 0$\\
34. &$L(\Delta_\rho[-\alpha-1, -\alpha-1], \Delta_\rho[-\alpha, -\alpha], \Delta_\rho[-\alpha+1, -\alpha]; \pi_{sc})$ & $\alpha - \frac{1}{2}$ & $\alpha \geq \frac{3}{2}$ \\
34. &$L(\Delta_\rho[-\alpha-2,-\alpha-3],\Delta_\rho[-\alpha-1,-\alpha-1],\Delta_\rho[-\alpha,-\alpha]; \pi_{sc})$ & $\alpha + \frac{5}{2}$ & $\alpha = \frac{1}{2}, 1$\\
35. &$L(\Delta_\rho[-\alpha+1, -\alpha]; T_{I,1}^{\alpha+1}(T_{I,1}^{\alpha}(\pi_{sc})))$ & $\alpha - \frac{1}{2}$ & $\alpha \geq 1$\\
35. &$L(\Delta_\rho[-\alpha-1, -\alpha-2]; T_{I,1}^{\alpha+1}(T_{I,1}^{\alpha}(\pi_{sc})))$ & $\alpha + \frac{3}{2}$ & $\alpha  \geq \frac{1}{2}$\\
36. &$L(\Delta_\rho[-\alpha -2, -\alpha-2], \Delta_\rho[-\alpha+1, -\alpha]; T_{I,1}^{\alpha}(\pi_{sc}))$ & $\alpha - \frac{3}{2}$ & $\alpha \geq 2$\\
36. &$L(\Delta_\rho[-\alpha-1, -\alpha-1], \Delta_\rho[-\alpha+1, -\alpha]; T_{I,1}^{\alpha}(\pi_{sc}))$ & $\alpha + \frac{1}{2}$ & $\alpha \geq \frac{3}{2}$\\
\hline
37. &$L(\Delta_\rho[-\alpha, -\alpha-1]; T_{I,1}^{\alpha-1}(T_{I,1}^{\alpha}(\pi_{sc})))$ & $\alpha + \frac{1}{2}$ & $\alpha \geq \frac{3}{2}$ \\
38. & $L(\Delta_\rho[-\alpha+2, -\alpha+1], \Delta_\rho-\alpha+1, -\alpha+1]; T_{I,1}^{\alpha}(\pi_{sc}))$ & $\alpha - \frac{3}{2}$ & $\alpha \geq 2$ \\
38. & $L(\Delta_\rho[-\alpha+1, -\alpha+1], \Delta_\rho[-\alpha, -\alpha-1]; T_{I,1}^{\alpha}(\pi_{sc}))$ & $\alpha + \frac{1}{2}$ & $\alpha \geq \frac{3}{2}$ \\
39. & $L(\Delta_\rho[-\alpha+1, -\alpha+1], \Delta_\rho[-\alpha+1, -\alpha], \Delta_\rho[-\alpha, -\alpha]; \pi_{sc})$ & $\alpha - \frac{1}{2}$ & $\alpha \geq \frac{3}{2}$ \\
39. & $L(\Delta_\rho[-\alpha+1, -\alpha+1], \Delta_\rho[-\alpha, -\alpha], \Delta_\rho[-\alpha, -\alpha-1]; \pi_{sc})$ & $\alpha + \frac{1}{2}$ &  $\alpha \geq \frac{3}{2}$ \\
\hline
40. & $L(\Delta_\rho[-1,-1],\Delta_\rho[-1,-2]; T_{IV,3}(\pi_{sc}))$ & $\frac{3}{2}$ & $\alpha = 1$\\
41. & $L(\Delta_\rho[-1,-2]; T_{V,2}^{+}(T_{I,1}^{1}(\pi_{sc})))$ & $\frac{3}{2}$ & $\alpha = 1$ \\
42. & $L(\Delta_\rho[0,-1],\Delta_\rho[-1,-2]; \pi_{sc})$ & $\frac{3}{2}$ & $\alpha = 1$ \\
43. & $L(\Delta_\rho[-1,-2]; T_{V,2}^{-}(T_{I,1}^{1}(\pi_{sc})))$ & $\frac{3}{2}$ & $\alpha = 1$ \\
\hline
44. & $L(\Delta_\rho[-\frac{3}{2},-\frac{5}{2}], \Delta_\rho[-\frac{1}{2}, -\frac{1}{2}], \Delta_\rho[-\frac{1}{2}, -\frac{1}{2}]; \pi_{sc})$ & $2$ & $\alpha = \frac{1}{2}$ \\
45. &$L(\Delta_\rho[-\frac{3}{2}, -\frac{5}{2}]; T_{I,2}^{\frac{1}{2}}(\pi_{sc}))$ & $2$& $\alpha = \frac{1}{2}$ \\
46. &$L(\Delta_\rho[-\frac{1}{2}, -\frac{1}{2}], \Delta_\rho[-\frac{3}{2}, -\frac{5}{2}]; T_{I,1}^{\frac{1}{2}}(\pi_{sc}))$ & $2$ & $\alpha = \frac{1}{2}$ \\
47. & $L(\Delta_\rho[-\frac{3}{2}, -\frac{5}{2}]; T_{III,2}^{\frac{1}{2}}(\pi_{sc})$ & $2$ & $\alpha = \frac{1}{2}$ \\
\hline
48. & $L(\Delta_\rho[0,-1],\Delta_\rho[-1,-2]; \pi_{sc})$ & $\frac{3}{2}$ & $\alpha = 0$\\
49. & $L(\Delta_\rho[-1,-1],\Delta_\rho[-1,-2]; T_{V,2}^{\pm}(\pi_{sc}))$ & $\frac{3}{2}$ & $\alpha = 0$ \\
50. & $L(\Delta_\rho[-1,-2]; T_{I,1}^{1}(T_{V,2}^{\mp}(\pi_{sc})))$ & $\frac{3}{2}$ & $\alpha = 0$
\end{longtable}
\end{center}

Proposition \ref{1dcompleseries} and the two tables above summarize all possible unitarizable one-parameter complementary series of the form $u_{\rho}(a,b) \lvert \cdot \rvert^{x} \rtimes \pi_A$, where $\pi_A \in \Pi_{A,gp}(G_m)$. To finish up our classification of one-parameter complementary series, we now consider all other families that are not of this form. To begin, we consider all other families induced from supercuspidal representations of smaller rank $G_m$ in Propositions \ref{1dcompscextra} to \ref{1dcompscextra,4} below. 

\begin{prop}\label{1dcompscextra}
     For all $x \in \mathbb{R}$, let \[\pi_{1,x} = L(\Delta_\rho[0,0], \Delta_\rho[3,1])\lvert \cdot \vert^x \rtimes \pi_{sc},\]
    \[\pi_{2,x} = L(\Delta_\rho[1,0], \Delta_\rho[3,2])\lvert \cdot \vert^x \rtimes \pi_{sc},\]
    \[\pi_{3,x} = L(\Delta_\rho[2,0], \Delta_\rho[3,3])\lvert \cdot \vert^x \rtimes \pi_{sc}.\]
    Then for $i = 1,2,3$, $\pi_{i,x}$ is reducible if and only if $x \in \{\pm(\alpha -3), \pm(\alpha-2), \pm(\alpha-1), \pm \alpha\}$. When $\pi_{i,x}$ is irreducible, it is unitarizable if and only if $-\alpha < x < \alpha-3$.
\end{prop}
\begin{proof}
    First let $i = 1$. The reducibility statement follows from Theorem \ref{Tadirred1}. By Proposition \ref{subquotientlist}, for $\alpha > 2$, the representation $L(\Delta_\rho[-\alpha+1, -\alpha-1], \Delta_\rho[-\alpha+2, -\alpha+2]; \pi_{sc})$ is a non-unitarizable subquotient of $\pi_{1,\alpha-2}$ and 
    \[L(\Delta_\rho[-\alpha, -\alpha-2], \Delta_\rho[-\alpha+1, -\alpha+1]; \pi_{sc})\]
    is a non-unitarizable subquotient of $\pi_{1,\alpha-1}$. 
    
    This proves that there is no complementary series when $x > \alpha -3$. By Propositions \ref{crnk4Artcritlist} and \ref{subquotientlist}, all subquotients of $\Pi_{|\alpha-3|, \alpha-2, \alpha-1, \alpha}$ are of Arthur type and hence unitarizable for $\alpha > 2$. Similarly, one can show that all subquotients of $\Pi_{x,x+1,x+2,x+3}$ are unitarizable, when $x = -\alpha, -(\alpha-1), -(\alpha-2), -(\alpha-3)$. This proves the statement for $\alpha > 2$. 

    For $\alpha = 2$, reducibility occurs at $0$ so there is no unitary complementary series for $x \geq 0$. When $x < 0$, the region $-2 < x < -1$ is unitarizable since when $x = -1$, all subquotients of $\Pi_{(-1,0,1,2)} \cong \Pi_{(0,1,1,2)}$ are unitarizable by Proposition \ref{subquotientlist}. 
    
    For $\alpha = 0,1$, reducibility occurs at $0$ so there is no complementary series beyond $0$. When $\alpha = 1$, 
    \[L(\Delta_\rho[-1,-1],\Delta_\rho[0,-2]; \pi_{sc})\]
    is a non-unitarizable subquotient of $\pi_{1,-1}$. When $\alpha = 0$, 
    \[L(\Delta_\rho[-3,-3],\Delta_\rho[0,-2]; \pi_{sc})\] is a non-unitarizable subquotient of $\pi_{1,-3}$, and 
    \[L(\Delta_\rho[-1,-1],\Delta_\rho[0,-2]; \pi_{sc})\] is a non-unitarizable subquotient of $\pi_{1,-1}$.

    Now it remains to consider the cases $\alpha = \frac{1}{2}, \frac{3}{2}, \frac{5}{2}$.  When $\alpha = \frac{1}{2}$,  the representation 
    \[L(\Delta_\rho[-\frac{1}{2}, -\frac{5}{2}], \Delta_\rho[-\frac{1}{2}, -\frac{1}{2}]; \pi_{sc})\]
    is a non-unitarizable subquotient of $\pi_{1, \frac{3}{2}}$ and $\pi_{1,-\frac{1}{2}}$.  and 
    \[L(\Delta_\rho[-\frac{3}{2}, -\frac{3}{2}], \Delta_\rho[\frac{1}{2}, -\frac{3}{2}]; \pi_{sc})\]
    is a non-unitarizable subquotient of $\pi_{1,\frac{1}{2}}$ and $\pi_{1,-\frac{3}{2}}$. 
    
    For $\alpha = \frac{3}{2}$, 
    \[L(\Delta_\rho[-\frac{1}{2}, -\frac{5}{2}], \Delta_\rho[-\frac{1}{2}, -\frac{1}{2}]; \pi_{sc})\]
    is a non-unitarizable subquotient of $\pi_{1,-\frac{1}{2}}$ and 
    \[L(\Delta_\rho[-\frac{3}{2}, -\frac{7}{2}], \Delta_\rho[-\frac{1}{2}, -\frac{1}{2}]; \pi_{sc})\]
    is a non-unitarizable subquotient of $\pi_{1,\frac{3}{2}}$.
    
    For $\alpha = \frac{5}{2}$, the representation
    \[L(\Delta_\rho[-\frac{5}{2}, -\frac{9}{2}], \Delta_\rho[-\frac{3}{2}, -\frac{3}{2}]; \pi_{sc})\]
    is a non-unitarizable subquotient of $\pi_{1,\frac{3}{2}}$; the representation
    \[L(\Delta_\rho[-\frac{3}{2}, -\frac{7}{2}], \Delta_\rho[-\frac{1}{2}, -\frac{1}{2}]; \pi_{sc})\]
    is a non-unitarizable subquotient of $\pi_{1,\frac{1}{2}}$. Same as before, all subquotients of $\Pi_{x,x+1,x+2,x+3}$ are unitarizable for $x = -\frac{5}{2}, -\frac{3}{2}, -\frac{1}{2}$, by Propositions \ref{subquotientlist} and \ref{crnk4Artcritlist}. This proves the claim for $i = 1$. The proof for $i = 2,3$ are similar.
\end{proof}

\begin{prop}\label{1dcompscextra,1}
    For $x \in \mathbb{R}$, let 
    \[\pi_{1,x} = L(\Delta_\rho[0,0],\Delta_\rho[1,1], \Delta_\rho[3,2])\lvert \cdot \rvert^{x} \rtimes \pi_{sc},\] 
    \[\pi_{2,x} = L(\Delta_\rho[0,0],\Delta_\rho[2,1], \Delta_\rho[3,3])\lvert \cdot \rvert^{x} \rtimes \pi_{sc},\] 
    \[\pi_{3,x} = L(\Delta_\rho[1,0],\Delta_\rho[2,2], \Delta_\rho[3,3]) \lvert \cdot \rvert^{x} \rtimes \pi_{sc}.\] 
    Then for $i = 1,2,3$, $\pi_{i,x}$ is reducible if and only if $x \in \{\pm(\alpha -2), \pm(\alpha-1), \pm\alpha, \pm(\alpha+1)\}$. When $\pi_{i,x}$ is irreducible, it is unitarizable if and only if $-\alpha \leq x < \alpha-3$. 
\end{prop}
\begin{proof}
    The proof is similar to that of Proposition \ref{1dcompscextra}, which we omit.
\end{proof}

\begin{prop}\label{1dcompscextra,2}
    For $x \geq 0$, let 
    \[\pi_{1,x} = L(\Delta_\rho[0,0], \Delta_\rho[1,0], \Delta_\rho[2,2])\lvert \cdot \rvert^{x} \rtimes \pi_{sc},\]
    \[\pi_{2,x} = L(\Delta_\rho[0,0], \Delta_\rho[2,0]) \lvert \cdot \rvert^{x} \rtimes \pi_{sc}.\]
    \[\pi_{3,x} = L(\Delta_\rho[0,0], \Delta_\rho[0,0], \Delta_\rho[2,1]) \lvert \cdot \rvert^{x} \rtimes \pi_{sc}.\]
    Then for $i = 1,2,3$, $\pi_{i,x}$ are irreducible when $x \in \{|\alpha-2|, |\alpha-1|, \alpha\}$. When irreducible, $\pi_{1,x}, \pi_{2,x}$ are non-unitarizable for any $x$. 
\end{prop}
\begin{proof}
    The reducibility statement follows from Theorem \ref{Tadirred1}. Let $i =1$. For all $\alpha$, the representation 
    \[L(\Delta_\rho[-\alpha-4, -\alpha-4], \Delta_\rho[-\alpha-2,-\alpha-3], \Delta_\rho[-\alpha-2, -\alpha-2]; \pi_{sc})\]
    is a non-unitarizable subquotient of $\pi_{1,|\alpha-2|}$, and 
    \[L(\Delta_\rho[-\alpha-3, -\alpha-3], \Delta_\rho[-\alpha-1,-\alpha-2], \Delta_\rho[-\alpha-1, -\alpha-1]; \pi_{sc})\] is a non-unitarizable subquotient of $\pi_{1,|\alpha-1|}$. This proves the claim for $\pi_{1,x}$. 
    Similarly, we have that for all $\alpha$, the representation 
    \[L(\Delta_\rho[-\alpha-2, -\alpha-4], \Delta_\rho[-\alpha-2, -\alpha-2]; \pi_{sc})\]
    is a non-unitarizable subquotient of $\pi_{2, |\alpha-1|}$ and 
    \[L(\Delta_\rho[-\alpha, -\alpha-2], \Delta_\rho[-\alpha, -\alpha]; \pi_{sc})\]
    is a non-unitarizable subquotient of $\pi_{2,\alpha}$. This proves the claim for $i = 1$. The proof is similar for $i = 2,3$.  
\end{proof}

\begin{prop}\label{1dcompscextra,3}
    For $x \geq 0$, let 
    \[\pi_{1,x} = L(\Delta_\rho[0,0], \Delta_\rho[2,1], \Delta_\rho[2,2])\lvert \cdot \rvert^{x} \rtimes \pi_{sc},\]
    \[\pi_{2,x} = L(\Delta_\rho[2,0],\Delta_\rho[2,2])\rtimes \pi_{sc},\]
    \[\pi_{3,x} = L(\Delta_\rho[1,0],\Delta_\rho[2,2],\Delta_\rho[2,2])\rtimes \pi_{sc}.\]
    Then for $i = 1,2,3$, $\pi_{i,x}$ are irreducible when $x \in \{|\alpha-2|, |\alpha-1|, \alpha\}$. When irreducible, $\pi_{i,x}$ are non-unitarizable for any $x$. 
\end{prop}
\begin{proof}
    The proof is similar to that of Proposition \ref{1dcompscextra,2}, which we omit. 
\end{proof}

\begin{prop}\label{1dcompscextra,4}
    For $x \geq 0$, let 
    \[\pi_{1,x} = L(\Delta_\rho[1,0], \Delta_\rho[1,1], \Delta_\rho[2,2])\lvert \cdot \rvert^{x} \rtimes \pi_{sc},\]
    \[\pi_{2,x} = L(\Delta_\rho[2,0], \Delta_\rho[1,1]) \lvert \cdot \rvert^{x} \rtimes \pi_{sc},\]
    \[\pi_{3,x} = L(\Delta_\rho[0,0], \Delta_\rho[1,1], \Delta_\rho[2,1])\lvert \cdot \rvert^{x} \rtimes \pi_{sc}.\]
    Then for $i = 1,2,3$, $\pi_{i,x}$ are irreducible if and only if $x \in \{|\alpha-2|, |\alpha-1|, \alpha\}$. When irreducible, $\pi_{i,x}$ are non-unitarizable for any $x$. 
\end{prop}
\begin{proof}
    The proof is similar to that of Proposition \ref{1dcompscextra,2}, which we omit. 
\end{proof}
In Propositions \ref{1dcompextra3+1,1} to \ref{1dcompextra3+1,3} below, we consider the unitarizability of $1$-parameter families induced from a critical unitarizable representation of a smaller rank $G_m$ that is not supercuspidal. In particular, we only need to consider the one-parameter families induced from a corank 1 critical unitarizable representation of smaller rank $G_m$, since all other cases are covered in Proposition \ref{1dcompleseries}. 

\begin{prop}\label{1dcompextra3+1,1}
    For $\alpha > 0$, let 
    \[\pi = T_{I,1}^{\alpha}(\pi_{sc}),\]
    \[\hat{\pi} = L(\Delta_\rho[-\alpha, -\alpha]; \pi_{sc}).\] 
    For $x \geq 0$, consider the family
    \[\tau_{x} = L(\Delta_\rho[x,x], \Delta_\rho[x+2, x+1]) \rtimes \pi,\]
    \[\hat{\tau}_{x} = L(\Delta_\rho[x,x], \Delta_\rho[x+2, x+1]) \rtimes \hat{\pi}.\]
    Then $\tau_{x}$ and $\hat{\tau}_{x}$ are not unitarizable for any $x \geq 0$. 
\end{prop}
\begin{proof}
    By Aubert duality it suffices to show it for $\tau_{x}$. By Proposition \ref{nontemp2C1}, $\tau_{x}$ always contain a non-unitarizable subquotient, unless $\epsilon_\rho \leq x \leq \alpha -4$. Therefore it suffices to show that for $\epsilon_\rho \leq x \leq \alpha-4$, $\tau_{x}$ is irreducible. 

    Fix $\epsilon_\rho \leq x \leq \alpha-4$, and let
    \[X_1 = \{\pm\alpha\},\] 
    \[X_2 = \{\pm x, \pm(x+1), \pm(x+2)\}.\]
    Then $X_1 \sqcup X_2$ forms a regular partition of the set 
    \[X = \{\pm\alpha, \pm x, \pm(x+1), \pm(x+2)\}.\]
    Let \[\beta_{1} = \rho\lvert \cdot \rvert^{\alpha}, \ \ \beta_{2} = L(\Delta_\rho[x,x], \Delta_\rho[x+2, x+1]),\]
    \[\gamma_{1} = T_{I,1}^{\alpha}(\pi_{sc}), \ \ \gamma_{2} = L(\Delta_\rho[-x-1,-x-2], \Delta_\rho[-x,-x]; \pi_{sc}),\]
    where $0_{F}$ denotes the trivial representation of $\GL_{0}(F)$. 
    Then by Theorem \ref{thm Jantzen}, we have that $\tau_{x}$ is irreducible if and only if both $\beta_1 \rtimes \gamma_2$ and $\beta_2 \rtimes \gamma_1$ are irreducible. For $\epsilon+1 \leq x \leq \alpha+1$, both representations above are irreducible by Theorem \ref{Tadirred1}, so we are done. 
    
\end{proof}

\begin{prop}\label{1dcompextra3+1,2}
    For $\alpha > 0$, let 
    \[\pi = T_{I,1}^{\alpha}(\pi_{sc}),\]
    \[\hat{\pi} = L(\Delta_\rho[-\alpha, -\alpha]; \pi_{sc}).\]
    For $x \geq 0$, consider the family
    \[\tau_{x} = L(\Delta_\rho[x+1, x], \Delta_\rho[x+2, x+2]) \rtimes \pi,\]
    \[\hat{\tau}_{x} = L(\Delta_\rho[x+1,x], \Delta_\rho[x+2, x+2]) \rtimes \hat{\pi}.\]
    Then $\tau_{x}$ and $\hat{\tau}_{x}$ are not unitarizable for any $x \geq 0$. 
\end{prop}
\begin{proof}
    By Aubert duality, we only need to prove it for $\tau_{x}$ a Let $\alpha > 2$, then the representation 
    \[L(\Delta_\rho[-\alpha+2, -\alpha]; T_{I,1}^{\alpha}(\pi_{sc}))\]
    is a non-unitarizable subquotient of $\tau_{\alpha-2}$ (this is proved using the same method as in the proof of Proposition \ref{1dcompleseries}).
    \[L(\Delta_\rho[-\alpha, -\alpha], \Delta_\rho[-\alpha, -\alpha], \Delta_\rho[-\alpha+1, -\alpha+1]; T_{I,1}^{\alpha}(\pi_{sc}))\]
    is a non-unitarizable subquotient of $\tau_{\alpha-1}$, and 
    \[L(\Delta_\rho[-\alpha-1, -\alpha-1], \Delta_\rho[-\alpha, -\alpha], \Delta_\rho[-\alpha, -\alpha]; T_{I,1}^{\alpha}(\pi_{sc}))\] is a non-unitarizable subquotient of $\tau_{\alpha}$. For $\alpha > 2$, this shows that there is no complementary series. 

    For $\alpha = 0,1,2$ we get reducibility at $0$ so there is no complementary series to consider. Now let $\alpha = \frac{1}{2}$, using the same method as before, one can show that the representation 
    \[L(\Delta_\rho[-\frac{3}{2}, -\frac{3}{2}]; T_{I,3}^{\frac{1}{2}}(\pi_{sc}))\]
    is a non-unitarizable subquotient of $\tau_{\frac{1}{2}}$. Let $\alpha = \frac{3}{2}$. The representation 
    \[L(\Delta_\rho[-\frac{3}{2}, -\frac{3}{2}], \Delta_\rho[-\frac{3}{2}, -\frac{3}{2}], \Delta_\rho[-\frac{1}{2}, -\frac{1}{2}]; T_{I,1}^{\frac{3}{2}}(\pi_{sc}))\]
    is a non-unitarizable subquotient of $\tau_{\frac{1}{2}}$, and 
    \[L(\Delta_\rho[-\frac{5}{2}, -\frac{5}{2}], \Delta_\rho[-\frac{3}{2}, -\frac{3}{2}], \Delta_\rho[-\frac{3}{2}, -\frac{3}{2}]; T_{I,1}^{\frac{3}{2}}(\pi_{sc}))\]
    is a non-unitarizable subquotient of $\tau_{\frac{3}{2}}$. This proves the claim. 
\end{proof}

\begin{prop}\label{1dcompextra3+1,3}
    For $\alpha = 0$, let 
    \[\pi^{\pm} = T_{V,2}^{\pm}(\pi_{sc}), \] and let 
    \[\tau_{1,x}^{\pm} = L(\Delta_\rho[x,x],\Delta_\rho[x+2,x+1]) \rtimes \pi^{\pm},\] 
    \[\tau_{2,x}^{\pm} = L(\Delta_\rho[x+1,x],\Delta_\rho[x+2,x+2]) \rtimes \pi^{\pm}.\] Then, $\tau_{1,x}^{\pm}$ and $\tau_{2,x}^{\pm}$ are not unitarizable for any $x \geq 0$.  
\end{prop}
\begin{proof}
    This follows directly from Propositions \ref{nontemp2C3} and \ref{nontemp2D3}. 
\end{proof}

\subsection{Other one-parameter families}

To conclude this section, we consider the unitarizability of one-parameter families induced from non-unitarizable representations. Since there are no non-unitarizable, critical type representations of corank $1$, we first consider one-parameter families induced from a non-unitarizable, corank 2 critical representation of smaller rank $G_m$, from Propositions \ref{1dcompnonunit1} to \ref{1dcompnonunit2} below.

\begin{prop}\label{1dcompnonunit1}
    Let 
    \[\pi = L(\Delta_\rho[-\alpha, -\alpha]; T_{I,1}^{\alpha}(\pi_{sc})),\]
    \[\hat{\pi} = L(\Delta_\rho[-\alpha, -\alpha], \Delta_\rho[-\alpha, -\alpha]; \pi_{sc}).\]
    For $x \geq 0$, consider the one-parameter family
     \[\tau_{x} = \Delta_\rho[x,x+1] \rtimes \pi,\]
    \[\hat{\tau}_{x} = \Delta_\rho[x,x+1] \rtimes \hat{\pi}.\]
    Then $\tau_{x}$ and $\hat{\tau}_{x}$ are reducible if and only if $x \in \{|\alpha-2|, |\alpha-1|, \alpha, \alpha+1\}$. When they are irreducible, $\tau_{x}$ and $\hat{\tau}_{x}$ are always non-unitarizable. 
\end{prop}
\begin{proof}
    The proof is identical to Proposition \ref{1dcompextra3+1,1}, which we omit. 
\end{proof}

\begin{prop}\label{1dcompnonunit2}
    Let 
    \[\pi = L(\Delta_\rho[-\alpha-1, -\alpha-1]; T_{I,1}^{\alpha}(\pi_{sc})),\]
    \[\hat{\pi} = L(\Delta_\rho[-\alpha, -\alpha-1]; \pi_{sc}).\]
    For $x \geq 0$, consider the one-parameter family
     \[\tau_{x} = \Delta_\rho[x,x+1] \rtimes \pi,\]
    \[\hat{\tau}_{x} = \Delta_\rho[x,x+1] \rtimes \hat{\pi}.\]
    Then $\tau_{x}$ and $\hat{\tau}_{x}$ are reducible if and only if $x \in \{|\alpha-2|, |\alpha-1|, \alpha, \alpha+1, \alpha+2\}$. When they are irreducible, $\tau_{x}$ and $\hat{\tau}_{x}$ are always non-unitarizable.
\end{prop}
\begin{proof}
    The proof is similar to that of Proposition \ref{1dcompnonunit1}, which we omit. 
\end{proof}

To conclude this subsection, we turn our attention to one-parameter families induced from a non-unitarizable, corank 3 critical type representation of smaller rank $G_m$ (listed as in \cite{Tad23}).

\begin{prop}\label{1dcompnonunitcrnk3}
    Let 
    \[\pi_{x} = \rho\lvert \cdot \rvert^{x} \rtimes \pi_{nu},\]
    where $\pi_{nu}$ is non-unitarizable corank 3 representation of critical type. Then $\pi_{x}$ is non-unitarizable when it is irreducible. 
\end{prop}
\begin{proof}
    The proof is similar to that of Proposition \ref{1dcompleseries}. The relevant data are summarized in Table \ref{tab:redpnt2} and Table \ref{tab: nonunitsub2} below. 
\end{proof}

\setlength{\tabcolsep}{4pt} % Adjust column spacing if needed
\renewcommand{\arraystretch}{1.2} % Adjust row spacing

\begin{center}
{\footnotesize
\begin{longtable}{|@{}c@{}|l@{}|l@{}|l@{}|@{}c@{}|}
\caption{Reducibility data for one parameter families $\pi_x$ induced from non-unitarizable representations} \label{tab:redpnt2}\\
\hline
\text{\small{\rm N}}$^\circ$ & $\pi_x$ & Reducibility points & Cases \\
\hline
\endfirsthead

\hline
\text{\small{\rm N}}$^\circ$ & $\Pi_x$ & Reducibility points & Cases \\
\hline
\endhead

\hline
\endfoot

\hline
\endlastfoot
1. & $\rho\lvert \cdot \rvert^{x} \rtimes L(\Delta_\rho[-\alpha-2, -\alpha-2], \Delta_\rho[-\alpha-1, -\alpha-1]; T_{I,1}^{\alpha}(\pi_{sc}))$ & $|\alpha-1|, \alpha, \alpha+3$ &$\alpha \geq \frac{1}{2}$ \\
1'. & $\rho\lvert \cdot \rvert^{x} \rtimes L(\Delta_\rho[-\alpha, -\alpha-2]; \pi_{sc})$& $|\alpha-1|, \alpha, \alpha+3$&$\alpha \geq \frac{1}{2}$ \\
2. & $\rho\lvert \cdot \rvert^{x} \rtimes L(\Delta_\rho[-\alpha-2, -\alpha-2], \Delta_\rho[-\alpha, -\alpha-1]; \pi_{sc})$ & $|\alpha-1|, \alpha, \alpha+1, \alpha+3$ &$\alpha \geq \frac{1}{2}$ \\
2'. & $\rho\lvert \cdot \rvert^{x} \rtimes L(\Delta_\rho[-\alpha-1, -\alpha-2], T_{I,1}^{\alpha}(\pi_{sc}))$ &$|\alpha-1|, \alpha, \alpha+1, \alpha+3$ &$\alpha \geq \frac{1}{2}$ \\
3. & $\rho\lvert \cdot \rvert^{x} \rtimes L(\Delta_\rho[-\alpha-2, -\alpha-2]; T_{I,1}^{\alpha+1}(T_{I,1}^{\alpha}(\pi_{sc}))))$ & $|\alpha-1|, \alpha+1, \alpha+3$&$\alpha \geq \frac{1}{2}$ \\
3'. & $\rho\lvert \cdot \rvert^{x} \rtimes L(\Delta_\rho[-\alpha-1, -\alpha-2], \Delta_\rho[-\alpha, -\alpha]; \pi_{sc})$ &$|\alpha-1|, \alpha+1, \alpha+3$ &$\alpha \geq \frac{1}{2}$ \\
4. & $\rho\lvert \cdot \rvert^{x} \rtimes L(\Delta_\rho[-\alpha-1, -\alpha-1], \Delta_\rho[-\alpha-1, -\alpha-1], \Delta_\rho[-\alpha, -\alpha]; \pi_{sc})$& $|\alpha-1|, \alpha, \alpha+2$&$\alpha \geq \frac{1}{2}$ \\
4'. & $\rho\lvert \cdot \rvert^{x} \rtimes L(\Delta_\rho[-\alpha-1, -\alpha-1], T_{I,1}^{\alpha+1}(T_{I,1}^{\alpha}(\pi_{sc})))$& $|\alpha-1|, \alpha, \alpha+2$&$\alpha \geq \frac{1}{2}$ \\
5. & $\rho\lvert \cdot \rvert^{x} \rtimes L(\Delta_\rho[-\alpha-1, -\alpha-1], \Delta_\rho[-\alpha-1, -\alpha-1]; T_{I,1}^{\alpha}(\pi_{sc}))$& $|\alpha-1|, \alpha, \alpha+2$&$\alpha \geq \frac{1}{2}$ \\
5'. & $\rho\lvert \cdot \rvert^{x} \rtimes L(\Delta_\rho[-\alpha-1, -\alpha-1], \Delta_\rho[-\alpha, -\alpha-1];\pi_{sc})$& $|\alpha-1|, \alpha, \alpha+2$&$\alpha \geq \frac{1}{2}$ \\

\hline
6. & $\rho\lvert \cdot \rvert^{x} \rtimes L(\Delta_\rho[-\alpha-1, -\alpha-1], \Delta_\rho[-\alpha, -\alpha],\Delta_\rho[-\alpha, -\alpha]; \pi_{sc})$& $\alpha-1, \alpha+1, \alpha+2$&$\alpha \geq 1$ \\
6'. & $\rho\lvert \cdot \rvert^{x} \rtimes L(\Delta_\rho[-\alpha, -\alpha]; T_{I,1}^{\alpha+1}(T_{I,1}^{\alpha}(\pi_{sc})))$& $\alpha-1, \alpha+1, \alpha+2$&$\alpha \geq 1$ \\
7. & $\rho\lvert \cdot \rvert^{x} \rtimes L(\Delta_\rho[-\alpha-1, -\alpha-1], \Delta_\rho[-\alpha, -\alpha]; T_{I,1}^{\alpha}(\pi_{sc}))$& $\alpha-1, \alpha+1, \alpha+2$&$\alpha \geq 1$ \\
7'. & $\rho\lvert \cdot \rvert^{x} \rtimes L(\Delta_\rho[-\alpha, -\alpha-1], \Delta_\rho[-\alpha, -\alpha]; \pi_{sc})$& $\alpha-1, \alpha+1, \alpha+2$&$\alpha \geq 1$ \\

8. & $\rho\lvert \cdot \rvert^{x} \rtimes L(\Delta_\rho[-\alpha, -\alpha-1]; T_{I,1}^{\alpha}(\pi_{sc}))$& $\alpha-1, \alpha+2$ &$\alpha \geq 1$ \\
9. & $\rho\lvert \cdot \rvert^{x} \rtimes L(\Delta_\rho[-\alpha, -\alpha], \Delta_\rho[-\alpha, -\alpha], \Delta_\rho[-\alpha, -\alpha]; \pi_{sc})$&$\alpha-1, \alpha+1$ &$\alpha \geq 1$ \\
9'. & $\rho\lvert \cdot \rvert^{x} \rtimes L(\Delta_\rho[-\alpha, -\alpha], \Delta_\rho[-\alpha, -\alpha]; T_{I,1}^{\alpha}(\pi_{sc}))$&$\alpha-1, \alpha+1$ &$\alpha \geq 1$ \\
\hline
10. & $\rho\lvert \cdot \rvert^{x} \rtimes L(\Delta_\rho[-\alpha-1, -\alpha-1], \Delta_\rho[-\alpha+1, -\alpha+1]; T_{I,1}^{\alpha}(\pi_{sc}))$&$|\alpha-2|, \alpha, \alpha+2$ &$\alpha > 1$ \\
10'. & $\rho\lvert \cdot \rvert^{x} \rtimes L(\Delta_\rho[-\alpha, -\alpha-1], \Delta_\rho[-\alpha+1, -\alpha+1]; \pi_{sc})$& $|\alpha-2|, \alpha, \alpha+2$&$\alpha > 1$ \\
11. & $\rho\lvert \cdot \rvert^{x} \rtimes L(\Delta_\rho[-\alpha-1, -\alpha-1], T_{I,1}^{\alpha-1}( T_{I,1}^{\alpha}(\pi_{sc})))$&$|\alpha-2|, \alpha, \alpha+2$ &$\alpha > 1$ \\
11'. & $\rho\lvert \cdot \rvert^{x} \rtimes L(\Delta_\rho[-\alpha+1, -\alpha-1]; \pi_{sc})$& $|\alpha-2|, \alpha, \alpha+2$&$\alpha > 1$ \\
12. & $\rho\lvert \cdot \rvert^{x} \rtimes L(\Delta_\rho[-\alpha, -\alpha], \Delta_\rho[-\alpha, -\alpha], \Delta_\rho[-\alpha+1, -\alpha+1]; \pi_{sc})$& $|\alpha-2|, \alpha-1, \alpha+1$&$\alpha > 1$ \\
12'. & $\rho\lvert \cdot \rvert^{x} \rtimes L(\Delta_\rho[-\alpha+1, -\alpha]; T_{I,1}^{\alpha}(\pi_{sc}))$&$|\alpha-2|, \alpha-1, \alpha+1$ &$\alpha > 1$ \\
13. & $\rho\lvert \cdot \rvert^{x} \rtimes L(\Delta_\rho[-\alpha, -\alpha], \Delta_\rho[-\alpha+1, -\alpha]; \pi_{sc})$&$|\alpha-2|, \alpha-1, \alpha+1$ &$\alpha > 1$ \\
13'. & $\rho\lvert \cdot \rvert^{x} \rtimes L(\Delta_\rho[-\alpha, -\alpha]; T_{I,1}^{\alpha-1}(T_{I,1}^{\alpha}(\pi_{sc})))$& $|\alpha-2|, \alpha-1, \alpha+1$&$\alpha > 1$ \\

14. & $\rho\lvert \cdot \rvert^{x} \rtimes L(\Delta_\rho[-\alpha, -\alpha], \Delta_\rho[-\alpha+1, -\alpha+1], \Delta_\rho[-\alpha+1, -\alpha+1]; \pi_{sc})$& $|\alpha-2|, \alpha, \alpha+1$&$\alpha > 1$ \\
14'. & $\rho\lvert \cdot \rvert^{x} \rtimes L(\Delta_\rho[-\alpha+1, -\alpha+1], \Delta_\rho[-\alpha+1, -\alpha+1]; T_{I,1}^{\alpha}(\pi_{sc}))$& $|\alpha-2|, \alpha, \alpha+1$&$\alpha > 1$ \\
15. & $\rho\lvert \cdot \rvert^{x} \rtimes L(\Delta_\rho[-\alpha+1, -\alpha], \Delta_\rho[-\alpha+1, -\alpha+1];\pi_{sc})$& $|\alpha-2|, \alpha, \alpha+1$&$\alpha > 1$ \\
15'. & $\rho\lvert \cdot \rvert^{x} \rtimes L(\Delta_\rho[-\alpha+1, -\alpha+1], T_{I,1}^{\alpha-1}(T_{I,1}^{\alpha}(\pi_{sc})))$&$|\alpha-2|, \alpha, \alpha+1$ &$\alpha > 1$ \\
\hline
16. & $\rho\lvert \cdot \rvert^{x} \rtimes L(\Delta_\rho[-2,-2]; T_{V,2}^{\pm}(T_{I,1}^{1}(\pi_{sc})))$& $0,1,3$&$\alpha = 1$ \\
17. & $\rho\lvert \cdot \rvert^{x} \rtimes L(\Delta_\rho[-1,-2]; T_{IV,3}(\pi_{sc}))$& $0,1,3$&$\alpha = 1$ \\
18. & $\rho\lvert \cdot \rvert^{x} \rtimes L(\Delta_\rho[0,-2]; \pi_{sc})$& $0,1,3$&$\alpha = 1$ \\
\hline
19. & $\rho\lvert \cdot \rvert^{x} \rtimes L(\Delta_\rho[-\frac{3}{2}, -\frac{3}{2}]; T_{I,2}^{\frac{1}{2}}(\pi_{sc}))$&$\frac{1}{2}, \frac{3}{2}, \frac{5}{2}$ &$\alpha = \frac{1}{2}$ \\
19'. & $\rho\lvert \cdot \rvert^{x} \rtimes L(\Delta_\rho[-\frac{1}{2}, -\frac{3}{2}], \Delta_\rho[-\frac{1}{2}, -\frac{1}{2}]; \pi_{sc})$& $\frac{1}{2}, \frac{3}{2}, \frac{5}{2}$&$\alpha = \frac{1}{2}$ \\
\hline
20. & $\rho\lvert \cdot \rvert^{x} \rtimes L(\Delta_\rho[-2,-2]; T_{I,1}^{1}(T_{V,2}^{\pm}(\pi_{sc})))$& $0,1,3$ &$\alpha = 0$ \\
20'. & $\rho\lvert \cdot \rvert^{x} \rtimes L(\Delta_\rho[-1,-2], T_{V,2}^{\mp}(\pi_{sc}))$&$0,1,3$  &$\alpha = 0$ \\
21. & $\rho\lvert \cdot \rvert^{x} \rtimes L(\Delta_\rho[-2,-2], \Delta_\rho[0,-1]; \pi_{sc})$&$0,1,3$  &$\alpha = 0$ \\

21.' & $\rho\lvert \cdot \rvert^{x} \rtimes L(\Delta_\rho[0,-2], \pi_{sc})$& $0,1,3$ &$\alpha = 0$ \\
22. & $\rho\lvert \cdot \rvert^{x} \rtimes L(\Delta_\rho[-1,-1],\Delta_\rho[0,-1]; \pi_{sc})$& $0,2$&$\alpha = 0$ \\
\end{longtable}
}
\end{center}

\setlength{\tabcolsep}{4pt} % Adjust column spacing
\renewcommand{\arraystretch}{1.2} % Adjust row spacing

\begin{center}
{\footnotesize
\begin{longtable}{|@{}c@{}|l@{}|l@{}|l@{}|@{}c@{}|}
\caption{Non-unitarizable irreducible subquotients at the reducibility point $x_1$} \label{tab: nonunitsub2} \\
\hline
\text{{\rm N}}$^\circ$ & $\pi$ & $x_1$ & Cases \\
\hline
\endfirsthead

\hline
\text{\small{\rm N}}$^\circ$ & $\pi$ & $x_1$ & Cases \\
\hline
\endhead

\hline
\endfoot

\hline
\endlastfoot
1. & $L(\Delta_\rho[-\alpha-2, -\alpha-2], \Delta_\rho[-\alpha, -\alpha], \Delta_\rho[-\alpha+1, -\alpha+1]; T_{I,1}^{\alpha}(\pi_{sc}))$ & $|\alpha-1|$ & $\alpha \geq \frac{1}{2}$ \\
1. & $L(\Delta_\rho[-\alpha-2, -\alpha-2], \Delta_\rho[-\alpha, -\alpha-1]; T_{I,1}^{\alpha}(\pi_{sc}))$ & $\alpha$ & $\alpha \geq \frac{1}{2}$ \\
2.  & $L(\Delta_\rho[-\alpha-2, -\alpha-2], \Delta_\rho[-\alpha, -\alpha-1], \Delta_\rho[-\alpha+1, -\alpha+1]; \pi_{sc})$ & $|\alpha-1|$ & $\alpha \geq \frac{1}{2}$ \\
2.  & $L(\Delta_\rho[-\alpha-2, -\alpha-2], \Delta_\rho[-\alpha-1, -\alpha-1], \Delta_\rho[-\alpha, -\alpha-1]; \pi_{sc})$ & $\alpha+1$ & $\alpha \geq \frac{1}{2}$ \\\
3.  & $L(\Delta_\rho[-\alpha-2, -\alpha-2], \Delta_\rho[-\alpha+1, -\alpha+1], T_{I,1}^{\alpha+1}(T_{I,1}^{\alpha}(\pi_{sc})))$ & $|\alpha-1|$ & $\alpha \geq \frac{1}{2}$ \\
3.  & $L(\Delta_\rho[-\alpha-2, -\alpha-2], \Delta_\rho[-\alpha-1, -\alpha-1], T_{I,1}^{\alpha+1}(T_{I,1}^{\alpha}(\pi_{sc})))$ & $\alpha+1$ & $\alpha \geq \frac{1}{2}$ \\
4.  & $L(\Delta_\rho[-\alpha-2, -\alpha-2], \Delta_\rho[-\alpha-1, -\alpha-1],\Delta_\rho[-\alpha, -\alpha], \Delta_\rho[-\alpha+1, -\alpha+1];  \pi_{sc})$ & $|\alpha-1|$ & $\alpha \geq \frac{1}{2}$ \\
4.  & $L(\Delta_\rho[-\alpha-2, -\alpha-2], \Delta_\rho[-\alpha-1, -\alpha-1],\Delta_\rho[-\alpha, -\alpha], \Delta_\rho[-\alpha, -\alpha];  \pi_{sc})$ & $\alpha$ & $\alpha \geq \frac{1}{2}$ \\
5.  & $L(\Delta_\rho[-\alpha-1, -\alpha-1], \Delta_\rho[-\alpha-1, -\alpha-1],\Delta_\rho[-\alpha+1, -\alpha+1]; T_{I,1}^{\alpha}(\pi_{sc}))$ & $|\alpha-1|$ & $\alpha \geq \frac{1}{2}$ \\
5.  & $L(\Delta_\rho[-\alpha-1, -\alpha-1], \Delta_\rho[-\alpha-1, -\alpha-1],\Delta_\rho[-\alpha, -\alpha]; T_{I,1}^{\alpha}(\pi_{sc}))$ & $\alpha$ & $\alpha \geq \frac{1}{2}$ \\
\hline
6.  & $L(\Delta_\rho[-\alpha-1, -\alpha-1], \Delta_\rho[-\alpha, -\alpha],\Delta_\rho[-\alpha, -\alpha], \Delta_\rho[-\alpha+1, -\alpha+1]; \pi_{sc}))$ & $|\alpha-1|$ & $\alpha \geq 1$ \\
6.  & $L(\Delta_\rho[-\alpha-1, -\alpha-1], \Delta_\rho[-\alpha-1, -\alpha-1],\Delta_\rho[-\alpha, -\alpha], \Delta_\rho[-\alpha, -\alpha]; \pi_{sc}))$ & $\alpha+1$ & $\alpha \geq 1$ \\
7.  & $L(\Delta_\rho[-\alpha-1, -\alpha-1],  \Delta_\rho[-\alpha, -\alpha]; T_{I,1}^{\alpha-1}(T_{I,1}^{\alpha}(\pi_{sc})))$ & $\alpha-1$ & $\alpha \geq 1$ \\
7.  & $L(\Delta_\rho[-\alpha-1, -\alpha-1],  \Delta_\rho[-\alpha-1, -\alpha-1], \Delta_\rho[-\alpha, -\alpha]; T_{I,1}^{\alpha}(\pi_{sc}))$ & $\alpha+1$ & $\alpha \geq 1$ \\
8. & $L(\Delta_\rho[-\alpha, -\alpha-1],  \Delta_\rho[-\alpha+1, -\alpha+1];  T_{I,1}^{\alpha}(\pi_{sc}))$ & $\alpha-1$ & $\alpha \geq 1$ \\
9. & $L(\Delta_\rho[-\alpha, -\alpha], \Delta_\rho[-\alpha, -\alpha], \Delta_\rho[-\alpha, -\alpha]  \Delta_\rho[-\alpha+1, -\alpha+1];  \pi_{sc})$ & $\alpha-1$ & $\alpha \geq 1$ \\
10. & $L(\Delta_\rho[-\alpha-1, -\alpha-1], \Delta_\rho[-\alpha+1, -\alpha+1], \Delta_\rho[-\alpha+2, -\alpha+2];  T_{I,1}^{\alpha}( \pi_{sc}))$ & $|\alpha-2|$ & $\alpha > 1$ \\
10. & $L(\Delta_\rho[-\alpha-1, -\alpha-1], \Delta_\rho[-\alpha+1, -\alpha];  T_{I,1}^{\alpha}(\pi_{sc}))$ & $\alpha$ & $\alpha > 1$ \\
11. & $L(\Delta_\rho[-\alpha-1, -\alpha-1], \Delta_\rho[-\alpha+2, -\alpha+2];  T_{I,1}^{\alpha-1}(T_{I,1}^{\alpha}(\pi_{sc})))$& $|\alpha-2|$ & $\alpha > 1$ \\
11. & $L(\Delta_\rho[-\alpha-1, -\alpha-1], \Delta_\rho[-\alpha, -\alpha];  T_{I,1}^{\alpha-1}(T_{I,1}^{\alpha}(\pi_{sc})))$& $\alpha$ & $\alpha > 1$ \\
12. & $L(\Delta_\rho[-\alpha, -\alpha], \Delta_\rho[-\alpha, -\alpha],\Delta_\rho[-\alpha+1, -\alpha+1],\Delta_\rho[-\alpha+2, -\alpha+2];  \pi_{sc})$& $|\alpha-2|$ & $\alpha > 1$ \\
12. & $L(\Delta_\rho[-\alpha, -\alpha], \Delta_\rho[-\alpha, -\alpha],\Delta_\rho[-\alpha+1, -\alpha+1],\Delta_\rho[-\alpha+1, -\alpha+1];  \pi_{sc})$& $\alpha-1$ & $\alpha > 1$ \\
13. & $L(\Delta_\rho[-\alpha, -\alpha], \Delta_\rho[-\alpha+1, -\alpha],\Delta_\rho[-\alpha+2, -\alpha+2];  \pi_{sc})$& $|\alpha-2|$ & $\alpha > 1$ \\
13. & $L(\Delta_\rho[-\alpha, -\alpha], \Delta_\rho[-\alpha+1, -\alpha],\Delta_\rho[-\alpha+1, -\alpha+1];  \pi_{sc})$& $\alpha-1$ & $\alpha > 1$ \\
14. & $L(\Delta_\rho[-\alpha, -\alpha], \Delta_\rho[-\alpha+1, -\alpha+1],\Delta_\rho[-\alpha+1, -\alpha+1],\Delta_\rho[-\alpha+2, -\alpha+2];  \pi_{sc})$& $|\alpha-2|$ & $\alpha > 1$ \\
14. & $L(\Delta_\rho[-\alpha, -\alpha], \Delta_\rho[-\alpha, -\alpha],\Delta_\rho[-\alpha+1, -\alpha+1],\Delta_\rho[-\alpha+1, -\alpha+1];  \pi_{sc})$& $\alpha$ & $\alpha > 1$ \\
15. & $L(\Delta_\rho[-\alpha+1, -\alpha], \Delta_\rho[-\alpha+1, -\alpha+1],\Delta_\rho[-\alpha+2, -\alpha+2];  \pi_{sc})$& $|\alpha-2|$ & $\alpha > 1$ \\
15. & $L(\Delta_\rho[-\alpha, -\alpha], \Delta_\rho[-\alpha+1, -\alpha],\Delta_\rho[-\alpha+1, -\alpha+1];  \pi_{sc})$& $\alpha$ & $\alpha > 1$ \\
16. & $L(\Delta_\rho[-2,-2]; T_{I,1}^{2}(T_{V,2}^{\pm}(\pi_{sc})))$ & $1$ &$\alpha = 1$ \\
17. & $L(\Delta_\rho[-1,-2], \Delta_\rho[-1,-1]; T_{IV,3}(\pi_{sc})$ & $1$ &$\alpha = 1$ \\
18. & $L(\Delta_\rho[-1,-1], \Delta_\rho[0,-2]; \pi_{sc})$ & $1$ &$\alpha = 1$ \\
19. & $L(\Delta_\rho[-\frac{3}{2}, -\frac{3}{2}]; T_{I,3}^{\frac{1}{2}}(\pi_{sc}))$ & $\frac{1}{2}$ & $\alpha = \frac{1}{2}$\\
19. & $L(\Delta_\rho[-\frac{3}{2}, -\frac{3}{2}],\Delta_\rho[-\frac{3}{2}, -\frac{3}{2}]; T_{I,2}^{\frac{1}{2}}(\pi_{sc}))$ & $\frac{3}{2}$ & $\alpha = \frac{1}{2}$\\
20. & $L(\Delta_\rho[-2,-2]; T_{I,2}^{1}(T_{V,2}^{\pm}(\pi_{sc})))$ & $1$ & $\alpha = 0$\\
21. & $L(\Delta_\rho[-2,-2], \Delta_\rho[-1,-1], \Delta_\rho[0,-1]; \pi_{sc})$ & $1$ & $\alpha = 0$\\
\end{longtable}
}
\end{center}

It remains to consider two-parameter unitarizable families (see the next section).

\section{Two-parameter unitarizable families }\label{2-parameterseries}
\subsection{Two-parameter complementary series}
In this section, we  classify all possible two-parameter unitarizable families of irreducible, unitarizable representations of corank $4$. 
For this, we need the exhaustive list of all critical type, unitarizable, irreducible representations of corank up to $2$, which can be found in \cite{Tad23}.
To begin our classification of two-parameter families, we consider induced representations of the form 
\[\pi_{(x,y)} = \rho\lvert \cdot \rvert^{x} \times \rho\lvert \cdot \rvert^{y} \rtimes \pi_{u},\]
from Propositions \ref{2dcompseries1} to \ref{2dcompseries7}, where $\pi_{u}$ is critical unitarizable of corank $2$. 

\begin{prop} \label{2dcompseries1}
    Let $\alpha \geq \frac{1}{2}$, and let 
    \[\pi = L(\Delta_\rho[-\alpha-1, -\alpha-1], \Delta_\rho[-\alpha, -\alpha]; \pi_{sc}), \ \ \hat{\pi} = T_{I,1}^{\alpha+1}(T_{I,1}^{\alpha}(\pi_{sc})).\]
    Then for $0 \leq x \leq y$, the two-parameter families
    \[\pi_{(x,y)} = \rho\lvert \cdot \rvert^{x} \times \rho \lvert \cdot \rvert^{y} \rtimes \pi,\]
    and 
    \[\hat{\pi}_{(x,y)} = \rho\lvert \cdot \rvert^{x} \times \rho \lvert \cdot \rvert^{y} \rtimes \hat{\pi}\] are unitary in the following regions: 
    \begin{gather*}
        x+y < 1, \ \ y < \alpha -1, \ \ (\alpha > \frac{3}{2}); \\
        x+1 < y < \alpha - 1, \ \ (\alpha > 2); \\
        y < \frac{1}{2}, \ \ (\alpha = \frac{1}{2}).
    \end{gather*}
\end{prop}

\begin{proof}
    The families $\rho\lvert \cdot \rvert^{x} \rtimes \pi$ and $\rho\lvert \cdot \rvert^{x} \rtimes \hat{\pi}$ are reducible at $|\alpha -1|$ and $\alpha+2$. The other reducibility hyperplanes are 
    \[x+y = 1, \ \ |x-y| = 1.\]
    By Aubert duality we only need to consider the family $\pi_{(x,y)}$. The Figures \ref{fig: l3}, \ref{fig: l2}, \ref{fig: l32}, \ref{fig: l1}, \ref{fig: l12} show the unitarizability of the two families in the region $\R^2_{+} = \{(x,y) \in \R^2: x,y \geq 0\}$ for the cases $\alpha = 3, \alpha = 2, \alpha = \frac{3}{2}, \alpha = 1,  \alpha = \frac{1}{2}$ respectively, as examples.

     When there is only one-parameter, we know exactly the irreducible subquotients of $\rho\lvert \cdot \rvert^x \rtimes \pi$. By induction, we can thus obtain the list of irreducible subquotients of $\pi_{(x,y)}$ (not counting multiplicites) when $(x,y)$ is a critical point. Thus we are able to classify the critical points as either strongly unitary, strongly non-unitary, or neither (denoted as black, white, and light gray balls respectively). This proves that list of unitary connected components we gave is exhaustive for $\alpha > 2$. 

     If a connected component contains a non-strongly unitary point on its boundary, then it's non-unitary. To prove our claim, it suffices to show that the region $C_1$ in Figure \ref{fig: l3} is unitary and the region $C_2$ in Figure \ref{fig: l12} is non-unitary. 
     
     First we prove that the region $C_1$ is unitary. Consider the line 
    $L:\{(x,y) \in C_1: x = 0\}$, which is nonempty and one-dimensional. By applying Step $(3$-$2)$ of Algorithm \ref{alg A bar}, one can see that the region $C_1$ is unitary if and only if any point on $L$ is strongly unitary. This reduces the problem to the unitarizability of the one-dimensional complementary series 
     \[ \rho \lvert \cdot \rvert^{y} \rtimes \pi .\]
     By \cite[Proposition $8.1$]{Tad23}, this is unitarizable for $1 \leq y < \alpha-1$. The conclusion follows. 

    Similarly, when $\alpha = \frac{1}{2}$, the slanted side of $C_2$ contains either the family
    \[\Delta_\rho[\frac{1}{2}, -\frac{1}{2}]\lvert \cdot \rvert^{x} \rtimes \pi,\]
    or 
    the family 
    \[L(\Delta_\rho[\frac{1}{2}, -\frac{1}{2}], \Delta_\rho[-\frac{1}{2}, \frac{1}{2}])\lvert \cdot \rvert^{x} \rtimes \pi,\]
    for $0 < x < \frac{1}{2}$. By Proposition \ref{1dcompleseries}, this is non-unitarizable. This proves that $C_2$ is non-unitary. 
\end{proof}

\begin{figure}[H]
\begin{tikzpicture} [thick, scale=0.7]
\draw[style=dotted,line width=1pt] (1,1) -- (12,1);
\draw[style=dotted,line width=1pt] (1,1) -- (1,12);
\draw[style=dotted,line width=1pt] (1,1) -- (12,12);
\draw[style=dashed,line width=1pt] (1,3) -- (3,1);
\draw[style=dashed,line width=1pt] (1,3) -- (10,12);
\draw[style=dashed,line width=1pt] (3,1) -- (12,10);
\draw[style=dashed,line width=1pt] (1,5) -- (12,5);
\draw[style=dashed,line width=1pt] (5,1) -- (5,12);
\draw[style=dashed,line width=1pt] (10,1) -- (10,12);
\draw[style=dashed,line width=1pt] (1,10) -- (12,10);
\draw (12,1) node[right] {$y=0$};
\draw (1,12) node[above] {$x=0$};
\draw (12.7,12) node[above] {$y=x$};
\draw (5,12) node[above] {$x=\alpha-1$};
\draw (8.7,12) node[above] {$y = x+1$};
\draw (12,5) node[right] {$y=\alpha-1$};
\draw (12,9) node[right] {$y=x-1$};
\draw (10.8,12) node[above] {$x = \alpha +2$};
\draw (12,10) node[right] {$ y=\alpha +2$};
\draw (2.4,2.4) node[above] {$y+x=1$};
\draw (10.6,10.3) node[above] {$(\alpha+2,\alpha+2)$};
\draw (5,5) node[above] {$(\alpha-1,\alpha-1)$};
\draw (7,10.3) node[above] {$(\alpha+1, \alpha+2)$};
\draw (3.7,10.3) node[above] {$(\alpha-1,\alpha+2)$};
\draw (1.6,1.5) node {$C_0$};
\draw (1.6,4.3) node {$C_1$};
\draw (4,7.5) node[above] {$(\alpha-1,\alpha)$};
\fill [pattern=north west lines, pattern color=gray] (1,1) -- (3,1) -- (1,3);
\fill [pattern=north west lines, pattern color=gray] (5,1) -- (3,1) -- (5,3);
\fill [pattern=north west lines, pattern color=gray] (1,3) -- (1,5) -- (3,5);
\draw[line width=2pt] (5,1) -- (5,3);
\draw[line width=2pt] (3,1) -- (5,3);
\draw[line width=2pt] (1,3) -- (3,1);
\draw[line width=2pt] (1,3) -- (3,5);
\draw[line width=2pt] (1,5) -- (3,5);
\shade[shading=ball,ball color=black] (3,5) circle (.09);
\shade[shading=ball,ball color=black] (5,3) circle (.09);
\shade[shading=ball,ball color=lightgray] (7,5) circle (.12);
\shade[shading=ball,ball color=lightgray] (5,7) circle (.12);
\shade[shading=ball,ball color=lightgray] (10,5) circle (.12);
\shade[shading=ball,ball color=lightgray] (5,10) circle (.12);
\shade[shading=ball,ball color=white] (5,5) circle (.19);
\shade[shading=ball,ball color=white] (10,10) circle (.19);
\shade[shading=ball,ball color=white] (10,8) circle (.19);
\shade[shading=ball,ball color=white] (8,10) circle (.19);
\end{tikzpicture}
\caption{Unitarizability for $\pi_{(x,y)}$ (case $\alpha = 3$)} \label{fig: l3}
\end{figure}

\begin{center}
\begin{figure}[H]
\begin{tikzpicture} [thick, scale=0.7]
\draw[style=dotted,line width=1pt] (1,1) -- (12,1);
\draw[style=dotted,line width=1pt] (1,1) -- (1,12);
\draw[style=dotted,line width=1pt] (1,1) -- (12,12);
\draw[style=dashed,line width=1pt] (1,3) -- (3,1);
\draw[style=dashed,line width=1pt] (1,3) -- (10,12);
\draw[style=dashed,line width=1pt] (3,1) -- (12,10);
\draw[style=dashed,line width=1pt] (1,3) -- (12,3);
\draw[style=dashed,line width=1pt] (3,1) -- (3,12);
\draw[style=dashed,line width=1pt] (9,1) -- (9,12);
\draw[style=dashed,line width=1pt] (1,9) -- (12,9);
\draw (12,1) node[right] {$y=0$};
\draw (.8,12) node[above] {$x=0$};
\draw (3.3,12) node[above] {$x=1$};
\draw (6.7,10) node[above] {$y=x+1$};
\draw (12,3) node[right] {$y=1$};
\draw (12,9) node[right] {$y=4$};
\draw (9,12) node[above] {$x=4$};
\draw (10,8) node[right] {$\ y=x-1$};
\draw (2,2) node[above] {$y+x=1$};
\draw (2.4,3) node[above] {$(1,1)$};=
\draw (3,5.1) node[left] {$(1,2)$};
\draw (2.3,10) node[below] {$(1,4)$};
\draw (7.2,9.5) node[left] {$(3,4)$ \ \ };
\draw (9,10) node[below] {$(4,4)$};
\draw (9,12) node[below] {$(4,5)$};
\draw (2,1.5) node[left] {$C_0$};
\fill [pattern=north west lines, pattern color=gray] (1,1) -- (3,1) -- (1,3);
\draw[line width=2pt] (1,3) -- (3,1);

\shade[shading=ball,ball color=black] (3,1) circle (.09);
\shade[shading=ball,ball color=black] (1,3) circle (.09);
\shade[shading=ball,ball color=lightgray] (5,3) circle (.14);
\shade[shading=ball,ball color=lightgray] (3,5) circle (.14);
\shade[shading=ball,ball color=lightgray] (9,11) circle (.14);
\shade[shading=ball,ball color=lightgray] (11,9) circle (.14);
\shade[shading=ball,ball color=lightgray] (3,9) circle (.14);
\shade[shading=ball,ball color=lightgray] (9,3) circle (.14);
\shade[shading=ball,ball color=white] (9,9) circle (.22);
\shade[shading=ball,ball color=white] (3,3) circle (.22);
\shade[shading=ball,ball color=white] (9,7) circle (.22);
\shade[shading=ball,ball color=white] (7,9) circle (.22);
 \end{tikzpicture}
 \caption{Unitarizability for $\pi_{(x,y)}$ (case $\alpha=2$)} \label{fig: l2}
\end{figure}
\end{center}

\begin{figure}[H]
\begin{tikzpicture} [thick, scale=0.55]
\draw[style=dotted,line width=1pt] (0,0) -- (15,0);
\draw[style=dotted,line width=1pt] (0,0) -- (0,15);
\draw[style=dotted,line width=1pt] (0,0) -- (15,15);
\draw[style=dashed,line width=1pt] (0,4) -- (4,0);
\draw[style=dashed,line width=1pt] (0,4) -- (11,15);
\draw[style=dashed,line width=1pt] (4,0) -- (15,11);
\draw[style=dashed,line width=1pt] (0,2) -- (15,2);
\draw[style=dashed,line width=1pt] (2,0) -- (2,15);
\draw[style=dashed,line width=1pt] (14,0) -- (14,14);
\draw[style=dashed,line width=1pt] (0,14) -- (14,14);
\draw (15,0) node[right] {$y=0$};
\draw (-0.3,15.2) node[above] {$x=0$};
\draw (2.3,15) node[above] {$x=\frac12$};
\draw (14.4,15) node[above] {$x=\frac72$};
\draw (15,2) node[right] {$y=\frac12$};
\draw (15,14) node[right] {$y=\frac52$};
\draw (10.8,15) node[above] {$y=x+1$};
\draw (15,11) node[right] {$\ y=x-1$};
\draw (3.4,0.4) node[above] {$y+x=1$};
\draw (3.4,2) node[above] {$(\frac12,\frac12)$};
\draw (2.8,5) node[above] {$(\frac12,\frac52)$};
\draw (11.2,14) node[above] {$(\frac52,\frac72)$};
\draw (2.8,14) node[below] {$(\frac12,\frac72)$};
\draw (14,15.2) node[below] {$(\frac72,\frac72)$};
\shade[shading=ball,ball color=lightgray] (2,2) circle (.15);
\shade[shading=ball,ball color=white] (14,10) circle (.18);
\shade[shading=ball,ball color=white] (10,14) circle (.18);
\shade[shading=ball,ball color=white] (6,2) circle (.18);
\shade[shading=ball,ball color=white] (2,6) circle (.18);
\shade[shading=ball,ball color=lightgray] (2,14) circle (.18);
\shade[shading=ball,ball color=lightgray] (14,2) circle (.18);
\shade[shading=ball,ball color=white] (14,14) circle (.3);
\end{tikzpicture}
\caption{Unitarizability for $\pi_{(x,y)}$ (case $\alpha=\frac32$)} \label{fig: l32}
\end{figure}

\begin{figure}[H]
\begin{tikzpicture} [thick, scale=0.55]
\draw[style=dashed,line width=1pt] (1,1) -- (18,1);
\draw[style=dashed,line width=1pt] (1,1) -- (1,18);
\draw[style=dotted,line width=1pt] (1,1) -- (18,18);
\draw[style=dashed,line width=1pt] (1,5) -- (5,1);
\draw[style=dashed,line width=1pt] (1,5) -- (14,18);
\draw[style=dashed,line width=1pt] (5,1) -- (18,14);
\draw[style=dashed,line width=1pt] (13,1) -- (13,18);
\draw[style=dashed,line width=1pt] (1,13) -- (18,13);
\draw (14,18) node[right] {$y=0$};
\draw (1,18) node[above] {$x=0$};
\draw (14,10) node[right] {$y=x-1$};
\draw (18,13) node[right] {$y=3$};
\draw (8.5,13.9) node[above] {$\ \ \ \ \ \ y=x+1$};
\draw (13,18) node[above] {$x=3$\ \ \ \ \ \ };
\draw (9.3,13) node[below] {$\ \ \ \ \ (2,3)$};
\draw (1,1) node[below] {$(0,0)$};
\draw (1,5) node[below] {$(0,1)$};
\draw (13,13) node[below] {$(3,3)$};
\draw (13,16) node[below] {$(3,4)$};
\draw (3,2) node[right] {$x+y=1$};
\shade[shading=ball,ball color=lightgray] (5,1) circle (.15);
\shade[shading=ball,ball color=lightgray] (1,1) circle (.15);
\shade[shading=ball,ball color=lightgray] (1,5) circle (.15);
\shade[shading=ball,ball color=lightgray] (13,17) circle (.24);
\shade[shading=ball,ball color=lightgray] (17,13) circle (.24);
\shade[shading=ball,ball color=white] (13,9) circle (.4);
\shade[shading=ball,ball color=white] (9,13) circle (.4);
\shade[shading=ball,ball color=white] (13,13) circle (.4);
\end{tikzpicture}
\caption{Unitarizability for $\pi_{(x,y)}$ (case $\alpha=1$)} \label{fig: l1}
\end{figure}

\begin{figure}[H]
\begin{tikzpicture} [thick, scale=0.6]
\draw[style=dotted,line width=1pt] (1,1) -- (15,1);
\draw[style=dotted,line width=1pt] (1,1) -- (1,15);
\draw[style=dotted,line width=1pt] (1,1) -- (15,15);
\draw[style=dashed,line width=1pt] (1,5) -- (5,1);
\draw[style=dashed,line width=1pt] (1,5) -- (11,15);
\draw[style=dashed,line width=1pt] (5,1) -- (15,11);
\draw[style=dashed,line width=1pt] (3,1) -- (3,15);
\draw[style=dashed,line width=1pt] (11,1) -- (11,15);
\draw[style=dashed,line width=1pt] (1,3) -- (15,3);
\draw[style=dashed,line width=1pt] (1,11) -- (15,11);
\draw (15,1) node[right] {$y=0$};
\draw (1,15.2) node[above] {$x=0$};
\draw (12,8) node[right] {$y=x-1$};
\draw (1,4.7) node[right] {$x+y=1$};
\draw (15,3) node[right] {$y=\frac12$};
\draw (15,11) node[right] {$y=\frac52$};
\draw (10.7,15) node[above] {$x=\frac52$};
\draw (3.4,15) node[above] {$x=\frac12$};
\draw (9.4,12) node[above] {$y=x+1$};
\draw (4,3) node[above] {$(\frac12,\frac12)$};
\draw (11.9,11) node[below] {$(\frac52,\frac52)$};
\draw (3.9,7) node[below] {$(\frac12,\frac32)$};
\draw (3.9,11) node[below] {$(\frac12,\frac52)$};
\draw (5.7,11.6) node[below] {$(\frac32,\frac52)$};
\draw (2,2.5) node[below] {$C_0$};
\draw (1.7,4.1) node[below] {$C_{2}$};
\fill [pattern=north west lines, pattern color=gray] (1,1) -- (1,3) -- (3,3) -- (3,1);
\draw[line width=2pt] (3,1) -- (3,3);
\draw[line width=2pt] (1,3) -- (3,3);
\draw[line width=1pt] (3,3) -- (5,1);
\draw[line width=1pt] (3,3) -- (1,5);
\draw[line width=1pt] (5,1) -- (7,3);
\draw[line width=1pt] (1,5) -- (3,7);
\draw[line width=1pt] (1,7) -- (3,7);
\draw[line width=1pt] (7,1) -- (7,3);
\draw[line width=1pt] (3,3) -- (7,3);
\draw[line width=1pt] (3,3) -- (3,7);
\shade[shading=ball,ball color=black] (3,3) circle (.15);
\shade[shading=ball,ball color=lightgray] (7,3) circle (.2);
\shade[shading=ball,ball color=lightgray] (3,7) circle (.2);
\shade[shading=ball,ball color=white] (11,3) circle (.2);
\shade[shading=ball,ball color=white] (3,11) circle (.2);
\shade[shading=ball,ball color=lightgray] (11,7) circle (.2);
\shade[shading=ball,ball color=lightgray] (7,11) circle (.2);
\shade[shading=ball,ball color=white] (11,11) circle (.33);
\end{tikzpicture}
\caption{Unitarizability for $\pi_{(x,y)}$ (case $\alpha=\frac12$)} \label{fig: l12}
\end{figure}

The proofs of Propositions \ref{2dcompseries2} to \ref{2dcompseries7} below are similar to that of Proposition \ref{2dcompseries1}, which we omit. 

\begin{prop}\label{2dcompseries2}
    Let $\alpha \geq \frac{3}{2}$, and let \[\pi = L(\Delta_\rho[-\alpha+1, -\alpha]; \pi_{sc}), \ \ \hat{\pi} = T_{I,1}^{\alpha-1}(T_{I,1}^{\alpha}(\pi_{sc})).\]
    Then for $0 \leq x \leq y$, the two-parameter families
    \[\pi_{(x,y)} = \rho\lvert \cdot \rvert^{x} \times \rho \lvert \cdot \rvert^{y} \rtimes \pi,\]
    and 
    \[\hat{\pi}_{(x,y)} = \rho\lvert \cdot \rvert^{x} \times \rho \lvert \cdot \rvert^{y} \rtimes \hat{\pi}\] are unitary in the following regions: 
    \begin{gather*}
        x+y < 1, \ \ (\alpha \geq \frac{3}{2}).
    \end{gather*}
\end{prop}

\begin{prop}\label{2dcompseries3}
    Let $\alpha \geq \frac{3}{2}$, and let \[\pi = L(\Delta_\rho[-\alpha+1, -\alpha+1]; T_{I,1}^{\alpha}(\pi_{sc})), \ \ \hat{\pi} = L(\Delta_\rho[-\alpha+1, -\alpha+1], \Delta_\rho[-\alpha, -\alpha]; \pi_{sc})\]
    Then for $0 \leq x \leq y$, the two-parameter families 
    \[\pi_{(x,y)} = \rho\lvert \cdot \rvert^{x} \times \rho \lvert \cdot \rvert^{y} \rtimes \pi,\]
    and 
    \[\hat{\pi}_{(x,y)} = \rho\lvert \cdot \rvert^{x} \times \rho \lvert \cdot \rvert^{y} \rtimes \hat{\pi},\] are unitary in the following regions: 
    \begin{gather*}
        y < \frac{1}{2}, \ \ (\alpha = \frac{3}{2}) \\
        x+y < 1, \ \ (\alpha >2).
    \end{gather*}
\end{prop}

\begin{prop}\label{2dcompseries4}
    Let $\alpha = \frac{1}{2}$, and let \[\pi = L(\Delta_\rho[-\frac{1}{2}, -\frac{1}{2}], \Delta_\rho[-\frac{1}{2}, -\frac{1}{2}]; \pi_{sc}), \ \ \hat{\pi} = T_{I,2}^{\frac{1}{2}}(\pi_{sc}).\]
    Then for $0 \leq x \leq y$, the two-parameter families
    \[\pi_{(x,y)} = \rho\lvert \cdot \rvert^{x} \times \rho \lvert \cdot \rvert^{y} \rtimes \pi,\]
    and 
    \[\hat{\pi}_{(x,y)} = \rho\lvert \cdot \rvert^{x} \times \rho \lvert \cdot \rvert^{y} \rtimes \hat{\pi},\] are unitary in the following regions: 
    \begin{gather*}
        y < \frac{1}{2}.
    \end{gather*}
\end{prop}

\begin{prop}\label{2dcompseries5}
    Let $\alpha = \frac{1}{2}$, and let \[\pi = L(\Delta_\rho[-\frac{1}{2}, -\frac{1}{2}], T_{I,1}^{\frac{1}{2}}(\pi_{sc})), \ \ \hat{\pi} = T_{III,2}^{\frac{1}{2}}(\pi_{sc})).\]
    Then for $0 \leq x \leq y$, the two-parameter families
    \[\pi_{(x,y)} = \rho\lvert \cdot \rvert^{x} \times \rho \lvert \cdot \rvert^{y} \rtimes \pi,\]
    and 
    \[\hat{\pi}_{(x,y)} = \rho\lvert \cdot \rvert^{x} \times \rho \lvert \cdot \rvert^{y} \rtimes \hat{\pi}\] are unitary in the following regions: 
    \begin{gather*}
        y < \frac{1}{2}.
    \end{gather*}
\end{prop}

\begin{prop}\label{2dcompseries6}
    Let $\alpha = 0$, and let \[\pi^{\pm} = L(\Delta_\rho[-1,-1]; T_{V,2}^{\pm}(\pi_{sc})), \ \ \hat{\pi}^{\pm} = T_{I,1}^{1}(T_{V,2}^{\mp}(\pi_{sc})).\]
    Then for $0 \leq x \leq y$, the two-parameter families
    \[\pi_{(x,y)}^{\pm} = \rho\lvert \cdot \rvert^{x} \times \rho \lvert \cdot \rvert^{y} \rtimes \pi^{\pm}\]
    and 
    \[\hat{\pi}_{(x,y)}^{\pm} = \rho\lvert \cdot \rvert^{x} \times \rho \lvert \cdot \rvert^{y} \rtimes \hat{\pi}^{\pm}\] are unitary in the following regions: 
    \begin{gather*}
        x+y < 1.
    \end{gather*}
\end{prop}

\begin{prop}\label{2dcompseries7}
    Let $\alpha = 0$, and let \[\pi = L(\Delta_\rho[0,-1]; \pi_{sc}).\]
    Then for $0 \leq x \leq y$, the two-parameter families
    \[\pi_{(x,y)} = \rho\lvert \cdot \rvert^{x} \times \rho \lvert \cdot \rvert^{y} \rtimes \pi\]
    are unitary in the following regions: 
    \begin{gather*}
        x+y < 1.
    \end{gather*}
\end{prop}

Propositions \ref{2dcompseries1} to \ref{2dcompseries7} gives all possible two-parameter complementary series of the form $\rho\lvert \cdot \rvert^{x} \times \rho\lvert \cdot \rvert^{y} \rtimes \pi_{u}$, where $\pi_{u}$ is unitarizable of corank $2$. 
Next, we consider two-parameter families of the form

\[u_{\rho}(a_1, b_1) \lvert \cdot \rvert^{x} \times u_{\rho}(a_2, b_2) \lvert \cdot \rvert^{y} \rtimes \pi_{A},\]
where $\pi_{A}$ is critical unitarizable of corank 1, from Propositions \ref{2dcompseries3+1,1} to \ref{2dcompseries3+1,2} below. 

\begin{prop}\label{2dcompseries3+1,1}
    For $x,y \geq 0$, let 
    \[\pi_{(x,y)} = \Delta_\rho[\frac{1}{2}, -\frac{1}{2}]\lvert \cdot \rvert^{x} \times \rho\lvert \cdot \rvert^{y} \rtimes T_{I,1}^{\alpha}(\pi_{sc}),\]
    and 
    \[\hat{\pi}_{(x,y)} = L(\Delta_\rho[-\frac{1}{2}, -\frac{1}{2}], \Delta_\rho[\frac{1}{2}, \frac{1}{2}])\lvert \cdot \rvert^{x} \times \rho\lvert \cdot \rvert^{y} \rtimes L(\Delta_\rho[-\alpha, -\alpha]; \pi_{sc}).\] Then for $x,y \geq 0$, $\pi_{(x,y)}$ and $\hat{\pi}_{(x,y)}$ are unitary in the following regions: 
    \begin{gather*}
        x+y < \frac{3}{2}, \ \ x < \alpha - \frac{3}{2}, \ \ y < \alpha -1, \ \ (\alpha \geq 2) \\
        y-x > 2, \ \ y < \alpha -1, \ \ (\alpha \geq \frac{7}{2}) \\
        x-y > \frac{3}{2}, \ \ x < \alpha - \frac{3}{2}, \ \ (\alpha \geq \frac{7}{2}) \\
        x < 1, \ \ y < \frac{1}{2}, \ \ (\alpha = \frac{1}{2}) \\
        x < \frac{1}{2}, \ \ y < 1 \ \ (\alpha = 0).
    \end{gather*}
\end{prop}
\begin{proof}
    In this case, the reducibility lines are 
    \[|x-y| = \frac{3}{2}, \ \ x+y = \frac{3}{2}, \ \ y = |\alpha-1|, \alpha+1, \ \ x = |\alpha \pm \frac{3}{2}|, \alpha + \frac{1}{2}.\]
    Similar to before, it suffices to consider the regions bounded by $y< |\alpha-1|, x < \min(|\alpha - \frac{3}{2}|, \alpha + \frac{1}{2})$. 
    The rest of the proof is similar to that of Proposition \ref{2dcompseriessc3}, which we omit. 
\end{proof}

\begin{prop}\label{2dcompseries3+1,2}
    For $x,y \geq 0$, let 
    \[\pi_{(x,y)} = \Delta_\rho[\frac{1}{2}, -\frac{1}{2}]\lvert \cdot \rvert^{x} \times \rho\lvert \cdot \rvert^{y} \rtimes L(\Delta_\rho[-\alpha, -\alpha]; \pi_{sc}),\]
    and 
    \[\hat{\pi}_{(x,y)} = L(\Delta_\rho[-\frac{1}{2}, -\frac{1}{2}], \Delta_\rho[\frac{1}{2}, \frac{1}{2}])\lvert \cdot \rvert^{x} \times \rho\lvert \cdot \rvert^{y} \rtimes T_{I,1}^{\alpha}(\pi_{sc}). \] Then for $x,y \geq 0$, $\pi_{(x,y)}$ and $\hat{\pi}_{(x,y)}$ are unitary in the following regions: 
    \begin{gather*}
        x+y < \frac{3}{2}, \ \ x < |\alpha - \frac{3}{2}|, \ \ y < |\alpha -1|, \ \ (\alpha \geq 2) \\
        y-x > 2, \ \ y < \alpha -1, \ \ (\alpha \geq 3) \\
        x-y > \frac{3}{2}, \ \ x < \alpha - \frac{3}{2}, \ \ (\alpha \geq \frac{7}{2}) \\
        x < \frac{1}{2}, \ \ y < 1 \ \ (\alpha = 0).
    \end{gather*}
\end{prop}
\begin{proof}
    The proof is similar to that of Proposition \ref{2dcompseries3+1,1}, which we omit. 
\end{proof}

At last, we consider the complementary series of the form 
\[u_{\rho}(a_1, b_1) \lvert \cdot \rvert^{x} \times u_{\rho}(a_2, b_2) \lvert \cdot \rvert^{y} \rtimes \pi_{sc},\]
where $\pi_{sc}$ is an irreducible supercuspidal representation of a small rank $G_m$, 
from Propositions \ref{2dcompseriessc1} to \ref{2dcompseriessc3} below. 

\begin{prop}\label{2dcompseriessc1}
    Let 
    \[\pi_{(x,y)} = \Delta_\rho[\frac{1}{2}, -\frac{1}{2}] \lvert \cdot \rvert^{x} \times \Delta_\rho[\frac{1}{2}, -\frac{1}{2}]\lvert \cdot \rvert^{y} \rtimes \pi_{sc},\] and 
    \[\hat{\pi}_{(x,y)} = L(\Delta_\rho[-\frac{1}{2}, -\frac{1}{2}], \Delta_\rho[\frac{1}{2}, \frac{1}{2}]) \lvert \cdot \rvert^{x} \times L(\Delta_\rho[-\frac{1}{2}, -\frac{1}{2}], \Delta_\rho[\frac{1}{2}, \frac{1}{2}])\lvert \cdot \rvert^{y} \rtimes \pi_{sc}.\]
    Then for $0 \leq x \leq y$, $\pi_{(x,y)}$ and $\hat{\pi}_{(x,y)}$ are unitary in the following regions: 
    \begin{gather*}
    y - x > 2, \ \ y < \alpha - \frac{1}{2},  \ \ (\alpha \geq 3) \\
    1 < y-x < 2, \ \ y < \alpha - \frac{1}{2}, \ \ (\alpha \geq \frac{5}{2}) \\
        x+y < 1, \ \ (\alpha \geq \frac{3}{2}) \\
        y < \frac{1}{2} \ \ (\alpha = 0, 1). 
    \end{gather*}
\end{prop}
\begin{proof}
    By Aubert duality, it suffices to show this for $\pi_{(x,y)}$. First note that in the region $0 \leq x \leq y$, the family
    \[ \Delta_\rho[\frac{1}{2}, -\frac{1}{2}] \lvert \cdot \rvert^{x} \times \Delta_\rho[\frac{1}{2}, -\frac{1}{2}]\lvert \cdot \rvert^{y}\]
    is reducible precisely when $|x-y| = 1,2$ or $x+y = 1$. 
    Let $\pi_{x} = \Delta_\rho[\frac{1}{2}, -\frac{1}{2}]\lvert \cdot \rvert^{x} \rtimes \pi_{sc}$, then $\pi_{x}$ is reducible when $x = |\alpha - \frac{1}{2}|, \alpha + \frac{1}{2}$. 
    When irreducible, the one-parameter family $\pi_{x}$ is unitarizable if and only if $0 \leq x < |\alpha - \frac{1}{2}| $.  By unitary induction, one can easily show that the regions above are unitary, by applying Step ($3$-$2$) of Algorithm \ref{alg A bar}. To show that the remaining regions are not unitary, we can simply consider the subquotient
    \[L(\Delta_\rho[-\alpha+1, -\alpha], \Delta_\rho[-\alpha+1, -\alpha]; \pi_{sc}),\]
    which is non-unitarizable for $\alpha \geq \frac{3}{2}$. Since 
    \begin{flalign*}
        L(\Delta_\rho[-\alpha+1, -\alpha], \Delta_\rho[-\alpha+1, -\alpha]; \pi_{sc}) &\leq S_{\rho\lvert \cdot \rvert^{\alpha-1}}^{(2)}(S_{\rho\lvert \cdot \rvert^{\alpha}}^{(2)}(D_{\rho\lvert \cdot \rvert^{\alpha}}^{(2)}(D_{\rho\lvert \cdot \rvert^{\alpha-1}}^{(2)})(\Delta_\rho[\alpha, \alpha-1]\\
        &\times \Delta_\rho[\alpha, \alpha-1] \rtimes \pi_{sc}))) \\
        & \cong \Pi_{(\alpha-1, \alpha-1)} ,
    \end{flalign*}
    the conclusion follows (see the notation in \cite[Section $2.2$]{HLL22}). 
    For $\alpha = \frac{1}{2}$, no regions are unitary since reducibility occurs at $0$. Finally for $\alpha = 0,1$, one can use Proposition \ref{1dcompleseries} to show that the region 
    \[x+y < 1, \ \ y > \frac{1}{2}\]
    is non-unitarizable (use case $N^{\circ} = 30$ in table \ref{tab: redpnt}). For $\alpha = 0$,  This concludes the proof. 
\end{proof}

\begin{prop}\label{2dcompseriessc2}
    Let 
    \[\pi_{(x,y)} = \Delta_\rho[\frac{1}{2}, -\frac{1}{2}] \lvert \cdot \rvert^{x} \times L(\Delta_\rho[-\frac{1}{2}, -\frac{1}{2}], \Delta_\rho[\frac{1}{2}, \frac{1}{2}])\lvert \cdot \rvert^{y} \rtimes \pi_{sc},\] and 
    \[\hat{\pi}_{(x,y)} = L(\Delta_\rho[-\frac{1}{2}, -\frac{1}{2}], \Delta_\rho[\frac{1}{2}, \frac{1}{2}]) \lvert \cdot \rvert^{x} \times \Delta_\rho[\frac{1}{2}, -\frac{1}{2}]\lvert \cdot \rvert^{y} \rtimes \pi_{sc}.\]
    Then for $0 \leq x \leq y$, $\pi_{(x,y)}$ and $\hat{\pi}_{(x,y)}$ are unitary in the following regions: 
    \begin{gather*}
    y - x > 2, \ \ y < \alpha - \frac{1}{2},  \ \ (\alpha \geq 3) \\
    1 < y-x < 2, \ \ y < \alpha - \frac{1}{2}, \ \  (\alpha \geq \frac{5}{2}) \\
        x+y < 1, \ \ (\alpha \geq \frac{3}{2}) \\
        y < \frac{1}{2}, \ \ (\alpha = 0, 1).
    \end{gather*}
\end{prop}
\begin{proof}
    The proof is similar to that of Proposition \ref{2dcompseriessc1}, which we omit. 
\end{proof}

\begin{prop}\label{2dcompseriessc3}
     Let 
    \[\pi_{(x,y)} = L(\Delta_\rho[-1,-1],\Delta_\rho[0,0], \Delta_\rho[1,1])\lvert \cdot \rvert^{x} \times  \rho\lvert \cdot \rvert^{y} \rtimes \pi_{sc},\] and 
    \[\hat{\pi}_{(x,y)} = \Delta_\rho[1,-1]\lvert \cdot \rvert^{x} \times \rho\lvert \cdot \rvert^{y} \rtimes \pi_{sc}.\]
    Then for $x,y \geq 0$, $\pi_{(x,y)}$ and $\hat{\pi}_{(x,y)}$ are unitary in the following regions: 
    \begin{gather*}
        x-y > 2, \ \ x < \alpha-1, \ \ (\alpha \geq \frac{7}{2}) \\
        x+y < 2, \ \ x < \alpha -1, \ \ (\alpha \geq 2)\\
        y-x > 2, \ \ y < \alpha, \ \ (\alpha \geq \frac{5}{2}) \\
       y < \frac{3}{2}, \ \ x < \frac{1}{2}, \ \ (\alpha = \frac{3}{2}) \\
       x,y < \frac{1}{2} \ \ (\alpha = \frac{1}{2}).
    \end{gather*}
\end{prop}
\begin{proof}
    Same as before, it suffices to consider $\pi_{(x,y)}$. First, note that in this case, the reducibility lines in the region $x,y \geq 0$ are 
    \[|x-y| = 2, \ \ x+y = 2, \ \ x = |\alpha-1|, \alpha, \alpha+1, \ \ y = \alpha. \]
    In particular, the connected components are asymmetric in the two coordinates. In the one-parameter families 
    \[L(\Delta_\rho[-1,-1],\Delta_\rho[0,0], \Delta_\rho[1,1])\lvert \cdot \rvert^{x} \rtimes \pi_{sc},\] and 
    \[\rho\lvert \cdot \rvert^{y} \rtimes \pi_{sc}\]
    respectively, there are no unitarizable points beyond the first reducibility point. 

    Thus, by applying Step ($3$-$2$) of Algorithm \ref{alg A bar}, it is straightforward to show that the regions listed are unitary. To show that other regions are not unitarizable, it suffices to consider the subquotient 
    \[L(\Delta_\rho[-\alpha, -\alpha], \Delta_\rho[-\alpha+2, -\alpha]; \pi_{sc}),\]
    which is a non-unitarizable subquotient of $\Pi_{(\alpha-1, \alpha)}$ for $\alpha \geq \frac{3}{2}$. This proves the claim for $\alpha \geq \frac{3}{2}$. 

    When $\alpha = 1$ or $\alpha = 0$, reducibility occurs at $x = 0$, so there are no unitary regions. When $\alpha = \frac{1}{2}$, the first reducibility point in $x$ and $y$ are both $\frac{1}{2}$, and the component 
    \[x , y < \frac{1}{2}\]
    contains the origin. This proves the claim. 
\end{proof}

This completes our consideration of the two-parameter complementary series.

\subsection{Other two-parameter families}
In this subsection, we show that there are no unitary regions in the 2-dimensional continuous families induced from non-unitarizable, critical type representations. 

\begin{prop}\label{2dcompseiresnonunit}
    Let 
    \[\pi_{(x,y)} = \rho\lvert \cdot \rvert^{x} \times \rho\lvert \cdot \rvert^{y} \rtimes \pi_{nu}, \]
    where $\pi_{nu}$ is critical non-unitarizable of corank $2$. Then when $\pi_{(x,y)}$ is irreducible, it is non-unitarizable.
\end{prop}
\begin{proof}
    This follows directly by considering the unitarizability of the one-parameter family 
    \[\rho\lvert \cdot \rvert^{x} \rtimes \pi_{nu},\] 
    when $\pi_{nu}$ is non-unitarizable, of critical type, and of corank $2$. It is proven by Tadi{\'c} in \cite{Tad23} that such families are non-unitarizable whenever they are irreducible. When $\pi_{(x,y)}$ is irreducible, $(x,y)$ lies in an irreducible component. If the component has nontrivial intersection with $x = 0$ or $y = 0$, then we can use Step ($3$-$2$) of Algorithm \ref{alg A bar} to reduce the problem of unitarizability to the one-parameter case, which is known. Otherwise, the connected component has a boundary point of the form $(x_0, y_0)$, where $x_0, y_0$ is a nonzero reducibility point of $\rho\lvert \cdot \rvert^{x} \rtimes \pi_{nu}$ and $\rho\lvert \cdot \rvert^{y} \rtimes \pi_{nu}$ respectively. In this case, one can first show that $\rho\lvert \cdot \rvert^{x} \rtimes \pi_{nu}$ has a non-unitarizable subquotient $\pi_{x_0}$ of corank $3$ at $x = x_0$ by the result of Tadi{\'c} in \cite[Proposition $8.1$]{Tad23}. Then, using the results of Proposition \ref{1dcompnonunitcrnk3}, one can show that a non-unitarizable subquotient of $\rho\lvert \cdot \rvert^{y_0} \rtimes \pi_{x_0}$ exists. 
    This proves the claim. 
\end{proof}

Since the lowest possible corank of a non-unitarizable representation  of critical type is $2$, the above proposition gives the only possible two-parameter families induced from a non-unitarizable representation of critical type. Therefore, this concludes Step $2$ of Algorithm \ref{alg A bar}, and we are now ready to state our final result in the next section. 

\section{Conclusion}\label{conclusion}
In this section, we  finally compile and state the main results of this paper. In Theorem \ref{finalconclusion}, we give the full unitary dual of corank $4$ (up to a $W$-orbit), and in Theorem \ref{uequalAbar}, we show that the full unitary dual is indeed equal to the set $\Pi_{\overline{A}}^{\lim}(G_n)$, when restricted to representations of corank $4$. Note that for presentation purposes, some unitarizable representations may be listed in Theorem \ref{finalconclusion}  more than once. In particular, the unitarizable representations in Proposition \ref{crnk4Artcritlist} may also appear in a complementary series. 

\begin{thm}\label{finalconclusion}
    The irreducible unitarizable subquotients of $\Pi_{\underline{x}}$, when $\underline{x} = (x_1, x_2, x_3, x_4) \in \R^{4}_{++}$, are precisely the following: 
    \begin{enumerate}
        \item $(\alpha \geq 1)$ All irreducible subquotients of $\Pi_{\underline{x}}$ when $\underline{x}$ lies in the closure of the connected components (\ref{ciu1}) to (\ref{ciu7}). 
        \item $(\alpha \geq 0)$ All irreducible subquotients of critical points $\underline{x}$ listed in Proposition \ref{crnk4Artcritlist}. 
        \item $(\alpha \geq 0)$ All one parameter complementary series of the form $u_{\rho}(a,b) \lvert \cdot \rvert^{x_1} \rtimes \pi_{A}$ listed in Table \ref{tab: redpnt}, for $0 \leq x_1 < \text{FRP}$, where FRP is the first nonzero reducibility point, given in the table. 
        \item $(\alpha = \frac{1}{2})$ All irreducible subquotients when $x_4 \leq \frac{1}{2}$. 
        \item $(\alpha \geq \frac{1}{2})$
        \begin{itemize}
            \item $\underline{x} = (x_1, x_2, \alpha, \alpha+1)$. The complementary series 
            \[\rho\lvert \cdot \rvert^{x_1} \times \rho\lvert \cdot \rvert^{x_2} \times L(\Delta_\rho[-\alpha-1, -\alpha-1], \Delta_\rho[-\alpha, -\alpha]; \pi_{sc}),\]
            and 
            \[\rho\lvert \cdot \rvert^{x_1} \times \rho\lvert \cdot \rvert^{x_2} \rtimes T_{I,1}^{\alpha+1}(T_{I,1}^{\alpha}(\pi_{sc})),\]
            for 
            \begin{gather*}
                x_1 + x_2 < 1, \ \ x_2 < \alpha -1, \ \ (\alpha > \frac{3}{2})\\
                x_1 +1 < x_2 < \alpha-1, \ \ (\alpha > 2) \\
                x_2 < \frac{1}{2} \ \ (\alpha = \frac{1}{2}).
            \end{gather*}
            \item $(\alpha \geq \frac{3}{2}), \underline{x} = (x_1, x_2, \alpha-1, \alpha)$. The complementary series  
            \[\rho\lvert \cdot \rvert^{x_1} \times \rho\lvert \cdot \rvert^{x_2} \times L(\Delta_\rho[-\alpha+1, -\alpha]; \pi_{sc}),\]
            and 
            \[\rho\lvert \cdot \rvert^{x_1} \times \rho\lvert \cdot \rvert^{x_2} \rtimes T_{I,1}^{\alpha-1}(T_{I,1}^{\alpha}(\pi_{sc})),\]
            for $x_1 + x_2 < 1$. 
            \item $(\alpha \geq \frac{3}{2}), \underline{x} = (x_1, x_2, \alpha-1, \alpha)$. The complementary series 
            \[\rho\lvert \cdot \rvert^{x_1} \times \rho\lvert \cdot \rvert^{x_2} \times L(\Delta_\rho[-\alpha+1, -\alpha+1]; T_{I,1}^{\alpha}(\pi_{sc})),\]
            and 
            \[\rho\lvert \cdot \rvert^{x_1} \times \rho\lvert \cdot \rvert^{x_2} \rtimes L(\Delta_\rho[-\alpha+1, -\alpha+1], \Delta_\rho[-\alpha, -\alpha]; \pi_{sc}),\]
            for \begin{gather*}
                x_2 < \frac{1}{2}, \ \ (\alpha = \frac{3}{2}) \\
                x_1 + x_2 < 1, \ \ (\alpha > 2).
            \end{gather*}
        \end{itemize}
        \item ($\alpha \geq 0$): 
        \begin{itemize}
            \item  The complementary series 
        \[\Delta_\rho[\frac{1}{2}, -\frac{1}{2}]\lvert \cdot \rvert^{x_1} \times \Delta_\rho[\frac{1}{2}, -\frac{1}{2}]\lvert \cdot \rvert^{x_2} \rtimes \pi_{sc},\] or its Aubert dual, when 
        \begin{gather*}
    x_2 - x_1 > 2, \ \ x_2 < \alpha - \frac{1}{2},  \ \ (\alpha \geq 3) \\
    1 < x_2 - x_1 < 2, \ \ x_2 < \alpha - \frac{1}{2}, \ \ (\alpha \geq \frac{5}{2}) \\
        x_1+x_2 < 1, \ \ (\alpha \geq \frac{3}{2}) \\
        x_2 < \frac{1}{2} \ \ (\alpha = 0, 1).
    \end{gather*}
    \item The complementary series 
    \[\Delta_\rho[\frac{1}{2}, -\frac{1}{2}] \lvert \cdot \rvert^{x_1}  \times L(\Delta_\rho[-\frac{1}{2}, -\frac{1}{2}], \Delta_\rho[\frac{1}{2}, \frac{1}{2}])\lvert \cdot \rvert^{x_2} \rtimes \pi_{sc},\] or its Aubert dual, when 
     \begin{gather*}
    x_2 - x_1 > 2, \ \ x_2 < \alpha - \frac{1}{2},  \ \ (\alpha \geq 3) \\
    1 < x_2-x_1 < 2, \ \ x_2 < \alpha - \frac{1}{2}, \ \ (\alpha \geq \frac{5}{2}) \\
        x_1+x_2 < 1, \ \ (\alpha \geq \frac{3}{2}) \\
        x_2 < \frac{1}{2} \ \ (\alpha = 0, 1).
    \end{gather*}
    \item The complementary series 
    \[L(\Delta_\rho[-1,-1], \Delta_\rho[0,0], \Delta_\rho[1,1])\lvert \cdot \rvert^{x_1} \times \rho\lvert \cdot \rvert^{x_2} \rtimes \pi_{sc}\] or its Aubert dual, when 
    \begin{gather*}
        x_1 - x_2 > 2, \ \ x_1 < \alpha-1, \ \ (\alpha \geq \frac{7}{2}) \\
        x_1+x_2 < 2, \ \ x_1 < \alpha -1, \ \ (\alpha \geq 2)\\
        x_2-x_1 > 2, \ \ x_2 < \alpha, \ \ (\alpha \geq 3) \\
       x_2 < \frac{3}{2}, \ \ x_1 < \frac{1}{2}, \ \ (\alpha = \frac{3}{2}) \\
       x_1,x_2 < \frac{1}{2} \ \ (\alpha = \frac{1}{2}).
    \end{gather*}
    \item The complementary series 
    \[\Delta_\rho[\frac{1}{2}, -\frac{1}{2}]\lvert \cdot \rvert^{x_1} \times \rho\lvert \cdot \rvert^{x_2} \rtimes T_{I,1}^{\alpha}(\pi_{sc}),\] when 
    \begin{gather*}
        x_1+x_2 < \frac{3}{2}, \ \ x_1 < \alpha - \frac{3}{2}, \ \ x_2 < \alpha -1, \ \ (\alpha \geq 2) \\
        x_2-x_1 > 2, \ \ x_2 < \alpha -1, \ \ (\alpha \geq 3) \\
        x_2-x_1 > \frac{3}{2}, \ \ x_1 < \alpha - \frac{3}{2}, \ \ (\alpha \geq \frac{7}{2}) \\
        x_1 < 1, \ \ x_2 < \frac{1}{2}, \ \ (\alpha = \frac{1}{2}) \\
        x_1 < \frac{1}{2}, \ \ x_2 < 1 \ \ (\alpha = 0).
    \end{gather*}
    \item The complementary series 
    \[\Delta_\rho[\frac{1}{2}, -\frac{1}{2}]\lvert \cdot \rvert^{x_1} \times \rho\lvert \cdot \rvert^{x_2} \rtimes L(\Delta_\rho[-\alpha, -\alpha]; \pi_{sc}),\] or its Aubert dual, when
    \begin{gather*}
        x_1+x_2 < \frac{3}{2}, \ \ x_1 < |\alpha - \frac{3}{2}|, \ \ x_2 < |\alpha -1|, \ \ (\alpha \geq 2) \\
        x_2-x_1 > 2, \ \ x_2 < \alpha -1, \ \ (\alpha \geq 3) \\
        x_2-x_1 > \frac{3}{2}, \ \ x_1 < \alpha - \frac{3}{2}, \ \ (\alpha \geq \frac{7}{2}) \\
        x_1< \frac{1}{2}, \ \ x_2 < 1, \ \ (\alpha = 0).
    \end{gather*}
        \end{itemize}
        \item $(\alpha = \frac{1}{2}),  \underline{x} = (x_1, x_2, \frac{1}{2}, \frac{1}{2})$. The complementary series 
        \[\rho\lvert \cdot \rvert^{x_1} \times \rho\lvert \cdot \rvert^{x_2} \rtimes \pi,\]
        where $\pi = L(\Delta_\rho[-\frac{1}{2}, -\frac{1}{2}], \Delta_\rho[-\frac{1}{2}, -\frac{1}{2}]; \pi_{sc}), T_{I,2}^{\frac{1}{2}}(\pi_{sc}), L(\Delta_\rho[-\frac{1}{2}, -\frac{1}{2}]; T_{I,1}^{\frac{1}{2}}(\pi_{sc})), T_{III,2}^{\frac{1}{2}}(\pi_{sc})$, when
        \[0 \leq x \leq y \leq \frac{1}{2}.\]
        \item $(\alpha = \frac{1}{2})$
        \begin{itemize}
            \item $\underline{x} = (x_1, \frac{1}{2}, \frac{1}{2}, \frac{3}{2})$ The complementary series 
        \[\rho\lvert \cdot \rvert^{x_1} \rtimes \pi,\]
        where $\pi = L(\Delta_\rho[-\frac{3}{2}, -\frac{3}{2}], \Delta_\rho[-\frac{1}{2}, -\frac{1}{2}], \Delta_\rho[-\frac{1}{2}, -\frac{1}{2}]; \pi_{sc}),  \\ L(\Delta_\rho[\frac{1}{2}, -\frac{3}{2}]; \pi_{sc}), L(\Delta_\rho[-\frac{1}{2}, -\frac{1}{2}]; T_{I,1}^{\frac{3}{2}}(T_{I,1}^{\frac{1}{2}}(\pi_{sc}))), L(\Delta_\rho[-\frac{1}{2}, -\frac{3}{2}]; T_{I,1}^{\frac{1}{2}}(\pi_{sc}))$ or their Aubert duals, when $x_1 < \frac{1}{2}$.
        \item $\underline{x} = (\frac{1}{2}, \frac{1}{2}, x_3, \frac{3}{2})$. The complementary series 
        \[\rho\lvert \cdot \rvert^{x_3} \rtimes \pi,\]
        for $\frac{1}{2} < x_3 < \frac{3}{2}$, where $\pi = L(\Delta_\rho[-\frac{3}{2}, -\frac{3}{2}], \Delta_\rho[-\frac{1}{2}, -\frac{1}{2}]; T_{I,1}^{\frac{1}{2}}(\pi_{sc}))$ or its Aubert dual.
        \item $\underline{x} = (x_1, \frac{1}{2}, \frac{1}{2}, \frac{1}{2})$ The complementary series 
        \[\rho\lvert \cdot \rvert^{x_1} \rtimes \pi,\]
        where $\pi = T_{I,3}^{\frac{1}{2}}(\pi_{sc}), L(\Delta_\rho[-\frac{1}{2}, -\frac{1}{2}], \Delta_\rho[-\frac{1}{2}, -\frac{1}{2}]; T_{I,1}^{\frac{1}{2}}(\pi_{sc}))$ or their Aubert duals, when $x_1 < \frac{1}{2}$.
        \item $\underline{x} = (\frac{1}{2}, \frac{1}{2}, \frac{1}{2}, x_4)$. The complementary series 
        \[\rho\lvert \cdot \rvert^{x_4} \rtimes L(\Delta_\rho[-\frac{1}{2}, -\frac{1}{2}]; T_{I,2}^{\frac{1}{2}}(\pi_{sc})),\]
        for $\frac{1}{2} < x_4 < \frac{3}{2}$. 

        \end{itemize}
    \item $(\alpha = 0)$ All irreducible subquotients when $\underline{x} = (0,0,x_3, x_4)$ with $x_3 + x_4  \leq 1$. 
    \item $(\alpha = 0), \underline{x} = (0,x_2, x_3, 1)$. 
        The complementary series 
        \[\rho\lvert \cdot \rvert^{x_2} \times \rho\lvert \cdot \rvert^{x_3} \rtimes \pi,\]
        for $x_2 + x_3 < 1$, where $\pi = L(\Delta_\rho[-1,-1]; T_{V,2}^{\pm}(\pi_{sc})), L(\Delta_\rho[0,-1]; \pi_{sc})$ or their Aubert duals. 
    \end{enumerate}
\end{thm}
\begin{proof}
    The unitarity of the above representations follow from previous propositions. It suffices to prove that the list is exhaustive, and by symmetry under the $W$-action, it suffices to consider the subquotients of $\Pi_{\underline{x}}$ when $\underline{x} \in \R^{4}_{++}$. If $\underline{x} \in \mathbb{R}^{4}_{++}$ is regular and unitarizable, then it must be contained in one of the connected components (\ref{ciu1}) to (\ref{ciu7}), as proven in Proposition \ref{exhaustconclusion}. 

    Otherwise, if $\underline{x} \in \mathbb{R}^{4}_{++}$ is not regular, then there are two possibilities. The first possibility occurs when the $W$-orbit of $\underline{x}$ intersects the image of a unitary$^{\pm}$ connected component of  $\R^3_{\reg,\shrt,++}$  (resp. $R^{3}_{\reg,\lev,++}$), under the affine isomorphism $\iota$ defined in Sections \ref{slantedsection} (resp. \ref{levelsection}). In this case, as stated in the proof of Proposition \ref{exhaustconclusion}, $\underline{x}$ lies in the closure of the connected components (\ref{ciu1}) to (\ref{ciu7}). 

    From now on, we assume the second possibility, which is when the $W$-orbit of $\underline{x}$ intersects neither $\iota(R^{3}_{\reg,\shrt,++})$, nor $\iota(\R^3_{\reg,\lev,++})$. In this case, if an irreducible representation $\pi$ is unitarizable of corank $4$, it must either be a subquotient of $\Pi_{\underline{x}}$ when $\underline{x}$ is an isolated point, or it lies in a continuous family of representations induced from a representation on the maximal Levi of $G_n$ which is of critical type. Any point $\underline{x} \in \R^4_{++}$ is either a regular point, or it is connected to a point of critical type by a one-parameter family which lies completely in a reducibility hyperplane. This implies that all isolated points are of critical type. Therefore, when $\underline{x}$ is isolated, the unitarizable subquotients of $\Pi_{\underline{x}}$ are listed in Proposition \ref{crnk4Artcritlist}. 

    When $\underline{x}$ is not an isolated point, it must be part of a continuous family of hermitian representations induced from a representation of critical type. By Proposition \ref{exhaustconclusion} and Remark \ref{n-1dcompseries}, it suffices to consider families with one or two parameters. 

    Suppose $\pi$ appears as part of a one-parameter complementary series induced from a unitarizable representation of critical type, then there are a few possibilities. First, suppose $\pi$ lies in a complementary series of the form $u_{\rho}(a,b) \lvert \cdot \rvert^{x} \rtimes \pi_{A}$, where $\pi_{A}$ is of critical type and of Arthur type, then it must be in one of the cases in table \ref{tab: redpnt}, which are all unitarizable from $0$ up to the first nonzero reducibility point. The endpoints of such complementary series which are unitarizable are already included in Proposition \ref{crnk4Artcritlist}. One can verify that all of the relevant complementary series are listed above. Otherwise, we have one of the cases below: 
    
    \begin{enumerate}
        \item $\pi$ lies in a continuous family of the form $\tau \lvert \cdot \rvert^{x} \rtimes \pi_{sc}$, where $\tau \neq u_{\rho}(a,b)$ for any $a,b$.
        \item $\pi \leq \tau \lvert \cdot \rvert^{x} \rtimes \pi_{A}$, where $\pi_{A}$ is of critical type and unitarizable of corank $1$, and $\tau \neq u_{\rho}(a,b)$ for any $a,b$.
    \end{enumerate} 
    
    Case $(1)$ is dealt with in Propositions \ref{1dcompscextra} to \ref{1dcompscextra,4}, and the relevant unitarizable subquotients are included. Case $(2)$ is dealt with in Propositions \ref{1dcompextra3+1,1} to \ref{1dcompextra3+1,3}, where we proved that in this case none of the subquotients are unitarizable. 
    
    Otherwise, $\pi$ is part of a hermitian family induced from a non-unitarizable representation of critical type. Since there are no non-unitarizable, critical type subquotients of corank $1$, the following two possibilities remain: 
    \begin{enumerate}
        \item $\pi \leq \rho \lvert \cdot \rvert^{x} \rtimes \pi_{nu}$, where $\pi_{nu}$ is non-unitarizable of corank $3$. 
        \item $\pi \leq \Delta_\rho[x+1,x] \rtimes \pi_{nu}$, where $\pi_{nu}$ is non-unitarizable of corank $2$. 
    \end{enumerate}

   Case $(1)$ is covered in Proposition \ref{1dcompnonunitcrnk3}, and case $(2)$ is covered by Propositions \ref{1dcompnonunit1} and \ref{1dcompnonunit2}. In each case,  we proved all subquotients are non-unitarizable. 
    
   Now, we consider the two-parameter families. First, consider those $2$-parameter complementary series. There are three possibilities: 
   \begin{enumerate}
       \item $\pi \leq \rho\lvert \cdot \rvert^{x} \times \rho\lvert \cdot \rvert^{y} \rtimes \pi_{A}$, where $\pi_{A}$ is unitarizable, of critical type, and of corank $2$. 
       \item $\pi \leq u_{\rho}(a,b) \lvert \cdot \rvert^{x} \times \rho\lvert \cdot \rvert^{y} \rtimes \pi_{A}$, where $\pi_{A}$ is unitarizable, of critical type, and of corank $1$.
       \item $\pi \leq u_{\rho}(a_1,b_1) \lvert \cdot \rvert^{x} \times u_{\rho}(a_2, b_2) \rho\lvert \cdot \rvert^{y} \rtimes \pi_{sc}$. 
   \end{enumerate}
   Case $(1)$ is covered by Propositions \ref{2dcompseries1} to \ref{2dcompseries7}. When $\alpha > 3$, one can check that these complementary series lie in the closure of (\ref{ciu4}), and when $2 < \alpha \leq 3$, they lie in the closure of (\ref{ciu5}). Case $(2)$ is covered by Propositions \ref{2dcompseries3+1,1} and \ref{2dcompseries3+1,2}, and case $(3)$ is covered by Propositions \ref{2dcompseriessc1} to \ref{2dcompseriessc3}. 
   
   Finally, consider those two-parameter hermitian families induced from non-unitarizable representations. Since there are no non-unitarizable, critical type representations of corank $1$, the only possibility we have are families of the form 
   \[\rho\lvert \cdot \rvert^{x} \times \rho\lvert \cdot \rvert^{y} \rtimes \pi_{nu},\]
   where $\pi_{nu}$ is critical type, non-unitarizable of corank $2$. By Proposition \ref{2dcompseiresnonunit}, none of the subquotients in families of this form are unitarizable. This concludes the proof. 
\end{proof}

\begin{remark}\label{isolated reps}
    By Theorem \ref{vogan thm}, one can see that the isolated representations in the list of Theorem \ref{finalconclusion} are exactly the representations listed in Proposition \ref{crnk4Artcritlist} which do not appear in complementary series or limits of complementary series.
\end{remark}

Now that we've constructed the full unitary dual for corank $4$ representations, we show that it is indeed the same as $\Pi_{\overline{A}}^{\lim}(G_n)$, as we expected in Conjecture \ref{unitary dual conjecture}. 

\begin{thm}\label{uequalAbar}
    If $\pi \in \Pi_{u}(G_n)$ is a unitary representation of corank $4$, then $\pi \in \Pi_{\overline{A}}^{\lim}(G_n)$. In particular, this implies that $\Pi_{u}(G_n) = \Pi_{\overline{A}}^{\lim}(G_n)$ when restricted to corank $4$ representations. 
\end{thm}
\begin{proof}
    First, we consider the standard case, where we take unitary representations of $G_n$ that is a subquotient of 
    \[\rho\lvert\cdot\rvert^{x_1} \times \rho\lvert\cdot\rvert^{x_2} \times \rho\lvert\cdot\rvert^{x_3} \times\rho\lvert\cdot\rvert^{x_4} \rtimes \pi_{sc}. \]
    In this case, it suffices to show that all representations given in Theorem \ref{finalconclusion} lie in $\Pi_{\overline{A}}^{\lim}(G_n)$. To prove this, we first note that the closure of all irreducible subquotients of $\Pi_{\underline{x}}$, where $\underline{x}$ lie in the closure of a unitary connected component in $\R^4_{++}$, lie in $\Pi_{\overline{A}}(G_n)$. 

    By definition, since all one-parameter complementary series listed in Proposition \ref{1dcompleseries} are of the form $u_{\rho}(a,b) \lvert \cdot \rvert^{x_1} \rtimes \Pi_{A}$, where $x_1$ lies in an irreducible region,  they are contained in $\Pi_{\overline{A}}^{\lim}(G_n)$. 

    Similarly, the two-parameter families listed in Theorem \ref{finalconclusion} are all of the form 
    \[u_{\rho}(a_1, b_1) \lvert \cdot \rvert^{x_1} \times u_{\rho}(a_2, b_2) \lvert \cdot \rvert^{x_2} \rtimes \pi_{A},\]
    where $(x_1, x_2)$ lies in an irreducible component $\R^2$ for the given parabolic induction. Therefore, they all lie in $\Pi_{\overline{A}}^{\lim}(G_n)$. 

    It remains to consider the mixed case, where $\pi \in \Pi_{u}(G_n)$ is an irreducible subquotient of 
    \[\rho_1\lvert\cdot\rvert^{x_1} \times \rho_2\lvert\cdot\rvert^{x_2} \times \rho_3\lvert\cdot\rvert^{x_3} \times\rho_4\lvert\cdot\rvert^{x_4} \rtimes \pi_{sc},\]
    for $\rho_i \in \mathcal{C}^{sd}$ not necessarily equal. Without loss of generality, let us assume that all $\rho_i$ are distinct (the case where they are not all distinct can be argued similarly). Let $X_{\rho_i} = \{\rho_i\lvert\cdot\rvert^x : x \in \R\}$ for $1 \leq i \leq 4$. 
    Then for each $i$, there exists a unique $\gamma \in \Irr(X_{\rho_i}; \pi_{sc})$ (see Definition \ref{defJantzendecomp}) such that for some $\beta$ supported in $\sqcup_{j \neq i} X_{\rho_j}$, we have 
    \[\pi \hookrightarrow \beta \rtimes \gamma.\]
    Denote this $\gamma$ by $X_{\rho_i}(\pi)$, then the correspondence $\pi \mapsto (X_{\rho_1}(\pi), X_{\rho_2}(\pi), X_{\rho_3}(\pi), X_{\rho_4}(\pi))$ determines a bijection by Theorem \ref{thm Jantzen}. By the results of \cite{HJLLZ24}, one sees that a representation of corank $\leq 3$ is non-unitarizable if and only if it lies inside a continuous family of irreducible, hermitian representations that contains a good parity subquotient which is not of Arthur type in its interior or closure. As Jantzen's decomposition preserves irreducibility, one can show that it also preserves unitarizability (i.e. $\pi$ is unitarizable if and only if $X_{\rho_i}(\pi)$ is unitarizable for all $1 \leq i \leq 4$), using a similar argument as in \cite[Chapter 9]{Tad23}. We omit the details here. This completes the proof. 
    \end{proof}

% We note that since our construction of the unitary dual follows precisely the steps in Algorithm \ref{alg A bar}, we may conclude that for representations of corank $4$, the unitary dual $\Pi_{u}(G_n) \subset \Pi_{\overline{A}}(G_n)$. By definition we always have, $\Pi_{\overline{A}}(G_n) \subset \Pi_{u}(G_n)$, which implies that 
% \[\Pi_{u}(G_n) = \Pi_{\overline{A}}(G_n),\]
% verifying conjecture ... in this case. 

Note that in Definition \ref{def closure A}, we defined two distinct sets $\Pi_{\overline{A}}^{\lim}(G_n)$ and $\Pi_{\overline{A}}^{\lim'}(G_n)$. By \cite[Theorem 5.2]{HJLLZ24}, we have $\Pi_{\overline{A}}^{\lim}(G_n) \subset \Pi_{\overline{A}}^{\lim'}(G_n) \subset \Pi_{u}(G_n)$. On the other hand, by Theorem \ref{uequalAbar}, in the case of corank $4$ representations,  $\Pi_{\overline{A}}^{\lim}(G_n) = \Pi_{u}(G_n)$. 
Therefore, in this case,  $\Pi_{\overline{A}}^{\lim}(G_n) = \Pi_{\overline{A}}^{\lim'}(G_n) = \Pi_{u}(G_n)$. Hence, we have that:

\begin{thm}\label{main conjecture holds}
    When $G_n = \Sp_{2n}$ or $\SO_{2n+1}$, Conjecture \ref{unitary dual conjecture} holds for representations of corank $4$. 
\end{thm}

\appendix
\section{\\List of representations of corank 4 that are of critical type and of Arthur type, sorted by cuspidal support} \label{artlist}

To explicitly construct the unitary dual for corank $4$ representations, we need a list of critical points $(x_1, x_2, x_3, x_4)$ ($0 \leq x_1 \leq x_2 \leq x_3 \leq x_4$) such that the corresponding induced representations
$$\Pi_{x_1, x_2, x_3, x_4} = \rho\lvert \cdot \rvert^{x_1} \times \rho\lvert \cdot \rvert^{x_2} \times \rho\lvert \cdot \rvert^{x_3} \times \rho\lvert \cdot \rvert^{x_4} \rtimes \pi_{sc}$$
only contain unitarizable subquotients, or, contains non-unitariable subquotients. To this end, in this appendix, we 
provide the complete list of corank 4 representations that are of critical type and of Arthur type, sorted by cuspidal support.
In the next appendix, we provide the opposite complete list of corank 4 representations that are of critical type but not of Arthur type.
These two lists help us determine whether a connected component is unitary.

\begin{prop} \label{crnk4Artcritlist}
    Let $\Pi_{x_1, x_2, x_3, x_4} = \rho\lvert \cdot \rvert^{x_1} \times \rho\lvert \cdot \rvert^{x_2} \times \rho\lvert \cdot \rvert^{x_3} \times \rho\lvert \cdot \rvert^{x_4} \rtimes \pi_{sc}$ for $0 \leq x_1 \leq x_2 \leq x_3 \leq x_4$. Then for any given critical point $(x_1, x_2, x_3, x_4)$, the following list contains all the irreducible subquotients of $\Pi_{x_1, x_2, x_3, x_4}$ which are of Arthur type. 
    \begin{enumerate}
        \item $(\alpha = 0)$: 
        \begin{enumerate}
        \item $(0,0,0,1)$: 
        \begin{itemize}
        \item $T_{I,1}^{1}(T_{V,6}^{\pm}(\pi_{sc}))$
            \item  $L(\Delta_\rho[-1,-1]; T_{V,6}^{\pm}(\pi_{sc}))$
            \item $L(\Delta_\rho[0,-1]; T_{V,4}^{\pm}(\pi_{sc}))$
        \end{itemize}
        \item $(0,0,1,1)$: 
        \begin{itemize}
        \item $T_{I,2}^{1}(T_{V,4}^{\pm}(\pi_{sc}))$
        \item $T_{II,3}^{1}(T_{V,4}^{\pm}(\pi_{sc}))$
            \item  $L(\Delta_\rho[-1,-1]; T_{I,1}^{1}(T_{V,4}^{\pm}(\pi_{sc})))$
            \item $L(\Delta_\rho[0,-1]; T_{I,1}^{1}(T_{V,2}^{\pm}(\pi_{sc})))$
            \item $L(\Delta_\rho[0,-1], \Delta_\rho[0,-1]; \pi_{sc})$  
            \item $L(\Delta_\rho[-1,-1], \Delta_\rho[-1,-1]; T_{V,4}^{\pm}(\pi_{sc}))$
            \item $L(\Delta_\rho[-1,-1], \Delta_\rho[0,-1]; T_{V,2}^{\pm}(\pi_{sc}))$
        \end{itemize}
        \item $(0,0,1,2)$: 
        \begin{itemize}
        \item $T_{I,1}^{2}(T_{I,1}^{1}(T_{V,4}^{\pm}(\pi_{sc})))$
            \item $L(\Delta_\rho[-2,-2], \Delta_\rho[-1,-1]; T_{V,4}^{\pm}(\pi_{sc}))$
        \end{itemize}
            \item $(0,1,1,2)$: 
            \begin{itemize}
            \item $T_{I,1}^{1}(T_{I,1}^{2}(T_{I,1}^{1}(T_{V,2}^{\pm}(\pi_{sc}))))$
                \item $L(\Delta_\rho[-1,-1]; T_{I,1}^{2}(T_{I,1}^{1}(T_{V,2}^{\pm}(\pi_{sc}))))$
                
                \item $L(\Delta_\rho[1,-2]; \pi_{sc})$ 
                \item $L(\Delta_\rho[-1,-2], \Delta_\rho[0,-1]; \pi_{sc})$  
                \item $L(\Delta_\rho[-2,-2], \Delta_\rho[-1,-1]; T_{V,2}^{\pm}(T_{I,1}^{1}(\pi_{sc})))$
                \item $L(\Delta_\rho[-2,-2], \Delta_\rho[-1,-1], \Delta_\rho[-1,-1]; T_{V,2}^{\pm}(\pi_{sc}))$
            \end{itemize}
            \item $(0,1,2,3)$: 
            \begin{itemize}
            \item $T_{I,1}^{3}(T_{I,1}^{2}(T_{I,1}^{1}(T_{V,2}^{\pm}(\pi_{sc}))))$
                \item $L(\Delta_\rho[-3,-3], \Delta_\rho[-2,-2], \Delta_\rho[-1,-1]; T_{V,2}^{\pm}(\pi_{sc}))$
            \end{itemize}
        \end{enumerate}
        \item $(\alpha = \frac{1}{2})$: 
        \begin{enumerate}
        \item $(\frac{1}{2}, \frac{1}{2}, \frac{1}{2}, \frac{1}{2})$:
        \begin{itemize}
        \item $T_{I,4}^{\frac{1}{2}}(\pi_{sc})$
        \item $T_{II,5}^{\frac{1}{2}}(\pi_{sc})$
        \item $L(\Delta_\rho[-\frac{1}{2}, -\frac{1}{2}]; T_{I,3}^{\frac{1}{2}}(\pi_{sc})$
        \item $L(\Delta_\rho[-\frac{1}{2}, -\frac{1}{2}],\Delta_\rho[-\frac{1}{2}, -\frac{1}{2}]; T_{II,3}^{\frac{1}{2}}(\pi_{sc}))$
            \item $L(\Delta_\rho[-\frac{1}{2}, -\frac{1}{2}],\Delta_\rho[-\frac{1}{2}, -\frac{1}{2}]; T_{I,2}^{\frac{1}{2}}(\pi_{sc}))$
            \item $L(\Delta_\rho[-\frac{1}{2}, -\frac{1}{2}],\Delta_\rho[-\frac{1}{2}, -\frac{1}{2}],\Delta_\rho[-\frac{1}{2}, -\frac{1}{2}]; T_{I,1}^{\frac{1}{2}}(\pi_{sc}))$
            \item $L(\Delta_\rho[-\frac{1}{2}, -\frac{1}{2}], \Delta_\rho[-\frac{1}{2}, -\frac{1}{2}], \Delta_\rho[-\frac{1}{2}, -\frac{1}{2}], \Delta_\rho[-\frac{1}{2}, -\frac{1}{2}]; \pi_{sc})$
            
        \end{itemize}
        \item $(\frac{1}{2}, \frac{1}{2}, \frac{1}{2}, \frac{3}{2})$: 
        \begin{itemize}
        \item $T_{I,2}^{\frac{1}{2}}(T_{I,1}^{\frac{3}{2}}(T_{I,1}^{\frac{1}{2}}(\pi_{sc})))$
        \item $T_{III,2}^{\frac{1}{2}}(T_{I,1}^{\frac{3}{2}}(T_{I,1}^{\frac{1}{2}}(\pi_{sc})))$
        \item $L(\Delta_\rho[-\frac{1}{2}, -\frac{1}{2}]; T_{I,1}^{\frac{1}{2}}(T_{I,1}^{\frac{3}{2}}(T_{I,1}^{\frac{1}{2}}(\pi_{sc}))))$
        \item $L(\Delta_\rho[-\frac{1}{2}, -\frac{1}{2}]; T_{I,1}^{\frac{3}{2}}(T_{II,3}^{\frac{1}{2}}(\pi_{sc}))))$
        \item $L(\Delta_\rho[\frac{1}{2}, -\frac{3}{2}]; T_{I,1}^{\frac{1}{2}}(\pi_{sc}))$ 
            \item $L(\Delta_\rho[-\frac{1}{2}, -\frac{3}{2}]; T_{I,2}^{\frac{1}{2}}(\pi_{sc}))$
            \item $L(\Delta_\rho[-\frac{1}{2}, -\frac{3}{2}]; T_{III,2}^{\frac{1}{2}}(\pi_{sc}))$
            \item $L(\Delta_\rho[-\frac{1}{2}, -\frac{1}{2}], \Delta_\rho[-\frac{1}{2}, -\frac{1}{2}]; T_{I,1}^{\frac{3}{2}}(T_{I,1}^{\frac{1}{2}}(\pi_{sc})))$
             \item $L(\Delta_\rho[-\frac{1}{2}, -\frac{1}{2}], \Delta_\rho[\frac{1}{2}, -\frac{3}{2}]; \pi_{sc})$  
             \item $L(\Delta_\rho[-\frac{1}{2}, -\frac{3}{2}], \Delta_\rho[-\frac{1}{2}, -\frac{1}{2}]; T_{I,1}^{\frac{1}{2}}(\pi_{sc}))$ 
            \item $L(\Delta_\rho[-\frac{3}{2}, -\frac{3}{2}], \Delta_\rho[-\frac{1}{2}, -\frac{1}{2}]; T_{I,2}^{\frac{1}{2}}(\pi_{sc}))$
             \item $L(\Delta_\rho[-\frac{3}{2}, -\frac{3}{2}], \Delta_\rho[-\frac{1}{2}, -\frac{1}{2}]; T_{III,2}^{\frac{1}{2}}(\pi_{sc}))$ 
            \item $L(\Delta_\rho[-\frac{3}{2}, -\frac{3}{2}],\Delta_\rho[-\frac{1}{2}, -\frac{1}{2}],\Delta_\rho[-\frac{1}{2}, -\frac{1}{2}] ; T_{I,1}^{\frac{1}{2}}(\pi_{sc}))$
            \item $L(\Delta_\rho[-\frac{3}{2}, -\frac{3}{2}], \Delta_\rho[-\frac{1}{2}, -\frac{1}{2}], \Delta_\rho[-\frac{1}{2}, -\frac{1}{2}], \Delta_\rho[-\frac{1}{2}, -\frac{1}{2}]; \pi_{sc})$
        \end{itemize}
         \item $(\frac{1}{2}, \frac{1}{2}, \frac{3}{2}, \frac{3}{2})$: 
         \begin{itemize}
         \item $T_{I,2}^{\frac{3}{2}}(T_{I,2}^{\frac{1}{2}}(\pi_{sc}))$ 
         \item $T_{I,2}^{\frac{3}{2}}(T_{III,2}^{\frac{1}{2}}(\pi_{sc}))$ 
         \item $L(\Delta_\rho[-\frac{3}{2}, -\frac{3}{2}]; T_{I,1}^{\frac{3}{2}}(T_{III,2}^{\frac{1}{2}}(\pi_{sc})))$  
             \item $L(\Delta_\rho[-\frac{3}{2}, -\frac{3}{2}], \Delta_\rho[-\frac{1}{2}, -\frac{1}{2}]; T_{I,1}^{\frac{3}{2}}(T_{I,1}^{\frac{1}{2}}(\pi_{sc})))$
             \item $L(\Delta_\rho[-\frac{3}{2}, -\frac{3}{2}], \Delta_\rho[-\frac{3}{2}, -\frac{3}{2}] , \Delta_\rho[-\frac{1}{2}, -\frac{1}{2}]; T_{I,1}^{\frac{1}{2}}(\pi_{sc}))  $
             \item $L(\Delta_\rho[-\frac{3}{2}, -\frac{3}{2}],\Delta_\rho[-\frac{3}{2}, -\frac{3}{2}],\Delta_\rho[-\frac{1}{2}, -\frac{1}{2}],\Delta_\rho[-\frac{1}{2}, -\frac{1}{2}]; \pi_{sc})$ 
         \end{itemize}
            \item $(\frac{1}{2}, \frac{1}{2}, \frac{3}{2}, \frac{5}{2})$: 
            \begin{itemize}
            \item $T_{I,1}^{\frac{5}{2}}(T_{I,1}^{\frac{3}{2}}(T_{I,2}^{\frac{1}{2}}(\pi_{sc}))$
            \item $T_{I,1}^{\frac{5}{2}}(T_{I,1}^{\frac{3}{2}}(T_{III,2}^{\frac{1}{2}}(\pi_{sc}))$
                \item $L(\Delta_\rho[-\frac{1}{2}, -\frac{1}{2}]; T_{I,1}^{\frac{5}{2}}(T_{I,1}^{\frac{3}{2}}(T_{I,1}^{\frac{1}{2}}(\pi_{sc}))))$
                \item $L(\Delta_\rho[-\frac{5}{2}, -\frac{5}{2}],\Delta_\rho[-\frac{3}{2}, -\frac{3}{2}],\Delta_\rho[-\frac{1}{2}, -\frac{1}{2}]; T_{I,1}^{\frac{1}{2}}(\pi_{sc}))  $
                \item $L(\Delta_\rho[-\frac{3}{2}, -\frac{5}{2}], \Delta_\rho[-\frac{1}{2}, -\frac{1}{2}] ,\Delta_\rho[-\frac{1}{2}, -\frac{1}{2}]; \pi_{sc})  $
                \item $L(\Delta_\rho[-\frac{5}{2}, -\frac{5}{2}],\Delta_\rho[-\frac{3}{2}, -\frac{3}{2}] ,\Delta_\rho[-\frac{1}{2}, -\frac{1}{2}],\Delta_\rho[-\frac{1}{2}, -\frac{1}{2}]; \pi_{sc})   $
            \end{itemize}
            \item $(\frac{1}{2}, \frac{3}{2}, \frac{5}{2}, \frac{7}{2})$: 
            \begin{itemize}
            \item $T_{I,1}^{\frac{7}{2}}(T_{I,1}^{\frac{5}{2}}(T_{I,1}^{\frac{3}{2}}(T_{I,1}^{\frac{1}{2}}(\pi_{sc}))))$
                \item $L(\Delta_\rho[-\frac{7}{2}, -\frac{7}{2}],\Delta_\rho[-\frac{5}{2}, -\frac{5}{2}] ,\Delta_\rho[-\frac{3}{2}, -\frac{3}{2}] ,\Delta_\rho[-\frac{1}{2}, -\frac{1}{2}]; \pi_{sc})  $
            \end{itemize}
        \end{enumerate}
        \item $(\alpha = 1)$: 
        \begin{enumerate}
        \item $(0,0,0,1)$ :  
        \begin{itemize}
        \item $T_{V,6}^{\pm}(T_{I,1}^{1}(\pi_{sc}))$
            \item $L(\Delta_\rho[-1,-1]; T_{IV,7}(\pi_{sc}))$
            \item $L(\Delta_\rho[0,-1]; T_{IV,5}(\pi_{sc}))$
        \end{itemize}
        \item $(0,0,1,1)$: 
        \begin{itemize}
        \item $T_{I,2}^{1}(T_{IV,5}(\pi_{sc}))$
        \item $T_{II,3}^{1}(T_{IV,5}(\pi_{sc}))$
            \item $L(\Delta_\rho[-1,-1]; T_{V,4}^{\pm}(T_{I,1}^{1}(\pi_{sc})))$
            \item $L(\Delta_\rho[0,-1]; T_{V,2}^{\pm}(T_{I,1}^{1}(\pi_{sc})))$
            \item $L(\Delta_\rho[0,-1], \Delta_\rho[0,-1]; \pi_{sc})$  
            \item $L(\Delta_\rho[-1,-1], \Delta_\rho[-1,-1]; T_{IV,5}(\pi_{sc}))$
            \item $L(\Delta_\rho[-1,-1], \Delta_\rho[0,-1]; T_{IV,3}(\pi_{sc}))$
        \end{itemize}
        \item $(0,0,1,2)$: 
        \begin{itemize}
            \item $T_{V,4}^{\pm}(T_{I,1}^{2}(T_{I,1}^{1}(\pi_{sc})))$
            \item $L(\Delta_\rho[-2,-2], \Delta_\rho[-1,-1]; T_{IV,5}(\pi_{sc}))$  
            \item $L(\Delta_\rho[-2,-2], \Delta_\rho[0,-1]; T_{IV,3}(\pi_{sc}))$ 
        \end{itemize}
        \item $(0,1,1,1)$: 
        \begin{itemize}
        \item $T_{I,3}^{1}(T_{IV,3}(\pi_{sc}))$
            \item $L(\Delta_\rho[-1,-1]; T_{I,2}^{1}(T_{IV,3}(\pi_{sc})))$
            \item $L(\Delta_\rho[-1,-1]; T_{III,3}^{1}(T_{IV,3}(\pi_{sc})))$
            \item $L(\Delta_\rho[-1,-1], \Delta_\rho[-1,-1]; T_{V,2}^{\pm}(T_{I,1}^{1}(\pi_{sc})))$
            \item $L(\Delta_\rho[-1,-1], \Delta_\rho[-1,-1], \Delta_\rho[-1,-1]; T_{IV,3}(\pi_{sc}))$
        \end{itemize}
            \item $(0,1,1,2)$: 
            \begin{itemize}
            \item $T_{V,2}^{\pm}(T_{I,1}^{1}(T_{I,1}^{2}(T_{I,1}^{1}(\pi_{sc}))))$
            \item $T_{I,1}^{2}(T_{III,2}^{1}(\pi_{sc}))$
                \item $L(\Delta_\rho[-1,-1]; T_{V,2}^{\pm}(T_{I,1}^{2}(T_{I,1}^{1}(\pi_{sc}))))$
                \item $L(\Delta_\rho[0,-1]; T_{I,1}^{2}(T_{I,1}^{1}(\pi_{sc})))$
                \item  $L(\Delta_\rho[-1,-2]; T_{V,2}^{\pm}(T_{I,1}^{1}(\pi_{sc})))$ with $\epsilon_{sc}(\rho \otimes S_3) = \mp 1$.
                \item $L(\Delta_\rho[0,-2]; T_{I,1}^{1}(\pi_{sc}))$
                \item $L(\Delta_\rho[1,-2]; \pi_{sc})$
                \item $L(\Delta_\rho[-1,-2], \Delta_\rho[0,-1]; \pi_{sc})$  
                \item $L(\Delta_\rho[-2,-2]; T_{II,3}^{1}(T_{IV,3}(\pi_{sc})))$
                \item $L(\Delta_\rho[-2,-2], \Delta_\rho[-1,-1]; T_{V,4}^{\pm}(T_{I,1}^{1}(\pi_{sc}))))$
                \item $L(\Delta_\rho[-2,-2], \Delta_\rho[-1,-1], \Delta_\rho[0,-1]; \pi_{sc})$  
                \item $L(\Delta_\rho[-2,-2], \Delta_\rho[-1,-1], \Delta_\rho[-1,-1]; T_{IV,3}(\pi_{sc}))$
            \end{itemize}
            \item $(0,1,2,3)$: 
            \begin{itemize}
            \item $T_{V,2}^{\pm}(T_{I,1}^{3}(T_{I,1}^{2}(T_{I,1}^{1}(\pi_{sc}))))$
                \item $L(\Delta_\rho[-3,-3], \Delta_\rho[-2,-2], \Delta_\rho[-1,-1]; T_{IV,3}(\pi_{sc}))$
                \item $L(\Delta_\rho[-3,-3], \Delta_\rho[-2,-2], \Delta_\rho[0,-1]; \pi_{sc})$
            \end{itemize}
            \item $(1,2,3,4)$: 
            \begin{itemize}
            \item $T_{I,1}^{4}(T_{I,1}^{3}(T_{I,1}^{2}(T_{I,1}^{1}(\pi_{sc}))))$
                \item $L(\Delta_\rho[-4,-4], \Delta_\rho[-3,-3], \Delta_\rho[-2,-2], \Delta_\rho[-1,-1]; \pi_{sc})$
            \end{itemize}
        \end{enumerate}
        \item $(\alpha = \frac{3}{2})$: 
        \begin{enumerate}
        \item $(\frac{1}{2}, \frac{1}{2}, \frac{1}{2}, \frac{3}{2})$: 
        \begin{itemize}
        \item $T_{I,3}^{\frac{1}{2}}(T_{I,1}^{\frac{3}{2}}(\pi_{sc}))$ 
            \item $L(\Delta_\rho[-\frac{1}{2}, -\frac{1}{2}]; T_{I,1}^{\frac{3}{2}}(T_{II,3}^{\frac{1}{2}}(\pi_{sc})))$
            \item $L(\Delta_\rho[-\frac{1}{2}, -\frac{1}{2}]; T_{I,2}^{\frac{1}{2}}(T_{I,1}^{\frac{3}{2}}(\pi_{sc})))$  
            \item $L(\Delta_\rho[-\frac{1}{2}, -\frac{3}{2}]; T_{II,3}^{\frac{1}{2}}(\pi_{sc}))$
            \item $L(\Delta_\rho[-\frac{1}{2}, -\frac{1}{2}], \Delta_\rho[\frac{1}{2}, -\frac{3}{2}]; \pi_{sc})$  
            \item $L(\Delta_\rho[-\frac{3}{2}, -\frac{3}{2}], \Delta_\rho[-\frac{1}{2}, -\frac{1}{2}]; T_{II,3}^{\frac{1}{2}}(\pi_{sc}))$
            \item $L(\Delta_\rho[-\frac{1}{2}, -\frac{1}{2}],\Delta_\rho[-\frac{1}{2}, -\frac{1}{2}]; T_{I,1}^{\frac{1}{2}}(T_{I,1}^{\frac{3}{2}}(\pi_{sc})))$
            \item $L(\Delta_\rho[-\frac{1}{2}, -\frac{1}{2}], \Delta_\rho[-\frac{1}{2}, -\frac{1}{2}], \Delta_\rho[-\frac{1}{2}, -\frac{1}{2}]; T_{I,1}^{\frac{3}{2}}(\pi_{sc}))$
            \item $L(\Delta_\rho[-\frac{1}{2}, -\frac{3}{2}], \Delta_\rho[-\frac{1}{2}, -\frac{1}{2}],\Delta_\rho[-\frac{1}{2}, -\frac{1}{2}]; \pi_{sc})$
            \item $L(\Delta_\rho[-\frac{3}{2}, -\frac{3}{2}], \Delta_\rho[-\frac{1}{2}, -\frac{1}{2}], \Delta_\rho[-\frac{1}{2}, -\frac{1}{2}], \Delta_\rho[-\frac{1}{2}, -\frac{1}{2}]; \pi_{sc})$
        \end{itemize}
        \item $(\frac{1}{2}, \frac{1}{2}, \frac{3}{2}, \frac{3}{2})$: 
        \begin{itemize}
        \item $T_{I,2}^{\frac{3}{2}}(T_{II,3}^{\frac{1}{2}}(\pi_{sc}))$
        \item $T_{III,2}^{\frac{3}{2}}(\pi_{sc})$
        \item $L(\Delta_\rho[\frac{1}{2}, -\frac{3}{2}]; T_{I,1}^{\frac{3}{2}}(\pi_{sc}))$  
        \item $L(\Delta_\rho[-\frac{1}{2}, -\frac{3}{2}]; T_{I,1}^{\frac{1}{2}}(T_{I,1}^{\frac{3}{2}}(\pi_{sc})))$
        \item $L(\Delta_\rho[-\frac{3}{2}, -\frac{3}{2}]; T_{I,1}^{\frac{3}{2}}(T_{II,3}^{\frac{1}{2}}(\pi_{sc})))$
            \item $L(\Delta_\rho[-\frac{3}{2}, -\frac{3}{2}], \Delta_\rho[-\frac{1}{2}, -\frac{1}{2}]; T_{I,1}^{\frac{1}{2}}(T_{I,1}^{\frac{3}{2}}(\pi_{sc})))$
            \item $L(\Delta_\rho[-\frac{3}{2}, -\frac{3}{2}], \Delta_\rho[\frac{1}{2}, -\frac{3}{2}];\pi_{sc})$ 
            \item $L(\Delta_\rho[-\frac{3}{2}, -\frac{3}{2}],\Delta_\rho[-\frac{1}{2}, -\frac{1}
            {2}],\Delta_\rho[-\frac{1}{2}, -\frac{1}{2}]; T_{I,1}^{\frac{3}{2}}(\pi_{sc}))$
            \item $L(\Delta_\rho[-\frac{3}{2}, -\frac{3}{2}],\Delta_\rho[-\frac{1}{2}, -\frac{3}{2}],\Delta_\rho[-\frac{1}{2}, -\frac{1}{2}]; \pi_{sc})$
            \item $L(\Delta_\rho[-\frac{3}{2}, -\frac{3}{2}], \Delta_\rho[-\frac{1}{2}, -\frac{1}{2}], \Delta_\rho[-\frac{1}{2}, -\frac{1}{2}]; T_{I,1}^{\frac{3}{2}}(\pi_{sc}))$
        \end{itemize}
        \item $(\frac{1}{2}, \frac{1}{2}, \frac{3}{2}, \frac{5}{2})$: 
        \begin{itemize}
        \item $T_{I,2}^{\frac{1}{2}}(T_{I,1}^{\frac{5}{2}}(T_{I,1}^{\frac{3}{2}}(\pi_{sc})))$
        \item $T_{II,3}^{\frac{1}{2}}(T_{I,1}^{\frac{5}{2}}(T_{I,1}^{\frac{3}{2}}(\pi_{sc})))$
        \item $L(\Delta_\rho[\frac{1}{2}, -\frac{5}{2}]; \pi_{sc})$ + 
            \item $L(\Delta_\rho[-\frac{5}{2}, -\frac{5}{2}],\Delta_\rho[-\frac{3}{2}, -\frac{3}{2}]; T_{II,3}^{\frac{1}{2}}(\pi_{sc}))$
            \item $L(\Delta_\rho[-\frac{5}{2}, -\frac{5}{2}], \Delta_\rho[\frac{1}{2}, -\frac{3}{2}]; \pi_{sc})$ 
            \item $L(\Delta_\rho[-\frac{1}{2}, -\frac{1}{2}],\Delta_\rho[-\frac{1}{2}, -\frac{1}{2}]; T_{I,1}^{\frac{5}{2}}(T_{I,1}^{\frac{3}{2}}(\pi_{sc})))$
            \item $L(\Delta_\rho[-\frac{3}{2}, -\frac{5}{2}], \Delta_\rho[-\frac{1}{2}, -\frac{1}{2}],\Delta_\rho[-\frac{1}{2}, -\frac{1}{2}]; \pi_{sc})$
            \item $L(\Delta_\rho[-\frac{5}{2}, -\frac{5}{2}], \Delta_\rho[-\frac{1}{2}, -\frac{3}{2}], \Delta_\rho[-\frac{1}{2}, -\frac{1}{2}]; \pi_{sc})$  
        \end{itemize}
        \item $(\frac{1}{2}, \frac{3}{2}, \frac{3}{2}, \frac{5}{2})$:
        \begin{itemize}
        \item $T_{I,1}^{\frac{3}{2}}(T_{I,1}^{\frac{1}{2}}(T_{I,1}^{\frac{5}{2}}(T_{I,1}^{\frac{3}{2}}(\pi_{sc}))))$
            
            \item $L(\Delta_\rho[-\frac{3}{2}, -\frac{3}{2}], \Delta_\rho[-\frac{1}{2}, -\frac{1}{2}]; T_{I,1}^{\frac{5}{2}}(T_{I,1}^{\frac{3}{2}}(\pi_{sc})))$
            \item $L(\Delta_\rho[-\frac{3}{2}, -\frac{5}{2}], \Delta_\rho[-\frac{1}{2}, -\frac{3}{2}]; \pi_{sc})$
            \item $L(\Delta_\rho[-\frac{5}{2}, -\frac{5}{2}],\Delta_\rho[-\frac{3}{2}, -\frac{3}{2}],\Delta_\rho[-\frac{1}{2}, -\frac{1}{2}]; T_{I,1}^{\frac{3}{2}}(\pi_{sc}))$
        \end{itemize}
            \item $(\frac{1}{2}, \frac{3}{2}, \frac{5}{2}, \frac{7}{2}):$
            \begin{itemize}
            \item $T_{I,1}^{\frac{1}{2}}(T_{I,1}^{\frac{7}{2}}(T_{I,1}^{\frac{5}{2}}(T_{I,1}^{\frac{3}{2}}(\pi_{sc}))))$
                \item $L(\Delta_\rho[-\frac{1}{2}, -\frac{1}{2}]; T_{I,1}^{\frac{7}{2}}(T_{I,1}^{\frac{5}{2}}(T_{I,1}^{\frac{3}{2}}(\pi_{sc}))))$
                \item $L(\Delta_\rho[-\frac{7}{2}, -\frac{7}{2}],\Delta_\rho[-\frac{5}{2}, -\frac{5}{2}],\Delta_\rho[-\frac{1}{2}, -\frac{3}{2}]; \pi_{sc})$
                \item $L(\Delta_\rho[-\frac{7}{2}, -\frac{7}{2}],\Delta_\rho[-\frac{5}{2}, -\frac{5}{2}],\Delta_\rho[-\frac{3}{2}, -\frac{3}{2}],\Delta_\rho[-\frac{1}{2}, -\frac{1}{2}]; \pi_{sc})$
            \end{itemize}
            \item $(\frac{3}{2}, \frac{5}{2}, \frac{7}{2}, \frac{9}{2}):$
            \begin{itemize}
            \item $T_{I,1}^{\frac{9}{2}}(T_{I,1}^{\frac{7}{2}}(T_{I,1}^{\frac{5}{2}}(T_{I,1}^{\frac{3}{2}}(\pi_{sc}))))$
                \item $L(\Delta_\rho[-\frac{9}{2}, -\frac{9}{2}],\Delta_\rho[-\frac{7}{2}, -\frac{7}{2}],\Delta_\rho[-\frac{5}{2}, -\frac{5}{2}],\Delta_\rho[-\frac{3}{2}, -\frac{3}{2}]; \pi_{sc})$
            \end{itemize}
        \end{enumerate}
        \item $(\alpha = 2)$: 
        \begin{enumerate}
        \item $(0,0,1,2)$: 
        \begin{itemize}
        \item $T_{V,4}^{\pm}(T_{I,1}^{2}(T_{I,1}^{1}(\pi_{sc})))$
            \item $L(\Delta_\rho[-1,-1]; T_{I,1}^{2}(T_{IV,5}(\pi_{sc})))$
            \item $L(\Delta_\rho[0,-1]; T_{IV,3}(T_{I,1}^{2}(\pi_{sc})))$
            \item $L(\Delta_\rho[0,-2]; T_{IV,3}(\pi_{sc}))$
            \item $L(\Delta_\rho[-1,-2]; T_{IV,5}(\pi_{sc}))$
            \item $L(\Delta_\rho[-2,-1]; T_{IV,5}(\pi_{sc}))$
            \item $L(\Delta_\rho[-2,-2], \Delta_\rho[0,-1]; T_{IV,3}(\pi_{sc}))$ 
            \item $L(\Delta_\rho[-2,-2], \Delta_\rho[-1,-1]; T_{IV,5}(\pi_{sc}))$ 
        \end{itemize}
        \item $(0,1,1,2)$: 
        \begin{itemize}
        \item $T_{I,2}^{1}(T_{IV,3}(T_{I,1}^{2}(\pi_{sc})))$
        \item $T_{II,3}^{1}(T_{IV,3}(T_{I,1}^{2}(\pi_{sc})))$
        \item $L(\Delta_\rho[0,-1]; T_{I,1}^{1}(T_{I,1}^{2}(\pi_{sc})))$
            \item $L(\Delta_\rho[-1,-1]; T_{V,2}^{\pm}(T_{I,1}^{1}(T_{I,1}^{2}(\pi_{sc}))))$
            \item $L(\Delta_\rho[-1,-1]; T_{I,1}^{1}(T_{IV,3}(T_{I,1}^{2}(\pi_{sc}))))$
            \item $L(\Delta_\rho[1,-2]; \pi_{sc})$  
            \item $L(\Delta_\rho[-2,-2]; T_{II,3}^{1}(T_{IV,3}(\pi_{sc})))$
            \item $L(\Delta_\rho[-1,-1], \Delta_\rho[-1,-1]; T_{V,2}^{\pm}(T_{I,1}^{1}(\pi_{sc})))$
            \item $L(\Delta_\rho[-1,-1], \Delta_\rho[0,-1]; T_{I,1}^{2}(\pi_{sc}))$
            \item $L(\Delta_\rho[-1,-1], \Delta_\rho[0,-2]; \pi_{sc})$ 
            \item $L(\Delta_\rho[-2,-2], \Delta_\rho[-1,-1], T_{IV,3}(T_{I,1}^{1}(\pi_{sc})))$
            \item $L(\Delta_\rho[-2,-2], \Delta_\rho[-1,-1], \Delta_\rho[0,-1]; \pi_{sc})$
            \item $L(\Delta_\rho[-2,-2], \Delta_\rho[-1,-1], \Delta_\rho[-1,-1]; T_{IV,3}(\pi_{sc}))$
        \end{itemize}
        \item $(0,1,2,2)$: 
        \begin{itemize}
            \item $L(\Delta_\rho[-2,-2], \Delta_\rho[-1,-1]; T_{IV,3}(T_{I,1}^{2}(\pi_{sc})))$
            \item $L(\Delta_\rho[-2,-2], \Delta_\rho[0,-1]; T_{I,1}^{2}(\pi_{sc}))$
        \end{itemize}
        \item $(0,1,2,3)$: 
        \begin{itemize}
        \item $T_{V,2}^{\pm}(T_{I,1}^{1}(T_{I,1}^{3}(T_{I,1}^{2}(\pi_{sc})))$
            \item $L(\Delta_\rho[-1,-1]; T_{IV,3}(T_{I,1}^{3}(T_{I,1}^{2}(\pi_{sc}))))$
            \item $L(\Delta_\rho[0,-1]; T_{I,1}^{3}(T_{I,1}^{2}(\pi_{sc})))$
            \item $L(\Delta_\rho[-3,-3], \Delta_\rho[0,-2]; \pi_{sc})$
            \item $L(\Delta_\rho[-3,-3], \Delta_\rho[-1,-2]; T_{IV,3}(\pi_{sc}))$
            \item $L(\Delta_\rho[-3,-3], \Delta_\rho[-2,-2], \Delta_\rho[-1,-1]; T_{IV,3}(\pi_{sc}))$
            \item $L(\Delta_\rho[-3,-3], \Delta_\rho[-2,-2], \Delta_\rho[0,-1]; \pi_{sc})$
        \end{itemize}
        \item $(1,2,2,3)$: 
        \begin{itemize}
        \item $T_{I,1}^{2}(T_{I,1}^{1}(T_{I,1}^{3}(T_{I,1}^{2}(\pi_{sc}))))$
        \item $L(\Delta_\rho[-2,-2], \Delta_\rho[-1,-1]; T_{I,1}^{3}(T_{I,1}^{2}(\pi_{sc})))$  
        \item $L(\Delta_\rho[-2,-3], \Delta_\rho[-1,-2]; \pi_{sc})$
            \item $L(\Delta_\rho[-3,-3], \Delta_\rho[-2,-2], \Delta_\rho[-1,-1]; T_{I,1}^{2}(\pi_{sc}))$
        \end{itemize}
            \item $(1,2,3,4)$: 
            \begin{itemize}
            \item $T_{I,1}^{1}(T_{I,1}^{4}(T_{I,1}^{3}(T_{I,1}^{2}(\pi_{sc}))))$
                \item $L(\Delta_\rho[-1,-1]; T_{I,1}^{4}(T_{I,1}^{3}(T_{I,1}^{2}(\pi_{sc}))))$
                \item $L(\Delta_\rho[-4,-4], \Delta_\rho[-3,-3], \Delta_\rho[-1,-2]; \pi_{sc})$
                \item $L(\Delta_\rho[-4,-4], \Delta_\rho[-3,-3], \Delta_\rho[-2,-2], \Delta_\rho[-1,-1]; \pi_{sc})$
            \end{itemize}
            \item $(2,3,4,5)$: 
            \begin{itemize}
            \item $T_{I,1}^{5}(T_{I,1}^{4}(T_{I,1}^{3}(T_{I,1}^{2}(\pi_{sc}))))$
                \item $L(\Delta_\rho[-5,-5], \Delta_\rho[-4,-4], \Delta_\rho[-3,-3], \Delta_\rho[-2,-2]; \pi_{sc})$
            \end{itemize}
        \end{enumerate}
        \item $(\alpha = \frac{5}{2})$: 
        \begin{enumerate}
        \item $(\frac{1}{2}, \frac{1}{2}, \frac{3}{2}, \frac{5}{2})$: 
        \begin{itemize}
        \item $T_{III,2}^{\frac{1}{2}}(T_{I,1}^{\frac{3}{2}}(T_{I,1}^{\frac{5}{2}}(\pi_{sc})))$ 
        \item $T_{I,2}^{\frac{1}{2}}(T_{I,1}^{\frac{3}{2}}(T_{I,1}^{\frac{5}{2}}(\pi_{sc})))$
        \item $L(\Delta_\rho[\frac{1}{2}, -\frac{3}{2}]; T_{I,1}^{\frac{5}{2}}(\pi_{sc}))$
        \item $L(\Delta_\rho[\frac{1}{2}, -\frac{5}{2}]; \pi_{sc})$
        \item $L(\Delta_\rho[-\frac{1}{2}, -\frac{1}{2}]; T_{I,1}^{\frac{1}{2}}(T_{I,1}^{\frac{3}{2}}(T_{I,1}^{\frac{5}{2}}(\pi_{sc}))))$
            \item $L(\Delta_\rho[-\frac{3}{2}, -\frac{3}{2}]; T_{I,1}^{\frac{5}{2}}(T_{II,3}^{\frac{1}{2}}(\pi_{sc}))$
            \item $L(\Delta_\rho[-\frac{3}{2}, -\frac{5}{2}]; T_{II,3}^{\frac{1}{2}}(\pi_{sc}))$
            \item $L(\Delta_\rho[-\frac{5}{2}, -\frac{5}{2}],\Delta_\rho[-\frac{3}{2}, -\frac{3}{2}]; T_{II,3}^{\frac{1}{2}}(\pi_{sc}))$
            \item $L(\Delta_\rho[-\frac{1}{2}, -\frac{1}{2}], \Delta_\rho[-\frac{1}{2}, -\frac{1}{2}]; T_{I,1}^{\frac{3}{2}}(T_{I,1}^{\frac{5}{2}}(\pi_{sc})))$
            \item $L(\Delta_\rho[-\frac{1}{2}, -\frac{3}{2}], \Delta_\rho[-\frac{1}{2}, -\frac{1}{2}]; T_{I,1}^{\frac{5}{2}}(\pi_{sc}))$
            \item $L(\Delta_\rho[-\frac{1}{2}, -\frac{5}{2}], \Delta_\rho[-\frac{1}{2}, -\frac{1}{2}]; \pi_{sc})$
            \item $L(\Delta_\rho[-\frac{5}{2}, -\frac{5}{2}], \Delta_\rho[\frac{1}{2}, -\frac{3}{2}]; \pi_{sc})$
            \item $L(\Delta_\rho[-\frac{3}{2}, -\frac{5}{2}], \Delta_\rho[-\frac{1}{2}, -\frac{1}{2}], \Delta_\rho[-\frac{1}{2}, -\frac{1}{2}]; \pi_{sc})$  
            \item $L(\Delta_\rho[-\frac{3}{2}, -\frac{3}{2}], \Delta_\rho[-\frac{1}{2}, -\frac{1}{2}],\Delta_\rho[-\frac{1}{2}, -\frac{1}{2}]; T_{I,1}^{\frac{5}{2}}(\pi_{sc}))$
            \item $L(\Delta_\rho[-\frac{5}{2}, -\frac{5}{2}], \Delta_\rho[-\frac{1}{2}, -\frac{3}{2}], \Delta_\rho[-\frac{1}{2}, -\frac{1}{2}]; \pi_{sc})$
            \item $L(\Delta_\rho[-\frac{5}{2}, -\frac{5}{2}],\Delta_\rho[-\frac{3}{2}, -\frac{3}{2}],\Delta_\rho[-\frac{1}{2}, -\frac{1}{2}],\Delta_\rho[-\frac{1}{2}, -\frac{1}{2}]; \pi_{sc})$
        \end{itemize}
        \item $(\frac{1}{2}, \frac{3}{2}, \frac{3}{2}, \frac{5}{2})$: 
        \begin{itemize}
        \item $L(\Delta_\rho[-\frac{3}{2}, -\frac{3}{2}], \Delta_\rho[-\frac{1}{2}, -\frac{1}{2}]; T_{I,1}^{\frac{3}{2}}(T_{I,1}^{\frac{5}{2}}(\pi_{sc})))$
        \item $L(\Delta_\rho[-\frac{3}{2}, -\frac{5}{2}], \Delta_\rho[-\frac{1}{2}, -\frac{3}{2}]; \pi_{sc})$
        \end{itemize}
        \item $(\frac{1}{2}, \frac{3}{2}, \frac{5}{2}, \frac{5}{2})$: 
        \begin{itemize}
        \item $L(\Delta_\rho[-\frac{5}{2}, -\frac{5}{2}],\Delta_\rho[-\frac{3}{2}, -\frac{3}{2}], \Delta_\rho[-\frac{1}{2}, -\frac{1}{2}]; T_{I,1}^{\frac{5}{2}}(\pi_{sc}))$
       \item $L(\Delta_\rho[-\frac{5}{2}, -\frac{5}{2}], \Delta_\rho[-\frac{1}{2}, -\frac{3}{2}]; T_{I,1}^{\frac{5}{2}}(\pi_{sc}))$
        \end{itemize}
        \item $(\frac{1}{2}, \frac{3}{2}, \frac{5}{2}, \frac{7}{2})$: 
        \begin{itemize}
        \item $T_{I,1}^{\frac{1}{2}}(T_{I,1}^{\frac{3}{2}}(T_{I,1}^{\frac{7}{2}}(T_{I,1}^{\frac{5}{2}}(\pi_{sc}))))$
            \item $L(\Delta_\rho[-\frac{1}{2}, -\frac{1}{2}]; T_{I,1}^{\frac{3}{2}}(T_{I,1}^{\frac{7}{2}}(T_{I,1}^{\frac{5}{2}}(\pi_{sc}))))$
            \item $L(\Delta_\rho[-\frac{1}{2}, -\frac{3}{2}]; T_{I,1}^{\frac{7}{2}}(T_{I,1}^{\frac{5}{2}}(\pi_{sc})))$
            \item $L(\Delta_\rho[-\frac{3}{2}, -\frac{3}{2}], \Delta_\rho[-\frac{1}{2}, -\frac{1}{2}]; T_{I,1}^{\frac{7}{2}}(T_{I,1}^{\frac{5}{2}}(\pi_{sc})))$
            \item $L(\Delta_\rho[-\frac{7}{2}, -\frac{7}{2}], \Delta_\rho[-\frac{1}{2}, -\frac{5}{2}]; \pi_{sc})$ 
            \item $L(\Delta_\rho[-\frac{7}{2}, -\frac{7}{2}],\Delta_\rho[-\frac{5}{2}, -\frac{5}{2}], \Delta_\rho[-\frac{1}{2}, -\frac{3}{2}]; \pi_{sc})$
            \item $L(\Delta_\rho[-\frac{7}{2}, -\frac{7}{2}], \Delta_\rho[-\frac{3}{2}, -\frac{5}{2}],\Delta_\rho[-\frac{1}{2}, -\frac{1}{2}]; \pi_{sc})$ 
            \item $L(\Delta_\rho[-\frac{7}{2}, -\frac{7}{2}], \Delta_\rho[-\frac{5}{2}, -\frac{5}{2}],\Delta_\rho[-\frac{3}{2}, -\frac{3}{2}],\Delta_\rho[-\frac{1}{2}, -\frac{1}{2}]; \pi_{sc})$
        \end{itemize}
        \item $(\frac{3}{2}, \frac{5}{2}, \frac{5}{2}, \frac{7}{2})$: 
        \begin{itemize}
        \item $T_{I,1}^{\frac{5}{2}}(T_{I,1}^{\frac{3}{2}}(T_{I,1}^{\frac{7}{2}}(T_{I,1}^{\frac{5}{2}}(\pi_{sc}))))$
            \item $L(\Delta_\rho[-\frac{5}{2}, -\frac{5}{2}], \Delta_\rho[-\frac{3}{2}, -\frac{3}{2}]; T_{I,1}^{\frac{7}{2}}(T_{I,1}^{\frac{5}{2}}(\pi_{sc})))$
            \item $L(\Delta_\rho[-\frac{5}{2}, -\frac{7}{2}], \Delta_\rho[-\frac{3}{2}, -\frac{5}{2}]; \pi_{sc})$
            \item $L(\Delta_\rho[-\frac{7}{2}, -\frac{7}{2}],\Delta_\rho[-\frac{5}{2}, -\frac{5}{2}],\Delta_\rho[-\frac{3}{2}, -\frac{3}{2}]; T_{I,1}^{\frac{5}{2}}(\pi_{sc}))$
        \end{itemize}
            \item $(\frac{3}{2}, \frac{5}{2}, \frac{7}{2}, \frac{9}{2})$: 
            \begin{itemize}
            \item $T_{I,1}^{\frac{3}{2}}(T_{I,1}^{\frac{9}{2}}(T_{I,1}^{\frac{7}{2}}(T_{I,1}^{\frac{5}{2}}(\pi_{sc}))))$
                \item $L(\Delta_\rho[-\frac{3}{2}, -\frac{3}{2}]; T_{I,1}^{\frac{9}{2}}(T_{I,1}^{\frac{7}{2}}(T_{I,1}^{\frac{5}{2}}(\pi_{sc}))))$
                \item $L(\Delta_\rho[-\frac{9}{2}, -\frac{9}{2}],\Delta_\rho[-\frac{7}{2}, -\frac{7}{2}],\Delta_\rho[-\frac{3}{2}, -\frac{5}{2}]; \pi_{sc})$
                \item $L(\Delta_\rho[-\frac{9}{2}, -\frac{9}{2}],\Delta_\rho[-\frac{7}{2}, -\frac{7}{2}],\Delta_\rho[-\frac{5}{2}, -\frac{5}{2}],\Delta_\rho[-\frac{3}{2}, -\frac{3}{2}]; \pi_{sc})$
            \end{itemize}
            \item $(\frac{5}{2}, \frac{7}{2}, \frac{9}{2}, \frac{11}{2})$: 
            \begin{itemize}
                \item $T_{I,1}^{\frac{11}{2}}(T_{I,1}^{\frac{9}{2}}(T_{I,1}^{\frac{7}{2}}(T_{I,1}^{\frac{5}{2}}(\pi_{sc}))))$
                \item $L(\Delta_\rho[-\frac{11}{2}, -\frac{11}{2}], \Delta_\rho[-\frac{9}{2}, -\frac{9}{2}], \Delta_\rho[-\frac{7}{2}, -\frac{7}{2}],\Delta_\rho[-\frac{5}{2}, -\frac{5}{2}]; \pi_{sc})$
            \end{itemize}
        \end{enumerate}
        \item $(\alpha = 3)$: 
        \begin{enumerate}
        \item $(0,1,2,3)$: 
        \begin{itemize}
        \item $T_{V,2}^{\pm}(T_{I,1}^{1}(T_{I,1}^{2}(T_{I,1}^{3}(\pi_{sc}))))$
            \item $L(\Delta_\rho[-1,-1]; T_{IV,3}(T_{I,1}^{2}(T_{I,1}^{3}(\pi_{sc}))))$
            \item $L(\Delta_\rho[-1,-2]; T_{IV,3}(T_{I,1}^{2}(\pi_{sc})))$
            \item $L(\Delta_\rho[0,-2]; T_{I,1}^{3}(\pi_{sc}))$
            \item $L(\Delta_\rho[-1,-3]; T_{IV,3}(\pi_{sc}))$
            \item $L(\Delta_\rho[0,-1]; T_{I,1}^{2}(T_{I,1}^{3}(\pi_{sc})))$  
            \item $L(\Delta_\rho[0,-3]; \pi_{sc})$
            \item $L(\Delta_\rho[-2,-2], \Delta_\rho[0,-1]; T_{I,1}^{3}(\pi_{sc}))$  
            \item $L(\Delta_\rho[-2,-3], \Delta_\rho[0,-1]; \pi_{sc})$ 
            \item $L(\Delta_\rho[-3,-3], \Delta_\rho[0,-2]; \pi_{sc})$  
            \item $L(\Delta_\rho[-2,-2], \Delta_\rho[-1,-1]; T_{IV,3}(T_{I,1}^{3}(\pi_{sc})))$
            \item $L(\Delta_\rho[-2,-3],\Delta_\rho[-1,-1]; T_{IV,3}(\pi_{sc}))$  
            \item $L(\Delta_\rho[-3,-3], \Delta_\rho[-1,-2]; T_{IV,3}(\pi_{sc}))$
            \item $L(\Delta_\rho[-3,-3], \Delta_\rho[-2,-2], \Delta_\rho[-1,-1]; T_{IV,3}(\pi_{sc}))$
            \item $L(\Delta_\rho[-3,-3], \Delta_\rho[-2,-2], \Delta_\rho[0,-1]; \pi_{sc})$
            
        \end{itemize}
        \item $(1,2,2,3)$: 
        \begin{itemize}
            \item $L(\Delta_\rho[-2,-2], \Delta_\rho[-1,-1]; T_{I,1}^{2}(T_{I,1}^{3}(\pi_{sc})))$
            \item $L(\Delta_\rho[-2,-3], \Delta_\rho[-1,-2]; \pi_{sc})$
        \end{itemize}
        \item $(1,2,3,3)$: 
        \begin{itemize}
            \item $L(\Delta_\rho[-3,-3], \Delta_\rho[-2,-2], \Delta_\rho[-1,-1]; T_{I,1}^{3}(\pi_{sc}))$
            \item $L(\Delta_\rho[-3,-3], \Delta_\rho[-1,-2]; T_{I,1}^{3}(\pi_{sc}))$ 
        \end{itemize}
        \item $(1,2,3,4)$: 
        \begin{itemize}
        \item $T_{I,1}^{1}(T_{I,1}^{2}(T_{I,1}^{4}(T_{I,1}^{3}(\pi_{sc}))))$
            \item $L(\Delta_\rho[-1,-1]; T_{I,1}^{2}(T_{I,1}^{4}(T_{I,1}^{3}(\pi_{sc}))))$
            \item $L(\Delta_\rho[-1,-2]; T_{I,1}^{4}(T_{I,1}^{3}(\pi_{sc})))$
            \item $L(\Delta_\rho[-2,-2], \Delta_\rho[-1,-1]; T_{I,1}^{4}(T_{I,1}^{3}(\pi_{sc})))$
            \item $L(\Delta_\rho[-4,-4], \Delta_\rho[-1,-3]; \pi_{sc})$  
            \item $L(\Delta_\rho[-4,-4], \Delta_\rho[-2,-3],\Delta_\rho[-1,-1]; \pi_{sc})$
            \item $L(\Delta_\rho[-4,-4], \Delta_\rho[-3,-3], \Delta_\rho[-1,-2]; \pi_{sc})$
            \item $L(\Delta_\rho[-4,-4], \Delta_\rho[-3,-3], \Delta_\rho[-2,-2], \Delta_\rho[-1,-1]; \pi_{sc})$
            
        \end{itemize}
        \item $(2,3,3,4)$: 
        \begin{itemize}
        \item $T_{I,1}^{3}(T_{I,1}^{4}(T_{I,1}^{3}(T_{I,1}^{2}(\pi_{sc}))))$
            \item $L(\Delta_\rho[-3,-3], \Delta_\rho[-2,-2]; T_{I,1}^{4}(T_{I,1}^{3}(\pi_{sc})))$
            \item $L(\Delta_\rho[-3,-4], \Delta_\rho[-2,-3]; \pi_{sc})$
            \item $L(\Delta_\rho[-4,-4], \Delta_\rho[-3,-3], \Delta_\rho[-2,-2]; T_{I,1}^{3}(\pi_{sc}))$
            
        \end{itemize}
            \item $(2,3,4,5)$: 
            \begin{itemize}
            \item $T_{I,1}^{2}(T_{I,1}^{5}(T_{I,1}^{4}(T_{I,1}^{3}(\pi_{sc}))))$
                \item $L(\Delta_\rho[-2,-2]; T_{I,1}^{5}(T_{I,1}^{4}(T_{I,1}^{3}(\pi_{sc}))))$
                \item $L(\Delta_\rho[-5,-5], \Delta_\rho[-4,-4], \Delta_\rho[-2,-3]; \pi_{sc})$
                \item $L(\Delta_\rho[-5,-5], \Delta_\rho[-4,-4], \Delta_\rho[-3,-3], \Delta_\rho[-2,-2]; \pi_{sc})$
                
            \end{itemize}
            \item $(3,4,5,6)$: 
            \begin{itemize}
            \item $T_{I,1}^{6}(T_{I,1}^{5}(T_{I,1}^{4}(T_{I,1}^{3}(\pi_{sc}))))$
                \item $L(\Delta_\rho[-6,-6], \Delta_\rho[-5,-5], \Delta_\rho[-4,-4], \Delta_\rho[-3,-3]; \pi_{sc})$
                
            \end{itemize}
        \end{enumerate}
        \item $(\alpha = \frac{7}{2})$: 
        \begin{enumerate}
        \item $(\frac{1}{2}, \frac{3}{2}, \frac{5}{2}, \frac{7}{2})$: 
        \begin{itemize}
        \item $T_{I,1}^{\frac{1}{2}}(T_{I,1}^{\frac{3}{2}}(T_{I,1}^{\frac{5}{2}}(T_{I,1}^{\frac{7}{2}}(\pi_{sc}))))$
            \item $L(\Delta_\rho[-\frac{1}{2}, -\frac{1}{2}]; T_{I,1}^{\frac{3}{2}}T_{I,1}^{\frac{5}{2}}(T_{I,1}^{\frac{7}{2}}(\pi_{sc})))$
            \item $L(\Delta_\rho[-\frac{1}{2}, -\frac{3}{2}]; T_{I,1}^{\frac{5}{2}}(T_{I,1}^{\frac{7}{2}}(\pi_{sc})))$
            \item $L(\Delta_\rho[-\frac{1}{2}, -\frac{5}{2}]; T_{I,1}^{\frac{7}{2}}(\pi_{sc}))$
            \item $L(\Delta_\rho[-\frac{1}{2}, -\frac{7}{2}]; \pi_{sc})$
            \item $L(\Delta_\rho[-\frac{3}{2}, -\frac{3}{2}],\Delta_\rho[-\frac{1}{2}, -\frac{1}{2}]; T_{I,1}^{\frac{5}{2}}(T_{I,1}^{\frac{7}{2}}(\pi_{sc})))$
            \item $L(\Delta_\rho[-\frac{3}{2}, -\frac{5}{2}], \Delta_\rho[-\frac{1}{2}, -\frac{1}{2}]; T_{I,1}^{\frac{7}{2}}(\pi_{sc}))$
            \item $L(\Delta_\rho[-\frac{3}{2}, -\frac{7}{2}], \Delta_\rho[-\frac{1}{2}, -\frac{1}{2}]; \pi_{sc})$
            \item $L(\Delta_\rho[-\frac{5}{2}, -\frac{5}{2}], \Delta_\rho[-\frac{1}{2}, -\frac{3}{2}]; T_{I,1}^{\frac{7}{2}}(\pi_{sc}))$
            \item $L(\Delta_\rho[-\frac{5}{2}, -\frac{7}{2}], \Delta_\rho[-\frac{1}{2}, -\frac{3}{2}]; \pi_{sc})$
            \item $L(\Delta_\rho[-\frac{7}{2}, -\frac{7}{2}],\Delta_\rho[-\frac{1}{2}, -\frac{5}{2}]; \pi_{sc})$
            \item $L(\Delta_\rho[-\frac{5}{2}, -\frac{5}{2}], \Delta_\rho[-\frac{3}{2}, -\frac{3}{2}],\Delta_\rho[-\frac{1}{2}, -\frac{1}{2}]; T_{I,1}^{\frac{7}{2}}(\pi_{sc}))$
            \item $L(\Delta_\rho[-\frac{5}{2}, -\frac{7}{2}], \Delta_\rho[-\frac{3}{2}, -\frac{3}{2}],\Delta_\rho[-\frac{1}{2}, -\frac{1}{2}]; \pi_{sc})$
            \item $L(\Delta_\rho[-\frac{7}{2}, -\frac{7}{2}], \Delta_\rho[-\frac{3}{2}, -\frac{5}{2}],\Delta_\rho[-\frac{1}{2}, -\frac{1}{2}]; \pi_{sc})$
            \item $L(\Delta_\rho[-\frac{7}{2}, -\frac{7}{2}],\Delta_\rho[-\frac{5}{2}, -\frac{5}{2}],\Delta_\rho[-\frac{1}{2}, -\frac{3}{2}]; \pi_{sc})$
            \item $L(\Delta_\rho[-\frac{7}{2}, -\frac{7}{2}],\Delta_\rho[-\frac{5}{2}, -\frac{5}{2}],\Delta_\rho[-\frac{3}{2}, -\frac{3}{2}],\Delta_\rho[-\frac{1}{2}, -\frac{1}{2}]; \pi_{sc})$
            
        \end{itemize}
        \item $(\frac{3}{2}, \frac{5}{2}, \frac{5}{2}, \frac{7}{2}):$
        \begin{itemize}
            \item $L(\Delta_\rho[-\frac{5}{2}, -\frac{5}{2}], \Delta_\rho[-\frac{3}{2}, -\frac{3}{2}]; T_{I,1}^{\frac{5}{2}}(T_{I,1}^{\frac{7}{2}}(\pi_{sc})))$
            \item $L(\Delta_\rho[-\frac{5}{2}, -\frac{7}{2}], \Delta_\rho[-\frac{3}{2}, -\frac{5}{2}]; \pi_{sc})$
        \end{itemize}
        \item $(\frac{3}{2}, \frac{5}{2}, \frac{7}{2}, \frac{7}{2})$: 
        \begin{itemize}
        \item $L(\Delta_\rho[-\frac{7}{2}, -\frac{7}{2}], \Delta_\rho[-\frac{3}{2}, -\frac{5}{2}]; T_{I,1}^{\frac{7}{2}}(\pi_{sc}))$  
            \item $L(\Delta_\rho[-\frac{7}{2}, -\frac{7}{2}], \Delta_\rho[-\frac{5}{2}, -\frac{5}{2}], \Delta_\rho[-\frac{3}{2}, -\frac{3}{2}]; T_{I,1}^{\frac{7}{2}}(\pi_{sc}))$
        \end{itemize}
        \item $(\frac{3}{2}, \frac{5}{2}, \frac{7}{2}, \frac{9}{2})$: 
        \begin{itemize}
            \item $T_{I,1}^{\frac{3}{2}}(T_{I,1}^{\frac{5}{2}}(T_{I,1}^{\frac{9}{2}}(T_{I,1}^{\frac{7}{2}}(\pi_{sc}))))$
            \item $L(\Delta_\rho[-\frac{3}{2}, -\frac{3}{2}]; T_{I,1}^{\frac{5}{2}}(T_{I,1}^{\frac{9}{2}}(T_{I,1}^{\frac{7}{2}}(\pi_{sc}))))$
            \item $L(\Delta_\rho[-\frac{3}{2}, -\frac{5}{2}]; T_{I,1}^{\frac{9}{2}}(T_{I,1}^{\frac{7}{2}}(\pi_{sc})))$
            \item $L(\Delta_\rho[-\frac{5}{2}, -\frac{5}{2}], \Delta_\rho[-\frac{3}{2}, -\frac{3}{2}]; T_{I,1}^{\frac{9}{2}}(T_{I,1}^{\frac{7}{2}}(\pi_{sc})))$
            \item $L(\Delta_\rho[-\frac{9}{2}, -\frac{9}{2}], \Delta_\rho[-\frac{3}{2}, -\frac{7}{2}]; \pi_{sc})$
            \item $L(\Delta_\rho[-\frac{9}{2}, -\frac{9}{2}], \Delta_\rho[-\frac{7}{2}, -\frac{7}{2}],\Delta_\rho[-\frac{3}{2}, -\frac{5}{2}]; \pi_{sc})$
            \item $L(\Delta_\rho[-\frac{9}{2}, -\frac{9}{2}], \Delta_\rho[-\frac{5}{2}, -\frac{7}{2}], \Delta_\rho[-\frac{3}{2}, -\frac{3}{2}]; \pi_{sc})$ 
            \item $L(\Delta_\rho[-\frac{9}{2}, -\frac{9}{2}], \Delta_\rho[-\frac{7}{2}, -\frac{7}{2}],\Delta_\rho[-\frac{5}{2}, -\frac{5}{2}],\Delta_\rho[-\frac{3}{2}, -\frac{3}{2}]; \pi_{sc})$
            
        \end{itemize}
        \item $(\frac{5}{2}, \frac{7}{2}, \frac{7}{2}, \frac{9}{2})$: 
        \begin{itemize}
        \item $T_{I,1}^{\frac{7}{2}}(T_{I,1}^{\frac{5}{2}}(T_{I,1}^{\frac{9}{2}}(T_{I,1}^{\frac{7}{2}}(\pi_{sc}))))$. 
            \item $L(\Delta_\rho[-\frac{7}{2}, -\frac{7}{2}], \Delta_\rho[-\frac{5}{2}, -\frac{5}{2}]; T_{I,1}^{\frac{9}{2}}(T_{I,1}^{\frac{7}{2}}(\pi_{sc})))$
            \item $L(\Delta_\rho[-\frac{7}{2}, -\frac{9}{2}], \Delta_\rho[-\frac{5}{2}, -\frac{7}{2}]; \pi_{sc})$
            \item $L(\Delta_\rho[-\frac{9}{2}, -\frac{9}{2}], \Delta_\rho[-\frac{7}{2}, -\frac{7}{2}],\Delta_\rho[-\frac{5}{2}, -\frac{5}{2}]; T_{I,1}^{\frac{7}{2}}(\pi_{sc}))$
            
        \end{itemize}
            \item $(\frac{5}{2}, \frac{7}{2}, \frac{9}{2}, \frac{11}{2})$: 
            \begin{itemize}
            \item $T_{I,1}^{\frac{5}{2}}(T_{I,1}^{\frac{11}{2}}(T_{I,1}^{\frac{9}{2}}(T_{I,1}^{\frac{7}{2}}(\pi_{sc}))))$
                \item $L(\Delta_\rho[-\frac{5}{2}, -\frac{5}{2}]; T_{I,1}^{\frac{11}{2}}(T_{I,1}^{\frac{9}{2}}(T_{I,1}^{\frac{7}{2}}(\pi_{sc}))))$
                \item $L(\Delta_\rho[-\frac{11}{2}, -\frac{11}{2}], \Delta_\rho[-\frac{9}{2}, -\frac{9}{2}], \Delta_\rho[-\frac{5}{2}, -\frac{7}{2}]; \pi_{sc})$
                \item $L(\Delta_\rho[-\frac{11}{2}, -\frac{11}{2}],\Delta_\rho[-\frac{9}{2}, -\frac{9}{2}],\Delta_\rho[-\frac{7}{2}, -\frac{7}{2}],\Delta_\rho[-\frac{5}{2}, -\frac{5}{2}]; \pi_{sc})$
                
            \end{itemize}
            \item $(\frac{7}{2}, \frac{9}{2}, \frac{11}{2}, \frac{13}{2})$: 
            \begin{itemize}
            \item $T_{I,1}^{\frac{7}{2}}(T_{I,1}^{\frac{9}{2}}(T_{I,1}^{\frac{11}{2}}(T_{I,1}^{\frac{13}{2}}(\pi_{sc}))))$
                \item $L(\Delta_\rho[-\frac{13}{2}, -\frac{13}{2}], \Delta_\rho[-\frac{11}{2}, -\frac{11}{2}], \Delta_\rho[-\frac{9}{2}, -\frac{9}{2}],\Delta_\rho[-\frac{7}{2}, -\frac{7}{2}]; \pi_{sc})$
            \end{itemize}
        \end{enumerate}
        \item $(\alpha \geq 4)$: 
        \begin{enumerate}
        \item $(\alpha -3, \alpha -2, \alpha -1, \alpha)$: 
        \begin{itemize}
        \item $T_{I,1}^{\alpha-3}(T_{I,1}^{\alpha-2}(T_{I,1}^{\alpha-1}(T_{I,1}^{\alpha}(\pi_{sc}))))$
        \item $L(\Delta_\rho[-\alpha+3, -\alpha+2]; T_{I,1}^{\alpha-1}(T_{I,1}^{\alpha}(\pi_{sc})))$  
            \item $L(\Delta_\rho[-\alpha +3, -\alpha +3]; T_{I,1}^{\alpha -2}(T_{I,1}^{\alpha -1}(T_{I,1}^{\alpha}(\pi_{sc}))))$
            \item $L(\Delta_\rho[-\alpha +3, -\alpha +1]; T_{I,1}^{\alpha}(\pi_{sc}))$
            \item $L(\Delta_\rho[-\alpha +3, -\alpha]; \pi_{sc})$
            \item $L(\Delta_\rho[-\alpha +2, -\alpha +2], \Delta_\rho[-\alpha +3, -\alpha +3]; T_{I,1}^{\alpha -1}(T_{I,1}^{\alpha}(\pi_{sc})))$
            \item $L(\Delta_\rho[-\alpha, -\alpha], \Delta_\rho[-\alpha+3, -\alpha+1]; \pi_{sc})$+ 
            \item $L(\Delta_\rho[-\alpha +2, -\alpha +1], \Delta_\rho[-\alpha +3, -\alpha +3]; T_{I,1}^{\alpha}(\pi_{sc}))$
            \item $L(\Delta_\rho[-\alpha +1, -\alpha +1], \Delta_\rho[-\alpha +3, -\alpha +2]; T_{I,1}^{\alpha}(\pi_{sc}))$
            \item $L(\Delta_\rho[-\alpha +1, -\alpha], \Delta_\rho[-\alpha +3, -\alpha +2]; \pi_{sc})$
            \item $L(\Delta_\rho[-\alpha +2, -\alpha], \Delta_\rho[-\alpha +3, -\alpha +3]; \pi_{sc})$
            \item $L(\Delta_\rho[-\alpha +1, -\alpha +1], \Delta_\rho[-\alpha +2, -\alpha +2], \Delta_\rho[-\alpha +3, -\alpha +3]; T_{I,1}^{\alpha}(\pi_{sc}))$
            \item $L(\Delta_\rho[-\alpha +1, -\alpha], \Delta_\rho[-\alpha +2, -\alpha +2], \Delta_\rho[-\alpha +3, -\alpha +3]; \pi_{sc})$
            \item $L(\Delta_\rho[-\alpha, -\alpha], \Delta_\rho[-\alpha +2, -\alpha +1], \Delta_\rho[-\alpha +3, -\alpha +3]; \pi_{sc})$
            \item $L(\Delta_\rho[-\alpha, -\alpha], \Delta_\rho[-\alpha +1, -\alpha +1], \Delta_\rho[-\alpha +3, -\alpha +2]; \pi_{sc})$
            \item $L(\Delta_\rho[-\alpha, -\alpha], \Delta_\rho[-\alpha +1, -\alpha +1], \Delta_\rho[-\alpha +2, -\alpha +2], \Delta_\rho[-\alpha +3, -\alpha +3]; \pi_{sc})$
        \end{itemize}
        \item $(\alpha -2, \alpha -1, \alpha -1, \alpha)$: 
        \begin{itemize}
            \item $L(\Delta_\rho[-\alpha +1, -\alpha +1], \Delta_\rho[-\alpha +2, -\alpha +2]; T_{I,1}^{\alpha -1}(T_{I,1}^{\alpha}(\pi_{sc})))$
            \item $L(\Delta_\rho[-\alpha +1, -\alpha], \Delta_\rho[-\alpha +2, -\alpha +1]; \pi_{sc})$
        \end{itemize}
        \item $(\alpha -2, \alpha -1, \alpha, \alpha)$: 
        \begin{itemize}
        \item $L(\Delta_\rho[-\alpha, -\alpha]; \Delta_\rho[-\alpha+2, -\alpha+1]; T_{I,1}^{\alpha}(\pi_{sc}))$
            \item $L(\Delta_\rho[-\alpha, -\alpha], \Delta_\rho[-\alpha +1, -\alpha +1], \Delta_\rho[-\alpha +2, -\alpha +2]; T_{I,1}^{\alpha}(\pi_{sc}))$
        \end{itemize}
        \item $(\alpha -2, \alpha-1, \alpha , \alpha +1)$: 
        \begin{itemize}
        \item $T_{I,1}^{\alpha-2}(T_{I,1}^{\alpha-1}(T_{I,1}^{\alpha+1}(T_{I,1}^{\alpha}(\pi_{sc}))))$ 
            \item $L(\Delta_\rho[-\alpha +2, -\alpha +2]; T_{I,1}^{\alpha -1}(T_{I,1}^{\alpha +1}(T_{I,1}^{\alpha}(\pi_{sc}))))$
            \item $L(\Delta_\rho[-\alpha +2, -\alpha +1]; T_{I,1}^{\alpha +1}(T_{I,1}^{\alpha}(\pi_{sc})))$
            \item $L(\Delta_\rho[-\alpha +1, -\alpha +1],\Delta_\rho[-\alpha +2, -\alpha +2]; T_{I,1}^{\alpha +1}(T_{I,1}^{\alpha}(\pi_{sc})))$
            \item $L(\Delta_\rho[-\alpha-1, -\alpha-1], \Delta_\rho[-\alpha+2, -\alpha]; \pi_{sc})$
            \item $L(\Delta_\rho[-\alpha-1, -\alpha-1], \Delta_\rho[-\alpha+1, -\alpha], \Delta_\rho[-\alpha+2, -\alpha+2]; \pi_{sc})$ 
            \item $L(\Delta_\rho[-\alpha -1, -\alpha -1], \Delta_\rho[-\alpha, -\alpha], \Delta_\rho[-\alpha +2, -\alpha +1]; \pi_{sc})$
            \item $L(\Delta_\rho[-\alpha -1, -\alpha -1], \Delta_\rho[-\alpha, -\alpha], \Delta_\rho[-\alpha +1, -\alpha +1], \Delta_\rho[-\alpha +2, -\alpha +2]; \pi_{sc})$
        \end{itemize}
        \item $(\alpha -1, \alpha, \alpha, \alpha +1)$: 
        \begin{itemize}
        \item $T_{I,1}^{\alpha}(T_{I,1}^{\alpha-1}(T_{I,1}^{\alpha+1}(T_{I,1}^{\alpha}(\pi_{sc}))))$
            \item $L(\Delta_\rho[-\alpha, -\alpha], \Delta_\rho[-\alpha +1, -\alpha +1]; T_{I,1}^{\alpha +1}(T_{I,1}^{\alpha}(\pi_{sc})))$
            \item $L(\Delta_\rho[-\alpha-1, -\alpha-1], \Delta_\rho[-\alpha+2, -\alpha]; \pi_{sc})$
            \item $L(\Delta_\rho[-\alpha -1, -\alpha -1], \Delta_\rho[-\alpha, -\alpha], \Delta\rho[-\alpha +1, -\alpha +1]; \pi_{sc})$
        \end{itemize}
            \item $(\alpha -1, \alpha, \alpha +1, \alpha +2)$: 
            \begin{itemize}
            \item $T_{I,1}^{\alpha-1}(T_{I,1}^{\alpha+2}(T_{I,1}^{\alpha+1}(T_{I,1}^{\alpha}(\pi_{sc}))))$
                \item $L(\Delta_\rho[-\alpha +1, -\alpha +1]; T_{I,1}^{\alpha +2}(T_{I,1}^{\alpha +1}(T_{I,1}^{\alpha}(\pi_{sc}))))$
                \item $L(\Delta_\rho[-\alpha -2, -\alpha -2], \Delta_\rho[-\alpha -1, -\alpha -1], \Delta_\rho[-\alpha +1, -\alpha]; \pi_{sc})$
                \item $L(\Delta_\rho[-\alpha -2, -\alpha -2], \Delta_\rho[-\alpha -1, -\alpha -1], \Delta_\rho[-\alpha, -\alpha], \Delta_\rho[-\alpha +1, -\alpha +1]; \pi_{sc})$
            \end{itemize}
            \item $(\alpha, \alpha +1, \alpha +2, \alpha +3): $
            \begin{itemize}
            \item $T_{I,1}^{\alpha+3}(T_{I,1}^{\alpha+2}(T_{I,1}^{\alpha+1}(T_{I,1}^{\alpha}(\pi_{sc}))))$ 
                \item $L(\Delta_\rho[-\alpha -3, -\alpha -3], \Delta_\rho[-\alpha -2, -\alpha -2], \Delta_\rho[-\alpha -1, -\alpha -1], \Delta_\rho[-\alpha, -\alpha]; \pi_{sc})$
            \end{itemize}
        \end{enumerate}
    \end{enumerate}
\end{prop}
\begin{proof} This follows directly from Propositions \ref{nontempA1} to \ref{nontemp4}. 
\end{proof}

\section{\\List of representations of corank 4 that are of critical type but not of Arthur type, sorted by cuspidal support} \label{nonartlist}

In this appendix, we provide the opposite list as the one given in the previous appendix, that is, the list of all corank 4 representations that are of critical type but not of Arthur type (see Proposition \ref{subquotientlist} below). Combining the two lists given in Propositions \ref{crnk4Artcritlist} and \ref{subquotientlist} together enables us to classify the unitarity of all critical points, a crucial step in the construction of the corank 4 unitary dual. 

As before, we sort our list of representations as subquotients of $\Pi_{x_1, x_2, x_3, x_4}:= \rho\lvert \cdot \rvert^{x_1} \times \rho\lvert \cdot \rvert^{x_2} \times \rho\lvert \cdot \rvert^{x_3} \times \rho\lvert \cdot \rvert^{x_4} \rtimes \pi_{sc}$. By Theorem \ref{unitiffArthur}, this list gives exactly the set of representations that are of critical type but not unitarizable.  Without loss of generality, we may assume $ 0\leq x_1 \leq x_2 \leq x_3 \leq x_4$. 

\begin{prop}\label{subquotientlist}
    Let $\Pi_{x_1, x_2, x_3, x_4} = \rho\lvert \cdot \rvert^{x_1} \times \rho\lvert \cdot \rvert^{x_2} \times \rho\lvert \cdot \rvert^{x_3} \times \rho\lvert \cdot \rvert^{x_4} \rtimes \pi_{sc}$ for $0 \leq x_1 \leq x_2 \leq x_3 \leq x_4$. Then for any critical point $(x_1, x_2, x_3, x_4)$, the following list contains all its irreducible subquotients that are not of Arthur type. 
    \begin{enumerate} 
        \item $(\alpha = 0)$:
        \begin{enumerate}
        \item $(0,0,1,2)$:
            \begin{itemize}
                \item $L(\Delta_{\rho}[-1,-2]; T_{V,4}^{\pm}(\pi_{sc}))$
                \item$L(\Delta_\rho[-2,-2]; T_{I,1}^{1}(T_{V,4}^{\pm}(\pi_{sc}))))$
                \item $L(\Delta_\rho[0,-2]; T_{V,2}^{\pm}(\pi_{sc}))$
                \item $L(\Delta_\rho[-2,-2], \Delta_\rho[0,-1]; T_{V,2}^{\pm}(\pi_{sc}))$
            \end{itemize}
            \item $(0,1,1,1)$: 
            \begin{itemize}
                \item$L(\Delta_\rho[-1,-1]; T_{I,2}^{1}(T_{V,2}^{\pm}(\pi_{sc})))$
                \item$L(\Delta_\rho[-1,-1], \Delta_\rho[-1,-1]; T_{I,1}^{1}(T_{V,2}^{\pm}(\pi_{sc})))$
                \item $L(\Delta_\rho[-1,-1], \Delta_\rho[-1,-1], \Delta_\rho[0,-1]; \pi_{sc})$
                \item $L(\Delta_\rho[-1,-1], \Delta_\rho[-1,-1], \Delta_\rho[-1,-1]; T_{V,2}^{\pm}(\pi_{sc}))$
            
            \end{itemize}
            \item $(0,1,1,2)$:  
            \begin{itemize}
                \item $L(\Delta_\rho[-1, -2]; T_{V,2}^{\pm}(T_{I,1}^{1}(\pi_{sc}))$
                \item$L(\Delta_\rho[-2,-2]; T_{I,2}^{1}(T_{V,2}^{\pm}(\pi_{sc}))))$
                \item $L(\Delta_\rho[-1,-1], \Delta_\rho[0,-2]; \pi_{sc})$  
                \item $L(\Delta_\rho[-1,-2], \Delta_\rho[-1,-1]; T_{V,2}^{\pm}(\pi_{sc}))$ 
                \item $L(\Delta_\rho[-2,-2], \Delta_\rho[-1,-1], \Delta_\rho[0,-1]; \pi_{sc})$
            \end{itemize}
            \item $(0,1,2,2)$: 
            \begin{itemize}
            \item $L(\Delta_\rho[-2,-2], \Delta_\rho[0,-2]; \pi_{sc})$ 
                \item $L(\Delta_\rho[-2,-2]; T_{I,1}^{2}(T_{I,1}^{1}(T_{V,2}^{\pm}(\pi_{sc}))))$
                \item $L(\Delta_\rho[-2,-2], \Delta_\rho[-1,-2]; T_{V,2}^{\pm}(\pi_{sc}))$ 
                \item $L(\Delta_\rho[-2,-2], \Delta_\rho[-2,-2]; T_{I,1}(T_{V,2}^{\pm}(\pi_{sc})))$
                \item $L(\Delta_\rho[-2,-2], \Delta_\rho[-2,-2], \Delta_\rho[-1,-1]; T_{V,2}^{\pm}(\pi_{sc}))$
                \item $L(\Delta_\rho[-2,-2], \Delta_\rho[-2,-2], \Delta_\rho[0,-1]; \pi_{sc})$
            \end{itemize}
            \item $(0,1,2,3)$: 
            \begin{itemize}
                \item$L(\Delta_\rho[0,-3]; \pi_{sc})$
                \item$L(\Delta_\rho[-1,-3]; T_{V,2}^{\pm}(\pi_{sc}))$
                \item $L(\Delta_\rho[-2,-3]; T_{I,1}^{1}(T_{V,2}^{\pm}(\pi_{sc})))$  
                \item $L(\Delta_\rho[-3, -3]; T_{I,1}^{2}(T_{I,1}^{1}(T_{V,2}^{\pm}(\pi_{sc}))))$
                \item $L(\Delta_\rho[-2,-3], \Delta_\rho[0,-1]; \pi_{sc})$ 
                \item $L(\Delta_\rho[-2,-3], \Delta_\rho[-1,-1]; T_{V,2}^{\pm}(\pi_{sc}))$  
                \item $L(\Delta_\rho[-3,-3], \Delta_\rho[0,-2]; \pi_{sc})$ 
                \item $L(\Delta_\rho[-3,-3], \Delta_\rho[-1,-2]; T_{V,2}^{\pm}(\pi_{sc}))$ 
                \item $L(\Delta_\rho[-3,-3], \Delta_\rho[-2,-2]; T_{I,1}^{1}(T_{V,2}^{\pm}(\pi_{sc})))$
                \item $L(\Delta_\rho[-3,-3], \Delta_\rho[-2,-2], \Delta_\rho[0,-1]; \pi_{sc})$
            \end{itemize}
        \end{enumerate}
        \item $(\alpha = \frac{1}{2})$ : 
        \begin{enumerate}
        \item $(\frac{1}{2}, \frac{1}{2}, \frac{1}{2}, \frac{3}{2})$:
        \begin{itemize}
        \item $L(\Delta_\rho[-\frac{3}{2}, -\frac{3}{2}]; T_{I,3}^{\frac{1}{2}}(\pi_{sc}))$ 
            \item $L(\Delta_\rho[-\frac{1}{2}, -\frac{3}{2}],\Delta_\rho[-\frac{1}{2}, -\frac{1}{2}],\Delta_\rho[-\frac{1}{2}, -\frac{1}{2}];  \pi_{sc})$
        \end{itemize}
        \item $(\frac{1}{2}, \frac{1}{2}, \frac{3}{2}, \frac{3}{2})$:
        \begin{itemize}
        \item $L(\Delta_\rho[-\frac{1}{2}, -\frac{3}{2}]; T_{I,1}^{\frac{3}{2}}(T_{I,1}^{\frac{1}{2}}(\pi_{sc})))$  
        \item $L(\Delta_\rho[-\frac{3}{2}, -\frac{3}{2}]; T_{I,1}^{\frac{1}{2}}(T_{I,1}^{\frac{3}{2}}(T_{I,1}^{\frac{1}{2}}(\pi_{sc}))))$ 
        \item $L(\Delta_\rho[-\frac{1}{2}, -\frac{3}{2}], \Delta_\rho[-\frac{1}{2}, -\frac{3}{2}]; \pi_{sc})$ 
            \item $L(\Delta_\rho[-\frac{3}{2}, -\frac{3}{2}], \Delta_\rho[-\frac{1}{2}, -\frac{3}{2}]; T_{I,1}^{\frac{1}{2}}(\pi_{sc}))$
            \item $L(\Delta_\rho[-\frac{3}{2}, -\frac{3}{2}], \Delta_\rho[\frac{1}{2}, -\frac{3}{2}]; \pi_{sc})$
            \item $L(\Delta_\rho[-\frac{3}{2}, -\frac{3}{2}], \Delta_\rho[-\frac{3}{2}, -\frac{3}{2}]; T_{I,2}^{\frac{1}{2}}(\pi_{sc}))$
            \item $L(\Delta_\rho[-\frac{3}{2}, -\frac{3}{2}], \Delta_\rho[-\frac{3}{2}, -\frac{3}{2}]; T_{III,2}^{\frac{1}{2}}(\pi_{sc}))$
            \item $L(\Delta_\rho[-\frac{3}{2}, -\frac{3}{2}], \Delta_\rho[-\frac{1}{2}, -\frac{3}{2}], \Delta_\rho[-\frac{1}{2}, -\frac{1}{2}]; \pi_{sc})$
        \end{itemize}
        \item $(\frac{1}{2}, \frac{1}{2}, \frac{3}{2}, \frac{5}{2})$: 
        \begin{itemize}
        \item $L(\Delta_\rho[-\frac{1}{2},-\frac{5}{2}]; T_{I,1}^{\frac{1}{2}}(\pi_{sc}))$
        \item $L(\Delta_\rho[-\frac{3}{2}, -\frac{5}{2}]; T_{I,2}^{\frac{1}{2}}(\pi_{sc}))$
            \item $L(\Delta_\rho[-\frac{5}{2}, -\frac{5}{2}]; T_{I,1}^{\frac{3}{2}}(T_{I,2}^{\frac{1}{2}}(\pi_{sc})))$
            \item $L(\Delta_\rho[-\frac{5}{2}, -\frac{5}{2}]; T_{I,1}^{\frac{3}{2}}(T_{III,2}^{\frac{1}{2}}(\pi_{sc}))))$ 
            \item $L(\Delta_\rho[-\frac{1}{2}, -\frac{5}{2}],\Delta_\rho[-\frac{1}{2}, -\frac{1}{2}]; \pi_{sc})$
            \item $L(\Delta_\rho[\frac{1}{2}, -\frac{5}{2}]; \pi_{sc}))$
             \item $L(\Delta_\rho[-\frac{3}{2}, -\frac{5}{2}], \Delta_\rho[-\frac{1}{2}, -\frac{1}{2}]; T_{I,1}^{\frac{1}{2}}(\pi_{sc}))$
             \item $L(\Delta_\rho[-\frac{5}{2}, -\frac{5}{2}], \Delta_\rho[-\frac{1}{2}, -\frac{1}{2}]; T_{I,1}^{\frac{3}{2}}(T_{I,1}^{\frac{1}{2}}(\pi_{sc})))$  
            \item $L(\Delta_\rho[-\frac{5}{2}, -\frac{5}{2}], \Delta_\rho[-\frac{3}{2}, -\frac{3}{2}]; T_{I,2}^{\frac{1}{2}}(\pi_{sc}))$
            \item $L(\Delta_\rho[-\frac{5}{2}, -\frac{5}{2}], \Delta_\rho[-\frac{3}{2}, -\frac{3}{2}]; T_{III,2}^{\frac{1}{2}}(\pi_{sc})))$         
             \item $L(\Delta_\rho[-\frac{5}{2}, -\frac{5}{2}], \Delta_\rho[-\frac{1}{2}, -\frac{3}{2}]; T_{I,1}^{\frac{1}{2}}(\pi_{sc}))$ with $\epsilon_{sc}(\rho \otimes S_2) = 1$
             \item $L(\Delta_\rho[-\frac{5}{2}, -\frac{5}{2}], \Delta_\rho[\frac{1}{2}, -\frac{3}{2}]; \pi_{sc})$
             \item $L(\Delta_\rho[-\frac{5}{2}, -\frac{5}{2}], \Delta_\rho[-\frac{1}{2}, -\frac{3}{2}], \Delta_\rho[-\frac{1}{2}, -\frac{1}{2}]; \pi_{sc})$
             \item $L(\Delta_\rho[-\frac{5}{2}, -\frac{5}{2}], \Delta_\rho[-\frac{3}{2}, -\frac{3}{2}], \Delta_\rho[-\frac{1}{2}, -\frac{1}{2}]; T_{I,1}^{\frac{1}{2}}(\pi_{sc}))$
        \end{itemize}
        \item $(\frac{1}{2}, \frac{3}{2}, \frac{3}{2}, \frac{3}{2})$: 
        \begin{itemize}
            \item $L(\Delta_\rho[-\frac{3}{2},-\frac{3}{2}], \Delta_\rho[-\frac{3}{2}, -\frac{3}{2}]; T_{I,1}^{\frac{3}{2}}(T_{I,1}^{\frac{1}{2}}(\pi_{sc})))$
            \item $L(\Delta_\rho[-\frac{3}{2}, -\frac{3}{2}],\Delta_\rho[-\frac{3}{2}, -\frac{3}{2}],\Delta_\rho[-\frac{1}{2}, -\frac{3}{2}]; \pi_{sc}))$
            \item $L(\Delta_\rho[-\frac{3}{2}, -\frac{3}{2}],\Delta_\rho[-\frac{3}{2}, -\frac{3}{2}],\Delta_\rho[-\frac{3}{2}, -\frac{3}{2}]; T_{I,1}^{\frac{1}{2}}(\pi_{sc})) $
            \item  $L(\Delta_\rho[-\frac{3}{2}, -\frac{3}{2}], \Delta_\rho[-\frac{3}{2}, -\frac{3}{2}],\Delta_\rho[-\frac{3}{2}, -\frac{3}{2}], \Delta_\rho[-\frac{1}{2}, -\frac{1}{2}]; \pi_{sc})$
            
        \end{itemize}
            \item $(\frac{1}{2}, \frac{3}{2}, \frac{3}{2}, \frac{5}{2})$: 
            \begin{itemize}
                \item $L(\Delta_\rho[-\frac{3}{2}, -\frac{3}{2}]; T_{I,1}^{\frac{5}{2}}(T_{I,1}^{\frac{3}{2}}(T_{I,1}^{\frac{1}{2}}(\pi_{sc}))))$
                \item $L(\Delta_\rho[-\frac{3}{2}, -\frac{5}{2}];T_{I,1}^{\frac{3}{2}}(T_{I,1}^{\frac{1}{2}}(\pi_{sc})))$ 
                \item $L(\Delta_\rho[-\frac{3}{2},-\frac{3}{2}], \Delta_\rho[-\frac{1}{2}, -\frac{5}{2}];\pi_{sc})$  
                \item $L(\Delta_\rho[-\frac{3}{2}, -\frac{5}{2}], \Delta_\rho[-\frac{1}{2}, -\frac{3}{2}]; \pi_{sc})$
                \item $L(\Delta_\rho[-\frac{3}{2}, -\frac{5}{2}], \Delta_\rho[-\frac{3}{2}, -\frac{3}{2}]; T_{I,1}^{\frac{1}{2}}(\pi_{sc}))$
                \item $L(\Delta_\rho[-\frac{5}{2},-\frac{5}{2}], \Delta_\rho[-\frac{3}{2}, -\frac{3}{2}];T_{I,1}^{\frac{3}{2}}(T_{I,1}^{\frac{1}{2}}(\pi_{sc})))$
                \item $L(\Delta_\rho[-\frac{3}{2}, -\frac{5}{2}], \Delta_\rho[-\frac{3}{2}, -\frac{3}{2}], \Delta_\rho[-\frac{1}{2}, -\frac{1}{2}];\pi_{sc})$  
                \item $L(\Delta_\rho[-\frac{5}{2}, -\frac{5}{2}], \Delta_\rho[-\frac{3}{2}, -\frac{3}{2}], \Delta_\rho[-\frac{1}{2}, -\frac{3}{2}]; \pi_{sc})$ 
                \item $L(\Delta_\rho[-\frac{5}{2}, -\frac{5}{2}], \Delta_\rho[-\frac{3}{2}, -\frac{3}{2}], \Delta_\rho[-\frac{3}{2}, -\frac{3}{2}]; T_{I,1}^{\frac{1}{2}}(\pi_{sc}))$  
                \item $L(\Delta_\rho[-\frac{5}{2}, -\frac{5}{2}], \Delta_\rho[-\frac{3}{2}, -\frac{3}{2}], \Delta_\rho[-\frac{3}{2}, -\frac{3}{2}], \Delta_\rho[-\frac{1}{2}, -\frac{1}{2}]; \pi_{sc})$
            \end{itemize}
            \item $(\frac{1}{2}, \frac{3}{2}, \frac{5}{2}, \frac{5}{2})$: 
            \begin{itemize}
                \item $L(\Delta_\rho[-\frac{5}{2}, -\frac{5}{2}]; T_{I,1}^{\frac{5}{2}}(T_{I,1}^{\frac{3}{2}}(T_{I,1}^{\frac{1}{2}}(\pi_{sc}))))$
                \item $L(\Delta_\rho[-\frac{5}{2}, -\frac{5}{2}], \Delta_\rho[-\frac{1}{2}, -\frac{5}{2}]; \pi_{sc})$  
                \item $L(\Delta_\rho[-\frac{5}{2}, -\frac{5}{2}], \Delta_\rho[-\frac{3}{2}, -\frac{5}{2}]; T_{I,1}^{\frac{1}{2}}(\pi_{sc}))$
                \item $L(\Delta_\rho[-\frac{5}{2}, -\frac{5}{2}], \Delta_\rho[-\frac{5}{2}, -\frac{5}{2}]; T_{I,1}^{\frac{3}{2}}(T_{I,1}^{\frac{1}{2}}(\pi_{sc})))$
                \item $L(\Delta_\rho[-\frac{5}{2}, -\frac{5}{2}], \Delta_\rho[-\frac{3}{2}, -\frac{5}{2}], \Delta_\rho[-\frac{1}{2}, -\frac{1}{2}]; \pi_{sc})$
                \item $L(\Delta_\rho[-\frac{5}{2}, -\frac{5}{2}], \Delta_\rho[-\frac{5}{2}, -\frac{5}{2}], \Delta_\rho[-\frac{3}{2}, -\frac{3}{2}], T_{I,1}^{\frac{1}{2}}(\pi_{sc}))$
                \item $L(\Delta_\rho[-\frac{5}{2}, -\frac{5}{2}], \Delta_\rho[-\frac{5}{2}, -\frac{5}{2}], \Delta_\rho[-\frac{1}{2}, -\frac{3}{2}]; \pi_{sc})$
                \item $L(\Delta_\rho[-\frac{5}{2}, -\frac{5}{2}], \Delta_\rho[-\frac{5}{2}, -\frac{5}{2}], \Delta_\rho[-\frac{3}{2}, -\frac{3}{2}], \Delta_\rho[-\frac{1}{2}, -\frac{1}{2}]; \pi_{sc})$
            \end{itemize}
            \item $(\frac{1}{2}, \frac{3}{2}, \frac{5}{2}, \frac{7}{2})$: 
            \begin{itemize}
            \item $L(\Delta_\rho[-\frac{1}{2}, -\frac{7}{2}]; \pi_{sc})$  
            \item $L(\Delta_\rho[-\frac{3}{2}, -\frac{7}{2}]; T_{I,1}^{\frac{1}{2}}(\pi_{sc}))$
            \item $L(\Delta_\rho[-\frac{5}{2}, -\frac{7}{2}]; T_{I,1}^{\frac{3}{2}}(T_{I,1}^{\frac{1}{2}}(\pi_{sc})))$  
                \item $L(\Delta_\rho[-\frac{7}{2}, -\frac{7}{2}]; T_{I,1}^{\frac{5}{2}}(T_{I,1}^{\frac{3}{2}}(T_{I,1}^{\frac{1}{2}}(\pi_{sc}))))$
                \item $L(\Delta_\rho[-\frac{3}{2}, -\frac{7}{2}], \Delta_\rho[-\frac{1}{2},-\frac{1}{2}]; \pi_{sc})$
                \item $L(\Delta_\rho[-\frac{5}{2}, -\frac{7}{2}], \Delta_\rho[-\frac{1}{2}, -\frac{3}{2}]; \pi_{sc})$  
                \item $L(\Delta_\rho[-\frac{5}{2}, -\frac{7}{2}], \Delta_\rho[-\frac{3}{2}, -\frac{3}{2}]; T_{I,1}^{\frac{1}{2}}(\pi_{sc}))$
                \item $L(\Delta_\rho[-\frac{7}{2}, -\frac{7}{2}], \Delta_\rho[-\frac{1}{2}, -\frac{5}{2}]; \pi_{sc})$  
                 \item $L(\Delta_\rho[-\frac{7}{2}, -\frac{7}{2}], \Delta_\rho[-\frac{3}{2}, -\frac{5}{2}]; T_{I,1}^{\frac{1}{2}}(\pi_{sc}))$
                 \item $L(\Delta_\rho[-\frac{7}{2}, -\frac{7}{2}], \Delta_\rho[-\frac{5}{2}, -\frac{5}{2}]; T_{I,1}^{\frac{3}{2}}(T_{I,1}^{\frac{1}{2}}(\pi_{sc})))$  
                 \item $L(\Delta_\rho[-\frac{5}{2}, -\frac{7}{2}], \Delta_\rho[-\frac{3}{2}, -\frac{3}{2}], \Delta_\rho[-\frac{1}{2}, -\frac{1}{2}]; \pi_{sc})$
                 \item $L(\Delta_\rho[-\frac{7}{2}, -\frac{7}{2}], \Delta_\rho[-\frac{5}{2}, -\frac{5}{2}], \Delta_\rho[-\frac{3}{2}, -\frac{3}{2}]; T_{I,1}^{\frac{1}{2}}(\pi_{sc}))$
                 \item $L(\Delta_\rho[-\frac{7}{2}, -\frac{7}{2}], \Delta_\rho[-\frac{3}{2}, -\frac{5}{2}], \Delta_\rho[-\frac{1}{2}, -\frac{1}{2}]; \pi_{sc})$
                 \item $L(\Delta_\rho[-\frac{7}{2}, -\frac{7}{2}], \Delta_\rho[-\frac{5}{2}, -\frac{5}{2}], \Delta_\rho[-\frac{1}{2}, -\frac{3}{2}]; \pi_{sc})$
            \end{itemize}
        \end{enumerate}
        \item $(\alpha =1)$:
        \begin{enumerate}
        \item $(0,0,1,2)$:
        \begin{itemize}
            \item $L(\Delta_\rho[-2,-2]; T_{V,4}^{\pm}(T_{I,1}^{1}(\pi_{sc})))$
            \item $L(\Delta_\rho[0,-2]; T_{IV,3}(\pi_{sc}))$
            \item $L(\Delta_\rho[-1,-2]; T_{IV,5}(\pi_{sc}))$
            \item $L(\Delta_\rho[-2,-2], \Delta_\rho[-1,-1], \Delta_\rho[0,-1]; \pi_{sc})$
        \end{itemize}
        \item $(0,1,1,1)$: 
        \begin{itemize}
            \item $L(\Delta_\rho[-1,-1], \Delta_\rho[0,-1]; T_{I,1}^{1}(\pi_{sc}))$
            \item $L(\Delta_\rho[-1,-1], \Delta_\rho[-1,-1], \Delta_\rho[0,-1]; \pi_{sc})$
        \end{itemize}
        \item $(0,1,1,2)$: 
        \begin{itemize}
            \item $L(\Delta_\rho[-2,-2]; T_{I,2}^{1}(T_{IV,3}(\pi_{sc})))$
            \item $L(\Delta_\rho[-1,-1], \Delta_\rho[0,-2]; \pi_{sc})$ 
            \item $L(\Delta_\rho[-1,-2], \Delta_\rho[-1,-1], T_{V,2}^{\pm}(\pi_{sc}))$
            \item $L(\Delta_\rho[-2,-2], \Delta_\rho[0,-1]; T_{I,1}^{1}(\pi_{sc}))$

        \end{itemize}
            \item $(0,1,2,2)$: 
            \begin{itemize}
                \item $L(\Delta_\rho[-2,-2]; T_{V,2}^{\pm}(T_{I,1}^{2}(T_{I,1}^{1}(\pi_{sc}))))$
                \item $L(\Delta_\rho[-2,-2], \Delta_\rho[0,-2]; \pi_{sc})$  
                \item $L(\Delta_\rho[-2,-2], \Delta_\rho[-1,-2]; T_{IV,3}(\pi_{sc}))$
                \item $L(\Delta_\rho[-2,-2], \Delta_\rho[-2,-2]; T_{V,2}^{\pm}(T_{I,1}^{1}(\pi_{sc})))$
                \item $L(\Delta_\rho[-2,-2], \Delta_\rho[-2,-2], \Delta_\rho[-1,-1]; T_{IV,3}(\pi_{sc}))$
                \item $L(\Delta_\rho[-2,-2], \Delta_\rho[-2,-2], \Delta_\rho[0,-1]; \pi_{sc})$
            \end{itemize}
            \item $(0,1,2,3)$: 
            \begin{itemize}
            \item $L(\Delta_\rho[0,-3]; \pi_{sc})$
             \item $L(\Delta_\rho[-2,-3]; T_{V,2}^{\pm}(T_{I,1}^{1}(\pi_{sc})))$
                \item $L(\Delta_\rho[-3,-3]; T_{V,2}^{\pm}(T_{I,1}^{2}(T_{I,1}^{1}(\pi_{sc}))))$.
                \item $L(\Delta_\rho[-1,-3]; T_{IV,3}(\pi_{sc}))$
                \item $L(\Delta_\rho[-2,-3], \Delta_\rho[0,-1]; \pi_{sc})$ 
                \item $L(\Delta_\rho[-2,-3], \Delta_\rho[-1,-1]; T_{IV,3}(\pi_{sc}))$
                \item $L(\Delta_\rho[-2,-3], \Delta_\rho[-1,-1]; T_{V,2}^{\pm}(\pi_{sc}))$
                \item $L(\Delta_\rho[-3,-3], \Delta_\rho[0,-2]; \pi_{sc})$  
                 \item $L(\Delta_\rho[-3,-3], \Delta_\rho[-1,-2]; T_{IV,3}(\pi_{sc}))$
                 \item $L(\Delta_\rho[-3,-3], \Delta_\rho[-2,-2]; T_{IV,3}(T_{I,1}^{1}(\pi_{sc})))$ 
            \end{itemize}
            \item $(1,1,1,1)$: 
            \begin{itemize}
                \item $L(\Delta_\rho[-1,-1], \Delta_\rho[-1,-1], \Delta_\rho[-1,-1]; T_{I,1}^{1}(\pi_{sc}))$
                \item $L(\Delta_\rho[-1,-1], \Delta_\rho[-1,-1],\Delta_\rho[-1,-1],\Delta_\rho[-1,-1]; \pi_{sc})$
            \end{itemize}
            \item $(1,1,1,2)$: 
            \begin{itemize}
                \item $L(\Delta_\rho[-1,-1], \Delta_\rho[-1,-1]; T_{I,1}^{2}(T_{I,1}^{1}(\pi_{sc})))$
                \item $L(\Delta_\rho[-1,-2], \Delta_\rho[-1,-1]; T_{I,1}^{1}(\pi_{sc}))$
                 \item $L(\Delta_\rho[-2,-2], \Delta_\rho[-1,-1], \Delta_\rho[-1,-1]; T_{I,1}^{1}(\pi_{sc}))$
                 \item $L(\Delta_\rho[-1,-2], \Delta_\rho[-1,-1], \Delta_\rho[-1,-1]; \pi_{sc})$
                 \item $L(\Delta_\rho[-2,-2], \Delta_\rho[-1,-1], \Delta_\rho[-1,-1], \Delta_\rho[-1,-1]; \pi_{sc})$
            \end{itemize}
            \item $(1,1,2,2)$: 
            \begin{itemize}
                \item $L(\Delta_\rho[-1,-2]; T_{I,1}^{2}(T_{I,1}^{1}(\pi_{sc})))$
                \item $L(\Delta_\rho[-2,-2], \Delta_\rho[-1,-1]; T_{I,1}^{2}(T_{I,1}^{1}(\pi_{sc})))$
                \item $L(\Delta_\rho[-2,-2], \Delta_\rho[-1,-2]; T_{I,1}^{1}(\pi_{sc}))$
                \item $L(\Delta_\rho[-1,-2], \Delta_\rho[-1,-2]; \pi_{sc})$
                \item $L(\Delta_\rho[-2,-2], \Delta_\rho[-2,-2], \Delta_\rho[-1,-1]; T_{I,1}^{1}(\pi_{sc}))$
                \item $L(\Delta_\rho[-2,-2], \Delta_\rho[-1,-2],\Delta_\rho[-1,-1]; \pi_{sc})$
                \item $L(\Delta_\rho[-2,-2], \Delta_\rho[-2,-2], \Delta_\rho[-1,-1], \Delta_\rho[-1,-1]; \pi_{sc})$
            \end{itemize}
            \item $(1,1,2,3)$: 
            \begin{itemize}
            \item $L(\Delta_\rho[-1,-1]; T_{I,1}^{3}(T_{I,1}^{2}(T_{I,1}^{1}(\pi_{sc}))))$ 
                \item $L(\Delta_\rho[-1,-3]; T_{I,1}^{1}(\pi_{sc}))$
                \item $L(\Delta_\rho[-1,-3], \Delta_\rho[-1,-1]; \pi_{sc})$ 
                \item $L(\Delta_\rho[-2,-3], \Delta_\rho[-1,-1]; T_{I,1}^{1}(\pi_{sc}))$
                \item $L(\Delta_\rho[-3,-3], \Delta_\rho[-1,-1]; T_{I,1}^{2}(T_{I,1}^{1}(\pi_{sc})))$
                \item $L(\Delta_\rho[-3,-3], \Delta_\rho[-1,-2]; T_{I,1}^{1}(\pi_{sc}))$
                \item $L(\Delta_\rho[-2,-3], \Delta_\rho[-1,-1], \Delta_\rho[-1,-1]; \pi_{sc})$
                \item $L(\Delta_\rho[-3,-3], \Delta_\rho[-2,-2], \Delta_\rho[-1,-1]; T_{I,1}^{1}(\pi_{sc}))$
                
                \item $L(\Delta_\rho[-3,-3], \Delta_\rho[-1,-2],\Delta_\rho[-1,-1]; \pi_{sc})$
                \item $L(\Delta_\rho[-3,-3], \Delta_\rho[-2,-2], \Delta_\rho[-1,-1], \Delta_\rho[-1,-1]; \pi_{sc})$
            \end{itemize}
            \item $(1,2,2,2)$: 
            \begin{itemize}
            \item $L(\Delta_\rho[-2,-2], \Delta_\rho[-2,-2]; T_{I,1}^{2}(T_{I,1}^{1}(\pi_{sc})))$  
                \item $L(\Delta_\rho[-2,-2], \Delta_\rho[-2,-2], \Delta_\rho[-2,-2]; T_{I,1}^{1}(\pi_{sc}))$
                \item $L(\Delta_\rho[-2,-2], \Delta_\rho[-2,-2], \Delta_\rho[-1,-2]; \pi_{sc})$
                \item $L(\Delta_\rho[-2,-2], \Delta_\rho[-2,-2], \Delta_\rho[-2,-2], \Delta_\rho[-1,-1]; \pi_{sc})$
            \end{itemize}
            \item $(1,2,2,3)$: 
            \begin{itemize}
            \item $L(\Delta_\rho[-2,-2]; T_{I,1}^{3}(T_{I,1}^{2}(T_{I,1}^{1}(\pi_{sc}))))$  
                \item $L(\Delta_\rho[-2,-3]; T_{I,1}^{2}(T_{I,1}^{1}(\pi_{sc})))$
                \item $L(\Delta_\rho[-2,-2], \Delta_\rho[-1,-3]; \pi_{sc})$ 
                \item $L(\Delta_\rho[-2,-3], \Delta_\rho[-2,-2]; T_{I,1}^{1}(\pi_{sc}))$
                \item $L(\Delta_\rho[-2,-3], \Delta_\rho[-1,-2]; \pi_{sc})$
                \item $L(\Delta_\rho[-2,-3], \Delta_\rho[-2,-2], \Delta_\rho[-1,-1]; \pi_{sc})$  
                \item $L(\Delta_\rho[-3,-3], \Delta_\rho[-2,-2]; T_{I,1}^{2}(T_{I,1}^{1}(\pi_{sc})))$
                \item $L(\Delta_\rho[-3,-3], \Delta_\rho[-2,-2], \Delta_\rho[-2,-2]; T_{I,1}^{1}(\pi_{sc}))$
                \item $L(\Delta_\rho[-3,-3], \Delta_\rho[-2,-2], \Delta_\rho[-1,-2]; \pi_{sc})$
                \item $L(\Delta_\rho[-3,-3], \Delta_\rho[-2,-2], \Delta_\rho[-2,-2], \Delta_\rho[-1,-1]; \pi_{sc})$
            \end{itemize}
            \item $(1,2,3,3)$: 
            \begin{itemize}
            \item $L(\Delta_\rho[-3,-3]; T_{I,1}^{3}(T_{I,1}^{2}(T_{I,1}^{1}(\pi_{sc}))))$  
            \item $L(\Delta_\rho[-3,-3], \Delta_\rho[-1,-3]; \pi_{sc})$  
                \item $L(\Delta_\rho[-3,-3], \Delta_\rho[-3,-3]; T_{I,1}^{2}(T_{I,1}^{1}(\pi_{sc})))$
                \item $L(\Delta_\rho[-3,-3], \Delta_\rho[-2,-3]; T_{I,1}^{1}(\pi_{sc}))$
                \item $L(\Delta_\rho[-3,-3], \Delta_\rho[-3,-3], \Delta_\rho[-2,-2]; T_{I,1}^{1}(\pi_{sc}))$
                \item $L(\Delta_\rho[-3,-3], \Delta_\rho[-2,-3],\Delta_\rho[-1,-1]; \pi_{sc})$
                \item $L(\Delta_\rho[-3,-3], \Delta_\rho[-3,-3], \Delta_\rho[-1,-2]; \pi_{sc})$
                \item $L(\Delta_\rho[-3,-3], \Delta_\rho[-3,-3], \Delta_\rho[-2,-2], \Delta_\rho[-1,-1]; \pi_{sc})$
            \end{itemize}
            \item $(1,2,3,4)$: 
            \begin{itemize}
                \item $L(\Delta_\rho[-1,-4]; \pi_{sc})$            
                \item $L(\Delta_\rho[-2,-4]; T_{I,1}^{1}(\pi_{sc}))$
                \item $L(\Delta_\rho[-3,-4]; T_{I,1}^{2}(T_{I,1}^{1}(\pi_{sc})))$
                \item $L(\Delta_\rho[-4,-4]; T_{I,1}^{3}(T_{I,1}^{2}(T_{I,1}^{1}(\pi_{sc}))))$ 
                \item $L(\Delta_\rho[-2,-4], \Delta_\rho[-1,-1]; \pi_{sc})$
                \item $L(\Delta_\rho[-3,-4], \Delta_\rho[-1,-2]; \pi_{sc})$
                 \item $L(\Delta_\rho[-3,-4], \Delta_\rho[-2,-2]; T_{I,1}^{1}(\pi_{sc}))$
                 \item $L(\Delta_\rho[-4,-4], \Delta_\rho[-1,-3]; \pi_{sc})$  
                \item $L(\Delta_\rho[-4,-4], \Delta_\rho[-2,-3]; T_{I,1}^{1}(\pi_{sc}))$
                \item $L(\Delta_\rho[-4,-4], \Delta_\rho[-3,-3]; T_{I,1}^{2}(T_{I,1}^{1}(\pi_{sc})))$
                \item $L(\Delta_\rho[-3,-4], \Delta_\rho[-2,-2], \Delta_\rho[-1,-1]; \pi_{sc})$
                \item $L(\Delta_\rho[-4,-4], \Delta_\rho[-3,-3], \Delta_\rho[-2,-2]; T_{I,1}^{1}(\pi_{sc}))$
                \item $L(\Delta_\rho[-4,-4], \Delta_\rho[-2,-3],\Delta_\rho[-1,-1]; \pi_{sc})$
                \item $L(\Delta_\rho[-4,-4], \Delta_\rho[-3,-3], \Delta_\rho[-1,-2]; \pi_{sc})$
            \end{itemize}
        \end{enumerate}
        \item $(\alpha = \frac{3}{2})$: 
        \begin{enumerate}
        \item $(\frac{1}{2}, \frac{1}{2}, \frac{3}{2}, \frac{3}{2})$: 
        \begin{itemize}
        \item $L(\Delta_\rho[-\frac{3}{2}, -\frac{3}{2}]; T_{I,2}^{\frac{1}{2}}(T_{I,1}^{\frac{3}{2}}(\pi_{sc})))$  
            \item $L(\Delta_\rho[-\frac{3}{2}, -\frac{3}{2}], \Delta_\rho[-\frac{1}{2}, -\frac{1}{2}]; T_{I,1}^{\frac{1}{2}}(T_{I,1}^{\frac{3}{2}}(\pi_{sc})))$
            \item $L(\Delta_\rho[-\frac{1}{2}, -\frac{3}{2}], \Delta_\rho[-\frac{1}{2}, -\frac{1}{2}]; T_{I,1}^{\frac{3}{2}}(\pi_{sc}))$
            \item $L(\Delta_\rho[-\frac{1}{2}, -\frac{3}{2}], \Delta_\rho[-\frac{1}{2}, -\frac{3}{2}]; \pi_{sc})$  
            \item $L(\Delta_\rho[-\frac{3}{2}, -\frac{3}{2}], \Delta_\rho[-\frac{3}{2}, -\frac{3}{2}]; T_{II,3}^{\frac{1}{2}}(\pi_{sc}))$ 
            \item $L(\Delta_\rho[-\frac{3}{2}, -\frac{3}{2}], \Delta_\rho[-\frac{1}{2}, -\frac{1}{2}], \Delta_\rho[-\frac{1}{2}, -\frac{1}{2}]; T_{I,1}^{\frac{3}{2}}(\pi_{sc}))$
            \item $L(\Delta_\rho[-\frac{3}{2}, -\frac{3}{2}], \Delta_\rho[-\frac{1}{2}, -\frac{3}{2}], \Delta_\rho[-\frac{1}{2}, -\frac{1}{2}]; \pi_{sc})$
        \end{itemize}
            \item $(\frac{1}{2}, \frac{1}{2}, \frac{3}{2}, \frac{5}{2})$: 
            \begin{itemize}
                \item $L(\Delta_\rho[-\frac{3}{2}, -\frac{5}{2}]; T_{II,3}^{\frac{1}{2}}(\pi_{sc}))$
                \item $L(\Delta_\rho[-\frac{1}{2}, -\frac{1}{2}]; T_{I,1}^{\frac{1}{2}}(T_{I,1}^{\frac{5}{2}}(T_{I,1}^{\frac{3}{2}}(\pi_{sc}))))$
                \item $L(\Delta_\rho[-\frac{5}{2}, -\frac{5}{2}]; T_{I,1}^{\frac{3}{2}}(T_{II,3}^{\frac{1}{2}}(\pi_{sc})))$
                \item $L(\Delta_\rho[-\frac{5}{2}, -\frac{5}{2}]; T_{I,2}^{\frac{1}{2}}(T_{I,1}^{\frac{3}{2}}(\pi_{sc})))$  
                \item $L(\Delta_\rho[\frac{1}{2}, -\frac{5}{2}]; \pi_{sc})$  
                \item $L(\Delta_\rho[-\frac{1}{2}, -\frac{5}{2}], \Delta_\rho[-\frac{1}{2}, -\frac{1}{2}]; \pi_{sc})$  
                \item $L(\Delta_\rho[-\frac{5}{2}, -\frac{5}{2}], \Delta_\rho[-\frac{1}{2}, -\frac{1}{2}]; T_{I,1}^{\frac{1}{2}}(T_{I,1}^{\frac{3}{2}}(\pi_{sc})))$ 
                \item $L(\Delta_\rho[-\frac{5}{2}, -\frac{5}{2}], \Delta_\rho[-\frac{1}{2}, -\frac{1}{2}], \Delta_\rho[-\frac{1}{2}, -\frac{1}{2}]; T_{I,1}^{\frac{3}{2}}(\pi_{sc}))$
                \item $L(\Delta_\rho[-\frac{5}{2},-\frac{5}{2}], \Delta_\rho[-\frac{3}{2},-\frac{3}{2}], \Delta_\rho[-\frac{1}{2},-\frac{1}{2}], \Delta_\rho[-\frac{1}{2},-\frac{1}{2}]; \pi_{sc})$
            \end{itemize}
            \item $(\frac{1}{2}, \frac{3}{2}, \frac{3}{2}, \frac{3}{2})$: 
            \begin{itemize}
                \item $L(\Delta_\rho[-\frac{3}{2}, -\frac{3}{2}], \Delta_\rho[-\frac{3}{2}, -\frac{3}{2}]; T_{I,1}^{\frac{1}{2}}(T_{I,1}^{\frac{3}{2}}(\pi_{sc})))$
                \item $L(\Delta_\rho[-\frac{3}{2}, -\frac{3}{2}], \Delta_\rho[-\frac{1}{2}, -\frac{3}{2}]; T_{I,1}^{\frac{3}{2}}(\pi_{sc}))$
                \item $L(\Delta_\rho[-\frac{3}{2}, -\frac{3}{2}], \Delta_\rho[-\frac{3}{2}, -\frac{3}{2}], \Delta_\rho[-\frac{1}{2}, -\frac{1}{2}]; T_{I,1}^{\frac{3}{2}}(\pi_{sc}))$
                \item $L(\Delta_\rho[-\frac{3}{2}, -\frac{3}{2}], \Delta_\rho[-\frac{3}{2}, -\frac{3}{2}], \Delta_\rho[-\frac{1}{2}, -\frac{3}{2}]; \pi_{sc})$
                \item $L(\Delta_\rho[-\frac{3}{2},-\frac{3}{2}], \Delta_\rho[-\frac{3}{2},-\frac{3}{2}], \Delta_\rho[-\frac{3}{2},-\frac{3}{2}], \Delta_\rho[-\frac{1}{2},-\frac{1}{2}]; \pi_{sc})$
            \end{itemize}
            \item $(\frac{1}{2}, \frac{3}{2}, \frac{3}{2}, \frac{5}{2})$: 
            \begin{itemize}
            \item $L(\Delta_\rho[-\frac{1}{2}, -\frac{3}{2}]; T_{I,1}^{\frac{5}{2}}(T_{I,1}^{\frac{3}{2}}(\pi_{sc})))$  
                \item $L(\Delta_\rho[-\frac{1}{2}, -\frac{5}{2}]; T_{I,1}^{\frac{3}{2}}(\pi_{sc}))$
                \item $L(\Delta_\rho[-\frac{3}{2}, -\frac{5}{2}]; T_{I,1}^{\frac{1}{2}}(T_{I,1}^{\frac{5}{2}}(\pi_{sc})))$  
                \item $L(\Delta_\rho[-\frac{3}{2}, -\frac{3}{2}]; T_{I,1}^{\frac{1}{2}}(T_{I,1}^{\frac{5}{2}}(T_{I,1}^{\frac{3}{2}}(\pi_{sc}))))$
                \item $L(\Delta_\rho[-\frac{5}{2}, -\frac{5}{2}], \Delta_\rho[-\frac{3}{2}, -\frac{3}{2}]; T_{I,1}^{\frac{1}{2}}(T_{I,1}^{\frac{3}{2}}(\pi_{sc})))$
                \item $L(\Delta_\rho[-\frac{3}{2}, -\frac{3}{2}], \Delta_\rho[-\frac{1}{2}, -\frac{5}{2}]; \pi_{sc})$  
                \item $L(\Delta_\rho[-\frac{3}{2}, -\frac{5}{2}], \Delta_\rho[-\frac{1}{2}, -\frac{1}{2}]; T_{I,1}^{\frac{3}{2}}(\pi_{sc}))$
                \item $L(\Delta_\rho[-\frac{3}{2}, -\frac{5}{2}], \Delta_\rho[-\frac{3}{2}, -\frac{3}{2}], \Delta_\rho[-\frac{1}{2}, -\frac{1}{2}]; \pi_{sc})$
                \item $L(\Delta_\rho[-\frac{5}{2}, -\frac{5}{2}], \Delta_\rho[-\frac{1}{2}, -\frac{3}{2}]; T_{I,1}^{\frac{3}{2}}(\pi_{sc}))$
                \item $L(\Delta_\rho[-\frac{5}{2}, -\frac{5}{2}], \Delta_\rho[-\frac{3}{2}, -\frac{3}{2}], \Delta_\rho[-\frac{1}{2}, -\frac{3}{2}]; \pi_{sc})$
                \item $L(\Delta_\rho[-\frac{5}{2},-\frac{5}{2}], \Delta_\rho[-\frac{3}{2},-\frac{3}{2}], \Delta_\rho[-\frac{3}{2},-\frac{3}{2}], \Delta_\rho[-\frac{1}{2},-\frac{1}{2}]; \pi_{sc})$
            \end{itemize}
            \item $(\frac{1}{2}, \frac{3}{2}, \frac{5}{2}, \frac{5}{2})$: 
            \begin{itemize}
                 \item $L(\Delta_\rho[-\frac{5}{2}, -\frac{5}{2}]; T_{I,1}^{\frac{1}{2}}(T_{I,1}^{\frac{5}{2}}(T_{I,1}^{\frac{3}{2}}(\pi_{sc}))))$
                \item $L(\Delta_\rho[-\frac{5}{2}, -\frac{5}{2}], \Delta_\rho[-\frac{1}{2}, -\frac{1}{2}]; T_{I,1}^{\frac{5}{2}}(T_{I,1}^{\frac{3}{2}}(\pi_{sc})))$
                \item $L(\Delta_\rho[-\frac{5}{2}, -\frac{5}{2}], \Delta_\rho[-\frac{1}{2}, -\frac{5}{2}]; \pi_{sc})$ 
                \item $L(\Delta_\rho[-\frac{5}{2}, -\frac{5}{2}], \Delta_\rho[-\frac{5}{2}, -\frac{5}{2}]; T_{I,1}^{\frac{1}{2}}(T_{I,1}^{\frac{3}{2}}(\pi_{sc})))$ 
                \item $L(\Delta_\rho[-\frac{5}{2}, -\frac{5}{2}], \Delta_\rho[-\frac{5}{2}, -\frac{5}{2}], \Delta_\rho[-\frac{1}{2}, -\frac{1}{2}]; T_{I,1}^{\frac{3}{2}}(\pi_{sc}))$
                \item $L(\Delta_\rho[-\frac{5}{2}, -\frac{5}{2}], \Delta_\rho[-\frac{3}{2}, -\frac{5}{2}], \Delta_\rho[-\frac{1}{2}, -\frac{1}{2}]; \pi_{sc})$
                \item $L(\Delta_\rho[-\frac{5}{2}, -\frac{5}{2}], \Delta_\rho[-\frac{5}{2}, -\frac{5}{2}], \Delta_\rho[-\frac{1}{2}, -\frac{3}{2}]; \pi_{sc})$
                \item $L(\Delta_\rho[-\frac{5}{2},-\frac{5}{2}], \Delta_\rho[-\frac{5}{2},-\frac{5}{2}], \Delta_\rho[-\frac{3}{2},-\frac{3}{2}], \Delta_\rho[-\frac{1}{2},-\frac{1}{2}]; \pi_{sc})$
            \end{itemize}
            \item $(\frac{1}{2}, \frac{3}{2}, \frac{5}{2}, \frac{7}{2})$:
            \begin{itemize}
                \item $L(\Delta_\rho[-\frac{1}{2}, -\frac{7}{2}]; \pi_{sc})$
                \item $L(\Delta_\rho[-\frac{5}{2}, -\frac{7}{2}]; T_{I,1}^{\frac{3}{2}}(T_{I,1}^{\frac{1}{2}}(\pi_{sc})))$
                \item $L(\Delta_\rho[-\frac{7}{2}, -\frac{7}{2}]; T_{I,1}^{\frac{1}{2}}(T_{I,1}^{\frac{5}{2}}(T_{I,1}^{\frac{3}{2}}(\pi_{sc}))))$
                \item $L(\Delta_\rho[-\frac{3}{2}, -\frac{7}{2}], \Delta_\rho[-\frac{1}{2}, -\frac{1}{2}]; \pi_{sc})$
                \item $L(\Delta_\rho[-\frac{5}{2}, -\frac{7}{2}], \Delta_\rho[-\frac{1}{2}, -\frac{1}{2}]; T_{I,1}^{\frac{3}{2}}(\pi_{sc}))$  
                \item $L(\Delta_\rho[-\frac{5}{2}, -\frac{7}{2}], \Delta_\rho[-\frac{1}{2}, -\frac{3}{2}]; \pi_{sc})$   
                \item $L(\Delta_\rho[-\frac{7}{2}, -\frac{7}{2}], \Delta_\rho[-\frac{1}{2}, -\frac{1}{2}]; T_{I,1}^{\frac{5}{2}}(T_{I,1}^{\frac{3}{2}}(\pi_{sc})))$ 
                \item $L(\Delta_\rho[-\frac{7}{2}, -\frac{7}{2}], \Delta_\rho[-\frac{1}{2}, -\frac{5}{2}]; \pi_{sc})$  
                \item $L(\Delta_\rho[-\frac{7}{2}, -\frac{7}{2}], \Delta_\rho[-\frac{5}{2}, -\frac{5}{2}]; T_{I,1}^{\frac{1}{2}}(T_{I,1}^{\frac{3}{2}}(\pi_{sc})))$
                \item $L(\Delta_\rho[-\frac{5}{2}, -\frac{7}{2}], \Delta_\rho[-\frac{3}{2}, -\frac{3}{2}], \Delta_\rho[-\frac{1}{2}, -\frac{1}{2}]; \pi_{sc})$
                \item $L(\Delta_\rho[-\frac{7}{2}, -\frac{7}{2}], \Delta_\rho[-\frac{3}{2}, -\frac{5}{2}], \Delta_\rho[-\frac{1}{2}, -\frac{1}{2}]; \pi_{sc})$ 
                \item $L(\Delta_\rho[-\frac{7}{2}, -\frac{7}{2}], \Delta_\rho[-\frac{5}{2}, -\frac{5}{2}], \Delta_\rho[-\frac{1}{2}, -\frac{1}{2}]; T_{I,1}^{\frac{3}{2}}(\pi_{sc}))$  
                \item $L(\Delta_\rho[-\frac{7}{2}, -\frac{7}{2}], \Delta_\rho[-\frac{3}{2}, -\frac{5}{2}], \Delta_\rho[-\frac{1}{2}, -\frac{1}{2}]; \pi_{sc})$
            \end{itemize}
            \item $(\frac{3}{2}, \frac{3}{2}, \frac{3}{2}, \frac{3}{2})$: 
            \begin{itemize}
                \item $L(\Delta_\rho[-\frac{3}{2}, -\frac{3}{2}], \Delta_\rho[-\frac{3}{2}, -\frac{3}{2}], \Delta_\rho[-\frac{3}{2}, -\frac{3}{2}]; T_{I,1}^{\frac{3}{2}}(\pi_{sc}))$
                \item $L(\Delta_\rho[-\frac{3}{2},-\frac{3}{2}], \Delta_\rho[-\frac{3}{2},-\frac{3}{2}],\Delta_\rho[-\frac{3}{2},-\frac{3}{2}],\Delta_\rho[-\frac{3}{2},-\frac{3}{2}]; \pi_{sc})$
            \end{itemize}
            \item $(\frac{3}{2}, \frac{3}{2}, \frac{3}{2}, \frac{5}{2})$: 
            \begin{itemize}
                \item $L(\Delta_\rho[-\frac{3}{2}, -\frac{3}{2}], \Delta_\rho[-\frac{3}{2},-\frac{3}{2}]; T_{I,1}^{\frac{5}{2}}(T_{I,1}^{\frac{3}{2}}(\pi_{sc})))$
                \item $L(\Delta_\rho[-\frac{3}{2}, -\frac{5}{2}], \Delta_\rho[-\frac{3}{2}, -\frac{3}{2}]; T_{I,1}^{\frac{3}{2}}(\pi_{sc}))$
                \item $L(\Delta_\rho[-\frac{5}{2}, -\frac{5}{2}], \Delta_\rho[-\frac{3}{2}, -\frac{3}{2}], \Delta_\rho[-\frac{3}{2}, -\frac{3}{2}]; T_{I,1}^{\frac{3}{2}}(\pi_{sc}))$
                \item $L(\Delta_\rho[-\frac{3}{2}, -\frac{5}{2}], \Delta_\rho[-\frac{3}{2}, -\frac{3}{2}], \Delta_\rho[-\frac{3}{2}, -\frac{3}{2}]; \pi_{sc})$
                \item $L(\Delta_\rho[-\frac{5}{2},-\frac{5}{2}], \Delta_\rho[-\frac{3}{2},-\frac{3}{2}], \Delta_\rho[-\frac{3}{2},-\frac{3}{2}], \Delta_\rho[-\frac{3}{2},-\frac{3}{2}]; \pi_{sc})$
            \end{itemize}
            \item $(\frac{3}{2}, \frac{3}{2}, \frac{5}{2}, \frac{5}{2})$: 
            \begin{itemize}
            \item $L(\Delta_\rho[-\frac{3}{2}, -\frac{5}{2}]; T_{I,1}^{\frac{5}{2}}(T_{I,1}^{\frac{3}{2}}(\pi_{sc})))$  
                \item $L(\Delta_\rho[-\frac{5}{2}, -\frac{5}{2}],\Delta_\rho[-\frac{3}{2}, -\frac{3}{2}]T_{I,1}^{\frac{5}{2}}(T_{I,1}^{\frac{3}{2}}(\pi_{sc})))$
                \item $L(\Delta_\rho[-\frac{5}{2}, -\frac{5}{2}], \Delta_\rho[-\frac{3}{2}, -\frac{5}{2}]; T_{I,1}^{\frac{3}{2}}(\pi_{sc}))$
                \item $L(\Delta_\rho[-\frac{3}{2}, -\frac{5}{2}], \Delta_\rho[-\frac{3}{2}, -\frac{5}{2}]; \pi_{sc})$
                \item $L(\Delta_\rho[-\frac{5}{2}, -\frac{5}{2}], \Delta_\rho[-\frac{5}{2}, -\frac{5}{2}], \Delta_\rho[-\frac{3}{2}, -\frac{3}{2}]; T_{I,1}^{\frac{3}{2}}(\pi_{sc}))$
                \item $L(\Delta_\rho[-\frac{5}{2}, -\frac{5}{2}], \Delta_\rho[-\frac{3}{2}, -\frac{5}{2}], \Delta_\rho[-\frac{3}{2}, -\frac{3}{2}]; \pi_{sc})$
                \item $L(\Delta_\rho[-\frac{5}{2},-\frac{5}{2}], \Delta_\rho[-\frac{5}{2},-\frac{5}{2}], \Delta_\rho[-\frac{3}{2},-\frac{3}{2}], \Delta_\rho[-\frac{3}{2},-\frac{3}{2}]; \pi_{sc})$
            \end{itemize}
            \item $(\frac{3}{2}, \frac{3}{2}, \frac{5}{2}, \frac{7}{2})$: 
            \begin{itemize}
            \item $L(\Delta_\rho[-\frac{3}{2}, -\frac{3}{2}]; T_{I,1}^{\frac{7}{2}}(T_{I,1}^{\frac{5}{2}}(T_{I,1}^{\frac{3}{2}}(\pi_{sc}))))$
                \item $L(\Delta_\rho[-\frac{3}{2}, -\frac{7}{2}]; T_{I,1}^{\frac{3}{2}}(\pi_{sc}))$
                \item $L(\Delta_\rho[-\frac{3}{2}, -\frac{3}{2}]; T_{I,1}^{\frac{7}{2}}(T_{I,1}^{\frac{5}{2}}(T_{I,1}^[\frac{3}{2}(\pi_{sc}))))$
                \item $L(\Delta_\rho[-\frac{5}{2}, -\frac{7}{2}], \Delta_\rho[-\frac{3}{2}, -\frac{3}{2}]; T_{I,1}^{\frac{3}{2}}(\pi_{sc}))$
                \item $L(\Delta_\rho[-\frac{7}{2}, -\frac{7}{2}], \Delta_\rho[-\frac{3}{2}, -\frac{5}{2}]; T_{I,1}^{\frac{3}{2}}(\pi_{sc}))$
                \item $L(\Delta_\rho[-\frac{3}{2}, -\frac{7}{2}], \Delta_\rho[-\frac{3}{2}, -\frac{3}{2}]; \pi_{sc})$
                \item $L(\Delta_\rho[-\frac{5}{2}, -\frac{7}{2}], \Delta_\rho[-\frac{3}{2}, -\frac{3}{2}]; T_{I,1}^{\frac{3}{2}}(\pi_{sc}))$ 
                \item $L(\Delta_\rho[-\frac{7}{2}, -\frac{7}{2}], \Delta_\rho[-\frac{5}{2}, -\frac{5}{2}], \Delta_\rho[-\frac{3}{2}, -\frac{3}{2}]; T_{I,1}^{\frac{3}{2}}(\pi_{sc})$
                \item $L(\Delta_\rho[-\frac{7}{2}, -\frac{7}{2}], \Delta_\rho[-\frac{3}{2}, -\frac{5}{2}], \Delta_\rho[-\frac{3}{2}, -\frac{3}{2}]; \pi_{sc})$
                \item $L(\Delta_\rho[-\frac{7}{2},-\frac{7}{2}], \Delta_\rho[-\frac{5}{2},-\frac{5}{2}], \Delta_\rho[-\frac{3}{2},-\frac{3}{2}], \Delta_\rho[-\frac{3}{2},-\frac{3}{2}]; \pi_{sc})$
            \end{itemize}
            \item $(\frac{3}{2}, \frac{5}{2}, \frac{5}{2}, \frac{5}{2})$: 
            \begin{itemize}
                \item $L(\Delta_\rho[-\frac{5}{2}, -\frac{5}{2}], \Delta_\rho[-\frac{5}{2},-\frac{5}{2}]; T_{I,1}^{\frac{5}{2}}(T_{I,1}^{\frac{3}{2}}(\pi_{sc})))$
                \item $L(\Delta_\rho[-\frac{5}{2}, -\frac{5}{2}], \Delta_\rho[-\frac{5}{2}, -\frac{5}{2}], \Delta_\rho[-\frac{5}{2}, -\frac{5}{2}]; T_{I,1}^{\frac{3}{2}}(\pi_{sc}))$
                \item $L(\Delta_\rho[-\frac{5}{2}, -\frac{5}{2}], \Delta_\rho[-\frac{5}{2}, -\frac{5}{2}], \Delta_\rho[-\frac{3}{2}, -\frac{5}{2}]; \pi_{sc})$
                \item $L(\Delta_\rho[-\frac{5}{2},-\frac{5}{2}], \Delta_\rho[-\frac{5}{2},-\frac{5}{2}], \Delta_\rho[-\frac{5}{2},-\frac{5}{2}], \Delta_\rho[-\frac{3}{2},-\frac{3}{2}]; \pi_{sc})$
            \end{itemize}
            \item $(\frac{3}{2}, \frac{5}{2}, \frac{5}{2}, \frac{7}{2})$:
            \begin{itemize}
                 \item $L(\Delta_\rho[-\frac{5}{2}, -\frac{5}{2}]; T_{I,1}^{\frac{7}{2}}(T_{I,1}^{\frac{5}{2}}(T_{I,1}^{\frac{3}{2}}(\pi_{sc}))))$
                 \item $L(\Delta_\rho[-\frac{5}{2}, -\frac{7}{2}]; T_{I,1}^{\frac{3}{2}}(\pi_{sc}))$ + 
                 \item $L(\Delta_\rho[-\frac{7}{2}, -\frac{7}{2}],\Delta_\rho[-\frac{5}{2}, -\frac{5}{2}]T_{I,1}^{\frac{5}{2}}(T_{I,1}^{\frac{3}{2}}(\pi_{sc})))$
                 \item $L(\Delta_\rho[-\frac{5}{2}, -\frac{7}{2}], \Delta_\rho[-\frac{5}{2}, -\frac{5}{2}]; T_{I,1}^{\frac{3}{2}}(\pi_{sc}))$
                 \item $L(\Delta_\rho[-\frac{5}{2}, -\frac{7}{2}], \Delta_\rho[-\frac{3}{2}, -\frac{5}{2}]; \pi_{sc})$
                 \item $L(\Delta_\rho[-\frac{5}{2}, -\frac{5}{2}],\Delta_\rho[-\frac{3}{2}, -\frac{7}{2}]; \pi_{sc})$
                 \item $L(\Delta_\rho[-\frac{5}{2}, -\frac{7}{2}], \Delta_\rho[-\frac{5}{2}, -\frac{5}{2}], \Delta_\rho[-\frac{3}{2}, -\frac{3}{2}]; \pi_{sc})$ + 
                 \item $L(\Delta_\rho[-\frac{7}{2}, -\frac{7}{2}], \Delta_\rho[-\frac{5}{2}, -\frac{5}{2}], \Delta_\rho[-\frac{5}{2}, -\frac{5}{2}]; T_{I,1}^{\frac{3}{2}}(\pi_{sc}))$
                 \item $L(\Delta_\rho[-\frac{7}{2}, -\frac{7}{2}], \Delta_\rho[-\frac{5}{2}, -\frac{5}{2}], \Delta_\rho[-\frac{3}{2}, -\frac{5}{2}]; \pi_{sc})$
                 \item $L(\Delta_\rho[-\frac{7}{2},-\frac{7}{2}], \Delta_\rho[-\frac{5}{2},-\frac{5}{2}], \Delta_\rho[-\frac{5}{2},-\frac{5}{2}], \Delta_\rho[-\frac{3}{2},-\frac{3}{2}]; \pi_{sc})$
            \end{itemize}
            \item $(\frac{3}{2}, \frac{5}{2}, \frac{7}{2}, \frac{7}{2})$:
            \begin{itemize}
                 \item $L(\Delta_\rho[-\frac{7}{2}, -\frac{7}{2}]; T_{I,1}^{\frac{7}{2}}(T_{I,1}^{\frac{5}{2}}(T_{I,1}^{\frac{3}{2}}(\pi_{sc}))))$
                 \item $L(\Delta_\rho[-\frac{7}{2}, -\frac{7}{2}], \Delta_\rho[-\frac{7}{2},-\frac{7}{2}]; T_{I,1}^{\frac{5}{2}}(T_{I,1}^{\frac{3}{2}}(\pi_{sc})$
                 \item $L(\Delta_\rho[-\frac{7}{2}, -\frac{7}{2}], \Delta_\rho[-\frac{5}{2}, -\frac{7}{2}]; T_{I,1}^{\frac{3}{2}}(\pi_{sc}))$
                 \item $L(\Delta_\rho[-\frac{7}{2}, -\frac{7}{2}],\Delta_\rho[-\frac{3}{2}, -\frac{7}{2}]; \pi_{sc})$
                 \item $L(\Delta_\rho[-\frac{7}{2}, -\frac{7}{2}], \Delta_\rho[-\frac{7}{2}, -\frac{7}{2}], \Delta_\rho[-\frac{5}{2}, -\frac{5}{2}]; T_{I,1}^{\frac{3}{2}}(\pi_{sc}))$
                 \item $L(\Delta_\rho[-\frac{7}{2}, -\frac{7}{2}], \Delta_\rho[-\frac{5}{2}, -\frac{7}{2}], \Delta_\rho[-\frac{3}{2}, -\frac{3}{2}]; \pi_{sc})$
                 \item $L(\Delta_\rho[-\frac{7}{2}, -\frac{7}{2}], \Delta_\rho[-\frac{7}{2}, -\frac{7}{2}], \Delta_\rho[-\frac{3}{2}, -\frac{5}{2}]; \pi_{sc})$
                 \item $L(\Delta_\rho[-\frac{7}{2},-\frac{7}{2}], \Delta_\rho[-\frac{7}{2},-\frac{7}{2}], \Delta_\rho[-\frac{5}{2},-\frac{5}{2}], \Delta_\rho[-\frac{3}{2},-\frac{3}{2}]; \pi_{sc})$
            \end{itemize}
            \item $(\frac{3}{2}, \frac{5}{2}, \frac{7}{2}, \frac{9}{2})$: 
            \begin{itemize}
            \item $L(\Delta_\rho[-\frac{3}{2}, -\frac{9}{2}]; \pi_{sc})$
                \item $L(\Delta_\rho[-\frac{5}{2}, -\frac{9}{2}]; T_{I,1}^{\frac{3}{2}}(\pi_{sc}))$
                \item $L(\Delta_\rho[-\frac{7}{2}, -\frac{9}{2}]; T_{I,1}^{\frac{5}{2}}(T_{I,1}^{\frac{3}{2}}(\pi_{sc})))$  
                \item $L(\Delta_\rho[-\frac{9}{2}, -\frac{9}{2}]; T_{I,1}^{\frac{7}{2}}(T_{I,1}^{\frac{5}{2}}(T_{I,1}^[\frac{3}{2}(\pi_{sc}))))$
                \item $L(\Delta_\rho[-\frac{9}{2}, -\frac{9}{2}],\Delta_\rho[-\frac{7}{2}, -\frac{7}{2}];T_{I,1}^{\frac{5}{2}}(T_{I,1}^{\frac{3}{2}}(\pi_{sc})))$
                \item $L(\Delta_\rho[-\frac{7}{2}, -\frac{9}{2}], \Delta_\rho[-\frac{5}{2}, -\frac{5}{2}]; T_{I,1}^{\frac{3}{2}}(\pi_{sc}))$
                \item $L(\Delta_\rho[-\frac{9}{2}, -\frac{9}{2}], \Delta_\rho[-\frac{5}{2}, -\frac{7}{2}]; T_{I,1}^{\frac{3}{2}}(\pi_{sc}))$
                \item $L(\Delta_\rho[-\frac{7}{2}, -\frac{9}{2}], \Delta_\rho[-\frac{3}{2}, -\frac{5}{2}]; \pi_{sc})$
                \item $L(\Delta_\rho[-\frac{5}{2}, -\frac{9}{2}], \Delta_\rho[-\frac{3}{2}, -\frac{3}{2}]; \pi_{sc})$
                \item $L(\Delta_\rho[-\frac{9}{2}, -\frac{9}{2}],\Delta_\rho[-\frac{3}{2}, -\frac{7}{2}]; \pi_{sc})$
                \item $L(\Delta_\rho[-\frac{9}{2}, -\frac{9}{2}], \Delta_\rho[-\frac{7}{2}, -\frac{7}{2}], \Delta_\rho[-\frac{5}{2}, -\frac{5}{2}]; T_{I,1}^{\frac{3}{2}}(\pi_{sc}))$
                \item $L(\Delta_\rho[-\frac{7}{2}, -\frac{9}{2}], \Delta_\rho[-\frac{5}{2}, -\frac{5}{2}], \Delta_\rho[-\frac{3}{2}, -\frac{3}{2}]; \pi_{sc})$
                \item $L(\Delta_\rho[-\frac{9}{2}, -\frac{9}{2}], \Delta_\rho[-\frac{5}{2}, -\frac{7}{2}], \Delta_\rho[-\frac{3}{2}, -\frac{3}{2}]; \pi_{sc})$
                \item $L(\Delta_\rho[-\frac{9}{2}, -\frac{9}{2}], \Delta_\rho[-\frac{7}{2}, -\frac{7}{2}], \Delta_\rho[-\frac{3}{2}, -\frac{5}{2}]; \pi_{sc})$
             \end{itemize}
        \end{enumerate}
        \item $(\alpha = 2)$: 
        \begin{enumerate}
        \item $(0,1,2,2)$: 
        \begin{itemize}
        \item $L(\Delta_\rho[0,-2]; T_{I,1}^{2}(\pi_{sc}))$
        \item $L(\Delta_\rho[-1,-2]; T_{IV,3}(T_{I,1}^{2}(\pi_{sc})))$
            \item $L(\Delta_\rho[-2,-2]; T_{V,2}^{\pm}(T_{I,1}^{1}(T_{I,1}^{2}(\pi_{sc}))))$
            \item $L(\Delta_\rho[-2,-2], \Delta_\rho[0,-1]; T_{I,1}^{2}(\pi_{sc}))$
            \item $L(\Delta_\rho[-2,-2], \Delta_\rho[0,-2];\pi_{sc})$ 
            \item $L(\Delta_\rho[-2,-2], \Delta_\rho[-1,-1]; T_{IV,3}(T_{I,1}^{2}(\pi_{sc})))$  
            \item $L(\Delta_\rho[-2,-2], \Delta_\rho[-1,-2]; T_{IV,3}(\pi_{sc}))$
            \item $L(\Delta_\rho[-2,-2], \Delta_\rho[-2,-2], \Delta_\rho[-1,-1]; T_{IV,3}(\pi_{sc}))$
            \item $L(\Delta_\rho[-2,-2], \Delta_\rho[-2,-2], \Delta_\rho[0,-1]; \pi_{sc})$
        \end{itemize}
        \item $(0,1,2,3)$: 
        \begin{itemize}
        \item $L(\Delta_\rho[0,-1]; T_{I,1}^{1}(T_{I,1}^{2}(\pi_{sc})))$  
        \item $L(\Delta_\rho[0,-3]; \pi_{sc})$
         \item $L(\Delta_\rho[-1,-3]; T_{IV,3}(\pi_{sc}))$
            \item $L(\Delta_\rho[-3,-3]; T_{V,2}^{\pm}(T_{I,1}^{1}(T_{I,1}^{2}(\pi_{sc}))))$
            \item $L(\Delta_\rho[-2,-3], \Delta_\rho[0,-1]; \pi_{sc})$  
            \item $L(\Delta_\rho[-2,-3], \Delta_\rho[-1,-1]; T_{IV,3}(\pi_{sc}))$
            \item $L(\Delta_\rho[-3,-3], \Delta_\rho[0,-1]; T_{I,1}^{2}(\pi_{sc}))$  
            \item $L(\Delta_\rho[-3,-3], \Delta_\rho[-1,-1]; T_{IV,3}(T_{I,1}^{2}(\pi_{sc})))$

        \end{itemize}
        \item $(1,1,1,2)$: 
        \begin{itemize}
            \item $L(\Delta_\rho[-1,-1], \Delta_\rho[-1,-1]; T_{I,1}^{1}(T_{I,1}^{2}(\pi_{sc})))$
            \item $L(\Delta_\rho[-1,-1], \Delta_\rho[-1,-1], \Delta_\rho[-1,-1]; T_{I,1}^{2}(\pi_{sc}))$
            \item $L(\Delta_\rho[-1,-2], \Delta_\rho[-1,-1], \Delta_\rho[-1,-1]; \pi_{sc})$
            \item $L(\Delta_\rho[-2,-2], \Delta_\rho[-1,-1], \Delta_\rho[-1,-1], \Delta_\rho[-1,-1]; \pi_{sc})$
        \end{itemize}
        \item $(1,1,2,2)$: 
        \begin{itemize}
            \item $L(\Delta_\rho[-1,-2]; T_{I,1}^{1}(T_{I,1}^{2}(\pi_{sc})))$
            \item $L(\Delta_\rho[-2,-2], \Delta_\rho[-1,-1]; T_{I,1}^{1}(T_{I,1}^{2}(\pi_{sc})))$
            \item $L(\Delta_\rho[-1,-2], \Delta_\rho[-1,-1]; T_{I,1}^{2}(\pi_{sc}))$
            \item $L(\Delta_\rho[-1,-2], \Delta_\rho[-1,-2]; \pi_{sc})$
            \item $L(\Delta_\rho[-2,-2], \Delta_\rho[-1,-1], \Delta_\rho[-1,-1]; T_{I,1}^{2}(\pi_{sc}))$
            \item $L(\Delta_\rho[-2,-2], \Delta_\rho[-1,-2], \Delta_\rho[-1,-1]; \pi_{sc})$
            \item $L(\Delta_\rho[-2,-2], \Delta_\rho[-2,-2], \Delta_\rho[-1,-1], \Delta_\rho[-1,-1]; \pi_{sc})$
        \end{itemize}
        \item $(1,1,2,3)$: 
        \begin{itemize}
            \item $L(\Delta_\rho[-1,-1]; T_{I,1}^{1}(T_{I,1}^{3}(T_{I,1}^{2}(\pi_{sc}))))$
            \item $L(\Delta_\rho[-1,-1], \Delta_\rho[-1,-1]; T_{I,1}^{3}(T_{I,1}^{2}(\pi_{sc})))$
            \item $L(\Delta_\rho[-3,-3], \Delta_\rho[-1,-1]; T_{I,1}^{1}(T_{I,1}^{2}(\pi_{sc}))$  
            \item $L(\Delta_\rho[-2,-3], \Delta_\rho[-1,-1], \Delta_\rho[-1,-1]; \pi_{sc})$
            \item $L(\Delta_\rho[-3,-3], \Delta_\rho[-1,-1], \Delta_\rho[-1,-1]; T_{I,1}^{2}(\pi_{sc}))$
            \item $L(\Delta_\rho[-3,-3], \Delta_\rho[-1,-2], \Delta_\rho[-1,-1]; \pi_{sc})$  
            \item $L(\Delta_\rho[-3,-3], \Delta_\rho[-2,-2], \Delta_\rho[-1,-1], \Delta_\rho[-1,-1]; \pi_{sc})$
        \end{itemize}
        \item $(1,2,2,2)$: 
        \begin{itemize}
             \item $L(\Delta_\rho[-2,-2], \Delta_\rho[-2,-2]; T_{I,1}^{1}(T_{I,1}^{2}(\pi_{sc})))$
             \item $L(\Delta_\rho[-2,-2], \Delta_\rho[-1,-2]; T_{I,1}^{2}(\pi_{sc}))$
             \item $L(\Delta_\rho[-2,-2], \Delta_\rho[-2,-2], \Delta_\rho[-1,-1]; T_{I,1}^{2}(\pi_{sc}))$
             \item $L(\Delta_\rho[-2,-2], \Delta_\rho[-2,-2], \Delta_\rho[-1,-2]; \pi_{sc})$
             \item $L(\Delta_\rho[-2,-2], \Delta_\rho[-2,-2], \Delta_\rho[-2,-2], \Delta_\rho[-1,-1]; \pi_{sc})$
        \end{itemize}
        \item $(1,2,2,3)$: 
        \begin{itemize}
        \item $L(\Delta_\rho[-1,-2]; T_{I,1}^{3}(T_{I,1}^{2}(\pi_{sc})))$  
            \item $L(\Delta_\rho[-2,-3]; T_{I,1}^{2}(T_{I,1}^{1}(\pi_{sc})))$
            \item  $L(\Delta_\rho[-2,-2]; T_{I,1}^{1}(T_{I,1}^{3}(T_{I,1}^{2}(\pi_{sc}))))$
            \item $L(\Delta_\rho[-2,-3]; T_{I,1}^{1}(T_{I,1}^{2}(\pi_{sc})))$
            \item $L(\Delta_\rho[-1,-3]; T_{I,1}^{2}(\pi_{sc}))$
            \item $L(\Delta_\rho[-3,-3], \Delta_\rho[-2,-2]; T_{I,1}^{1}(T_{I,1}^{2}(\pi_{sc})))$
            \item $L(\Delta_\rho[-3,-3], \Delta_\rho[-1,-2]; T_{I,1}^{2}(\pi_{sc}))$
            \item $L(\Delta_\rho[-2,-3], \Delta_\rho[-1,-1]; T_{I,1}^{2}(\pi_{sc}))$
            \item $L(\Delta_\rho[-2,-2], \Delta_\rho[-1,-3]; \pi_{sc})$  
            \item $L(\Delta_\rho[-2,-3], \Delta_\rho[-2,-2], \Delta_\rho[-1,-1]; \pi_{sc})$
            \item $L(\Delta_\rho[-3,-3], \Delta_\rho[-2,-2], \Delta_\rho[-1,-2]; \pi_{sc})$
            \item $L(\Delta_\rho[-3,-3], \Delta_\rho[-2,-2], \Delta_\rho[-2,-2], \Delta_\rho[-1,-1]; \pi_{sc})$
        \end{itemize}
        \item $(1,2,3,3)$: 
        \begin{itemize}
            \item $L(\Delta_\rho[-3,-3]; T_{I,1}^{1}(T_{I,1}^{3}(T_{I,1}^{2}(\pi_{sc}))))$
            \item $L(\Delta_\rho[-3,-3], \Delta_\rho[-3,-3]; T_{I,1}^{1}(T_{I,1}^{2}(\pi_{sc})))$
            \item $L(\Delta_\rho[-3,-3], \Delta_\rho[-1,-1]; T_{I,1}^{3}(T_{I,1}^{2}(\pi_{sc})))$ 
            \item $L(\Delta_\rho[-3,-3], \Delta_\rho[-1,-3]; \pi_{sc})$ 
            \item $L(\Delta_\rho[-3,-3], \Delta_\rho[-3,-3], \Delta_\rho[-1,-1]; T_{I,1}^{2}(\pi_{sc}))$
            \item $L(\Delta_\rho[-3,-3], \Delta_\rho[-2,-3], \Delta_\rho[-1,-1]; \pi_{sc})$
            \item $L(\Delta_\rho[-3,-3], \Delta_\rho[-3,-3], \Delta_\rho[-1,-2]; \pi_{sc})$
            \item $L(\Delta_\rho[-3,-3], \Delta_\rho[-2,-2], \Delta_\rho[-1,-1], \Delta_\rho[-1,-1]; \pi_{sc})$
        \end{itemize}
        \item $(1,2,3,4)$: 
        \begin{itemize}
            \item $L(\Delta_\rho[-1,-4]; 
            \pi_{sc})$
            \item $L(\Delta_\rho[-4,-4]; T_{I,1}^{1}(T_{I,1}^{3}(T_{I,1}^{2}(\pi_{sc}))))$
            \item $L(\Delta_\rho[-3,-4]; T_{I,1}^{1}(T_{I,1}^{2}(\pi_{sc})))$
            \item $L(\Delta_\rho[-4,-4], \Delta_\rho[-3,-3]; T_{I,1}^{1}(T_{I,1}^{2}(\pi_{sc})))$
            \item $L(\Delta_\rho[-4,-4], \Delta_\rho[-1,-1]; T_{I,1}^{3}(T_{I,1}^{2}(\pi_{sc})))$  
            \item $L(\Delta_\rho[-4,-4], \Delta_\rho[-1,-3]; \pi_{sc})$  
            \item $L(\Delta_\rho[-3,-4], \Delta_\rho[-1,-2]; \pi_{sc})$
            \item $L(\Delta_\rho[-2,-4], \Delta_\rho[-1,-1]; \pi_{sc})$
            \item $L(\Delta_\rho[-3,-4], \Delta_\rho[-2,-2], \Delta_\rho[-1,-1]; \pi_{sc})$
            \item $L(\Delta_\rho[-4,-4], \Delta_\rho[-3,-3], \Delta\rho[-1,-1]; T_{I,1}^{2}(\pi_{sc}))$  
            \item $L(\Delta_\rho[-4,-4], \Delta_\rho[-2,-3], \Delta_\rho[-1,-1]; \pi_{sc})$
            \item $L(\Delta_\rho[-4,-4], \Delta_\rho[-3,-3], \Delta_\rho[-1,-2]; \pi_{sc})$
        \end{itemize}
        \item $(2,2,2,2)$: 
        \begin{itemize}
            \item $L(\Delta_\rho[-2,-2], \Delta_\rho[-2,-2], \Delta_\rho[-2,-2]; T_{I,1}^{2}(\pi_{sc}))$
            \item $L(\Delta_\rho[-2,-2], \Delta_\rho[-2,-2],\Delta_\rho[-2,-2],\Delta_\rho[-2,-2]; \pi_{sc})$
        \end{itemize}
        \item $(2,2,2,3)$: 
        \begin{itemize}
            \item $L(\Delta_\rho[-2,-2], \Delta_\rho[-2,-2]; T_{I,1}^{3}(T_{I,1}^{2}(\pi_{sc})))$
            \item $L(\Delta_\rho[-2,-3], \Delta_\rho[-2,-2]; T_{I,1}^{2}(\pi_{sc}))$
            \item $L(\Delta_\rho[-3,-3], \Delta_\rho[-2,-2], \Delta_\rho[-2,-2]; T_{I,1}^{2}(\pi_{sc}))$
            \item $L(\Delta_\rho[-2,-3], \Delta_\rho[-2,-2], \Delta_\rho[-2,-2]; \pi_{sc})$
            \item $L(\Delta_\rho[-3,-3], \Delta_\rho[-2,-2], \Delta_\rho[-2,-2], \Delta_\rho[-2,-2]; \pi_{sc})$
        \end{itemize}
        \item $(2,2,3,3)$:
        \begin{itemize}
        \item $L(\Delta_\rho[-2,-3]; T_{I,1}^{3}(T_{I,1}^{2}(\pi_{sc})))$ 
            \item $L(\Delta_\rho[-3,-3], \Delta_\rho[-2,-2]; T_{I,1}^{3}(T_{I,1}^{2}(\pi_{sc})))$
            \item $L(\Delta_\rho[-3,-3], \Delta_\rho[-2,-3]; T_{I,1}^{2}(\pi_{sc}))$
            \item $L(\Delta_\rho[-2,-3], \Delta_\rho[-2,-3]; \pi_{sc})$
            \item $L(\Delta_\rho[-3,-3], \Delta_\rho[-3,-3], \Delta_\rho[-2,-2]; T_{I,1}^{2}(\pi_{sc}))$
            \item $L(\Delta_\rho[-3,-3], \Delta_\rho[-2,-3], \Delta_\rho[-2,-2]; \pi_{sc})$
            \item $L(\Delta_\rho[-3,-3], \Delta_\rho[-3,-3], \Delta_\rho[-2,-2], \Delta_\rho[-2,-2]; \pi_{sc})$
        \end{itemize}
            \item $(2,2,3,4)$: 
            \begin{itemize}
                \item $L(\Delta_\rho[-2,-2]; T_{I,1}^{4}(T_{I,1}^{3}(T_{I,1}^{2}(\pi_{sc}))))$
                \item $L(\Delta_\rho[-2,-4]; T_{I,1}^{2}(\pi_{sc}))$
                \item $L(\Delta_\rho[-3,-4], \Delta_\rho[-2,-2]; T_{I,1}^{2}(\pi_{sc}))$
                \item $L(\Delta_\rho[-4,-4], \Delta_\rho[-2,-3]; T_{I,1}^{2}(\pi_{sc}))$
                \item $L(\Delta_\rho[-4,-4], \Delta_\rho[-2,-2]; T_{I,1}^{3}(T_{I,1}^{2}(\pi_{sc})))$ 
                \item $L(\Delta_\rho[-2,-4], \Delta_\rho[-2,-2]; \pi_{sc})$
                \item $L(\Delta_\rho[-4,-4], \Delta_\rho[-3,-3], \Delta_\rho[-2,-2]; T_{I,1}^{2}(\pi_{sc}))$
                \item $L(\Delta_\rho[-3,-4], \Delta_\rho[-2,-2], \Delta_\rho[-2,-2]; \pi_{sc})$
                \item $L(\Delta_\rho[-4,-4], \Delta_\rho[-2,-3], \Delta_\rho[-2,-2]; \pi_{sc})$
                \item $L(\Delta_\rho[-4,-4], \Delta_\rho[-3,-3], \Delta_\rho[-2,-2], \Delta_\rho[-2,-2]; \pi_{sc})$
            \end{itemize}
            \item $(2,3,3,3)$: 
            \begin{itemize}
                \item $L(\Delta_\rho[-3,-3], \Delta_\rho[-3,-3]; T_{I,1}^{3}(T_{I,1}^{2}(\pi_{sc})))$
                \item $L(\Delta_\rho[-3,-3], \Delta_\rho[-3,-3], \Delta_\rho[-2,-3]; \pi_{sc})$  
                \item $L(\Delta_\rho[-3,-3], \Delta_\rho[-3,-3], \Delta_\rho[-3,-3]; T_{I,1}^{2}(\pi_{sc}))$
                \item $L(\Delta_\rho[-3,-3], \Delta_\rho[-3,-3], \Delta_\rho[-3,-3], \Delta_\rho[-2,-2]; \pi_{sc})$
            \end{itemize}
            \item $(2,3,3,4)$: 
            \begin{itemize}
                \item $L(\Delta_\rho[-3,-3]; T_{I,1}^{4}(T_{I,1}^{3}(T_{I,1}^{2}(\pi_{sc}))))$
                \item $L(\Delta_\rho[-4,-4], \Delta_\rho[-3,-3]; T_{I,1}^{3}(T_{I,1}^{2}(\pi_{sc})))$
                \item $L(\Delta_\rho[-3,-4], \Delta_\rho[-3,-3]; T_{I,1}^{2}(\pi_{sc}))$
                \item $L(\Delta_\rho[-3,-4], \Delta_\rho[-2,-3]; \pi_{sc})$
                \item $L(\Delta_\rho[-2,-4], \Delta_\rho[-3,-3]; \pi_{sc})$
                \item $L(\Delta_\rho[-3,-3], \Delta_\rho[-2,-4]; \pi_{sc})$
                \item $L(\Delta_\rho[-4,-4], \Delta_\rho[-3,-3], \Delta_\rho[-3,-3]; T_{I,1}^{2}(\pi_{sc}))$
                \item $L(\Delta_\rho[-3,-4], \Delta_\rho[-3,-3], \Delta_\rho[-2,-2]; \pi_{sc})$ 
                \item $L(\Delta_\rho[-4,-4], \Delta_\rho[-3,-3], \Delta_\rho[-2,-3]; \pi_{sc})$ 
                \item $L(\Delta_\rho[-4,-4], \Delta_\rho[-3,-3], \Delta_\rho[-3,-3], \Delta_\rho[-2,-2]; \pi_{sc})$
            \end{itemize}
        \item $(2,3,4,4)$: 
        \begin{itemize}
            \item $L(\Delta_\rho[-4,-4]; T_{I,1}^{4}(T_{I,1}^{3}(T_{I,1}^{2}(\pi_{sc}))))$
            \item $L(\Delta_\rho[-4,-4], \Delta_\rho[-4,-4]; T_{I,1}^{3}(T_{I,1}^{2}(\pi_{sc})))$
            \item $L(\Delta_\rho[-4,-4], \Delta_\rho[-3,-4]; T_{I,1}^{2}(\pi_{sc}))$
            \item $L(\Delta_\rho[-4,-4], \Delta_\rho[-2,-4]; \pi_{sc})$
            \item $L(\Delta_\rho[-4,-4], \Delta_\rho[-4,-4], \Delta_\rho[-2,-3]; \pi_{sc})$  
            \item $L(\Delta_\rho[-4,-4], \Delta_\rho[-4,-4], \Delta_\rho[-3,-3]; T_{I,1}^{2}(\pi_{sc}))$
            \item $L(\Delta_\rho[-4,-4], \Delta_\rho[-3,-4], \Delta_\rho[-2,-2]; \pi_{sc})$
            \item $L(\Delta_\rho[-4,-4], \Delta_\rho[-4,-4], \Delta_\rho[-3,-3], \Delta_\rho[-2,-2]; \pi_{sc})$
        \end{itemize}
        \item $(2,3,4,5)$: 
        \begin{itemize}
            \item $L(\Delta_\rho[-5,-5]; T_{I,1}^{4}(T_{I,1}^{3}(T_{I,1}^{2}(\pi_{sc}))))$
            \item $L(\Delta_\rho[-4,-5]; T_{I,1}^{3}(T_{I,1}^{2}(\pi_{sc}))$ 
            \item $L(\Delta_\rho[-3,-5]; T_{I,1}^{2}(\pi_{sc}))$
            \item $L(\Delta_\rho[-2,-5]; \pi_{sc})$
            \item $L(\Delta_\rho[-5,-5], \Delta_\rho[-4,-4]; T_{I,1}^{3}(T_{I,1}^{2}(\pi_{sc})))$
            \item $L(\Delta_\rho[-4,-5], \Delta_\rho[-3,-3]; T_{I,1}^{2}(\pi_{sc}))$
            \item $L(\Delta_\rho[-5,-5], \Delta_\rho[-3,-4]; T_{I,1}^{2}(\pi_{sc}))$
            \item $L(\Delta_\rho[-4,-5], \Delta_\rho[-2,-3]; \pi_{sc})$
            \item $L(\Delta_\rho[-3,-5], \Delta_\rho[-2,-2]; \pi_{sc})$
            \item $L(\Delta_\rho[-5,-5], \Delta_\rho[-2,-4]; \pi_{sc})$
            \item $L(\Delta_\rho[-4,-5], \Delta_\rho[-3,-3], \Delta_\rho[-2,-2]; \pi_{sc})$
            \item $L(\Delta_\rho[-5,-5], \Delta_\rho[-4,-4], \Delta_\rho[-2,-3]' \pi_{sc})$  
            \item $L(\Delta_\rho[-5,-5], \Delta_\rho[-4,-4], \Delta_\rho[-3,-3]; T_{I,1}^{2}(\pi_{sc}))$
            \item $L(\Delta_\rho[-5,-5], \Delta_\rho[-3,-4], \Delta_\rho[-2,-2]; \pi_{sc})$
        \end{itemize}
        \end{enumerate}
    \item $(\alpha = \frac{5}{2})$: 
       \begin{enumerate}
       \item $(\frac{1}{2}, \frac{3}{2}, \frac{3}{2}, \frac{5}{2})$:
       \begin{itemize}
           \item $L(\Delta_\rho[-\frac{3}{2}, -\frac{3}{2}]; T_{I,1}^{\frac{1}{2}}(T_{I,1}^{\frac{3}{2}}(T_{I,1}^{\frac{5}{2}}(\pi_{sc}))))$  
        \item $L(\Delta_\rho[-\frac{1}{2}, -\frac{3}{2}]; T_{I,1}^{\frac{3}{2}}(T_{I,1}^{\frac{5}{2}}(\pi_{sc})))$
        \item $L(\Delta_\rho[-\frac{3}{2}, -\frac{3}{2}], \Delta_\rho[-\frac{1}{2}, -\frac{3}{2}]; T_{I,1}^{\frac{5}{2}}(\pi_{sc}))$
        \item $L(\Delta_\rho[-\frac{3}{2}, -\frac{3}{2}], \Delta_\rho[-\frac{1}{2}, -\frac{5}{2}]; \pi_{sc})$  
        \item $L(\Delta_\rho[-\frac{3}{2}, -\frac{3}{2}], \Delta_\rho[-\frac{3}{2}, -\frac{3}{2}], \Delta_\rho[-\frac{1}{2}, -\frac{1}{2}]; T_{I,1}^{\frac{5}{2}}(\pi_{sc}))$
            \item $L(\Delta_\rho[-\frac{5}{2}, -\frac{5}{2}], \Delta_\rho[-\frac{3}{2}, -\frac{3}{2}], \Delta_\rho[-\frac{1}{2}, -\frac{3}{2}]; \pi_{sc})$ 
            \item $L(\Delta_\rho[-\frac{5}{2}, -\frac{5}{2}], \Delta_\rho[-\frac{3}{2}, -\frac{3}{2}], \Delta_\rho[-\frac{3}{2}, -\frac{3}{2}], \Delta_\rho[-\frac{1}{2}, -\frac{1}{2}]; \pi_{sc})$ 
       \end{itemize}
       \item $(\frac{1}{2}, \frac{3}{2}, \frac{5}{2}, \frac{5}{2})$: 
       \begin{itemize}
            \item $L(\Delta_\rho[-\frac{1}{2}, -\frac{5}{2}]; T_{I,1}^{\frac{5}{2}}(\pi_{sc}))$ 
        \item $L(\Delta_\rho[-\frac{5}{2}, -\frac{5}{2}]; T_{I,1}^{\frac{1}{2}}(T_{I,1}^{\frac{3}{2}}(T_{I,1}^{\frac{5}{2}}(\pi_{sc}))))$  
        \item $L(\Delta_\rho[-\frac{3}{2}, -\frac{5}{2}], \Delta_\rho[-\frac{1}{2}, -\frac{1}{2}]; T_{I,1}^{\frac{5}{2}}(\pi_{sc}))$  
        \item $L(\Delta_\rho[-\frac{5}{2}, -\frac{5}{2}], \Delta_\rho[-\frac{1}{2}, -\frac{1}{2}]; T_{I,1}^{\frac{3}{2}}(T_{I,1}^{\frac{5}{2}}(\pi_{sc})))$  
        \item $L(\Delta_\rho[-\frac{5}{2}, -\frac{5}{2}], \Delta_\rho[-\frac{1}{2}, -\frac{5}{2}]; \pi_{sc})$ 
            \item $L(\Delta_\rho[-\frac{5}{2}, -\frac{5}{2}], \Delta_\rho[-\frac{5}{2}, -\frac{5}{2}], \Delta_\rho[-\frac{1}{2}, -\frac{3}{2}]; \pi_{sc})$  
            \item $L(\Delta_\rho[-\frac{5}{2}, -\frac{5}{2}], \Delta_\rho[-\frac{3}{2}, -\frac{5}{2}], \Delta_\rho[-\frac{1}{2}, -\frac{1}{2}]; \pi_{sc})$
            \item $L(\Delta_\rho[-\frac{1}{2}, -\frac{1}{2}], \Delta_\rho[-\frac{3}{2}, -\frac{3}{2}], \Delta_\rho[-\frac{5}{2}, -\frac{5}{2}], \Delta_\rho[-\frac{5}{2}, -\frac{5}{2}]; \pi_{sc})$  
       \end{itemize}
       \item $(\frac{1}{2}, \frac{3}{2}, \frac{5}{2}, \frac{7}{2})$: 
       \begin{itemize}
       \item $L(\Delta_\rho[-\frac{7}{2}, -\frac{7}{2}]; T_{I,1}^{\frac{1}{2}}(T_{I,1}^{\frac{3}{2}}(T_{I,1}^{\frac{5}{2}}(\pi_{sc}))))$ 
       \item $L(\Delta_\rho[-\frac{1}{2}, -\frac{7}{2}]; \pi_{sc})$  
       \item $L(\Delta_\rho[-\frac{5}{2}, -\frac{7}{2}], \Delta_\rho[-\frac{1}{2}, -\frac{3}{2}]; \pi_{sc})$  
       \item $L(\Delta_\rho[-\frac{7}{2}, -\frac{7}{2}], \Delta_\rho[-\frac{1}{2}, -\frac{3}{2}]; T_{I,1}^{\frac{5}{2}}(\pi_{sc}))$  
       \item $L(\Delta_\rho[-\frac{7}{2}, -\frac{7}{2}], \Delta_\rho[-\frac{1}{2}, -\frac{1}{2}]; T_{I,1}^{\frac{3}{2}}(T_{I,1}^{\frac{5}{2}}(\pi_{sc})))$  
           \item $L(\Delta_\rho[-\frac{3}{2}, \frac{7}{2}], \Delta_\rho[-\frac{1}{2}, -\frac{1}{2}]; \pi_{sc})$
           \item $L(\Delta_\rho[-\frac{7}{2}, -\frac{7}{2}],\Delta_\rho[-\frac{3}{2}, -\frac{3}{2}], \Delta_\rho[-\frac{1}{2}, -\frac{1}{2}];T_{I,1}^{\frac{5}{2}}(\pi_{sc}))$
           \item $L(\Delta_\rho[-\frac{5}{2}, -\frac{7}{2}],\Delta_\rho[-\frac{3}{2}, -\frac{3}{2}],\Delta_\rho[-\frac{1}{2}, -\frac{1}{2}]; \pi_{sc})$
       \end{itemize}
       \item $(\frac{3}{2}, \frac{3}{2}, \frac{3}{2}, \frac{5}{2})$: 
       \begin{itemize}
           \item $L(\Delta_\rho[-\frac{3}{2}, -\frac{3}{2}],\Delta_\rho[-\frac{3}{2}, -\frac{3}{2}]; T_{I,1}^{\frac{3}{2}}(T_{I,1}^{\frac{5}{2}}(\pi_{sc}))$
           \item $L(\Delta_\rho[-\frac{3}{2}, -\frac{3}{2}],\Delta_\rho[-\frac{3}{2}, -\frac{3}{2}],\Delta_\rho[-\frac{3}{2}, -\frac{3}{2}]; T_{I,1}^{\frac{5}{2}}(\pi_{sc}))$
           \item $L(\Delta_\rho[-\frac{3}{2}, -\frac{5}{2}], \Delta_\rho[-\frac{3}{2}, -\frac{3}{2}], \Delta_\rho[-\frac{3}{2}, -\frac{3}{2}]; \pi_{sc})$
           \item $L(\Delta_\rho[-\frac{5}{2},-\frac{5}{2}], \Delta_\rho[-\frac{3}{2},-\frac{3}{2}], \Delta_\rho[-\frac{3}{2},-\frac{3}{2}], \Delta_\rho[-\frac{3}{2},-\frac{3}{2}]; \pi_{sc})$
       \end{itemize}
       \item $(\frac{3}{2}, \frac{3}{2}, \frac{5}{2}, \frac{5}{2})$:
       \begin{itemize}
       \item $L(\Delta_\rho[-\frac{3}{2}, -\frac{5}{2}]; T_{I,1}^{\frac{3}{2}}(T_{I,1}^{\frac{5}{2}}(\pi_{sc})))$
           \item $L(\Delta_\rho[-\frac{3}{2}, -\frac{5}{2}],\Delta_\rho[-\frac{3}{2},-\frac{3}{2}]; T_{I,1}^{\frac{5}{2}}(\pi_{sc}))$
           \item $L(\Delta_\rho[-\frac{3}{2}, -\frac{5}{2}], \Delta_\rho[-\frac{3}{2}, -\frac{5}{2}]; \pi_{sc})$
           \item $L(\Delta_\rho[-\frac{5}{2}, -\frac{5}{2}], \Delta_\rho[-\frac{3}{2}, -\frac{3}{2}]; T_{I,1}^{\frac{5}{2}}(T_{I,1}^{\frac{3}{2}}(\pi_{sc})))$
           \item $L(\Delta_\rho[-\frac{5}{2}, -\frac{5}{2}], \Delta_\rho[-\frac{3}{2}, -\frac{3}{2}], \Delta_\rho[-\frac{3}{2}, -\frac{3}{2}]; T_{I,1}^{\frac{5}{2}}(\pi_{sc}))$
           \item $L(\Delta_\rho[-\frac{5}{2}, -\frac{5}{2}], \Delta_\rho[-\frac{3}{2}, -\frac{5}{2}], \Delta_\rho[-\frac{3}{2}, -\frac{3}{2}]; \pi_{sc})$
           \item $L(\Delta_\rho[-\frac{5}{2},-\frac{5}{2}], \Delta_\rho[-\frac{5}{2},-\frac{5}{2}], \Delta_\rho[-\frac{3}{2},-\frac{3}{2}], \Delta_\rho[-\frac{3}{2},-\frac{3}{2}]; \pi_{sc})$
       \end{itemize}
       \item $(\frac{3}{2}, \frac{3}{2}, \frac{5}{2}, \frac{7}{2})$: 
       \begin{itemize}
           \item $L(\Delta_\rho[-\frac{3}{2}, -\frac{3}{2}]; T_{I,1}^{\frac{3}{2}}(T_{I,1}^{\frac{7}{2}}(T_{I,1}^{\frac{5}{2}}(\pi_{sc}))))$
           \item $L(\Delta_\rho[-\frac{3}{2}, -\frac{3}{2}];T_{I,1}^{\frac{3}{2}}(T_{I,1}^{\frac{7}{2}}(T_{I,1}^{\frac{5}{2}}(\pi_{sc}))))$
           \item $L(\Delta_\rho[-\frac{7}{2}, -\frac{7}{2}], \Delta_\rho[-\frac{3}{2}, -\frac{3}{2}]; T_{I,1}^{\frac{3}{2}}(T_{I,1}^{\frac{5}{2}}(\pi_{sc})))$
           \item $L(\Delta_\rho[-\frac{3}{2}, -\frac{3}{2}],\Delta_\rho[-\frac{3}{2}, -\frac{3}{2}]; T_{I,1}^{\frac{7}{2}}(T_{I,1}^{\frac{5}{2}}(\pi_{sc})))$  
           \item $L(\Delta_\rho[-\frac{3}{2}, \frac{7}{2}], \Delta_\rho[-\frac{3}{2}, -\frac{3}{2}]; \pi_{sc})$
           \item $L(\Delta_\rho[-\frac{7}{2}, -\frac{7}{2}], \Delta_\rho[-\frac{3}{2}, -\frac{3}{2}], \Delta_\rho[-\frac{3}{2}, -\frac{3}{2}]; T_{I,1}^{\frac{5}{2}}(\pi_{sc}))$
           \item $L(\Delta_\rho[-\frac{5}{2}, -\frac{7}{2}],\Delta_\rho[-\frac{3}{2}, -\frac{3}{2}],\Delta_\rho[-\frac{3}{2}, -\frac{3}{2}]; \pi_{sc})$
           \item $L(\Delta_\rho[-\frac{7}{2}, -\frac{7}{2}], \Delta_\rho[-\frac{5}{2},-\frac{3}{2}], \Delta_\rho[-\frac{3}{2}, -\frac{3}{2}]; \pi_{sc})$ 
           \item $L(\Delta_\rho[-\frac{7}{2},-\frac{7}{2}], \Delta_\rho[-\frac{5}{2},-\frac{5}{2}], \Delta_\rho[-\frac{3}{2},-\frac{3}{2}], \Delta_\rho[-\frac{3}{2},-\frac{3}{2}]; \pi_{sc})$
       \end{itemize}
       \item $(\frac{3}{2}, \frac{5}{2}, \frac{5}{2}, \frac{5}{2})$: 
       \begin{itemize}
           \item $L(\Delta_\rho[-\frac{5}{2}, -\frac{5}{2}], \Delta_\rho[-\frac{3}{2},-\frac{5}{2}]; T_{I,1}^{\frac{5}{2}}(\pi_{sc}))$
           \item $L(\Delta_\rho[-\frac{5}{2}, -\frac{5}{2}], \Delta_\rho[-\frac{5}{2}, -\frac{5}{2}]; T_{I,1}^{\frac{3}{2}}(T_{I,1}^{\frac{5}{2}}(\pi_{sc})))$ 
           \item $L(\Delta_\rho[-\frac{5}{2}, -\frac{5}{2}],\Delta_\rho[-\frac{5}{2}, -\frac{5}{2}], \Delta_\rho[-\frac{3}{2}, -\frac{3}{2}]; T_{I,1}^{\frac{5}{2}}(\pi_{sc}))$
           \item $L(\Delta_\rho[-\frac{5}{2}, -\frac{5}{2}], \Delta_\rho[-\frac{5}{2}, -\frac{5}{2}], \Delta_\rho[-\frac{3}{2}, -\frac{5}{2}]; \pi_{sc})$
           \item $L(\Delta_\rho[-\frac{5}{2},-\frac{5}{2}], \Delta_\rho[-\frac{5}{2},-\frac{5}{2}], \Delta_\rho[-\frac{5}{2},-\frac{5}{2}], \Delta_\rho[-\frac{3}{2},-\frac{3}{2}]; \pi_{sc})$
       \end{itemize}
       \item $(\frac{3}{2}, \frac{5}{2}, \frac{5}{2}, \frac{7}{2})$: 
       \begin{itemize}
           \item $L(\Delta_\rho[-\frac{5}{2}, -\frac{5}{2}]; T_{I,1}^{\frac{3}{2}}(T_{I,1}^{\frac{7}{2}}(T_{I,1}^{\frac{5}{2}}(\pi_{sc}))))$
           \item $L(\Delta_\rho[-\frac{5}{2}, -\frac{5}{2}];T_{I,1}^{\frac{3}{2}}(T_{I,1}^{\frac{7}{2}}(T_{I,1}^{\frac{5}{2}}(\pi_{sc}))))$
           \item $L(\Delta_\rho[-\frac{3}{2}, -\frac{7}{2}]; T_{I,1}^{\frac{5}{2}}(\pi_{sc}))$
           \item $L(\Delta_\rho[-\frac{5}{2}, -\frac{5}{2}], \Delta_\rho[-\frac{3}{2}, -\frac{3}{2}]; T_{I,1}^{\frac{5}{2}}(T_{I,1}^{\frac{7}{2}}(\pi_{sc}))$
           \item $L(\Delta_\rho[-\frac{5}{2}, -\frac{7}{2}],\Delta_\rho[-\frac{3}{2},-\frac{3}{2}]; T_{I,1}^{\frac{5}{2}}(\pi_{sc}))$
           \item $L(\Delta_\rho[-\frac{7}{2}, -\frac{7}{2}], \Delta_\rho[-\frac{3}{2},-\frac{5}{2}]; T_{I,1}^{\frac{5}{2}}(\pi_{sc}))$
           
           \item $L(\Delta_\rho[-\frac{3}{2}, \frac{7}{2}], \Delta_\rho[-\frac{5}{2}, -\frac{5}{2}]; \pi_{sc})$
           \item $L(\Delta_\rho[-\frac{5}{2}, -\frac{5}{2}], \Delta_\rho[-\frac{3}{2}, -\frac{7}{2}]; \pi_{sc})$
           \item $L(\Delta_\rho[-\frac{5}{2}, -\frac{7}{2}],\Delta_\rho[-\frac{5}{2}, -\frac{5}{2}],\Delta_\rho[-\frac{3}{2}, -\frac{3}{2}]; T_{I,1}^{\frac{5}{2}}(\pi_{sc}))$
           \item $L(\Delta_\rho[-\frac{7}{2}, -\frac{7}{2}], \Delta_\rho[-\frac{5}{2}, -\frac{5}{2}], \Delta_\rho[-\frac{3}{2}, -\frac{5}{2}]; \pi_{sc})$
           \item $L(\Delta_\rho[-\frac{7}{2},-\frac{7}{2}], \Delta_\rho[-\frac{5}{2},-\frac{5}{2}], \Delta_\rho[-\frac{5}{2},-\frac{5}{2}], \Delta_\rho[-\frac{3}{2},-\frac{3}{2}]; \pi_{sc})$
       \end{itemize}
       \item $(\frac{3}{2}, \frac{5}{2}, \frac{7}{2}, \frac{7}{2})$: 
       \begin{itemize}
           \item $L(\Delta_\rho[-\frac{7}{2}, -\frac{7}{2}]; T_{I,1}^{\frac{3}{2}}(T_{I,1}^{\frac{7}{2}}(T_{I,1}^{\frac{5}{2}}(\pi_{sc}))))$
           \item $L(\Delta_\rho[-\frac{7}{2}, -\frac{7}{2}];T_{I,1}^{\frac{3}{2}}(T_{I,1}^{\frac{7}{2}}(T_{I,1}^{\frac{5}{2}}(\pi_{sc}))))$
           \item $L(\Delta_\rho[-\frac{7}{2},-\frac{7}{2}], \Delta_\rho[-\frac{7}{2},-\frac{7}{2}]; T_{I,1}^{\frac{3}{2}}(T_{I,1}^{\frac{5}{2}}(\pi_{sc})))$
           \item $L(\Delta_\rho[-\frac{7}{2}, \frac{7}{2}], \Delta_\rho[-\frac{3}{2}, -\frac{7}{2}]; \pi_{sc})$
           
           \item $L(\Delta_\rho[-\frac{7}{2}, -\frac{7}{2}],\Delta_\rho[-\frac{7}{2}, -\frac{7}{2}], \Delta_\rho[-\frac{3}{2}, -\frac{3}{2}]; T_{I,1}^{\frac{5}{2}}(\pi_{sc}))$
           \item $L(\Delta_\rho[-\frac{7}{2}, -\frac{7}{2}], \Delta_\rho[-\frac{5}{2}, -\frac{7}{2}], \Delta_\rho[-\frac{3}{2}, -\frac{3}{2}]; \pi_{sc})$
           \item $L(\Delta_\rho[-\frac{7}{2}, -\frac{7}{2}], \Delta_\rho[-\frac{7}{2}, -\frac{7}{2}], \Delta_\rho[-\frac{3}{2}, -\frac{5}{2}]; \pi_{sc})$
          
           \item $L(\Delta_\rho[-\frac{7}{2},-\frac{7}{2}], \Delta_\rho[-\frac{7}{2},-\frac{7}{2}], \Delta_\rho[-\frac{5}{2},-\frac{5}{2}], \Delta_\rho[-\frac{3}{2},-\frac{3}{2}]; \pi_{sc})$
       \end{itemize}
       \item $(\frac{3}{2}, \frac{5}{2}, \frac{7}{2}, \frac{9}{2})$:
       \begin{itemize}
           \item $L(\Delta_\rho[-\frac{9}{2}, -\frac{9}{2}];T_{I,1}^{\frac{3}{2}}(T_{I,1}^{\frac{7}{2}}(T_{I,1}^{\frac{5}{2}}(\pi_{sc}))))$
           \item $L(\Delta_\rho[-\frac{7}{2}, -\frac{9}{2}]; T_{I,1}^{\frac{3}{2}}(T_{I,1}^{\frac{5}{2}}(\pi_{sc})))$
           \item $L(\Delta_\rho[-\frac{3}{2}, -\frac{9}{2}]; \pi_{sc})$
           \item $L(\Delta_\rho[-\frac{7}{2}, -\frac{9}{2}], \Delta_\rho[-\frac{3}{2}, -\frac{5}{2}]; \pi_{sc})$
           \item $L(\Delta_\rho[-\frac{5}{2}, \frac{9}{2}], \Delta_\rho[-\frac{3}{2}, -\frac{3}{2}]; \pi_{sc})$
           \item $L(\Delta_\rho[-\frac{9}{2}, \frac{9}{2}], \Delta_\rho[-\frac{3}{2}, -\frac{7}{2}]; \pi_{sc})$
           \item $L(\Delta_\rho[-\frac{9}{2}, -\frac{9}{2}], \Delta_\rho[-\frac{3}{2}, -\frac{3}{2}]; T_{I,1}^{\frac{7}{2}}(T_{I,1}^{\frac{5}{2}}(\pi_{sc})))$  
           \item $L(\Delta_\rho[-\frac{7}{2}, -\frac{9}{2}],\Delta_\rho[-\frac{3}{2},-\frac{3}{2}]; T_{I,1}^{\frac{5}{2}}(\pi_{sc}))$
           \item $L(\Delta_\rho[-\frac{9}{2}, -\frac{9}{2}],\Delta_\rho[-\frac{7}{2}, -\frac{7}{2}]; T_{I,1}^{\frac{3}{2}}(T_{I,1}^{\frac{5}{2}}(\pi_{sc})))$  
            \item $L(\Delta_\rho[-\frac{9}{2}, -\frac{9}{2}], \Delta_\rho[-\frac{7}{2}, -\frac{7}{2}], \Delta_\rho[-\frac{3}{2}, -\frac{3}{2}]; T_{I,1}^{\frac{5}{2}}(\pi_{sc}))$ 
           \item $L(\Delta_\rho[-\frac{7}{2}, -\frac{9}{2}], \Delta_\rho[-\frac{5}{2}, -\frac{5}{2}], \Delta_\rho[-\frac{3}{2}, -\frac{3}{2}]; \pi_{sc})$
           \item $L(\Delta_\rho[-\frac{9}{2}, -\frac{9}{2}], \Delta_\rho[-\frac{5}{2}, -\frac{7}{2}], \Delta_\rho[-\frac{3}{2}, -\frac{3}{2}]; \pi_{sc})$
       \end{itemize}
       \item $(\frac{5}{2}, \frac{5}{2}, \frac{5}{2}, \frac{5}{2})$: 
       \begin{itemize}
           \item  $L(\Delta_\rho[-\frac{5}{2}, -\frac{5}{2}],\Delta_\rho[-\frac{5}{2}, -\frac{5}{2}],\Delta_\rho[-\frac{5}{2}, -\frac{5}{2}]; T_{I,1}^{\frac{5}{2}}(\pi_{sc}))$
           \item $L(\Delta_\rho[-\frac{5}{2}, -\frac{5}{2}], \Delta_\rho[-\frac{5}{2}, -\frac{5}{2} ], \Delta_\rho[-\frac{5}{2}, -\frac{5}{2} ], \Delta_\rho[-\frac{5}{2}, -\frac{5}{2} ]; \pi_{sc})$
       \end{itemize}
       \item $(\frac{5}{2}, \frac{5}{2}, \frac{5}{2}, \frac{7}{2})$: 
       \begin{itemize}
           \item  $L(\Delta_\rho[-\frac{5}{2}, -\frac{7}{2}],\Delta_\rho[-\frac{5}{2},-\frac{5}{2}]; T_{I,1}^{\frac{5}{2}}(\pi_{sc}))$
           \item $L(\Delta_\rho[-\frac{5}{2},-\frac{5}{2}],\Delta_\rho[-\frac{5}{2},-\frac{5}{2}]; T_{I,1}^{\frac{7}{2}}(T_{I,1}^{\frac{5}{2}}(\pi_{sc})))$  
           \item $L(\Delta_\rho[-\frac{7}{2}, -\frac{7}{2}], \Delta_\rho[-\frac{5}{2}, -\frac{5}{2}], \Delta_\rho[-\frac{5}{2}, -\frac{5}{2}]; T_{I,1}^{\frac{5}{2}}(\pi_{sc}))$
           \item $L(\Delta_\rho[-\frac{5}{2}, -\frac{7}{2}], \Delta_\rho[-\frac{5}{2}, -\frac{5}{2}], \Delta_\rho[-\frac{5}{2}, -\frac{5}{2}]; \pi_{sc})$
           \item $L(\Delta_\rho[-\frac{7}{2},-\frac{7}{2}], \Delta_\rho[-\frac{5}{2},-\frac{5}{2}], \Delta_\rho[-\frac{5}{2},-\frac{5}{2}], \Delta_\rho[-\frac{5}{2},-\frac{5}{2}]; \pi_{sc})$
       \end{itemize}
       \item $(\frac{5}{2}, \frac{5}{2}, \frac{7}{2}, \frac{7}{2})$:
       \begin{itemize}
           \item $L(\Delta_\rho[-\frac{7}{2}, -\frac{7}{2}], \Delta_\rho[-\frac{5}{2}, -\frac{5}{2}]; T_{I,1}^{\frac{7}{2}}(T_{I,1}^{\frac{5}{2}}(\pi_{sc})))$
           \item $L(\Delta_\rho[-\frac{7}{2}, -\frac{7}{2}], \Delta_\rho[-\frac{5}{2}, -\frac{5}{2}]; T_{I,1}^{\frac{5}{2}}(T_{I,1}^{\frac{7}{2}}(\pi_{sc}))$
           \item $L(\Delta_\rho[-\frac{7}{2}, -\frac{7}{2}], \Delta_\rho[-\frac{5}{2},-\frac{7}{2}]; T_{I,1}^{\frac{5}{2}}(\pi_{sc}))$
           \item $L(\Delta_\rho[-\frac{5}{2}, -\frac{7}{2}], \Delta_\rho[-\frac{5}{2}, -\frac{7}{2}]; \pi_{sc})$
           \item $L(\Delta_\rho[-\frac{7}{2}, -\frac{7}{2}],\Delta_\rho[-\frac{7}{2}, -\frac{7}{2}], \Delta_\rho[-\frac{5}{2}, -\frac{5}{2}]; T_{I,1}^{\frac{5}{2}}(\pi_{sc}))$
           \item $L(\Delta_\rho[-\frac{7}{2},-\frac{7}{2}], \Delta_\rho[-\frac{7}{2},-\frac{7}{2}], \Delta_\rho[-\frac{5}{2},-\frac{5}{2}], \Delta_\rho[-\frac{5}{2},-\frac{5}{2}]; \pi_{sc})$
            \item $L(\Delta_\rho[-\frac{7}{2},-\frac{7}{2}], \Delta_\rho[-\frac{7}{2},-\frac{7}{2}], \Delta_\rho[-\frac{5}{2},-\frac{5}{2}], \Delta_\rho[-\frac{5}{2},-\frac{5}{2}]; \pi_{sc})$
       \end{itemize}
        \item $(\frac{5}{2}, \frac{5}{2}, \frac{7}{2}, \frac{9}{2})$: 
        \begin{itemize}
            \item $L(\Delta_\rho[-\frac{5}{2}, -\frac{5}{2}]; T_{I,1}^{\frac{9}{2}}(T_{I,1}^{\frac{7}{2}}(T_{I,1}^{\frac{5}{2}}(\pi_{sc}))))$
            \item $L(\Delta_\rho[-\frac{5}{2}, -\frac{9}{2}]; T_{I,1}^{\frac{5}{2}}(\pi_{sc}))$
            \item $L(\Delta_\rho[-\frac{7}{2}, -\frac{9}{2}],\Delta_\rho[-\frac{5}{2},-\frac{5}{2}]; T_{I,1}^{\frac{5}{2}}(\pi_{sc}))$
            \item $L(\Delta_\rho[-\frac{9}{2}, -\frac{9}{2}], \Delta_\rho[-\frac{5}{2},-\frac{7}{2}]; T_{I,1}^{\frac{5}{2}}(\pi_{sc}))$
            \item $L(\Delta_\rho[-\frac{5}{2}, \frac{9}{2}], \Delta_\rho[-\frac{5}{2}, -\frac{5}{2}]; \pi_{sc})$
            \item $L(\Delta_\rho[-\frac{9}{2}, -\frac{9}{2}], \Delta_\rho[-\frac{5}{2}, -\frac{5}{2}]; T_{I,1}^{\frac{7}{2}}(T_{I,1}^{\frac{5}{2}}(\pi_{sc})))$  
            \item $L(\Delta_\rho[-\frac{9}{2},-\frac{9}{2}],\Delta_\rho[-\frac{7}{2}, -\frac{7}{2}], \Delta_\rho[-\frac{5}{2},-\frac{5}{2}]; T_{I,1}^{\frac{5}{2}}(\pi_{sc}))$
            \item $L(\Delta_\rho[-\frac{7}{2}, -\frac{9}{2}], \Delta_\rho[-\frac{5}{2}, -\frac{5}{2}], \Delta_\rho[-\frac{5}{2}, -\frac{5}{2}]; \pi_{sc})$
            \item $L(\Delta_\rho[-\frac{9}{2}, -\frac{9}{2}], \Delta_\rho[-\frac{5}{2}, -\frac{7}{2}], \Delta_\rho[-\frac{5}{2}, -\frac{5}{2}]; \pi_{sc})$
            \item $L(\Delta_\rho[-\frac{9}{2},-\frac{9}{2}], \Delta_\rho[-\frac{7}{2},-\frac{7}{2}], \Delta_\rho[-\frac{5}{2},-\frac{5}{2}], \Delta_\rho[-\frac{5}{2},-\frac{5}{2}]; \pi_{sc})$
        \end{itemize}
        \item $(\frac{5}{2}, \frac{7}{2}, \frac{7}{2}, \frac{7}{2})$
        \begin{itemize}
        \item $L(\Delta_\rho[-\frac{7}{2}, -\frac{7}{2}], \Delta_\rho[-\frac{7}{2}, -\frac{7}{2}]; T_{I,1}^{\frac{7}{2}}(T_{I,1}^{\frac{5}{2}}(\pi_{sc})))$ 
            \item $L(\Delta_\rho[-\frac{7}{2}, -\frac{7}{2}],\Delta_\rho[-\frac{7}{2}, -\frac{7}{2}],\Delta_\rho[-\frac{7}{2}, -\frac{7}{2}]; T_{I,1}^{\frac{5}{2}}(\pi_{sc}))$
            \item $L(\Delta_\rho[-\frac{7}{2}, -\frac{7}{2}], \Delta_\rho[-\frac{7}{2}, -\frac{7}{2}], \Delta_\rho[-\frac{5}{2}, -\frac{7}{2}]; \pi_{sc})$
             \item $L(\Delta_\rho[-\frac{7}{2},-\frac{7}{2}], \Delta_\rho[-\frac{7}{2},-\frac{7}{2}], \Delta_\rho[-\frac{7}{2},-\frac{7}{2}], \Delta_\rho[-\frac{5}{2},-\frac{5}{2}]; \pi_{sc})$
        \end{itemize}
        \item $(\frac{5}{2}, \frac{7}{2}, \frac{7}{2}, \frac{9}{2})$: 
        \begin{itemize}
            \item $L(\Delta_\rho[-\frac{7}{2}, -\frac{7}{2}]; T_{I,1}^{\frac{9}{2}}(T_{I,1}^{\frac{7}{2}}(T_{I,1}^{\frac{5}{2}}(\pi_{sc}))))$
            \item $L(\Delta_\rho[-\frac{9}{2}, -\frac{9}{2}], \Delta_\rho[-\frac{7}{2}, -\frac{7}{2}]; T_{I,1}^{\frac{7}{2}}(T_{I,1}^{\frac{5}{2}}(\pi_{sc})))$
            \item $L(\Delta_\rho[-\frac{9}{2}, -\frac{9}{2}], \Delta_\rho[-\frac{7}{2}, -\frac{7}{2}]; T_{I,1}^{\frac{5}{2}}(T_{I,1}^{\frac{7}{2}}(\pi_{sc}))$
            \item $L(\Delta_\rho[-\frac{7}{2}, -\frac{9}{2}],\Delta_\rho[-\frac{7}{2},-\frac{7}{2}]; T_{I,1}^{\frac{5}{2}}(\pi_{sc}))$
            \item $L(\Delta_\rho[-\frac{7}{2}, -\frac{9}{2}], \Delta_\rho[-\frac{5}{2}, -\frac{7}{2}]; \pi_{sc})$
            \item $L(\Delta_\rho[-\frac{5}{2}, \frac{9}{2}], \Delta_\rho[-\frac{7}{2}, -\frac{7}{2}]; \pi_{sc})$
            \item $L(\Delta_\rho[-\frac{7}{2}, \frac{7}{2}], \Delta_\rho[-\frac{5}{2}, -\frac{9}{2}]; \pi_{sc})$
            \item $L(\Delta_\rho[-\frac{9}{2}, -\frac{9}{2}], \Delta_\rho[-\frac{7}{2}, -\frac{7}{2}], \Delta_\rho[-\frac{7}{2}, -\frac{7}{2}]; T_{I,1}^{\frac{5}{2}}(\pi_{sc}))$
            \item $L(\Delta_\rho[-\frac{9}{2}, -\frac{9}{2}], \Delta_\rho[-\frac{7}{2}, -\frac{7}{2}], \Delta_\rho[-\frac{5}{2}, -\frac{7}{2}]; \pi_{sc})$
            \item $L(\Delta_\rho[-\frac{9}{2},-\frac{9}{2}], \Delta_\rho[-\frac{7}{2},-\frac{7}{2}], \Delta_\rho[-\frac{7}{2},-\frac{7}{2}], \Delta_\rho[-\frac{5}{2},-\frac{5}{2}]; \pi_{sc})$
        \end{itemize}
        \item $(\frac{5}{2}, \frac{7}{2}, \frac{9}{2}, \frac{9}{2})$: 
        \begin{itemize}
            \item $L(\Delta_\rho[-\frac{9}{2}, -\frac{9}{2}]; T_{I,1}^{\frac{9}{2}}(T_{I,1}^{\frac{7}{2}}(T_{I,1}^{\frac{5}{2}}(\pi_{sc}))))$
            \item $L(\Delta_\rho[-\frac{9}{2}, -\frac{9}{2}], \Delta_\rho[-\frac{9}{2}, -\frac{9}{2}]; T_{I,1}^{\frac{7}{2}}(T_{I,1}^{\frac{5}{2}}(\pi_{sc})))$  
            \item $L(\Delta_\rho[-\frac{9}{2}, -\frac{9}{2}], \Delta_\rho[-\frac{7}{2},-\frac{9}{2}]; T_{I,1}^{\frac{5}{2}}(\pi_{sc}))$
            \item $L(\Delta_\rho[-\frac{9}{2}, \frac{9}{2}], \Delta_\rho[-\frac{5}{2}, -\frac{9}{2}]; \pi_{sc})$
            \item $L(\Delta_\rho[-\frac{9}{2}, -\frac{9}{2}],\Delta_\rho[-\frac{9}{2}, -\frac{9}{2}], \Delta_\rho[-\frac{7}{2}, -\frac{7}{2}]; T_{I,1}^{\frac{5}{2}}(\pi_{sc}))$
            \item $L(\Delta_\rho[-\frac{9}{2}, -\frac{9}{2}], \Delta_\rho[-\frac{7}{2}, -\frac{9}{2}], \Delta_\rho[-\frac{5}{2}, -\frac{5}{2}]; \pi_{sc})$
            \item $L(\Delta_\rho[-\frac{9}{2}, -\frac{9}{2}], \Delta_\rho[-\frac{9}{2}, -\frac{9}{2}], \Delta_\rho[-\frac{5}{2}, -\frac{7}{2}]; \pi_{sc})$
            \item $L(\Delta_\rho[-\frac{9}{2},-\frac{9}{2}], \Delta_\rho[-\frac{9}{2},-\frac{9}{2}], \Delta_\rho[-\frac{7}{2},-\frac{7}{2}], \Delta_\rho[-\frac{5}{2},-\frac{5}{2}]; \pi_{sc})$
        \end{itemize}
        \item $(\frac{5}{2}, \frac{7}{2}, \frac{9}{2}, \frac{11}{2})$: 
        \begin{itemize}
            \item $L(\Delta_\rho[-\frac{11}{2}, -\frac{11}{2}]; T_{I,1}^{\frac{9}{2}}(T_{I,1}^{\frac{7}{2}}(T_{I,1}^{\frac{5}{2}}(\pi_{sc}))))$
            \item $L(\Delta_\rho[-\frac{9}{2}, -\frac{11}{2}]; T_{I,1}^{\frac{7}{2}}(T_{I,1}^{\frac{5}{2}}(\pi_{sc})))$  
            \item $L(\Delta_\rho[-\frac{7}{2}, -\frac{11}{2}]; T_{I,1}^{\frac{5}{2}}(\pi_{sc}))$
            \item $L(\Delta_\rho[-\frac{5}{2}, -\frac{11}{2}]; \pi_{sc})$
            \item $L(\Delta_\rho[-\frac{11}{2}, -\frac{11}{2}], \Delta_\rho[-\frac{9}{2}, -\frac{9}{2}]; T_{I,1}^{\frac{7}{2}}(T_{I,1}^{\frac{5}{2}}(\pi_{sc})))$
            \item $L(\Delta_\rho[-\frac{9}{2}, -\frac{11}{2}],\Delta_\rho[-\frac{7}{2},-\frac{7}{2}]; T_{I,1}^{\frac{5}{2}}(\pi_{sc}))$
            \item $L(\Delta_\rho[-\frac{11}{2}, -\frac{11}{2}], \Delta_\rho[-\frac{7}{2},-\frac{9}{2}]; T_{I,1}^{\frac{5}{2}}(\pi_{sc}))$
            \item $L(\Delta_\rho[-\frac{9}{2}, -\frac{11}{2}], \Delta_\rho[-\frac{5}{2}, -\frac{7}{2}]; \pi_{sc})$
            \item $L(\Delta_\rho[-\frac{7}{2}, \frac{11}{2}], \Delta_\rho[-\frac{5}{2}, -\frac{5}{2}]; \pi_{sc})$
            \item $L(\Delta_\rho[-\frac{11}{2}, \frac{11}{2}], \Delta_\rho[-\frac{5}{2}, -\frac{9}{2}]; \pi_{sc})$
            \item $L(\Delta_\rho[-\frac{11}{2},-\frac{11}{2}],\Delta_\rho[-\frac{9}{2}, -\frac{9}{2}], \Delta_\rho[-\frac{7}{2},-\frac{7}{2}]; T_{I,1}^{\frac{5}{2}}(\pi_{sc}))$
            \item $L(\Delta_\rho[-\frac{9}{2}, -\frac{11}{2}], \Delta_\rho[-\frac{7}{2}, -\frac{7}{2}], \Delta_\rho[-\frac{5}{2}, -\frac{5}{2}]; \pi_{sc})$
            \item $L(\Delta_\rho[-\frac{11}{2}, -\frac{11}{2}], \Delta_\rho[-\frac{7}{2}, -\frac{9}{2}], \Delta_\rho[-\frac{5}{2}, -\frac{5}{2}]; \pi_{sc})$
            \item $L(\Delta_\rho[-\frac{11}{2}, -\frac{11}{2}], \Delta_\rho[-\frac{9}{2}, -\frac{9}{2}], \Delta_\rho[-\frac{5}{2}, -\frac{7}{2}]; \pi_{sc})$
        \end{itemize}
       \end{enumerate}
       \item $(\alpha =3)$: 
       \begin{enumerate}
       \item $(1,1,2,3)$: 
       \begin{itemize}
       \item $L(\Delta_\rho[-1,-1];T_{I,1}^{1}(T_{I,1}^{2}(T_{I,1}^{3}(\pi_{sc}))))$
           \item $L(\Delta_\rho[-1,-1], \Delta_\rho[-1,-1]; T_{I,1}^{2}(T_{I,1}^{3}(\pi_{sc})))$
           \item $L(\Delta_\rho[-1,-2], \Delta_\rho[-1,-1]; T_{I,1}^{3}(\pi_{sc}))$
           \item $L(\Delta_\rho[-1,-3], \Delta_\rho[-1,-1];\pi_{sc})$
           \item $L(\Delta_\rho[-2,-2], \Delta_\rho[-1,-1], \Delta_\rho[-1,-1]; T_{I,1}^{3}(\pi_{sc}))$
           \item $L(\Delta_\rho[-2,-3], \Delta_\rho[-1,-1], \Delta_\rho[-1,-1]; \pi_{sc})$
           \item $L(\Delta_\rho[-3,-3], \Delta_\rho[-1,-2], \Delta_\rho[-1,-1]; \pi_{sc})$
           \item $L(\Delta_\rho[-3,-3], \Delta_\rho[-2,-2], \Delta_\rho[-1,-1], \Delta_\rho[-1,-1]; \pi_{sc})$
       \end{itemize}
       \item $(1,2,2,3)$: 
       \begin{itemize}
       \item $L(\Delta_\rho[-1,-2]; T_{I,1}^{2}(T_{I,1}^{3}(\pi_{sc})))$
           \item $L(\Delta_\rho[-2,-2]; T_{I,1}^{1}(T_{I,1}^{2}(T_{I,1}^{3}(\pi_{sc}))))$
           \item $L(\Delta_\rho[-2,-2], \Delta_\rho[-1,-2]; T_{I,1}^{3}(\pi_{sc}))$ 
           \item $L(\Delta_\rho[-2,-2],\Delta_\rho[-1,-3];\pi_{sc})$  
           \item $L(\Delta_\rho[-2,-2], \Delta_\rho[-2,-2], \Delta_\rho[-1,-1]; T_{I,1}^{3}(\pi_{sc}))$
           \item $L(\Delta_\rho[-2,-3], \Delta_\rho[-2,-2], \Delta_\rho[-1,-1]; \pi_{sc})$
           \item $L(\Delta_\rho[-3,-3],\Delta_\rho[-2,-2], \Delta_\rho[-1,-2]; \pi_{sc})$
           \item $L(\Delta_\rho[-3,-3], \Delta_\rho[-2,-2], \Delta_\rho[-2,-2], \Delta_\rho[-1,-1]; \pi_{sc})$
       \end{itemize}
       \item $(1,2,3,3)$: 
       \begin{itemize}
            \item $L(\Delta_\rho[-3,-3]; T_{I,1}^{1}(T_{I,1}^{2}(T_{I,1}^{3}(\pi_{sc}))))$
            \item $L(\Delta_\rho[-1,-3]; T_{I,1}^{3}(\pi_{sc}))$
            \item $L(\Delta_\rho[-3,-3],\Delta_\rho[-1,-1];T_{I,1}^{2}(T_{I,1}^{3}(\pi_{sc})))$  
            \item $L(\Delta_\rho[-2,-3], \Delta_\rho[-1,-1]; T_{I,1}^{3}(\pi_{sc}))$
            \item $L(\Delta_\rho[-3,-3],\Delta_\rho[-1,-3]; \pi_{sc})$
            \item $L(\Delta_\rho[-3,-3], \Delta_\rho[-2,-3], \Delta_\rho[-1,-1]; \pi_{sc})$
            \item $L(\Delta_\rho[-3,-3],\Delta_\rho[-3,-3], \Delta_\rho[-1,-2]; \pi_{sc})$
            \item $L(\Delta_\rho[-3,-3], \Delta_\rho[-3,-3], \Delta_\rho[-2,-2], \Delta_\rho[-1,-1]; \pi_{sc})$
       \end{itemize}
       \item $(1,2,3,4)$: 
       \begin{itemize}
           \item $L(\Delta_\rho[-4,-4]; T_{I,1}^{1}(T_{I,1}^{2}(T_{I,1}^{3}(\pi_{sc}))))$
           \item $L(\Delta_\rho[-1,-4]; \pi_{sc})$
           \item $L(\Delta_\rho[-3,-4], \Delta_\rho[-1,-2]; \pi_{sc})$
           \item $L(\Delta_\rho[-2,-4], \Delta_\rho[-1,-1]; \pi_{sc})$
           \item $L(\Delta_\rho[-4,-4], \Delta_\rho[-1,-2];T_{I,1}^{3}(\pi_{sc}))$  
           \item $L(\Delta_\rho[-4,-4], \Delta_\rho[-1,-1];T_{I,1}^{2}(T_{I,1}^{3}(\pi_{sc})))$ 
           \item $L(\Delta_\rho[-4,-4], \Delta_\rho[-2,-2],\Delta_\rho[-1,-1];T_{I,1}^{3}(\pi_{sc}))$ 
           \item $L(\Delta_\rho[-3,-4], \Delta_\rho[-2,-2], \Delta_\rho[-1,-1]; \pi_{sc})$
       \end{itemize}
       \item $(2,2,2,3)$: 
       \begin{itemize}
           \item  $L(\Delta_\rho[-2,-2], \Delta_\rho[-2,-2]; T_{I,1}^{2}(T_{I,1}^{3}(\pi_{sc})))$
           \item $L(\Delta_\rho[-2,-2], \Delta_\rho[-2,-2], \Delta_\rho[-2,-2]; T_{I,1}^{3}(\pi_{sc}))$
           \item $L(\Delta_\rho[-2,-3], \Delta_\rho[-2,-2], \Delta_\rho[-2,-2]; \pi_{sc})$
           \item $L(\Delta_\rho[-3,-3], \Delta_\rho[-2,-2], \Delta_\rho[-2,-2], \Delta_\rho[-2,-2]; \pi_{sc})$
       \end{itemize}
       \item $(2,2,3,3)$: 
       \begin{itemize}
       \item $L(\Delta_\rho[-2,-3];T_{I,1}^{2}(T_{I,1}^{3}(\pi_{sc})))$  
           \item $L(\Delta_\rho[-3,-3], \Delta_\rho[-2,-2]; T_{I,1}^{2}(T_{I,1}^{3}(\pi_{sc})))$
           \item $L(\Delta_\rho[-2,-3], \Delta_\rho[-2,-2]; T_{I,1}^{3}(\pi_{sc}))$
           \item $L(\Delta_\rho[-2,-3], \Delta_\rho[-2,-3]; \pi_{sc})$
           \item $L(\Delta_\rho[-3,-3], \Delta_\rho[-2,-2], \Delta_\rho[-2,-2]; T_{I,1}^{3}(\pi_{sc}))$
           \item $L(\Delta_\rho[-3,-3], \Delta_\rho[-2,-3], \Delta_\rho[-2,-2]; \pi_{sc})$
           \item $L(\Delta_\rho[-3,-3], \Delta_\rho[-3,-3], \Delta_\rho[-2,-2], \Delta_\rho[-2,-2]; \pi_{sc})$
       \end{itemize}
       \item $(2,2,3,4)$: 
       \begin{itemize}
           \item $L(\Delta_\rho[-2,-2]; T_{I,1}^{2}(T_{I,1}^{4}(T_{I,1}^{3}(\pi_{sc}))))$
           \item $L(\Delta_\rho[-2,-4], \Delta_\rho[-2,-2]; \pi_{sc})$
           \item $L(\Delta_\rho[-2,-2], \Delta_\rho[-2,-2];T_{I,1}^{4}(T_{I,1}^{3}(\pi_{sc})))$  
           \item $L(\Delta_\rho[-4,-4], \Delta_\rho[-2,-2]; T_{I,1}^{2}(T_{I,1}^{3}(\pi_{sc})))$
           \item $L(\Delta_\rho[-4,-4], \Delta_\rho[-2,-2], \Delta_\rho[-2,-2]; T_{I,1}^{3}(\pi_{sc}))$
           \item $L(\Delta_\rho[-3,-4], \Delta_\rho[-2,-2], \Delta_\rho[-2,-2]; \pi_{sc})$
           \item $L(\Delta_\rho[-4,-4], \Delta_\rho[-2,-3], \Delta_\rho[-2,-2];\pi_{sc})$ 
           \item $L(\Delta_\rho[-4,-4], \Delta_\rho[-3,-3], \Delta_\rho[-2,-2], \Delta_\rho[-2,-2]; \pi_{sc})$
       \end{itemize}
       \item $(2,3,3,3)$: 
       \begin{itemize}
       \item $L(\Delta_\rho[-3,-3], \Delta_\rho[-2,-3]; T_{I,1}^{3}(\pi_{sc}))$ + 
           \item  $L(\Delta_\rho[-3,-3], \Delta_\rho[-3,-3]; T_{I,1}^{2}(T_{I,1}^{3}(\pi_{sc})))$
           \item $L(\Delta_\rho[-3,-3], \Delta_\rho[-3,-3], \Delta_\rho[-2,-2]; T_{I,1}^{3}(\pi_{sc}))$
           \item $L(\Delta_\rho[-3,-3],\Delta_\rho[-3,-3], \Delta_\rho[-2,-3]; \pi_{sc})$
           \item $L(\Delta_\rho[-3,-3], \Delta_\rho[-3,-3], \Delta_\rho[-3,-3], \Delta_\rho[-2,-2]; \pi_{sc})$
       \end{itemize}
       \item $(2,3,3,4)$: 
       \begin{itemize}
           \item $L(\Delta_\rho[-3,-3]; T_{I,1}^{2}(T_{I,1}^{4}(T_{I,1}^{3}(\pi_{sc}))))$
           \item $L(\Delta_\rho[-2,-3];T_{I,1}^{4}(T_{I,1}^{3}(\pi_{sc})))$  
           \item $L(\Delta_\rho[-2,-4]; T_{I,1}^{3}(\pi_{sc}))$
           \item $L(\Delta_\rho[-3,-4];T_{I,1}^{2}(T_{I,1}^{3}(\pi_{sc})))$  
           \item $L(\Delta_\rho[-4,-4], \Delta_\rho[-3,-3]; T_{I,1}^{2}(T_{I,1}^{3}(\pi_{sc})))$
           \item $L(\Delta_\rho[-3,-2], \Delta_\rho[-2,-2]; T_{I,1}^{3}(\pi_{sc}))$
           \item $L(\Delta_\rho[-3,-4], \Delta_\rho[-2,-2];T_{I,1}^{3}(\pi_{sc}))$  
           \item $L(\Delta_\rho[-4,-4], \Delta_\rho[-2,-3];T_{I,1}^{3}(\pi_{sc}))$  
           \item $L(\Delta_\rho[-2,-4], \Delta_\rho[-3,-3]; \pi_{sc})$
           \item $L(\Delta_\rho[-3,-3], \Delta_\rho[-2,-4]; \pi_{sc})$
           \item $L(\Delta_\rho[-3,-4], \Delta_\rho[-3,-3], \Delta_\rho[-2,-2]; \pi_{sc})$
           \item $L(\Delta_\rho[-4,-4],\Delta_\rho[-3,-3], \Delta_\rho[-2,-3]; \pi_{sc})$
           \item $L(\Delta_\rho[-4,-4], \Delta_\rho[-3,-3], \Delta_\rho[-3,-3], \Delta_\rho[-2,-2]; \pi_{sc})$
       \end{itemize}
           \item $(2,3,4,4)$: 
           \begin{itemize}
               \item $L(\Delta_\rho[-4,-4]; T_{I,1}^{2}(T_{I,1}^{4}(T_{I,1}^{3}(\pi_{sc}))))$
               \item $L(\Delta_\rho[-4,-4], \Delta_\rho[-4,-4]; T_{I,1}^{2}(T_{I,1}^{3}(\pi_{sc})))$
               \item $L(\Delta_\rho[-4,-4], \Delta_\rho[-2,-2]; T_{I,1}^{4}(T_{I,1}^{3}(\pi_{sc})))$  
               \item $L(\Delta_\rho[-4,-4], \Delta_\rho[-2,-4]; \pi_{sc})$
               \item $L(\Delta_\rho[-4,-4], \Delta_\rho[-4,-4], \Delta_\rho[-2,-2]; T_{I,1}^{3}(\pi_{sc}))$
               \item $L(\Delta_\rho[-4,-4], \Delta_\rho[-3,-4], \Delta_\rho[-2,-2]; \pi_{sc})$
               \item $L(\Delta_\rho[-4,-4],\Delta_\rho[-4,-4], \Delta_\rho[-2,-3]; \pi_{sc})$
               \item $L(\Delta_\rho[-4,-4], \Delta_\rho[-4,-4], \Delta_\rho[-3,-3], \Delta_\rho[-2,-2]; \pi_{sc})$
           \end{itemize}
           \item $(2,3,4,5)$: 
           \begin{itemize}
           \item $L(\Delta_\rho[-4,-5]; T_{I,1}^{2}(T_{I,1}^{3}(\pi_{sc})))$+
           \item $L(\Delta_\rho[-2,-5]; \pi_{sc})$  
               \item $L(\Delta_\rho[-5,-5]; T_{I,1}^{2}(T_{I,1}^{4}(T_{I,1}^{3}(\pi_{sc}))))$
               \item $L(\Delta_\rho[-5,-5], \Delta_\rho[-4,-4]; T_{I,1}^{2}(T_{I,1}^{3}(\pi_{sc})))$
               \item $L(\Delta_\rho[-5,-5],\Delta_\rho[-2,-2]; T_{I,1}^{4}(T_{I,1}^{3}(\pi_{sc})))$  
               \item $L(\Delta_\rho[-4,-5], \Delta_\rho[-2,-2];T_{I,1}^{3}(\pi_{sc}))$  
               \item $L(\Delta_\rho[-4,-5], \Delta_\rho[-2,-3]; \pi_{sc})$
               \item $L(\Delta_\rho[-3,-5], \Delta_\rho[-2,-2]; \pi_{sc})$
               \item $L(\Delta_\rho[-5,-5], \Delta_\rho[-2,-4]; \pi_{sc})$
               \item $L(\Delta_\rho[-4,-5], \Delta_\rho[-3,-3], \Delta_\rho[-2,-2]; \pi_{sc})$
               \item $L(\Delta_\rho[-5,-5], \Delta_\rho[-4,-4], \Delta_\rho[-2,-2]; T_{I,1}^{3}(\pi_{sc}))$  
               \item $L(\Delta_\rho[-5,-5], \Delta_\rho[-3,-4], \Delta_\rho[-2,-2]; \pi_{sc})$
           \end{itemize}
           \item $(3,3,3,3)$: 
           \begin{itemize}
               \item $L(\Delta_\rho[-3,-3], \Delta_\rho[-3,-3], \Delta_\rho[-3,-3]; T_{I,1}^{3}(\pi_{sc}))$
               \item $L(\Delta_\rho[-3,-3], \Delta_\rho[-3,-3], \Delta_\rho[-3,-3], \Delta_\rho[-3,-3]; \pi_{sc})$
           \end{itemize}
           \item $(3,3,3,4)$: 
           \begin{itemize}
           \item $L(\Delta_\rho[-3,-3],\Delta_\rho[-3,-3];T_{I,1}^{4}(T_{I,1}^{3}(\pi_{sc})))$  
               \item $L(\Delta_\rho[-3,-4], \Delta_\rho[-3,-3]; T_{I,1}^{3}(\pi_{sc}))$
               \item $L(\Delta_\rho[-4,-4], \Delta_\rho[-3,-3], \Delta_\rho[-3,-3]; T_{I,1}^{3}(\pi_{sc}))$
               \item $L(\Delta_\rho[-3,-4], \Delta_\rho[-3,-3], \Delta_\rho[-3,-3]; \pi_{sc})$
               \item $L(\Delta_\rho[-4,-4], \Delta_\rho[-3,-3], \Delta_\rho[-3,-3], \Delta_\rho[-3,-3]; \pi_{sc})$
           \end{itemize}
           \item $(3,3,4,4)$: 
           \begin{itemize}
           \item $L(\Delta_\rho[-3,-4]; T_{I,1}^{4}(T_{I,1}^{3}(\pi_{sc})))$  
           \item $L(\Delta_\rho[-4,-4], \Delta_\rho[-3,-4]; T_{I,1}^{3}(\pi_{sc}))$
               \item $L(\Delta_\rho[-4,-4],\Delta_\rho[-3,-3]; T_{I,1}^{4}(T_{I,1}^{3}(\pi_{sc})))$
               \item $L(\Delta_\rho[-3,-4], \Delta_\rho[-3,-4]; \pi_{sc})$
               \item $L(\Delta_\rho[-4,-4], \Delta_\rho[-4,-4], \Delta_\rho[-3,-3]; T_{I,1}^{3}(\pi_{sc}))$
               \item $L(\Delta_\rho[-4,-4], \Delta_\rho[-3,-4], \Delta_\rho[-3,-3]; \pi_{sc})$
               \item $L(\Delta_\rho[-4,-4], \Delta_\rho[-4,-4], \Delta_\rho[-3,-3], \Delta_\rho[-3,-3]; \pi_{sc})$
           \end{itemize}
           \item $(3,3,4,5):$
           \begin{itemize}
               \item $L(\Delta_\rho[-3,-3]; T_{I,1}^{5}(T_{I,1}^{4}(T_{I,1}^{3}(\pi_{sc}))))$
               \item $L(\Delta_\rho[-3,-5]; T_{I,1}^{3}(\pi_{sc}))$
               \item $L(\Delta_\rho[-5,-5],\Delta_\rho[-3,-3]; T_{I,1}^{4}(T_{I,1}^{3}(\pi_{sc})))$ 
               \item $L(\Delta_\rho[-4,-5], \Delta_\rho[-3,-3]; T_{I,1}^{3}(\pi_{sc}))$
               \item $L(\Delta_\rho[-3,-5], \Delta_\rho[-3,-3]; \pi_{sc})$
               \item $L(\Delta_\rho[-5,-5],\Delta_\rho[-3,-4];T_{I,1}^{3}(\pi_{sc}))$  
               \item $L(\Delta_\rho[-5,-5], \Delta_\rho[-4,-4], \Delta_\rho[-3,-3]; T_{I,1}^{3}(\pi_{sc}))$
               \item $L(\Delta_\rho[-4,-5], \Delta_\rho[-3,-3], \Delta_\rho[-3,-3]; \pi_{sc})$
               \item $L(\Delta_\rho[-5,-5], \Delta_\rho[-3,-4], \Delta_\rho[-3,-3]; \pi_{sc})$
               \item $L(\Delta_\rho[-5,-5], \Delta_\rho[-4,-4], \Delta_\rho[-3,-3], \Delta_\rho[-3,-3]; \pi_{sc})$
           \end{itemize}
           \item $(3,4,4,4)$: 
           \begin{itemize}
           \item $L(\Delta_\rho[-4,-4],\Delta_\rho[-4,-4];T_{I,1}^{4}(T_{I,1}^{3}(\pi_{sc})))$ 
               \item $L(\Delta_\rho[-4,-4], \Delta_\rho[-4,-4], \Delta_\rho[-4,-4]; T_{I,1}^{3}(\pi_{sc}))$
               \item $L(\Delta_\rho[-4,-4],\Delta_\rho[-4,-4], \Delta_\rho[-3,-4]; \pi_{sc})$
               \item $L(\Delta_\rho[-4,-4], \Delta_\rho[-4,-4], \Delta_\rho[-4,-4], \Delta_\rho[-3,-3]; \pi_{sc})$
           \end{itemize}
           \item $(3,4,4,5)$: 
           \begin{itemize}
               \item $L(\Delta_\rho[-4,-4]; T_{I,1}^{5}(T_{I,1}^{3}(T_{I,1}^{3}(\pi_{sc}))))$
               \item $L(\Delta_\rho[-4,-5]; T_{I,1}^{4}(\pi_{sc}))$  
               \item $L(\Delta_\rho[-5,-5],\Delta_\rho[-4,-4]; T_{I,1}^{4}(T_{I,1}^{3}(\pi_{sc})))$
               \item $L(\Delta_\rho[-4,-5], \Delta_\rho[-4,-4]; T_{I,1}^{3}(\pi_{sc}))$
               \item $L(\Delta_\rho[-4,-5], \Delta_\rho[-3,-4]; \pi_{sc})$
               \item $L(\Delta_\rho[-3,-5], \Delta_\rho[-4,-4]; \pi_{sc})$
               \item $L(\Delta_\rho[-4,-4], \Delta_\rho[-3,-5]; \pi_{sc})$
               \item $L(\Delta_\rho[-4,-5],\Delta_\rho[-4,-4],\Delta_\rho[-3,-3];\pi_{sc})$ 
               \item $L(\Delta_\rho[-5,-5], 
               \Delta_\rho[-4,-4], \Delta_\rho[-4,-4]; T_{I,1}^{3}(\pi_{sc}))$
               \item $L(\Delta_\rho[-5,-5],\Delta_\rho[-4,-4], \Delta_\rho[-3,-4]; \pi_{sc})$
               \item $L(\Delta_\rho[-5,-5], \Delta_\rho[-4,-4], \Delta_\rho[-4,-4], \Delta_\rho[-3,-3]; \pi_{sc})$
           \end{itemize}
           \item $(3,4,5,5)$: 
           \begin{itemize}
               \item $L(\Delta_\rho[-5,-5]; T_{I,1}^{5}(T_{I,1}^{4}(T_{I,1}^{3}(\pi_{sc}))))$
               \item $L(\Delta_\rho[-5,-5], \Delta_\rho[-3,-5]; \pi_{sc})$
               \item $L(\Delta_\rho[-5,-5], \Delta_\rho[-4,-5]; T_{I,1}^{3}(\pi_{sc}))$
               \item $L(\Delta_\rho[-5,-5], \Delta_\rho[-5,-5]; T_{I,1}^{4}(T_{I,1}^{3}(\pi_{sc})))$ 
               \item $L(\Delta_\rho[-5,-5], \Delta_\rho[-5,-5], \Delta_\rho[-4,-4]; T_{I,1}^{3}(\pi_{sc}))$
               \item $L(\Delta_\rho[-5,-5], \Delta_\rho[-4,-5], \Delta_\rho[-3,-3]; \pi_{sc})$
               \item $L(\Delta_\rho[-5,-5],\Delta_\rho[-5,-5], \Delta_\rho[-3,-4]; \pi_{sc})$
               \item $L(\Delta_\rho[-5,-5], \Delta_\rho[-5,-5], \Delta_\rho[-4,-4], \Delta_\rho[-3,-3]; \pi_{sc})$
           \end{itemize}
           \item $(3,4,5,6)$: 
           \begin{itemize}
           \item $L(\Delta_\rho[-3,-6];\pi_{sc})$ 
               \item $L(\Delta_\rho[-6,-6]; T_{I,1}^{5}(T_{I,1}^{4}(T_{I,1}^{3}(\pi_{sc}))))$
               \item $L(\Delta_\rho[-4,-6]; T_{I,1}^{3}(\pi_{sc}))$
               \item $L(\Delta_\rho[-5,-6];T_{I,1}^{4}(T_{I,1}^{3}(\pi_{sc})))$ 
               \item $L(\Delta_\rho[-6,-6],\Delta_\rho[5,-5]; T_{I,1}^{4}(T_{I,1}^{3}(\pi_{sc})))$
               \item $L(\Delta_\rho[-5,-6], \Delta_\rho[-4,-4]; T_{I,1}^{3}(\pi_{sc}))$
               \item $L(\Delta_\rho[-5,-6], \Delta_\rho[-3,-4]; \pi_{sc})$
               \item $L(\Delta_\rho[-4,-6], \Delta_\rho[-3,-3]; \pi_{sc})$
               \item $L(\Delta_\rho[-6,-6], \Delta_\rho[-3,-5]; \pi_{sc})$
               \item $L(\Delta_\rho[-6,-6],\Delta_\rho[-4,-5];T_{I,1}^{3}(\pi_{sc}))$
               \item $L(\Delta_\rho[-6,-6], \Delta_\rho[-5,-5], \Delta_\rho[-4,-4]; T_{I,1}^{3}(\pi_{sc}))$
               \item $L(\Delta_\rho[-5,-6], \Delta_\rho[-4,-4], \Delta_\rho[-3,-3]; \pi_{sc})$
               \item $L(\Delta_\rho[-6,-6],\Delta_\rho[-4,-5], \Delta_\rho[-3,-3]; \pi_{sc})$ 
               \item $L(\Delta_\rho[-6,-6],\Delta_\rho[-5,-5], \Delta_\rho[-3,-4]; \pi_{sc})$
           \end{itemize}
       \end{enumerate}
       \item $(\alpha = \frac{7}{2})$:
       \begin{enumerate}
       \item $(\frac{3}{2}, \frac{3}{2}, \frac{5}{2}, \frac{7}{2}):$
       \begin{itemize}
           \item $L(\Delta_\rho[-\frac{3}{2}, -\frac{3}{2}]; T_{I,1}^{\frac{3}{2}}(T_{I,1}^{\frac{5}{2}}(T_{I,1}^{\frac{7}{2}}(\pi_{sc}))))$
           \item $L(\Delta_\rho[-\frac{3}{2}, -\frac{3}{2}], \Delta_\rho[-\frac{3}{2}, -\frac{3}{2}]; T_{I,1}^{\frac{5}{2}}(T_{I,1}^{\frac{7}{2}}(\pi_{sc})))$
           \item $L(\Delta_\rho[-\frac{3}{2}, -\frac{5}{2}], \Delta_\rho[-\frac{3}{2}, -\frac{3}{2}]; T_{I,1}^{\frac{7}{2}}(\pi_{sc}))$
           \item $L(\Delta_\rho[-\frac{3}{2}, -\frac{7}{2}], \Delta_\rho[-\frac{3}{2}, -\frac{3}{2}]; \pi_{sc})$
           \item $L(\Delta_\rho[-\frac{5}{2}, -\frac{5}{2}], \Delta_\rho[-\frac{3}{2}, -\frac{3}{2}], \Delta_\rho[-\frac{3}{2}, -\frac{3}{2}]; T_{I,1}^{\frac{7}{2}}(\pi_{sc}))$
           \item $L(\Delta_\rho[-\frac{5}{2}, -\frac{7}{2}], \Delta_\rho[-\frac{3}{2}, -\frac{3}{2}], \Delta_\rho[-\frac{3}{2}, -\frac{3}{2}];\pi_{sc})$
           \item $L(\Delta_\rho[-\frac{7}{2}, -\frac{7}{2}], \Delta_\rho[-\frac{3}{2}, -\frac{5}{2}], \Delta_\rho[-\frac{3}{2}, -\frac{3}{2}]; \pi_{sc})$
           \item $L(\Delta_\rho[-\frac{7}{2},-\frac{7}{2}], \Delta_\rho[-\frac{5}{2},-\frac{5}{2}], \Delta_\rho[-\frac{3}{2},-\frac{3}{2}], \Delta_\rho[-\frac{3}{2},-\frac{3}{2}]; \pi_{sc})$
       \end{itemize}
       \item $(\frac{3}{2}, \frac{5}{2}, \frac{5}{2}, \frac{7}{2}):$
       \begin{itemize}
           \item $L(\Delta_\rho[-\frac{5}{2}, -\frac{5}{2}]; T_{I,1}^{\frac{3}{2}}(T_{I,1}^{\frac{5}{2}}(T_{I,1}^{\frac{7}{2}}(\pi_{sc}))))$
           \item $L(\Delta_\rho[-\frac{3}{2}, -\frac{5}{2}]; T_{I,1}^{\frac{5}{2}}(T_{I,1}^{\frac{7}{2}}(\pi_{sc})))$  
           \item $L(\Delta_\rho[-\frac{5}{2}, -\frac{5}{2}],\Delta_\rho[-\frac{3}{2}, -\frac{7}{2}]; \pi_{sc})$
           \item $L(\Delta_\rho[-\frac{5}{2}, -\frac{5}{2}], \Delta_\rho[-\frac{3}{2}, -\frac{5}{2}]; T_{I,1}^{\frac{7}{2}}(\pi_{sc}))$ 
           \item $L(\Delta_\rho[-\frac{5}{2}, -\frac{7}{2}], \Delta_\rho[-\frac{5}{2}, -\frac{5}{2}], \Delta_\rho[-\frac{3}{2}, -\frac{3}{2}]; \pi_{sc})$
           \item $L(\Delta_\rho[-\frac{5}{2}, -\frac{5}{2}], \Delta_\rho[-\frac{5}{2}, -\frac{5}{2}], \Delta_\rho[-\frac{3}{2}, -\frac{3}{2}]; T_{I,1}^{\frac{7}{2}}(\pi_{sc}))$
           \item $L(\Delta_\rho[-\frac{7}{2}, -\frac{7}{2}], \Delta_\rho[-\frac{5}{2}, -\frac{5}{2}], \Delta_\rho[-\frac{3}{2}, -\frac{5}{2}]; \pi_{sc})$
           \item $L(\Delta_\rho[-\frac{7}{2},-\frac{7}{2}], \Delta_\rho[-\frac{5}{2},-\frac{5}{2}], \Delta_\rho[-\frac{5}{2},-\frac{5}{2}], \Delta_\rho[-\frac{3}{2},-\frac{3}{2}]; \pi_{sc})$
       \end{itemize}
       \item $(\frac{3}{2}, \frac{5}{2}, \frac{7}{2}, \frac{7}{2}):$
       \begin{itemize}
       \item $L(\Delta_\rho[-\frac{7}{2}, -\frac{7}{2}]; T_{I,1}^{\frac{3}{2}}(T_{I,1}^{\frac{5}{2}}(T_{I,1}^{\frac{7}{2}}(\pi_{sc}))))$ 
       
           \item $L(\Delta_\rho[-\frac{3}{2}, -\frac{7}{2}]; T_{I,1}^{\frac{7}{2}}(\pi_{sc}))$
           \item $L(\Delta_\rho[-\frac{7}{2}, -\frac{7}{2}], \Delta_\rho[-\frac{3}{2}, -\frac{3}{2}]; T_{I,1}^{\frac{5}{2}}(T_{I,1}^{\frac{7}{2}}(\pi_{sc})))$ 
           \item $L(\Delta_\rho[-\frac{5}{2}, -\frac{7}{2}], \Delta_\rho[-\frac{3}{2}, -\frac{3}{2}]; T_{I,1}^{\frac{7}{2}}(\pi_{sc}))$
           \item $L(\Delta_\rho[-\frac{7}{2}, -\frac{7}{2}],\Delta_\rho[-\frac{3}{2}, -\frac{7}{2}]; \pi_{sc})$
           \item $L(\Delta_\rho[-\frac{7}{2}, -\frac{7}{2}], \Delta_\rho[-\frac{5}{2}, -\frac{7}{2}], \Delta_\rho[-\frac{3}{2}, -\frac{3}{2}]; \pi_{sc})$
           \item $L(\Delta_\rho[-\frac{7}{2}, -\frac{7}{2}], \Delta_\rho[-\frac{7}{2}, -\frac{7}{2}], \Delta_\rho[-\frac{3}{2}, -\frac{5}{2}]; \pi_{sc})$
           \item $L(\Delta_\rho[-\frac{7}{2},-\frac{7}{2}], \Delta_\rho[-\frac{7}{2},-\frac{7}{2}], \Delta_\rho[-\frac{5}{2},-\frac{5}{2}], \Delta_\rho[-\frac{3}{2},-\frac{3}{2}]; \pi_{sc})$
       \end{itemize}
       \item $(\frac{3}{2}, \frac{5}{2}, \frac{7}{2}, \frac{9}{2}):$
       \begin{itemize}
           \item $L(\Delta_\rho[-\frac{9}{2}, -\frac{9}{2}]; T_{I,1}^{\frac{3}{2}}(T_{I,1}^{\frac{5}{2}}(T_{I,1}^{\frac{7}{2}}(\pi_{sc}))))$
           \item $L(\Delta_\rho[-\frac{3}{2}, -\frac{9}{2}]; \pi_{sc})$
           \item $L(\Delta_\rho[-\frac{7}{2}, -\frac{9}{2}], \Delta_\rho[-\frac{3}{2}, -\frac{5}{2}]; \pi_{sc})$
           \item $L(\Delta_\rho[-\frac{5}{2}, -\frac{9}{2}], \Delta_\rho[-\frac{3}{2}, -\frac{3}{2}]; \pi_{sc})$
           \item $L(\Delta_\rho[-\frac{9}{2}, -\frac{9}{2}], \Delta_\rho[-\frac{5}{2}, -\frac{5}{2}], \Delta_\rho[-\frac{3}{2}, -\frac{3}{2}];T_{I,1}^{\frac{7}{2}}(\pi_{sc}))$ 
           \item $L(\Delta_\rho[-\frac{9}{2}, -\frac{9}{2}], \Delta_\rho[-\frac{3}{2}, -\frac{3}{2}]; T_{I,1}^{\frac{5}{2}}(T_{I,1}^{\frac{7}{2}}(\pi_{sc})))$+
           \item $L(\Delta_\rho[-\frac{9}{2}, -\frac{9}{2}], \Delta_\rho[-\frac{3}{2},-\frac{5}{2}];T_{I,1}^{\frac{7}{2}}(\pi_{sc}))$ 
           \item $L(\Delta_\rho[-\frac{7}{2}, -\frac{9}{2}], \Delta_\rho[-\frac{5}{2}, -\frac{5}{2}], \Delta_\rho[-\frac{3}{2}, -\frac{3}{2}]; \pi_{sc})$
       \end{itemize}
       \item $(\frac{5}{2}, \frac{5}{2}, \frac{5}{2}, \frac{7}{2})$: 
       \begin{itemize}
           \item $L(\Delta_\rho[-\frac{5}{2}, -\frac{5}{2}], \Delta_\rho[-\frac{5}{2}, -\frac{5}{2}]; T_{I,1}^{\frac{5}{2}}(T_{I,1}^{\frac{7}{2}}(\pi_{sc})))$
           \item $L(\Delta_\rho[-\frac{5}{2}, -\frac{5}{2}], \Delta_\rho[-\frac{5}{2}, -\frac{5}{2}], \Delta_\rho[-\frac{5}{2},-\frac{5}{2}], T_{I,1}^{\frac{7}{2}}(\pi_{sc}))$
           \item $L(\Delta_\rho[-\frac{5}{2}, -\frac{7}{2}], \Delta_\rho[-\frac{5}{2}, -\frac{5}{2}], \Delta_\rho[-\frac{5}{2}, -\frac{5}{2}]; \pi_{sc})$
           \item $L(\Delta_\rho[-\frac{7}{2},-\frac{7}{2}], \Delta_\rho[-\frac{5}{2},-\frac{5}{2}], \Delta_\rho[-\frac{5}{2},-\frac{5}{2}], \Delta_\rho[-\frac{5}{2},-\frac{5}{2}]; \pi_{sc})$
       \end{itemize}
       \item $(\frac{5}{2}, \frac{5}{2}, \frac{7}{2}, \frac{7}{2})$: 
       \begin{itemize}
           \item $L(\Delta_\rho[-\frac{7}{2}, -\frac{7}{2}], \Delta_\rho[-\frac{5}{2}, -\frac{5}{2}]; T_{I,1}^{\frac{5}{2}}(T_{I,1}^{\frac{7}{2}}(\pi_{sc})))$
           \item $L(\Delta_\rho[-\frac{5}{2}, -\frac{7}{2}], \Delta_\rho[-\frac{5}{2}, -\frac{5}{2}]; T_{I,1}^{\frac{7}{2}}(\pi_{sc}))$
           \item $L(\Delta_\rho[-\frac{5}{2}, -\frac{7}{2}], \Delta_\rho[-\frac{5}{2}, -\frac{7}{2}]; \pi_{sc})$
           \item $L(\Delta_\rho[-\frac{7}{2}, -\frac{7}{2}], \Delta_\rho[-\frac{5}{2}, -\frac{5}{2}], \Delta_\rho[-\frac{5}{2}, -\frac{5}{2}]; T_{I,1}^{\frac{7}{2}}(\pi_{sc}))$
           \item $L(\Delta_\rho[-\frac{7}{2}, -\frac{7}{2}], \Delta_\rho[-\frac{5}{2}, -\frac{7}{2}], \Delta_\rho[-\frac{5}{2}, -\frac{5}{2}]; \pi_{sc})$
           \item $L(\Delta_\rho[-\frac{7}{2},-\frac{7}{2}], \Delta_\rho[-\frac{7}{2},-\frac{7}{2}], \Delta_\rho[-\frac{5}{2},-\frac{5}{2}], \Delta_\rho[-\frac{5}{2},-\frac{5}{2}]; \pi_{sc})$
           \item $L(\Delta_\rho[-\frac{7}{2},-\frac{7}{2}], \Delta_\rho[-\frac{7}{2},-\frac{7}{2}], \Delta_\rho[-\frac{5}{2},-\frac{5}{2}], \Delta_\rho[-\frac{5}{2},-\frac{5}{2}]; \pi_{sc})$
       \end{itemize}
       \item $(\frac{5}{2}, \frac{5}{2}, \frac{7}{2}, \frac{9}{2})$: 
       \begin{itemize}
           \item $L(\Delta_\rho[-\frac{5}{2}, -\frac{5}{2}]; T_{I,1}^{\frac{5}{2}}(T_{I,1}^{\frac{9}{2}}(T_{I,1}^{\frac{7}{2}}(\pi_{sc}))))$
           \item $L(\Delta_\rho[-\frac{5}{2}, -\frac{5}{2}], \Delta_\rho[-\frac{5}{2}, -\frac{5}{2}]; T_{I,1}^{\frac{9}{2}}(T_{I,1}^{\frac{7}{2}}(\pi_{sc})))$
           \item $L(\Delta_\rho[-\frac{9}{2}, -\frac{9}{2}], \Delta_\rho[-\frac{5}{2}, -\frac{5}{2}]; T_{I,1}^{\frac{5}{2}}(T_{I,1}^{\frac{7}{2}}(\pi_{sc})))$  
           \item $L(\Delta_\rho[-\frac{5}{2}, -\frac{9}{2}], \Delta_\rho[-\frac{5}{2}, -\frac{5}{2}]; \pi_{sc})$
           \item $L(\Delta_\rho[-\frac{9}{2}, -\frac{9}{2}], \Delta_\rho[-\frac{5}{2}, -\frac{5}{2}], \Delta_\rho[-\frac{5}{2}, -\frac{5}{2}]; T_{I,1}^{\frac{7}{2}}(\pi_{sc}))$
           \item $L(\Delta_\rho[-\frac{7}{2}, -\frac{9}{2}], \Delta_\rho[-\frac{5}{2}, -\frac{5}{2}], \Delta_\rho[-\frac{5}{2}, -\frac{5}{2}]; \pi_{sc})$
           \item $L(\Delta_\rho[-\frac{9}{2}, -\frac{9}{2}], \Delta_\rho[-\frac{5}{2}, -\frac{7}{2}], \Delta_\rho[-\frac{5}{2}, -\frac{5}{2}]; \pi_{sc})$ 
           \item $L(\Delta_\rho[-\frac{9}{2}, -\frac{9}{2}], \Delta_\rho[-\frac{7}{2}, -\frac{7}{2}], \Delta_\rho[-\frac{5}{2}, -\frac{5}{2}], \Delta_\rho[-\frac{5}{2}, -\frac{5}{2}]; \pi_{sc})$
       \end{itemize}
       \item $(\frac{5}{2}, \frac{7}{2}, \frac{7}{2}, \frac{7}{2})$: 
       \begin{itemize}
           \item  $L(\Delta_\rho[-\frac{7}{2}, -\frac{7}{2}], \Delta_\rho[-\frac{7}{2}, -\frac{7}{2}]; T_{I,1}^{\frac{5}{2}}(T_{I,1}^{\frac{7}{2}}(\pi_{sc})))$
           \item $L(\Delta_\rho[-\frac{7}{2}, -\frac{7}{2}], \Delta_\rho[-\frac{7}{2}, -\frac{7}{2}];T_{I,1}^{\frac{7}{2}}(T_{I,1}^{\frac{5}{2}}(\pi_{sc})))$ 
           \item $L(\Delta_\rho[-\frac{7}{2}, -\frac{7}{2}], \Delta_\rho[-\frac{7}{2}, -\frac{7}{2}], \Delta_\rho[-\frac{5}{2}, -\frac{5}{2}]; T_{I,1}^{\frac{7}{2}}(\pi_{sc}))$
           \item $L(\Delta_\rho[-\frac{7}{2}, -\frac{7}{2}], \Delta_\rho[-\frac{7}{2}, -\frac{7}{2}], \Delta_\rho[-\frac{5}{2}, -\frac{7}{2}]; \pi_{sc})$
           \item $L(\Delta_\rho[-\frac{7}{2},-\frac{7}{2}], \Delta_\rho[-\frac{7}{2},-\frac{7}{2}], \Delta_\rho[-\frac{7}{2},-\frac{7}{2}], \Delta_\rho[-\frac{5}{2},-\frac{5}{2}]; \pi_{sc})$
       \end{itemize}
       \item $(\frac{5}{2}, \frac{7}{2}, \frac{7}{2}, \frac{9}{2})$: 
       \begin{itemize}
           \item $L(\Delta_\rho[-\frac{7}{2}, -\frac{7}{2}]; T_{I,1}^{\frac{5}{2}}(T_{I,1}^{\frac{9}{2}}(T_{I,1}^{\frac{7}{2}}(\pi_{sc}))))$
           \item $L(\Delta_\rho[-\frac{5}{2}, -\frac{9}{2}]; T_{I,1}^{\frac{7}{2}}(\pi_{sc}))$
           \item $L(\Delta_\rho[-\frac{7}{2}, -\frac{9}{2}];T_{I,1}^{\frac{5}{2}}(T_{I,1}^{\frac{7}{2}}(\pi_{sc})))$ 
           \item $L(\Delta_\rho[-\frac{5}{2}, -\frac{7}{2}]; T_{I,1}^{\frac{9}{2}}(T_{I,1}^{\frac{7}{2}}(\pi_{sc})))$ 
           \item $L(\Delta_\rho[-\frac{9}{2}, -\frac{9}{2}], \Delta_\rho[-\frac{7}{2}, -\frac{7}{2}]; T_{I,1}^{\frac{5}{2}}(T_{I,1}^{\frac{7}{2}}
           (\pi_{sc})))$
           \item $L(\Delta_\rho[-\frac{9}{2}, -\frac{9}{2}], \Delta_\rho[-\frac{7}{2}, -\frac{7}{2}]; T_{I,1}^{\frac{7}{2}}(T_{I,1}^{\frac{5}{2}}(\pi_{sc})))$ 
           \item $L(\Delta_\rho[-\frac{7}{2}, -\frac{9}{2}], \Delta_\rho[-\frac{5}{2}, -\frac{5}{2}]; T_{I,1}^{\frac{7}{2}}(\pi_{sc}))$
           \item $L(\Delta_\rho[-\frac{5}{2}, -\frac{9}{2}], \Delta_\rho[-\frac{7}{2}, -\frac{7}{2}]; \pi_{sc})$
           \item $L(\Delta_\rho[-\frac{7}{2}, -\frac{7}{2}],\Delta_\rho[-\frac{5}{2}, -\frac{9}{2}]; \pi_{sc})$
           \item $L(\Delta_\rho[-\frac{7}{2}, -\frac{9}{2}], \Delta_\rho[-\frac{7}{2}, -\frac{7}{2}], \Delta_\rho[-\frac{5}{2}, -\frac{5}{2}]; \pi_{sc})$
           \item $L(\Delta_\rho[-\frac{9}{2}, -\frac{9}{2}], \Delta_\rho[-\frac{7}{2}, -\frac{7}{2}], \Delta_\rho[-\frac{5}{2}, -\frac{7}{2}]; \pi_{sc})$
           \item $L(\Delta_\rho[-\frac{9}{2},-\frac{9}{2}], \Delta_\rho[-\frac{7}{2},-\frac{7}{2}], \Delta_\rho[-\frac{7}{2},-\frac{7}{2}], \Delta_\rho[-\frac{5}{2},-\frac{5}{2}]; \pi_{sc})$
        \end{itemize}
     \item $(\frac{5}{2}, \frac{7}{2}, \frac{9}{2}, \frac{9}{2})$: 
       \begin{itemize}
           \item $L(\Delta_\rho[-\frac{9}{2}, -\frac{9}{2}]; T_{I,1}^{\frac{5}{2}}(T_{I,1}^{\frac{9}{2}}(T_{I,1}^{\frac{7}{2}}(\pi_{sc}))))$
           \item $L(\Delta_\rho[-\frac{9}{2}, -\frac{9}{2}], \Delta_\rho[-\frac{9}{2}, -\frac{9}{2}]; T_{I,1}^{\frac{5}{2}}(T_{I,1}^{\frac{7}{2}}(\pi_{sc})))$
           \item $L(\Delta_\rho[-\frac{9}{2}, -\frac{9}{2}], \Delta_\rho[-\frac{5}{2}, -\frac{5}{2}]; T_{I,1}^{\frac{9}{2}}(T_{I,1}^{\frac{7}{2}}(\pi_{sc})))$
           \item $L(\Delta_\rho[-\frac{9}{2}, -\frac{9}{2}],\Delta_\rho[-\frac{5}{2}, -\frac{9}{2}]; \pi_{sc})$
           \item $L(\Delta_\rho[-\frac{9}{2}, -\frac{9}{2}], \Delta_\rho[-\frac{9}{2}, -\frac{9}{2}], \Delta_\rho[-\frac{5}{2}, -\frac{5}{2}]; T_{I,1}^{\frac{7}{2}}(\pi_{sc}))$
           \item $L(\Delta_\rho[-\frac{9}{2}, -\frac{9}{2}], \Delta_\rho[-\frac{7}{2}, -\frac{9}{2}], \Delta_\rho[-\frac{5}{2}, -\frac{5}{2}]; \pi_{sc})$
           \item $L(\Delta_\rho[-\frac{9}{2}, -\frac{9}{2}], \Delta_\rho[-\frac{9}{2}, -\frac{9}{2}], \Delta_\rho[-\frac{5}{2}, -\frac{7}{2}]; \pi_{sc})$
           \item $L(\Delta_\rho[-\frac{9}{2},-\frac{9}{2}], \Delta_\rho[-\frac{9}{2},-\frac{9}{2}], \Delta_\rho[-\frac{7}{2},-\frac{7}{2}], \Delta_\rho[-\frac{5}{2},-\frac{5}{2}]; \pi_{sc})$
        \end{itemize}
    \item $(\frac{5}{2}, \frac{7}{2}, \frac{9}{2}, \frac{11}{2})$: 
       \begin{itemize}
           \item $L(\Delta_\rho[-\frac{11}{2}, -\frac{11}{2}]; T_{I,1}^{\frac{5}{2}}(T_{I,1}^{\frac{9}{2}}(T_{I,1}^{\frac{7}{2}}(\pi_{sc}))))$
           \item $L(\Delta_\rho[-\frac{9}{2}, -\frac{11}{2}]; T_{I,1}^{\frac{5}{2}}(T_{I,1}^{\frac{7}{2}}(\pi_{sc})))$ +
           \item $L(\Delta_\rho[-\frac{5}{2}, -\frac{11}{2}]; \pi_{sc})$
           \item $L(\Delta_\rho[-\frac{11}{2}, -\frac{11}{2}], \Delta_\rho[-\frac{9}{2}, -\frac{9}{2}]; T_{I,1}^{\frac{5}{2}}(T_{I,1}^{\frac{7}{2}}(\pi_{sc})))$
           \item $L(\Delta_\rho[-\frac{11}{2}, -\frac{11}{2}], \Delta_\rho[-\frac{5}{2}, -\frac{5}{2}]; T_{I,1}^{\frac{9}{2}}(T_{I,1}^{\frac{7}{2}}(\pi_{sc}))$ 
           \item $L(\Delta_\rho[-\frac{9}{2},-\frac{11}{2}], \Delta_\rho[-\frac{5}{2},-\frac{5}{2}];T_{I,1}^{\frac{7}{2}}(\pi_{sc}))$ 
           \item $L(\Delta_\rho[-\frac{9}{2}, -\frac{11}{2}], \Delta_\rho[-\frac{5}{2}, -\frac{7}{2}]; \pi_{sc})$
           \item $L(\Delta_\rho[-\frac{7}{2}, -\frac{11}{2}], \Delta_\rho[-\frac{5}{2}, -\frac{5}{2}]; \pi_{sc})$
           \item $L(\Delta_\rho[-\frac{11}{2}, -\frac{11}{2}],\Delta_\rho[-\frac{5}{2}, -\frac{9}{2}]; \pi_{sc})$
           \item $L(\Delta_\rho[-\frac{11}{2},-\frac{11}{2}],\Delta_\rho[-\frac{9}{2},-\frac{9}{2}],\Delta_\rho[-\frac{5}{2},-\frac{5}{2}];T_{I,1}^{\frac{7}{2}}(\pi_{sc}))$
           \item $L(\Delta_\rho[-\frac{9}{2}, -\frac{11}{2}], \Delta_\rho[-\frac{7}{2}, -\frac{7}{2}], \Delta_\rho[-\frac{5}{2}, -\frac{5}{2}]; \pi_{sc})$
           \item $L(\Delta_\rho[-\frac{11}{2}, -\frac{11}{2}], \Delta_\rho[-\frac{7}{2}, -\frac{9}{2}], \Delta_\rho[-\frac{5}{2}, -\frac{5}{2}]; \pi_{sc})$
        \end{itemize}
        \item $(\frac{7}{2}, \frac{7}{2}, \frac{7}{2}, \frac{7}{2})$: 
        \begin{itemize}
            \item  $L(\Delta_\rho[-\frac{7}{2}, -\frac{7}{2}], \Delta_\rho[-\frac{7}{2}, -\frac{7}{2}], \Delta_\rho[-\frac{7}{2},-\frac{7}{2}], T_{I,1}^{\frac{7}{2}}(\pi_{sc}))$
            \item $L(\Delta_\rho[-\frac{7}{2},-\frac{7}{2}], \Delta_\rho[-\frac{7}{2},-\frac{7}{2}], \Delta_\rho[-\frac{7}{2},-\frac{7}{2}], \Delta_\rho[-\frac{7}{2},-\frac{7}{2}]; \pi_{sc})$
        \end{itemize}
        \item $(\frac{7}{2}, \frac{7}{2}, \frac{7}{2}, \frac{9}{2})$: 
        \begin{itemize}
            \item  $L(\Delta_\rho[-\frac{7}{2}, -\frac{7}{2}], \Delta_\rho[-\frac{7}{2}, -\frac{7}{2}]; T_{I,1}^{\frac{9}{2}}(T_{I,1}^{\frac{7}{2}}(\pi_{sc})))$
            \item $L(\Delta_\rho[-\frac{7}{2}, -\frac{9}{2}], \Delta_\rho[-\frac{7}{2}, -\frac{7}{2}]; T_{I,1}^{\frac{7}{2}}(\pi_{sc}))$
            \item $L(\Delta_\rho[-\frac{9}{2}, -\frac{9}{2}], \Delta_\rho[-\frac{7}{2}, -\frac{7}{2}], \Delta_\rho[-\frac{7}{2}, -\frac{7}{2}]; T_{I,1}^{\frac{7}{2}}(\pi_{sc}))$
            \item $L(\Delta_\rho[-\frac{7}{2}, -\frac{9}{2}], \Delta_\rho[-\frac{7}{2}, -\frac{7}{2}], \Delta_\rho[-\frac{7}{2}, -\frac{7}{2}]; \pi_{sc})$
            \item $L(\Delta_\rho[-\frac{9}{2},-\frac{9}{2}], \Delta_\rho[-\frac{7}{2},-\frac{7}{2}], \Delta_\rho[-\frac{7}{2},-\frac{7}{2}], \Delta_\rho[-\frac{7}{2},-\frac{7}{2}]; \pi_{sc})$
        \end{itemize}
        \item $(\frac{7}{2}, \frac{7}{2}, \frac{9}{2}, \frac{9}{2}):$
        \begin{itemize}
        \item $L(\Delta_\rho[-\frac{7}{2}, -\frac{9}{2}],\Delta_\rho[-\frac{7}{2}, -\frac{9}{2}]; \pi_{sc})$ 
            \item $L(\Delta_\rho[-\frac{9}{2}, -\frac{9}{2}], \Delta_\rho[-\frac{7}{2}, -\frac{7}{2}]; T_{I,1}^{\frac{9}{2}}(T_{I,1}^{\frac{7}{2}}(\pi_{sc})))$
            \item $L(\Delta_\rho[-\frac{7}{2}, -\frac{9}{2}], \Delta_\rho[-\frac{7}{2}, -\frac{9}{2}]; \pi_{sc})$
            \item $L(\Delta_\rho[-\frac{9}{2}, -\frac{9}{2}], \Delta_\rho[-\frac{7}{2}, -\frac{9}{2}]; T_{I,1}^{\frac{7}{2}}(\pi_{sc}))$
            \item $L(\Delta_\rho[-\frac{9}{2}, -\frac{9}{2}], \Delta_\rho[-\frac{9}{2}, -\frac{9}{2}], \Delta_\rho[-\frac{7}{2}, -\frac{7}{2}]; T_{I,1}^{\frac{7}{2}}(\pi_{sc}))$
            \item $L(\Delta_\rho[-\frac{9}{2}, -\frac{9}{2}], \Delta_\rho[-\frac{7}{2}, -\frac{9}{2}], \Delta_\rho[-\frac{7}{2}, -\frac{7}{2}]; \pi_{sc})$
            \item $L(\Delta_\rho[-\frac{9}{2},-\frac{9}{2}], \Delta_\rho[-\frac{9}{2},-\frac{9}{2}], \Delta_\rho[-\frac{7}{2},-\frac{7}{2}], \Delta_\rho[-\frac{7}{2},-\frac{7}{2}]; \pi_{sc})$
        \end{itemize}
           \item $(\frac{7}{2}, \frac{7}{2}, \frac{9}{2}, \frac{11}{2}):$
           \begin{itemize}
               \item $L(\Delta_\rho[-\frac{7}{2}, -\frac{7}{2}]; T_{I,1}^{\frac{11}{2}}(T_{I,1}^{\frac{9}{2}}(T_{I,1}^{\frac{7}{2}}(\pi_{sc}))))$
               \item $L(\Delta_\rho[-\frac{7}{2}, -\frac{11}{2}]; T_{I,1}^{\frac{7}{2}}(\pi_{sc}))$
               \item $L(\Delta_\rho[-\frac{9}{2}, -\frac{11}{2}], \Delta_\rho[-\frac{7}{2}, -\frac{7}{2}]; T_{I,1}^{\frac{7}{2}}(\pi_{sc}))$
               \item $L(\Delta_\rho[-\frac{7}{2}, -\frac{11}{2}], \Delta_\rho[-\frac{7}{2}, -\frac{7}{2}]; \pi_{sc})$
               \item $L(\Delta_\rho[-\frac{11}{2}, -\frac{11}{2}], \Delta_\rho[-\frac{7}{2},-\frac{7}{2}]; T_{I,1}^{\frac{9}{2}}(T_{I,1}^{\frac{7}{2}}(\pi_{sc})))$ 
               \item $L(\Delta_\rho[-\frac{11}{2}, -\frac{11}{2}], \Delta_\rho[-\frac{7}{2}, -\frac{9}{2}], T_{I,1}^{\frac{7}{2}}(\pi_{sc}))$ +
               \item $L(\Delta_\rho[-\frac{11}{2}, -\frac{11}{2}], \Delta_\rho[-\frac{9}{2}, -\frac{9}{2}], \Delta_\rho[-\frac{7}{2}, -\frac{7}{2}]; T_{I,1}^{\frac{7}{2}}(\pi_{sc}))$
               \item $L(\Delta_\rho[-\frac{9}{2}, -\frac{11}{2}], \Delta_\rho[-\frac{7}{2}, -\frac{7}{2}], \Delta_\rho[-\frac{7}{2}, -\frac{7}{2}];\pi_{sc})$
               \item $L(\Delta_\rho[-\frac{11}{2}, -\frac{11}{2}], \Delta_\rho[-\frac{7}{2}, -\frac{9}{2}], \Delta_\rho[-\frac{7}{2}, -\frac{7}{2}]; \pi_{sc})$
               \item $L(\Delta_\rho[-\frac{11}{2},-\frac{11}{2}], \Delta_\rho[-\frac{9}{2},-\frac{9}{2}], \Delta_\rho[-\frac{7}{2},-\frac{7}{2}], \Delta_\rho[-\frac{7}{2},-\frac{7}{2}]; \pi_{sc})$
           \end{itemize}
           \item $(\frac{7}{2}, \frac{9}{2}, \frac{9}{2}, \frac{9}{2})$: 
           \begin{itemize}
               \item $L(\Delta_\rho[-\frac{9}{2}, -\frac{9}{2}], \Delta_\rho[-\frac{9}{2}, -\frac{9}{2}]; T_{I,1}^{\frac{9}{2}}(T_{I,1}^{\frac{7}{2}}(\pi_{sc})))$
               \item  $L(\Delta_\rho[-\frac{9}{2}, -\frac{9}{2}], \Delta_\rho[-\frac{9}{2}, -\frac{9}{2}], \Delta_\rho[-\frac{9}{2},-\frac{9}{2}], T_{I,1}^{\frac{7}{2}}(\pi_{sc}))$
               \item $L(\Delta_\rho[-\frac{9}{2}, -\frac{9}{2}], \Delta_\rho[-\frac{9}{2}, -\frac{9}{2}], \Delta_\rho[-\frac{7}{2}, -\frac{9}{2}]; \pi_{sc})$
               \item $L(\Delta_\rho[-\frac{9}{2},-\frac{9}{2}], \Delta_\rho[-\frac{9}{2},-\frac{9}{2}], \Delta_\rho[-\frac{9}{2},-\frac{9}{2}], \Delta_\rho[-\frac{7}{2},-\frac{7}{2}]; \pi_{sc})$
           \end{itemize}
           \item $(\frac{7}{2}, \frac{9}{2}, \frac{9}{2}, \frac{11}{2}):$
           \begin{itemize}
               \item $L(\Delta_\rho[-\frac{9}{2}, -\frac{9}{2}]; T_{I,1}^{\frac{11}{2}}(T_{I,1}^{\frac{9}{2}}(T_{I,1}^{\frac{7}{2}}(\pi_{sc}))))$
               \item $L(\Delta_\rho[-\frac{9}{2}, -\frac{11}{2}]; T_{I,1}^{\frac{9}{2}}(T_{I,1}^{\frac{7}{2}}(\pi_{sc})))$ 
               \item $L(\Delta_\rho[-\frac{11}{2}, -\frac{11}{2}], \Delta_\rho[-\frac{9}{2}, -\frac{9}{2}]; T_{I,1}^{\frac{9}{2}}(T_{I,1}^{\frac{7}{2}}(\pi_{sc})))$
               \item $L(\Delta_\rho[-\frac{9}{2}, -\frac{11}{2}], \Delta_\rho[-\frac{9}{2}, -\frac{9}{2}]; T_{I,1}^{\frac{7}{2}}(\pi_{sc}))$
               \item $L(\Delta_\rho[-\frac{9}{2}, -\frac{11}{2}], \Delta_\rho[-\frac{7}{2}, -\frac{9}{2}]; \pi_{sc})$
               \item $L(\Delta_\rho[-\frac{7}{2}, -\frac{11}{2}], \Delta_\rho[-\frac{9}{2}, -\frac{9}{2}]; \pi_{sc})$
               \item $L(\Delta_\rho[-\frac{9}{2}, -\frac{9}{2}],\Delta_\rho[-\frac{7}{2}, -\frac{11}{2}]; \pi_{sc})$
               \item $L(\Delta_\rho[-\frac{9}{2}, -\frac{11}{2}], \Delta_\rho[-\frac{9}{2}, -\frac{9}{2}], \Delta_\rho[-\frac{7}{2}, -\frac{7}{2}];\pi_{sc})$ 
               \item $L(\Delta_\rho[-\frac{11}{2}, -\frac{11}{2}], \Delta_\rho[-\frac{9}{2}, -\frac{9}{2}], \Delta_\rho[-\frac{9}{2}, -\frac{9}{2}]; T_{I,1}^{\frac{7}{2}}(\pi_{sc}))$
               \item $L(\Delta_\rho[-\frac{11}{2}, -\frac{11}{2}], \Delta_\rho[-\frac{9}{2}, -\frac{9}{2}], \Delta_\rho[-\frac{7}{2}, -\frac{9}{2}]; \pi_{sc})$
               \item $L(\Delta_\rho[-\frac{11}{2},-\frac{11}{2}], \Delta_\rho[-\frac{9}{2},-\frac{9}{2}], \Delta_\rho[-\frac{9}{2},-\frac{9}{2}], \Delta_\rho[-\frac{7}{2},-\frac{7}{2}]; \pi_{sc})$
           \end{itemize}
           \item $(\frac{7}{2}, \frac{9}{2}, \frac{11}{2}, \frac{11}{2}):$
           \begin{itemize}
               \item $L(\Delta_\rho[-\frac{11}{2}, -\frac{11}{2}]; T_{I,1}^{\frac{11}{2}}(T_{I,1}^{\frac{9}{2}}(T_{I,1}^{\frac{7}{2}}(\pi_{sc}))))$
               \item $L(\Delta_\rho[-\frac{11}{2}, -\frac{11}{2}], \Delta_\rho[-\frac{11}{2}, -\frac{11}{2}]; T_{I,1}^{\frac{9}{2}}(T_{I,1}^{\frac{7}{2}}(\pi_{sc})))$
               \item $L(\Delta_\rho[-\frac{11}{2}, -\frac{11}{2}],\Delta_\rho[-\frac{7}{2}, -\frac{11}{2}]; \pi_{sc})$
               \item $L(\Delta_\rho[-\frac{11}{2},-\frac{11}{2}], \Delta_\rho[-\frac{9}{2}, -\frac{11}{2}]; T_{I,1}^{\frac{7}{2}}(\pi_{sc}))$ 
               \item $L(\Delta_\rho[-\frac{11}{2}, -\frac{11}{2}], \Delta_\rho[-\frac{11}{2}, -\frac{11}{2}], \Delta_\rho[-\frac{9}{2}, -\frac{9}{2}]; T_{I,1}^{\frac{7}{2}}(\pi_{sc}))$
               \item $L(\Delta_\rho[-\frac{11}{2}, -\frac{11}{2}], \Delta_\rho[-\frac{9}{2}, -\frac{11}{2}], \Delta_\rho[-\frac{7}{2}, -\frac{7}{2}]; \pi_{sc})$
               \item $L(\Delta_\rho[-\frac{11}{2}, -\frac{11}{2}], \Delta_\rho[-\frac{11}{2}, -\frac{11}{2}], \Delta_\rho[-\frac{7}{2}, -\frac{9}{2}]; \pi_{sc})$
               \item $L(\Delta_\rho[-\frac{11}{2},-\frac{11}{2}], \Delta_\rho[-\frac{11}{2},-\frac{11}{2}], \Delta_\rho[-\frac{9}{2},-\frac{9}{2}], \Delta_\rho[-\frac{7}{2},-\frac{7}{2}]; \pi_{sc})$
           \end{itemize}
           \item $(\frac{7}{2}, \frac{9}{2}, \frac{11}{2}, \frac{13}{2}):$
           \begin{itemize}
               \item $L(\Delta_\rho[-\frac{13}{2}, -\frac{13}{2}]; T_{I,1}^{\frac{11}{2}}(T_{I,1}^{\frac{9}{2}}(T_{I,1}^{\frac{7}{2}}(\pi_{sc}))))$
               \item $L(\Delta_\rho[-\frac{9}{2}, -\frac{13}{2}]; T_{I,1}^{\frac{7}{2}}(\pi_{sc}))$
               \item $L(\Delta_\rho[-\frac{7}{2}, -\frac{13}{2}]; \pi_{sc})$
               \item $L(\Delta_\rho[-\frac{11}{2},-\frac{13}{2}]; T_{I,1}^{\frac{9}{2}}(T_{I,1}^{\frac{7}{2}}(\pi_{sc})))$
               \item $L(\Delta_\rho[-\frac{13}{2}, -\frac{13}{2}], \Delta_\rho[-\frac{11}{2}, -\frac{11}{2}]; T_{I,1}^{\frac{9}{2}}(T_{I,1}^{\frac{7}{2}}(\pi_{sc})))$
               \item $L(\Delta_\rho[-\frac{11}{2}, -\frac{13}{2}], \Delta_\rho[-\frac{9}{2}, -\frac{9}{2}]; T_{I,1}^{\frac{7}{2}}(\pi_{sc}))$
               \item $L(\Delta_\rho[-\frac{13}{2}, -\frac{13}{2}], \Delta_\rho[-\frac{9}{2}, -\frac{11}{2}], T_{I,1}^{\frac{7}{2}}(\pi_{sc}))$ 
               \item $L(\Delta_\rho[-\frac{11}{2}, -\frac{13}{2}], \Delta_\rho[-\frac{7}{2}, -\frac{9}{2}]; \pi_{sc})$
               \item $L(\Delta_\rho[-\frac{9}{2}, -\frac{13}{2}], \Delta_\rho[-\frac{7}{2}, -\frac{7}{2}]; \pi_{sc})$
               \item $L(\Delta_\rho[-\frac{13}{2}, -\frac{13}{2}],\Delta_\rho[-\frac{7}{2}, -\frac{11}{2}]; \pi_{sc})$
               \item $L(\Delta_\rho[-\frac{11}{2}, -\frac{13}{2}], \Delta_\rho[-\frac{9}{2}, -\frac{9}{2}], \Delta_\rho[-\frac{7}{2}, -\frac{7}{2}]; \pi_{sc})$ 
               \item $L(\Delta_\rho[-\frac{13}{2}, -\frac{13}{2}], \Delta_\rho[-\frac{11}{2}, -\frac{11}{2}], \Delta_\rho[-\frac{9}{2}, -\frac{9}{2}]; T_{I,1}^{\frac{7}{2}}(\pi_{sc}))$
               \item $L(\Delta_\rho[-\frac{13}{2}, -\frac{13}{2}], \Delta_\rho[-\frac{11}{2}, -\frac{11}{2}], \Delta_\rho[-\frac{7}{2}, -\frac{9}{2}]; \pi_{sc})$
               \item $L(\Delta_\rho[-\frac{13}{2},-\frac{13}{2}],\Delta_\rho[-\frac{9}{2},-\frac{11}{2}], \Delta_\rho[-\frac{7}{2}, -\frac{7}{2}];\pi_{sc})$
           \end{itemize}
       \end{enumerate}
       \item $(\alpha \geq 4)$: 
       \begin{enumerate}
       \item $(\alpha -2, \alpha -2, \alpha -1, \alpha)$: 
       \begin{itemize}
           \item $L(\Delta_\rho[-\alpha +2, -\alpha +2]; T_{I,1}^{\alpha -2}(T_{I,1}^{\alpha -1}(T_{I,1}^{\alpha}(\pi_{sc}))))$
           \item $L(\Delta_\rho[-\alpha +2 ,-\alpha +2], \Delta_\rho[-\alpha +2, -\alpha +2]; T_{I,1}^{\alpha -1}(T_{I,1}^{\alpha}(\pi_{sc})))$
           \item $L(\Delta_\rho[-\alpha +2, -\alpha +1], \Delta_\rho[-\alpha +2, -\alpha +2]; T_{I,1}^{\alpha}(\pi_{sc}))$
           \item $L(\Delta_\rho[-\alpha +2, -\alpha], \Delta_\rho[-\alpha +2, -\alpha+2]; \pi_{sc})$
           \item $L(\Delta_\rho[-\alpha +1, -\alpha +1], \Delta_\rho[-\alpha +2, -\alpha +2], \Delta_\rho[-\alpha+2, -\alpha+2]; T_{I,1}^{\alpha}(\pi_{sc}))$
           \item $L(\Delta_\rho[-\alpha+1 , -\alpha ], \Delta_\rho[-\alpha+2 , -\alpha +2], \Delta_\rho[-\alpha +2, -\alpha +2]; \pi_{sc})$
           \item $L(\Delta_\rho[-\alpha , -\alpha ], \Delta_\rho[-\alpha+2, -\alpha +1], \Delta_\rho[-\alpha+2, -\alpha+2]; \pi_{sc})$
           \item  $L(\Delta_\rho[-\alpha, -\alpha ], \Delta_\rho[-\alpha +1, -\alpha +1], \Delta_\rho[-\alpha +2, \alpha +2], \Delta_\rho[-\alpha+2, -\alpha+2]; \pi_{sc})$
       \end{itemize}
       \item $(\alpha -2, \alpha -1, \alpha -1, \alpha)$: 
       \begin{itemize}
       \item $L(\Delta_\rho[-\alpha+2, -\alpha+1]; T_{I,1}^{\alpha-1}(T_{I,1}^{\alpha}(\pi_{sc})))$  
           \item $L(\Delta_\rho[-\alpha +1, -\alpha +1]; T_{I,1}^{\alpha -2}(T_{I,1}^{\alpha -1}(T_{I,1}^{\alpha}(\pi_{sc}))))$
           \item $L(\Delta_\rho[-\alpha+1, -\alpha+1], \Delta_\rho[-\alpha+2, -\alpha+1]; T_{I,1}^{\alpha}(\pi_{sc}))$  
           \item $L(\Delta_\rho[-\alpha +1, -\alpha +1], \Delta_\rho[-\alpha +2, -\alpha]; \pi_{sc})$
           \item $L(\Delta_\rho[-\alpha +1, -\alpha +1], \Delta_\rho[-\alpha +1, -\alpha +1], \Delta_\rho[-\alpha+2, -\alpha+2]; T_{I,1}^{\alpha}(\pi_{sc}))$
           \item $L(\Delta_\rho[-\alpha+1 , -\alpha ], \Delta_\rho[-\alpha+1 , -\alpha +1], \Delta_\rho[-\alpha+2 , -\alpha +2]; \pi_{sc})$
           \item $L(\Delta_\rho[-\alpha , -\alpha ], \Delta_\rho[-\alpha+1, -\alpha +1], \Delta_\rho[-\alpha+2, -\alpha+1]; \pi_{sc})$
           \item  $L(\Delta_\rho[-\alpha , -\alpha ], \Delta_\rho[-\alpha +1, -\alpha +1], \Delta_\rho[-\alpha +1, \alpha +1], \Delta_\rho[-\alpha+2, -\alpha+2]; \pi_{sc})$
       \end{itemize}
       \item $(\alpha -2, \alpha -1, \alpha , \alpha)$: 
       \begin{itemize}
           \item $L(\Delta_\rho[-\alpha , -\alpha ]; T_{I,1}^{\alpha -2}(T_{I,1}^{\alpha -1}(T_{I,1}^{\alpha}(\pi_{sc}))))$
           \item $L(\Delta_\rho[-\alpha +2, -\alpha]; T_{I,1}^{\alpha}(\pi_{sc}))$
           \item $L(\Delta_\rho[-\alpha +1, -\alpha], \Delta_\rho[-\alpha +2, -\alpha +2]; T_{I,1}^{\alpha}(\pi_{sc}))$
           \item $L(\Delta_\rho[-\alpha, -\alpha], \Delta_\rho[-\alpha+2, -\alpha+2]; T_{I,1}^{\alpha-1}(T_{I,1}^{\alpha}(\pi_{sc})))$ 
           \item $L(\Delta_\rho[-\alpha , -\alpha ], \Delta_\rho[-\alpha +2, -\alpha]; \pi_{sc})$
           \item $L(\Delta_\rho[-\alpha , -\alpha ], \Delta_\rho[-\alpha+1, -\alpha ], \Delta_\rho[-\alpha+2, -\alpha+2]; \pi_{sc})$
           \item $L(\Delta_\rho[-\alpha , -\alpha ], \Delta_\rho[-\alpha, -\alpha ], \Delta_\rho[-\alpha+2, -\alpha+1]; \pi_{sc})$
           \item  $L(\Delta_\rho[-\alpha , -\alpha ], \Delta_\rho[-\alpha , -\alpha ], \Delta_\rho[-\alpha +1, \alpha +1], \Delta_\rho[-\alpha+2, -\alpha+2]; \pi_{sc})$
       \end{itemize}
       \item $(\alpha -2, \alpha -1, \alpha , \alpha+1)$: 
       \begin{itemize}
       \item $L(\Delta_\rho[-\alpha+2, -\alpha-1]; \pi_{sc})$ 
           \item $L(\Delta_\rho[-\alpha -1, -\alpha -1]; T_{I,1}^{\alpha -2}(T_{I,1}^{\alpha -1}(T_{I,1}^{\alpha}(\pi_{sc}))))$
           \item $L(\Delta_\rho[-\alpha+1, -\alpha-1], \Delta_\rho[-\alpha+2, -\alpha+2]; \pi_{sc})$  
           \item $L(\Delta_\rho[-\alpha , -\alpha-1], \Delta_\rho[-\alpha +2, -\alpha+1]; \pi_{sc})$
           \item $L(\Delta_\rho[-\alpha-1, -\alpha-1], \Delta_\rho[-\alpha+2, -\alpha+2]; T_{I,1}^{\alpha-1}(T_{I,1}^{\alpha}(\pi_{sc})))$  
           \item $L(\Delta_\rho[-\alpha-1, -\alpha-1], \Delta_\rho[-\alpha+2, -\alpha+1]; T_{I,1}^{\alpha}(\pi_{sc}))$ 
           \item $L(\Delta_\rho[-\alpha , -\alpha -1], \Delta_\rho[-\alpha+1 , -\alpha+1 ], \Delta_\rho[-\alpha +2, -\alpha+2 ]; \pi_{sc})$
           \item $L(\Delta_\rho[-\alpha-1, -\alpha-1], \Delta_\rho[-\alpha+1, -\alpha+1], \Delta_\rho[-\alpha+2, -\alpha+2]; T_{I,1}^{\alpha}(\pi_{sc}))$  
       \end{itemize}
       \item $(\alpha -1, \alpha -1, \alpha -1, \alpha)$: 
       \begin{itemize}
           \item  $L(\Delta_\rho[-\alpha +1 ,-\alpha +1], \Delta_\rho[-\alpha +1, -\alpha +1]; T_{I,1}^{\alpha -1}(T_{I,1}^{\alpha}(\pi_{sc})))$
           \item $L(\Delta_\rho[-\alpha +1, -\alpha +1], \Delta_\rho[-\alpha +1, -\alpha +1], \Delta_\rho[-\alpha +1, -\alpha +1], T_{I,1}^{\alpha}(\pi_{sc}))$
           \item $L(\Delta_\rho[-\alpha+1 , -\alpha ], \Delta_\rho[-\alpha+1 , -\alpha +1], \Delta_\rho[-\alpha +1, -\alpha +1]; \pi_{sc})$
           \item $L(\Delta_\rho[-\alpha, -\alpha ], \Delta_\rho[-\alpha +1, -\alpha +1], \Delta_\rho[-\alpha +1, -\alpha +1], \Delta_\rho[-\alpha +1, -\alpha +1]; \pi_{sc})$
       \end{itemize}
       \item $(\alpha -1, \alpha -1, \alpha, \alpha)$: 
       \begin{itemize}
       \item $L(\Delta_\rho[-\alpha+1, -\alpha]; T_{I,1}^{\alpha-1}(T_{I,1}^{\alpha}(\pi_{sc})))$ 
           \item $L(\Delta_\rho[-\alpha, -\alpha], \Delta_\rho[-\alpha +1, -\alpha +1]; T_{I,1}^{\alpha -1}(T_{I,1}^{\alpha}(\pi_{sc})))$
           \item $L(\Delta_\rho[-\alpha +1, -\alpha ], \Delta_\rho[-\alpha +1, -\alpha +1]; T_{I,1}^{\alpha}(\pi_{sc}))$
           \item $L(\Delta_\rho[-\alpha +1, -\alpha], \Delta_\rho[-\alpha +1, -\alpha]; \pi_{sc})$
           \item $L(\Delta_\rho[-\alpha , -\alpha ], \Delta_\rho[-\alpha +1, -\alpha +1], \Delta_\rho[-\alpha+1, -\alpha+1]; T_{I,1}^{\alpha}(\pi_{sc}))$
           \item $L(\Delta_\rho[-\alpha , -\alpha ], \Delta_\rho[-\alpha+1, -\alpha ], \Delta_\rho[-\alpha+1, -\alpha+1]; \pi_{sc})$
           \item $L(\Delta_\rho[-\alpha, -\alpha ], \Delta_\rho[-\alpha , -\alpha ], \Delta_\rho[-\alpha +1, -\alpha +1], \Delta_\rho[-\alpha +1, -\alpha +1]; \pi_{sc})$
       \end{itemize}
       \item $(\alpha -1, \alpha -1, \alpha, \alpha +1)$: 
       \begin{itemize}
           \item $L(\Delta_\rho[-\alpha +1, -\alpha +1]; T_{I,1}^{\alpha -1}(T_{I,1}^{\alpha +1}(T_{I,1}^{\alpha}(\pi_{sc}))))$
           \item $L(\Delta_\rho[-\alpha-1, -\alpha-1], \Delta_\rho[-\alpha+1, -\alpha+1]; T_{I,1}^{\alpha-1}(T_{I,1}^{\alpha}(\pi_{sc})))$ 
           \item $L(\Delta_\rho[-\alpha +1, -\alpha +1], \Delta_\rho[-\alpha +1, -\alpha +1]; T_{I,1}^{\alpha +1}(T_{I,1}^{\alpha}(\pi_{sc})))$
           \item $L(\Delta_\rho[-\alpha +1, -\alpha-1], \Delta_\rho[-\alpha +1, -\alpha+1]; \pi_{sc})$
           \item $L(\Delta_\rho[-\alpha -1, -\alpha -1], \Delta_\rho[-\alpha +1, -\alpha +1], \Delta_\rho[-\alpha+1, -\alpha+1]; T_{I,1}^{\alpha}(\pi_{sc}))$
           \item $L(\Delta_\rho[-\alpha-1, -\alpha-1], \Delta_\rho[-\alpha+1, -\alpha], \Delta_\rho[-\alpha+1, -\alpha+1]; \pi_{sc})$ 
           \item $L(\Delta_\rho[-\alpha , -\alpha -1], \Delta_\rho[-\alpha+1 , -\alpha+1 ], \Delta_\rho[-\alpha +1, -\alpha+1 ]; \pi_{sc})$
           \item  $L(\Delta_\rho[-\alpha -1, -\alpha -1], \Delta_\rho[-\alpha , -\alpha ], \Delta_\rho[-\alpha +1, \alpha +1], \Delta_\rho[-\alpha+1, -\alpha+1]; \pi_{sc})$
       \end{itemize}
       \item $(\alpha -1, \alpha, \alpha, \alpha)$: 
       \begin{itemize}
       \item $L(\Delta_\rho[-\alpha, -\alpha], \Delta_\rho[-\alpha+1, -\alpha]; T_{I,1}^{\alpha}(\pi_{sc}))$ 
           \item $L(\Delta_\rho[-\alpha  ,-\alpha ], \Delta_\rho[-\alpha , -\alpha ]; T_{I,1}^{\alpha -1}(T_{I,1}^{\alpha}(\pi_{sc})))$
           \item $L(\Delta_\rho[-\alpha , -\alpha ], \Delta_\rho[-\alpha , -\alpha ], \Delta_\rho[-\alpha+1, -\alpha+1]; T_{I,1}^{\alpha}(\pi_{sc}))$
           \item $L(\Delta_\rho[-\alpha , -\alpha ], \Delta_\rho[-\alpha, -\alpha ], \Delta_\rho[-\alpha+1, -\alpha]; \pi_{sc})$
           \item $L(\Delta_\rho[-\alpha, -\alpha ], \Delta_\rho[-\alpha , -\alpha ], \Delta_\rho[-\alpha , -\alpha ], \Delta_\rho[-\alpha +1, -\alpha +1]; \pi_{sc})$
       \end{itemize}
       \item $(\alpha -1, \alpha , \alpha, \alpha +1)$: 
       \begin{itemize}
       \item $L(\Delta_\rho[-\alpha+1, -\alpha]; T_{I,1}^{\alpha+1}(T_{I,1}^{\alpha}(\pi_{sc})))$  
       \item $L(\Delta_\rho[-\alpha+1, -\alpha-1]; T_{I,1}^{\alpha}(\pi_{sc}))$  
           \item $L(\Delta_\rho[-\alpha , -\alpha ]; T_{I,1}^{\alpha -1}(T_{I,1}^{\alpha +1}(T_{I,1}^{\alpha}(\pi_{sc}))))$
           \item $L(\Delta_\rho[-\alpha, -\alpha-1]; T_{I,1}^{\alpha-1}(T_{I,1}^{\alpha}(\pi_{sc})))$ 
           \item $L(\Delta_\rho[-\alpha-1, -\alpha-1], \Delta_\rho[-\alpha , -\alpha ]; T_{I,1}^{\alpha -1}(T_{I,1}^{\alpha}(\pi_{sc})))$
           \item $L(\Delta_\rho[-\alpha-1, -\alpha-1], \Delta_\rho[-\alpha+1, -\alpha]; T_{I,1}^{\alpha}(\pi_{sc}))$  
           \item $L(\Delta_\rho[-\alpha , -\alpha -1], \Delta_\rho[-\alpha +1, -\alpha +1]; T_{I,1}^{\alpha}(\pi_{sc}))$
           \item $L(\Delta_\rho[-\alpha , -\alpha-1], \Delta_\rho[-\alpha +1, -\alpha]; \pi_{sc})$
           \item $L(\Delta_\rho[-\alpha , -\alpha ], \Delta_\rho[-\alpha +1, -\alpha-1]; \pi_{sc})$
           \item $L(\Delta_\rho[-\alpha , -\alpha -1], \Delta_\rho[-\alpha , -\alpha ], \Delta_\rho[-\alpha+1 , -\alpha+1 ]; \pi_{sc})$
           \item $L(\Delta_\rho[-\alpha -1, -\alpha -1], \Delta_\rho[-\alpha, -\alpha ], \Delta_\rho[-\alpha+1, -\alpha]; \pi_{sc})$
       \end{itemize}
       \item $(\alpha -1, \alpha , \alpha+1, \alpha +1)$: 
       \begin{itemize}
           \item $L(\Delta_\rho[-\alpha -1, -\alpha -1]; T_{I,1}^{\alpha -1}(T_{I,1}^{\alpha +1}(T_{I,1}^{\alpha}(\pi_{sc}))))$
           \item $L(\Delta_\rho[-\alpha -1 ,-\alpha -1], \Delta_\rho[-\alpha -1, -\alpha -1]; T_{I,1}^{\alpha -1}(T_{I,1}^{\alpha}(\pi_{sc})))$
           \item $L(\Delta_\rho[-\alpha -1, -\alpha -1], \Delta_\rho[-\alpha +1, -\alpha-1]; \pi_{sc})$
           \item $L(\Delta_\rho[-\alpha-1, -\alpha-1], \Delta_\rho[-\alpha+1, -\alpha+1], T_{I,1}^{\alpha+1}(T_{I,1}^{\alpha}(\pi_{sc})))$
           \item $L(\Delta_\rho[-\alpha -1, -\alpha -1], \Delta_\rho[-\alpha -1, -\alpha -1], \Delta_\rho[-\alpha+1, -\alpha+1]; T_{I,1}^{\alpha}(\pi_{sc}))$
           \item $L(\Delta_\rho[-\alpha -1, -\alpha -1], \Delta_\rho[-\alpha, -\alpha -1], \Delta_\rho[-\alpha+1, -\alpha+1]; \pi_{sc})$
           \item $L(\Delta_\rho[-\alpha -1, -\alpha -1], \Delta_\rho[-\alpha-1, -\alpha -1], \Delta_\rho[-\alpha+1, -\alpha]; \pi_{sc})$
           \item  $L(\Delta_\rho[-\alpha -1, -\alpha -1], \Delta_\rho[-\alpha -1, -\alpha -1], \Delta_\rho[-\alpha , \alpha ], \Delta_\rho[-\alpha+1, -\alpha+1]; \pi_{sc})$
       \end{itemize}
       \item $(\alpha -1, \alpha , \alpha+1, \alpha +2)$: 
       \begin{itemize}
           \item $L(\Delta_\rho[-\alpha -2, -\alpha -2]; T_{I,1}^{\alpha -1}(T_{I,1}^{\alpha +1}(T_{I,1}^{\alpha}(\pi_{sc}))))$
           \item $L(\Delta_\rho[-\alpha-1, -\alpha-2]; T_{I,1}^{\alpha}(T_{I,1}^{\alpha-1}(\pi_{sc})))$ 
           \item $L(\Delta_\rho[-\alpha+1, -\alpha-2]; \pi_{sc})$  
           \item $L(\Delta_\rho[-\alpha-2, -\alpha-2], \Delta_\rho[-\alpha -1, -\alpha -1]; T_{I,1}^{\alpha -1}(T_{I,1}^{\alpha}(\pi_{sc})))$
           \item $L(\Delta_\rho[-\alpha-2, -\alpha-2], \Delta_\rho[-\alpha+1, -\alpha+1]; T_{I,1}^{\alpha+1}(T_{I,1}^{\alpha}(\pi_{sc})))$ 
           \item $L(\Delta_\rho[-\alpha-1, -\alpha-2], \Delta_\rho[-\alpha+1, -\alpha+1]; T_{I,1}^{\alpha}(\pi_{sc}))$  
           \item $L(\Delta_\rho[-\alpha-1, -\alpha-2], \Delta_\rho[-\alpha+1, -\alpha]; \pi_{sc})$  
             \item $L(\Delta_\rho[-\alpha , -\alpha-2], \Delta_\rho[-\alpha+1 , -\alpha+1]; \pi_{sc})$
             \item $L(\Delta_\rho[-\alpha -2, -\alpha -2], \Delta_\rho[-\alpha +1, -\alpha-1]; \pi_{sc})$
             \item $L(\Delta_\rho[-\alpha -1, -\alpha -2], \Delta_\rho[-\alpha , -\alpha ], \Delta_\rho[-\alpha+1, -\alpha+1 ]; \pi_{sc})$
             \item $L(\Delta_\rho[-\alpha-2, -\alpha-2], \Delta_\rho[-\alpha-1, -\alpha-1], \Delta_\rho[-\alpha+1, -\alpha+1]; T_{I,1}^{\alpha}(\pi_{sc}))$  
             \item $L(\Delta_\rho[-\alpha -2, -\alpha -2], \Delta_\rho[-\alpha, -\alpha -1], \Delta_\rho[-\alpha+1, -\alpha+1]; \pi_{sc})$
       \end{itemize}
       \item $(\alpha, \alpha, \alpha, \alpha)$: 
       \begin{itemize}
           \item $L(\Delta_\rho[-\alpha , -\alpha ], \Delta_\rho[-\alpha , -\alpha ], \Delta_\rho[-\alpha , -\alpha ], T_{I,1}^{\alpha}(\pi_{sc}))$
           \item $L(\Delta_\rho[-\alpha, -\alpha], \Delta_\rho[-\alpha , -\alpha ], \Delta_\rho[-\alpha , -\alpha ], \Delta_\rho[-\alpha , -\alpha ]; \pi_{sc})$
       \end{itemize}
       \item $(\alpha, \alpha, \alpha, \alpha +1)$: 
       \begin{itemize}
           \item  $L(\Delta_\rho[-\alpha , -\alpha ], \Delta_\rho[-\alpha , -\alpha ]; T_{I,1}^{\alpha +1}(T_{I,1}^{\alpha}(\pi_{sc})))$
           \item $L(\Delta_\rho[-\alpha , -\alpha -1], \Delta_\rho[-\alpha , -\alpha ]; T_{I,1}^{\alpha}(\pi_{sc}))$
           \item $L(\Delta_\rho[-\alpha -1, -\alpha -1], \Delta_\rho[-\alpha , -\alpha ], \Delta_\rho[-\alpha, -\alpha]; T_{I,1}^{\alpha}(\pi_{sc}))$
           \item $L(\Delta_\rho[-\alpha , -\alpha -1], \Delta_\rho[-\alpha , -\alpha ], \Delta_\rho[-\alpha , -\alpha ];\pi_{sc})$
           \item $L(\Delta_\rho[-\alpha-1, -\alpha-1], \Delta_\rho[-\alpha , -\alpha ], \Delta_\rho[-\alpha , -\alpha ], \Delta_\rho[-\alpha , -\alpha ]; \pi_{sc})$
       \end{itemize}
       \item $(\alpha, \alpha, \alpha+1, \alpha +1)$:
       \begin{itemize}
       \item $L(\Delta_\rho[-\alpha, -\alpha-1]; T_{I,1}^{\alpha+1}(T_{I,1}^{\alpha}(\pi_{sc})))$  
           \item $L(\Delta_\rho[-\alpha -1, -\alpha -1], \Delta_\rho[-\alpha, -\alpha]; T_{I,1}^{\alpha +1}(T_{I,1}^{\alpha}(\pi_{sc})))$
           \item $L(\Delta_\rho[-\alpha-1, -\alpha-1], \Delta_\rho[-\alpha, -\alpha-1]; T_{I,1}^{\alpha}(\pi_{sc}))$ 
           \item $L(\Delta_\rho[-\alpha , -\alpha-1], \Delta_\rho[-\alpha , -\alpha-1]; \pi_{sc})$
           \item $L(\Delta_\rho[-\alpha -1, -\alpha -1], \Delta_\rho[-\alpha -1, -\alpha -1], \Delta_\rho[-\alpha, -\alpha]; T_{I,1}^{\alpha}(\pi_{sc}))$
           \item $L(\Delta_\rho[-\alpha -1, -\alpha -1], \Delta_\rho[-\alpha, -\alpha -1], \Delta_\rho[-\alpha, -\alpha]; \pi_{sc})$
           \item $L(\Delta_\rho[-\alpha-1, -\alpha-1], \Delta_\rho[-\alpha-1 , -\alpha -1], \Delta_\rho[-\alpha , -\alpha ], \Delta_\rho[-\alpha , -\alpha ]; \pi_{sc})$
       \end{itemize}
           \item $(\alpha, \alpha, \alpha +1, \alpha +2)$: 
           \begin{itemize}
               \item $L(\Delta_\rho[-\alpha, -\alpha]; T_{I,1}^{\alpha +2}(T_{I,1}^{\alpha +1}(T_{I,1}^{\alpha}(\pi_{sc}))))$
               \item $L(\Delta_\rho[-\alpha , -\alpha-2]; T_{I,1}^{\alpha}(\pi_{sc}))$
               \item $L(\Delta_\rho[-\alpha -1, -\alpha -2], \Delta_\rho[-\alpha , -\alpha ]; T_{I,1}^{\alpha}(\pi_{sc}))$
                \item $L(\Delta_\rho[-\alpha , -\alpha-2], \Delta_\rho[-\alpha , -\alpha]; \pi_{sc})$
                \item $L(\Delta_\rho[-\alpha-2, -\alpha-2], \Delta_\rho[-\alpha, -\alpha]; T_{I,1}^{\alpha+1}(T_{I,1}^{\alpha}(\pi_{sc})))$  
                \item $L(\Delta_\rho[-\alpha-2, -\alpha-2], \Delta_\rho[-\alpha, -\alpha-1]; T_{I,1}^{\alpha}(\pi_{sc}))$  
                \item $L(\Delta_\rho[-\alpha -2, -\alpha -2], \Delta_\rho[-\alpha -1, -\alpha -1], \Delta_\rho[-\alpha, -\alpha]; T_{I,1}^{\alpha}(\pi_{sc}))$
                \item $L(\Delta_\rho[-\alpha-1 , -\alpha -2], \Delta_\rho[-\alpha , -\alpha ], \Delta_\rho[-\alpha , -\alpha ], \pi_{sc})$
                \item $L(\Delta_\rho[-\alpha -2, -\alpha -2], \Delta_\rho[-\alpha, -\alpha -1], \Delta_\rho[-\alpha, -\alpha]; \pi_{sc})$
                \item  $L(\Delta_\rho[-\alpha -2, -\alpha -2], \Delta_\rho[-\alpha -1, -\alpha -1], \Delta_\rho[-\alpha, -\alpha ], \Delta_\rho[-\alpha, -\alpha]; \pi_{sc})$
           \end{itemize}
           \item $(\alpha, \alpha +1, \alpha +1, \alpha +1)$:
           \begin{itemize}
               \item  $L(\Delta_\rho[-\alpha -1, -\alpha -1], \Delta_\rho[-\alpha -1, -\alpha -1]; T_{I,1}^{\alpha +1}(T_{I,1}^{\alpha}(\pi_{sc})))$
               \item $L(\Delta_\rho[-\alpha -1, -\alpha -1], \Delta_\rho[-\alpha -1, -\alpha -1], \Delta_\rho[-\alpha -1, -\alpha -1]; T_{I,1}^{\alpha}(\pi_{sc}))$
               \item $L(\Delta_\rho[-\alpha -1, -\alpha -1], \Delta_\rho[-\alpha-1, -\alpha -1], \Delta_\rho[-\alpha, -\alpha-1]; \pi_{sc})$
               \item $L(\Delta_\rho[-\alpha-1, -\alpha-1], \Delta_\rho[-\alpha-1 , -\alpha -1], \Delta_\rho[-\alpha-1 , -\alpha-1 ], \Delta_\rho[-\alpha , -\alpha ]; \pi_{sc})$
           \end{itemize}
           \item $(\alpha, \alpha+1, \alpha +1, \alpha +2)$: 
           \begin{itemize}
               \item $L(\Delta_\rho[-\alpha -1, -\alpha - 1]; T_{I,1}^{\alpha +2}(T_{I,1}^{\alpha +1}(T_{I,1}^{\alpha}(\pi_{sc}))))$
               \item $L(\Delta_\rho[-\alpha-1, -\alpha-2]; T_{I,1}^{\alpha+1}(T_{I,1}^{\alpha}(\pi_{sc})))$ 
               \item $L(\Delta_\rho[-\alpha -1, -\alpha -2], \Delta_\rho[-\alpha -1, -\alpha -1]; T_{I,1}^{\alpha}(\pi_{sc}))$
               \item $L(\Delta_\rho[-\alpha -1, -\alpha-2], \Delta_\rho[-\alpha , -\alpha-1]; \pi_{sc})$
                \item $L(\Delta_\rho[-\alpha -1, -\alpha -1], \Delta_\rho[-\alpha , -\alpha-2]; \pi_{sc})$
                
                 \item $L(\Delta_\rho[-\alpha -2, -\alpha -2], \Delta_\rho[-\alpha-1, -\alpha-1]; T_{I,1}^{\alpha +1}(T_{I,1}^{\alpha}(\pi_{sc})))$
                 \item $L(\Delta_\rho[-\alpha-1, -\alpha-2], \Delta_\rho[-\alpha-1, -\alpha-1], \Delta_\rho[-\alpha, -\alpha]; \pi_{sc})$
                 \item $L(\Delta_\rho[-\alpha -2, -\alpha -2], \Delta_\rho[-\alpha-1, -\alpha -1], \Delta_\rho[-\alpha, -\alpha-1]; \pi_{sc})$
                \item $L(\Delta_\rho[-\alpha -2, -\alpha -2], \Delta_\rho[-\alpha -1, -\alpha -1], \Delta_\rho[-\alpha-1, -\alpha-1]; T_{I,1}^{\alpha}(\pi_{sc}))$
                 \item  $L(\Delta_\rho[-\alpha -2, -\alpha -2], \Delta_\rho[-\alpha -1, -\alpha -1], \Delta_\rho[-\alpha -1, \alpha -1], \Delta_\rho[-\alpha, -\alpha]; \pi_{sc})$
           \end{itemize}
           \item $(\alpha, \alpha+1, \alpha +2, \alpha +2)$: 
           \begin{itemize}
               \item $L(\Delta_\rho[-\alpha,-2 -\alpha-2]; T_{I,1}^{\alpha +2}(T_{I,1}^{\alpha +1}(T_{I,1}^{\alpha}(\pi_{sc}))))$
               \item $L(\Delta_\rho[-\alpha-2, -\alpha -2], \Delta_\rho[-\alpha -2, -\alpha -2]; T_{I,1}^{\alpha +1}(T_{I,1}^{\alpha}(\pi_{sc})))$
               \item $L(\Delta_\rho[-\alpha -2, -\alpha-2], \Delta_\rho[-\alpha-1, -\alpha-2]; T_{I,1}^{\alpha}(\pi_{sc}))$ 
               \item $L(\Delta_\rho[-\alpha -2, -\alpha -2], \Delta_\rho[-\alpha , -\alpha-2]; \pi_{sc})$
               \item $L(\Delta_\rho[-\alpha -2, -\alpha -2], \Delta_\rho[-\alpha -2, -\alpha -2], \Delta_\rho[-\alpha-1, -\alpha-1]; T_{I,1}^{\alpha}(\pi_{sc}))$
               \item $L(\Delta_\rho[-\alpha -2, -\alpha -2], \Delta_\rho[-\alpha-1, -\alpha -2], \Delta_\rho[-\alpha, -\alpha]; \pi_{sc})$
                \item $L(\Delta_\rho[-\alpha -2, -\alpha -2], \Delta_\rho[-\alpha-2, -\alpha -2], \Delta_\rho[-\alpha, -\alpha-1]; \pi_{sc})$
                \item  $L(\Delta_\rho[-\alpha -2, -\alpha -2], \Delta_\rho[-\alpha -2, -\alpha -2], \Delta_\rho[-\alpha -1, \alpha -1], \Delta_\rho[-\alpha, -\alpha]; \pi_{sc})$
           \end{itemize}
           \item $(\alpha, \alpha+1, \alpha +2, \alpha +3)$: 
           \begin{itemize}
               \item $L(\Delta_\rho[-\alpha-3, -\alpha-3]; T_{I,1}^{\alpha +2}(T_{I,1}^{\alpha +1}(T_{I,1}^{\alpha}(\pi_{sc}))))$
               \item $L(\Delta_\rho[-\alpha-2, -\alpha-3]; T_{I,1}^{\alpha+1}(T_{I,1}^{\alpha}(\pi_{sc})))$
               \item $L(\Delta_\rho[-\alpha -1, -\alpha -3]; T_{I,1}^{\alpha}(\pi_{sc}))$
               \item $L(\Delta_\rho[-\alpha -3, -\alpha -3], \Delta_\rho[-\alpha-2, -\alpha-2]; T_{I,1}^{\alpha +1}(T_{I,1}^{\alpha}(\pi_{sc})))$
               \item $L(\Delta_\rho[-\alpha, -\alpha-3]; \pi_{sc})$ 
               \item $L(\Delta_\rho[-\alpha-3, -\alpha-3], \Delta_\rho[-\alpha -1, -\alpha -2]; T_{I,1}^{\alpha}(\pi_{sc}))$
               \item $L(\Delta_\rho[-\alpha -2, -\alpha -3], \Delta_\rho[-\alpha -1, -\alpha -1]; T_{I,1}^{\alpha}(\pi_{sc}))$
               \item $L(\Delta_\rho[-\alpha -2, -\alpha-3], \Delta_\rho[-\alpha , -\alpha-1]; \pi_{sc})$
                \item $L(\Delta_\rho[-\alpha -1, -\alpha-3], \Delta_\rho[-\alpha , -\alpha]; \pi_{sc})$
                \item $L(\Delta_\rho[-\alpha -3, -\alpha -3], \Delta_\rho[-\alpha , -\alpha-2]; \pi_{sc})$
                \item $L(\Delta_\rho[-\alpha -3, -\alpha -3], \Delta_\rho[-\alpha -2, -\alpha -2], \Delta_\rho[-\alpha-1, -\alpha-1]; T_{I,1}^{\alpha}(\pi_{sc}))$
                \item $L(\Delta_\rho[-\alpha-2 , -\alpha -3], \Delta_\rho[-\alpha -1, -\alpha-1 ], \Delta_\rho[-\alpha , -\alpha ];\pi_{sc})$
                \item $L(\Delta_\rho[-\alpha -3, -\alpha-3], \Delta_\rho[-\alpha -1, -\alpha -2], \Delta_\rho[-\alpha, -\alpha]; \pi_{sc})$
                \item $L(\Delta_\rho[-\alpha -3, -\alpha -3], \Delta_\rho[-\alpha-2, -\alpha -2], \Delta_\rho[-\alpha, -\alpha-1]; \pi_{sc})$
           \end{itemize}
       \end{enumerate}
    \end{enumerate}
\end{prop}

\begin{proof} The list in this proposition is also directly implied by Propositions \ref{nontempA1} to \ref{nontemp4}. 
\end{proof}

\end{document}